\documentclass[10pt,reqno]{article}
\usepackage{amsfonts,latexsym,amsmath, amssymb, mathrsfs, verbatim}
\usepackage[nobysame]{amsrefs}
\usepackage{url,color}
\usepackage{pifont}
\usepackage{upgreek}
\usepackage{subfigure}
\usepackage{fancyhdr}
\usepackage[linktocpage]{hyperref}
\usepackage{calligra}
\usepackage{marvosym}
\usepackage[percent]{overpic}
\usepackage{pict2e}
\usepackage{tensor}
\usepackage{amsthm}
\usepackage[hypcap=false]{caption}
\usepackage{xcolor}
\usepackage{wasysym}
\usepackage{accents}
\usepackage{enumitem}
\usepackage{tikz}

\usetikzlibrary{shapes.geometric, arrows}
\usetikzlibrary{mindmap}
\setlist[enumerate]{label*=\arabic*.}
\numberwithin{equation}{section}

\topmargin - .23 in

\newtheorem{theorem}{Theorem}[section]
\newtheorem{proposition}{Proposition}[section]
\newtheorem{lemma}[proposition]{Lemma}
\newtheorem{corollary}[proposition]{Corollary}

\theoremstyle{definition}
\newtheorem{definition}{Definition}[section]
\newtheorem{remark}[proposition]{Remark}

\DeclareMathAlphabet{\mathcalligra}{T1}{calligra}{m}{n}
\DeclareFontShape{T1}{calligra}{m}{n}{<->s*[2.2]callig15}{}

\newcommand{\ent}{s}
\newcommand{\vars}{\dot{\ent}}
\newcommand{\varp}{\dot{\pressure}}

\newcommand{\pressure}{p}

\newcommand{\varv}{\dot{v}}
\newcommand{\density}{\varrho}

\newcommand{\enth}{H}
\newcommand{\lnenth}{h}
\newcommand{\varh}{\dot{\lnenth}}
\newcommand{\temp}{\theta}

\newcommand{\vort}{\upomega}

\newcommand{\speed}{c}

\newcommand{\variables}{\Psi}

\newcommand{\uvariable}{\underline{\variables}}

\newcommand{\norm}[1]{\left\lVert#1\right\rVert}
\newcommand{\onenorm}[4]{\norm{#1}_{L_{#2}^{#3}(#4)}}
\newcommand{\twonorms}[6]{\norm{#1}_{L_{#2}^{#3}L_{#4}^{#5}(#6)}}

\newcommand{\threenorms}[8]{\norm{#1}_{L_{#2}^{#3}L_{#4}^{#5}L_{#6}^{#7}(#8)}}

\newcommand{\abs}[1]{\left\vert#1\right\vert}

\newcommand{\sgabs}[1]{\left\vert#1\right\vert_{\gsphere}}

\newcommand{\curl}{\mbox{\upshape curl}}
\newcommand{\antisymmetic}{\epsilon}
\newcommand{\tensorproduct}{\widehat{\otimes}}

\newcommand{\sphere}{\theta}

\newcommand{\region}{\mathcal{M}}
\newcommand{\hypersur}{\mathcal{H}}
\newcommand{\hypersurbar}{\underline{\mathcal{H}}}

\newcommand{\inihypersur}{\mathring{\hypersur}}
\newcommand{\inihypersurbar}{\mathring{\hypersurbar}}

\newcommand{\suubarprime}{S_{u,\ubarprime}}
\newcommand{\suubar}{S_{u,\ubar}}
\newcommand{\suprimeubar}{S_{\uprime,\ubar}}
\newcommand{\Stwo}{\mathbb{S}^2}

\newcommand{\gensmoothfunction}{\mathrm{f}}

\newcommand{\Ricfour}[2]{\mathbf{Ric}_{#1 #2}}
\newcommand{\angRicfour}{\mathbf{Ric} \mkern-20mu / \mkern 11mu}
\newcommand{\hangRicfour}{\widehat{\mathbf{Ric}} \mkern-20mu / \mkern 11mu}
\newcommand{\Riemfour}[4]{\mathbf{R}_{#1 #2 #3 #4}}
\newcommand{\Riem}{\mathbf{R}}
\newcommand{\uRiem}{\underline{\Riem}}
\newcommand{\upRiem}{\overline{\Riem}}

\newcommand{\EMT}{T}
\newcommand{\Scalar}{\mathbf{R}}
\newcommand{\Weyl}{\mathbf{W}}
\newcommand{\Weyldual}{\mkern-2mu ^{\star}\mathbf{W}}
\newcommand{\AlphaAB}{\alpha}
\newcommand{\CAlphaAB}{\check{\alpha}}
\newcommand{\FAlphaAB}{\AlphaAB^{F}}
\newcommand{\HodgeLeftDual}{\mkern 0mu ^{\star}}
\newcommand{\Alphabar}{\underline{\AlphaAB}}
\newcommand{\CAlphabar}{\check{\underline{\AlphaAB}}}
\newcommand{\FAlphabar}{\Alphabar^{F}}
\newcommand{\BetaA}{\beta}
\newcommand{\CBetaA}{\check{\beta}}
\newcommand{\FBetaA}{\BetaA^{F}}
\newcommand{\Betabar}{\underline{\BetaA}}
\newcommand{\CBetabar}{\check{\underline{\BetaA}}}
\newcommand{\FBetabar}{\Betabar^{F}}
\newcommand{\Rho}{\rho}
\newcommand{\CRho}{\check{\rho}}
\newcommand{\FRho}{\Rho^{F}}
\newcommand{\Rhodual}{\sigma}
\newcommand{\FRhodual}{\Rhodual^{F}}
\newcommand{\CRhodual}{\check{\Rhodual}}
\newcommand{\DWeyl}{\mathbf{J}}
\newcommand{\DWeylStar}{\mathbf{J}^{\star}}
\newcommand{\DWeylBianchi}{\widetilde{\mathbf{J}}}
\newcommand{\DWeylBianchiStar}{\widetilde{\mathbf{J}}^{\star}}
\newcommand{\GaussCurvature}{\mathbf{K}}
\newcommand{\spherenormal}{\mathbf{N}}

\newcommand{\lapse}{\Omega}

\newcommand{\ueta}{\underline{\upeta}}

\newcommand{\hchi}{\hat{\upchi}}
\newcommand{\uchi}{\underline{\upchi}}
\newcommand{\huchi}{\hat{\uchi}}

\newcommand{\upsi}{\underline{\psi}}
\newcommand{\uPsi}{\underline{\Psi}}
\newcommand{\unpsi}{\overline{\psi}}

\newcommand{\BAsize}{\Lambda}
\newcommand{\BAfluid}{\Delta}
\newcommand{\iniBAsize}{\mathring{\BAsize}}
\newcommand{\iniBAfluid}{\mathring{\BAfluid}}
\newcommand{\varJ}{\dot{J}}
\newcommand{\Jen}[1]{^{(#1)} \mkern-3mu \mathbf{J}}

\newcommand{\deform}[1]{{^{(#1)} \mkern-1mu \pmb{\pi}}}

\newcommand{\Fluid}[2]{{^{(#1)} \mkern-1mu \mathfrak{#2}}}

\newcommand{\gacoustical}{\gfour_{Acou}}
\newcommand{\gfour}{\mathbf{g}}

\newcommand{\gsphere}{g \mkern-8.5mu / }

\newcommand{\Chfour}{\pmb{\Gamma}}

\newcommand{\volume}{\upsilon}

\newcommand{\gvol}{\varpi_{\gfour}}

\newcommand{\spherevol}{\varpi_{\gsphere}}

\newcommand{\sphereproject}{{\Pi \mkern-12mu / } \, }

\newcommand{\vproject}{\mathbf{\Pi}}

\newcommand{\mytr}{{\mbox{\upshape tr}}}

\newcommand{\gtr}{\mytr_{\gsphere}}

\newcommand{\ms}{\hspace{0.25cm}}

\newcommand{\Timelike}{\mathbf{T}}

\newcommand{\Lunit}{L}

\newcommand{\uLunit}{\underline{L}}

\newcommand{\Lgeo}{\Lunit_{(Geo)}}
\newcommand{\uLgeo}{\uLunit_{(Geo)}}
\newcommand{\uomega}{\underline{\omega}}

\newcommand{\p}{\partial}
\newcommand{\pfour}{\pmb\partial}

\newcommand{\ubar}{\underline{u}}
\newcommand{\ubarprime}{\underline{u}^\prime}
\newcommand{\uprime}{u^\prime}

\newcommand{\uts}{u_{\textit{T}}}
\newcommand{\uinf}{u_{\infty}}

\newcommand{\Dfour}{\mathbf{D}}
\newcommand{\WDfour}{\tilde{\Dfour}}
\newcommand{\angv}{{v \mkern-10mu / }}
\newcommand{\angD}{ {\Dfour \mkern-14mu / \,}}
\newcommand{\WangD}{ {\tilde{\Dfour} \mkern-14mu / \,}}

\newcommand{\angV}{ {V \mkern-14mu / }}

\newcommand{\angnabla}{ {\nabla \mkern-14mu / \,} }

\newcommand{\angdiv}{\mbox{\upshape{div} $\mkern-17mu /$\,}}
\newcommand{\angcurl}{\mbox{\upshape{curl} $\mkern-17mu /$\,}}

\newcommand{\rgeo}{\widetilde{r}}

\newcommand{\mass}{\upmu}
\newcommand{\umass}{\underline{\upmu}}

\newcommand{\Lie}{\mathcal{L}}

\newcommand{\angLie}{ { \mathcal{L} \mkern-10mu / } }

\newcommand{\diff}{\mathrm d}

\newcommand{\average}[1]{\overline{#1}}

\newcommand{\Energy}{\mathbb{F}}

\newcommand{\uEnergy}{\underline{\mathbb{F}}}

\title{Formation of Trapped Surfaces for the Einstein-Euler System}

\date{\today}      
\begin{document}
	\author{Xinliang An\thanks{National University of Singapore, Singapore. matax@nus.edu.sg} and 
		Sifan Yu\thanks{National University of Singapore, Singapore. sifanyu@nus.edu.sg}}
	\maketitle
	\begin{abstract}
		
		We prove the formation of trapped surfaces for the Einstein--Euler system in $3+1$ dimensions beyond spherical symmetry. We consider smooth characteristic initial data in a short-pulse regime with rapid fluid injection. Starting from initial data containing no trapped surfaces, we obtain two sufficient criteria for trapped surface formation, corresponding to the shear and fluid contributions, respectively. In particular, for the fluid-driven criterion, we impose no positive lower bound on the shear flux. To handle the multiple propagation speeds associated with gravitational waves, sound waves, and fluid transport, our analysis combines three ingredients: a coercive fluid energy current whose null fluxes reflect the unequal velocity components, curvature renormalizations adapted to the matter terms, and weighted geometric estimates controlling derivatives in all spacetime directions. To ensure that fluid injection alone can lead to trapped surface formation, we use particle-number conservation in a designed spacetime region to control the focusing process and obtain the desired sphere-integrated fluxes.

	\end{abstract}
	\tableofcontents
	\newpage
	
	\def\Ric{\mathrm{Ric}}\def\bfg{\mathbf{g}}\def\MM{\mathcal{M}}\def\Cref{\ref}\def\cref{\ref}

	\def\ub{\underline{u}}\def\Cb{\underline{C}}\def\Lb{\underline{L}}\def\Lh{\hat{L}}\def\Lbh{\hat{\Lb}}\def\phib{\underline{\phi}}\def\Phib{\underline{\Phi}}\def\Db{\underline{D}}\def\Dh{\hat{D}}\def\Dbh{\hat{\Db}}\def\omb{\underline{\omega}}\def\omh{\hat{\omega}}\def\ombh{\hat{\omb}}\def\Pb{\underline{P}}\def\chib{\underline{\chi}}\def\chih{\hat{\chi}}\def\chibh{\hat{\chib}}

	\def\alb{\underline{\alpha}}\def\zeb{\underline{\zeta}}\def\beb{\underline{\beta}}\def\etb{\underline{\eta}}\def\Mb{\underline{M}}\def\oth{\hat{\otimes}}
	
	\gdef\dxd{\displaystyle}
	\gdef\x{\xi}
	\def\al{\aligned}
	\def\eal{\endaligned}
	\def\be{\begin{equation}}
		\def\ee{\end{equation}}

	\def\a {\alpha}
	\def\b {\beta}
	\def\ab {\alphab}
	\def\bb {\betab}
	\def\nab {\nabla}

	\def\ub {\underline{u}}
	\def\th {\theta}
	\def\Lb {\underline{L}}
	\def\Hb {\underline{H}}
	\def\chib {\underline{\chi}}
	\def\chih {\hat{\chi}}
	\def\chibh {\hat{\underline{\chi}}}
	\def\omegab {\underline{\omega}}
	\def\etab {\underline{\eta}}
	\def\betab {\underline{\beta}}
	\def\alphab {\underline{\alpha}}
	\def\Psib {\underline{\Psi}}
	\def\hot{\widehat{\otimes}}
	\def\Phib {\underline{\Phi}}
	\def\thb {\underline{\theta}}
	\def\t {\tilde}
	\def\st {\tilde{s}}
	\def\h {\hat}

	\def\rhoc{\check{\rho}}
	\def\sigmac{\check{\sigma}}
	\def\Psic{\check{\Psi}}
	\def\kappab{\underline{\kappa}}
	\def\betabc {\check{\underline{\beta}}}
	
	\def\d {\delta}
	\def\f {\frac}
	\def\l {\bigg(}
	\def\r {\bigg)}
	\def\S {S_{u,\underline{u}}}
	\def\o{\omega}
	\def\O{\Omega}\
	\def\be{\begin{equation}\begin{split}}
			\def\en{\end{split}\end{equation}}
	\def\at{a^{\frac{1}{2}}}
	\def\af{a^{\frac{1}{4}}}
	\def\od{\omega^{\dagger}}
	\def\ombd{\underline{\omega}^{\dagger}}
	\def\K{K-\frac{1}{|u|^2}}
	\def\ut{\frac{1}{|u|^2}}
	\def\s{\frac{\delta a^{\frac{1}{2}}}{|u|}}
	\def\Kb{K-\frac{1}{(u+\underline{u})^2}}
	\def\ut{\frac{1}{|u|^2}}
	\def\s{\frac{\delta a^{\frac{1}{2}}}{|u|}}
	\def\Kb{K-\frac{1}{(u+\underline{u})^2}}
	\def\Kt{K-\frac{1}{|u|^2}-\frac14 \nab^A\phi\nab_A\phi} 
	\def\bf{b^{\frac{1}{4}}}
	\def\bt{b^{\frac{1}{2}}}
	\def\lo{\lambda_1}
	\def\lt{\lambda_2}
	\def\phib{\bar{\phi}}
	\def\bR{\bar{R}}
	\def\tR{\tilde{R}}
	\def\S{S_{u,\ub}}
	
	\def\de{\delta}
	\def\ls{\lesssim}
	\def\om{\omega}
	\def\Om{\Omega}
	\def\O{\Omega} 
	\def\nab{\nabla}
	\def\tp{\widetilde{\phi}}
	\def\tO{\tilde{O}}
	\def\cp{\check{\phi}}
	\def\ms{\mu^{*}}\def\F{\big(1+\frac{\ub^{\frac12}a^{\frac14}}{|u|^{\frac12}\O} \big)}
	\def\ho{\hat{\omega}}
	
	\def\e{\epsilon}\def\et{\frac{\epsilon}{2}}\def\ef{\frac{\epsilon}{4}}\def\LH{L^2(H_u)}\def\LHb{L^2(\underline{H}_{\underline{u}})}\def\M{\mathcal}\def\TM{\tilde{\mathcal}}\def\q{\underline{\psi}\hspace{1pt}}\def\Li{_{L^{\infty}(S_{u,\underline{u}})}}\def\Lt{_{L^{2}(S)}}\def\da{\delta^{-\frac{\epsilon}{2}}}\def\db{\delta^{1-\frac{\epsilon}{2}}}\def\D{\Delta}


	\renewcommand{\div}{\mbox{div }}
	\def\curl{\mbox{curl }}\def\trchb{\mbox{tr} \chib}\def\trch{\mbox{tr}\chi}
	\def\tr{\mbox{tr}}
	\def\Ls{{\mathcal L} \mkern-10mu /\,}\def\eps{{\epsilon} \mkern-8mu /\,}

	\def\xib{\underline{\xi}}\def\psib{\underline{\psi}}\def\rhob{\underline{\rho}}\def\thetab{\underline{\theta}}\def\gammab{\underline{\gamma}}\def\nub{\underline{\nu}}\def\lb{\underline{l}}\def\mub{\underline{\mu}}\def\Xib{\underline{\Xi}}\def\Thetab{\underline{\Theta}}\def\Lambdab{\underline{\Lambda}}\def\vphb{\underline{\varphi}}
	
	\def\ih{\hat{i}}
	
	\def\tcL{\widetilde{\mathscr{L}}}
	
	\def\sRic{Ric\mkern-19mu /\,\,\,\,}\def\sL{{\cal L}\mkern-10mu /}\def\sLh{\hat{\sL}}\def\sg{g\mkern-9mu /}\def\seps{\epsilon\mkern-8mu /}\def\sd{d\mkern-10mu /}\def\sR{R\mkern-10mu /}\def\snab{\nabla\mkern-13mu /}\def\sdiv{\mbox{div}\mkern-19mu /\,\,\,\,}\def\scurl{\mbox{curl}\mkern-19mu /\,\,\,\,}\def\slap{\mbox{$\triangle  \mkern-13mu / \,$}}\def\sGamma{\Gamma\mkern-10mu /}\def\somega{\omega\mkern-10mu /}\def\somb{\omb\mkern-10mu /}\def\spi{\pi\mkern-10mu /}\def\sJ{J\mkern-10mu /}\renewcommand{\sp}{p\mkern-9mu /}
	\def\su{u\mkern-8mu /}

	\section{Introduction}
	\label{S:Introduction}

	Let $(\region,\gfour)$ be a four-dimensional Lorentzian spacetime, with covariant derivative $\Dfour$, Ricci tensor $\Ricfour{\mu}{\nu}$, and scalar curvature $\Scalar$. We study the following Einstein--Euler equations
\begin{subequations}\label{1.0}
	\begin{align}
		\label{1.0.1}\Ricfour{\mu}{\nu}-\frac12\Scalar\gfour_{\mu\nu}&=\EMT_{\mu\nu},
		&\EMT_{\mu\nu}&:=(\density+\pressure)v_\mu v_\nu+\pressure\gfour_{\mu\nu},\\
		\label{1.0.2}\Dfour_\mu\EMT^{\mu\nu}&=0,\\
		\label{1.0.3}\Dfour_\mu I^\mu&=0,&
		I^\mu&:=nv^\mu.
	\end{align}
\end{subequations}
Here $v$ is the future-directed fluid velocity, normalized by $\gfour(v,v)=-1$; $\density$, $\pressure$, and $n$ are the proper energy density, pressure, and particle number density. Greek indices are spacetime indices. Writing $\ent$ for the entropy per particle and $\enth=(\density+\pressure)/n$ for the enthalpy per particle, the equation of state has the form $\pressure=\pressure(\enth,\ent)$. The sound speed satisfies $\speed^2:=\left.(\partial\pressure/\partial\density)\right|_{\ent}\leq\speed_0^2<1$, as in Section~\ref{SS:deffluid}. Moreover, we assume that the fluid is away from vacuum\footnote{It is well-known that the fluid part of the equation is degenerate if density $\density=0$. In that case, one should consider a free-boundary problem with degeneracy at the boundary. Little is understood in this subject; see \cite{DisconziSpeck2027} for a local well-posedness result in spherical symmetry.}; that is, we assume that there is a constant background enthalpy $\average{\enth}$ and we define the logarithmic enthalpy $\lnenth:=\ln(\enth/\average{\enth})$.

		A central problem in gravitational collapse is to understand how a closed trapped surface can form from initial data that contain no such surface. A trapped surface is a closed spacelike two-surface whose two future null expansions are everywhere negative: its area decreases along both families of future-directed orthogonal light rays. Under the null convergence condition and the global hypothesis of a noncompact Cauchy hypersurface, Penrose's theorem \cite{Penrose1965} connects the existence of such a surface to future null geodesic incompleteness. The evolutionary problem is to produce this geometric condition through the Einstein equations while controlling the solution up to the time of trapping.
		
		In this paper, we study this problem for the Einstein--Euler system in $3+1$ dimensions, beyond spherical symmetry. Our main theorem gives two sufficient conditions for trapped-surface formation from characteristic initial data. In the first, the required focusing is supplied by the outgoing shear. In the second, it is supplied by the perfect-fluid contribution to the outgoing null Raychaudhuri equation. The fluid alternative assumes no corresponding positive lower bound on the shear flux; the shear remains part of the coupled evolution. To the authors' knowledge, this is the first result on trapped surface formation beyond spherical symmetry driven by two distinct mechanisms in a system with multiple characteristic speeds. The argument combines estimates for a rapidly injected fluid with a propagation mechanism for its focusing flux.

	\subsection{Background and the Einstein--Euler Problem}

		Christodoulou \cite{Christodoulou2009} established the dynamical formation of trapped surfaces in vacuum without symmetry assumptions. His short-pulse construction concentrates incoming gravitational radiation in a short interval while preserving a hierarchy of estimates for the spacetime geometry. The integrated shear can then overcome the initially positive outgoing expansion. The essential analytic task is to propagate the large initial data far enough for this focusing to occur. Klainerman and Rodnianski \cite{KlainermanRodnianski2012} gave a proof using scale-invariant estimates and fewer curvature derivatives. An and Luk \cite{AnLuk2017} obtained a scale-critical formation criterion, and An \cite{An2022} developed an approach based on signatures for decay rates. Klainerman, Luk, and Rodnianski \cite{KlainermanLukRodnianski2014} showed how an angularly localized concentration of radiation can form a trapped surface by combining null evolution with a deformation of the final surface. These works identify complementary aspects of the vacuum problem: the size, propagation, and angular distribution of the focusing flux. We also mention the works \cites{An2020,An2022,An2025,AnAthanasiou2022,AnHan2021,AnHe2024,AnHe2025,AnWu2026} of the first author on extending trapped surface formation to broader physical settings.
		
		For a self-gravitating perfect fluid, spherical symmetry already permits nontrivial dynamical collapse. Earlier research dates back to the 1930s; see \cite{Tolman1934,Lemaitre1933}. Assuming a pressureless Einstein--dust model in spherical symmetry in an FLRW background, Oppenheimer and Snyder \cite{OppenheimerSnyder1939} proved the dynamical formation of singularities (later known as black holes); see also \cite{OppenheimerVolkoff1939,Tolman1939} for indications of gravitational collapse of a neutron star with a mass above the TOV limit. Burtscher and LeFloch \cite{BurtscherLeFloch2014} proved trapped-surface formation from initially untrapped data for spherically symmetric Einstein--Euler solutions of bounded variation.
	
		In a different direction, Christodoulou \cite{Christodoulou1984} constructed spherically symmetric naked singularities for the Einstein-dust, and Guo, Had\v{z}i\'{c}, and Jang \cite{GuoHadzicJang2023} for the Einstein--Euler system.
		These results concern spherical evolution. The present problem requires control of angular derivatives, shear, and Weyl curvature together with the fluid, and a focusing lower bound that holds along every generator of the final sphere.
		
		Matter focusing results also appear for nonsymmetric Einstein--matter system with one wave-speed, for example \cites{Yu2013,AnAthanasiou2022} for the Einstein-Maxwell system. For fluid (Euler) matter, an additional difficulty is the different wave speeds between gravitational and acoustic propagations. The spacetime geometry is organized by metric-null hypersurfaces, whereas the fluid evolves along its timelike flow and acoustic characteristics. Thus the leading-order positive fluid-matter term in the Raychaudhuri equation cannot simply be propagated by metric-null characteristic equations as for the shear tensors. One needs estimates adapted to the fluid velocity, control of differentiated matter sources in the curvature equations, and a separate conservation argument to retain enough matter flux in the focusing region.
	\subsection{Geometric Setting and Main Result}
	\label{SS:IntroMain}

		We use optical functions $u$ and $\ubar$ for $\gfour$ and their level hypersurfaces $\hypersur_u$ and $\hypersurbar_{\ubar}$. Their intersections $S_{u,\ubar}$ are spacelike two-spheres, with induced metric $\gsphere$. Let $e_4$ and $e_3$ be their outgoing and incoming future null normals, normalized by $\gfour(e_3,e_4)=-2$. The positive null lapse $\lapse$ is chosen so that $e_4=\lapse^{-1}\partial_{\ubar}$. We write $\upchi$ and $\uchi$ for the corresponding null second fundamental forms, $\gtr\upchi$ and $\gtr\uchi$ for their traces with respect to $\gsphere$, and $\hchi$ and $\huchi$ for their trace-free parts. Thus $S_{u,\ubar}$ is trapped precisely when
	\begin{equation}\label{EQ:IntroTrapped}
		{\gtr\upchi<0\qquad\hbox{and}\qquad\gtr\uchi<0
			\quad\hbox{everywhere on }S_{u,\ubar}.}
	\end{equation}

		The initial hypersurfaces are $\hypersur_{\uinf}$ and $\hypersurbar_0$, where $\uinf<\uts<0$ and $\uts$ is the final incoming-coordinate value. The short interval on $\hypersur_{\uinf}$ has length $\delta>0$. We use $\vartheta=(\sphere^1,\sphere^2)\in\Stwo$ for the transported angular label. For the perfect fluid, the scalar Ricci component relevant to outgoing focusing is
	\begin{equation}\label{EQ:IntroMatter}
		{
			2\FAlphaAB=\Ricfour{4}{4}=\EMT_{44}
			=(\density+\pressure)(v_4)^2,
			\qquad v_4=\gfour(v,e_4).
		}
	\end{equation}
		In particular, $\FAlphaAB$ is a scalar matter component, distinct from the Weyl curvature tensor component $\AlphaAB$. Its normalization accounts for the factor $2$ in the focusing equation below.
		
		We consider the fast-injection regime of Section~\ref{SS:Norms}, with short scale $\delta>0$ and weighted initial-data bounds \eqref{Data} at differentiability level $N\geq6$. In the characteristic-data setting of Section~\ref{SS:InitialData}, Theorem~\ref{TH:Main} gives the following criterion.
	
	\begin{theorem}[Main result; cf. Theorem~\ref{TH:Main}]
		\label{TH:IntroMain}
		For the Einstein--Euler equations \eqref{1.0}, consider characteristic initial data on $\hypersur_{\uinf}\cup\hypersurbar_0$ in the fast-injection regime described above. Assume
		\begin{itemize}[leftmargin=*,itemsep=5pt]
			\item the weighted initial-data bounds \eqref{Data} on $u=\uinf$ and $\ubar=0$;
			\item the common upper bound along $u=\uinf$,
			\begin{equation}\label{EQ:IntroUpper}
				\int_0^\delta |\uinf|^2
				\bigl(|\hchi|^2+2\FAlphaAB\bigr)(\uinf,\ubarprime,\vartheta)
				\,\diff\ubarprime<2|\uinf|
				\quad\hbox{for every }\vartheta\in\Stwo.
			\end{equation}
		\end{itemize}
		Then the initial hypersurfaces contain no trapped surface in the initial region under consideration. 
        
        Moreover, we have:
		\begin{itemize}[leftmargin=*,itemsep=5pt]
			\item \textup{(Shear focusing.)} If
			\begin{equation}\label{EQ:IntroShear}
				\int_0^\delta |\uinf|^2|\hchi|^2(\uinf,\ubarprime,\vartheta)
				\,\diff\ubarprime>5|\uts|
				\quad\hbox{for every }\vartheta\in\Stwo,
			\end{equation}
			then $S_{\uts,\delta}$ is a trapped surface.
			\item \textup{(Fluid focusing.)} If
			\begin{equation}\label{EQ:IntroFluid}
				\int_0^\delta |\uinf|^2\FAlphaAB(\uinf,\ubarprime,\vartheta)
				\,\diff\ubarprime>4|\uts|
				\quad\hbox{for every }\vartheta\in\Stwo,
			\end{equation}
			then $S_{\uts,\delta_T}$ is a trapped surface, where
			$\delta_T:=\bigl(1+2(\uts-\uinf)\bigr)\delta$.
		\end{itemize}
	\end{theorem}
	
	\noindent\begin{minipage}{\textwidth}
		\centering
		\begin{minipage}[t]{0.48\textwidth}
			\centering
			\textbf{Shear focusing}\par\smallskip
			\begin{tikzpicture}[scale=0.82,font=\small,line join=round]
				\path[use as bounding box] (0,-2.1) rectangle (5.7,3.8);
				\coordinate (A) at (3,-1);
				\coordinate (B) at (4,0);
				\coordinate (C) at (1,1);
				\coordinate (D) at (2,2);
				\fill[black!18] (A)--(B)--(D)--(C)--cycle;
				\draw[dashed,black!55] (2.65,-1.35)--(4.35,0.35);
				\draw[dashed,black!55] (0.65,0.65)--(2.35,2.35);
				\draw[dashed,black!55] (0.65,1.35)--(3.35,-1.35);
				\draw[dashed,black!55] (1.65,2.35)--(4.35,-0.35);
				\draw[thick] (A)--(B)--(D)--(C)--cycle;
				\path (A)--node[midway,sloped,below=4pt] {$\hypersur_{\uinf}\ (u=\uinf)$} (B);
				\path (C)--node[midway,sloped,above=4pt] {$\hypersur_{\uts}$} (D);
				\path (C)--node[midway,sloped,below=4pt] {$\hypersurbar_0\ (\ubar=0)$} (A);
				\path (D)--node[midway,sloped,above=4pt] {$\hypersurbar_\delta\ (\ubar=\delta)$} (B);
				\fill (D) circle (1.7pt);
				\node[above=9pt] at (D) {$S_{\uts,\delta}$};
			\end{tikzpicture}
		\end{minipage}\hfill
		\begin{minipage}[t]{0.48\textwidth}
			\centering
			\textbf{Fluid focusing}\par\smallskip
			\begin{tikzpicture}[scale=0.82,font=\small,line join=round]
				\path[use as bounding box] (0,-2.1) rectangle (5.7,3.8);
				\coordinate (A) at (3,-1);
				\coordinate (B) at (4,0);
				\coordinate (C) at (1,1);
				\coordinate (D) at (2,2);
				\coordinate (T) at (2.6,2.6);
				\fill[black!18] (A)--(B)--(T)--(C)--cycle;
				\draw[dashed,black!55] (2.65,-1.35)--(4.85,0.85);
				\draw[dashed,black!55] (0.65,0.65)--(2.95,2.95);
				\draw[dashed,black!55] (0.65,1.35)--(3.35,-1.35);
				\draw[dashed,black!55] (T)--(4.6,0.6);
				\draw[dashed,black!60] (D)--node[midway,sloped,below=2pt] {$\ubar=\delta$} (B);
				\draw[thick] (A)--(B)--(T)--(C)--cycle;
				\path (A)--node[midway,sloped,below=4pt] {$\hypersur_{\uinf}\ (u=\uinf)$} (B);
				\path (C)--node[midway,sloped,above=4pt] {$\hypersur_{\uts}$} (T);
				\path (C)--node[midway,sloped,below=4pt] {$\hypersurbar_0\ (\ubar=0)$} (A);
				\path (T)--node[midway,sloped,above=4pt] {$\hypersurbar_{\delta_T}\ (\ubar=\delta_T)$} (4.6,0.6);
				\path (B)--node[midway,sloped,above=3pt] {$\mathcal T$} (T);
				\fill (T) circle (1.7pt);
				\node[above=9pt] at (T) {$S_{\uts,\delta_T}$};
			\end{tikzpicture}
		\end{minipage}
		
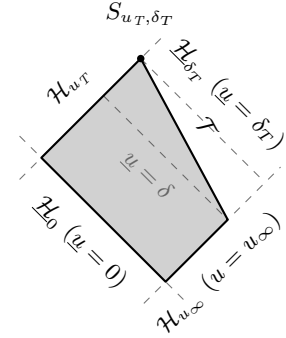
\captionof{figure}{Schematic regions for the two alternatives in Theorem~\ref{TH:IntroMain}. The dots represent the final trapped spheres. On the right, the timelike boundary $\mathcal T$ is given by $\ubar=\bigl(1+2(u-\uinf)\bigr)\delta$, so its endpoint has $\ubar=\delta_T$.}
		\label{FIG:IntroFormation}
	\end{minipage}
	\par\medskip

		All three integrals are taken along individual null generators, with coordinate measure $\diff\ubarprime$; they are not integrals over a sphere. The common bound controls the initial total focusing, whereas the two alternatives prescribe a lower bound that becomes sufficient at the later value $\uts$. The constants and the final spheres differ because the two fluxes are propagated by different arguments. In particular, the fluid lower bound contains $\FAlphaAB$, while the common upper bound contains $2\FAlphaAB$. No spherical symmetry is imposed, although the fluid estimates impose quantitative control of angular variation.
	
	\subsection{Physical Heuristics}\label{SS:PhysicalHeuristics}
In this section, we discuss the physical heuristics underlying trapped surface formation. Inspired by Christodoulou \cite{Christodoulou2009}, our proof of trapped surface formation relies on a suitable choice of initial data prescribed by a short-pulse ansatz and a delicate argument for propagating estimates using the vectorfield method. We begin by discussing the mechanism of trapped surface formation under the short-pulse ansatz for several models: the Einstein vacuum, Einstein--Maxwell, and Einstein--scalar field equations. Then, in Section \ref{SSS:HeuristicsEuler}, we introduce the corresponding mechanism for the Einstein--Euler system. Unlike the aforementioned systems, the Einstein--Euler system features multiple propagation speeds associated with light, sound, and fluid flow. Consequently, it is challenging to formulate an ansatz\footnote{See Section \ref{SS:ShortPulseAnsatz} for the short-pulse ansatz used in this article, which is adapted to the Einstein--Euler system.} that is compatible with the short-pulse ansatz in \cite{Christodoulou2009} and to propagate estimates for both the spacetime geometry and the fluid variables that remain consistent with this ansatz. We discuss our ansatz and estimates with more details in Section \ref{SS:IntroIngredients}.
		
		\subsubsection{Einstein vacuum scenario}\label{SSS:HeuristicsVacuum}
		
		We begin with the finite-radius form of Christodoulou's short-pulse theorem. It is obtained from the normalization $u_\infty=-1$ by rescaling the area radius. In particular, the free outgoing data are prescribed on the initial outgoing null hypersurface $\hypersur_{\uinf}$, whereas the data on the initial incoming null hypersurface $\hypersurbar_0$ are Minkowskian.
		
		\begin{theorem}[Christodoulou \cite{Christodoulou2009}]
			\label{TH:HeuristicsChristodoulou}
			Fix $\uinf<u_*<0$ and $B>0$. Prescribe Minkowskian characteristic data on $\hypersurbar_0$ for
			\[
			\uinf\leq u<0,
			\]
			and prescribe the outgoing shear on $\hypersur_{\uinf}$ for
			\[
			0\leq\ubar\leq\delta.
			\]
			There exists
			\[
			\delta_0=\delta_0(B,\uinf,u_*)>0
			\]
			such that, whenever $0<\delta\leq\delta_0$ and
			\begin{equation}\label{EQ:HeuristicsChristodoulouUpper}
				\sum_{i=0}^{5}\sum_{j=0}^{3}
				\delta^{\frac12+j}|\uinf|^i
				\bigl\|(\partial_{\ubar})^j\angnabla^i\hchi\bigr\|_
				{L^\infty_{\ubar}L^2(S_{\uinf,\ubar})}
				\leq B,
			\end{equation}
			the vacuum Einstein equations admit a regular development in
			\[
			\uinf\leq u\leq u_*,
			\qquad
			0\leq\ubar\leq\delta.
			\]
			If, in addition,
			\begin{equation}\label{EQ:HeuristicsChristodoulouLower}
				\inf_{\vartheta\in\Stwo}|\uinf|^2
				\int_0^\delta
				|\hchi|^2(\uinf,\ubarprime,\vartheta)\,\diff\ubarprime
				\geq M_*>2|u_*|,
			\end{equation}
			then $S_{u_*,\delta}$ is a trapped surface, provided $\delta_0$ is chosen sufficiently small depending also on the strict gap $M_*-2|u_*|$.
		\end{theorem}
		
		We now explain the focusing mechanism behind this theorem and its matter analogues. The outgoing coordinate ranges over $0\leq\ubar\leq\delta$. We decompose $\upchi_{AB}$ and $\uchi_{AB}$ into their trace and trace-free parts:
		\begin{equation}\label{EQ:HeuristicsDecomposition}
			\upchi_{AB}
			=\frac12(\gtr\upchi)\gsphere_{AB}+\hchi_{AB},
			\qquad
			\uchi_{AB}
			=\frac12(\gtr\uchi)\gsphere_{AB}+\huchi_{AB}.
		\end{equation}
		Here $\gsphere_{AB}$ is the induced metric on $S_{u,\ubar}$. A two-sphere $S_{u,\ubar}$ is trapped if
		\[
		\gtr\upchi<0
		\qquad\hbox{and}\qquad
		\gtr\uchi<0
		\]
		hold pointwise on $S_{u,\ubar}$.
		
		For the Minkowskian data on $\hypersurbar_0$, each $S_{u,0}$ is a standard two-sphere of radius $|u|$, and
		\begin{equation}\label{EQ:HeuristicsMinkowski}
			\lapse(u,0)=1,
			\qquad
			\gtr\upchi(u,0)=\frac{2}{|u|},
			\qquad
			\gtr\uchi(u,0)=-\frac{2}{|u|}.
		\end{equation}
		Thus the incoming expansion is negative on $\hypersurbar_0$. The exact incoming Raychaudhuri equation and the identity
		\[
		e_3\log\lapse=-2\uomega
		\]
		give
		\begin{equation}\label{EQ:HeuristicsIncomingVacuum}
			e_3(\lapse^{-1}\gtr\uchi)
			=
			-\lapse^{-1}
			\left(
			\frac12(\gtr\uchi)^2+|\huchi|^2
			\right)
			\leq0.
		\end{equation}
		The short-pulse estimates keep $\gtr\uchi$ negative on $\hypersur_{\uinf}$, while \eqref{EQ:HeuristicsIncomingVacuum} propagates this sign throughout the development. It therefore remains only to make the outgoing expansion negative.
		
		Using $\Ricfour{4}{4}=0$, the transport equations governing the outgoing expansion and the outgoing shear are
		\begin{align}
			e_4\gtr\upchi+\frac12(\gtr\upchi)^2
			&=
			-|\hchi|^2-2\omega\gtr\upchi,
			\label{EQ:HeuristicsVacuumRaychaudhuri}\\
			\angD_3\hchi+\frac12(\gtr\uchi)\hchi
			&=
			-\frac12(\gtr\upchi)\huchi
			+\angnabla\tensorproduct\upeta
			+2\uomega\hchi
			+\upeta\tensorproduct\upeta.
			\label{EQ:HeuristicsVacuumShear}
		\end{align}
		Since
		$
		e_4=\lapse^{-1}\partial_{\ubar},
e_4\log\lapse=-2\omega
		$,	the first equation becomes
		\[
		\partial_{\ubar}(\lapse^{-1}\gtr\upchi)
		=
		-\frac12(\gtr\upchi)^2-|\hchi|^2
		\leq-|\hchi|^2.
		\]
		Integrating along a fixed outgoing generator yields
		\begin{equation}\label{EQ:HeuristicsVacuumIntegrated}
			\begin{aligned}
				(\lapse^{-1}\gtr\upchi)(u,\ubar,\vartheta)
				&\leq
				\gtr\upchi(u,0)
				-
				\int_0^{\ubar}
				|\hchi|^2(u,\ubarprime,\vartheta)\,\diff\ubarprime\\
				&=
				\frac{2}{|u|}
				-
				\int_0^{\ubar}
				|\hchi|^2(u,\ubarprime,\vartheta)\,\diff\ubarprime.
			\end{aligned}
		\end{equation}
		
		The leading part of \eqref{EQ:HeuristicsVacuumShear}, together with
		\[
		\gtr\uchi=-\frac{2}{|u|}+\mbox{lower order terms},
		\]
		transports $|u|\hchi$ approximately unchanged in the incoming direction. Using
		$e_3
		=
		\lapse^{-1}
		\left(\partial_u+\Lambda^A\partial_{\theta^A}
		\right)$ and applying $e_3\log\lapse=-2\uomega$, we further obtain
		\begin{equation}\label{EQ:HeuristicsVacuumWeightedShear}
			\begin{aligned}
				\partial_u\bigl(|u|^2\lapse^2|\hchi|^2\bigr)
				={}&
				2|u|^2\lapse^3
				\left(
				-\frac12(\gtr\upchi)\huchi
				+\angnabla\tensorproduct\upeta
				+\upeta\tensorproduct\upeta
				\right)\cdot\hchi\\
				&-|u|^2\Lambda^A
				\partial_{\theta^A}(\lapse^2|\hchi|^2)\\
				&-|u|^2\lapse^3
				\left(
				\gtr\uchi+\frac{2}{|u|}
				\right)|\hchi|^2
				+2|u|(\lapse-1)\lapse^2|\hchi|^2.
			\end{aligned}
		\end{equation}
		The initial shear is chosen so that
		\begin{equation}\label{EQ:HeuristicsChristodoulouSize}
			\|\hchi\|_{L^\infty(S_{\uinf,\ubar})}
			\approx
			\frac{\delta^{-1/2}}{|\uinf|},
			\qquad
			|\uinf|^2
			\int_0^\delta
			|\hchi|^2(\uinf,\ubarprime,\vartheta)
			\,\diff\ubarprime
			\geq M_*>2|u_*|
		\end{equation}
		for every $\vartheta\in\Stwo$. The interval length $\delta$ compensates for the squared amplitude $\delta^{-1}$, since
		$\delta(\delta^{-1/2})^2=1$.
        
	\vspace{1mm}
		The short-pulse estimates control the incoming shear, torsion, lapse, angular shift, and angular derivatives in \eqref{EQ:HeuristicsVacuumWeightedShear}. Since $\lapse=1$ on $\hypersur_{\uinf}$, integration in $u$ and $\ubar$ yields
		\begin{equation}\label{EQ:HeuristicsVacuumFluxComparison}
			\left|
			|u|^2
			\int_0^\delta
			|\hchi|^2(u,\ubarprime,\vartheta)\,\diff\ubarprime
			-
			|\uinf|^2
			\int_0^\delta
			|\hchi|^2(\uinf,\ubarprime,\vartheta)\,\diff\ubarprime
			\right|
			\leq C\delta^{1/2}.
		\end{equation}
		Here $C$ depends on the fixed initial bounds and on the interval $\uinf\leq u\leq u_*$,
		but not on $\delta$. Combining \eqref{EQ:HeuristicsVacuumIntegrated} and \eqref{EQ:HeuristicsVacuumFluxComparison} on $\hypersur_{u_*}$, for sufficiently small $\delta$ we obtain
		\begin{equation}\label{EQ:HeuristicsChristodoulouConclusion}
			\begin{aligned}
				(\lapse^{-1}\gtr\upchi)(u_*,\delta,\vartheta)
				&\leq
				\frac{2}{|u_*|}
				-
				\frac{|\uinf|^2}{|u_*|^2}
				\int_0^\delta
				|\hchi|^2(\uinf,\ubarprime,\vartheta)
				\,\diff\ubarprime
				+C\delta^{1/2}\\
				&\leq
				\frac{2}{|u_*|}
				-
				\frac{M_*}{|u_*|^2}
				+C\delta^{1/2}
				<0.
			\end{aligned}
		\end{equation}
		Note that the strict inequality follows from
		$
		M_*>2|u_*|
		$. For example, $u_*=-1/4$ and $M_*=1$ produce the familiar leading comparison
		$8-16<0$.

        \vspace{1mm}
		Since $\lapse>0$ and $\gtr\uchi<0$, the sphere $S_{u_*,\delta}$ is trapped. The choice
		\[
		|\hchi|
		\approx
		\frac{\delta^{-1/2}}{|\uinf|}
		\]
		on $\hypersur_{\uinf}$ also determines the hierarchy in $\delta$ and $|u|$ for the remaining geometric quantities. The rigorous proof consists in controlling every error term in \eqref{EQ:HeuristicsVacuumWeightedShear} and extending the regular development from $\hypersur_{\uinf}$ to $\hypersur_{u_*}$.
		
		\subsubsection{Einstein--Maxwell and Einstein--scalar field scenarios}
		\label{SSS:HeuristicsMatter}
		
		The same mechanism extends to the Einstein--Maxwell and Einstein--scalar field systems: the corresponding null matter component contributes an additional nonnegative focusing flux to the outgoing Raychaudhuri equation.
		
		For the Einstein--Maxwell system,
		\[
		\Ricfour{\mu}{\nu}
		=
		F_{\mu\lambda}F_\nu{}^\lambda
		-
		\frac14
		\gfour_{\mu\nu}
		F_{\alpha\beta}F^{\alpha\beta},
		\qquad
		\Dfour^\mu F_{\mu\nu}=0,
		\qquad
		\Dfour_{[\mu}F_{\nu\lambda]}=0.
		\]
		The outgoing electromagnetic component is the sphere-tangent one-form
		\[
		(\alpha_F)_A=F(e_A,e_4),
		\]
		and
		\[
		\Ricfour{4}{4}=|\alpha_F|^2.
		\]
		Hence the outgoing Raychaudhuri equation becomes
		\begin{equation}\label{EQ:HeuristicsMaxwellRaychaudhuri}
			\partial_{\ubar}(\lapse^{-1}\gtr\upchi)
			=
			-\frac12(\gtr\upchi)^2
			-|\hchi|^2
			-|\alpha_F|^2.
		\end{equation}
		Thus $\alpha_F$ enters with the same favorable sign as $\hchi$. Moreover, $\alpha_F$ is transported by the metric-null characteristic transport equation as well (as for the shear tensor). For the hyperbolic estimates, we can thus treat the gravitational field and the Maxwell field in similarly the same framework.  
        
        \vspace{1mm}
        
        If we further impose that the outgoing data on $\hypersur_{\uinf}$ satisfy
		\begin{equation}\label{EQ:HeuristicsMaxwellData}
			\left(
			|\hchi|^2+|\alpha_F|^2
			\right)^{1/2}
			\approx
			\frac{\delta^{-1/2}}{|\uinf|},
			\qquad
			|\uinf|^2
			\int_0^\delta
			\left(
			|\hchi|^2+|\alpha_F|^2
			\right)
			(\uinf,\ubarprime,\vartheta)
			\,\diff\ubarprime
			\geq M_*>2|u_*|,
		\end{equation}
		then the short-pulse transport estimates of An--Athanasiou \cite{AnAthanasiou2022} propagate the combined weighted flux from $\hypersur_{\uinf}$ to $\hypersur_{u_*}$. Integrating \eqref{EQ:HeuristicsMaxwellRaychaudhuri} therefore gives the analogue of \eqref{EQ:HeuristicsChristodoulouConclusion}, with $|\hchi|^2$ replaced by
		\[
		|\hchi|^2+|\alpha_F|^2,
		\]
		and yields the trapped sphere $S_{u_*,\delta}$.
        
		\vspace{2mm}
		For the Einstein--scalar field system
		\begin{equation}\label{EQ:HeuristicsScalarSystem}
			\Ricfour{\mu}{\nu}
			=
			\Dfour_\mu\phi\,\Dfour_\nu\phi,
			\qquad
			\Dfour^\mu\Dfour_\mu\phi=0,
		\end{equation}
		the corresponding Raychaudhuri equation becomes
		\begin{equation}\label{EQ:HeuristicsScalarRaychaudhuri}
			\partial_{\ubar}(\lapse^{-1}\gtr\upchi)
			=
			-\frac12(\gtr\upchi)^2
			-|\hchi|^2
			-(e_4\phi)^2.
		\end{equation}
		Here $e_4\phi$ plays the same focusing role as $\hchi$ and it is again transported by the metric-null characteristic transport equation(as for the shear tensor). In the corresponding characteristic problem with Minkowskian data on $\hypersurbar_0$, we can prescribe the below initial data on $\hypersur_{\uinf}$
		\begin{equation}\label{EQ:HeuristicsScalarData}
			\left(
			|\hchi|^2+(e_4\phi)^2
			\right)^{1/2}
			\approx
			\frac{\delta^{-1/2}}{|\uinf|},
			\qquad
			|\uinf|^2
			\int_0^\delta
			\left(
			|\hchi|^2+(e_4\phi)^2
			\right)
			(\uinf,\ubarprime,\vartheta)
			\,\diff\ubarprime
			\geq M_*>2|u_*|.
		\end{equation}
		The wave equation propagates $|u|e_4\phi$ in the same leading-order manner in which the shear equation propagates $|u|\hchi$. The regular-development estimates for the special $\Omega=1$ and $\mu^*=0$ case established by An \cite{An2025} then give the corresponding weighted-flux comparison. Through \eqref{EQ:HeuristicsScalarRaychaudhuri}, we then know that $S_{u_*,\delta}$ is the formed trapped surface.
        
        \vspace{1mm}
      In both systems, for the leading-order terms, the concentrated matter flux are treated in the same way as for the shear flux, since the matter field and the gravitational field share the same metric-induced wave speed. And both the shear tensor and the matter field can independently form the trapped surfaces.  

			\subsubsection{Einstein--Euler scenario}
		\label{SSS:HeuristicsEuler}
		
		For the Einstein--Euler system, the outgoing Raychaudhuri equation is
		\begin{equation}\label{EQ:HeuristicsEulerRaychaudhuri}
			\partial_{\ubar}(\lapse^{-1}\gtr\upchi)
			=
			-\frac12(\gtr\upchi)^2
			-|\hchi|^2
			-(\density+\pressure)v_4^2.
		\end{equation}
		It contains two nonnegative focusing densities: the gravitational contribution $|\hchi|^2$ and the fluid contribution
		\[
		(\density+\pressure)v_4^2=:2\FAlphaAB.
		\]
		These lead to the two alternative independent trapped surface formation criterion in Theorem~\ref{TH:IntroMain}.

\vspace{2mm}
        
		\medskip
		\noindent\textbf{Focusing by gravitational radiation.}
		The transport equation for the outgoing shear is
		\begin{equation}\label{EQ:HeuristicsEulerShear}
			\angD_3\hchi+\frac12(\gtr\uchi)\hchi
			=
			-\frac12(\gtr\upchi)\huchi
			+\angnabla\tensorproduct\upeta
			+2\uomega\hchi
			+\upeta\tensorproduct\upeta
			+\FRhodual.
		\end{equation}
		The additional angular Ricci term is of the expression
		\[
		\FRhodual_{AB}
		=
		(\density+\pressure)
		\left(
		v_Av_B
		-\frac12|\angv|^2\gsphere_{AB}
		\right).
		\]
		Note that our large null velocity component $v_4$ does not appear in this angular trace-free source. Since the angular fluid velocity is perturbative, the leading-order estimates for \eqref{EQ:HeuristicsEulerShear} preserve up lower order terms, and it is the same relation between $|u|\hchi$ and $|\uinf|\hchi(\uinf,\ubar,\vartheta)$ as in the vacuum argument. Consequently, the strict initial lower bound \eqref{EQ:IntroShear} implies
		\begin{align}\label{ES:IntrohchiBound}
			\int_0^\delta
		|\uts|^2
		|\hchi|^2(\uts,\ubarprime,\vartheta)
		\,\diff\ubarprime
		>4|\uts|
	\end{align}
		for every $\vartheta\in\Stwo$. Because the fluid term in \eqref{EQ:HeuristicsEulerRaychaudhuri} is nonpositive, it can be discarded due to the favorable sign.  One key difference between here and the Einstein vacuum scenarios is that the initial incoming characteristic hypersurface is not Minkowskian, since the fluid is (assumed to be) away from vacuum. Using transport equations along the incoming null geodesics (e.g. \eqref{e3gtrupchi}), we derive estimates for the data\footnote{See Section \ref{SS:CharData} for the construction and estimates for the initial data.} on $\hypersurbar_0$. In particular, we bound the initially positive outgoing expansion $\gtr\upchi$ by
		\begin{align}\label{ES:IntroupchiBound}
			\abs{\gtr\upchi(\uts,0,\vartheta)-\frac{2}{|\uts|}}
			\leq\frac{1}{|\uts|^2}+\iniBAsize\delta a^{1/2},
		\end{align}
		where $\iniBAsize\delta a^{1/2}<<\frac{1}{\abs{\uts}}$. Combining \eqref{EQ:HeuristicsEulerRaychaudhuri}, \eqref{ES:IntrohchiBound} and \eqref{ES:IntroupchiBound}, following the similar approach as in Section \ref{SSS:HeuristicsVacuum}, for $\vartheta\in \mathbb{S}^2$ we obtain:
		\begin{equation}\label{EQ:HeuristicsEulerShearNegative}
			(\lapse^{-1}\gtr\upchi)
			(\uts,\delta,\vartheta)
			<
			\frac{3}{|\uts|}
			-
			\frac{4}{|\uts|}
			+
			\iniBAsize\delta a^{1/2}
			<0.
		\end{equation}
Hence, $S_{u_T, \delta}$ is the formed trapped surface. 

\vspace{2mm}
        
		\medskip
		\noindent\textbf{Focusing by the fluid.}
		In our second trapped surface formation mechanism, the shear tensor could be small, and the fluid alone supplies the necessary focusing. In this paper, our fast-injection hierarchy is designed to be
			\begin{align}\label{EQ:IntroFast0}
			|v_4|\sim&~ b^{1/2},& |v_3|\sim&~ b^{-1/2},& \abs{v_A,\lnenth,\ent}\sim&~ \BAfluid,
		\end{align}
		where $v_I=\gfour(v,e_I)$, $b\sim \delta^{-1}$ is a large parameter that describes the fast-inject speed of the fluid and $\BAfluid$ is a small parameter that control the rest of the fluid variables.
		
		The essential difference from gravitational, electromagnetic, and scalar radiation is that fluid particles follow timelike, rather than null, trajectories. Consequently, they need not remain within a region of fixed outgoing width. To capture the transported particles, we introduce the expanding lateral boundary
		\[
		\ubar
		=
		\bigl(1+2(u-\uinf)\bigr)\delta,
		\qquad
		\uinf\leq u\leq\uts.
		\]
		Its unit outer normal and the particle flux through it are
		\[
		\mathbf S
		=
		\frac{2\delta e_4-e_3}{\sqrt{8\delta}},
		\qquad
		\gfour(nv,\mathbf S)
		=
		\frac{n(2\delta v_4-v_3)}{\sqrt{8\delta}}.
		\]
		The ansatz \eqref{EQ:IntroFast0} implies a condition
		\[
		2\delta(-v_4)\geq-v_3,
		\]
		which guaranties that the outward velocity of the fluid does not exceed the slope of this boundary. Thus, no particle flux is lost through the lateral boundary in the unfavorable direction. Applying particle-number conservation,
	\begin{align}
		\Dfour_\mu I^\mu&=0,&
        I^\mu&=nv^\mu
	\end{align}
		to region as in the fluid-focusing case in Figure \ref{FIG:IntroFormation} gives
		\begin{equation}\label{EQ:HeuristicsEulerParticleFlux}
			\int_0^{(1+2(\uts-\uinf))\delta}
			\int_{S_{\uts,\ubar}}
			n(-v_4)\lapse
			\,\diff\spherevol\,\diff\ubar
			\geq
			\int_0^\delta
			\int_{S_{\uinf,\ubar}}
			n(-v_4)\lapse
			\,\diff\spherevol\,\diff\ubar.
		\end{equation}
	Estimate \eqref{EQ:HeuristicsEulerParticleFlux} is initially a comparison of fluxes integrated over the entire sphere; by itself, it does not yet give the pointwise-in-angle lower bound required by Raychaudhuri's equation. To derive the favorable lower bound, we first estimate the angular derivative that controls the oscillation from one generator to another. Then, using $\density+\pressure=n\enth$ and the estimates \eqref{EQ:IntroFast0}, in Lemma~\ref{LE:EstimatesforuFAlpha} we obtain the following lower bound along every final null generator:
		\begin{equation}\label{EQ:HeuristicsEulerGeneratorFlux}
			\begin{aligned}
				|\uts|^2
				\int_0^{(1+2(\uts-\uinf))\delta}
				(\density+\pressure)v_4^2
				(\uts,\ubarprime,\vartheta)
				\,\diff\ubarprime
				&\geq
				\frac12|\uinf|^2
				\int_0^\delta
				(\density+\pressure)v_4^2
				(\uinf,\ubarprime,\vartheta)
				\,\diff\ubarprime\\
				&>4|\uts|.
			\end{aligned}
		\end{equation}
		Recalling \eqref{EQ:IntroMatter}, the last inequality is precisely the initial fluid hypothesis \eqref{EQ:IntroFluid}.
		Substituting \eqref{EQ:HeuristicsEulerGeneratorFlux} into \eqref{EQ:HeuristicsEulerRaychaudhuri} and discarding the favorable shear term gives
		\begin{equation}\label{EQ:HeuristicsEulerFluidNegative}
			\begin{aligned}
				&(\lapse^{-1}\gtr\upchi)
				\bigl(
				\uts,
				(1+2(\uts-\uinf))\delta,
				\vartheta
				\bigr)\\
				&\qquad\leq
				\frac{3}{|\uts|}
				+\iniBAsize\delta a^{1/2}
				-
				\int_0^{(1+2(\uts-\uinf))\delta}
				(\density+\pressure)v_4^2
				(\uts,\ubarprime,\vartheta)
				\,\diff\ubarprime\\
				&\qquad<
				\frac{3}{|\uts|}
				-
				\frac{4}{|\uts|}
				+
				\iniBAsize\delta a^{1/2}
				<0.
			\end{aligned}
		\end{equation}
For the incoming null expansion scalar $\gtr\uchi$, the incoming Raychaudhuri equation implies
\begin{align}
			e_3(\lapse^{-1}\gtr\uchi)
		=
		-\lapse^{-1}
		\left(
		\frac12(\gtr\uchi)^2
		+|\huchi|^2
		+(\density+\pressure)v_3^2
		\right)
		\leq0.
\end{align}
		Therefore the incoming expansion remains negative. Combining this fact with \eqref{EQ:HeuristicsEulerFluidNegative} yields the trapped surfaces asserted in the fluid focusing case in Theorem~\ref{TH:IntroMain}.

	\subsection{New Ingredients}\label{SS:IntroIngredients}
	In the above discussion, we note that there are two important ingredients in the trapped surface formation mechanism: \textbf{1)} precise control of all geometric quantities and matter consistent with the short-pulse ansatz and \textbf{2)} discover the amplifying effect for the terms on the RHS of Raychaudhuri equation \eqref{EQ:HeuristicsEulerRaychaudhuri} which leads to the dynamical formation of a trapped surface on $\hypersur_{\uts}$, while there is the absence of trapped surfaces on $\hypersur_{\uinf}$. In this section, we will discuss the proofs addressing these two main challenges in more detail.

	\subsubsection{The two focusing mechanisms}
	\label{SS:IntroFocusing}
		The geometric mechanism is expressed by the outgoing null Raychaudhuri equation. In the present coordinates, its exact lapse-corrected form is
	\begin{equation}\label{EQ:IntroRaychaudhuri}
		{
			\partial_{\ubar}\bigl(\lapse^{-1}\gtr\upchi\bigr)
			=-\frac12(\gtr\upchi)^2-|\hchi|^2-2\FAlphaAB.
		}
	\end{equation}
		This follows from \eqref{EQ:e4upchi} and the lapse identity \eqref{EQ:omegauomega}; the lapse derivative cancels the term involving the connection coefficient $\omega$. The two nonnegative focusing contributions are the squared shear and the null matter density \eqref{EQ:IntroMatter}. At fixed $u$ and $\vartheta$, integration to any $\ubar_*\geq0$ within the controlled region gives
	\begin{equation}\label{EQ:IntroIntegrated}
		{\begin{aligned}
				\bigl(\lapse^{-1}\gtr\upchi\bigr)(u,\ubar_*,\vartheta)
				\leq{}&\bigl(\lapse^{-1}\gtr\upchi\bigr)(u,0,\vartheta)\\
				&-\int_0^{\ubar_*}\bigl(|\hchi|^2+2\FAlphaAB\bigr)
				(u,\ubarprime,\vartheta)\,\diff\ubarprime.
		\end{aligned}}
	\end{equation}
		The sign of the right-hand side is not determined by the initial lower bound alone. One must first retain that lower bound as the solution evolves from $\uinf$ to $\uts$, while controlling the initial value of the outgoing expansion on $\hypersurbar_0$. In the shear alternative, the transport estimate of Lemma~\ref{LE:Estimatesforuhchi} compares the appropriately weighted shear at the initial and final incoming-coordinate values. Together with the lapse estimate, it preserves enough of the initial shear flux for \eqref{EQ:IntroIntegrated} to become negative at $\ubar_* =\delta$.
	
		The fluid alternative needs a different comparison. The quantity $(\density+\pressure)v_4^2$ evolves through the Euler equations, and matter can cross the side of a region bounded by metric-null hypersurfaces. Particle-current conservation is therefore applied on an expanding timelike wedge. This first gives a comparison of sphere-integrated fluxes. Control of the enthalpy, velocity, sphere geometry, and angular oscillation is then needed to obtain the lower bound along every outgoing generator. These are distinct analytic steps: conservation supplies the total flux, and angular control distributes that information over the final sphere.

	\subsubsection{A fluid energy current adapted to null geometry}
	\label{SS:IntroFluidCurrent}
	
	The fluid estimates must be adapted to the double-null foliation used to measure gravitational focusing. In the nondegenerate subluminal range $0<\speed<1$, the sound cones lie inside the metric light cones, so a metric-null hypersurface is spacelike for the acoustic metric. We construct coercive fluid fluxes through these same hypersurfaces, with weights that reflect the strongly unequal null components of the velocity.
	
	In terms of the short scale $\delta$, the fast-injection hierarchy takes the form
	\begin{align}\label{EQ:IntroFast}
		|v_4|\sim&~ b^{1/2},& |v_3|\sim&~ b^{-1/2},& \abs{v_A,\lnenth,\ent}\sim&~ \BAfluid,
	\end{align}
	where $v_I=\gfour(v,e_I)$, $b\sim \delta^{-1}$ is a large parameter that describes the fast-injection speed of the fluid, and $\BAfluid$ is a small parameter that controls the rest\footnote{$v_3$ has extra smallness simply due to normalization $\gfour(v,v)=-1$.} of the fluid variables. Specifically, the angular velocity components, the logarithmic enthalpy $\lnenth=\log(\enth/\average{H})$, and the entropy are perturbative in the weighted hierarchy. Here $\average{H}>0$ is the fixed background enthalpy. We commute with $\delta\Dfour_4$, $|u|\Dfour_3$, $|u|\Dfour_A$ and $\{e_A\}$ with $A=1,2$ is an orthonormal frame tangent to $S_{u,\ubar}$. We note that these weights account for the short outgoing interval and the radial and angular scales of the geometry.

    \vspace{1mm}
	Let $(\varh,\vars,\varv)$ denote the variations of $(\lnenth,\ent,v)$ generated by these weighted derivatives, and let $\temp$ denote the temperature. The current introduced in Definition~\ref{DE:EnergyCurrent} is
	\begin{equation}\label{EQ:IntroCurrent}
		{\begin{aligned}
				\varJ^\mu={}&\frac{v^\mu}{\speed^2}
				\left(\varh^2+\frac{q^2\temp^2\speed^2}{\enth^2}\vars^2
				-2\frac{\temp}{\enth}\varh\vars\right)\\
				&+2\varv^\mu\left(\varh-\frac{\temp}{\enth}\vars\right)
				+v^\mu\gfour(\varv,\varv),
		\end{aligned}}
	\end{equation}
	where $q$ is the large constant chosen in Lemma~\ref{LE:Coerciveness}. As in Proposition \ref{PR:DiveofFluidCurrent}, the key property of $\varJ^\mu$ is that $\Dfour_\mu\varJ^\mu$ does not contain top-order terms; there is a precise cancellation due to the structure of the Euler equations\footnote{See Proposition \ref{PR:EQVariations} for the commuted Euler equations.}. This observation is applied frequently in compressible Euler and Einstein-Euler equations; see \cite{Speck2009,Speck2012,RodnianskiSpeck2013,Christodoulou2007,Song2025,LiLiZhu2026}.

    \vspace{1mm}
	One new observation in this article is that there are natural asymmetric coercive weights in the energy estimates when we apply the energy current method in double-null foliations. In particular:
	\begin{subequations}\label{ES:IntroCoercive}
		\begin{align}
			-\varJ_4\geq &C\cdot\left[(-v_4)(\vars^2+\varh^2+\varv_A\varv_A)+(-v_3)\varv_4^2\right],\\
			-\varJ_3\geq& C \cdot\left[(-v_3)(\vars^2+\varh^2+\varv_A\varv_A)+(-v_4)\varv_3^2\right].
		\end{align}
	\end{subequations}
	Lemma~\ref{LE:Coerciveness} quantifies these bounds. This asymmetry allows the energy estimates to retain the distinct sizes for fluid components throughout the coupled evolution. In particular, \eqref{ES:IntroCoercive} and \eqref{EQ:IntroCurrent} are compatible in the sense that, via energy estimates and initial data assumptions, one can show (see Proposition \ref{PR:HigherOrderEnergyEstimates}) that, heuristically:
	\begin{align}\label{ES:introvarv}
		|\varv_4|\sim&~ \BAfluid b^{1/2},& |\varv_3|\sim&~ \BAfluid b^{-1/2},& \abs{\varv_A,\varh,\vars}\sim&~ \BAfluid.
	\end{align}
	
	We also note that in the above description since we commute with spacetime covariant derivatives in all directions\footnote{For trapped surface formation for the Eintein vacuum equations (and all other single wave speed models), to close the estimates, one only needs to commute the equations with the angular derivative $\angnabla$, which significantly reduces the workload and difficulties.}, additional difficulties arise when large gravitational terms and fluid term enter into the bulk term $\int\Dfour_\mu\varJ^\mu$. We carefully estimate the size of every term in the commuted equations in Section \ref{SS:EstimatesFluidterms} to treat these terms.
	
	\subsubsection{Curvature renormalization and the coupled estimates}
	\label{SS:IntroCoupledEstimates}
	
		The geometric and fluid estimates must be closed at compatible derivative orders. The Bianchi equations differentiate the matter Ricci tensor, while the commuted Euler equations require derivatives of the geometry. We organize this interaction by modifying the curvature variables with matter terms. For example, Definition~\ref{DE:Renormalization} introduces
	\begin{subequations}\label{EQ:IntroRenormalization}
		\begin{align}
			\CAlphaAB:=&~\AlphaAB-\FAlphaAB\gsphere,&
			\CAlphabar:=&~\Alphabar-\FAlphabar\gsphere,\\
			\CBetaA:=&~\BetaA-\FBetaA,&
			\CBetabar:=&~\Betabar+\FBetabar,\\
			\CRho:=&~\Rho-\FRho-\frac{1}{3}\Scalar-\frac{1}{2}\hchi\cdot\huchi,&
			\CRhodual:=&~\Rhodual+\frac{1}{2}\huchi\wedge\hchi.
		\end{align}
	\end{subequations}
	Here $\AlphaAB,\Alphabar,\BetaA,\Betabar,\Rho,\Rhodual$ are Weyl curvature components as in Definition \ref{DE:CurvatureComponents}, the superscript $F$ denotes matter Ricci components, $\Scalar$ is the spacetime scalar curvature, and the dot contracts sphere-tangent indices. These combinations absorb selected null derivatives of matter into the curvature variables and incorporate the quadratic shear correction in the Gauss curvature\footnote{As in \eqref{EQ:Renormalization2}, the Gauss curvature of the $\suubar$-induced metric $\gsphere$ satisfies $\GaussCurvature:=-\CRho-\frac{1}{4}\gtr\upchi\gtr\uchi$.}. The Ricci curvature components are estimated using the fluid energy estimates of Section~\ref{SS:IntroFluidCurrent}. We note that the closer the Weyl curvature components are in size to their Ricci curvature counterparts, the more important renormalization becomes. In the present work, $\AlphaAB\sim a^{\frac{1}{4}}\FAlphaAB$, $\BetaA\sim\BAfluid\FBetaA$, and $\Betabar\sim\BAfluid\FBetabar$. Thus, analogous estimates for $\AlphaAB,\BetaA,\Betabar$ hold even without the renormalization \eqref{EQ:IntroRenormalization}. Renormalization is more important for $\Alphabar,\Rho,\Rhodual$, which lie at the critical level. We expect this comparison between the Weyl and Ricci curvature components to play an even greater role in establishing related results\footnote{For example, for trapped surface formation in the Einstein-scalar field and Einstein-Maxwell equations, the renormalization \eqref{EQ:IntroRenormalization} is important for $\BetaA$ since it is at the critical level there; see \cite{AnAthanasiou2022,Yu2011,An2025}.} under more critical assumptions on the initial geometric and fluid data.
	
	At top order, the $N-1$ curvature commutations in Proposition~\ref{PR:topL2curvature} are matched with order-$N$ matter estimates. Transport and Hodge estimates for the Ricci coefficients, including the renormalized mass-aspect quantities in Proposition~\ref{PR:topL2RicciCoefficients}, recover the geometric derivatives needed by the fluid commutators. The resulting bootstrap links fluid energy, curvature flux, and sphere geometry while respecting both the derivative count and the null weights. The smallness in the weighted hierarchy controls the errors, even as the integrated shear or matter flux remains large enough to produce focusing. We emphasize that we commute the null structure equations and renormalized null Bianchi equations with derivatives in all directions, compared to the Einstein vacuum equations and other single wave speed models, where only $\angnabla$ is commuted. This is because the fluid energy estimates require control over derivatives (of both the fluid and the geometry) in every direction, since one considers the spacetime divergence of the fluid energy current $\Dfour_\mu\varJ^\mu$. The readers are referred to Section \ref{SS:IntroFluidCurrent} for a more detailed discussion of the estimates of the fluid variables.
	
	\subsubsection{From particle conservation to focusing at every angle}
	\label{SS:IntroParticleFocusing}
	
	The fluid-driven trapped-surface formation argument requires a positive lower bound for $\FAlphaAB$ that survives evolution. We present our heuristic view below.
	To start with, the ansatz \eqref{EQ:IntroFast}, the relations $\enth=\average{H}e^\lnenth$ and $\density+\pressure=n\enth$ inform us that, for small log enthalpy $\lnenth$, we have $\density+\pressure\sim n$. We apply this property to the conservation of mass equation \eqref{1.0.3}. Specifically, we apply the divergence theorem to the fluid current $I^\mu=nv^\mu$ in the region $\mathcal{D}_{\uinf,\delta}^{\uts,2\delta}$ bounded by $\hypersur_{\uinf},\hypersur_{\uts},\hypersurbar_{0}$ and $\mathcal{T}:=\{\ubar=\left(1+2(u-\uinf)\right)\delta,u\in[\uinf,\uts]\}$ as in Figure \ref{FIG:IntroFormation}. That is,
	\begin{align}
		\mathcal{D}_{\uinf,\delta}^{\uts,2\delta}=\left\{(u,\ubar,\sphere^1,\sphere^2)|u\in[\uinf,\uts], \ubar\in\left[0,\left(1+2(u-\uinf)\right)\delta\right]\right\}.
	\end{align}
	In this region, we note that the fluid flows transversally into $\mathcal{D}_{\uinf,\delta}^{\uts,2\delta}$ via $\hypersur_{\uinf},\hypersurbar_{0}$ and $\mathcal{T}$, despite that $\mathcal{T}$ is a timelike manifold (generated by $\frac{1}{\sqrt{8\delta}}(2\delta e_4+e_3)$). This is due to the fact that the fluid is focusing at a relatively large speed\footnote{See previous Section \ref{SS:IntroFluidCurrent} for an introduction to the control of the fluid variables.}. Moreover, since the region $\mathcal{D}_{\uinf,\delta}^{\uts,2\delta}$ is still local in $\ubar$, the estimates in our semi-glboal existence theorem still holds. Through the above argument, we show that, on spherical average, along $\hypersur_{\uts}$, we have $\density+\pressure\sim n\sim\frac{1}{\abs{u_T}^2}$, which has an amplified effect compared to $\density+\pressure\sim n\sim\frac{1}{\abs{\uinf}^2}$ along $\hypersur_{\uinf}$.

    \vspace{1mm}
	Our next step is to pass from this spherical average to a lower bound along every generator, as required by \eqref{EQ:IntroIntegrated}. Here we observe that our established angular estimates control the oscillation of the integrated flux $|u|^2\int\FAlphaAB\,\diff\ubar$. Specifically, \eqref{ES:introvarv} implies:
    \begin{align}
        \abs{\FAlphaAB(u,\ubar,\vartheta)-\FAlphaAB(u,\ubar,\vartheta')}&\leq \BAfluid \abs{\uinf} \quad  \mbox{for all} \quad (u,\ubar,\vartheta),(u,\ubar,\vartheta')\in\mathcal{D}_{\uinf,\delta}^{\uts,2\delta}.
    \end{align}
    Once the desired averaged comparison is derived on $\mathcal{D}_{\uinf,\delta}^{\uts,2\delta}$, the angular and sphere-area low-order terms can be absorbed and we then arrive at the comparison formula used in Lemma~\ref{LE:EstimatesforuFAlpha}:
	\begin{equation}\label{EQ:IntroFluxTarget}
		{\begin{aligned}
				|\uts|^2\int_0^{\delta_T}\FAlphaAB(\uts,\ubarprime,\vartheta)
				\,\diff\ubarprime
				\geq\frac12|\uinf|^2\int_0^\delta
				\FAlphaAB(\uinf,\ubarprime,\vartheta)\,\diff\ubarprime.
		\end{aligned}}
	\end{equation}

	\noindent  The above explanations supply the link between particle conservation and fluid-focusing mechanism at every angle, while allowing general non-spherically symmetric spacetimes.

    \vspace{1mm}

	For either focusing mechanism, the propagated lower bound is finally inserted into \eqref{EQ:IntroIntegrated}. The lapse estimate and the initial expansion bound \eqref{AS:SmallonHbar} control the positive initial contribution, so the outgoing expansion becomes negative on the corresponding target sphere. Meanwhile, the zeroth-order geometric estimate \eqref{ES:LinftyRicciCoefficients4} keeps the incoming expansion close to $-2/|u|$ and strictly negative. The two bounds together give a closed trapped sphere in the constructed spacetime.
	
	\subsection{Other Related Results}
	As mentioned above, after Christodoulou's monumental work \cite{Christodoulou2009}, a systematic approach based on introducing signatures for the short pulse was developed by Klainerman-Rodnianski in \cite{KlainermanRodnianski2012} to simplify and extend \cite{Christodoulou2009} in a finite region. Yu in \cite{Yu2011} extended Klainerman-Rodnianski's approach to the Einstein-Maxwell system. The author in \cite{An2012} introduced a new signature $s_2$, called the signature for decay rates, and extended Klainerman-Rodnianski's results to past null infinity. Based on the hyperbolic estimates in \cite{Christodoulou2009} or \cite{KlainermanRodnianski2012}, Klainerman-Luk-Rodnianski in \cite{KlainermanLukRodnianski2014} found a fully anisotropic mechanism for trapped surface formation. 
	
	In \cite{AnLuk2017} the author and Luk designed a new set of short-pulse hierarchies and proved the first scale-critical trapped surface formation criterion for the Einstein vacuum equations. Using only $s_2$, the signature for decay rates, the author in \cite{An2022} gave a short proof of scale-critical trapped surface formation in the far-field regime and related the short-pulse hierarchies to the peeling properties of gravitational waves. In \cite{AnAthanasiou2022}, the author and Athanasiou further extended the new approach in \cite{An2022} to the Einstein-Maxwell system. With spherically symmetric singular initial data prescribed along the incoming null hypersurface, Li and Liu in \cite{LiLiu2017} studied the Einstein-scalar field system and proved an almost scale-critical trapped surface formation criterion. For connections between trapped surface formation and low-regularity solutions, interested readers are also referred to \cites{LukRodnianski2017,LukRodnianski2020} by Luk-Rodnianski and \cite{LukMoschidis2023} by Luk-Moschidis. Very recently, the methods in \cite{An2022} by the first author and in \cite{AnAthanasiou2022} by the first author and Athanasiou have also been extended to the Einstein-Yang-Mills system, the Einstein-scalar field system, the Einstein-Vlasov system, the Einstein-Maxwell-charged scalar field system, the Einstein-Maxwell-Klein-Gordon system, by Athanasiou-Mondal-Yau \cite{AthanasiouMondalYau2023}, Zhao-Hilditch-Kroon \cite{ZhaoHilditchKroon2023}, Athanasiou-Mondal-Yau \cite{AthanasiouMondalYau2025}, Shen-Wan \cite{ShenWan20251}, Giorgi-Shen-Wan \cite{GiorgiShenWan2025}, respectively. The above evolutionary results all employ double-null foliation. We also note the result \cite{ChenKlainerman2024TS} by Chen-Klainerman using geodesic foliation. For constructing Cauchy initial data leading to trapped surface formation, interested readers are referred to \cite{LiYu2015} by Li-Yu, \cite{ShenWan20252} by Shen-Wan via gluing arguments, and \cite{ChenKlainerman2025} by Chen-Klainerman via solving the constraint equations. We also mention the recent encouraging progress of the Einstein-Euler system: Guo-Hadzic-Jang \cite{GuoHadzicJang2023} on constructing the first naked singularity solution for this system, Miao-Shahshahani \cite{MiaoShahshahani2024} on free-boundary for hard-phase model, Song \cite{Song2025} on weak null singularities, Li-Li-Zhu \cite{LiLiZhu2026} on trapped surface formation with shear-tensor dominant dynamics, and Anderson-Luk \cite{AndersonLuk} on shock formation for this system without symmetry. 
	
	\subsection{Organization of the Paper}

		Section~2 fixes the geometric and fluid conventions and records the null structure, Bianchi, and Euler equations. Section~3 introduces the weighted norms, constructs the characteristic data, and states the main theorem and bootstrap assumptions. Section~4 records the commutation formulas, and Section~\ref{S:basic} establishes the preliminary metric, transport, Sobolev, and elliptic estimates. Section~6 proves the fluid energy estimates. Sections~7 and~8 treat the Ricci coefficients below and at top order. Section~\ref{S:Curvature} proves the renormalized Weyl-curvature estimates, and Section~\ref{S:Formation} gives the two trapped-surface formation arguments. The appendix collects the notations.
	
	\subsection{Notations}\label{SS:Notations}
	We use the following notations throughout the paper:
	\begin{itemize}
		\item Greek ``space-time" indices take on the values $0,1,2,3$, while the Latin ``spatial" indices take on the values $1,2,3$. 
		\item Geometric coordinates indices $I,J$ take values $1,2,3,4$, and $A,B,C$ take values $1,2$.
		\item We use the Einstein summation convention.
		\item $\antisymmetic_{\alpha\beta\gamma\delta}$ and $\antisymmetic_{ijk}$ denotes the fully antisymmetric symbol normalized by $\antisymmetic_{0123}=1$ and $\antisymmetic_{123}=1$.
		\item We raise and lower indices for tensorfield with respect to spacetime metric $\gfour$ and its inverse. For example, $V^\alpha=(\gfour^{-1})^{\alpha\beta}V_\beta$, $V_\beta=\gfour_{\alpha\beta}V^\alpha$.
		\item We denote the schematic spatial partial derivatives by $\p$ and the schematic space-time partial derivatives by $\pfour$.
		\item We denote $A\lesssim B$ as $A\leq C\cdot B$ for some universal constant $C$.
	\end{itemize}
	We refer readers to Appendix \ref{AP:notations} for more notations that we use throughout the article.
	
	\subsection{Acknowledgments}
	XA is supported by MOE Tier 1 grant A-8002933-00-00.
	
	\section{Settings}
	
	\subsection{Definitions of the Spacetime}\label{SS:despacetime}
	Let $(\region,\gfour)$ be a $3+1$-dimensional spacetime. Let the speed of light be $1$.
	$\Dfour$ denotes the Levi-Civita connection with respect to $\gfour$. We let $\Riemfour{\alpha}{\beta}{\gamma}{\delta}$ denote the Riemann curvature tensor of $\gfour$, the Ricci curvature tensor $\Ricfour{\alpha}{\beta}:=(\gfour^{-1})^{\gamma\delta}\Riemfour{\gamma}{\alpha}{\delta}{\beta}$, and the scalar curvature $\Scalar:=(\gfour^{-1})^{\alpha\beta}\Ricfour{\alpha}{\beta}$. We use the notation that $\langle\Dfour_X\Dfour_YW-\Dfour_Y\Dfour_XW,Z\rangle=\Riemfour{Z}{W}{X}{Y}+\langle\Dfour_{[X,Y]}W,Z\rangle$ where $X,Y,W,Z$ are vectorfields, $[\cdot,\cdot]$ is the Lie bracket and $\langle\cdot,\cdot\rangle:=\gfour(\cdot,\cdot)$. 
	
	If $\upxi$ is a space-time tensor, then $\sphereproject\upxi$ denotes its $\gfour$-orthogonal projection onto $\suubar$. If both $V$ and $\upxi$ are $\suubar$-tangent, we denote the covariant derivative on $\suubar$ with respect to the induced metric $\gsphere$ as $\angnabla_V\upxi$. Then we define $\angD_V\upxi:=\sphereproject\Dfour_V\upxi$ where $V$ is a vector and $\Dfour_V\upxi$ is the covariant derivative of $\upxi$ in the $V$ direction. Note that $\angD_V\upxi:=\angnabla_V\upxi$ when both $V$ and $\upxi$ are $\suubar$-tangent.
	\subsection{Definitions of the Fluid Variables and Related Quantities}\label{SS:deffluid}
	\subsubsection{The Basic Fluid Variables} \label{SSS:thebasicfluidvariables}
	
	The fluid velocity $v^\alpha$ is a future-directed four-vector normalized by $\gfour(v,v)=-1$. $\pressure$ denotes the pressure, $\density$ denotes the proper energy density, $n$ denotes the proper number density, $\ent$ denotes the entropy per particle, $\temp$ denotes the temperature, and 
	\begin{equation}\label{DE:enthalpy}
		\enth=(\density+\pressure)/n
	\end{equation}
	is the enthalpy per particle. Thermodynamics supplies the following laws:
	\begin{align}
		\enth&=\left.\frac{\p\density}{\p n}\right\vert_\ent,&
		\temp&=\left.\frac{1}{n}\frac{\p\density}{\p\ent}\right\vert_n,&
		\diff\enth&=\frac{\diff \pressure}{n}+\temp\diff\ent,
	\end{align}
	where $\left.\frac{\p}{\p n}\right\vert_\ent$ denotes partial differentiation with respect to $n$ at fixed $\ent$ and $\left.\frac{\p}{\p \ent}\right\vert_n$ denotes partial differentiation with respect to $\ent$ at fixed $n$.
	
	We assume an equation of state of the form $\pressure=\pressure(\enth,\ent)$. The speed of sound $\speed$ is defined as follows:
	\begin{equation}
		\speed:=\sqrt{\left.\frac{\p \pressure}{\p\density}\right\vert_\ent}.
	\end{equation}
	We assume that $\speed^2\leq \speed_0^2<1$.
	\subsubsection{$v$-orthogonal Vorticity}\label{SSS:v-ortho}
	\begin{definition}[The $v$-orthogonal vorticity of a one form]\label{vorthogonalvorticity}
		Given a space-time one-form $V$, we define the corresponding $v$-orthogonal (with respect to the Minkowski metric) vorticity vectorfield as follows:
		\begin{equation}
			\text{vort}^\alpha(V):=-\antisymmetic^{\alpha\beta\gamma\delta}v_\beta\Dfour_\gamma V_\delta.
		\end{equation}
		
	\end{definition}
	
	\begin{definition}[The $v$-orthogonal vorticity vector field]\label{vorticity}
		We define the vorticity vector field $\vort^\alpha$ as follows:
		\begin{align}
			\vort^\alpha:=\text{vort}^\alpha(\enth v)=-\antisymmetic^{\alpha\beta\gamma\delta}v_\beta\Dfour_\gamma(\enth v_\delta).
		\end{align}
	\end{definition}
	\begin{definition}[Logarithmic enthalpy]\label{DE:loglnenth}
		Let $\average{H}>0$ be a fixed constant value of the background enthalpy. We define the logarithmic enthalpy $\lnenth$ as follows:
		\begin{equation}\label{DE:logenthalpy}
			\lnenth:=\ln(\enth/\average{H}).
		\end{equation}
	\end{definition}
	\subsection{The Einstein-Euler Equations}\label{SS:EinsteinEuler}
	We consider the Einstein-Euler equations as follows:
	\begin{subequations}\label{EQ:EE}
		\begin{align}
			\Ricfour{\mu}{\nu}-\frac{1}{2}\Scalar\gfour_{\mu\nu}=&~\EMT_{\mu\nu},\\
			\label{EQ:DivEMT}\Dfour_{\mu}\EMT^{\mu\nu}=&~0,\\
			\label{EQ:DivCurrent}\Dfour_{\mu}I^\mu=&~0,
		\end{align}
	\end{subequations}
	where $I^\mu:=nv^\mu$ is the fluid current and $\EMT_{\mu\nu}$ is the energy-momentum tensor defined as follows:
	\begin{align}
		\EMT_{\mu\nu}:=(\density+\pressure)v_{\mu}v_{\nu}+\pressure\gfour_{\mu\nu}.
	\end{align}
	\subsection{Null Frames}\label{SS:NullFrame}
	First, we define the eikonal functions $u$ and $\ubar$ to be the solutions of the following eikonal equations, respectively:
	\begin{align}\label{EQ:eikonalequations}
		(\gfour^{-1})^{\alpha\beta}\p_\alpha u\p_\beta u&=0,&
		(\gfour^{-1})^{\alpha\beta}\p_\alpha \ubar\p_\beta \ubar&=0.
	\end{align}
	
	We employ double null coordinates $(u,\ubar,\theta^1,\theta^2)$. 
	We define our geometric null frame and its associated vectorfields based on $u$ and $\ubar$ as follows:\\
	
	\textbf{Geometric time $\uptau$}
	\begin{align}\label{DE:Deftau}
		\uptau:=u+\ubar.
	\end{align}
	
	\textbf{Geometric Radius $\rgeo$}
	\begin{align}\label{DE:Defr}
		\rgeo:=\ubar-u.
	\end{align}
	
	\textbf{Geodesic vectorfields $\Lgeo^\alpha$ and $\uLgeo^\alpha$}
	\begin{align}\label{D:geo}
		\Lgeo^\alpha&:=-2(\gfour^{-1})^{\alpha\beta}\p_\beta u,&
		\uLgeo^\alpha&:=-2(\gfour^{-1})^{\alpha\beta}\p_\beta \ubar.
	\end{align}
	
	\textbf{Null lapse $\lapse$}\label{DE:lapse}
	\begin{align}\label{D:nullapse}
		2\lapse^{-2}:=-\gfour(\Lgeo,\uLgeo).
	\end{align}
	
	\textbf{Principal null vectorfields $e_3$ and $e_4$}
	\begin{align}\label{D:e3e4}
		e_4&:=\lapse\Lgeo,&
		e_3&:=\lapse\uLgeo.
	\end{align}
	
	\textbf{Rescaled null vectorfields $\Lunit$ and $\uLunit$}
	\begin{align}\label{D:LunituLunit}
		\Lunit&:=\lapse^2\Lgeo,&
		\uLunit&:=\lapse^2\uLgeo,
	\end{align}
	
	We denote the level sets of $u$ and $\ubar$ by $\hypersur_u$ and $\hypersurbar_{\ubar}$, respectively. Notice that $\hypersur_u$ and $\hypersurbar_{\ubar}$ are null hypersurfaces with principal null generators $e_4$ and $e_3$, respectively.
	We denote the $2$-sphere $\suubar$ as follows:
	\begin{align}
		\suubar:=\hypersur_u\bigcup\hypersurbar_{\ubar}.
	\end{align}
	
	We denote $\region$ the space-time region, where $\region=\bigcup\limits_{u}\hypersur_{u}=\bigcup\limits_{\ubar}\hypersurbar_{\ubar}$. 
	We define $\gsphere$ to be the inverse of $\gsphere^{-1}$, and $\gsphere^{-1}$ as follows:
	\begin{align}\label{EQ:decompositionofmetric}
		(\gsphere^{-1})^{\alpha\beta}:=(\gfour^{-1})^{\alpha\beta}+\frac{1}{2}e_4^\alpha e_3^\beta+\frac{1}{2}e_3^\alpha e_4^\beta.
	\end{align}
	
	It is easy to check that $\gsphere$ is the induced metric of $\gfour$ on $\suubar$. We fix a pair of locally defined unit orthonormal spherical vectorfields\footnote{In the rest of the paper, we automatically sum over $A$ if there are two $A$'s in the expression.} on $\suubar$ by $\{e_A\}_{A=1,2}$ such that $(\gsphere^{-1})^{\alpha\beta}=\sum\limits_{A=1,2}e_A^\alpha e_A^\beta$. We call $\{e_I\}_{I=1}^4$ a null frame for the geometry. 
	
	By the construction above, we have the following identities:
	\begin{subequations}\label{EQ:innerproduct}
		\begin{align}
			\gfour(e_4,e_4)=&~0,&\gfour(e_3,e_3)=&~0,\\
			\gfour(e_3,e_4)=&~-2,&\gfour(e_A,e_B)=&~\gsphere_{AB},\\
			\gfour(e_A,e_4)=&~0,&\gfour(e_A,e_3)=&~0.
		\end{align}
	\end{subequations}
	
	In these coordinates, we have:
	\begin{align}\label{DE:Coordinatese3e4}
		e_3=&~\lapse^{-1}\left(\frac{\p}{\p u}+\Lambda^A\frac{\p}{\p\theta^A}\right), &
		e_4=&~\lapse^{-1}\frac{\p}{\p\ubar},
	\end{align}
	for some $\Lambda^A$ such that $\Lambda^A=0$ on $\hypersurbar_0$. We choose $u,\ubar$ such that $\lapse=1$ on $\hypersur_{\uinf}\cup\hypersurbar_{0}$. We can express the spacetime metric $\gfour$ relative to the geometric coordinates $(u,\ubar,\sphere^1,\sphere^2)$:
	\begin{align}\label{EQ:gfouringeometriccoordinates}
		\gfour=-2\lapse^2\left(\diff u\otimes\diff\ubar+\diff \ubar\otimes\diff u\right)+\gsphere_{AB}\left(\diff\sphere^A-\Lambda^A\diff \ubar\right)\otimes\left(\diff\sphere^B-\Lambda^B\diff\ubar\right).
	\end{align}
	
	We define the geometric components of $v$ as follows:
	\begin{align}\label{EQ:devi}
		v_I:=&~\gfour(v,e_I),&
		I=&~1,2,3,4.
	\end{align}
	By the fact that $v$ is a future-directed timelike vectorfield, $v_3,v_4<0$.
	Then we have the decomposition of $v$ as follows:
	\begin{subequations}
		\begin{align}
			v=&~-v_\Timelike\Timelike+v_\spherenormal\spherenormal+(\gsphere^{-1})^{AB}v_Ae_B,\\
			v=&~-\frac{1}{2}v_{4}e_3-\frac{1}{2}v_{3}e_4+(\gsphere^{-1})^{AB}v_Ae_B,
		\end{align}
	\end{subequations}
	where $v_\xi:=\gfour(v,\xi)$ for any vectorfield $\xi$.
	The normalization $\gfour(v,v)=-1$ is equivalent to the following:
	\begin{subequations}
		\begin{align}
			-v_\Timelike^2+v_\spherenormal^2+\gsphere^{AB}v_Av_B=&~-1,\\
			-v_{3}v_{4}+\gsphere^{AB}v_Av_B=&~-1
		\end{align}
	\end{subequations}
	
	\subsection{Geometric Quantities and Equations}
	\subsubsection{Notations}\label{SSS:basicnotations}
	We set $\phi^{(1)}\cdot\phi^{(2)}$ to stand for a contraction of the tensor product of $\phi^{(1)}$ and $\phi^{(2)}$ with respect to the metric $\gsphere$ on $\suubar$. We also denote the following:
	\begin{subequations}
		\begin{align}
			(\phi^{(1)}\tensorproduct\phi^{(2)})_{AB}:=&~\phi^{(1)}_A\phi^{(2)}_B+\phi^{(1)}_B\phi^{(2)}_A-\gsphere_{AB}\left(\phi^{(1)}\cdot\phi^{(2)}\right),& &\text{ for one-forms }\phi^{(1)}_A, \phi^{(2)}_A\\
			\phi^{(1)}\wedge\phi^{(2)}:=&~\antisymmetic^{AB}\left(\gsphere^{-1}\right)^{CD}\phi^{(1)}_{AC}\phi^{(2)}_{BD},& &\text{ for symmetric two-tensors } \phi^{(1)}_{AB}, \phi^{(2)}_{AB}
		\end{align}
	\end{subequations}
	\begin{definition}[Spherical projection operators]\label{DE:sphereprojection}
		We define spherical projection operator $\sphereproject$ as follows:
		\begin{align}
			\sphereproject^\alpha_\beta:=\delta^\alpha_\beta+\frac{1}{2}e_4^\alpha(e_3)_\beta+\frac{1}{2}e_3^\alpha(e_4)_\beta
		\end{align}
		Then, for any vectorfield/one-form $V$, we define $\angV:=\sphereproject V$ to be the spherical projection of $V$.
		We define the spherical projection for tensorfields in a similar fashion.
	\end{definition}
	\subsubsection{Connection coefficients}\label{SS:connectioncoefficients}
	With the construction of the eikonal functions, in this section we define eikonal quantities, including connection coefficients and curvature components.
	\begin{definition}[Deformation tensors]\label{DE:DeformationTensor}
		We define the deformation tensor of a vectorfield $V$ as follows:
		\begin{align}\label{D:deformationtensor}
			\deform{V}_{\alpha\beta}:=\Dfour_\alpha V_{\beta}+\Dfour_{\beta}V_{\alpha}
		\end{align}
	\end{definition}
	\begin{definition}[Connection coefficients]\label{D:Defconnectioncoeff}
		Given a 2-sphere $\suubar$ and an arbitrary frame $\{e_A\}_{A=1,2}$ tangent to it (in particular, for $\{e_A\}_{A=1,2}$ defined in Section \ref{SS:NullFrame}), we define the connection coefficients relative to the null frame in Section \ref{SS:NullFrame} as follows:
		\begin{subequations}
			\begin{align}
				\upchi_{AB}&=\gfour(\Dfour_Ae_4,e_B),&
				\uchi_{AB}&=\gfour(\Dfour_Ae_3,e_B),\\
				\upeta_A&=-\frac{1}{2}\gfour(\Dfour_3e_A,e_4),&
				\ueta_A&=-\frac{1}{2}\gfour(\Dfour_4e_A,e_3),\\
				\omega&=-\frac{1}{4}\gfour(\Dfour_4e_3,e_4),&
				\uomega&=-\frac{1}{4}\gfour(\Dfour_3e_4,e_3),\\
				\upzeta_A&=\frac{1}{2}\gfour(\Dfour_Ae_4,e_3).
			\end{align}
		\end{subequations}
		Moreover, we define the trace and traceless part of a $\suubar$-tangent two-tensor $\xi$ as follows:
		\begin{align}
			\gtr\xi:=&~(\gsphere^{-1})^{AB}\xi_{AB},&
			\hat{\xi}:=&~\xi-\frac{1}{2}\gtr\xi \gsphere.
		\end{align}
	\end{definition}	
	\begin{proposition}[Connection coefficients and relationship between tensors]\label{PR:PropConnection}
		By Definition \ref{D:Defconnectioncoeff}, the covariant derivatives of the geometric frame can be expressed as follows:
		\begin{subequations}\label{E:propconnection}
			\begin{align}
				\label{EQ:DLL}\Dfour_{4}e_4&=-2\omega e_4,&
				\Dfour_{3}e_3&=-2\uomega e_3,\\
				\Dfour_Ae_4&=\upchi_{AB}e_B-\upzeta_Ae_4,&
				\Dfour_Ae_3&=\uchi_{AB}e_B+\upzeta_Ae_3,\\
				\Dfour_{3}e_4&=2\upeta_A e_A+2\uomega e_4,&
				\label{EQ:DLbarL}\Dfour_{4}e_3&=2\ueta_Ae_A+2\omega e_3,\\		
				\Dfour_{4} e_A&=\angD_{4} e_A+\ueta_Ae_4,&
				\Dfour_{3}e_A&=\angD_{3}e_A+\upeta_Ae_3,\\
				\Dfour_Ae_B&=\angnabla_Ae_B+\frac{1}{2}\uchi_{AB}e_4+\frac{1}{2}\upchi_{AB}e_3,\\
				\Dfour_\Timelike\Timelike&=\angnabla_A(\ln\lapse) e_A+\frac{1}{2}(\uomega-\omega)e_4+\frac{1}{2}(\omega-\uomega)e_3,\\
				\Dfour_\spherenormal\spherenormal&=-\angnabla_A(\ln\lapse) e_A-\frac{1}{2}(\omega+\uomega)e_4+\frac{1}{2}(\omega+\uomega)e_3.
			\end{align}
		\end{subequations}
		Moreover,
		\begin{subequations}
			\begin{align}
				\label{EQ:omegauomega}\omega&=-\frac{1}{2}e_4\ln \lapse,&
				\uomega&=-\frac{1}{2}e_3\ln \lapse,\\
				\label{EQ:etaupeta}\upeta_A&=\upzeta_A+\angnabla_A\ln\lapse,&
				\ueta_A&=-\upzeta_A+\angnabla_A\ln\lapse.&
			\end{align}
		\end{subequations}
		In addition, 
		\begin{align}\label{EQ:Lieg}
			\angLie_4\gfour&=\angLie_4\gsphere=2\upchi,&
			\angLie_3\gfour&=\angLie_3\gsphere=2\uchi,
		\end{align}
		where $\Lie$ is the standard Lie derivative and $\angLie:=\sphereproject\Lie$.
	\end{proposition}
	\subsubsection{Curvature components}
	\begin{definition}[Standard Weyl tensor]\label{DE:Weyl}
		We define the Weyl tensor $\Weyl_{\alpha\beta\gamma\delta}$ as follows:
		\begin{align}
			\Weyl_{\alpha\beta\gamma\delta}:=\Riemfour{\alpha}{\beta}{\gamma}{\delta}+\frac{1}{2}\left(\Ricfour{\alpha}{\delta}\gfour_{\beta\gamma}+\Ricfour{\beta}{\gamma}\gfour_{\alpha\delta}-\Ricfour{\alpha}{\gamma}\gfour_{\beta\delta}-\Ricfour{\beta}{\delta}\gfour_{\alpha\gamma}\right)+\frac{1}{6}\Scalar\left(\gfour_{\alpha\gamma}\gfour_{\beta\delta}-\gfour_{\alpha\delta}\gfour_{\beta\gamma}\right).
		\end{align}
	\end{definition}
	\begin{definition}[Hodge dual]\label{DE:HodgeDual}	
		Let $\antisymmetic_{\alpha\beta\gamma\delta}$ be the volume element of $\region$. We define the spacetime (left) Hodge dual as follows:
		\begin{align}
			\HodgeLeftDual\xi_{\alpha\beta\gamma\delta}:=&~\antisymmetic_{\alpha\beta}^{\ \ \mu\nu}\xi_{\mu\nu\gamma\delta}.
		\end{align}
		We define the (left) Hodge dual on $\suubar$-tangent vectorfields and symmetric traceless 2-tensors as follows:
		\begin{subequations}
			\begin{align}
				\HodgeLeftDual\xi_A:=&~\antisymmetic_{AB}\xi_B,\\
				\HodgeLeftDual\xi_{AB}:=&~\antisymmetic_{AC}\xi_{BC},
			\end{align}
		\end{subequations}
		where $\antisymmetic_{AB}$ is the volume element of $\suubar$.
	\end{definition}
	\begin{definition}[Projected Ricci Curvature]\label{DE:ProjectedRicci}
		We recall the spherical projection operator $\sphereproject$ in Definition \ref{DE:sphereprojection}. Then, we denote \begin{subequations}
			\begin{align}
				\angRicfour_{\alpha\beta}:=&~\sphereproject_\alpha^{\mu}\sphereproject_\beta^{\nu}\Ricfour{\mu}{\nu},\\ \angRicfour_{4\alpha}:=&~\sphereproject_\alpha^{\mu}\Ricfour{4}{\mu}&
				\angRicfour_{3\alpha}:=&~\sphereproject_\alpha^{\mu}\Ricfour{3}{\mu}
			\end{align}
		\end{subequations}
	\end{definition}
	\begin{definition}[Curvature components]\label{DE:CurvatureComponents}
		We denote by $\GaussCurvature$ the Gaussian curvature. We define the following null curvature components:
		\begin{subequations}
			\begin{align}
				\AlphaAB_{AB}=&~\Weyl_{A4B4},& \Alphabar_{AB}=&~\Weyl_{A3B3},\\
				\BetaA_A=&~\frac{1}{2}\Weyl_{A434},& \Betabar_A=&~\frac{1}{2}\Weyl_{A334},\\
				\Rho=&~\frac{1}{4}\Weyl_{4343}, &\Rhodual=&~\frac{1}{4}\Weyldual_{4343}.
			\end{align}
		\end{subequations}
		Moreover, we define the Ricci curvature tensor components as follows:
		\begin{subequations}
			\begin{align}
				\FAlphaAB:=&~\frac{1}{2}\Ricfour{4}{4},&
				\FAlphabar:=&~\frac{1}{2}\Ricfour{3}{3},\\
				\FBetaA_A:=&~\frac{1}{2}\angRicfour_{4A},&
				\FBetabar_A:=&~\frac{1}{2}\angRicfour_{3A},\\
				\FRho:=&~\frac{1}{2}\Ricfour{3}{4},&
				\FRhodual_{AB}:=&~\hangRicfour_{AB}.
			\end{align}
		\end{subequations}
	\end{definition}

	\begin{proposition}[Preliminary results for Weyl tensor components]\label{PR:Weyl}
		By straightforward computations, we have the following expressions of the Weyl tensor components:
		\begin{subequations}
			\begin{align}
				\Weyl_{AB34}=&~2\Rhodual\antisymmetic_{AB},\\
				\Weyl_{ABC3}=&~\antisymmetic_{AB}\HodgeLeftDual\Betabar_{C},&
				\Weyl_{ABC4}=&~-\antisymmetic_{AB}\HodgeLeftDual\BetaA_{C},\\
				\Weyl_{A3B4}=&~-\Rho\gsphere_{AB}+\Rhodual\antisymmetic_{AB},&
				\Weyl_{ABCD}=&~\antisymmetic_{AB}\antisymmetic_{CD}\Rho,
			\end{align}
		\end{subequations}
		In addition, for the spacetime dual of $W$, we have:
		\begin{subequations}
			\begin{align}
				\HodgeLeftDual\Weyl_{A3B3}=&~\HodgeLeftDual\Alphabar_{AB},&
				\HodgeLeftDual\Weyl_{A4B4}=&~-\HodgeLeftDual\AlphaAB_{AB},\\
				\HodgeLeftDual\Weyl_{A334}=&~2\HodgeLeftDual\Betabar_{A},&
				\HodgeLeftDual\Weyl_{A434}=&~-2\HodgeLeftDual\BetaA_{A},\\
				\HodgeLeftDual\Weyl_{3434}=&~4\Rhodual,&
				\HodgeLeftDual\Weyl_{AB34}=&~-2\Rho\antisymmetic_{AB},\\
				\HodgeLeftDual\Weyl_{ABC3}=&~-\antisymmetic_{AB}\Betabar_{C},&
				\HodgeLeftDual\Weyl_{ABC4}=&~-\antisymmetic_{AB}\BetaA_{C},\\
				\HodgeLeftDual\Weyl_{A3B4}=&~-\Rhodual\gsphere_{AB}-\Rho\antisymmetic_{AB},&
				\HodgeLeftDual\Weyl_{ABCD}=&~\antisymmetic_{AB}\antisymmetic_{CD}\Rhodual.
			\end{align}
		\end{subequations}
	\end{proposition}
	\begin{proposition}[Preliminary results for Ricci curvature components]\label{PR:Ricci0}
		We compute the Ricci curvature components in the null frame as follows:
		\begin{subequations}
			\begin{align}
				\Scalar=&~-\Ricfour{3}{4}+\gtr\angRicfour=\density-3\pressure,\\
				2\FAlphaAB=&~\Ricfour{4}{4}=\EMT_{44}=(\density+\pressure)(v_4)^2,\\
				2\FAlphabar=&~\Ricfour{3}{3}=\EMT_{33}=(\density+\pressure)(v_3)^2,\\
				2\FBetaA=&~\Ricfour{4}{A}=\EMT_{4A}=(\density+\pressure)v_4v_A,\\
				2\FBetabar=&~\Ricfour{3}{A}=\EMT_{3A}=(\density+\pressure)v_3v_A,\\
				\frac{1}{2}(\Scalar+2\FRho)\gsphere_{AB}+\FRhodual_{AB}=&~\Ricfour{A}{B}=\EMT_{AB}+\frac{1}{2}\Scalar\gsphere_{AB}=(\density+\pressure)v_Av_B+\frac{1}{2}(\density-\pressure)\gsphere_{AB},\\
				2\FRho=&~\Ricfour{3}{4}=\EMT_{34}-\Scalar=(\density+\pressure)v_3v_4-(\density-\pressure)\\
				\notag=&~\gtr\angRicfour-\Scalar=(\density+\pressure)\angv\cdot\angv+2\pressure.
			\end{align}
		\end{subequations}
		Moreover, the Weyl tensor components satisfy the following identities:
		\begin{subequations}
			\begin{align}
				\Weyl_{A4B4}=&~\Riemfour{4}{A}{4}{B}-\frac{1}{2}\Ricfour{4}{4}\gsphere_{AB},&\Weyl_{A3B3}=&~\Riemfour{3}{A}{3}{B}-\frac{1}{2}\Ricfour{3}{3}\gsphere_{AB},\\
				\Weyl_{A434}=&~\Riemfour{A}{4}{3}{4}-\Ricfour{A}{4},&
				\Weyl_{A334}=&~\Riemfour{A}{3}{3}{4}+\Ricfour{A}{3},\\
				\Weyl_{4343}=&~\Riemfour{4}{3}{4}{3}-2\Ricfour{3}{4}-\frac{2}{3}\Scalar=2\Riemfour{A}{3}{4}{A}-\frac{2}{3}\Scalar,&
				\Weyl_{AB34}=&~\Riemfour{A}{B}{3}{4},\\
				\Weyl_{ABC3}=&~\Riemfour{A}{B}{C}{3}+\frac{1}{2}\left(\Ricfour{A}{3}\gsphere_{BC}-\Ricfour{B}{3}\gsphere_{AC}\right),&
				\Weyl_{ABC4}=&~\Riemfour{A}{B}{C}{4}+\frac{1}{2}\left(\Ricfour{A}{4}\gsphere_{BC}-\Ricfour{B}{4}\gsphere_{AC}\right),\\
				\Weyl_{A3B4}=&~\Riemfour{A}{3}{B}{4}+\frac{1}{2}\left(2\Ricfour{A}{B}-\Ricfour{3}{4}\gsphere_{AB}\right)-\frac{1}{3}\Scalar\gsphere_{AB}, &
				\HodgeLeftDual\Weyl_{3434}=&~\HodgeLeftDual\Riem_{3434},
			\end{align}
			\begin{align}
				\Weyl_{ABCD}=&~\Riemfour{A}{B}{C}{D}+\frac{1}{2}\left(\Ricfour{A}{D}\gsphere_{BC}+\Ricfour{B}{C}\gsphere_{AD}-\Ricfour{A}{C}\gsphere_{BD}-\Ricfour{B}{D}\gsphere_{AC}\right)+\frac{1}{6}\Scalar(\gsphere_{AC}\gsphere_{BD}-\gsphere_{AD}\gsphere_{BC})
			\end{align}
		\end{subequations}
	\end{proposition}

	\begin{proposition}[Null structure equations verified by geometric quantities]\label{PR:NullStructure1}
		By the geometric construction, the following equations hold:\\
		
		\noindent	\textbf{Transport equations for $\upchi$}
		\begin{subequations}
			\begin{align}
				\label{EQ:e4upchi}e_4\gtr\upchi+\frac{1}{2}(\gtr\upchi)^2=&~-\sgabs{\hchi}^2-2\omega\gtr\upchi-2\FAlphaAB,\\
				\label{EQ:angDhchi}\angD_4\hchi+(\gtr\upchi)\hchi=&~-2\omega\hchi-\AlphaAB,\\	
				\label{e3gtrupchi}e_3\gtr\upchi+\frac{1}{2}\gtr\uchi\gtr\upchi=&~2\angdiv\upeta+2\uomega\gtr\upchi-\hchi\cdot\huchi+2\sgabs{\upeta}^2\\
				\notag&+2\Rho+\frac{1}{3}\Scalar,\\
				\label{EQ:angD3hchi}\angD_3\hchi+\frac{1}{2}(\gtr\uchi)\hchi=&~-\frac{1}{2}\gtr\upchi\huchi+\angnabla\tensorproduct\upeta+2\uomega\hchi+\upeta\tensorproduct\upeta+\FRhodual.
			\end{align}
		\end{subequations}
		\textbf{Transport equations for $\uchi$}
		\begin{subequations}
			\begin{align}
				\label{e3uchi}e_3\gtr\uchi+\frac{1}{2}(\gtr\uchi)^2=&~-\sgabs{\huchi}^2-2\uomega\gtr\uchi-2\FAlphabar,\\
				\label{e4huchi}\angD_3\huchi+(\gtr\uchi)\huchi=&~-2\uomega\huchi-\Alphabar,\\
				\label{Luchi}e_4\gtr\uchi+\frac{1}{2}\gtr\upchi\gtr\uchi=&~2\angdiv\ueta+2\omega\gtr\uchi-\hchi\cdot\huchi+2\sgabs{\ueta}^2+2\Rho+\frac{1}{3}\Scalar,\\
				\label{Lhuchi}\angD_4\huchi+\frac{1}{2}(\gtr\upchi)\huchi=&~-\frac{1}{2}\gtr\uchi\hchi+\angnabla\tensorproduct\ueta+2\omega\huchi+\ueta\tensorproduct\ueta+\FRhodual.
			\end{align}
		\end{subequations}
		\textbf{Div-curl systems for $\hchi,\huchi$}
		\begin{subequations}
			\begin{align}
				\label{angdivhchi}\angdiv\hchi+\hchi\cdot\upzeta=&~\frac{1}{2}(\angnabla\gtr\upchi+\upzeta\gtr\upchi)-\BetaA,\\
				\label{angdivhuchi}\angdiv\huchi-\huchi\cdot\upzeta=&~\frac{1}{2}(\angnabla\gtr\uchi-\upzeta\gtr\uchi)+\Betabar.
			\end{align}
		\end{subequations}
		\textbf{Equations for $\upeta,\ueta,\upzeta$}
		\begin{subequations}
			\begin{align}
				\label{Lupeta}\angD_4\upeta&=(\ueta-\upeta)\cdot\chi-\BetaA-\FBetaA,\\
				\label{uLueta}\angD_3\ueta&=(\upeta-\ueta)\cdot\uchi+\Betabar-\FBetabar,\\
				\label{Lzeta}\angD_4\upzeta&=2\angnabla\omega+2\omega(\ueta+\upzeta)+\upchi\cdot(\ueta-\upzeta)-\BetaA-\FBetaA,\\
				\label{uLzeta}\angD_3\upzeta&=-2\angnabla\uomega+2\uomega(\upzeta-\upeta)-\uchi\cdot(\ueta+\upzeta)-\Betabar+\FBetabar,\\
				\label{curleta}\angcurl\upeta&=-\angcurl\ueta=\frac{1}{2}\huchi\wedge\hchi+\Rhodual.	
			\end{align}
		\end{subequations}
		\textbf{Equations for $\omega,\uomega,\GaussCurvature$}
		\begin{subequations}
			\begin{align}
				\label{Luomega}\angD_4\uomega=&~2\omega\uomega+2\sgabs{\upzeta}^2-\frac{1}{2}\sgabs{\ueta}^2+\frac{1}{2}\Rho+\frac{1}{2}\FRho+\frac{1}{12}\Scalar,\\
				\label{Lbaromega}\angD_3\omega=&~2\omega\uomega+2\sgabs{\upzeta}^2-\frac{1}{2}\sgabs{\upeta}^2+\frac{1}{2}\Rho+\frac{1}{2}\FRho+\frac{1}{12}\Scalar,\\
				\label{KandRho}\GaussCurvature=&~-\Rho+\frac{1}{2}\hchi\cdot\huchi-\frac{1}{4}\gtr\upchi\gtr\uchi+\FRho+\frac{1}{3}\Scalar.
			\end{align}
		\end{subequations}
	\end{proposition}
	\subsection{Null Bianchi Equations}\label{SS:Bianchi}
	\begin{definition}[First derivative of the Weyl tensor]\label{DE:DerivativeofWeyl}
		We define the divergence of the Weyl tensor $\DWeyl_{\beta\gamma\delta}$ and its right Hodge dual as follows:
		\begin{subequations}
			\begin{align}
				\DWeyl_{\beta\gamma\delta}:=&~\Dfour^{\alpha}\Weyl_{\alpha\beta\gamma\delta},&
				\DWeylStar_{\beta\gamma\delta}:=&~\Dfour^{\alpha}\Weyldual_{\alpha\beta\gamma\delta}
			\end{align}
		\end{subequations}
		We note that $C_{\beta\gamma\delta}:=2\DWeyl_{\beta\gamma\delta}$ is the so-called Cotton tensor.
		By direct computation, we have:
		\begin{align}
			\DWeyl_{\beta\gamma\delta}=&~-\frac{1}{2}\Dfour_\delta\Ricfour{\beta}{\gamma}+\frac{1}{2}\Dfour_\gamma\Ricfour{\beta}{\delta}-\frac{1}{12}\left(\Dfour_\gamma\Scalar\gfour_{\beta\delta}-\Dfour_\delta\Scalar\gfour_{\beta\gamma}\right),\\
			\DWeylStar_{\beta\gamma\delta}=&~\DWeyl_{\beta\mu\nu}\antisymmetic^{\mu\nu}_{\ \ \gamma\delta}.
		\end{align}
		By the Bianchi identities, we have:
		\begin{subequations}
			\begin{align}
				\DWeylBianchi_{\sigma\gamma\delta\alpha\beta}:=&~\Dfour_{[\sigma}\Weyl_{\gamma\delta]\alpha\beta}=\antisymmetic^\mu_{\ \sigma\gamma\delta}\DWeylStar_{\mu\alpha\beta},&
				\DWeylBianchiStar_{\sigma\gamma\delta\alpha\beta}:=&~\Dfour_{[\sigma}\Weyldual_{\gamma\delta]\alpha\beta}=-\antisymmetic^\mu_{\ \sigma\gamma\delta}\DWeyl_{\mu\alpha\beta}.
			\end{align}
		\end{subequations}
	\end{definition}
	
	\begin{lemma}\cite[(7.3.7). Null Bianchi identities for Riemann curvature tensor]{ChristodoulouKlainerman1993}.\label{LE:NullBianchi0}
		\begin{subequations}
			\begin{align}
				\Dfour_4\Weyl_{A3B3}=&~-\frac{1}{2}(\Dfour_B\Weyl_{A334}+\Dfour_A\Weyl_{B334})+\frac{1}{2}\Dfour_3(\Weyl_{A3B4}+\Weyl_{B3A4})+\frac{1}{2}(\DWeylBianchi_{A3B34}+\DWeylBianchi_{B3A34}),\\
				\Dfour_3\Weyl_{A4B4}=&~\frac{1}{2}(\Dfour_B\Weyl_{A434}+\Dfour_A\Weyl_{B434})+\frac{1}{2}\Dfour_4(\Weyl_{A3B4}+\Weyl_{B3A4})+\frac{1}{2}(\DWeylBianchi_{A4B43}+\DWeylBianchi_{B4A43}),\\
				\Dfour_4\Weyl_{A434}=&~2\Dfour^B\Weyl_{A4B4}-2\DWeyl_{4A4},\\
				\Dfour_3\Weyl_{A434}=&~-2\Dfour^B\Weyl_{B3A4}+2\DWeyl_{3A4},\\
				\Dfour_4\Weyl_{A334}=&~2\Dfour^B\Weyl_{A3B4}-2\DWeyl_{4A3},\\
				\Dfour_3\Weyl_{A334}=&~-2\Dfour^B\Weyl_{A3B3}+2\DWeyl_{3A3},\\
				\Dfour_4\Weyl_{3434}=&~2\Dfour^B\Weyl_{B434}-2\DWeyl_{434},\\
				\Dfour_3\Weyl_{3434}=&~-2\Dfour^B\Weyl_{B334}+2\DWeyl_{334},\\
				\Dfour_4\HodgeLeftDual\Weyl_{3434}=&~2\Dfour^B\HodgeLeftDual\Weyl_{B434}-2\DWeylStar_{434},\\
				\Dfour_3\HodgeLeftDual\Weyl_{3434}=&~-2\Dfour^B\HodgeLeftDual\Weyl_{B334}+2\DWeylStar_{334}.
			\end{align}
		\end{subequations}
	\end{lemma}
	\begin{lemma}\cite[(7.3.8)-(7.3.9). Derivatives of Weyl curvature tensor components]{ChristodoulouKlainerman1993}\label{LE:DofWeyl}.
		\begin{subequations}
			\begin{align}
				\Dfour_3\Weyl_{A4B4}=&~\angD_3\AlphaAB_{AB}-4\uomega\AlphaAB_{AB}-4\left(\upeta\tensorproduct\BetaA\right)_{AB},\\
				\Dfour_4\Weyl_{A3B3}=&~\angD_4\Alphabar_{AB}-4\omega\Alphabar_{AB}+4\left(\ueta\tensorproduct\Betabar\right)_{AB},\\
				\Dfour_4\Weyl_{A434}=&~2\angD_{4}\BetaA_A+4\omega\BetaA_A-2\ueta_B\AlphaAB_{AB},\\
				\Dfour_3\Weyl_{A434}=&~2\angD_{3}\BetaA_A-4\uomega\BetaA_A-6\upeta_A\Rho-6\HodgeLeftDual\upeta_A\Rhodual,\\
				\Dfour_4\Weyl_{A334}=&~2\angD_{4}\Betabar_A-4\omega\Betabar_A+6\ueta_A\Rho-6\HodgeLeftDual\ueta_A\Rhodual,\\
				\Dfour_3\Weyl_{A334}=&~2\angD_{3}\Betabar_A-4\uomega\Betabar_A-2\upeta_B\Alphabar_{AB}.
			\end{align}
		\end{subequations}
		Moreover, we have:
		\begin{subequations}
			\begin{align}
				\frac{1}{2}\Dfour_3(\Weyl_{A3B4}+\Weyl_{B3A4})=&~\left(-\angD_3\Rho-2\upeta\cdot\Betabar\right)\gsphere_{AB},\\
				\frac{1}{2}\Dfour_4(\Weyl_{A3B4}+\Weyl_{B3A4})=&~\left(-\angD_4\Rho+2\ueta\cdot\BetaA\right)\gsphere_{AB},\\	
				\Dfour_4\Weyl_{3434}=&~4\angD_4\Rho-8\ueta\cdot\BetaA,\\
				\Dfour_3\Weyl_{3434}=&~4\angD_3\Rho+8\upeta\cdot\Betabar,\\
				\Dfour_4\HodgeLeftDual\Weyl_{3434}=&~4\angD_4\Rhodual+8\ueta\cdot\HodgeLeftDual\BetaA,\\
				\Dfour_3\HodgeLeftDual\Weyl_{3434}=&~4\angD_3\Rhodual+8\upeta\cdot\HodgeLeftDual\Betabar.
			\end{align}
		\end{subequations}
		On the other hand,
		\begin{subequations}\label{EQ:DerivativeoWeyl4}
			\begin{align}
				\Dfour^B\Weyl_{A4B4}=&~\angnabla^B\AlphaAB_{AB}-2\gtr\upchi\BetaA_A+2\upzeta_B\AlphaAB_{AB},\\
				\Dfour^B\Weyl_{A3B3}=&~\angnabla^B\Alphabar_{AB}+2\gtr\uchi\Betabar_A-2\upzeta_B\Alphabar_{AB}.
			\end{align}
		\end{subequations}
		Using the fact that $u_{AC}v_{CB}+u_{BC}v_{CA}=(u\cdot v)\gsphere_{AB}$ for traceless symmetric forms $u,v$ tangent to $\suubar$, we have:
		\begin{subequations}\label{EQ:DerivativeoWeyl5}
			\begin{align}
				\label{EQ:DerivativeoWeyl5.1}\frac{1}{2}\left(\Dfour_B\Weyl_{A434}+\Dfour_A\Weyl_{B434}\right)=&~\angnabla_A\BetaA_B+\angnabla_B\BetaA_A-\frac{1}{2}\gtr\uchi\AlphaAB_{AB}-3\upchi_{AB}\Rho-3\HodgeLeftDual\hchi_{AB}\Rhodual-\frac{1}{2}(\AlphaAB\cdot\huchi)\gsphere_{AB}+(\upzeta_A\BetaA_B+\upzeta_B\BetaA_A),\\
				\label{EQ:DerivativeoWeyl5.2}\frac{1}{2}\left(\Dfour_B\Weyl_{A334}+\Dfour_A\Weyl_{B334}\right)=&~\angnabla_A\Betabar_B+\angnabla_B\Betabar_A+\frac{1}{2}\gtr\upchi\Alphabar_{AB}+3\uchi_{AB}\Rho-3\HodgeLeftDual\huchi_{AB}\Rhodual+\frac{1}{2}(\Alphabar\cdot\hchi)\gsphere_{AB}-(\upzeta_A\Betabar_B+\upzeta_B\Betabar_A).
			\end{align}
		\end{subequations}
		\begin{subequations}\label{EQ:DerivativeoWeyl6}
			\begin{align}	
				\Dfour^B\Weyl_{A3B4}=&~-\angnabla_A\Rho+\antisymmetic_{AB}\angnabla_B\Rhodual-\gtr\upchi\Betabar_A+2\huchi_{AB}\BetaA_B,\\
				\Dfour^B\Weyl_{B3A4}=&~-\angnabla_A\Rho+\antisymmetic_{BA}\angnabla_B\Rhodual+\gtr\uchi\BetaA_A-2\hchi_{AB}\Betabar_B.
			\end{align}
		\end{subequations}
		Contracting \eqref{EQ:DerivativeoWeyl5.1} and \eqref{EQ:DerivativeoWeyl5.2}, we obtain:
		\begin{subequations}\label{EQ:DerivativeoWeyl7}
			\begin{align}
				\label{EQ:DerivativeoWeyl7.1}\Dfour^A\Weyl_{A434}=&~2\angdiv\BetaA-3\gtr\upchi\Rho-\AlphaAB\cdot\uchi+2\upzeta\cdot\BetaA\\
				\label{EQ:DerivativeoWeyl7.2}\Dfour^A\Weyl_{A334}=&~2\angdiv\Betabar+3\gtr\uchi\Rho+\Alphabar\cdot\upchi-2\upzeta\cdot\Betabar.
			\end{align}
		\end{subequations}
		Taking the Hodge dual of \eqref{EQ:DerivativeoWeyl7}, we have:
		\begin{subequations}\label{EQ:DerivativeoWeyl8}
			\begin{align}
				\label{EQ:DerivativeoWeyl8.1}\Dfour^A\HodgeLeftDual\Weyl_{A434}=&~-2\angdiv\HodgeLeftDual\BetaA-3\gtr\upchi\Rhodual+\HodgeLeftDual\AlphaAB\cdot\huchi-2\upzeta\cdot\HodgeLeftDual\BetaA\\
				\label{EQ:DerivativeoWeyl8.2}\Dfour^A\HodgeLeftDual\Weyl_{A334}=&~2\angdiv\HodgeLeftDual\Betabar+3\gtr\upchi\Rhodual-\HodgeLeftDual\Alphabar\cdot\hchi+2\upzeta\cdot\HodgeLeftDual\Betabar.
			\end{align}
		\end{subequations}
	\end{lemma}
	\begin{lemma}[Expression of the Cotton tensor components]\label{LE:CottonTensor}
		\begin{subequations}
			\begin{align}
				\DWeyl_{4A4}=&~-\frac{1}{2}\Dfour_4\Ricfour{4}{A}+\frac{1}{2}\Dfour_A\Ricfour{4}{4}\\
				\notag=&~-\angD_{4}\FBetaA_A+\angnabla_A\FAlphaAB-2\omega\FBetaA_A-2\upchi_{AB}\FBetaA_B+2\upzeta_A\FAlphaAB  {+\ueta_A\FAlphaAB},\\
				\DWeyl_{3A4}=&~-\frac{1}{2}\Dfour_4\Ricfour{3}{A}+\frac{1}{2}\Dfour_A\Ricfour{3}{4}+\frac{1}{6}\Dfour_A\Scalar\\
				\notag=&~-\angD_{4}\FBetabar_A+\angnabla_A\FRho+\frac{1}{6}\angnabla_A\Scalar+\frac{1}{2}(\Scalar+  {4}\FRho)\ueta_A+\FRhodual_{AB}\ueta_B+2\omega\FBetabar_A-\uchi_{AB}\FBetaA_B-\upchi_{AB}\FBetabar_B\\
				\DWeyl_{4A3}=&~-\frac{1}{2}\Dfour_3\Ricfour{4}{A}+\frac{1}{2}\Dfour_A\Ricfour{3}{4}+\frac{1}{6}\Dfour_A\Scalar\\
				\notag=&~-\angD_{3}\FBetaA_A+\angnabla_A\FRho+\frac{1}{6}\angnabla_A\Scalar+\frac{1}{2}(\Scalar+  {4}\FRho)\upeta_A+\FRhodual_{AB}\upeta_B+2\uomega\FBetaA_A-\uchi_{AB}\FBetaA_B-\upchi_{AB}\FBetabar_B\\
				\DWeyl_{3A3}=&~-\frac{1}{2}\Dfour_3\Ricfour{3}{A}+\frac{1}{2}\Dfour_A\Ricfour{3}{3}\\
				\notag=&~-\angD_{3}\FBetabar_A+\angnabla_A\FAlphabar-2\uomega\FBetabar_A-2\uchi_{AB}\FBetabar_B-2\upzeta_A\FAlphabar  {+\upeta_A\FAlphabar},\\
				\DWeyl_{434}=&~-\frac{1}{2}\Dfour_4\Ricfour{3}{4}+\frac{1}{2}\Dfour_3\Ricfour{4}{4}-\frac{1}{6}\Dfour_4\Scalar\\
				\notag=&~-\angD_{4}\FRho+\angD_{3}\FAlphaAB-\frac{1}{6}\angD_4\Scalar+2\ueta\cdot\FBetaA-4\upeta\cdot\FBetaA-4\uomega\FAlphaAB\\
				\DWeyl_{334}=&~-\frac{1}{2}\Dfour_4\Ricfour{3}{3}+\frac{1}{2}\Dfour_3\Ricfour{3}{4}+\frac{1}{6}\Dfour_3\Scalar\\
				\notag=&~-\angD_{4}\FAlphabar+\angD_{3}\FRho  {+}\frac{1}{6}\angD_3\Scalar+4\ueta\cdot\FBetabar+4\omega\FAlphabar-2\upeta\cdot\FBetabar,\\
				\DWeylStar_{434}=&~\antisymmetic^{AB}\Dfour_A\Ricfour{4}{B}\\
				\notag=&~2\angdiv\HodgeLeftDual\FBetaA-\hchi\wedge\FRhodual+2\upzeta\cdot\HodgeLeftDual\FBetaA,\\
				\DWeylStar_{334}=&~\antisymmetic^{AB}\Dfour_A\Ricfour{3}{B}\\
				\notag=&~2\angdiv\HodgeLeftDual\FBetabar-\huchi\wedge\FRhodual-2\upzeta\cdot\HodgeLeftDual\FBetabar,\\
				\DWeylBianchi_{A3B34}=&~-\DWeylBianchi_{B3A34},\\
				\DWeylBianchi_{A4B43}=&~-\DWeylBianchi_{B4A43}.
			\end{align}
		\end{subequations}
		Using the contracted Bianchi equations, we have more expressions for some of the components as follows:
		\begin{subequations}
			\begin{align}
				\DWeyl_{3A4}=&~\frac{1}{2}\Dfour_3\Ricfour{4}{A}-\Dfour_B\Ricfour{A}{B}+\frac{1}{2}\Dfour_A\Ricfour{3}{4}+\frac{2}{3}\Dfour_A\Scalar\\
				\notag=&~\angD_{ 3}\FBetaA_A-\angnabla_B\FRhodual_{AB}+\frac{1}{6}\angnabla_A\Scalar-\frac{1}{2}\upeta_A(\Scalar+ {4}\FRho)-\upeta_B\FRhodual_{AB}-2\uomega\FBetaA_A {+(\gtr\uchi)\FBetaA_A+(\gtr\upchi)\FBetabar_A},\\
				\DWeyl_{4A3}=&~\frac{1}{2}\Dfour_4\Ricfour{3}{A}-\Dfour_B\Ricfour{A}{B}+\frac{1}{2}\Dfour_A\Ricfour{3}{4}+\frac{2}{3}\Dfour_A\Scalar\\
				\notag=&~\angD_{ 4}\FBetabar_A-\angnabla_B\FRhodual_{AB}+\frac{1}{6}\angnabla_A\Scalar-\frac{1}{2}\ueta_A(\Scalar+ {4}\FRho)-\ueta_B\FRhodual_{AB}-2\omega\FBetabar_A {+(\gtr\uchi)\FBetaA_A+(\gtr\upchi)\FBetabar_A},\\
				\DWeyl_{434}=&~-\Dfour_4\Ricfour{3}{4}+\Dfour_A\Ricfour{A}{4}-\frac{2}{3}\Dfour_4\Scalar\\
				\notag=&~-2\angD_{4}\FRho+2\angdiv\FBetaA-\frac{2}{3}\angD_4\Scalar+4\ueta\cdot\FBetaA-\frac{1}{2}(\Scalar+ {4}\FRho)\gtr\upchi-\hchi\cdot\FRhodual+2\upzeta\cdot\FBetaA {-(\gtr\uchi)\FAlphaAB},\\
				\DWeyl_{334}=&~\Dfour_3\Ricfour{3}{4}-\Dfour_A\Ricfour{A}{3}+\frac{2}{3}\Dfour_3\Scalar\\
				\notag=&~2\angD_{3}\FRho-2\angdiv\FBetabar+\frac{2}{3}\angD_3\Scalar-4\upeta\cdot\FBetabar+\frac{1}{2}(\Scalar+ {4}\FRho)\gtr\uchi+\huchi\cdot\FRhodual+2\upzeta\cdot\FBetabar {+(\gtr\upchi)\FAlphabar},\\
				\DWeyl_{434}=&~\Dfour_3\Ricfour{4}{4}-\Dfour_A\Ricfour{A}{4}+\frac{1}{3}\Dfour_4\Scalar\\
				\notag=&~2\angD_{3}\FAlphaAB-2\angdiv\FBetaA+\frac{1}{3}\angD_4\Scalar-8\uomega\FAlphaAB-8\upeta {\cdot}\FBetaA+\frac{1}{2}(\Scalar+ {4}\FRho)\gtr\upchi+\hchi\cdot\FRhodual-2\upzeta\cdot\FBetaA {+(\gtr\uchi)\FAlphaAB},\\
				\DWeyl_{334}=&~-\Dfour_4\Ricfour{3}{3}+\Dfour_A\Ricfour{A}{3}-\frac{1}{3}\Dfour_3\Scalar\\
				\notag=&~-2\angD_{4}\FAlphabar+2\angdiv\FBetabar-\frac{1}{3}\angD_3\Scalar+8\omega\FAlphabar+8\ueta {\cdot}\FBetabar-\frac{1}{2}(\Scalar+ {4}\FRho)\gtr\uchi-\huchi\cdot\FRhodual-2\upzeta\cdot\FBetabar {-(\gtr\upchi)\FAlphabar}.
			\end{align}
		\end{subequations}
	\end{lemma}
	\begin{proposition}[Null Bianchi identities]\label{PR:NullBianchi0}
		Following the similar approach as in \cite[Chapter 7]{ChristodoulouKlainerman1993}, combining Lemma \ref{LE:NullBianchi0}, Lemma \ref{LE:DofWeyl} and Lemma \ref{LE:CottonTensor}, we have the following null Bianchi equations:
		\begin{subequations}\label{EQ:NullBianchi0}
			\begin{align}
				\label{EQ:D4Rho1}\angD_4\Rho+\frac{3}{2}\gtr\upchi\Rho=&~\angdiv\BetaA-\frac{1}{2}\huchi\cdot\AlphaAB+\upzeta\cdot\BetaA+2\ueta\cdot\BetaA\\
				\notag&+\angD_{4}\FRho-\angdiv\FBetaA+\frac{1}{3}\Dfour_4\Scalar-2\ueta\cdot\FBetaA+\frac{1}{4}(\Scalar+ {4}\FRho)\gtr\upchi+\frac{1}{2}\hchi\cdot\FRhodual-\upzeta\cdot\FBetaA {+\frac12(\gtr\uchi)\FAlphaAB},\\
				\label{EQ:D3Rho1}\angD_3\Rho+\frac{3}{2}\gtr\uchi\Rho=&~-\angdiv\Betabar-\frac{1}{2}\hchi\cdot\Alphabar+\upzeta\cdot\Betabar-2\upeta\cdot\Betabar\\
				\notag&+\angD_{3}\FRho-\angdiv\FBetabar+\frac{1}{3}\Dfour_3\Scalar-2\upeta\cdot\FBetabar+\frac{1}{4}(\Scalar+ {4}\FRho)\gtr\uchi+
				\frac{1}{2}\huchi\cdot\FRhodual+\upzeta\cdot\FBetabar {+\frac12(\gtr\upchi)\FAlphabar},\\
				\label{EQ:D4Rhodual1}\angD_4\Rhodual+\frac{3}{2}\gtr\upchi\Rhodual=&~-\angdiv\HodgeLeftDual\BetaA+\frac{1}{2}
				\huchi\wedge\AlphaAB-\upzeta\cdot\HodgeLeftDual\BetaA-2\ueta\cdot\HodgeLeftDual\BetaA\\
				\notag&-\angdiv\HodgeLeftDual\FBetaA+\frac{1}{2}\hchi\wedge\FRhodual-\upzeta\cdot\HodgeLeftDual\FBetaA,\\
				\label{EQ:D3Rhodual1}\angD_3\Rhodual+\frac{3}{2}\gtr\uchi\Rhodual=&~-\angdiv\HodgeLeftDual\Betabar-\frac{1}{2}
				\hchi\wedge\Alphabar {+}\upzeta\cdot\HodgeLeftDual\Betabar-2\upeta\cdot\HodgeLeftDual\Betabar\\
				\notag&+\angdiv\HodgeLeftDual\FBetabar-\frac{1}{2}\huchi\wedge\FRhodual-\upzeta\cdot\HodgeLeftDual\FBetabar.
			\end{align}
		\end{subequations}
		\begin{subequations}
			\begin{align}
				\angD_{4}\BetaA+2\gtr\upchi\BetaA=&~\angdiv\AlphaAB-2\omega\BetaA+\ueta\cdot\AlphaAB+2\upzeta\cdot\AlphaAB\\
				\notag&+\angD_{4}\FBetaA-\angnabla\FAlphaAB+2\omega\FBetaA+2\upchi\cdot\FBetaA-2\upzeta\FAlphaAB {-\ueta\FAlphaAB},\\
				\angD_{3}\BetaA+\gtr\uchi\BetaA=&~\angnabla\Rho+\HodgeLeftDual\angnabla\Rhodual+2\uomega\BetaA+2\hchi\cdot\Betabar+3\upeta\Rho+3\HodgeLeftDual\upeta\Rhodual\\
				\notag&+\angD_{ 3}\FBetaA-\angdiv\FRhodual+\frac{1}{6}\angnabla\Scalar-\frac{1}{2}\upeta(\Scalar+ {4}\FRho)-\upeta\cdot\FRhodual-2\uomega\FBetaA {+(\gtr\uchi)\FBetaA+(\gtr\upchi)\FBetabar},\\
				\label{EQ:D4Betabar1}\angD_{4}\Betabar+\gtr\upchi\Betabar=&~-\angnabla\Rho+\HodgeLeftDual\angnabla\Rhodual+2\omega\Betabar-3\ueta\Rho+3\HodgeLeftDual\ueta\Rhodual+2\huchi\cdot\BetaA\\
				\notag&-\angD_{ 4}\FBetabar+\angdiv\FRhodual-\frac{1}{6}\angnabla\Scalar+\frac{1}{2}\ueta(\Scalar+ {4}\FRho)+\ueta\cdot\FRhodual+2\omega\FBetabar {-(\gtr\uchi)\FBetaA-(\gtr\upchi)\FBetabar},\\
				\angD_{3}\Betabar+2\gtr\uchi\Betabar=&~-\angdiv\Alphabar+2\upzeta\cdot\Alphabar-2\uomega\Betabar-\upeta\cdot\Alphabar\\
				\notag&-\angD_{3}\FBetabar+\angnabla\FAlphabar-2\uomega\FBetabar-2\uchi\cdot\FBetabar-2\upzeta\FAlphabar {+\upeta\FAlphabar},
			\end{align}
		\end{subequations}
		\begin{subequations}\label{EQ:DAlphaDAlphabar}
			\begin{align}
				\angD_4\Alphabar+\frac{1}{2}\gtr\upchi\Alphabar=&~-\angnabla\tensorproduct\Betabar+4\omega\Alphabar+(\upzeta-4\ueta)\tensorproduct\Betabar-3\huchi\Rho+3\HodgeLeftDual\huchi\Rhodual\\
				\notag&+\left(\angD_{4}\FAlphabar-\angdiv\FBetabar+\frac{1}{6}\angD_3\Scalar-4\omega\FAlphabar-4\ueta\FBetabar+\frac{1}{4}(\Scalar+2\FRho)\gtr\uchi+\frac{1}{2}\huchi\cdot\FRhodual+\upzeta\cdot\FBetabar\right)\gsphere,\\
				\angD_3\AlphaAB+\frac{1}{2}\gtr\uchi\AlphaAB=&~\angnabla\tensorproduct\BetaA+4\uomega\AlphaAB+\left(\upzeta+4\upeta\right)\tensorproduct\BetaA-3\hchi\Rho-3\HodgeLeftDual\hchi\Rhodual\\
				\notag&+\left(\angD_{3}\FAlphaAB-\angdiv\FBetaA+\frac{1}{6}\angD_4\Scalar-4\uomega\FAlphaAB-4\upeta\FBetaA+\frac{1}{4}(\Scalar+2\FRho)\gtr\upchi+\frac{1}{2}\hchi\cdot\FRhodual-\upzeta\cdot\FBetaA\right)\gsphere.
			\end{align}
		\end{subequations}
	\end{proposition}
	\begin{definition}[Renormalization of Weyl tensor components]\label{DE:Renormalization}
		We renormalize the Weyl tensor components as follows:
		\begin{subequations}\label{EQ:Renormalization}
			\begin{align}
				\CAlphaAB:=&~\AlphaAB-\FAlphaAB\gsphere,&
				\CAlphabar:=&~\Alphabar-\FAlphabar\gsphere,\\
				\CBetaA:=&~\BetaA-\FBetaA,&
				\CBetabar:=&~\Betabar+\FBetabar,\\
				\CRho:=&~\Rho-\FRho-\frac{1}{3}\Scalar-\frac{1}{2}\hchi\cdot\huchi,&
				\CRhodual:=&~\Rhodual+\frac{1}{2}\huchi\wedge\hchi.
			\end{align}
		\end{subequations}
		Then the Gaussian curvature of $\gsphere$ satisfies:
		\begin{align}\label{EQ:Renormalization2}
			\GaussCurvature:=&~-\CRho-\frac{1}{4}\gtr\upchi\gtr\uchi.
		\end{align}
	\end{definition}
	\begin{proposition}[Renormalized null Bianchi equations]\label{PR:ReNullBianchi}
		Under the renormalization in Definition \ref{DE:Renormalization}, the null Bianchi equations \eqref{EQ:NullBianchi0}-\eqref{EQ:DAlphaDAlphabar} read:
		\begin{subequations}
			\begin{align}
				\label{EQ:D4Rho2}\angD_4\CRho+\frac{3}{2}\gtr\upchi\CRho=&~\angdiv\CBetaA-\frac{1}{4}\left( {2}\FRho+\Scalar\right)\gtr\upchi {+\frac12(\gtr\uchi)\FAlphaAB}\\
				\notag&+\frac{1}{2}\left(\frac{1}{2}\gtr\uchi\hchi-\angnabla\tensorproduct\ueta-\ueta\tensorproduct\ueta\right)\cdot\hchi+\upzeta\cdot\CBetaA+2\ueta\cdot\CBetaA\\
				\label{EQ:D3Rho2}\angD_{3}\CRho+\frac{3}{2}\gtr\uchi\CRho=&~-\angdiv\CBetabar-\frac{1}{4}\left( {2}\FRho+\Scalar\right)\gtr\uchi {+\frac12(\gtr\upchi)\FAlphabar}\\
				\notag&+\frac{1}{2}\left(\frac{1}{2}\gtr\upchi\huchi-\angnabla\tensorproduct\upeta-\upeta\tensorproduct\upeta\right)\cdot\huchi+\upzeta\cdot\CBetabar-2\upeta\cdot\CBetabar,\\
				\label{EQ:D4Rhodual2}\angD_4\CRhodual+\frac{3}{2}\gtr\upchi\CRhodual=&~-\angdiv\HodgeLeftDual\CBetaA-2\angdiv\HodgeLeftDual\FBetaA-\frac{1}{2}\hchi\wedge(\angnabla\tensorproduct\ueta+\ueta\tensorproduct\ueta)-\upzeta\cdot(\HodgeLeftDual\BetaA+\HodgeLeftDual\FBetaA)-2\ueta\cdot\HodgeLeftDual\BetaA,\\
				\label{EQ:D3Rhodual2}\angD_3\CRhodual+\frac{3}{2}\gtr\uchi\CRhodual=&~-\angdiv\HodgeLeftDual\CBetabar+2\angdiv\HodgeLeftDual\FBetabar+\frac{1}{2}
				\huchi\wedge(\angnabla\tensorproduct\upeta+\upeta\tensorproduct\upeta) {+}\upzeta\cdot(\HodgeLeftDual\Betabar-\HodgeLeftDual\FBetabar)-2\upeta\cdot\HodgeLeftDual\Betabar,
			\end{align}
		\end{subequations}
		\begin{subequations}
			\begin{align}
				\label{EQ:D4CBeta}\angD_{4}\CBetaA+2\gtr\upchi\CBetaA=&~\angdiv\CAlphaAB-2\gtr\upchi\FBetaA-2\omega\CBetaA+\ueta\cdot {\CAlphaAB}+2\upzeta\cdot\CAlphaAB+2\upchi\cdot\FBetaA,\\
				\label{EQ:D3CBetabar}\angD_{3}\CBetabar+2\gtr\uchi\CBetabar=&~-\angdiv\CAlphabar+2\gtr\uchi\FBetabar+2\upzeta\cdot\CAlphabar-2\uomega\CBetabar-\upeta\cdot {\CAlphabar}-2\uchi\cdot\FBetabar,\\
				\label{EQ:D3CBeta}\angD_{3}\CBetaA+\gtr\uchi\CBetaA=&~\angnabla\CRho+\HodgeLeftDual\angnabla\CRhodual+\frac{1}{2}\angnabla(\hchi\cdot\huchi) {+\frac{1}{2}\angnabla(\hchi\wedge\huchi)}+2\uomega\CBetaA+2\hchi\cdot\Betabar+3\upeta\Rho+3\HodgeLeftDual\upeta\Rhodual\\
				\notag&+\angnabla\FRho-\angdiv\FRhodual+\frac{1}{2}\angnabla\Scalar {-\frac{1}{2}\upeta(\Scalar+4\FRho)-\upeta\cdot\FRhodual+(\gtr\upchi)\FBetabar},\\
				\label{EQ:D4CBetabar}\angD_{4}\CBetabar+\gtr\upchi\CBetabar=&~-\angnabla\CRho+\angnabla\CRhodual-\frac{1}{2}\angnabla(\hchi\cdot\huchi) {+\frac{1}{2}\HodgeLeftDual\angnabla(\hchi\wedge\huchi)}+2\omega\CBetabar-3\ueta\Rho+3\HodgeLeftDual\ueta\Rhodual+2\huchi\cdot\BetaA\\
				\notag&-\angnabla\FRho+\angdiv\FRhodual-\frac{1}{2}\angnabla\Scalar {+\frac{1}{2}\ueta(\Scalar+4\FRho)+\ueta\cdot\FRhodual-(\gtr\uchi)\FBetaA},
			\end{align}
		\end{subequations}
		\begin{subequations}
			\begin{align}
				\label{EQ:D4CAlphabar}\angD_4\CAlphabar+\frac{1}{2}\gtr\upchi\CAlphabar=&~-\angnabla\tensorproduct\CBetabar+4\omega\CAlphabar+(\upzeta-4\ueta)\tensorproduct\CBetabar-3\huchi\Rho+3\HodgeLeftDual\huchi\Rhodual\\
				\notag&+\angnabla\tensorproduct\FBetabar+\left(-\angdiv\FBetabar+\frac{1}{6}\angD_3\Scalar-\frac{1}{2}\gtr\upchi\FAlphabar+\frac{1}{4}(\Scalar+2\FRho)\gtr\uchi+\frac{1}{2}\huchi\cdot\FRhodual\right)\gsphere,\\
				\label{EQ:D3CAlpha}\angD_3\CAlphaAB+\frac{1}{2}\gtr\uchi\CAlphaAB=&~\angnabla\tensorproduct\CBetaA+4\uomega\CAlphaAB+\left(\upzeta+4\upeta\right)\tensorproduct\CBetaA-3\hchi\Rho-3\HodgeLeftDual\hchi\Rhodual\\
				\notag&+\angnabla\tensorproduct\FBetaA+\left(-\angdiv\FBetaA+\frac{1}{6}\angD_4\Scalar-\frac{1}{2}\gtr\uchi\FAlphaAB+\frac{1}{4}(\Scalar+2\FRho)\gtr\upchi+\frac{1}{2}\hchi\cdot\FRhodual\right)\gsphere.
			\end{align}
		\end{subequations}
		In addition, we have the following equations for $\GaussCurvature=-\CRho-\frac{1}{4}\gtr\upchi\gtr\uchi$:
		\begin{align}
			\label{EQ:D4Rho3}\angD_{4}\GaussCurvature+\gtr\upchi\GaussCurvature=&~\angdiv\CBetaA-\upzeta\cdot\CBetaA-2\ueta\cdot\CBetaA+\frac{1}{2}\hchi\cdot(\angnabla\tensorproduct\ueta)+\frac{1}{2}\hchi\cdot(\ueta\tensorproduct\ueta)-\frac{1}{2}\gtr\upchi\angdiv\ueta-\frac{1}{2}\gtr\upchi\sgabs{\ueta}^2\\
			\notag&+\frac{1}{2}\FRho\gtr\upchi+\frac{1}{2}\gtr\uchi\FAlphaAB,\\
			\label{EQ:D3Rho3}\angD_{3}\GaussCurvature+\gtr\uchi\GaussCurvature=&~\angdiv\CBetabar-\upzeta\cdot\CBetabar+2\upeta\cdot\CBetabar+\frac{1}{2}\huchi\cdot(\angnabla\tensorproduct\upeta)+\frac{1}{2}\huchi\cdot(\upeta\tensorproduct\upeta)-\frac{1}{2}\gtr\uchi\angdiv\upeta-\frac{1}{2}\gtr\uchi\sgabs{\upeta}^2\\
			\notag&+\frac{1}{2}\FRho\gtr\uchi+\frac{1}{2}\gtr\upchi\FAlphabar.
		\end{align}
	\end{proposition}

	\subsection{First-order Euler Equations}\label{SS:EulerEquations}
	\begin{definition}[Acoustical metric and its inverse]\label{DE:acousticunrescaled}
		We define the acoustical metric $(\gacoustical)_{\alpha\beta}$ and its inverse $(\gacoustical^{-1})^{\alpha\beta}$ as follows:
		\begin{subequations}
			\begin{align}
				(\gacoustical)_{\alpha\beta}&:=\speed^{-2}\gfour_{\alpha\beta}+(\speed^{-2}-1)v_\alpha v_\beta,\\
				(\gacoustical^{-1})^{\alpha\beta}
				&:=\speed^2(\gfour^{-1})^{\alpha\beta}+(\speed^2-1)v^\alpha v^\beta.
			\end{align}
		\end{subequations}	
	\end{definition}
	
	Now, we provide the first-order Euler equations in the framework of Einstein-Euler system Section \ref{SS:EinsteinEuler}.
	\begin{proposition}\cite[First-order Euler equations]{DisconziSpeck2019,Christodoulou2007}.
		Contracting equation (\ref{EQ:DivEMT}) with $v_\nu$ and $\delta^\alpha_\nu+v^\alpha v_\nu$, respectively, we have:
		\begin{subequations}\label{EQ:FOE1}
			\begin{align}
				v^\mu\Dfour_\mu\pressure+\speed^2(\pressure+\density)\Dfour_\mu v^\mu=&~0,\\
				(\density+\pressure)v^\mu\Dfour_\mu v^\alpha+\vproject^{\alpha\mu}\Dfour_\mu\pressure=&~0,\\
				v^\mu\Dfour_\mu\ent=&~0.
			\end{align}
		\end{subequations}
		where $\vproject^{\mu\nu}:=(\gfour^{-1})^{\mu\nu}+v^\mu v^\nu$, or, equivalently:
		\begin{subequations}\label{EQ:FOE2}
			\begin{align}
				\label{EQ:FOE2.1}v^\mu\Dfour_\mu\lnenth+\speed^2\Dfour_\mu v^\mu=&~0,\\
				\label{EQ:FOE2.2}v^\mu\Dfour_\mu v^\alpha+\vproject^{\alpha\mu}\Dfour_\mu\lnenth-\frac{\temp}{\enth}\Dfour^\alpha\ent=&~0,\\
				\label{EQ:FOE2.3}v^\mu\Dfour_\mu\ent=&~0.
			\end{align}
		\end{subequations}
		
	\end{proposition}
	
	\section{Initial Data, Main Theorem and Bootstrap Assumptions}
	\subsection{Schematic Notation}\label{SS:SchematicNotation}
	Throughout this article, we use the following notation:
	\begin{subequations}\label{DE:notation}
		\begin{align}
			\psi\in&\{\hchi,\gtr\upchi,\omega\},\\
			\upsi\in&\{\upeta,\ueta,\huchi,\gtr\uchi+\frac{2}{\abs{u}},\uomega\},&
			\unpsi\in&\{\upeta,\huchi,\gtr\uchi+\frac{2}{\abs{u}},\uomega\},\\
			\upRiem\in&\{\BetaA,\Rho,\Rhodual,\Betabar,\Alphabar\},&
			\uRiem\in&\{\CRho,\CRhodual,\Betabar,\Alphabar\},\\
			\variables\in&\{v_4,v_3,v_A,\lnenth,\ent\},&
			\uvariable\in&\{\angv,\lnenth,\ent\},
		\end{align}
	\end{subequations}
	where $\angv\in\{v_1,v_2\}$ denotes the sphere-projected velocity field.
	\subsection{Choice of Parameters}\label{SS:ChoiceofParameters}
	In this section, we list our choice of parameters. 
	\begin{itemize}
		\item $a$: The size of the gravitational pulse.
		\item $b$: The size of the fluid focusing speed.
		\item $\delta$: The interval (in $\ubar$) that the gravitational pulse is injected.
		\item $\BAfluid$: The size of the perturbative fluid variables $\angv,\lnenth,\ent$.
		\item $N\geq 6$: Regularity needed to close the estimates.
	\end{itemize}

	In this article, we let\footnote{$a\simeq b$ means $C^{-1}b\leq a\leq C b$ for some universal constant $C$.} $a\simeq b$ be large numbers, and $\delta\simeq a^{-1}$ be a small number.
	Let $0<\iniBAfluid<a^{-\frac{1}{5}}$. Let $\iniBAsize>1$ be a large number such that $1<\iniBAsize<a^{\frac{1}{200}}$. Then we have $0<\iniBAfluid<a^{-\frac{1}{5}}<\iniBAsize^{-40}<1$. 
	
	Let $\iniBAfluid<\BAfluid<a^{-\frac{1}{10}}$. Let $\iniBAsize,\BAsize>1$ be a large number such that $\iniBAsize<\BAsize<a^{\frac{1}{100}}$. Then we have $\iniBAfluid<\BAfluid<a^{-\frac{1}{10}}<\BAsize^{-10}<\iniBAsize^{-10}<1$.
	
	\subsection{Norms}\label{SS:Norms}
	\begin{definition}[Weighted-mixed operators]\label{DE:operators}
		Let $\angD$ be as in Section \ref{SS:despacetime}. Let $F$ be a tensorfield and $\vec{m_i}$'s be multi-indices consisting of $3-$tuples. For $\vec{m}=(i,j,k)$ with $\abs{\vec{m}}\leq N$, we define
		\begin{subequations}\label{EQ:Operators}
			\begin{align}
				\Dfour^{\vec{m}}F:=&~\Dfour_4^i\Dfour_3^j\Dfour_{A_i\cdots A_k}^kF,\\
				\angD^{\vec{m}}F:=&~\angD_4^i\angD_3^j\angD_{A_i\cdots A_k}^kF,\\
				\WDfour^{\vec{m}}F:=&~(\delta\Dfour_4)^i(\abs{u}^k\Dfour_3)^j(\abs{u}\Dfour_{A_i\cdots A_k}^k)F,\\
				\WangD^{\vec{m}}F:=&~(\delta\angD_4)^i(\abs{u}^k\angD_3)^j(\abs{u}\angnabla_{A_i\cdots A_k}^k)F.
			\end{align}
		\end{subequations}
	\end{definition}
	Let $N\geq 6$. See Section \ref{SS:ChoiceofParameters} for choice of parameters $a,b,\delta$. For $0\leq i+j+k\leq N$, let $\vec{m}=(i,j,k)$ be defined as in Definition \ref{DE:operators}. Letting $\mathcal{D}:=[\uinf,\uts]\times[0,5\delta]$ be a bounded region for $(u,\ubar)$, and $\sphereproject$ as in Section \ref{SS:despacetime}, we define the following norms that will be used throughout the article:\\
	
	\noindent\textbf{Curvature components:}\\
	
	For $0\leq i+j+k\leq \frac{N-1}{2}$:
	\begin{align}
		\mathcal{R}_{\mathcal{D},i,j,k,\infty}:=&~\sup\limits_{(u,\ubar)\in\mathcal{D}}\left(\frac{\abs{u}}{\delta^{-1}a^{\frac{1}{4}}}\norm{\WangD^{\vec{m}}\AlphaAB}_{L^\infty(\suubar)}+\frac{\abs{u}^2}{a^{\frac{1}{2}}}\norm{\WangD^{\vec{m}}\BetaA}_{L^\infty(\suubar)}\right.\\
		\notag&\left.+\frac{\abs{u}^3}{\delta a^{\frac{1}{2}}}\norm{\WangD^{\vec{m}}(\CRho,\CRhodual,\Betabar)}_{L^2(\suubar)}+\frac{\abs{u}^3}{\delta^2 a}\norm{\WangD^{\vec{m}}\Alphabar}_{L^2(\suubar)}\right).
	\end{align}
	
	For $i+j+k= N-1$ with $k\geq 1$:
	\begin{align}	
		\tilde{\mathcal{R}}_{\mathcal{D},N-1,2}:=&~\sup\limits_{(u,\ubar)\in\mathcal{D}}\left(\frac{\abs{u}^{\frac{1}{2}}}{ \delta^{-\frac{1}{2}}a^{\frac{1}{4}}}\norm{\WangD^{\vec{m}}\AlphaAB}_{L^2(\hypersur_{u})}+\frac{1}{ \delta^{-\frac{1}{2}}a^{\frac{1}{4}}}\norm{\abs{u}^{\frac{1}{2}}\WangD^{\vec{m}}\BetaA}_{L^2(\hypersurbar_{\ubar})}\right.\\
		\notag&\left.+\frac{\abs{u}}{ \delta^{\frac{1}{2}} a^{\frac{1}{2}}}\norm{\WangD^{\vec{m}}\BetaA}_{L^2(\hypersur_{u})}+\frac{1}{\delta^{\frac{1}{2}} a^{\frac{1}{2}}}\norm{\abs{u}\WangD^{\vec{m}}(\CRho,\CRhodual)}_{L^2(\hypersurbar_{\ubar})}\right.\\
		\notag&\left.+\frac{\abs{u}^3}{\delta^{\frac{3}{2}} a}\norm{\WangD^{\vec{m}}(\CRho,\CRhodual,\Betabar)}_{L^2(\hypersur_{u})}+\frac{1}{\delta^{\frac{3}{2}} a}\norm{\abs{u}^3\WangD^{\vec{m}}\Betabar}_{L^2(\hypersurbar_{\ubar})}+\frac{1}{\delta^{2} a}\norm{\abs{u}^3\WangD^{\vec{m}}\Alphabar}_{L^2(\hypersurbar_{\ubar})}\right).
	\end{align}
	
	For $0\leq i+j+k\leq N-2$:
	\begin{align}	
		\mathcal{R}_{\mathcal{D},i,j,k,2}:=&~\sup\limits_{(u,\ubar)\in\mathcal{D}}\left(\frac{1}{\delta^{-1}a^{\frac{1}{4}}}\norm{\WangD^{\vec{m}}\AlphaAB}_{L^2(\suubar)}+\frac{\abs{u}}{a^{\frac{1}{2}}}\norm{\WangD^{\vec{m}}\BetaA}_{L^2(\suubar)}\right.\\
		\notag&\left.+\frac{\abs{u}^2}{\delta a^{\frac{1}{2}}}\norm{\WangD^{\vec{m}}(\CRho,\CRhodual,\Betabar)}_{L^2(\suubar)}+\frac{\abs{u}^2}{\delta^2 a}\norm{\WangD^{\vec{m}}\Alphabar}_{L^2(\suubar)}\right).
	\end{align}
	We denote $\mathcal{R}:=\sup_{i+j+k\leq N-1}\mathcal{R}_{i,j,k,2}+\sup_{l+m+n\leq (N-1)/2}\mathcal{R}_{l,m,n,\infty}+\tilde{\mathcal{R}}_{N-1,2}$.
	
	\noindent\textbf{Ricci coefficients:}\\
	
	For $i+j+k=N$ and $k\geq 2$:
	\begin{align}
		\tilde{\mathcal{O}}_{\mathcal{D},N,2}:=&~\sup\limits_{(u,\ubar)\in\mathcal{D}}\left(\frac{1}{\delta^{\frac{1}{2}}a^{\frac{1}{2}}}\norm{\WangD^{\vec{m}}(\gtr\upchi,\hchi)}_{L^2(\hypersur_{u})}+\frac{1}{\delta^{\frac{1}{2}}a^{\frac{1}{2}}}\norm{\WangD^{\vec{m}}\omega}_{L^2(\hypersur_{u})}+\frac{1}{\delta^{\frac{1}{2}}a^{\frac{1}{2}}}\norm{\WangD^{\vec{m}}\upeta}_{L^2(\hypersurbar_{\ubar})}\right.\\
		\notag&\left.+\frac{1}{\delta^{\frac{3}{2}}a}\norm{\abs{u}\WangD^{\vec{m}}(\upeta,\ueta)}_{L^2(\hypersur_{u})}+\frac{\abs{u}^{\frac{1}{2}}}{\delta a^{\frac{1}{2}}}\norm{\WangD^{\vec{m}}(\gtr\uchi,\huchi)}_{L^2(\hypersurbar_{\ubar})}\right),
	\end{align}
	
	For $0\leq i+j+k\leq \frac{N-1}{2}$:
	\begin{align}
		\mathcal{O}_{\mathcal{D},i,j,k,\infty}:=&~\sup\limits_{(u,\ubar)\in\mathcal{D}}\left(\frac{\abs{u}}{a^{1/2}}\norm{\WangD^{\vec{m}}\psi}_{L^\infty(\suubar)}+\frac{\abs{u}^2}{\delta a^{1/2}}\norm{\WangD^{\vec{m}}\upsi}_{L^\infty(\suubar)}\right).
	\end{align}
	
	For $i+j+k\leq N-1$ with $k\geq 2$:
	\begin{align}
		\mathcal{O}_{\mathcal{D},i,j,k,2}:=&~\sup\limits_{(u,\ubar)\in\mathcal{D}}\left(\frac{1}{a^{1/2}}\norm{\WangD^{\vec{m}}\psi}_{L^2(\suubar)}+\frac{\abs{u}}{\delta a^{1/2}}\norm{\WangD^{\vec{m}}\upsi}_{L^2(\suubar)}\right).
	\end{align}
	\noindent\textbf{Fluid variables:}\\
	
	For $0\leq i+j+k\leq \frac{N-1}{2}$:
	\begin{align}
		\mathcal{E}_{\mathcal{D},i,j,k,\infty}:=&~\sup\limits_{(u,\ubar)\in\mathcal{D}}\left(\frac{1}{b^{\frac{1}{2}}}\norm{\WDfour^{\vec{m}}(v_4)}_{L^\infty(\suubar)}+\frac{1}{b^{-\frac{1}{2}}}\norm{\WDfour^{\vec{m}}(v_3)}_{L^\infty(\suubar)}\right.\\
		\notag&\left.+\abs{u}\norm{\WDfour^{\vec{m}}(v_A)}_{L^\infty(\suubar)}+\abs{u}\norm{\WDfour^{\vec{m}}(\lnenth,\ent)}_{L^\infty(\suubar)}\right).	
	\end{align}
	
	For $i+j+k\leq N-1$ with $k\geq 2$:
	\begin{align}
		\mathcal{E}_{\mathcal{D},i,j,k,2}:=&~\sup\limits_{(u,\ubar)\in\mathcal{D}}\left(\frac{1}{b^{\frac{1}{2}}\abs{u}}\norm{\WDfour^{\vec{m}}(v_4)}_{L^2(\suubar)}+\frac{1}{b^{-\frac{1}{2}}\abs{u}}\norm{\WDfour^{\vec{m}}(v_3)}_{L^2(\suubar)}\right.\\
		\notag&\left.+\norm{\WDfour^{\vec{m}}(v_A)}_{L^2(\suubar)}+\norm{\WDfour^{\vec{m}}(\lnenth,\ent)}_{L^2(\suubar)}\right).
	\end{align}
	
	For $i+j+k=N$ with $k\geq 2$:
	\begin{align}
		\tilde{\mathcal{E}}_{\mathcal{D},N,2}:=&~\sup\limits_{(u,\ubar)\in\mathcal{D}}\left(\norm{\WDfour^{\vec{m}}(v_4)}_{L^2(\hypersur_{u})}+\frac{1}{b^{-\frac{1}{2}}}\norm{\WDfour^{\vec{m}}(v_3)}_{L^2(\hypersurbar_{\ubar})}\right.\\
		\notag&\left.+\frac{1}{b^{-\frac{1}{2}}}\norm{\WDfour^{\vec{m}}(v_A)}_{L^2(\hypersur_{u})}+\frac{1}{b^{-\frac{1}{2}}}\norm{\WDfour^{\vec{m}}(\lnenth,\ent)}_{L^2(\hypersur_{u})}+\norm{\WDfour^{\vec{m}}(v_A)}_{L^2(\hypersurbar_{\ubar})}+\norm{\WDfour^{\vec{m}}(\lnenth,\ent)}_{L^2(\hypersurbar_{\ubar})}\right).
	\end{align}
	
	\subsection{Short Pulse Ansatz}\label{SS:ShortPulseAnsatz}
	Inspired by \cite{Christodoulou2009}, we present our short-pulse ansatz for the Einstein-Euler system that dynamically forms trapped surface under either gravitational focusing or fluid focusing. That is, in this article, we construct initial data and propagate all quantities to spacetime region $\region:=[\uinf,\uts]\times[0,(1+2(\uts-\uinf))\delta]\times\Stwo$ so that quantities are consistent with our ansatz. We note that the size of quantities should be viewed as an upper bound, while its lower bound depends on the initial data assumptions under different focusing mechanisms; see Section \ref{SS:InitialData}-\ref{SS:MainTheorem}. 
	
	Recall Section \ref{SS:ChoiceofParameters} for parameters $a,b,\BAfluid,\delta$, Definition \ref{D:Defconnectioncoeff} for connection coefficients, Definition \ref{DE:CurvatureComponents} for curvature tensor components. We now present the ansatz in below:
	\subsubsection{The short pulse ansatz for the fast-injected fluid variables}
	 Recall \eqref{EQ:devi} for the definition of $v_I$. 
	We present the following ansatz for the quantities related to fluid variables:
	\begin{subequations}
		\begin{align}
			v_4\sim&~ -b^{1/2},& v_3\sim&~ -b^{-1/2},& v_A,\lnenth,\ent\sim&~ \BAfluid,
		\end{align}
	\end{subequations}
	\subsubsection{Short pulse ansatz on the gravitational quantities}
	We present the following ansatz for the quantities related to metric $\gfour$:
	\begin{subequations}
		\begin{align}
			\upchi,\omega\sim& \frac{a^{\frac{1}{2}}}{\abs{u}},&
			\ueta,\upeta\sim& \frac{\delta a^{\frac{1}{2}}}{\abs{u}^2},&
			\huchi\sim& \frac{\delta a^{\frac{1}{2}}}{\abs{u}^2},&
			\uomega\sim& \frac{\delta a^{\frac{1}{2}}}{\abs{u}^2},&
			\gtr\uchi+\frac{2}{\abs{u}}\sim& \frac{\delta a^{\frac{1}{2}}}{\abs{u}^2},\\
			\AlphaAB\sim& \frac{\delta^{-1}a^{\frac{1}{4}}}{\abs{u}},&
			\BetaA\sim& \frac{a^{\frac{1}{2}}}{\abs{u}^2},&
			\CRho,\CRhodual\sim& \frac{\delta a^{\frac{1}{2}}}{\abs{u}^3},&
			\Betabar\sim& \frac{\delta a^{\frac{1}{2}}}{\abs{u}^3},&
			\Alphabar\sim& \frac{\delta a^{\frac{1}{2}}}{\abs{u}^3}.
		\end{align}
	\end{subequations}
	\subsection{Characteristic Data}\label{SS:CharData}
	Inspired by \cite{Luk2012,Christodoulou2009}, in this section, we discuss the prescription of the characteristic initial data on $\hypersur_{\uinf}\cup\hypersurbar_{0}$. We impose the ansatz that the induced metric $\gsphere$ on $\hypersur_{\uinf}\cup\hypersurbar_{0}$ is of the following form:
	\begin{align}\label{EQ:CharData1}
		\gsphere=\Phi^2\hat{\gsphere}.
	\end{align}
	First of all, recall $\lapse=1$ on $\hypersur_{\uinf}\cup\hypersurbar_{0}$ in Section \ref{SS:NullFrame}. Let $b$ be a large number.
	We make the following assumptions\footnote{We note that these are not the only assumptions that one could make, especially for the gravity terms. We make these assumptions to conveniently estimate the favorable size of the geometric quantities on the characteristic initial hypersurfaces.} on the initial sphere $S_{\uinf,0}=\hypersur_{\uinf}\cup\hypersurbar_{0}$: 
	\begin{subequations}\label{EQ:CharData1.5}
		\begin{align}\label{EQ:CharData1.5.1}
			\lnenth,\ent,v_A&=0,&
			v_4&=b^{\frac{1}{2}},&
			v_3&=b^{-\frac{1}{2}},\\
			\Phi|_{S_{\uinf,0}}&=1,& \left.\frac{\p\Phi}{\p\ubar}\right|_{S_{\uinf,0}}&=\frac{1}{\abs{\uinf}},&
			\left.\frac{\p\Phi}{\p u}\right|_{S_{\uinf,0}}=-\frac{1}{\abs{\uinf}},\\
			\label{EQ:CharData1.5.2}\hchi|_{S_{\uinf,0}} \sim a^{\frac{1}{4}}, \upzeta,\huchi,\omega,\uomega|_{S_{\uinf,0}}&=0, &
			\AlphaAB,\Alphabar,\BetaA,\Betabar|_{S_{\uinf,0}}&=0, & \text{in shear focusing case},\\
			\label{EQ:CharData1.5.3}\upzeta,\hchi,\huchi,\omega,\uomega|_{S_{\uinf,0}}&=0&
			\AlphaAB,\Alphabar,\BetaA,\Betabar|_{S_{\uinf,0}}&=0, & \text{in fluid focusing case}.
		\end{align}
	\end{subequations} 
	In above, \eqref{EQ:CharData1.5.2} and \eqref{EQ:CharData1.5.3} correspond to assumption 1 and assumption 2, respectively, in Theorem \ref{TH:Main}.
	In coordinates \eqref{DE:Coordinatese3e4}, on $\hypersur_{\uinf}$ (where $\lapse=1,\Lambda^A=0$), \eqref{EQ:Lieg} reads:
	\begin{align}\label{EQ:CharData2}
		\frac{\p}{\p\ubar}\gsphere&=2\upchi=2\hchi+\gtr\upchi\cdot\gsphere.
	\end{align}
	Meanwhile, by \eqref{EQ:CharData1},
	\begin{align}\label{EQ:CharData3}
		\frac{\p}{\p\ubar}\gsphere&=\Phi^2\frac{\p}{\p\ubar}\hat{\gsphere}+2\Phi\frac{\p\Phi}{\p\ubar}\cdot\hat{\gsphere}.
	\end{align}
	Since $\frac{\p}{\p\ubar}\hat{\gsphere}$ is traceless, we have:
	\begin{align}\label{EQ:CharData4}
		\hchi&=\frac{1}{2}\Phi^2\frac{\p}{\p\ubar}\hat{\gsphere},&
		\gtr\upchi&=\frac{2}{\Phi}\frac{\p\Phi}{\p\ubar}.
	\end{align}
	Then, in coordinates \eqref{DE:Coordinatese3e4}, \eqref{EQ:e4upchi} reads
	\begin{align}\label{EQ:CharData5}
		\frac{\p^2\Phi}{\p\ubar^2}+\frac{1}{2}\sgabs{\hchi}^2+\FAlphaAB=0,
	\end{align}
	where 
	\begin{align}\label{EQ:CharData6}
		\sgabs{\hchi}^2=\frac{1}{4}\left(\hat{\gsphere}^{-1}\right)^{AC}\cdot\left(\hat{\gsphere}^{-1}\right)^{BD}\cdot\frac{\p}{\p\ubar}\hat{\gsphere}_{AB}\cdot\frac{\p}{\p\ubar}\hat{\gsphere}_{CD}.
	\end{align}
	Then, we call the smooth data set $(\hchi,\gtr\upchi,\lnenth,\ent,v_4,v_A)$ on $\hypersur_{\uinf}$ regular if $(\hchi,\gtr\upchi)$ is of the form \eqref{EQ:CharData4} with constraint\ equation \eqref{EQ:CharData5} satisfied with initial data $\Phi|_{S_{\uinf,0}}=1$ and $\left.\frac{\p\Phi}{\p\ubar}\right|_{S_{\uinf,0}}=\frac{1}{\abs{\uinf}}$ as in \eqref{EQ:CharData1.5.1}. We note that given $\sgabs{\hchi}^2$ and $\FAlphaAB$, the ODE \eqref{EQ:CharData5} for $\Phi$ is well-posed (at least locally).
	
	Similarly, we call the smooth data set $(\huchi,\gtr\uchi,\lnenth,\ent,v_3,v_A)$ on $\hypersurbar_{0}$ regular if 
	\begin{align}\label{EQ:CharData9}
		\huchi&=\frac{1}{2}\Phi^2\frac{\p}{\p u}\hat{\gsphere},&
		\gtr\uchi&=\frac{2}{\Phi}\frac{\p\Phi}{\p u},
	\end{align}
	and the following constraint equation is satisfied with initial data $\Phi|_{S_{\uinf,0}}=1$ and $\left.\frac{\p\Phi}{\p u}\right|_{S_{\uinf,0}}=-\frac{1}{\abs{\uinf}}$ as in \eqref{EQ:CharData1.5.1}:
	\begin{align}\label{EQ:CharData7}
		\frac{\p^2\Phi}{\p u^2}+\frac{1}{2}\sgabs{\huchi}^2+\FAlphabar=0,
	\end{align}
	where 
	\begin{align}\label{EQ:CharData8}
		\sgabs{\huchi}^2=\frac{1}{4}\left(\hat{\gsphere}^{-1}\right)^{AC}\cdot\left(\hat{\gsphere}^{-1}\right)^{BD}\cdot\frac{\p}{\p u}\hat{\gsphere}_{AB}\cdot\frac{\p}{\p u}\hat{\gsphere}_{CD}.
	\end{align}
	
	We note that, fluid variables are allowed to be prescribed freely on $\hypersur_{\uinf}\cup\hypersurbar_{0}$, since $\gfour$-null hypersurfaces are spacelike with respect to the acoustic geometry\footnote{For the acoustical metric in Definition \ref{DE:acousticunrescaled}, by straightforward computation, $\gfour_{Acou}(e_3,e_3),\gfour_{Acou}(e_4,e_4)>0$. We note that in this article, the proof does not rely on the acoustical metric. Hence we only provide its very definition without proving that $\gfour_{Acou}$ defined in Definition \ref{DE:acousticunrescaled} indeed governs the dynamics of sound waves.} simply due to the fact that the speed of sound is strictly smaller than the speed of light. With prescribed fluid variables (hence Ricci curvature components), we refer readers to \cite[Chapter 2]{Christodoulou2009} and \cite[Section 2-3]{Luk2012} for detailed data prescription/constructions.
	\subsubsection{Estimate of the characteristic data on $\hypersur_{\uinf}$}\label{SSS:DataonH}
	In this section, we impose the so-called short pulse data on $\hypersur_{\uinf}$ and estimate the size of the geometric quantities. 
	Let $a$ be a large number associated with the pulse size of the gravitational field, $b$ be a large number associated with the pulse size of the fluid, and $\delta$ be a small number such that $\delta a, \delta b\leq \abs{\uts}\leq \abs{\uinf}$ corresponding to the region (length of $\ubar$) where pulse is injected. Let $0<\iniBAfluid<a^{-\frac{1}{5}}$ be a positive number. Let $\vec{m}:=(i,j,k)$ with $1\leq \abs{\vec{m}}\leq N$ and $\WangD^{\vec{m}}$ be defined as in Definition \ref{DE:operators}. We assume that there is a regular data set $(\hchi,\gtr\upchi,\lnenth,\ent,v_4,v_A)$ on $\hypersur_{\uinf}$ such that 
	\begin{subequations}\label{EQ:DataonH1}
		\begin{align}
			\abs{\lnenth,\ent, v_A}&\leq \iniBAfluid,&
			\abs{\WangD^{\vec{m}}\left(\lnenth,\ent, v_A\right)}&\leq \iniBAfluid,\\
			b^{\frac{1}{2}}(1-\iniBAfluid)\leq -v_4&\leq b^{\frac{1}{2}}(1+\iniBAfluid),&
			\abs{\WangD^{\vec{m}}(v_4)}&\leq \iniBAfluid b^{\frac{1}{2}}.
		\end{align}
	\end{subequations}
	Moreover, $\gfour(v,v)=-1$ allows us to assume that
	\begin{align}
		b^{-\frac{1}{2}}(1-\iniBAfluid)\leq-v_3&\leq b^{-\frac{1}{2}}(1+\iniBAfluid),&
		\abs{\WangD^{\vec{m}}(v_3)}&\leq \iniBAfluid b^{-\frac{1}{2}}.
	\end{align}
	We make a technical assumption\footnote{This assumption is equivalent to asking for a close-to-homogeneous matter distribution on $\hypersur_{\uinf}$ with total amount of substance $\leq 1$.} on the equation of state on $\hypersur_{\uinf}$, $\pressure=\pressure(\enth,\ent)$ such that 
	\begin{align}
		n=&~\left.\frac{\p \pressure}{\p\enth}\right\vert_\ent(\uinf,\ubar,\sphere^1,\sphere^2)\lesssim\frac{1}{\abs{\uinf}^2},&
		\text{for all }& (\ubar,\sphere^1,\sphere^2)\in[0,(1+2(\uts-\uinf))\delta]\times\Stwo.
	\end{align}  Then, \eqref{EQ:DataonH1} implies
	\begin{align}\label{EQ:DataonH2}
		\abs{\left(\pressure(\lnenth,\ent),\density(\lnenth,\ent)\right)}&\leq \frac{1}{\abs{\uinf}^2},&
		\abs{\WangD^{\vec{m}}\left(\pressure(\lnenth,\ent),\density(\lnenth,\ent)\right)}&\lesssim \iniBAfluid\frac{1}{\abs{\uinf}^2}.
	\end{align}
	Then, we assume that $(\hchi,\gtr\upchi)$ satisfies \eqref{EQ:CharData4}-\eqref{EQ:CharData6} and 		
	\begin{align}\label{EQ:DataonH3}	
		\abs{\WangD^{\vec{m}}\hchi}&\lesssim \frac{a^{\frac{1}{2}}}{\abs{\uinf}},&
		\abs{\WangD^{\vec{m}}\left(\gtr\upchi-\frac{2}{\abs{\uinf}}\right)}&\lesssim \frac{1}{\abs{\uinf}^2}.
	\end{align}
	To be specific, we assume $\sgabs{\hchi}^2\lesssim \frac{a}{\abs{\uinf}^2}$ and \eqref{EQ:DataonH3} is a consequence of \eqref{EQ:CharData4}, \eqref{EQ:CharData5} and \eqref{EQ:DataonH2}.
	
	Under the above assumptions, in the rest of the subsection, we estimate the geometric quantities on $\hypersur_{\uinf}$.
	
	First of all, $\omega=0$ simply due to \eqref{EQ:omegauomega} and $\lapse=1$ on $\hypersur_{\uinf}$. 
	
	Using Proposition \ref{PR:Ricci0}, we have the following estimates for the Ricci curvature tensor components defined in Definition \ref{DE:CurvatureComponents}:
	\begin{subequations}\label{EQ:DataonH2.5}
		\begin{align}
			\label{EQ:DataonH2.5.1}\abs{\FAlphaAB}&\lesssim\frac{b}{\abs{\uinf}^2},&
			\abs{\WangD^{\vec{m}}\FAlphaAB}&\lesssim \iniBAfluid \frac{b}{\abs{\uinf}^2},&
			\abs{\WangD^{\vec{m}}\FBetaA}&\lesssim\iniBAfluid \frac{b^{\frac{1}{2}}}{\abs{\uinf}^2},\\
			\label{EQ:DataonH2.5.2}\abs{\FAlphabar}&\lesssim \frac{b^{-1}}{\abs{\uinf}^2},&
			\abs{\WangD^{\vec{m}}\FAlphabar}&\lesssim \iniBAfluid \frac{b^{-1}}{\abs{\uinf}^2},&
			\abs{\WangD^{\vec{m}}\FBetabar}&\leq \iniBAfluid \frac{b^{-\frac{1}{2}}}{\abs{\uinf}^2},\\
			\label{EQ:DataonH2.5.3}\abs{\FRho,\FRhodual,\Scalar}&\leq  \frac{1}{\abs{\uinf}^2},&
			\abs{\WangD^{\vec{m}}(\FRho,\FRhodual,\Scalar)}&\leq \iniBAfluid \frac{1}{\abs{\uinf}^2}.
		\end{align}
	\end{subequations}
	
	Also, the Gauss curvature $\GaussCurvature$ of $\gsphere$ satisfies 
	\begin{align}\label{EQ:DataonH4}
		\abs{\WangD^{\vec{m}}\left(\GaussCurvature-\frac{1}{\abs{\uinf}^2}\right)}\lesssim \frac{\delta a^{\frac{1}{2}}}{\abs{\uinf}^3}.
	\end{align}
	
	Then, integrating \eqref{Lzeta} over $\ubar$ and substituting for $\BetaA$ by \eqref{angdivhchi}, using \eqref{EQ:DataonH2.5.1} for $\FBetaA$ and the initial condition \eqref{EQ:CharData1.5.2}, we derive
	\begin{align}\label{EQ:DataonH5}
		\abs{\WangD^{\vec{m}}\upzeta}\lesssim\frac{\delta a^{\frac{1}{2}}}{\abs{\uinf}^2}+\iniBAfluid\frac{\delta b^{\frac{1}{2}}}{\abs{\uinf}^2}.
	\end{align}
	We note that by \eqref{EQ:etaupeta} and $\lapse=1$ on $\hypersur_{\uinf}$, $\upzeta=\upeta=-\ueta$.
	
	Integrating \eqref{Luchi} over $\ubar$ and substituting\footnote{See \eqref{EQ:Renormalization}-\eqref{EQ:Renormalization2} for the derivation of \eqref{EQ:DataonH6}.}
	\begin{align}\label{EQ:DataonH6}
		2\rho-\hchi\cdot\huchi=-2\GaussCurvature-\frac{1}{2}\gtr\upchi\gtr\uchi+2\FRho+\frac{2}{3}\Scalar,
	\end{align} 
	using \eqref{EQ:CharData9}, \eqref{EQ:DataonH2.5.3} for $\FRho,\Scalar$, and initial condition \eqref{EQ:CharData1.5.1}, we derive
	\begin{align}\label{EQ:DataonH7}
		\abs{\WangD^{\vec{m}}\left(\gtr\uchi+\frac{2}{\abs{\uinf}}\right)}\lesssim\frac{\delta}{\abs{\uinf}^2}+\iniBAfluid\frac{\delta}{\abs{\uinf}^2}.
	\end{align}
	
	Integrating \eqref{Lhuchi} over $\ubar$, using \eqref{EQ:CharData9}, \eqref{EQ:DataonH2.5.3} for $\FRhodual$, and the initial condition \eqref{EQ:CharData1.5.1}, we derive
	\begin{align}\label{EQ:DataonH8}
		\abs{\WangD^{\vec{m}}\huchi}\lesssim\frac{\delta a^{\frac{1}{2}}}{\abs{\uinf}^2}+\iniBAfluid\frac{\delta}{\abs{\uinf}^2}.
	\end{align}
	
	Now, we estimate the rest of the Weyl curvature components. Using data \eqref{EQ:CharData1.5}, previously proved results \eqref{EQ:DataonH3}, \eqref{EQ:DataonH5}, by \eqref{EQ:angDhchi} and \eqref{angdivhchi}, we have
	\begin{subequations}\label{EQ:DataonH9}
		\begin{align}
			\abs{\WangD^{\vec{m}}\AlphaAB}&\lesssim\frac{ a^{\frac{5}{4}}}{\abs{\uinf}},&
			\abs{\WangD^{\vec{m}}\CAlphaAB}&\lesssim\frac{ a^{\frac{5}{4}}}{\abs{\uinf}}+\iniBAfluid\frac{b}{\abs{\uinf}^2},\\
			\abs{\WangD^{\vec{m}}\BetaA}&\lesssim\frac{ a^{\frac{1}{2}}}{\abs{\uinf}^2},&
			\abs{\WangD^{\vec{m}}\CBetaA}&\lesssim\frac{ a^{\frac{1}{2}}}{\abs{\uinf}^2}+\iniBAfluid\frac{1}{\abs{\uinf}^2}.
		\end{align}
	\end{subequations}
	
	By the expressions of $\Rho,\CRho$ in \eqref{EQ:DataonH6}, \eqref{EQ:Renormalization} and previously proved results \eqref{EQ:DataonH3}, \eqref{EQ:DataonH4}, \eqref{EQ:DataonH7}, \eqref{EQ:DataonH8}, \eqref{EQ:DataonH2.5.3}, we derive:
	\begin{align}\label{EQ:DataonH10}
		\abs{\WangD^{\vec{m}}\Rho}&\lesssim\frac{1}{\abs{\uinf}^3}+\iniBAfluid\frac{1}{\abs{\uinf}^2},&
		\abs{\WangD^{\vec{m}}\CRho}&\lesssim\frac{\delta a^{\frac{1}{2}}}{\abs{\uinf}^3}+\iniBAfluid\frac{\delta}{\abs{\uinf}^3}.
	\end{align}
	
	By \eqref{curleta}, \eqref{EQ:Renormalization}, and previously proved results \eqref{EQ:DataonH3}, \eqref{EQ:DataonH5}, \eqref{EQ:DataonH7}, \eqref{EQ:DataonH8}, \eqref{EQ:DataonH2.5.3}, we obtain:
	\begin{align}\label{EQ:DataonH11}
		\abs{\WangD^{\vec{m}}\Rhodual}&\lesssim\frac{1}{\abs{\uinf}^3}+\iniBAfluid\frac{\delta b^{\frac{1}{2}}}{\abs{\uinf}^3},&
		\abs{\WangD^{\vec{m}}\CRhodual}&\lesssim\frac{\delta a^{\frac{1}{2}}}{\abs{\uinf}^3}+\iniBAfluid\frac{\delta b^{\frac{1}{2}}}{\abs{\uinf}^3}.
	\end{align}
	
	Integrating \eqref{EQ:D4CBetabar} over $\ubar$, using previously proved results \eqref{EQ:DataonH2.5}-\eqref{EQ:DataonH8}, and the initial condition \eqref{EQ:CharData1.5.2}, we derive:
	\begin{align}\label{EQ:DataonH12}
		\abs{\WangD^{\vec{m}}\CBetabar}&\lesssim\frac{\delta }{\abs{\uinf}^3}+\iniBAfluid\frac{\delta}{\abs{\uinf}^3},&
		\abs{\WangD^{\vec{m}}\Betabar}&\lesssim\frac{\delta }{\abs{\uinf}^3}+\iniBAfluid\frac{ b^{-\frac{1}{2}}}{\abs{\uinf}^3}.
	\end{align}
	
	Integrating \eqref{EQ:D4CAlphabar} over $\ubar$, using previously proved results \eqref{EQ:DataonH2.5}-\eqref{EQ:DataonH12}, and the initial condition \eqref{EQ:CharData1.5.2}, we have:
	\begin{align}\label{EQ:DataonH13}
		\abs{\WangD^{\vec{m}}\CAlphabar}&\lesssim\frac{\delta^{\frac{3}{2}} }{\abs{\uinf}^5}+\iniBAfluid\frac{\delta }{\abs{\uinf}^3},&
		\abs{\WangD^{\vec{m}}\Alphabar}&\lesssim\frac{\delta^{\frac{3}{2}} }{\abs{\uinf}^5}+\iniBAfluid\frac{\delta }{\abs{\uinf}^3}.
	\end{align}
	
	For the remaining connection coefficient $\uomega$, integrating \eqref{Luomega} over $\ubar$, using previously proved results \eqref{EQ:DataonH2.5}, \eqref{EQ:DataonH5}, and \eqref{EQ:DataonH10}, and the initial condition \eqref{EQ:CharData1.5.2}, we derive:
	\begin{align}\label{EQ:DataonH14}
		\abs{\WangD^{\vec{m}}\uomega}&\lesssim\frac{\delta }{\abs{\uinf}^3}+\iniBAfluid\frac{\delta }{\abs{\uinf}^2}.
	\end{align}
	
	\subsubsection{Estimate of the characteristic data on $\hypersurbar_{0}$}\label{SSS:DataonHbar}
	In this section, we impose the characteristic data on the incoming initial null hypersurface $\hypersurbar_{0}$ and estimate the size of the geometric quantities. Compared to the Einstein vacuum, Einstein-scalar-field, and Einstein-Maxwell in \cite{AnLuk2017,Christodoulou2009,AnAthanasiou2022,ShenWan20251,Yu2013,Yu2011}, where the spacetime on $\hypersurbar_{0}$ is assumed to be Minkowskian, $\pressure,\density,\enth$ are non-trivial in our setting to avoid degeneracy of the fluid part of the Einstein-Euler equations \eqref{EQ:FOE1}. Therefore, we need to carefully prescribe initial data on $\hypersurbar_{0}$ and show the data constructed/prescribed is suitable for forming a trapped surface dynamically. Moreover, from the structural equations, we can rule out the possibility that the trapped surfaces can form on $\hypersurbar_{0}$ within our short pulse ansatz.
	
	Let $a,b,\delta,\iniBAfluid$ be the parameters as in Section \ref{SSS:DataonH}. Let $\vec{m}:=(i,j,k)$ with $1\leq\abs{\vec{m}}\leq N$ and $\WangD^{\vec{m}}$ be defined as in Definition \ref{DE:operators}. We assume that there is a regular data set $(\huchi,\gtr\uchi,\lnenth,\ent,v_3,v_A)$ on $\hypersurbar_{0}$ such that the following bounds hold for fluid variables on $\hypersurbar_{0}$:
	\begin{subequations}\label{EQ:DataonHbar1}
		\begin{align}
			\abs{\lnenth,\ent, v_A}&\leq \iniBAfluid,&
			\abs{\WangD^{\vec{m}}\left(\lnenth,\ent, v_A\right)}&\leq \iniBAfluid,\\
			b^{-\frac{1}{2}}(1-\iniBAfluid)\leq -v_3&\leq b^{-\frac{1}{2}}(1+\iniBAfluid),&
			\abs{\WangD^{\vec{m}}(v_3)}&\leq \iniBAfluid b^{-\frac{1}{2}}.
		\end{align}
	\end{subequations}
	In addition, $\gfour(v,v)=-1$ allows us to assume that
	\begin{align}
		b^{\frac{1}{2}}(1-\iniBAfluid)\leq -v_4&\leq b^{\frac{1}{2}}(1+\iniBAfluid),&
		\abs{\WangD^{\vec{m}}(v_4)}&\leq \iniBAfluid b^{\frac{1}{2}}.
	\end{align}
	Moreover, we assume that the corner compatibility condition is satisfied. That is, the fluid equations \eqref{EQ:FOE2} are satisfied on $S_{\uinf,0}$ at all orders. In fact, in the construction of data in Section \ref{SSS:DataonH}, we prescribed $e_3$-derivatives of all the geometric quantities (including fluid variables) at all orders\footnote{In our analysis on trapped surface formation, we only need to differentiate geometric quantities $N$ times. However, the current techniques in the proof of local well-posedness of the characteristic initial value problem require the initial data being smooth (due to the Whitney extension); see \cite{Rendall1990}. Hence, we requires our regular data set and all geometric quantities to be smooth.} and the prescribed data satisfy \eqref{EQ:FOE2} at all orders on $S_{\uinf,0}$. Here, we only require the imposed data on $\hypersurbar_{0}$ to inherit the properties (due to the data construction on $\hypersur_{\uinf}$) on $S_{\uinf,0}$.
	
	Similar to Section \ref{SSS:DataonH}, we make the assumption on the equation of state on $\hypersurbar_{0}$, $\pressure=\pressure(\enth,\ent)$ such that $n=\left.\frac{\p \pressure}{\p\enth}\right\vert_\ent(u,0,\sphere^1,\sphere^2)\lesssim\frac{1}{\abs{u}^2}\leq\frac{1}{\abs{\uts}^2}$ for all $(\sphere^1,\sphere^2)\in\Stwo$. Then, \eqref{EQ:DataonH1} implies
	\begin{align}\label{EQ:DataonHbar2}
		\abs{\left(\pressure(\lnenth,\ent),\density(\lnenth,\ent)\right)}&\leq \frac{1}{\abs{u}^2},&
		\abs{\WangD^{\vec{m}}\left(\pressure(\lnenth,\ent),\density(\lnenth,\ent)\right)}&\lesssim \iniBAfluid\frac{1}{\abs{u}^2}.
	\end{align}
	Then, we assume that $(\huchi,\gtr\uchi)$ satisfies \eqref{EQ:CharData9}-\eqref{EQ:CharData8} and 		
	\begin{align}\label{EQ:DataonHbar3}	
		\abs{\WangD^{\vec{m}}\huchi}&\leq \frac{\delta}{\abs{u}^2},&
		\abs{\WangD^{\vec{m}}\left(\gtr\uchi+\frac{2}{\abs{u}}\right)}&\lesssim \frac{\delta^2}{\abs{u}^2}+\iniBAfluid\frac{b^{-1}}{\abs{u}^2}.
	\end{align}
	To be specific, we assume $\sgabs{\huchi}^2\lesssim \frac{\delta^2} {\abs{u}^4}$ and \eqref{EQ:DataonHbar3} is a consequence of $\abs{u}\leq\abs{\uinf}$, \eqref{EQ:CharData9}, \eqref{EQ:CharData7} and \eqref{EQ:DataonHbar2}.
	
	Under the above assumptions, in the rest of the subsection, we estimate the geometric quantities on $\hypersurbar_{0}$.
	
	First of all, $\uomega=0$ simply due to \eqref{EQ:omegauomega} and $\lapse=1$ on $\hypersurbar_{0}$. 
	
	Using Proposition \ref{PR:Ricci0}, we have the following estimates for the Ricci curvature tensor components defined in Definition \ref{DE:CurvatureComponents}:
	\begin{subequations}\label{EQ:DataonHbar2.5}
		\begin{align}
			\label{EQ:DataonHbar2.5.1}\abs{\FAlphaAB}&\leq\frac{b}{\abs{u}^2},&
			\abs{\WangD^{\vec{m}}\FAlphaAB}&\leq \iniBAfluid \frac{b}{\abs{u}^2},&
			\abs{\WangD^{\vec{m}}\FBetaA}&\leq \iniBAfluid \frac{b^{\frac{1}{2}}}{\abs{u}^2},\\
			\label{EQ:DataonHbar2.5.2}\abs{\FAlphabar}&\leq  \frac{b^{-1}}{\abs{u}^2},&
			\abs{\WangD^{\vec{m}}\FAlphabar}&\leq \iniBAfluid \frac{b^{-1}}{\abs{u}^2},&
			\abs{\WangD^{\vec{m}}\FBetabar}&\leq \iniBAfluid \frac{b^{-\frac{1}{2}}}{\abs{u}^2},\\
			\label{EQ:DataonHbar2.5.3}\abs{\FRho,\FRhodual,\Scalar}&\leq \frac{1}{\abs{u}^2},&
			\abs{\WangD^{\vec{m}}(\FRho,\FRhodual,\Scalar)}&\leq \iniBAfluid \frac{1}{\abs{u}^2}.
		\end{align}
	\end{subequations}
	
	Also, the Gauss curvature $\GaussCurvature$ of $\gsphere$ satisfies 
	\begin{align}\label{EQ:DataonHbar4}
		\abs{\WangD^{\vec{m}}\left(\GaussCurvature-\frac{1}{\abs{u}^2}\right)}\lesssim \frac{\delta a^{\frac{1}{2}}}{\abs{u}^3}.
	\end{align}
	
	Then, integrating \eqref{uLzeta} over $u$ and substituting for $\Betabar$ by \eqref{angdivhuchi}, using \eqref{EQ:DataonHbar2.5.1} for $\FBetabar$ and the initial condition \eqref{EQ:CharData1.5.2}, we derive
	\begin{align}\label{EQ:DataonHbar5}
		\abs{\WangD^{\vec{m}}\upzeta}\lesssim\frac{\delta a^{\frac{1}{2}}}{\abs{u}^2}+\iniBAfluid\frac{ b^{-\frac{1}{2}}}{\abs{u}}.
	\end{align}
	We note that by \eqref{EQ:etaupeta} and $\lapse=1$ on $\hypersurbar_{0}$, $\upzeta=\upeta=-\ueta$.
	
	Integrating \eqref{e3gtrupchi} over $u$ and using \eqref{EQ:DataonH6}, \eqref{EQ:CharData9}, \eqref{EQ:DataonHbar2.5.3} for $\FRho,\Scalar$, and the  initial condition \eqref{EQ:CharData1.5.1}, we derive
	\begin{align}\label{EQ:DataonHbar7}
		\abs{\WangD^{\vec{m}}\left(\gtr\upchi-\frac{2}{\abs{u}}\right)}\leq\frac{1}{\abs{u}^2}+\frac{\delta a^{\frac{1}{2}}}{\abs{u}^2}+\frac{ b^{-\frac{1}{2}}}{\abs{u}}.
	\end{align}
	We note that the largest term, the first term, on the RHS of \eqref{EQ:DataonHbar7} is contributed by the term $2\Rho+\frac{1}{3}\Scalar+\cdots=2\FRho+\Scalar=\gtr\angRicfour=(\density+\pressure)\angv\cdot\angv+\density-\pressure$ (see Proposition \ref{PR:Ricci0}), where $\cdots$ represents the terms with smaller size (of size $\lesssim\delta^{\frac{1}{2}}$). We emphasize two important aspects regarding the term $(\density+\pressure)\angv\cdot\angv$: 1) This term is positive since $\density,\pressure>0$, $\gsphere$ is Riemannian, and $\density-\pressure>0$ due to the speed of sound is less than the speed of light $\speed^2:=\frac{\p\pressure}{\p\density}<1$ and $\pressure|_{\density=0}=0$. Hence under the assumptions of data, \textbf{$\gtr\upchi$ is strictly monotonically increasing along $u$, which suggests that $\gtr\upchi$ remains positive on the incoming initial hypersurface $\hypersurbar_{0}$.} 2) $\gtr\upchi$ does not grow too large on $\hypersurbar_{0}$ ($\gtr\upchi<\frac{3}{\abs{u}}$ for all $u\in[\uinf,\uts]$). Otherwise, it will be difficult for space-time to form a trapped surface at the preferred location $S_{\uts,a(\delta)}$, since in our analysis (as in Section \ref{S:Formation}), $\gtr\upchi|_{S_{\uts,a(\delta)}}\sim \gtr\upchi|_{S_{\uts,0}}-\frac{4}{\abs{\uts}}$.
	
	Integrating \eqref{EQ:angD3hchi} over $u$, using \eqref{EQ:CharData4}, \eqref{EQ:DataonHbar2.5.3} for $\FRhodual$, and the initial condition \eqref{EQ:CharData1.5.1}, we derive
	\begin{align}\label{EQ:DataonHbar8}
		\abs{\WangD^{\vec{m}}\hchi}\lesssim\frac{\delta a^{\frac{1}{2}}}{\abs{u}^2}+\iniBAfluid\frac{1}{\abs{u}}.
	\end{align}
	
	Now, we estimate the rest of the Weyl curvature components. Using data \eqref{EQ:CharData1.5}, previously proved results \eqref{EQ:DataonHbar3}, \eqref{EQ:DataonHbar5}, by \eqref{e4huchi} and \eqref{angdivhuchi}, we have
	\begin{subequations}\label{EQ:DataonHbar9}
		\begin{align}
			\abs{\WangD^{\vec{m}}\Alphabar}&\lesssim\frac{\delta  }{\abs{u}^3},&
			\abs{\WangD^{\vec{m}}\CAlphabar}&\lesssim\frac{\delta  }{\abs{u}^3}+\iniBAfluid\frac{b^{-1}}{\abs{u}^2},\\
			\abs{\WangD^{\vec{m}}\Betabar}&\lesssim\frac{\delta  a^{\frac{1}{2}}}{\abs{u}^3},&
			\abs{\WangD^{\vec{m}}\CBetabar}&\lesssim\frac{\delta a^{\frac{1}{2}}}{\abs{u}^3}+\iniBAfluid\frac{b^{-\frac{1}{2}}}{\abs{u}^2}.
		\end{align}
	\end{subequations}
	
	By the expressions for $\Rho,\CRho$ in \eqref{EQ:DataonH6}, \eqref{EQ:Renormalization} and previously proved results \eqref{EQ:DataonHbar3}, \eqref{EQ:DataonHbar4}, \eqref{EQ:DataonHbar7}, \eqref{EQ:DataonHbar8}, \eqref{EQ:DataonHbar2.5.3}, we derive:
	\begin{align}\label{EQ:DataonHbar10}
		\abs{\WangD^{\vec{m}}\Rho}&\lesssim\frac{\delta a^{\frac{1}{2}}}{\abs{u}^3}+\iniBAfluid\frac{1}{\abs{u}^2},&
		\abs{\WangD^{\vec{m}}\CRho}&\lesssim\frac{\delta a^{\frac{1}{2}}}{\abs{\uinf}^3}+\iniBAfluid\frac{b^{-\frac{1}{2}}}{\abs{u}^3}.
	\end{align}
	
	By \eqref{curleta}, \eqref{EQ:Renormalization}, and previously proved results \eqref{EQ:DataonHbar3}, \eqref{EQ:DataonHbar5}, \eqref{EQ:DataonHbar7}, \eqref{EQ:DataonHbar8}, \eqref{EQ:DataonHbar2.5.3}, we obtain:
	\begin{align}\label{EQ:DataonHbar11}
		\abs{\WangD^{\vec{m}}\Rhodual}&\lesssim\frac{\delta}{\abs{u}^3}+\iniBAfluid\frac{ b^{-\frac{1}{2}}}{\abs{u}^3},&
		\abs{\WangD^{\vec{m}}\CRhodual}&\lesssim\frac{\delta }{\abs{\uinf}^3}+\iniBAfluid\frac{ b^{-\frac{1}{2}}}{\abs{u}^3}.
	\end{align}
	
	Integrating \eqref{EQ:D3CBeta} over $u$, using previously proved results \eqref{EQ:DataonHbar2.5}-\eqref{EQ:DataonHbar8}, and the initial condition \eqref{EQ:CharData1.5.2}, we derive:
	\begin{align}\label{EQ:DataonHbar12}
		\abs{\WangD^{\vec{m}}\CBetaA}&\lesssim\frac{\delta }{\abs{u}^3}+\iniBAfluid\frac{\delta a^{\frac{1}{2}} b^{\frac{1}{2}}}{\abs{u}^2}+\iniBAfluid\frac{1}{\abs{u}},&
		\abs{\WangD^{\vec{m}}\BetaA}&\lesssim\frac{\delta }{\abs{u}^3}+\iniBAfluid\frac{ b^{\frac{1}{2}}}{\abs{u}^3}.
	\end{align}
	
	Integrating \eqref{EQ:D3CAlpha} over $u$, using previously proved results \eqref{EQ:DataonHbar2.5}-\eqref{EQ:DataonHbar12}, and the initial condition \eqref{EQ:CharData1.5.2}, we have:
	\begin{align}\label{EQ:DataonHbar13}
		\abs{\WangD^{\vec{m}}\CAlphaAB}&\lesssim\frac{\delta} {\abs{u}^3}+\iniBAfluid\frac{\delta^{-1} }{\abs{u}},&
		\abs{\WangD^{\vec{m}}\AlphaAB}&\lesssim\frac{\delta }{\abs{u}^3}+\iniBAfluid\frac{\delta^{-1} }{\abs{u}}+\iniBAfluid\frac{b }{\abs{u}^2}.
	\end{align}
	
	For the remaining connection coefficient $\omega$, integrating \eqref{Lbaromega} over $u$, using previously proved results \eqref{EQ:DataonHbar2.5}, \eqref{EQ:DataonHbar5}, and \eqref{EQ:DataonHbar10}, and the initial condition \eqref{EQ:CharData1.5.2}, we derive:
	\begin{align}\label{EQ:DataonHbar14}
		\abs{\WangD^{\vec{m}}\omega}&\lesssim\frac{\delta a^{\frac{1}{2}} }{\abs{u}^2}+\iniBAfluid\frac{1 }{\abs{u}}.
	\end{align}
	\subsection{Initial Data Assumptions}\label{SS:InitialData}
	Recall Section \ref{SS:ChoiceofParameters} for the choice of parameters $\iniBAfluid,\iniBAsize$. Let $\mathring{\mathcal{D}}_{\inihypersur}:=\{\uinf\}\times[0,\delta]$, $\mathring{\mathcal{D}}_{\inihypersurbar}:=[\uinf,\uts]\times\{0\}$, and $\mathring{\mathcal{D}}:=\mathring{\mathcal{D}}_{\inihypersur}\cup\mathring{\mathcal{D}}_{\inihypersurbar}$. We make the following initial data assumptions where $\mathcal{R},\mathcal{O},\mathcal{E}$ are defined in Section \ref{SS:Norms}:
	\begin{subequations}\label{Data}
		\begin{align}
			\label{Data11}\mathcal{R}_{\mathring{\mathcal{D}},l,m,n,\infty}\leq&~\iniBAsize,&
			0\leq l+m+n\leq&~ \frac{N-1}{2},\\
			\label{Data12}\mathcal{R}_{\mathring{\mathcal{D}},i,j,k,2}\leq&~\iniBAsize,&
			0\leq i+j+k\leq&~ N-1,\\
			\label{Data21}\mathcal{O}_{\mathring{\mathcal{D}},l,m,n,\infty}\leq&~\iniBAsize,&
			0\leq l+m+n\leq&~ \frac{N}{2},\\
			\label{Data22}\mathcal{O}_{\mathring{\mathcal{D}},i,j,k,2}\leq&~\iniBAsize,&
			0\leq i+j+k\leq&~ N,\\
			\label{Data31}\mathcal{E}_{\mathring{\mathcal{D}},l,m,n,\infty}\leq&~\iniBAfluid,&
			1\leq l+m+n\leq&~ \frac{N-1}{2},\\
			\label{Data32}\mathcal{E}_{\mathring{\mathcal{D}},i,j,k,2}\leq&~\iniBAfluid,&
			1\leq i+j+k\leq&~ N,\\
			\label{Data4}(1-\iniBAfluid)\leq \abs{b^{\frac{1}{2}}v_3,b^{-\frac{1}{2}}v_4}\leq&~(1+\iniBAfluid),&
			\abs{\angv,\lnenth,\ent}\leq\iniBAfluid.
		\end{align}
	\end{subequations}
	Moreover, on $\hypersurbar_0$, we assume that 
	\begin{align}\label{AS:SmallonHbar}
		\norm{\gtr\upchi-\frac{2}{\abs{u}}}_{L^\infty(S_{u,0})}&\leq \frac{1}{\abs{u}}+\iniBAsize \delta a^{\frac{1}{2}}.
	\end{align}
	\begin{remark}
		Due to the construction of the characteristic data and structural equations in Proposition \ref{PR:NullStructure1}, Proposition \ref{PR:ReNullBianchi}, and Section \ref{SS:EulerEquations}, the data set that satisfies \eqref{Data}-\eqref{AS:SmallonHbar} exists. See \cite[Chapter 2]{Christodoulou2009} for detailed arguments on $\hypersur_{\uinf}$. 
	\end{remark}
	\subsection{Main Theorem}\label{SS:MainTheorem}
	\begin{theorem}[Trapped surface formation]\label{TH:Main}
		We assume that the data satisfy \eqref{Data} and: 
		\begin{align}\label{AS:LargenessofData}
			\int_{0}^{\delta}\abs{\uinf}^2(\abs{\hchi}^2+2\FAlphaAB)(\uinf,\ubarprime,\sphere^1,\sphere^2)\diff\ubarprime&< 2\abs{\uinf},&
			\text{for all } &(\sphere^1,\sphere^2)\in\Stwo.
		\end{align}
		Then, there is no trapped surface on $\hypersur_{\uinf}\cup\hypersurbar_{0}$. Moreover,
		\begin{enumerate}[label=\arabic*]
			\item If data satisfy:
			\begin{align}\label{AS:LargenessofData1}
				\int_{0}^{\delta}\abs{\uinf}^2\abs{\hchi}^2(\uinf,\ubarprime,\sphere^1,\sphere^2)\diff\ubarprime&>5\abs{\uts},&
				\text{for all } &(\sphere^1,\sphere^2)\in\Stwo,
			\end{align}
			then the $2$-sphere $S_{\uts,\delta}$ is a trapped surface. That is, $\gtr\upchi,\gtr\uchi<0$ on $S_{\uts,\delta}$.
			\item If data satisfy:
			\begin{align}
				\label{AS:LargenessofData2}\int_{0}^{\delta}\abs{\uinf}^2\FAlphaAB(\uinf,\ubarprime,\sphere^1,\sphere^2)\diff\ubarprime&>4\abs{\uts},&
				\text{for all } &(\sphere^1,\sphere^2)\in\Stwo,
			\end{align}
			then the $2$-sphere $S_{\uts,(1+2(\uts-\uinf))\delta}$ is a trapped surface. That is, $\gtr\upchi,\gtr\uchi<0$ on $S_{\uts,(1+2(\uts-\uinf))\delta}$.
		\end{enumerate}
		
	\end{theorem}
	\begin{remark}
		We note that by Proposition \ref{PR:Ricci0}, $\FAlphaAB$ is strictly positive. 
	\end{remark}
	\subsection{Bootstrap Assumptions}\label{SS:BA}
	Recall Section \ref{SS:ChoiceofParameters} for the choice of parameters $\iniBAfluid,\BAfluid,\iniBAsize,\BAsize$. Let $\mathcal{D}:=[\uinf,\uts]\times[0,(1+2(\uts-\uinf))\delta]$. We make the following bootstrap assumptions where $\mathcal{R},\mathcal{O},\mathcal{E}$ are defined in Section \ref{SS:Norms}:
	\begin{subequations}\label{BA}
		\begin{align}
			\label{BA11}\mathcal{R}_{\mathcal{D},l,m,n,\infty}\leq&~\BAsize,&
			0\leq l+m+n\leq&~ \frac{N-1}{2},\\
			\label{BA12}\mathcal{R}_{\mathcal{D},i,j,k,2}\leq&~\BAsize,&
			0\leq i+j+k\leq&~ N-2,\\
			\label{BA13}\tilde{\mathcal{R}}_{\mathcal{D},N-1,2}\leq&~\BAsize,&
			i+j+k=&~ N-1 \text{ and } k\geq 1,\\
			\label{BA21}\mathcal{O}_{\mathcal{D},l,m,n,\infty}\leq&~\BAsize,&
			0\leq  l+m+n\leq&~ \frac{N}{2},\\
			\label{BA22}\mathcal{O}_{\mathcal{D},i,j,k,2}\leq&~\BAsize,&
			0\leq i+j+k\leq&~ N-2 \text{ or } i+j+k=~ N-1 \text{ with } k\geq 1\\
			\label{BA23}\tilde{\mathcal{O}}_{\mathcal{D},N,2}\leq&~\BAsize,&
			i+j+k=&~ N \text{ and } k\geq 2,\\
			\label{BA31}\mathcal{E}_{\mathcal{D},l,m,n,\infty}\leq&~\BAfluid,&
			1\leq  l+m+n\leq&~ \frac{N}{2},\\
			\label{BA32}\mathcal{E}_{\mathcal{D},i,j,k,2}\leq&~\BAfluid,&
			1\leq i+j+k\leq&~ N-2 \text{ or } i+j+k=~ N-1 \text{ with }k\geq 1,\\
			\label{BA33}\tilde{\mathcal{E}}_{\mathcal{D},N,2}\leq&~\BAfluid,&
			1\leq i+j+k\leq&~ N \text{ with }k\geq 1,\\
			\label{BA4}(1-\BAfluid)\leq -b^{\frac{1}{2}}v_3,-b^{-\frac{1}{2}}v_4\leq&~(1+\BAfluid),&
			\abs{\angv,\lnenth,\ent}\leq&~\BAfluid.
		\end{align}
	\end{subequations}
	For convenience, we drop the ``$\mathcal{D}$" notation in \eqref{BA} in the rest of the article unless we emphasize that we consider the quantities in a different region.
	\section{Commutator Formula}
	\begin{proposition}\cite[Lemma 3.3, commutator formula for null frames]{KlainermanRodnianski2012}.\label{PR:commutator}
		Let $f$ be a scalar function. Then:
		\begin{subequations}\label{EQ:Commutator1}
			\begin{align}
				[\Dfour_3,\Dfour_4]f=&~-2\omega\Dfour_3f+2\uomega\Dfour_4f+4\upzeta\cdot\angnabla f,\\
				[\Dfour_4,\Dfour_A]f=&~\frac{1}{2}(\upeta_A+\ueta_A)\Dfour_4f-\upchi_{AB}\cdot\Dfour_B f,\\
				[\Dfour_3,\Dfour_A]f=&~\frac{1}{2}(\upeta_A+\ueta_A)\Dfour_3f-\uchi_{AB}\cdot\Dfour_B f,\\
				[\angD_3,\angD_4]f=&~-2\omega\angD_3f+2\uomega\angD_4f+4\upzeta\cdot\angnabla f,\\
				[\angD_4,\angnabla]f=&~\frac{1}{2}(\upeta+\ueta)\Dfour_4f-\upchi\cdot\angnabla f,\\
				[\angD_3,\angnabla]f=&~\frac{1}{2}(\upeta+\ueta)\Dfour_3f-\uchi\cdot\angnabla f.
			\end{align}
		\end{subequations}
		For any one-form $V$, we have:
		\begin{subequations}\label{EQ:Commutator2.1}
			\begin{align}
				[\Dfour_3,\Dfour_4]V_\mu=&~-2\omega\Dfour_3V_\mu+2\uomega\Dfour_4V_\mu+4\upzeta_C\Dfour_C V_\mu+\Riem_{\mu\alpha34}V^\alpha,\\
				[\Dfour_4,\Dfour_A]V_\mu=&~-\upchi_{AC}\Dfour_C V_\mu+\frac{1}{2}(\upeta_A+\ueta_A)\Dfour_4V_\mu+\Riem_{\mu\alpha4A}V^\alpha,\\
				[\Dfour_3,\Dfour_A]V_\mu=&~-\uchi_{AC}\Dfour_C V_\mu+\frac{1}{2}(\upeta_A+\ueta_A)\Dfour_3V_\mu+\Riem_{\mu\alpha3A}V^\alpha.
			\end{align}
		\end{subequations}
		Let $V$ be a one-form tangent to $\suubar$. We have:
		\begin{subequations}\label{EQ:Commutator2}
			\begin{align}
				[\angD_3,\angD_4]V_B=&~-2\omega\angD_3V_B+2\uomega\angD_4V_B+4\upzeta_C\angnabla_C V_B+2\Rhodual\antisymmetic_{BC}V_C+2\ueta_B\upeta_CV_C-2\upeta_B\ueta_CV_C,\\
				[\angD_4,\angnabla_A]V_B=&~-\upchi_{AC}\angnabla_C V_B+\antisymmetic_{AC}\HodgeLeftDual\BetaA_BV_C+\frac{1}{2}(\upeta_A+\ueta_A)\angD_4V_B-\upchi_{AC}\ueta_B V_C+\upchi_{AB}\ueta\cdot V\\
				\notag&+\FBetaA_AV_B-\FBetaA_BV_A,\\
				[\angD_3,\angnabla_A]V_B=&~-\uchi_{AC}\angnabla_C V_B+\antisymmetic_{AC}\HodgeLeftDual\Betabar_BV_C+\frac{1}{2}(\upeta_A+\ueta_A)\angD_3V_B-\uchi_{AC}\upeta_B V_C+\uchi_{AB}\upeta\cdot V\\
				\notag&+\FBetabar_AV_B-\FBetabar_BV_A.
			\end{align}
		\end{subequations}
	\end{proposition}
	\begin{proposition}[Commuting transport equations with $\Dfour$]\label{PR:CommuteDijk}
		Let $A=1,2$ and let $\vec{m_a}:=(i_a,j_a,k_a)$'s be multi-indices consisting of $3-$tuples. For any $\binom{a}{b}$-tensorfield $F$, we have the following schematic commutator formula:
		\begin{align}\label{EQ:CommutatorD1}
			\Dfour_3 \Dfour^i_4\Dfour_3^j\Dfour_A^{k}F_{\nu_1\cdots\nu_b}^{\mu_1\cdots\mu_a}=&~\sum\limits_{\vec{m_1}+\vec{m_2}+\vec{m_3}=(i-i_0,j,k-k_0)}\Dfour^{\vec{m_1}}(\upeta+\ueta)^{k_0}\cdot\Dfour^{\vec{m_2}}\omega^{i_0}\cdot\Dfour^{\vec{m_3}}\Dfour_3F_{\nu_1\cdots\nu_b}^{\mu_1\cdots\mu_a}\\
			\notag&+~\sum\limits_{\substack{\vec{m_1}+\vec{m_2}+\vec{m_3}+\vec{m_4}=(i-i_0,j,k-k_0)\\\vec{m_4}\neq(0,0,0)}}\Dfour^{\vec{m_1}}(\upeta+\ueta)^{k_0}\cdot\Dfour^{\vec{m_2}}\omega^{i_0}\cdot\Dfour^{\vec{m_3}}(\uchi,\uomega,\upzeta)\cdot\Dfour^{\vec{m_4}}F_{\nu_1\cdots\nu_b}^{\mu_1\cdots\mu_a}\\
			\notag&+\sum\limits_{\vec{m_1}+\vec{m_2}+\vec{m_3}+\vec{m_4}=(i-i_0,j,k-k_0-1)}\Dfour^{\vec{m_1}}(\upeta+\ueta)^{k_0}\cdot\Dfour^{\vec{m_2}}\omega^{i_0}\cdot\Dfour^{\vec{m_3}}\Riem^{\mu_i}_{\mu_03A}\cdot\Dfour^{\vec{m_4}}F_{\nu_1\cdots\nu_b}^{\hat{\mu}_i\mu_0}\\
			\notag&+\sum\limits_{\vec{m_1}+\vec{m_2}+\vec{m_3}+\vec{m_4}=(i-i_0,j,k-k_0-1)}\Dfour^{\vec{m_1}}(\upeta+\ueta)^{k_0}\cdot\Dfour^{\vec{m_2}}\omega^{i_0}\cdot\Dfour^{\vec{m_3}}\Riem^{\nu_0}_{\nu_i3A}\cdot\Dfour^{\vec{m_4}}F_{\hat{\nu_i}\nu_0}^{\mu_1\cdots\mu_a}\\
			\notag&+\sum\limits_{\vec{m_1}+\vec{m_2}+\vec{m_3}+\vec{m_4}=(i-i_0-1,j,k-k_0)}\Dfour^{\vec{m_1}}(\upeta+\ueta)^{k_0}\cdot\Dfour^{\vec{m_2}}\omega^{i_0}\cdot\Dfour^{\vec{m_3}}\Riem^{\mu_i}_{\mu_034}\cdot\Dfour^{\vec{m_4}}F_{\nu_1\cdots\nu_b}^{\hat{\mu}_i\mu_0}\\
			\notag&+\sum\limits_{\vec{m_1}+\vec{m_2}+\vec{m_3}+\vec{m_4}=(i-i_0-1,j,k-k_0)}\Dfour^{\vec{m_1}}(\upeta+\ueta)^{k_0}\cdot\Dfour^{\vec{m_2}}\omega^{i_0}\cdot\Dfour^{\vec{m_3}}\Riem^{\nu_0}_{\nu_i34}\cdot\Dfour^{\vec{m_4}}F_{\hat{\nu_i}\nu_0}^{\mu_1\cdots\mu_a},
		\end{align}
		where
		\begin{subequations}
			\begin{align}
				\Dfour^{i_1}\Psi^{i_2}:=&~\sum\limits_{j_1+\cdots+j_{i_2}=i_1}\Dfour^{j_1}\Psi\cdot\Dfour^{j_{i_2}}\Psi,\\
				F_{\hat{\nu_i}\nu_0}^{\mu_1\cdots\mu_a}:=&~F_{\nu_1\cdots\nu_{i-1}\nu_{i+1}\cdots\nu_b\nu_0}^{\mu_1\cdots\mu_a}
			\end{align}
		\end{subequations}
		Similarly, 
		\begin{align}\label{EQ:CommutatorD2}
			\Dfour_4 \Dfour^i_4\Dfour_3^j\Dfour_A^{k}F_{\nu_1\cdots\nu_b}^{\mu_1\cdots\mu_a}=&~\sum\limits_{\vec{m_1}+\vec{m_2}+\vec{m_3}=(i,j-j_0,k-k_0)}\Dfour^{\vec{m_1}}(\upeta+\ueta)^{k_0}\cdot\Dfour^{\vec{m_2}}\uomega^{j_0}\cdot\Dfour^{\vec{m_3}}\Dfour_4F_{\nu_1\cdots\nu_b}^{\mu_1\cdots\mu_a}\\
			\notag&+~\sum\limits_{\substack{\vec{m_1}+\vec{m_2}+\vec{m_3}+\vec{m_4}=(i-i_0,j,k-k_0)\\\vec{m_4}\neq(0,0,0)}}\Dfour^{\vec{m_1}}(\upeta+\ueta)^{k_0}\cdot\Dfour^{\vec{m_2}}\uomega^{j_0}\cdot\Dfour^{\vec{m_3}}(\upchi,\omega)\cdot\Dfour^{\vec{m_4}}F_{\nu_1\cdots\nu_b}^{\mu_1\cdots\mu_a}\\
			\notag&+\sum\limits_{\vec{m_1}+\vec{m_2}+\vec{m_3}+\vec{m_4}=(i,j-j_0,k-k_0-1)}\Dfour^{\vec{m_1}}(\upeta+\ueta)^{k_0}\cdot\Dfour^{\vec{m_2}}\uomega^{j_0}\cdot\Dfour^{\vec{m_3}}\Riem^{\mu_i}_{\mu_04A}\cdot\Dfour^{\vec{m_4}}F_{\nu_1\cdots\nu_b}^{\hat{\mu}_i\mu_0}\\
			\notag&+\sum\limits_{\vec{m_1}+\vec{m_2}+\vec{m_3}+\vec{m_4}=(i,j-j_0,k-k_0-1)}\Dfour^{\vec{m_1}}(\upeta+\ueta)^{k_0}\cdot\Dfour^{\vec{m_2}}\uomega^{j_0}\cdot\Dfour^{\vec{m_3}}\Riem^{\nu_0}_{\nu_i4A}\cdot\Dfour^{\vec{m_4}}F_{\hat{\nu_i}\nu_0}^{\mu_1\cdots\mu_a}\\
			\notag&+\sum\limits_{\vec{m_1}+\vec{m_2}+\vec{m_3}+\vec{m_4}=(i,j-j_0-1,k-k_0)}\Dfour^{\vec{m_1}}(\upeta+\ueta)^{k_0}\cdot\Dfour^{\vec{m_2}}\uomega^{j_0}\cdot\Dfour^{\vec{m_3}}\Riem^{\mu_i}_{\mu_034}\cdot\Dfour^{\vec{m_4}}F_{\nu_1\cdots\nu_b}^{\hat{\mu}_i\mu_0}\\
			\notag&+\sum\limits_{\vec{m_1}+\vec{m_2}+\vec{m_3}+\vec{m_4}=(i,j-j_0-1,k-k_0)}\Dfour^{\vec{m_1}}(\upeta+\ueta)^{k_0}\cdot\Dfour^{\vec{m_2}}\uomega^{j_0}\cdot\Dfour^{\vec{m_3}}\Riem^{\nu_0}_{\nu_i34}\cdot\Dfour^{\vec{m_4}}F_{\hat{\nu_i}\nu_0}^{\mu_1\cdots\mu_a}.
		\end{align}
		For $B=1,2$,
		\begin{align}\label{EQ:CommutatorD3}
			\Dfour_B \Dfour^i_4\Dfour_3^j\Dfour_A^{k}F_{\nu_1\cdots\nu_b}^{\mu_1\cdots\mu_a}=&~\sum\limits_{\vec{m_1}+\vec{m_2}+\vec{m_3}=(i-i_0,j-j_0,k)}\angD^{\vec{m_1}}\upchi^{i_0}\cdot\angD^{\vec{m_2}}\uchi^{ j_0}\cdot\Dfour^{\vec{m_3}}\Dfour_BF_{\nu_1\cdots\nu_b}^{\mu_1\cdots\mu_a}\\
			\notag&+~\sum\limits_{\substack{\vec{m_1}+\vec{m_2}+\vec{m_3}+\vec{m_4}=(i-i_0,j,k-k_0)\\\vec{m_4}\neq(0,0,0)}}\angD^{\vec{m_1}}\upchi^{i_0}\cdot\angD^{\vec{m_2}}\uchi^{ j_0}\cdot\Dfour^{\vec{m_3}}(\upeta_B+\ueta_B)\cdot\Dfour^{\vec{m_4}}F_{\nu_1\cdots\nu_b}^{\mu_1\cdots\mu_a}\\
			\notag&+\sum\limits_{\substack{\vec{m_1}+\vec{m_2}+\vec{m_3}+\vec{m_4}=(i-i_0,j-j_0-1,k)\\}}\angD^{\vec{m_1}}\upchi^{i_0}\cdot\angD^{\vec{m_2}}\uchi^{ j_0}\cdot\Dfour^{\vec{m_3}}\Riem^{\mu_i}_{\mu_03B}\cdot\Dfour^{\vec{m_4}}F_{\nu_1\cdots\nu_b}^{\hat{\mu}_i\mu_0}\\
			\notag&+\sum\limits_{\substack{\vec{m_1}+\vec{m_2}+\vec{m_3}+\vec{m_4}=(i-i_0,j-j_0-1,k)\\}}\angD^{\vec{m_1}}\upchi^{i_0}\cdot\angD^{\vec{m_2}}\uchi^{ j_0}\cdot\Dfour^{\vec{m_3}}\Riem^{\nu_0}_{\nu_i3B}\cdot\Dfour^{\vec{m_4}}F_{\hat{\nu_i}\nu_0}^{\mu_1\cdots\mu_a}\\
			\notag&+\sum\limits_{\substack{\vec{m_1}+\vec{m_2}+\vec{m_3}+\vec{m_4}=(i-i_0-1,j-j_0,k)\\}}\angD^{\vec{m_1}}\upchi^{i_0}\cdot\angD^{\vec{m_2}}\uchi^{ j_0}\cdot\Dfour^{\vec{m_3}}\Riem^{\mu_i}_{\mu_04B}\cdot\Dfour^{\vec{m_4}}F_{\nu_1\cdots\nu_b}^{\hat{\mu}_i\mu_0}\\
			\notag&+\sum\limits_{\substack{\vec{m_1}+\vec{m_2}+\vec{m_3}+\vec{m_4}=(i-i_0-1,j-j_0,k)\\}}\angD^{\vec{m_1}}\upchi^{i_0}\cdot\angD^{\vec{m_2}}\uchi^{ j_0}\cdot\Dfour^{\vec{m_3}}\Riem^{\nu_0}_{\nu_i4B}\cdot\Dfour^{\vec{m_4}}F_{\hat{\nu_i}\nu_0}^{\mu_1\cdots\mu_a}.
		\end{align}
	\end{proposition}
	\begin{proposition}[Commuting transport equations with $\angD$].\label{PR:Commuteijk}
		Let $F$ be a $\suubar$-tangent tensorfield. Then we have the following transport equation in schematic form:
		\begin{align}\label{EQ:angi}
			\angD_4\angnabla^iF=&~\sum\limits_{i_1+i_2+i_3=i}\angnabla^{i_1}(\upeta+\ueta)^{i_2}\angnabla^{i_3}\angD_4F\\
			\notag&+\sum\limits_{\substack{i_1+i_2+i_3+i_4=i\\i_3\neq i}}\angnabla^{i_1}(\upeta+\ueta)^{i_2}\angnabla^{i_3}\upchi\angnabla^{i_4}F\\
			\notag&+\sum\limits_{i_1+i_2+i_3+i_4=i-1}\angnabla^{i_1}(\upeta+\ueta)^{i_2}\angnabla^{i_3}\BetaA\angnabla^{i_4}F,\\
			\notag&+\sum\limits_{i_1+i_2+i_3+i_4=i-1}\angnabla^{i_1}(\upeta+\ueta)^{i_2}\angnabla^{i_3}\FBetaA\angnabla^{i_4}F,\\
			\angD_3\angnabla^iF=&~\sum\limits_{i_1+i_2+i_3=i}\angnabla^{i_1}(\upeta+\ueta)^{i_2}\angnabla^{i_3}\angD_3F\\
			\notag&+\sum\limits_{\substack{i_1+i_2+i_3+i_4=i\\i_3\neq i}}\angnabla^{i_1}(\upeta+\ueta)^{i_2}\angnabla^{i_3}\uchi\angnabla^{i_4}F\\
			\notag&+\sum\limits_{i_1+i_2+i_3+i_4=i-1}\angnabla^{i_1}(\upeta+\ueta)^{i_2}\angnabla^{i_3}\Betabar\angnabla^{i_4}F\\
			\notag&+\sum\limits_{i_1+i_2+i_3+i_4=i-1}\angnabla^{i_1}(\upeta+\ueta)^{i_2}\angnabla^{i_3}\FBetabar\angnabla^{i_4}F,
		\end{align}
	\end{proposition}
	where \begin{align}
		\Dfour^{i_1}\Psi^{i_2}:=\sum\limits_{j_1+\cdots+j_{i_2}=i_1}\Dfour^{j_1}\Psi\cdot\cdot\Dfour^{j_{i_2}}\Psi.
	\end{align}
Moreover, let $\vec{m_a}:=(i_a,j_a,k_a)$'s be multi-indices consisting of $3-$tuples. Then we write the commuted geometric equations as follows:

\begin{align}\label{EQ:ang3ijk}
	&\angD_3\angD_4^{i}\angD_3^{j}\angnabla^{k}F\\
	\notag=&~\sum\limits_{\vec{m_1}+\vec{m_2}+\vec{m_3}=(i-i_0,j,k-k_0)}\angD^{\vec{m_1}}(\upeta+\ueta)^{k_0}\cdot\angD^{\vec{m_2}}\omega^{i_0}\cdot\angD^{\vec{m_3}}\angD_3F\\
	\notag&+\sum\limits_{\substack{\vec{m_1}+\vec{m_2}+\vec{m_3}+\vec{m_4}=(i-i_0,j,k-k_0)\\\vec{m_1},\vec{m_2},\vec{m_3}\neq(i-i_0,j,k-k_0)}}\angD^{\vec{m_1}}(\upeta+\ueta)^{k_0}\cdot\angD^{\vec{m_2}}\omega^{i_0}\cdot\angD^{\vec{m_3}}\uchi\cdot\angD^{\vec{m_4}}F\\
	\notag&+\sum\limits_{\vec{m_1}+\vec{m_2}+\vec{m_3}+\vec{m_4}=(i-i_0,j,k-k_0-1)}\angD^{\vec{m_1}}(\upeta+\ueta)^{k_0}\cdot\angD^{\vec{m_2}}\omega^{i_0}\cdot\angD^{\vec{m_3}}\Betabar\cdot\angD^{\vec{m_4}}F\\
	\notag&+\sum\limits_{\substack{\vec{m_1}+\vec{m_2}+\vec{m_3}+\vec{m_4}=(i-i_0,j,k-k_0)\\\vec{m_4}\neq(0,0,0)}}\angD^{\vec{m_1}}(\upeta+\ueta)^{k_0}\cdot\angD^{\vec{m_2}}\omega^{i_0}\cdot\angD^{\vec{m_3}}\uomega\cdot\angD^{\vec{m_4}}F\\
	\notag&+\sum\limits_{\substack{\vec{m_1}+\vec{m_2}+\vec{m_3}+\vec{m_4}=(i-i_0-1,j,k-k_0+1)\\\vec{m_4}\neq(0,0,0)}}\angD^{\vec{m_1}}(\upeta+\ueta)^{k_0}\cdot\angD^{\vec{m_2}}\omega^{i_0}\cdot\angD^{\vec{m_3}}\upzeta\cdot\angD^{\vec{m_4}}F\\
	\notag&+\sum\limits_{\vec{m_1}+\vec{m_2}+\vec{m_3}+\vec{m_4}=(i-i_0-1,j,k-k_0)}\angD^{\vec{m_1}}(\upeta+\ueta)^{k_0}\cdot\angD^{\vec{m_2}}\omega^{i_0}\cdot\angD^{\vec{m_3}}\Rhodual\cdot\angD^{\vec{m_4}}F\\
	\notag&+\sum\limits_{\vec{m_1}+\vec{m_2}+\vec{m_3}+\vec{m_4}+\vec{m_5}=(i-i_0-1,j,k-k_0)}\angD^{\vec{m_1}}(\upeta+\ueta)^{k_0}\cdot\angD^{\vec{m_2}}\omega^{i_0}\cdot\angD^{\vec{m_3}}\upeta\cdot\angD^{\vec{m_4}}\ueta\cdot\angD^{\vec{m_5}}F\\
	\notag&+\sum\limits_{\vec{m_1}+\vec{m_2}+\vec{m_3}+\vec{m_4}=(i-i_0,j,k-k_0-1)}\angD^{\vec{m_1}}(\upeta+\ueta)^{k_0}\cdot\angD^{\vec{m_2}}\omega^{i_0}\cdot\angD^{\vec{m_3}}\FBetabar\cdot\angD^{\vec{m_4}}F.	
\end{align}	
In particular, a slightly more precise expression of the second term on the RHS of \eqref{EQ:ang3ijk} is as follows:
\begin{align}\label{EQ:ang3ijk2}
	&\sum\limits_{\vec{m_1}+\vec{m_2}+\vec{m_3}+\vec{m_4}=(i-i_0,j,k-k_0)}\angD^{\vec{m_1}}(\upeta+\ueta)^{k_0}\cdot\angD^{\vec{m_2}}\omega^{i_0}\cdot\angD^{\vec{m_3}}\uchi\cdot\angD^{\vec{m_4}}F\\
	\notag=&~k\uchi\cdot\angD^{\vec{m}_5}F+\sum\limits_{\substack{\vec{m_6}+\vec{m_7}+\vec{m_8}+\vec{m_9}=(i-i_0,j,k-k_0)\\\abs{\vec{m}_9}<\abs{\vec{m}_5}}}\angD^{\vec{m_6}}(\upeta+\ueta)^{k_0}\cdot\angD^{\vec{m_7}}\omega^{i_0}\cdot\angD^{\vec{m_8}}\uchi\cdot\angD^{\vec{m_9}}F,
\end{align}
where $\vec{m}_5=(i,j,k)$.
Similarly, we have:
\begin{align}\label{EQ:ang4ijk}
	&\angD_4\angD_4^{i}\angD_3^{j}\angnabla^{k}F\\
	\notag=&~\sum\limits_{\vec{m_1}+\vec{m_2}+\vec{m_3}=(i,j- j_0,k-k_0)}\angD^{\vec{m_1}}(\upeta+\ueta)^{k_0}\cdot\angD^{\vec{m_2}}\uomega^{ j_0}\cdot\angD^{\vec{m_3}}\angD_4F\\
	\notag&+\sum\limits_{\substack{\vec{m_1}+\vec{m_2}+\vec{m_3}+\vec{m_4}=(i-i_0,j,k-k_0)\\\vec{m_1},\vec{m_2},\vec{m_3}\neq(i-i_0,j,k-k_0)}}\angD^{\vec{m_1}}(\upeta+\ueta)^{k_0}\cdot\angD^{\vec{m_2}}\uomega^{ j_0}\cdot\angD^{\vec{m_3}}\upchi\cdot\angD^{\vec{m_4}}F\\
	\notag&+\sum\limits_{\vec{m_1}+\vec{m_2}+\vec{m_3}+\vec{m_4}=(i,j- j_0,k-k_0-1)}\angD^{\vec{m_1}}(\upeta+\ueta)^{k_0}\cdot\angD^{\vec{m_2}}\uomega^{ j_0}\cdot\angD^{\vec{m_3}}\BetaA\cdot\angD^{\vec{m_4}}F\\
	\notag&+\sum\limits_{\substack{\vec{m_1}+\vec{m_2}+\vec{m_3}+\vec{m_4}=(i-i_0,j,k-k_0)\\\vec{m_4}\neq(0,0,0)}}\angD^{\vec{m_1}}(\upeta+\ueta)^{k_0}\cdot\angD^{\vec{m_2}}\uomega^{ j_0}\cdot\angD^{\vec{m_3}}\omega\cdot\angD^{\vec{m_4}}F\\
	\notag&+\sum\limits_{\substack{\vec{m_1}+\vec{m_2}+\vec{m_3}+\vec{m_4}=(i,j- j_0-1,k-k_0+1)\\\vec{m_4}\neq(0,0,0)}}\angD^{\vec{m_1}}(\upeta+\ueta)^{k_0}\cdot\angD^{\vec{m_2}}\uomega^{ j_0}\cdot\angD^{\vec{m_3}}\upzeta\cdot\angD^{\vec{m_4}}F\\
	\notag&+\sum\limits_{\vec{m_1}+\vec{m_2}+\vec{m_3}+\vec{m_4}=(i,j- j_0-1,k-k_0)}\angD^{\vec{m_1}}(\upeta+\ueta)^{k_0}\cdot\angD^{\vec{m_2}}\uomega^{ j_0}\cdot\angD^{\vec{m_3}}\Rhodual\cdot\angD^{\vec{m_4}}F\\
	\notag&+\sum\limits_{\vec{m_1}+\vec{m_2}+\vec{m_3}+\vec{m_4}+\vec{m_5}=(i,j- j_0-1,k-k_0)}\angD^{\vec{m_1}}(\upeta+\ueta)^{k_0}\cdot\angD^{\vec{m_2}}\uomega^{ j_0}\cdot\angD^{\vec{m_3}}\upeta\cdot\angD^{\vec{m_4}}\ueta\cdot\angD^{\vec{m_5}}F\\
	\notag&+\sum\limits_{\vec{m_1}+\vec{m_2}+\vec{m_3}+\vec{m_4}=(i,j- j_0,k-k_0-1)}\angD^{\vec{m_1}}(\upeta+\ueta)^{k_0}\cdot\angD^{\vec{m_2}}\uomega^{ j_0}\cdot\angD^{\vec{m_3}}\FBetaA\cdot\angD^{\vec{m_4}}F,
\end{align}

\begin{align}\label{EQ:angAijk}
	&\angnabla\angD_4^{i}\angD_3^{j}\angnabla^{k}F\\
	\notag=&~\sum\limits_{\vec{m_1}+\vec{m_2}+\vec{m_3}=(i-i_0,j- j_0,k)}\angD^{\vec{m_1}}\upchi^{i_0}\cdot\angD^{\vec{m_2}}\uchi^{ j_0}\cdot\angD^{\vec{m_3}}\angnabla F\\
	\notag&+\sum\limits_{\substack{\vec{m_1}+\vec{m_2}+\vec{m_3}+\vec{m_4}=(i-i_0,j,k-k_0)\\\vec{m_1},\vec{m_2},\vec{m_3}\neq(i-i_0,j,k-k_0)}}\angD^{\vec{m_1}}\upchi^{i_0}\cdot\angD^{\vec{m_2}}\uchi^{ j_0}\cdot\angD^{\vec{m_3}}(\upeta+\ueta)\cdot\angD^{\vec{m_4}}F\\
	\notag&+\sum\limits_{\vec{m_1}+\vec{m_2}+\vec{m_3}+\vec{m_4}=(i-i_0-1,j- j_0,k)}\angD^{\vec{m_1}}\upchi^{i_0}\cdot\angD^{\vec{m_2}}\uchi^{ j_0}\cdot\angD^{\vec{m_3}}\BetaA\cdot\angD^{\vec{m_4}}F\\
	\notag&+\sum\limits_{\vec{m_1}+\vec{m_2}+\vec{m_3}+\vec{m_4}=(i-i_0,j- j_0-1,k)}\angD^{\vec{m_1}}\upchi^{i_0}\cdot\angD^{\vec{m_2}}\uchi^{ j_0}\cdot\angD^{\vec{m_3}}\Betabar\cdot\angD^{\vec{m_4}}F\\
	\notag&+\sum\limits_{\vec{m_1}+\vec{m_2}+\vec{m_3}+\vec{m_4}=(i-i_0-1,j- j_0,k)}\angD^{\vec{m_1}}\upchi^{i_0}\cdot\angD^{\vec{m_2}}\uchi^{ j_0}\cdot\angD^{\vec{m_3}}\FBetaA\cdot\angD^{\vec{m_4}}F\\
	\notag&+\sum\limits_{\vec{m_1}+\vec{m_2}+\vec{m_3}+\vec{m_4}=(i-i_0,j- j_0-1,k)}\angD^{\vec{m_1}}\upchi^{i_0}\cdot\angD^{\vec{m_2}}\uchi^{ j_0}\cdot\angD^{\vec{m_3}}\FBetabar\cdot\angD^{\vec{m_4}}F.
\end{align}

\begin{proof}[Discussion of the proof of Proposition \ref{PR:Commuteijk}]
	Using \eqref{EQ:Commutator2} repeatedly, \eqref{EQ:angi}-\eqref{EQ:angAijk} are obtained by induction arguments. See \cite[Proposition 12]{LukRodnianski2015} for a detailed derivation of \eqref{EQ:angi}; the proofs of \eqref{EQ:ang3ijk}-\eqref{EQ:angAijk} are in a similar fashion.
\end{proof}
\section{The Preliminary Estimates}\label{S:basic}
In our analysis, we aim to integrate along null geodesics. However, by Section \ref{SS:NullFrame}, Lie derivatives $\Lie_{\Lunit}$ and $\Lie_{\uLunit}$ differ slightly from $\Dfour_{4},\Dfour_{3}$. In the next lemma, we provide the distinctions between $\Lie_{\Lunit}$, $\Lie_{\uLunit}$ and $\Dfour_4,\Dfour_{ 3}$.
\begin{lemma}\cite[(3.1.53), Difference between $\Lie_{\Lunit}$ and $\Lie_{\uLunit}$ and $\Dfour_{4},\Dfour_{3}$]{KlainermanNicolo2003}.\label{LE:LieandCovariant}
	Let $\Lie_{V}$ be the standard Lie derivative with respect to vectorfield $V$. For any scalar function $f$, $\suubar$-tensorfield $\xi_0\in\{\upeta,\ueta,\upzeta\},\xi_2\in\{\hchi,\huchi\}$, the following identities hold:
	\begin{subequations}
		\begin{align}
			\Lie_{\Lunit}f=&~\lapse\Dfour_{ 4}f,&
			\Lie_{\uLunit}f=&~\lapse\Dfour_{ 3}f,\\
			\Lie_{\Lunit}\xi_0=&~\lapse\left(\angD_{ 4}\xi_0+\upchi\cdot\xi_0\right),&
			\Lie_{\uLunit}\xi_0=&~\lapse\left(\angD_{ 3}\xi_0+\uchi\cdot\xi_0\right),\\
			\Lie_{\Lunit}\xi_2=&~\lapse\left(\angD_{ 4}\xi_2+2\upchi\cdot\xi_2\right),&
			\Lie_{\uLunit}\xi_2=&~\lapse\left(\angD_{ 3}\xi_2+2\uchi\cdot\xi_2\right).
		\end{align}
	\end{subequations}
\end{lemma}
\subsection{Estimates for Metric Components}\label{SS:metric}
In below, we provide the transport lemma for integrals on the sphere $\suubar$.
\begin{lemma}\cite[Lemma 3.1.3]{KlainermanNicolo2003}\label{LE:e4SuubarInt}.
	For $\suubar$ and a scalar function $f$, we have the following equations:
	\begin{subequations}
		\begin{align}
			\label{EQ:LieL}\Lie_{\Lunit}\left(\int_{\suubar}f\diff\spherevol\right)=&~\int_{\suubar}\left(\Lie_{\Lunit}f+\lapse\gtr\upchi f\right)\diff\spherevol,\\
			\label{EQ:LieuL}\Lie_{\uLunit}\left(\int_{\suubar}f\diff\spherevol\right)=&~\int_{\suubar}\left(\Lie_{\uLunit}f+\lapse\gtr\uchi f\right)\diff\spherevol,
		\end{align}
	\end{subequations}
	where $\Lie_{V}$ is the standard Lie derivative with respect to vectorfield $V$.
\end{lemma}
Now, we control $\Omega$ under the bootstrap assumptions:
\begin{proposition}[Control of lapse function $\lapse$]\label{PR:lapse}
	Under the initial data assumptions in Section \ref{SS:InitialData} and bootstrap assumptions \eqref{BA}, we have:
	\begin{align}
		\norm{\lapse^{-1}-1}_{L^\infty(\suubar)}\lesssim\frac{\delta a^{\frac{1}{2}}}{\abs{u}}\BAsize.
	\end{align}
\end{proposition}
\begin{proof}[Proof of Proposition \ref{PR:lapse}]
	By bootstrap assumptions \eqref{BA}, we integrate the equation \eqref{EQ:omegauomega} along $\hypersur_{u}$ and use the fact that $\lapse=1$ on $\hypersurbar_{0}$ to obtain:
	\begin{align}
		\norm{\lapse^{-1}-1}_{L^\infty(\suubar)}\lesssim\int_0^{\ubar}\norm{\omega}_{L^\infty(\suubarprime)}\diff\ubarprime\lesssim\frac{\delta a^{\frac{1}{2}}}{\abs{u}}\BAsize.
	\end{align}
\end{proof}
\begin{proposition}[Control of spherical induced metric $\gsphere$]\label{PR:gsphere}
	Under the initial data assumptions in Section \ref{SS:InitialData} and bootstrap assumptions \eqref{BA}, for metric $\gsphere$ on $\suubar$, we have:
	\begin{subequations}
		\begin{align}
			\label{ES:detgsphere}c^\prime\leq\det\gsphere\leq C^\prime,\\
			\label{ES:gsphere}\abs{\gsphere_{AB}}+\abs{(\gsphere^{-1})^{AB}}\leq C^\prime.
		\end{align}
	\end{subequations}
	Moreover, fixing $u,\sphere^1,\sphere^2$, let $\lambda_1(\ubar),\lambda_{2}(\ubar)$ be the larger and smaller eigenvalue of $\gsphere^{-1}(u,0,\sphere^1,\sphere^2)\gsphere(u,\ubar,\sphere^1,\sphere^2)$. Then:
	\begin{align}
		\label{ES:eigenvalue}\abs{\lambda_{1}(\ubar)-1}+\abs{\lambda_2(\ubar)-1}\lesssim\frac{\delta a^{\frac{1}{2}}}{\abs{u}}\BAsize.
	\end{align}
\end{proposition}
\begin{proof}[Proof of Proposition \ref{PR:gsphere}]
	We consider the first variation formula:
	\begin{align}\label{EQ:angLieLunit}
		\angLie_{\Lunit}\gsphere=2\lapse\upchi.
	\end{align}
	\eqref{EQ:angLieLunit} implies
	\begin{align}
		\frac{\p}{\p\ubar}\ln(\det\gsphere)=2\lapse\gtr\upchi.
	\end{align}
	By bootstrap assumptions \eqref{BA} and Proposition \eqref{PR:lapse}, we obtain:
	\begin{align}\label{ES:volume1}
		\abs{\frac{\det\gsphere(u,\ubar,\sphere^1,\sphere^2)}{\det\gsphere(u,0,\sphere^1,\sphere^2)}-1}\lesssim\frac{\delta a^{\frac{1}{2}}}{\abs{u}}\BAsize.
	\end{align}
	By initial data assumption in Section \ref{SS:InitialData}, we derive
	\eqref{ES:detgsphere}.
	\eqref{ES:gsphere} follows from initial data assumption in Section \ref{SS:InitialData}, \eqref{ES:detgsphere}, \eqref{ES:eigenvalue}, the fact that $\lambda_1\leq\sup\gsphere_{AB}$ and:
	\begin{align}
		\abs{\gsphere_{AB}(u,\ubar,\sphere^1,\sphere^2)-\gsphere_{AB}(u,0,\sphere^1,\sphere^2)}\leq\int_{0}^{\ubar}\abs{\upchi_{AB}}\diff\ubarprime\leq \lambda_1\delta\norm{\upchi}_{L^\infty(\suubar)}\lesssim \frac{\delta^2 a}{\abs{u}^2}\BAsize.
	\end{align}
	Now we prove \eqref{ES:eigenvalue}. We define $\volume(\ubar):=\sqrt{\frac{\lambda_1(\ubar)}{\lambda_2(\ubar)}}$. Following the derivation of (5.93) in \cite{Christodoulou2009}, by \eqref{EQ:angLieLunit}, we have:
	\begin{align}
		\volume(\ubar)\leq 1+\int_{0}^{\ubar}\sgabs{\lapse\hchi}\volume\diff\ubarprime.
	\end{align}
	Via Gr\"onwall's inequality, we derive:
	\begin{align}\label{ES:volume2}
		\abs{\volume-1}\lesssim\frac{\delta a^{\frac{1}{2}}}{\abs{u}}\BAsize.
	\end{align}
	\eqref{ES:eigenvalue} follows from \eqref{ES:volume1},\eqref{ES:volume2} and the fact that  $\volume=\lambda_1\cdot\frac{\sqrt{\det\gsphere(u,0,\sphere^1,\sphere^2)}}{\sqrt{\det\gsphere(u,\ubar,\sphere^1,\sphere^2)}}$
\end{proof}
\begin{corollary}[Control of area of sphere $\suubar$]\label{CO:Area}
	Under the initial data assumptions in Section \ref{SS:InitialData} and bootstrap assumptions \eqref{BA}, for the area on $\suubar$, we have:
	\begin{align}
		\sup\limits_{u,\ubar}\abs{\mbox{\upshape Area}(\suubar)-\mbox{\upshape Area}(S_{u,0})}\lesssim\frac{\delta a^{\frac{1}{2}}}{\abs{u}}\BAsize.
	\end{align}
\end{corollary}
\subsection{Transport Lemma, Sobolev Embeddings and Elliptic Estimates}\label{SS:TransportSobolevElliptic}
\begin{proposition}\cite[Proposition 3.5]{AnAthanasiou2022}\label{PR:transport1}.
	Let $F$ be a $\suubar$-tangent tensorfield. Under the bootstrap assumptions \eqref{BA}, we have:
	\begin{subequations}
		\begin{align}
			\norm{F}_{L^2(\suubar)}&\lesssim\norm{F}_{L^2(\suubarprime)}+\int_{\ubarprime}^{\ubar}\norm{\angD_4F}_{L^2(S_{u,\ubar^{\prime\prime}})}\diff\ubar^{\prime\prime},\\
			\norm{F}_{L^2(\suubar)}&\lesssim\norm{F}_{L^2(\suprimeubar)}+\int_{\uprime}^{u}\norm{\angD_3F}_{L^2(S_{u^{\prime\prime},\ubar})}\diff u^{\prime\prime}.
		\end{align}
	\end{subequations}
\end{proposition}
\begin{proposition}\cite[Proposition 3.7]{AnAthanasiou2022}\label{PR:transport2}.
	Let $F,\mathfrak{F}$ be $\suubar$-tangent tensorfields satisfying the following transport equation:
	\begin{align}
		\angD_{3}F+\lambda_0\gtr\uchi F=&~\mathfrak{F}.
	\end{align} 
	Denoting $\tilde{\lambda_0}=2\lambda_0-1$, under the bootstrap assumptions \eqref{BA}, for metric $\gsphere$ on $\suubar$, we have:
	\begin{align}
		\abs{u}^{\tilde{\lambda_0}}\norm{F}_{L^2(\suubar)}\lesssim\abs{\uinf}^{\tilde{\lambda_0}}\norm{F}_{L^2(S_{\uinf,\ubar})}+\int_{\uinf}^{u}\abs{\uprime}^{\tilde{\lambda_0}}\norm{\mathfrak{F}}_{L^2(\suprimeubar)}\diff\uprime
	\end{align}
\end{proposition}
\begin{proposition}\cite[Proposition 5.6. Control of isoperimetric constant]{AnLuk2017}.\label{PR:Isoperimetric}
	Let isoperimetric constant $I(\suubar)$ be as follows:
	\begin{align}
		I(\suubar):=\sup\limits_{\substack{U\subset\suubar\\\p U\in C^1}}\frac{\min\{\mbox{\upshape Area}(U),\mbox{\upshape Area}(U^c)\}}{(\mbox{\upshape Perimeter}(\p U))^2}.
	\end{align}
	Then, Under the initial data assumptions in Section \ref{SS:InitialData} and bootstrap assumptions \eqref{BA}, we have:
	\begin{align}
		I(\suubar)\leq\frac{1}{\pi}.
	\end{align}
\end{proposition}
\begin{proposition}\cite[Lemma 5.1-Lemma 5.2]{Christodoulou2009}\label{PR:Sobolevinequality}.
	For any compact connected Riemannian 2-manifold $(S,\gamma)$, denoting $F$ the $S$-tangent tensorfield, the following Sobolev inequalities hold for $2<p<\infty$:
	\begin{subequations}
		\begin{align}
			(\mbox{\upshape Area}(S))^{-\frac{1}{p}}\norm{F}_{L^p(S)}&\leq C_p\sqrt{\max\{I(S),1\}}\cdot\left(\norm{\angnabla F}_{L^2(S)}+(\mbox{\upshape Area}(S))^{-\frac{1}{2}}\norm{F}_{L^2(S)}\right),\\
			\norm{F}_{L^\infty(S)}&\leq C_p\sqrt{\max\{I(S),1\}}\cdot(\mbox{\upshape Area}(S))^{\frac{1}{2}-\frac{1}{p}}\left(\norm{\angnabla F}_{L^p(S)}+(\mbox{\upshape Area}(S))^{-\frac{1}{2}}\norm{F}_{L^p(S)}\right).
		\end{align}
	\end{subequations}
\end{proposition}
\begin{proposition}[$L^2-L^\infty$ Sobolev embedding]\label{PR:SobolevEmbedding}
	Let $F$ be a $\suubar$-tangent tensorfield. Under the initial data assumptions in Section \ref{SS:InitialData} and bootstrap assumptions \eqref{BA}, combining Corollary \ref{CO:Area}, Proposition \ref{PR:Isoperimetric}-\ref{PR:Sobolevinequality}, we have the following Sobolev inequality:
	\begin{align}
		\norm{F}_{L^\infty(\suubar)}\lesssim\sum\limits_{i=0,1,2}\norm{\abs{u}^{i-1}\angnabla^iF}_{L^2(\suubar)}.
	\end{align}
\end{proposition}
\begin{lemma}\cite[Chapter 2, $L^2$-estimates for div-curl system on $\suubar$]{ChristodoulouKlainerman1993}.\label{LE:ellipticES}
	On an arbitrary compact Riemannian manifold $(\suubar,\gsphere)$, with $\GaussCurvature$ the Gauss curvature of $\suubar$, we have, for vectorfield $\xi$:
	\begin{align}\label{E:DivCurlonsuubar}
		\int_{\suubar}\sgabs{\angnabla\xi}^2+\GaussCurvature\sgabs{\xi}^2\diff\spherevol=\int_{\suubar}\sgabs{\angdiv\xi}^2+\sgabs{\angcurl\xi}^2\diff\spherevol.
	\end{align}
	Moreover, if $\xi$ is a symmetric traceless 2-tensor, then
	\begin{align}\label{E:Divonsuubar}
		\int_{\suubar}\sgabs{\angnabla\xi}^2+2\GaussCurvature\sgabs{\xi}^2\diff\spherevol=2\int_{\suubar}\sgabs{\angdiv\xi}^2\diff\spherevol.
	\end{align}
\end{lemma}
\begin{lemma}[Elliptic estimate for symmetric type $\binom{0}{2}$ tensorfield on $\suubar$]\label{LE:EllipticEstimateonSuubar}
	Let $\xi$ be a symmetric type $\binom{0}{2}$ tensorfield on $\suubar$. Then we have the following elliptic estimate:
	\begin{align}\label{ES:EllipticEstimateonSuubar}
		\onenorm{\angnabla\xi}{\gsphere}{2}{\suubar}\leq2\onenorm{\angdiv\xi}{\gsphere}{2}{\suubar}+\onenorm{\angnabla\gtr\xi}{\gsphere}{2}{\suubar}+2\onenorm{\GaussCurvature}{\gsphere}{\infty}{\suubar}\onenorm{\xi}{\gsphere}{2}{\suubar},
	\end{align}
	where $\hat{\xi}$ is the $g$-traceless part of $\xi$, $\GaussCurvature$ is the Gauss curvature of $\suubar$.
\end{lemma}
\begin{proof}[Proof of Lemma \ref{LE:EllipticEstimateonSuubar}]
	We decompose $\xi$ into trace part and traceless part as follows:
	\begin{align}\label{EQ:ProofEllipticEstimates1}
		\xi=&~\hat{\xi}+\frac{1}{2}\gtr\xi\gsphere.
	\end{align}
	Taking $\angnabla$ and $\angdiv$ on \eqref{EQ:ProofEllipticEstimates1}, we derive:
	\begin{subequations}\label{EQ:ProofEllipticEstimates2}
		\begin{align}
			\angnabla\xi=&~\angnabla\hat{\xi}+\frac{1}{2}\angnabla\gtr\xi\gsphere,\\
			\angdiv\xi=&~\angdiv\hat{\xi}+\frac{1}{2}\angnabla\gtr\xi\cdot\gsphere.
		\end{align}
	\end{subequations}
	By \eqref{E:Divonsuubar}, we have:
	\begin{align}\label{ES:EllipticEstimateonSuubar1.1}
		\onenorm{\angnabla\hat{\xi}}{\gsphere}{2}{\suubar}\leq2\onenorm{\angdiv\hat{\xi}}{\gsphere}{2}{\suubar}+2\onenorm{\GaussCurvature}{\gsphere}{\infty}{\suubar}\onenorm{\hat{\xi}}{\gsphere}{2}{\suubar}.
	\end{align}
	Desired result \eqref{ES:EllipticEstimateonSuubar} follows from combining \eqref{ES:EllipticEstimateonSuubar1.1} and \eqref{EQ:ProofEllipticEstimates2}.
\end{proof}

\subsection{Comparison Lemmas, Basic estimates for $\psi,\upsi$}
\begin{lemma}[Estimates for the Christoffel symbols compared to velocity components]\label{LE:EstimatesofTensorDifference}
	Let $A,B=1,2$ and $I,J=1,2,3,4$. Let $\Dfour\in\{\Dfour_4,\Dfour_3,\Dfour_A\},\angD\in\{\angD_4,\angD_3,\angD_A\}$ be defined as in Section \ref{SS:despacetime}. Let velocity components $v_I$ be as in \eqref{EQ:devi}. We have the following estimates for the $v$:
	\begin{subequations}\label{ES:EstimatesforChristoffelSymbols}
		\begin{align}
			\norm{\Dfour_I (v_J)-\Dfour_I v_J}_{L^2(\suubar)}\lesssim&\norm{\Dfour_I  v_J}_{L^2(\suubar)},&
			(I,J)\neq&(A,3),(A,B)\\
			\norm{\Dfour ( v_I)-\Dfour  v_I}_{L^\infty(\suubar)}\lesssim&\norm{\Dfour  v_I}_{L^\infty(\suubar)},&
			(I,J)\neq&(A,3),(A,B)\\
			\norm{\angD  v_A-\Dfour  v_A}_{L^2(\suubar)}\lesssim&\norm{\Dfour  v_A}_{L^2(\suubar)},&
			(I,J)\neq&(A,3),(A,B)\\
			\norm{\angD  v_A-\Dfour  v_A}_{L^\infty(\suubar)}\lesssim&\norm{\Dfour  v_A}_{L^\infty(\suubar)}&
			(I,J)\neq&(A,3),(A,B).
		\end{align}
	\end{subequations}
\end{lemma}
\begin{proof}[Proof of Lemma \ref{LE:EstimatesofTensorDifference}]
	
	\eqref{ES:EstimatesforChristoffelSymbols} follows directly from
	Proposition \ref{PR:PropConnection}, bootstrap assumptions \eqref{BA21}-\eqref{BA22}, \eqref{BA31}-\eqref{BA32}. 
\end{proof}
In below, we provide a useful lemma for estimating the products of derivatives of $\upsi$.
\begin{lemma}[Bounds for products of derivatives of $\upsi$]\label{LE:BoundsforProducts}
	Let $\upsi$ be as in Section \ref{SS:SchematicNotation} and $\WangD^{\vec{m}}$ be as in Definition \ref{DE:operators}. Under the bootstrap assumptions \eqref{BA}, we have the estimates for the derivatives of $\upsi$:
	\begin{align}\label{ES:BoundsforProducts1}
		\sum\limits_{\abs{\vec{m_1}}\leq N-1}\norm{\abs{u}^l\WangD^{\vec{m_1}}\upsi^l}_{L^2(\suubar)}+\sum\limits_{\abs{\vec{m_2}}\leq N-3}\norm{\abs{u}^{l+1}\WangD^{\vec{m_2}}\upsi^l}_{L^\infty(\suubar)}\lesssim\delta a^{\frac{1}{2}}\BAsize.
	\end{align}
	Moreover, for mixed products, we have:
	\begin{subequations}\label{ES:BoundsforProducts2}
		\begin{align}\label{ES:BoundsforProducts21}
			&\sum\limits_{\substack{\vec{m}_1+\vec{m}_2=(i-i_0,j,k-k_0)\\i+j+k\leq N-1}}\delta^{i_0}\cdot\twonorms{\abs{u}^{k_0}\WangD^{\vec{m}_1}\psi^{i_0+l}\cdot\WangD^{\vec{m}_2}\upsi^{k_0+n}}{u}{1}{\sphere}{2}{\hypersurbar_{\ubar}}
			\lesssim\left(\BAsize a^{\frac{1}{2}}\right)^{l}\cdot\left(\BAsize\delta a^{\frac{1}{2}}\right)^{i_0+k_0+n}\cdot\abs{u}^{-k_0-i_0-l-2n+2},\\
			\label{ES:BoundsforProducts22}&\sum\limits_{\substack{\vec{m}_1+\vec{m}_2=(i-i_0,j,k-k_0)\\\abs{\vec{m}_1+\vec{m}_2}\leq \frac{N}{2}}}\delta^{i_0}\cdot\twonorms{\abs{u}^{k_0}\WangD^{\vec{m}_1}\psi^{i_0+l}\cdot\WangD^{\vec{m}_2}\upsi^{k_0+n}}{u}{\infty}{\sphere}{\infty}{\hypersurbar_{\ubar}}
			\lesssim\left(\BAsize a^{\frac{1}{2}}\right)^{l}\cdot\left(\BAsize\delta a^{\frac{1}{2}}\right)^{i_0+k_0+n}\cdot\abs{u}^{-k_0-i_0-l-2n}.
		\end{align}
	\end{subequations}
	Similarly, on $\hypersur_{u}$:
	\begin{subequations}\label{ES:BoundsforProducts3}
		\begin{align}\label{ES:BoundsforProducts31}
			&\sum\limits_{\substack{\vec{m}_1+\vec{m}_2=(i-i_0,j,k-k_0)\\i+j+k\leq N-1}}\delta^{i_0}\abs{u}^{k_0}\cdot\twonorms{\WangD^{\vec{m}_1}\psi^{i_0+l}\cdot\WangD^{\vec{m}_2}\upsi^{k_0+n}}{\ubar}{1}{\sphere}{2}{\hypersur_{u}}
			\lesssim\delta \left(\BAsize a^{\frac{1}{2}}\right)^{l}\cdot\left(\BAsize\delta a^{\frac{1}{2}}\right)^{i_0+k_0+n}\cdot\abs{u}^{-k_0-i_0-l-2n+1},\\
			\label{ES:BoundsforProducts32}&\sum\limits_{\substack{\vec{m}_1+\vec{m}_2=(i-i_0,j,k-k_0)\\\abs{\vec{m}_1+\vec{m}_2}\leq \frac{N}{2}}}\delta^{i_0}\abs{u}^{k_0}\cdot\twonorms{\WangD^{\vec{m}_1}\psi^{i_0+l}\cdot\WangD^{\vec{m}_2}\upsi^{k_0+n}}{\ubar}{\infty}{\sphere}{\infty}{\hypersur_{u}}
			\lesssim\left(\BAsize a^{\frac{1}{2}}\right)^{l}\cdot\left(\BAsize\delta a^{\frac{1}{2}}\right)^{i_0+k_0+n}\cdot\abs{u}^{-k_0-i_0-l-2n}.
		\end{align}
	\end{subequations}
\end{lemma}
\begin{proof}[Proof of Lemma \ref{LE:BoundsforProducts}]
	In \eqref{ES:BoundsforProducts1}, we prove the $L^2(\suubar)$-estimates only while the proof of  $L^\infty(\suubar)$-estimates is in a similar fashion. By bootstrap assumptions \eqref{BA}, we have:
	\begin{align}
		\sum\limits_{\abs{\vec{m}_1}\leq N-1}\norm{\abs{u}^l\WangD^{\vec{m}_1}\upsi^l}_{L^2(\suubar)}\lesssim&~\left(\sum\limits_{\abs{\vec{m}_{11}}\leq N-1}\norm{\abs{u}\WangD^{\vec{m}_{11}}\upsi}_{L^2(\suubar)}\right)\left(\sum\limits_{\abs{\vec{m}_{12}}\leq \frac{N}{2}}\norm{\abs{u}\WangD^{\vec{m}_{12}}\upsi}^l_{L^\infty(\suubar)}\right)\\
		\notag\lesssim&~\sum\limits_{\abs{\vec{m}_{11}}\leq N-1}\norm{\abs{u}\WangD^{\vec{m}_{11}}\upsi}_{L^2(\suubar)}\\
		\notag\lesssim&~\delta a^{\frac{1}{2}}\BAsize.
	\end{align}
	Now we prove \eqref{ES:BoundsforProducts21}; the proof of \eqref{ES:BoundsforProducts22} is in a similar fashion. By bootstrap assumptions \eqref{BA},
	\begin{align}
		&\sum\limits_{\vec{m}_1+\vec{m}_2=(i-i_0,j,k-k_0)}\twonorms{\WangD^{\vec{m}_1}\psi^{i_0+1}\cdot\WangD^{\vec{m}_2}\upsi^{k_0+1}}{u}{1}{\sphere}{2}{\hypersurbar_{\ubar}}\\
		\notag\lesssim&~\sum\limits_{\substack{\vec{m}_1+\vec{m}_2+\vec{m}_{11}=(i-i_0,j,k-k_0)\\\abs{\vec{m}_{11}}>\frac{N}{2}}}\twonorms{\WangD^{\vec{m}_{11}}\psi}{u}{1}{\sphere}{2}{\hypersurbar_{\ubar}}\cdot\left(\twonorms{\WangD^{\vec{m}_1}\psi}{u}{\infty}{\sphere}{\infty}{\hypersurbar_{\ubar}}^{i_0}\cdot\twonorms{\WangD^{\vec{m}_2}\upsi}{u}{\infty}{\sphere}{\infty}{\hypersurbar_{\ubar}}^{k_0+1}\right)\\
		\notag&+\sum\limits_{\substack{\vec{m}_1+\vec{m}_2+\vec{m}_{21}=(i-i_0,j,k-k_0)\\\abs{\vec{m}_{21}}>\frac{N}{2}}}\twonorms{\WangD^{\vec{m}_{21}}\upsi}{u}{1}{\sphere}{2}{\hypersurbar_{\ubar}}\cdot\left(\twonorms{\WangD^{\vec{m}_1}\psi}{u}{\infty}{\sphere}{\infty}{\hypersurbar_{\ubar}}^{i_0+1}\cdot\twonorms{\WangD^{\vec{m}_2}\upsi}{u}{\infty}{\sphere}{\infty}{\hypersurbar_{\ubar}}^{k_0}\right)\\
		\notag\lesssim&~\BAsize^{i_0+k_0+2}\left( a^{\frac{1}{2}}\right)^{i_0+1}\left(\frac{\delta a^{\frac{1}{2}}}{\abs{u}^2}\right)^{k_0+1}\abs{u}^2\\
		\notag\lesssim&~\BAsize^{i_0+k_0+2}\delta^{k_0+1}\cdot a^{\frac{i_0+k_0+2}{2}}\cdot\abs{u}^{-2k_0}.
	\end{align}
\end{proof}
\section{Energy Estimates for Fluid Variables}
\subsection{Cone Flux Energy and Energy Current}
\begin{definition}\label{DE:ConeFluxEnergy}
	We define the flux energy $\Energy[\varphi](u),\tilde{\Energy}[\varphi](u)$ on $\hypersur_u$ and $\uEnergy[\varphi](\ubar),\tilde{\uEnergy}[\varphi](\ubar)$ on $\hypersurbar_{\ubar}$ as follows:
	\begin{subequations}
		\begin{align}
			\Energy[\varphi](u):=&~\int_{\hypersur_u}\varphi^2\diff\ubar\diff\spherevol,&
			\uEnergy[\varphi](\ubar):=&~\int_{\hypersurbar_{\ubar}}\varphi^2\diff u\diff\spherevol,\\
			\tilde{\Energy}_{(k)}[\varphi](u):=&~\sum\limits_{j=0}^{k}\Energy[\pfour^j\varphi](u),&
			\tilde{\uEnergy}_{(k)}[\varphi](\ubar):=&~\sum\limits_{j=0}^{k}\uEnergy[\pfour^j\varphi](\ubar).
		\end{align}
	\end{subequations}
\end{definition}

\begin{definition}[Energy Current] \label{DE:EnergyCurrent}
	Let $(\varh,\vars,\varv^a):=(\delta \Dfour_4)^{i}(\abs{u}\Dfour_3)^j(\abs{u}\Dfour_A)^k(\lnenth,\ent,v^a)$, $a=1,2,3$, be the variation of $(\lnenth,\ent,v^a)$ generated by $(\delta \Dfour_4)^{i}(\abs{u}\Dfour_3)^j(\abs{u}\Dfour_A)^k$. 
	For geometric coordinates indices $I,J=1,2,3,4$, we use the notation 
	\begin{align}
		\varv_I:=(\delta \Dfour_4)^{i}(\abs{u}\Dfour_3)^j(\abs{u}\Dfour_A)^k (v_I),
	\end{align}
	where we define $\varv^0$ such that 
	\begin{align}\label{EQ:gvarvv=0}
		(\gfour^{-1})^{IJ}\varv_I v_J=&0.
	\end{align}
	We define the energy current $\varJ$ as follows:
	\begin{align}\label{EQ:EnergyCurrent}
		\varJ^\mu:=&~\frac{v^\mu}{\speed^2}\left(\varh^2+\frac{q^2\temp^2\speed^2}{\enth^2}\vars^2-2\frac{\temp}{\enth}\varh\vars\right)+(-\varv_4e_3^\mu-\varv_3e_4^\mu+2\varv_B e_B^\mu)(\varh-\frac{\temp}{\enth}\vars)+v^\mu\left(-\varv_3\varv_4+\varv_B\varv_B\right),
	\end{align}
	where $B=1,2$, $q$ is a large constant determined by $\speed<1$ in Lemma \ref{LE:Coerciveness}.
\end{definition}
\begin{remark}
	We note that under bootstrap assumptions \eqref{BA}, $\varv_I$ has the same size to $\gfour(\varv,e_I)$ except for $\Dfour_Av_B,\Dfour_Av_3$ as we showed in Lemma \ref{LE:EstimatesofTensorDifference}. Hence in the rest of the paper, we abuse the notation of $\Dfour_Jv_I$ and $\Dfour_J(v_I)$ for $(I,J)\neq (A,3),(A,B)$, while we distinguish $\Dfour_Av_B$ and $\Dfour_A(v_B)$, $\Dfour_Av_3$ and $\Dfour_A(v_3)$ carefully. 
\end{remark}
\subsection{Commuted Euler Equations}\label{SS:CommutedEulerEq}
\begin{proposition}[Equations satisfied by variations]\label{PR:EQVariations}
	Let $\WDfour$ be as in Definition \ref{DE:operators}. For solutions $(\lnenth,\ent,v^\mu)$ of \eqref{EQ:FOE2}, its variation $(\varh,\vars,\varv^\mu)$ with respect to operator $\WDfour^{\vec{m}}:=(\delta \Dfour_4)^{i}(\abs{u}\Dfour_3)^j(\abs{u}\Dfour_A)^k$ defined in Definition \ref{DE:EnergyCurrent} satisfies the following equations of variation:
	\begin{subequations}\label{EQ:FOE3}
		\begin{align}
			\label{EQ:FOE3.1}v^\mu\Dfour_\mu\varh+\frac{\speed^2}{2}(-\Dfour_4 \varv_3-\Dfour_3 \varv_4+2\Dfour_B \varv_B)=&~\mathfrak{F}_{(i,j,k)},\\
			\label{EQ:FOE3.2.1}v^\mu\Dfour_\mu \varv_4+\vproject^{\mu}_4\Dfour_\mu\varh-\frac{\temp}{\enth}\Dfour_4\vars=&~(\mathfrak{G}_{(i,j,k)})_4,\\
			\label{EQ:FOE3.2.2}v^\mu\Dfour_\mu \varv_3+\vproject^{\mu}_3\Dfour_\mu\varh-\frac{\temp}{\enth}\Dfour_3\vars=&~(\mathfrak{G}_{(i,j,k)})_3,\\
			\label{EQ:FOE3.2.3}v^\mu\Dfour_\mu (\varv_B)+\vproject^{\mu}_B\Dfour_\mu\varh-\frac{\temp}{\enth}\Dfour_B\vars=&~(\mathfrak{G}_{(i,j,k)})_B,\\
			\label{EQ:FOE3.3}v^\mu\Dfour_\mu\vars=&~\mathfrak{H}_{(i,j,k)}.
		\end{align}
	\end{subequations}
	where $\mathfrak{F}_{(i,j,k)}:=\Fluid{\varh}{F}+\Fluid{\varv}{F},\mathfrak{G}_{(i,j,k)}:=\Fluid{\varv}{G}+\Fluid{\varh}{G}+\Fluid{\vars}{G}$ ($\mathfrak{F}_{(0,0,0)}=\mathfrak{G}^\alpha_{(0,0,0)}=\mathfrak{H}_{(0,0,0)}=0$) and $\mathfrak{H}$ are inhomogeneous terms that arise from commuting \eqref{EQ:FOE1} with variation generator operator $\WDfour^{\vec{m}}:=(\delta \Dfour_4)^i(\abs{u}\Dfour_3)^j(\abs{u}\Dfour_A)^k$. Schematically, the components of $\mathfrak{F}_{(i,j,k)},\mathfrak{G}_{(i,j,k)}, \mathfrak{H}_{(i,j,k)}$ read as follows:
	\begin{align}\label{EQ:Fvarh}
		\Fluid{\varh}{F}=&\sum\limits_{\substack{\vec{m_1}+\vec{m_2}+\vec{m_3}+\vec{m_4}=(i-i_0,j,k-k_0)\\\abs{\vec{m_3}}\geq1
		}}\delta^{i_0}\abs{u}^{k_0}\WDfour^{\vec{m_1}}(\upeta+\ueta)^{k_0}\cdot\WDfour^{\vec{m_2}}\omega^{i_0}\cdot\WDfour^{\vec{m_3}}v_4\cdot\WDfour^{\vec{m_4}}\Dfour_3\lnenth\\
		\notag&+\sum\limits_{\substack{\vec{m_1}+\vec{m_2}+\vec{m_3}+\vec{m_4}+\vec{m_5}=(i-i_0,j,k-k_0)\\\vec{m_1},\vec{m_2},\vec{m_4}\neq(i-i_0,j,k-k_0)}}\delta^{i_0}\abs{u}^{k_0}\WDfour^{\vec{m_1}}(\upeta+\ueta)^{k_0}\cdot\WDfour^{\vec{m_2}}\omega^{i_0}\cdot\WDfour^{\vec{m_3}}v_4\cdot\WDfour^{\vec{m_4}}(\uchi,\uomega,\upzeta)\cdot\WDfour^{\vec{m_5}}\lnenth\\
		\notag&+\sum\limits_{\substack{\vec{m_1}+\vec{m_2}+\vec{m_3}+\vec{m_4}=(i,j-j_0,k-k_0)\\\abs{\vec{m_3}}\geq 1
		}}\abs{u}^{j_0+k_0}\WDfour^{\vec{m_1}}(\upeta+\ueta)^{k_0}\cdot\WDfour^{\vec{m_2}}\uomega^{j_0}\cdot\WDfour^{\vec{m_3}}(v_3)\cdot\WDfour^{\vec{m_4}}\Dfour_4\lnenth\\
		\notag&+\sum\limits_{\substack{\vec{m_1}+\vec{m_2}+\vec{m_3}+\vec{m_4}+\vec{m_5}=(i,j-j_0,k-k_0)\\\vec{m_1},\vec{m_2},\vec{m_4}\neq(i,j-j_0,k-k_0)}}\abs{u}^{j_0+k_0}\WDfour^{\vec{m_1}}(\upeta+\ueta)^{k_0}\cdot\WDfour^{\vec{m_2}}\uomega^{j_0}\cdot\WDfour^{\vec{m_3}}(v_3)\cdot\WDfour^{\vec{m_4}}(\upchi,\omega)\cdot\WDfour^{\vec{m_5}}\lnenth\\
		\notag&+\sum\limits_{\substack{\vec{m_1}+\vec{m_2}+\vec{m_3}+\vec{m_4}=(i-i_0,j-j_0,k)\\\abs{\vec{m_3}}\geq 1
		}}\delta^{i_0}\abs{u}^{j_0}\WDfour^{\vec{m_1}}\upchi^{i_0}\cdot\WDfour^{\vec{m_2}}\uchi^{j_0}\cdot\WDfour^{\vec{m_3}} (v_B)\cdot\WDfour^{\vec{m_4}}\Dfour_B\lnenth\\
		\notag&+\sum\limits_{\substack{\vec{m_1}+\vec{m_2}+\vec{m_3}+\vec{m_4}+\vec{m_5}=(i-i_0,j-j_0,k)\\\vec{m_1},\vec{m_2},\vec{m_4}\neq(i-i_0,j-j_0,k)}}\delta^{i_0}\abs{u}^{j_0}\WDfour^{\vec{m_1}}\upchi^{i_0}\cdot\WDfour^{\vec{m_2}}\uchi^{j_0}\cdot\WDfour^{\vec{m_3}} (v_B)\cdot\WDfour^{\vec{m_4}}(\upeta+\ueta)\cdot\WDfour^{\vec{m_5}}\lnenth.
	\end{align}

	\begin{align}\label{EQ:varvF}
		\Fluid{\varv}{F}=&\sum\limits_{\substack{\vec{m_1}+\vec{m_2}+\vec{m_3}+\vec{m_4}=(i-i_0,j,k-k_0)\\\abs{\vec{m_3}}\geq 1
		}}\delta^{i_0}\abs{u}^{k_0}\WDfour^{\vec{m_1}}(\upeta+\ueta)^{k_0}\cdot\WDfour^{\vec{m_2}}\omega^{i_0}\cdot\WDfour^{\vec{m_3}}\speed^2\cdot\WDfour^{\vec{m_4}}\Dfour_3v_4\\
		\notag&+\sum\limits_{\substack{\vec{m_1}+\vec{m_2}+\vec{m_3}+\vec{m_4}+\vec{m_5}=(i-i_0,j,k-k_0)\\\vec{m_1},\vec{m_2},\vec{m_4}\neq(i-i_0,j,k-k_0)}}\delta^{i_0}\abs{u}^{k_0}\WDfour^{\vec{m_1}}(\upeta+\ueta)^{k_0}\cdot\WDfour^{\vec{m_2}}\omega^{i_0}\cdot\WDfour^{\vec{m_3}}\speed^2\cdot\WDfour^{\vec{m_4}}(\uchi,\uomega,\upzeta)\cdot\WDfour^{\vec{m_5}}v_4\\
		\notag&+\sum\limits_{\vec{m_1}+\vec{m_2}+\vec{m_3}+\vec{m_4}+\vec{m_5}=(i-i_0,j,k-k_0-1)}\delta^{i_0}\abs{u}^{k_0+1}\WDfour^{\vec{m_1}}(\upeta+\ueta)^{k_0}\cdot\WDfour^{\vec{m_2}}\omega^{i_0}\cdot\WDfour^{\vec{m_3}}\speed^2\cdot\WDfour^{\vec{m_4}}\Riem^{I}_{43A}\cdot\WDfour^{\vec{m_5}}(v_{I})\\
		\notag&+\sum\limits_{\vec{m_1}+\vec{m_2}+\vec{m_3}+\vec{m_4}+\vec{m_5}=(i-i_0-1,j,k-k_0)}\delta^{i_0+1}\abs{u}^{k_0}\WDfour^{\vec{m_1}}(\upeta+\ueta)^{k_0}\cdot\WDfour^{\vec{m_2}}\omega^{i_0}\cdot\WDfour^{\vec{m_3}}\speed^2\cdot\WDfour^{\vec{m_4}}\Riem^{I}_{434}\cdot\WDfour^{\vec{m_5}}(v_I)\\
		\notag&+\sum\limits_{\substack{\vec{m_1}+\vec{m_2}+\vec{m_3}+\vec{m_4}=(i,j-j_0,k-k_0)\\\abs{\vec{m_3}}\geq 1
		}}\abs{u}^{j_0+k_0}\WDfour^{\vec{m_1}}(\upeta+\ueta)^{k_0}\cdot\WDfour^{\vec{m_2}}\uomega^{j_0}\cdot\WDfour^{\vec{m_3}}\speed^2\cdot\WDfour^{\vec{m_4}}\Dfour_4v_3\\
		\notag&+\sum\limits_{\substack{\vec{m_1}+\vec{m_2}+\vec{m_3}+\vec{m_4}+\vec{m_5}=(i,j-j_0,k-k_0)\\\vec{m_1},\vec{m_2},\vec{m_4}\neq(i,j-j_0,k-k_0)}}\abs{u}^{j_0+k_0}\WDfour^{\vec{m_1}}(\upeta+\ueta)^{k_0}\cdot\WDfour^{\vec{m_2}}\uomega^{j_0}\cdot\WDfour^{\vec{m_3}}\speed^2\cdot\WDfour^{\vec{m_4}}(\upchi,\omega)\cdot\WDfour^{\vec{m_5}}(v_3)\\
		\notag&+\sum\limits_{\vec{m_1}+\vec{m_2}+\vec{m_3}+\vec{m_4}+\vec{m_5}=(i,j-j_0,k-k_0-1)}\abs{u}^{j_0+k_0+1}\WDfour^{\vec{m_1}}(\upeta+\ueta)^{k_0}\cdot\WDfour^{\vec{m_2}}\uomega^{j_0}\cdot\WDfour^{\vec{m_3}}\speed^2\cdot\WDfour^{\vec{m_4}}\Riem^{I}_{34A}\cdot\WDfour^{\vec{m_5}}(v_{I})\\
		\notag&+\sum\limits_{\vec{m_1}+\vec{m_2}+\vec{m_3}+\vec{m_4}+\vec{m_5}=(i,j-j_0-1,k-k_0)}\abs{u}^{j_0+k_0+1}\WDfour^{\vec{m_1}}(\upeta+\ueta)^{k_0}\cdot\Dfour^{\vec{m_2}}\uomega^{j_0}\cdot\WDfour^{\vec{m_3}}\speed^2\cdot\WDfour^{\vec{m_4}}\Riem^{I}_{334}\cdot\WDfour^{\vec{m_5}}(v_I)\\
		\notag&+\sum\limits_{\substack{\vec{m_1}+\vec{m_2}+\vec{m_3}+\vec{m_4}=(i-i_0,j-j_0,k)\\\abs{\vec{m_3}}\geq 1
		}}\delta^{i_0}\abs{u}^{j_0}\WDfour^{\vec{m_1}}\upchi^{i_0}\cdot\WDfour^{\vec{m_2}}\uchi^{j_0}\cdot\WDfour^{\vec{m_3}}\speed^2\cdot\WDfour^{\vec{m_4}}\Dfour_Bv_B\\
		\notag&+\sum\limits_{\substack{\vec{m_1}+\vec{m_2}+\vec{m_3}+\vec{m_4}+\vec{m_5}=(i-i_0,j-j_0,k)\\\vec{m_1},\vec{m_2},\vec{m_4}\neq(i-i_0,j-j_0,k)}}\delta^{i_0}\abs{u}^{j_0}\WDfour^{\vec{m_1}}\upchi^{i_0}\cdot\WDfour^{\vec{m_2}}\uchi^{j_0}\cdot\WDfour^{\vec{m_3}}\speed^2\cdot\WDfour^{\vec{m_4}}(\upeta_B+\ueta_B)\cdot\WDfour^{\vec{m_5}}(v_B)\\
		\notag&+\sum\limits_{\vec{m_1}+\vec{m_2}+\vec{m_3}+\vec{m_4}+\vec{m_5}=(i-i_0,j-j_0-1,k)}\delta^{i_0}\abs{u}^{j_0+1}\WDfour^{\vec{m_1}}\upchi^{i_0}\cdot\WDfour^{\vec{m_2}}\uchi^{ j_0}\cdot\WDfour^{\vec{m_3}}\speed^2\cdot\WDfour^{\vec{m_4}}\Riem^{I}_{B3B}\cdot\WDfour^{\vec{m_5}}(v_I).\\
		\notag&+\sum\limits_{\vec{m_1}+\vec{m_2}+\vec{m_3}+\vec{m_4}+\vec{m_5}=(i-i_0-1,j-j_0,k)}\delta^{i_0+1}\abs{u}^{j_0}\WDfour^{\vec{m_1}}\upchi^{i_0}\cdot\WDfour^{\vec{m_2}}\uchi^{ j_0}\cdot\WDfour^{\vec{m_3}}\speed^2\cdot\WDfour^{\vec{m_4}}\Riem^{I}_{B4B}\cdot\WDfour^{\vec{m_5}}(v_I).
	\end{align}
	
	Letting $I=1,2,3,4$, we have:
	\begin{align}\label{EQ:varvG}
		\Fluid{\varv}{G}_I=&\sum\limits_{\substack{\vec{m_1}+\vec{m_2}+\vec{m_3}+\vec{m_4}=(i-i_0,j,k-k_0)\\\abs{\vec{m_3}}\geq 1
		}}\delta^{i_0}\abs{u}^{k_0}\WDfour^{\vec{m_1}}(\upeta+\ueta)^{k_0}\cdot\WDfour^{\vec{m_2}}\omega^{i_0}\cdot\WDfour^{\vec{m_3}}v_4\cdot\WDfour^{\vec{m_4}}\Dfour_3v_{I}\\
		\notag&+\sum\limits_{\substack{\vec{m_1}+\vec{m_2}+\vec{m_3}+\vec{m_4}+\vec{m_5}=(i-i_0,j,k-k_0)\\\vec{m_1},\vec{m_2},\vec{m_4}\neq(i-i_0,j,k-k_0)}}\delta^{i_0}\abs{u}^{k_0}\WDfour^{\vec{m_1}}(\upeta+\ueta)^{k_0}\cdot\WDfour^{\vec{m_2}}\omega^{i_0}\cdot\WDfour^{\vec{m_3}}v_4\cdot\WDfour^{\vec{m_4}}(\uchi,\uomega,\upzeta)\cdot\WDfour^{\vec{m_5}}(v_{I})\\
		\notag&+\sum\limits_{\vec{m_1}+\vec{m_2}+\vec{m_3}+\vec{m_4}+\vec{m_5}=(i-i_0,j,k-k_0-1)}\delta^{i_0}\abs{u}^{k_0+1}\WDfour^{\vec{m_1}}(\upeta+\ueta)^{k_0}\cdot\WDfour^{\vec{m_2}}\omega^{i_0}\cdot\WDfour^{\vec{m_3}}v_4\cdot\WDfour^{\vec{m_4}}\Riem^{J}_{I 3A}\cdot\WDfour^{\vec{m_5}}(v_{J})\\
		\notag&+\sum\limits_{\vec{m_1}+\vec{m_2}+\vec{m_3}+\vec{m_4}+\vec{m_5}=(i-i_0-1,j,k-k_0)}\delta^{i_0+1}\abs{u}^{k_0}\WDfour^{\vec{m_1}}(\upeta+\ueta)^{k_0}\cdot\WDfour^{\vec{m_2}}\omega^{i_0}\cdot\WDfour^{\vec{m_3}}v_4\cdot\WDfour^{\vec{m_4}}\Riem^{J}_{I 34}\cdot\WDfour^{\vec{m_5}}(v_J)\\
		\notag&+\sum\limits_{\substack{\vec{m_1}+\vec{m_2}+\vec{m_3}+\vec{m_4}=(i,j-j_0,k-k_0)\\\abs{\vec{m_3}}\geq 1
		}}\abs{u}^{j_0+k_0}\WDfour^{\vec{m_1}}(\upeta+\ueta)^{k_0}\cdot\WDfour^{\vec{m_2}}\uomega^{j_0}\cdot\WDfour^{\vec{m_3}}(v_3)\cdot\WDfour^{\vec{m_4}}\Dfour_4v_{I}\\
		\notag&+\sum\limits_{\substack{\vec{m_1}+\vec{m_2}+\vec{m_3}+\vec{m_4}+\vec{m_5}=(i,j-j_0,k-k_0)\\\vec{m_1},\vec{m_2},\vec{m_4}\neq(i,j-j_0,k-k_0)}}\abs{u}^{j_0+k_0}\WDfour^{\vec{m_1}}(\upeta+\ueta)^{k_0}\cdot\WDfour^{\vec{m_2}}\uomega^{j_0}\cdot\WDfour^{\vec{m_3}}(v_3)\cdot\WDfour^{\vec{m_4}}(\upchi,\omega)\cdot\WDfour^{\vec{m_5}}(v_{I})\\
		\notag&+\sum\limits_{\vec{m_1}+\vec{m_2}+\vec{m_3}+\vec{m_4}+\vec{m_5}=(i,j-j_0,k-k_0-1)}\abs{u}^{j_0+k_0+1}\WDfour^{\vec{m_1}}(\upeta+\ueta)^{k_0}\cdot\WDfour^{\vec{m_2}}\uomega^{j_0}\cdot\WDfour^{\vec{m_3}}(v_3)\cdot\WDfour^{\vec{m_4}}\Riem^{J}_{I4A}\cdot\WDfour^{\vec{m_5}}(v_{J})\\
		\notag&+\sum\limits_{\vec{m_1}+\vec{m_2}+\vec{m_3}+\vec{m_4}+\vec{m_5}=(i,j-j_0-1,k-k_0)}\abs{u}^{j_0+k_0+1}\WDfour^{\vec{m_1}}(\upeta+\ueta)^{k_0}\cdot\Dfour^{\vec{m_2}}\uomega^{j_0}\cdot\WDfour^{\vec{m_3}}(v_3)\cdot\WDfour^{\vec{m_4}}\Riem^{J}_{I34}\cdot\WDfour^{\vec{m_5}}(v_J)\\
		\notag&+\sum\limits_{\substack{\vec{m_1}+\vec{m_2}+\vec{m_3}+\vec{m_4}=(i-i_0,j-j_0,k)\\\abs{\vec{m_3}}\geq 1
		}}\delta^{i_0}\abs{u}^{j_0}\WDfour^{\vec{m_1}}\upchi^{i_0}\cdot\WDfour^{\vec{m_2}}\uchi^{j_0}\cdot\WDfour^{\vec{m_3}} (v_B)\cdot\WDfour^{\vec{m_4}}\Dfour_Bv_{I}\\
		\notag&+\sum\limits_{\substack{\vec{m_1}+\vec{m_2}+\vec{m_3}+\vec{m_4}+\vec{m_5}=(i-i_0,j-j_0,k)\\\vec{m_1},\vec{m_2},\vec{m_4}\neq(i-i_0,j-j_0,k)}}\delta^{i_0}\abs{u}^{j_0}\WDfour^{\vec{m_1}}\upchi^{i_0}\cdot\WDfour^{\vec{m_2}}\uchi^{j_0}\cdot\WDfour^{\vec{m_3}} (v_B)\cdot\WDfour^{\vec{m_4}}(\upeta_B+\ueta_B)\cdot\WDfour^{\vec{m_5}}(v_{I})\\
		\notag&+\sum\limits_{\vec{m_1}+\vec{m_2}+\vec{m_3}+\vec{m_4}+\vec{m_5}=(i-i_0,j-j_0-1,k)}\delta^{i_0}\abs{u}^{j_0+1}\WDfour^{\vec{m_1}}\upchi^{i_0}\cdot\WDfour^{\vec{m_2}}\uchi^{ j_0}\cdot\WDfour^{\vec{m_3}} (v_B)\cdot\WDfour^{\vec{m_4}}\Riem^{J}_{I3B}\cdot\WDfour^{\vec{m_5}}(v_J).\\
		\notag&+\sum\limits_{\vec{m_1}+\vec{m_2}+\vec{m_3}+\vec{m_4}+\vec{m_5}=(i-i_0-1,j-j_0,k)}\delta^{i_0+1}\abs{u}^{j_0}\WDfour^{\vec{m_1}}\upchi^{i_0}\cdot\WDfour^{\vec{m_2}}\uchi^{ j_0}\cdot\WDfour^{\vec{m_3}} (v_B)\cdot\WDfour^{\vec{m_4}}\Riem^{J}_{I4B}\cdot\WDfour^{\vec{m_5}}(v_J).
	\end{align}
	\begin{align}\label{EQ:varhG}
		\Fluid{\varh}{G}_{I}=&\sum\limits_{\substack{\vec{m_1}+\vec{m_2}+\vec{m_3}=(i-i_0,j,k-k_0)\\\abs{\vec{m_3}}<\abs{\vec{m}}
		}}\delta^{i_0}\abs{u}^{k_0}\WDfour^{\vec{m_1}}(\upeta+\ueta)^{k_0}\cdot\WDfour^{\vec{m_2}}\omega^{i_0}\cdot\WDfour^{\vec{m_3}}\Dfour_3\lnenth\cdot \gfour(e_4,e_I)\\
		\notag&+\sum\limits_{\substack{\vec{m_1}+\vec{m_2}+\vec{m_3}+\vec{m_4}=(i-i_0,j,k-k_0)\\\vec{m_1},\vec{m_2},\vec{m_3}\neq(i-i_0,j,k-k_0)}}\delta^{i_0}\abs{u}^{k_0}\WDfour^{\vec{m_1}}(\upeta+\ueta)^{k_0}\cdot\WDfour^{\vec{m_2}}\omega^{i_0}\cdot\WDfour^{\vec{m_3}}(\uchi,\uomega,\upzeta)\cdot\WDfour^{\vec{m_4}}\lnenth\cdot \gfour(e_4,e_I)\\
		\notag&+\sum\limits_{\substack{\vec{m_1}+\vec{m_2}+\vec{m_3}=(i,j-j_0,k-k_0)\\\abs{\vec{m_3}}<\abs{\vec{m}}
		}}\abs{u}^{j_0+k_0}\WDfour^{\vec{m_1}}(\upeta+\ueta)^{k_0}\cdot\WDfour^{\vec{m_2}}\uomega^{j_0}\cdot\WDfour^{\vec{m_3}}\Dfour_4\lnenth\cdot \gfour(e_3,e_I)\\
		\notag&+\sum\limits_{\substack{\vec{m_1}+\vec{m_2}+\vec{m_3}+\vec{m_4}=(i,j-j_0,k-k_0)\\\vec{m_1},\vec{m_2},\vec{m_3}\neq(i,j-j_0,k-k_0)}}\abs{u}^{j_0+k_0}\WDfour^{\vec{m_1}}(\upeta+\ueta)^{k_0}\cdot\WDfour^{\vec{m_2}}\uomega^{j_0}\cdot\WDfour^{\vec{m_3}}(\upchi,\omega)\cdot\WDfour^{\vec{m_4}}\lnenth\cdot \gfour(e_3,e_I)\\
		\notag&+\sum\limits_{\substack{\vec{m_1}+\vec{m_2}+\vec{m_3}=(i-i_0,j-j_0,k)\\\abs{\vec{m_3}}<\abs{\vec{m}}
		}}\delta^{i_0}\abs{u}^{j_0}\WDfour^{\vec{m_1}}\upchi^{i_0}\cdot\WDfour^{\vec{m_2}}\uchi^{j_0}\cdot\WDfour^{\vec{m_3}}\Dfour_B\lnenth\cdot \gfour(e_B,e_I)\\
		\notag&+\sum\limits_{\substack{\vec{m_1}+\vec{m_2}+\vec{m_3}+\vec{m_4}=(i-i_0,j-j_0,k)\\\vec{m_1},\vec{m_2},\vec{m_3}\neq(i-i_0,j-j_0,k)}}\delta^{i_0}\abs{u}^{j_0}\WDfour^{\vec{m_1}}\upchi^{i_0}\cdot\WDfour^{\vec{m_2}}\uchi^{j_0}\cdot\WDfour^{\vec{m_3}}(\upeta_B+\ueta_B)\cdot\WDfour^{\vec{m_4}}\lnenth\cdot \gfour(e_B,e_I)\\
		\notag&+\sum\limits_{\substack{\vec{m_1}+\vec{m_2}+\vec{m_3}+\vec{m_4}+\vec{m_5}=(i-i_0,j,k-k_0)\\\abs{\vec{m_3}}+\abs{\vec{m_5}}\geq 1
		}}\delta^{i_0}\abs{u}^{k_0}\WDfour^{\vec{m_1}}(\upeta+\ueta)^{k_0}\cdot\WDfour^{\vec{m_2}}\omega^{i_0}\cdot\WDfour^{\vec{m_3}}v_4\cdot\WDfour^{\vec{m_4}}\Dfour_3\lnenth\cdot\WDfour^{\vec{m_5}}(v_{I})\\
		\notag&+\sum\limits_{\substack{\vec{m_1}+\vec{m_2}+\vec{m_3}+\vec{m_4}+\vec{m_5}+\vec{m_6}=(i-i_0,j,k-k_0)\\\vec{m_1},\vec{m_2},\vec{m_4}\neq(i-i_0,j,k-k_0)}}\delta^{i_0}\abs{u}^{k_0}\WDfour^{\vec{m_1}}(\upeta+\ueta)^{k_0}\cdot\WDfour^{\vec{m_2}}\omega^{i_0}\cdot\WDfour^{\vec{m_3}}v_4\cdot\WDfour^{\vec{m_4}}(\uchi,\uomega,\upzeta)\cdot\WDfour^{\vec{m_5}}\lnenth\cdot\WDfour^{\vec{m_6}}(v_{I})\\
		\notag&+\sum\limits_{\substack{\vec{m_1}+\vec{m_2}+\vec{m_3}+\vec{m_4}+\vec{m_5}=(i,j-j_0,k-k_0)\\\abs{\vec{m_3}}+\abs{\vec{m_5}}\geq 1
		}}\abs{u}^{j_0+k_0}\WDfour^{\vec{m_1}}(\upeta+\ueta)^{k_0}\cdot\WDfour^{\vec{m_2}}\uomega^{j_0}\cdot\WDfour^{\vec{m_3}}(v_3)\cdot\WDfour^{\vec{m_4}}\Dfour_4\lnenth\cdot\WDfour^{\vec{m_5}}(v_{I})\\
		\notag&+\sum\limits_{\substack{\vec{m_1}+\vec{m_2}+\vec{m_3}+\vec{m_4}+\vec{m_5}+\vec{m_6}=(i,j-j_0,k-k_0)\\\vec{m_1},\vec{m_2},\vec{m_4}\neq(i,j-j_0,k-k_0)}}\abs{u}^{j_0+k_0}\WDfour^{\vec{m_1}}(\upeta+\ueta)^{k_0}\cdot\WDfour^{\vec{m_2}}\uomega^{j_0}\cdot\WDfour^{\vec{m_3}}(v_3)\cdot\WDfour^{\vec{m_4}}(\upchi,\omega)\cdot\WDfour^{\vec{m_5}}\lnenth\cdot\WDfour^{\vec{m_6}}(v_{I})\\
		\notag&+\sum\limits_{\substack{\vec{m_1}+\vec{m_2}+\vec{m_3}+\vec{m_4}+\vec{m_5}=(i-i_0,j-j_0,k)\\\abs{\vec{m_3}}+\abs{\vec{m_5}}\geq 1
		}}\delta^{i_0}\abs{u}^{j_0}\WDfour^{\vec{m_1}}\upchi^{i_0}\cdot\WDfour^{\vec{m_2}}\uchi^{j_0}\cdot\WDfour^{\vec{m_3}} (v_B)\cdot\WDfour^{\vec{m_4}}\Dfour_B\lnenth\cdot\WDfour^{\vec{m_5}}(v_{I})\\
		\notag&+\sum\limits_{\substack{\vec{m_1}+\vec{m_2}+\vec{m_3}+\vec{m_4}+\vec{m_5}+\vec{m_6}=(i-i_0,j-j_0,k)\\\vec{m_1},\vec{m_2},\vec{m_4}\neq(i-i_0,j-j_0,k)}}\delta^{i_0}\abs{u}^{j_0}\WDfour^{\vec{m_1}}\upchi^{i_0}\cdot\WDfour^{\vec{m_2}}\uchi^{j_0}\cdot\WDfour^{\vec{m_3}} (v_B)\cdot\WDfour^{\vec{m_4}}(\upeta_B+\ueta_B)\cdot\WDfour^{\vec{m_5}}\lnenth\cdot\WDfour^{\vec{m_6}}(v_{I}).
	\end{align}
	\begin{align}\label{EQ:varsG}
		\Fluid{\vars}{G}_{I}=&\sum\limits_{\substack{\vec{m_1}+\vec{m_2}+\vec{m_3}+\vec{m_4}=(i-i_0,j,k-k_0)\\\abs{\vec{m_3}}\geq 1
		}}\delta^{i_0}\abs{u}^{k_0}\WDfour^{\vec{m_1}}(\upeta+\ueta)^{k_0}\cdot\WDfour^{\vec{m_2}}\omega^{i_0}\cdot\WDfour^{\vec{m_3}}\frac{\temp}{\enth}\cdot\WDfour^{\vec{m_4}}\Dfour_3\ent\cdot \gfour(e_4,e_I)\\
		\notag&+\sum\limits_{\substack{\vec{m_1}+\vec{m_2}+\vec{m_3}+\vec{m_4}+\vec{m_5}=(i-i_0,j,k-k_0)\\\vec{m_1},\vec{m_2},\vec{m_4}\neq(i-i_0,j,k-k_0)}}\delta^{i_0}\abs{u}^{k_0}\WDfour^{\vec{m_1}}(\upeta+\ueta)^{k_0}\cdot\WDfour^{\vec{m_2}}\omega^{i_0}\cdot\WDfour^{\vec{m_3}}\frac{\temp}{\enth}\cdot\WDfour^{\vec{m_4}}(\uchi,\uomega,\upzeta)\cdot\WDfour^{\vec{m_5}}\ent\cdot \gfour(e_4,e_I)\\
		\notag&+\sum\limits_{\substack{\vec{m_1}+\vec{m_2}+\vec{m_3}+\vec{m_4}=(i,j-j_0,k-k_0)\\\abs{\vec{m_3}}\geq 1
		}}\abs{u}^{j_0+k_0}\WDfour^{\vec{m_1}}(\upeta+\ueta)^{k_0}\cdot\WDfour^{\vec{m_2}}\uomega^{j_0}\cdot\WDfour^{\vec{m_3}}\frac{\temp}{\enth}\cdot\WDfour^{\vec{m_4}}\Dfour_4\ent\cdot \gfour(e_3,e_I)\\
		\notag&+\sum\limits_{\substack{\vec{m_1}+\vec{m_2}+\vec{m_3}+\vec{m_4}+\vec{m_5}=(i,j-j_0,k-k_0)\\\vec{m_1},\vec{m_2},\vec{m_4}\neq(i,j-j_0,k-k_0)}}\abs{u}^{j_0+k_0}\WDfour^{\vec{m_1}}(\upeta+\ueta)^{k_0}\cdot\WDfour^{\vec{m_2}}\uomega^{j_0}\cdot\WDfour^{\vec{m_3}}\frac{\temp}{\enth}\cdot\WDfour^{\vec{m_4}}(\upchi,\omega)\cdot\WDfour^{\vec{m_5}}\ent\cdot \gfour(e_3,e_I)\\
		\notag&+\sum\limits_{\substack{\vec{m_1}+\vec{m_2}+\vec{m_3}+\vec{m_4}=(i-i_0,j-j_0,k)\\\abs{\vec{m_3}}\geq 1
		}}\delta^{i_0}\abs{u}^{j_0}\WDfour^{\vec{m_1}}\upchi^{i_0}\cdot\WDfour^{\vec{m_2}}\uchi^{j_0}\cdot\WDfour^{\vec{m_3}}\frac{\temp}{\enth}\cdot\WDfour^{\vec{m_4}}\Dfour_B\ent\cdot \gfour(e_B,e_I)\\
		\notag&+\sum\limits_{\substack{\vec{m_1}+\vec{m_2}+\vec{m_3}+\vec{m_4}+\vec{m_5}=(i-i_0,j-j_0,k)\\\vec{m_1},\vec{m_2},\vec{m_4}\neq(i-i_0,j-j_0,k)}}\delta^{i_0}\abs{u}^{j_0}\WDfour^{\vec{m_1}}\upchi^{i_0}\cdot\WDfour^{\vec{m_2}}\uchi^{j_0}\cdot\WDfour^{\vec{m_3}}\frac{\temp}{\enth}\cdot\WDfour^{\vec{m_4}}(\upeta_B+\ueta_B)\cdot\WDfour^{\vec{m_5}}\ent\cdot \gfour(e_B,e_I).
	\end{align}
	\begin{align}\label{EQ:H}
		\mathfrak{H}_{(i,j,k)}=&\sum\limits_{\substack{\vec{m_1}+\vec{m_2}+\vec{m_3}+\vec{m_4}=(i-i_0,j,k-k_0)\\\abs{\vec{m_3}}\geq 1
		}}\delta^{i_0}\abs{u}^{k_0}\WDfour^{\vec{m_1}}(\upeta+\ueta)^{k_0}\cdot\WDfour^{\vec{m_2}}\omega^{i_0}\cdot\WDfour^{\vec{m_3}}v_4\cdot\WDfour^{\vec{m_4}}\Dfour_3\ent\\
		\notag&+\sum\limits_{\substack{\vec{m_1}+\vec{m_2}+\vec{m_3}+\vec{m_4}+\vec{m_5}=(i-i_0,j,k-k_0)\\\vec{m_1},\vec{m_2},\vec{m_4}\neq(i-i_0,j,k-k_0)}}\delta^{i_0}\abs{u}^{k_0}\WDfour^{\vec{m_1}}(\upeta+\ueta)^{k_0}\cdot\WDfour^{\vec{m_2}}\omega^{i_0}\cdot\WDfour^{\vec{m_3}}v_4\cdot\WDfour^{\vec{m_4}}(\uchi,\uomega,\upzeta)\cdot\WDfour^{\vec{m_5}}\ent\\
		\notag&+\sum\limits_{\substack{\vec{m_1}+\vec{m_2}+\vec{m_3}+\vec{m_4}=(i,j-j_0,k-k_0)\\\abs{\vec{m_3}}\geq 1
		}}\abs{u}^{j_0+k_0}\WDfour^{\vec{m_1}}(\upeta+\ueta)^{k_0}\cdot\WDfour^{\vec{m_2}}\uomega^{j_0}\cdot\WDfour^{\vec{m_3}}(v_3)\cdot\WDfour^{\vec{m_4}}\Dfour_4\ent\\
		\notag&+\sum\limits_{\substack{\vec{m_1}+\vec{m_2}+\vec{m_3}+\vec{m_4}+\vec{m_5}=(i,j-j_0,k-k_0)\\\vec{m_1},\vec{m_2},\vec{m_4}\neq(i,j-j_0,k-k_0)}}\abs{u}^{j_0+k_0}\WDfour^{\vec{m_1}}(\upeta+\ueta)^{k_0}\cdot\WDfour^{\vec{m_2}}\uomega^{j_0}\cdot\WDfour^{\vec{m_3}}(v_3)\cdot\WDfour^{\vec{m_4}}(\upchi,\omega)\cdot\WDfour^{\vec{m_5}}\ent\\
		\notag&+\sum\limits_{\substack{\vec{m_1}+\vec{m_2}+\vec{m_3}+\vec{m_4}=(i-i_0,j-j_0,k)\\\abs{\vec{m_3}}\geq 1
		}}\delta^{i_0}\abs{u}^{j_0}\WDfour^{\vec{m_1}}\upchi^{i_0}\cdot\WDfour^{\vec{m_2}}\uchi^{j_0}\cdot\WDfour^{\vec{m_3}} (v_B)\cdot\WDfour^{\vec{m_4}}\Dfour_B\ent\\
		\notag&+\sum\limits_{\substack{\vec{m_1}+\vec{m_2}+\vec{m_3}+\vec{m_4}+\vec{m_5}=(i-i_0,j-j_0,k)\\\vec{m_1},\vec{m_2},\vec{m_4}\neq(i-i_0,j-j_0,k)}}\delta^{i_0}\abs{u}^{j_0}\WDfour^{\vec{m_1}}\upchi^{i_0}\cdot\WDfour^{\vec{m_2}}\uchi^{j_0}\cdot\WDfour^{\vec{m_3}} (v_B)\cdot\WDfour^{\vec{m_4}}(\upeta_B+\ueta_B)\cdot\WDfour^{\vec{m_5}}\ent.
	\end{align}
	
\end{proposition}
\begin{proof}[Proof of Proposition \ref{PR:EQVariations}]
	\eqref{EQ:FOE3} are derived by straightforward computations. To be specific, we commute $\WDfour^{\vec{m}}$ with \eqref{EQ:FOE2} where $\vec{m}=(i,j,k)$, keep the top-order term (in the form of $\Dfour^{\vec{n}}$ with $\abs{\vec{n}}=\abs{\vec{m}+1}$) on the left, and use Proposition \ref{PR:CommuteDijk} to express the commutator terms. In particular, we note that by the construction of $\varv^0$ in Definition \ref{DE:EnergyCurrent},
	we take $\Dfour_\mu$ derivative of \eqref{EQ:gvarvv=0} to obtain:
	\begin{align}
		\Dfour_\mu\varv^0 =&~-\frac{1}{\gfour_{0\beta}v^\beta}\left(\gfour_{\alpha\beta}\varv^\alpha\Dfour_\mu v^\beta+\gfour_{a\beta}(\Dfour_\mu\varv^a)v^\beta\right).
	\end{align}
	Meanwhile, by $\gfour(v,v)=-1$, we have
	\begin{align}
		\Dfour_\mu\left(\WDfour^{\vec{m}}v^0\right) =&~-\frac{1}{\gfour_{0\beta}v^\beta}\gfour_{ab}(\Dfour_\mu\WDfour^{\vec{m}}v^a)v^b\\
		\notag&+\sum\limits_{\substack{\vec{m_1}+\vec{m_2}=\vec{m}\\\vec{m_1},\vec{m_2}\neq\vec{m}}}\gfour_{\alpha\beta}\Dfour_\mu\WDfour^{\vec{m_1}}v^\alpha\WDfour^{\vec{m_2}}v^\beta.
	\end{align}
	Therefore, $\Dfour_\mu\varv^0-\Dfour_\mu\left(\WDfour^{\vec{m}}v^0\right)$ only consists of lower-order terms (no $(\abs{\vec{m}}+1)$-derivative terms), which could be absorbed by inhomogeneous terms $\mathfrak{F},\mathfrak{G}$. We omit listing every lower-order terms provided by $\Dfour_\mu\varv^0-\Dfour_\mu\left(\WDfour^{\vec{m}}v^0\right)$ in the equations \eqref{EQ:varvF}-\eqref{EQ:varvG} as there are harmless small terms since the size of $\varv^0-\WDfour^{\vec{m}}v^0$ is $\BAfluid^2$.
\end{proof}
\subsection{Estimates for the Inhomogeneous Terms of Equation of Variations}\label{SS:EstimatesFluidterms}
\begin{proposition}[Below-top-order estimates of inhomogenuous terms of equation of variations]\label{PR:EstimateofFGH}
	Let $\BAsize,\BAfluid$ be as in Section \ref{SS:BA}. Let $\mathfrak{F}:=\mathfrak{F}_{(i,j,k)},\mathfrak{G}:=\mathfrak{G}_{(i,j,k)},\mathfrak{H}:=\mathfrak{H}_{(i,j,k)}$ be defined as in Proposition \ref{PR:EQVariations} with $i+j+k\leq N-2$ or $i+j+k\leq N-1$ with $k\geq 1$. Then under the bootstrap assumptions, we have the following estimates\footnote{On the RHS of \eqref{ES:EstimateFGH}, the first term represents the bound purely for fluid terms, and the second term with $a$ appearing represents the bound for the gravitational impact. We use the same notation in the rest of this section.}:
	\begin{subequations}\label{ES:EstimateFGH}
		\begin{align}
			\label{ES:EstimateFGH1}\norm{\mathfrak{F}_{(i,j,k)}}_{L^2(\suubar)}\lesssim&~\BAfluid^2\delta^{-1}b^{-\frac{1}{2}}\abs{u}+\BAsize\BAfluid a^{\frac{1}{2}} b^{-\frac{1}{2}},\\
			\label{ES:EstimateFGH2}\norm{(\mathfrak{G}^\alpha_{(i,j,k)})_3}_{L^2(\suubar)}\lesssim&~\BAfluid^2\delta^{-1}b^{-1}\abs{u}+\BAsize\BAfluid a^{\frac{1}{2}} b^{-1},\\
			\label{ES:EstimateFGH3}\norm{(\mathfrak{G}^\alpha_{(i,j,k)})_B}_{L^2(\suubar)}\lesssim&~\BAfluid^2\delta^{-1}b^{-\frac{1}{2}}\abs{u}+\BAsize\delta^{-1}a^{\frac{3}{8}} b^{-1}\abs{u},\\
			\label{ES:EstimateFGH4}\norm{(\mathfrak{G}^\alpha_{(i,j,k)})_4}_{L^2(\suubar)}\lesssim&~\BAfluid^2\delta^{-1}\abs{u}+\BAsize\BAfluid\delta^{-1}a^{\frac{1}{2}} b^{-1}\abs{u},\\
			\label{ES:EstimateFGH5}\norm{\mathfrak{H}_{(i,j,k)}}_{L^2(\suubar)}\lesssim&~\BAfluid^2\delta^{-1}b^{-\frac{1}{2}}\abs{u}+\BAsize\BAfluid a^{\frac{1}{2}} b^{-\frac{1}{2}}.
		\end{align}
	\end{subequations}
\end{proposition}
\subsubsection{Proof of Proposition \ref{PR:EstimateofFGH}}\label{SSS:ProofofBulk}
In this section, we silently use preliminary results Lemma \ref{LE:EstimatesofTensorDifference} - Lemma \ref{LE:BoundsforProducts}, and bootstrap assumptions \eqref{BA} for bounding every individual terms without mentioning. Let $\vec{m}:=(i,j,k)$ be any non-negative $3$-tuples such that $\abs{\vec{m}}\leq N-2$ or $\abs{\vec{m}}\leq N-1$ with $k\geq 1$. Let $\WDfour^{\vec{m}}$ be as in Definition \ref{DE:operators}.
\begin{proof}[Proof of \eqref{ES:EstimateFGH1}]
	As in Proposition \ref{PR:EQVariations}, $\mathfrak{F}_{(i,j,k)}:=\Fluid{\varh}{F}+\Fluid{\varv}{F}$. We first estimate terms in $\Fluid{\varh}{F}$. For the first term on the RHS of \eqref{EQ:Fvarh}, we have:
	\begin{align}\label{PV:varhF1}
		&\sum\limits_{\substack{\vec{m_1}+\vec{m_2}+\vec{m_3}+\vec{m_4}=(i-i_0,j,k-k_0)\\\abs{\vec{m_3}}\geq1
		}}\delta^{i_0}\abs{u}^{k_0}\norm{\WDfour^{\vec{m_1}}(\upeta+\ueta)^{k_0}\cdot\WDfour^{\vec{m_2}}\omega^{i_0}\cdot\WDfour^{\vec{m_3}}v_4\cdot\WDfour^{\vec{m_4}}\Dfour_3\lnenth}_{L^2(\suubar)}\\
		\notag\lesssim&\sum\limits_{\substack{\vec{m_1}+\vec{m_2}+\vec{m_3}+\vec{m_4}=(i-i_0,j,k-k_0)\\\abs{\vec{m_1}+\vec{m_2}}> \frac{N-1}{2}
		}}\delta^{i_0}\abs{u}^{k_0}\norm{\WDfour^{\vec{m_1}}\upsi^{k_0}\cdot\WDfour^{\vec{m_2}}\psi^{i_0}}_{L^2(\suubar)}\cdot\norm{\WDfour^{\vec{m_3}}v_4}_{L^\infty(\suubar)}\cdot\norm{\WDfour^{\vec{m_4}}\Dfour_3\lnenth}_{L^\infty(\suubar)}\\
		\notag&+\sum\limits_{\substack{\vec{m_1}+\vec{m_2}+\vec{m_3}+\vec{m_4}=(i-i_0,j,k-k_0)\\\abs{\vec{m_3}}> \frac{N-1}{2}
		}}\delta^{i_0}\abs{u}^{k_0}\norm{\WDfour^{\vec{m_1}}\upsi^{k_0}\cdot\WDfour^{\vec{m_2}}\psi^{i_0}}_{L^\infty(\suubar)}\cdot\norm{\WDfour^{\vec{m_3}}v_4}_{L^2(\suubar)}\cdot\norm{\WDfour^{\vec{m_4}}\Dfour_3\lnenth}_{L^\infty(\suubar)}\\
		\notag&+\sum\limits_{\substack{\vec{m_1}+\vec{m_2}+\vec{m_3}+\vec{m_4}=(i-i_0,j,k-k_0)\\\abs{\vec{m_4}+1}> \frac{N-1}{2}
		}}\delta^{i_0}\abs{u}^{k_0}\norm{\WDfour^{\vec{m_1}}\upsi^{k_0}\cdot\WDfour^{\vec{m_2}}\psi^{i_0}}_{L^\infty(\suubar)}\cdot\norm{\WDfour^{\vec{m_3}}v_4}_{L^\infty(\suubar)}\cdot\norm{\WDfour^{\vec{m_4}}\Dfour_3\lnenth}_{L^2(\suubar)}\\
		\notag\lesssim&\left(\BAsize\frac{\delta a^{\frac{1}{2}}}{\abs{u}}\right)^{i_0+k_0}\cdot\BAfluid^2 b^{\frac{1}{2}}.
	\end{align}
	For the second term on the RHS of \eqref{EQ:Fvarh}:
	\begin{align}\label{PV:varhF2}
		&\sum\limits_{\substack{\vec{m_1}+\vec{m_2}+\vec{m_3}+\vec{m_4}+\vec{m_5}=(i-i_0,j,k-k_0)\\\vec{m_1},\vec{m_2},\vec{m_4}\neq(i-i_0,j,k-k_0)}}\delta^{i_0}\abs{u}^{k_0}\norm{\WDfour^{\vec{m_1}}(\upeta+\ueta)^{k_0}\cdot\WDfour^{\vec{m_2}}\omega^{i_0}\cdot\WDfour^{\vec{m_3}}v_4\cdot\WDfour^{\vec{m_4}}(\uchi,\uomega,\upzeta)\cdot\WDfour^{\vec{m_5}}\lnenth}_{L^2(\suubar)}\\
		\notag\lesssim&\sum\limits_{\substack{\vec{m_1}+\vec{m_2}+\vec{m_3}+\vec{m_4}=(i-i_0,j,k-k_0)\\\abs{\vec{m_1}+\vec{m_2}}> \frac{N-1}{2}
		}}\delta^{i_0}\abs{u}^{k_0}\norm{\WDfour^{\vec{m_1}}\upsi^{k_0+1}\cdot\WDfour^{\vec{m_2}}\psi^{i_0}}_{L^2(\suubar)}\cdot\norm{\WDfour^{\vec{m_3}}v_4}_{L^\infty(\suubar)}\cdot\norm{\WDfour^{\vec{m_4}}\lnenth}_{L^\infty(\suubar)}\\
		\notag&+\sum\limits_{\substack{\vec{m_1}+\vec{m_2}+\vec{m_3}+\vec{m_4}=(i-i_0,j,k-k_0)\\\abs{\vec{m_3}}> \frac{N-1}{2}
		}}\delta^{i_0}\abs{u}^{k_0}\norm{\WDfour^{\vec{m_1}}\upsi^{k_0+1}\cdot\WDfour^{\vec{m_2}}\psi^{i_0}}_{L^\infty(\suubar)}\cdot\norm{\WDfour^{\vec{m_3}}v_4}_{L^2(\suubar)}\cdot\norm{\WDfour^{\vec{m_4}}\lnenth}_{L^\infty(\suubar)}\\
		\notag&+\sum\limits_{\substack{\vec{m_1}+\vec{m_2}+\vec{m_3}+\vec{m_4}=(i-i_0,j,k-k_0)\\\abs{\vec{m_4}}> \frac{N-1}{2}
		}}\delta^{i_0}\abs{u}^{k_0}\norm{\WDfour^{\vec{m_1}}\upsi^{k_0+1}\cdot\WDfour^{\vec{m_2}}\psi^{i_0}}_{L^\infty(\suubar)}\cdot\norm{\WDfour^{\vec{m_3}}v_4}_{L^\infty(\suubar)}\cdot\norm{\WDfour^{\vec{m_4}}\lnenth}_{L^2(\suubar)}\\
		\notag\lesssim&\left(\BAsize\frac{\delta a^{\frac{1}{2}}}{\abs{u}}\right)^{i_0+k_0+1}\cdot\BAfluid b^{\frac{1}{2}}.
	\end{align}
	The third and fourth terms of \eqref{EQ:Fvarh} can be estimates in a similar fashion:
	\begin{align}\label{PV:varhF3}
		&\sum\limits_{\substack{\vec{m_1}+\vec{m_2}+\vec{m_3}+\vec{m_4}=(i,j-j_0,k-k_0)\\\abs{\vec{m_3}}\geq 1
		}}\abs{u}^{j_0+k_0}\norm{\WDfour^{\vec{m_1}}(\upeta+\ueta)^{k_0}\cdot\WDfour^{\vec{m_2}}\uomega^{j_0}\cdot\WDfour^{\vec{m_3}}(v_3)\cdot\WDfour^{\vec{m_4}}\Dfour_4\lnenth}_{L^2(\suubar)}\\
		\notag\lesssim&\sum\limits_{\substack{\vec{m_1}+\vec{m_2}+\vec{m_3}=(i,j-j_0,k-k_0)\\\abs{\vec{m_1}}> \frac{N-1}{2}
		}}\abs{u}^{j_0+k_0}\norm{\WDfour^{\vec{m_1}}\upsi^{j_0+k_0}}_{L^2(\suubar)}\cdot\norm{\WDfour^{\vec{m_2}}(v_3)}_{L^\infty(\suubar)}\cdot\norm{\WDfour^{\vec{m_3}}\Dfour_4\lnenth}_{L^\infty(\suubar)}\\
		\notag&+\sum\limits_{\substack{\vec{m_1}+\vec{m_2}+\vec{m_3}=(i,j-j_0,k-k_0)\\\abs{\vec{m_2}}> \frac{N-1}{2}
		}}\abs{u}^{j_0+k_0}\norm{\WDfour^{\vec{m_1}}\upsi^{j_0+k_0}}_{L^\infty(\suubar)}\cdot\norm{\WDfour^{\vec{m_2}}(v_3)}_{L^2(\suubar)}\cdot\norm{\WDfour^{\vec{m_3}}\Dfour_4\lnenth}_{L^\infty(\suubar)}\\
		\notag&+\sum\limits_{\substack{\vec{m_1}+\vec{m_2}+\vec{m_3}=(i,j-j_0,k-k_0)\\\abs{\vec{m_3}+1}> \frac{N-1}{2}
		}}\abs{u}^{j_0+k_0}\norm{\WDfour^{\vec{m_1}}\upsi^{j_0+k_0}}_{L^\infty(\suubar)}\cdot\norm{\WDfour^{\vec{m_2}}(v_3)}_{L^\infty(\suubar)}\cdot\norm{\WDfour^{\vec{m_3}}\Dfour_4\lnenth}_{L^2(\suubar)}\\
		\notag\lesssim&\left(\BAsize\frac{\delta a^{\frac{1}{2}}}{\abs{u}}\right)^{i_0+k_0}\cdot\BAfluid^2 b^{-\frac{1}{2}}\cdot\delta^{-1}\cdot\abs{u}.
	\end{align}
	\begin{align}\label{PV:varhF4}
		&\sum\limits_{\substack{\vec{m_1}+\vec{m_2}+\vec{m_3}+\vec{m_4}+\vec{m_5}=(i,j-j_0,k-k_0)\\\vec{m_1},\vec{m_2},\vec{m_4}\neq(i,j-j_0,k-k_0)}}\abs{u}^{j_0+k_0}\norm{\WDfour^{\vec{m_1}}(\upeta+\ueta)^{k_0}\cdot\WDfour^{\vec{m_2}}\uomega^{j_0}\cdot\WDfour^{\vec{m_3}}(v_3)\cdot\WDfour^{\vec{m_4}}(\upchi,\omega)\cdot\WDfour^{\vec{m_5}}\lnenth}_{L^2(\suubar)}\\
		\notag\lesssim&\sum\limits_{\substack{\vec{m_1}+\vec{m_2}+\vec{m_3}+\vec{m_4}=(i,j-j_0,k-k_0)\\\abs{\vec{m_1}+\vec{m_2}}> \frac{N-1}{2}
		}}\abs{u}^{j_0+k_0}\norm{\WDfour^{\vec{m_1}}\upsi^{j_0+k_0}\cdot\WDfour^{\vec{m_2}}\psi}_{L^2(\suubar)}\cdot\norm{\WDfour^{\vec{m_3}}(v_3)}_{L^\infty(\suubar)}\cdot\norm{\WDfour^{\vec{m_4}}\lnenth}_{L^\infty(\suubar)}\\
		\notag&+\sum\limits_{\substack{\vec{m_1}+\vec{m_2}+\vec{m_3}+\vec{m_4}=(i,j-j_0,k-k_0)\\\abs{\vec{m_3}}> \frac{N-1}{2}
		}}\abs{u}^{j_0+k_0}\norm{\WDfour^{\vec{m_1}}\upsi^{j_0+k_0}\cdot\WDfour^{\vec{m_2}}\psi}_{L^\infty(\suubar)}\cdot\norm{\WDfour^{\vec{m_3}}(v_3)}_{L^2(\suubar)}\cdot\norm{\WDfour^{\vec{m_4}}\lnenth}_{L^\infty(\suubar)}\\
		\notag&+\sum\limits_{\substack{\vec{m_1}+\vec{m_2}+\vec{m_3}+\vec{m_4}=(i,j-j_0,k-k_0)\\\abs{\vec{m_4}}> \frac{N-1}{2}
		}}\abs{u}^{j_0+k_0}\norm{\WDfour^{\vec{m_1}}\upsi^{j_0+k_0}\cdot\WDfour^{\vec{m_2}}\psi}_{L^\infty(\suubar)}\cdot\norm{\WDfour^{\vec{m_3}}(v_3)}_{L^\infty(\suubar)}\cdot\norm{\WDfour^{\vec{m_4}}\lnenth}_{L^2(\suubar)}\\
		\notag\lesssim&\left(\BAsize\frac{\delta a^{\frac{1}{2}}}{\abs{u}}\right)^{i_0+k_0}\cdot \BAsize a^{\frac{1}{2}}\cdot\BAfluid b^{-\frac{1}{2}}.
	\end{align}
	For the last two terms in $\Fluid{\varh}{F}$, we have:
	\begin{align}\label{PV:varhF5}
		&\sum\limits_{\substack{\vec{m_1}+\vec{m_2}+\vec{m_3}+\vec{m_4}=(i-i_0,j-j_0,k)\\\abs{\vec{m_3}}\geq 1
		}}\delta^{i_0}\abs{u}^{j_0}\norm{\WDfour^{\vec{m_1}}\upchi^{i_0}\cdot\WDfour^{\vec{m_2}}\uchi^{j_0}\cdot\WDfour^{\vec{m_3}} (v_B)\cdot\WDfour^{\vec{m_4}}\Dfour_B\lnenth}_{L^2(\suubar)}\\
		\notag\lesssim&\sum\limits_{\substack{\vec{m_1}+\vec{m_2}+\vec{m_3}+\vec{m_4}=(i-i_0,j-j_0,k)\\\abs{\vec{m_1}+\vec{m_2}}> \frac{N-1}{2}
		}}\delta^{i_0}\abs{u}^{j_0}\norm{\WDfour^{\vec{m_1}}\upsi^{j_0}\cdot\WDfour^{\vec{m_2}}\psi^{i_0}}_{L^2(\suubar)}\cdot\norm{\WDfour^{\vec{m_3}} (v_B)}_{L^\infty(\suubar)}\cdot\norm{\WDfour^{\vec{m_4}}\Dfour_B\lnenth}_{L^\infty(\suubar)}\\
		\notag&+\sum\limits_{\substack{\vec{m_1}+\vec{m_2}+\vec{m_3}+\vec{m_4}=(i-i_0,j-j_0,k)\\\abs{\vec{m_3}}> \frac{N-1}{2}
		}}\delta^{i_0}\abs{u}^{j_0}\norm{\WDfour^{\vec{m_1}}\upsi^{j_0}\cdot\WDfour^{\vec{m_2}}\psi^{i_0}}_{L^\infty(\suubar)}\cdot\norm{\WDfour^{\vec{m_3}} (v_B)}_{L^2(\suubar)}\cdot\norm{\WDfour^{\vec{m_4}}\Dfour_B\lnenth}_{L^\infty(\suubar)}\\
		\notag&+\sum\limits_{\substack{\vec{m_1}+\vec{m_2}+\vec{m_3}+\vec{m_4}=(i-i_0,j-j_0,k)\\\abs{\vec{m_4}+1}> \frac{N-1}{2}
		}}\delta^{i_0}\abs{u}^{j_0}\norm{\WDfour^{\vec{m_1}}\upsi^{j_0}\cdot\WDfour^{\vec{m_2}}\psi^{i_0}}_{L^\infty(\suubar)}\cdot\norm{\WDfour^{\vec{m_3}} (v_B)}_{L^\infty(\suubar)}\cdot\norm{\WDfour^{\vec{m_4}}\Dfour_B\lnenth}_{L^2(\suubar)}\\
		\notag\lesssim&\left(\BAsize\frac{\delta a^{\frac{1}{2}}}{\abs{u}}\right)^{i_0+j_0}\cdot\BAfluid^2 b^{\frac{1}{2}}.
	\end{align}
	\begin{align}\label{PV:varhF6}
		&\sum\limits_{\substack{\vec{m_1}+\vec{m_2}+\vec{m_3}+\vec{m_4}+\vec{m_5}=(i-i_0,j-j_0,k)\\\vec{m_1},\vec{m_2},\vec{m_4}\neq(i-i_0,j-j_0,k)}}\delta^{i_0}\abs{u}^{j_0}\norm{\WDfour^{\vec{m_1}}\upchi^{i_0}\cdot\WDfour^{\vec{m_2}}\uchi^{j_0}\cdot\WDfour^{\vec{m_3}} (v_B)\cdot\WDfour^{\vec{m_4}}(\upeta_B+\ueta_B)\cdot\WDfour^{\vec{m_5}}\lnenth}_{L^2(\suubar)}\\
		\notag\lesssim&\sum\limits_{\substack{\vec{m_1}+\vec{m_2}+\vec{m_3}+\vec{m_4}=(i-i_0,j-j_0,k)\\\abs{\vec{m_1}+\vec{m_2}}> \frac{N-1}{2}
		}}\delta^{i_0}\abs{u}^{j_0}\norm{\WDfour^{\vec{m_1}}\upsi^{j_0+1}\cdot\WDfour^{\vec{m_2}}\psi^{i_0}}_{L^2(\suubar)}\cdot\norm{\WDfour^{\vec{m_3}} (v_B)}_{L^\infty(\suubar)}\cdot\norm{\WDfour^{\vec{m_4}}\lnenth}_{L^\infty(\suubar)}\\
		\notag&+\sum\limits_{\substack{\vec{m_1}+\vec{m_2}+\vec{m_3}+\vec{m_4}=(i-i_0,j-j_0,k)\\\abs{\vec{m_3}}> \frac{N-1}{2}
		}}\delta^{i_0}\abs{u}^{j_0}\norm{\WDfour^{\vec{m_1}}\upsi^{j_0+1}\cdot\WDfour^{\vec{m_2}}\psi^{i_0}}_{L^\infty(\suubar)}\cdot\norm{\WDfour^{\vec{m_3}} (v_B)}_{L^2(\suubar)}\cdot\norm{\WDfour^{\vec{m_4}}\lnenth}_{L^\infty(\suubar)}\\
		\notag&+\sum\limits_{\substack{\vec{m_1}+\vec{m_2}+\vec{m_3}+\vec{m_4}=(i-i_0,j-j_0,k)\\\abs{\vec{m_4}}> \frac{N-1}{2}
		}}\delta^{i_0}\abs{u}^{j_0}\norm{\WDfour^{\vec{m_1}}\upsi^{j_0+1}\cdot\WDfour^{\vec{m_2}}\psi^{i_0}}_{L^\infty(\suubar)}\cdot\norm{\WDfour^{\vec{m_3}} (v_B)}_{L^\infty(\suubar)}\cdot\norm{\WDfour^{\vec{m_4}}\lnenth}_{L^2(\suubar)}\\
		\notag\lesssim&\left(\BAsize\frac{\delta a^{\frac{1}{2}}}{\abs{u}}\right)^{i_0+j_0+1}\cdot\BAfluid^2b^{\frac{1}{2}}.
	\end{align}
	Combining \eqref{PV:varhF1}-\eqref{PV:varhF6}, we obtain the estimate for $\Fluid{\varh}{F}$:
	\begin{align}\label{PV:varhF}
		\norm{\Fluid{\varh}{F}}_{L^2(\suubar)}\lesssim&~\BAfluid^2\delta^{-1}b^{-\frac{1}{2}}\abs{u}+\BAsize\BAfluid a^{\frac{1}{2}} b^{-\frac{1}{2}}.
	\end{align}
	
	Now, we estimate terms in $\Fluid{\varv}{F}$. For the first term on the RHS of \eqref{EQ:varvF}, we have:
	\begin{align}\label{PV:varvF1}
		&\sum\limits_{\substack{\vec{m_1}+\vec{m_2}+\vec{m_3}+\vec{m_4}=(i-i_0,j,k-k_0)\\\abs{\vec{m_3}}\geq 1
		}}\delta^{i_0}\abs{u}^{k_0}\norm{\WDfour^{\vec{m_1}}(\upeta+\ueta)^{k_0}\cdot\WDfour^{\vec{m_2}}\omega^{i_0}\cdot\WDfour^{\vec{m_3}}\speed^2\cdot\WDfour^{\vec{m_4}}\Dfour_3v_4}_{L^2(\suubar)}\\
		\notag\lesssim&\sum\limits_{\substack{\vec{m_1}+\vec{m_2}+\vec{m_3}+\vec{m_4}=(i-i_0,j,k-k_0)\\\abs{\vec{m_1}+\vec{m_2}}> \frac{N-1}{2}
		}}\delta^{i_0}\abs{u}^{k_0}\norm{\WDfour^{\vec{m_1}}\upsi^{k_0}\cdot\WDfour^{\vec{m_2}}\psi^{i_0}}_{L^2(\suubar)}\cdot\norm{\WDfour^{\vec{m_3}}\speed^2}_{L^\infty(\suubar)}\cdot\norm{\WDfour^{\vec{m_4}}\Dfour_3v_4}_{L^\infty(\suubar)}\\
		\notag&+\sum\limits_{\substack{\vec{m_1}+\vec{m_2}+\vec{m_3}+\vec{m_4}=(i-i_0,j,k-k_0)\\\abs{\vec{m_3}}> \frac{N-1}{2}
		}}\delta^{i_0}\abs{u}^{k_0}\norm{\WDfour^{\vec{m_1}}\upsi^{k_0}\cdot\WDfour^{\vec{m_2}}\psi^{i_0}}_{L^\infty(\suubar)}\cdot\norm{\WDfour^{\vec{m_3}}\speed^2}_{L^2(\suubar)}\cdot\norm{\WDfour^{\vec{m_4}}\Dfour_3v_4}_{L^\infty(\suubar)}\\
		\notag&+\sum\limits_{\substack{\vec{m_1}+\vec{m_2}+\vec{m_3}+\vec{m_4}=(i-i_0,j,k-k_0)\\\abs{\vec{m_4}+1}> \frac{N-1}{2}
		}}\delta^{i_0}\abs{u}^{k_0}\norm{\WDfour^{\vec{m_1}}\upsi^{k_0}\cdot\WDfour^{\vec{m_2}}\psi^{i_0}}_{L^\infty(\suubar)}\cdot\norm{\WDfour^{\vec{m_3}}\speed^2}_{L^\infty(\suubar)}\cdot\norm{\WDfour^{\vec{m_4}}\Dfour_3v_4}_{L^2(\suubar)}\\
		\notag\lesssim&\left(\BAsize\frac{\delta a^{\frac{1}{2}}}{\abs{u}}\right)^{i_0+k_0}\cdot\BAfluid^2 b^{\frac{1}{2}}.
	\end{align}
	For the second term on the RHS of \eqref{EQ:varvF}:
	\begin{align}\label{PV:varvF2}
		&\sum\limits_{\substack{\vec{m_1}+\vec{m_2}+\vec{m_3}+\vec{m_4}+\vec{m_5}=(i-i_0,j,k-k_0)\\\vec{m_1},\vec{m_2},\vec{m_4}\neq(i-i_0,j,k-k_0)}}\delta^{i_0}\abs{u}^{k_0}\norm{\WDfour^{\vec{m_1}}(\upeta+\ueta)^{k_0}\cdot\WDfour^{\vec{m_2}}\omega^{i_0}\cdot\WDfour^{\vec{m_3}}\speed^2\cdot\WDfour^{\vec{m_4}}(\uchi,\uomega,\upzeta)\cdot\WDfour^{\vec{m_5}}v_4}_{L^2(\suubar)}\\
		\notag\lesssim&\sum\limits_{\substack{\vec{m_1}+\vec{m_2}+\vec{m_3}+\vec{m_4}+\vec{m_5}=(i-i_0,j,k-k_0)\\\abs{\vec{m_1}+\vec{m_2}}> \frac{N-1}{2}
		}}\delta^{i_0}\abs{u}^{k_0}\norm{\WDfour^{\vec{m_1}}\upsi^{k_0+1}\cdot\WDfour^{\vec{m_2}}\psi^{i_0}}_{L^2(\suubar)}\cdot\norm{\WDfour^{\vec{m_3}}\speed^2}_{L^\infty(\suubar)}\cdot\norm{\WDfour^{\vec{m_4}}v_4}_{L^\infty(\suubar)}\\
		\notag&+\sum\limits_{\substack{\vec{m_1}+\vec{m_2}+\vec{m_3}+\vec{m_4}+\vec{m_5}=(i-i_0,j,k-k_0)\\\abs{\vec{m_3}}> \frac{N-1}{2}
		}}\delta^{i_0}\abs{u}^{k_0}\norm{\WDfour^{\vec{m_1}}\upsi^{k_0+1}\cdot\WDfour^{\vec{m_2}}\psi^{i_0}}_{L^\infty(\suubar)}\cdot\norm{\WDfour^{\vec{m_3}}\speed^2}_{L^2(\suubar)}\cdot\norm{\WDfour^{\vec{m_4}}v_4}_{L^\infty(\suubar)}\\
		\notag&+\sum\limits_{\substack{\vec{m_1}+\vec{m_2}+\vec{m_3}+\vec{m_4}+\vec{m_5}=(i-i_0,j,k-k_0)\\\abs{\vec{m_4}}> \frac{N-1}{2}
		}}\delta^{i_0}\abs{u}^{k_0}\norm{\WDfour^{\vec{m_1}}\upsi^{k_0+1}\cdot\WDfour^{\vec{m_2}}\psi^{i_0}}_{L^\infty(\suubar)}\cdot\norm{\WDfour^{\vec{m_3}}\speed^2}_{L^\infty(\suubar)}\cdot\norm{\WDfour^{\vec{m_4}}v_4}_{L^2(\suubar)}\\
		\notag\lesssim&\left(\BAsize\frac{\delta a^{\frac{1}{2}}}{\abs{u}}\right)^{i_0+k_0+1}\cdot\BAfluid b^{\frac{1}{2}}.
	\end{align}
	We bound the third term of \eqref{EQ:varvF} using Proposition \ref{PR:Weyl} and Proposition \ref{PR:Ricci0}:
	\begin{align}\label{PV:varvF3}
		&\sum\limits_{\vec{m_1}+\vec{m_2}+\vec{m_3}+\vec{m_4}+\vec{m_5}=(i-i_0,j,k-k_0-1)}\delta^{i_0}\abs{u}^{k_0+1}\norm{\WDfour^{\vec{m_1}}(\upeta+\ueta)^{k_0}\cdot\WDfour^{\vec{m_2}}\omega^{i_0}\cdot\WDfour^{\vec{m_3}}\speed^2\cdot\WDfour^{\vec{m_4}}\Riem^{I}_{43A}\cdot\WDfour^{\vec{m_5}}(v_{I})}_{L^2(\suubar)}\\
		\notag\lesssim&\sum\limits_{\vec{m_1}+\vec{m_2}+\vec{m_3}+\vec{m_4}+\vec{m_5}=(i-i_0,j,k-k_0-1)}\delta^{i_0}\abs{u}^{k_0+1}\norm{\WDfour^{\vec{m_1}}(\upeta+\ueta)^{k_0}\cdot\WDfour^{\vec{m_2}}\omega^{i_0}\cdot\WDfour^{\vec{m_3}}\speed^2\cdot\WDfour^{\vec{m_4}}(\Betabar,\FBetabar)\cdot\WDfour^{\vec{m_5}}v_{4}}_{L^2(\suubar)}\\
		\notag&+\sum\limits_{\vec{m_1}+\vec{m_2}+\vec{m_3}+\vec{m_4}+\vec{m_5}=(i-i_0,j,k-k_0-1)}\delta^{i_0}\abs{u}^{k_0+1}\norm{\WDfour^{\vec{m_1}}(\upeta+\ueta)^{k_0}\cdot\WDfour^{\vec{m_2}}\omega^{i_0}\cdot\WDfour^{\vec{m_3}}\speed^2\cdot\WDfour^{\vec{m_4}}(\Rho,\Rhodual,\FRho,\FRhodual,\Scalar)\cdot\WDfour^{\vec{m_5}}(\angv)}_{L^2(\suubar)}\\
		\notag\lesssim&\sum\limits_{\substack{\vec{m_1}+\vec{m_2}+\vec{m_3}+\vec{m_4}+\vec{m_5}=(i-i_0,j,k-k_0-1)\\\abs{\vec{m_1}+\vec{m_2}\vec{m_3}}> \frac{N-1}{2}
		}}\delta^{i_0}\abs{u}^{k_0+1}\norm{\WDfour^{\vec{m_1}}\upsi^{k_0}\cdot\WDfour^{\vec{m_2}}\psi^{i_0}\cdot\WDfour^{\vec{m_3}}\speed^2}_{L^\infty(\suubar)}\\
		\notag&\cdot\left(\norm{\WDfour^{\vec{m_4}}(\Betabar,\FBetabar)}_{L^2(\suubar)}\cdot\norm{\WDfour^{\vec{m_5}}v_4}_{L^\infty(\suubar)}+\norm{\WDfour^{\vec{m_4}}(\Rho,\Rhodual,\FRho,\FRhodual,\Scalar)}_{L^2(\suubar)}\cdot\norm{\WDfour^{\vec{m_5}}(\angv)}_{L^\infty(\suubar)}\right)\\
		\notag&+\sum\limits_{\substack{\vec{m_1}+\vec{m_2}+\vec{m_3}+\vec{m_4}+\vec{m_5}=(i-i_0,j,k-k_0-1)\\\abs{\vec{m_4}}> \frac{N-1}{2}
		}}\delta^{i_0}\abs{u}^{k_0+1}\norm{\WDfour^{\vec{m_1}}\upsi^{k_0}\cdot\WDfour^{\vec{m_2}}\psi^{i_0}\cdot\WDfour^{\vec{m_3}}\speed^2}_{L^\infty(\suubar)}\\
		\notag&\cdot\left(\norm{\WDfour^{\vec{m_4}}(\Betabar,\FBetabar)}_{L^2(\suubar)}\cdot\norm{\WDfour^{\vec{m_5}}v_4}_{L^\infty(\suubar)}+\norm{\WDfour^{\vec{m_4}}(\Rho,\Rhodual,\FRho,\FRhodual,\Scalar)}_{L^2(\suubar)}\cdot\norm{\WDfour^{\vec{m_5}}(\angv)}_{L^\infty(\suubar)}\right)\\
		\notag&+\sum\limits_{\substack{\vec{m_1}+\vec{m_2}+\vec{m_3}+\vec{m_4}+\vec{m_5}=(i-i_0,j,k-k_0-1)\\\abs{\vec{m_5}}> \frac{N-1}{2}
		}}\delta^{i_0}\abs{u}^{k_0+1}\norm{\WDfour^{\vec{m_1}}\upsi^{k_0}\cdot\WDfour^{\vec{m_2}}\psi^{i_0}\cdot\WDfour^{\vec{m_3}}\speed^2}_{L^\infty(\suubar)}\\
		\notag&\cdot\left(\norm{\WDfour^{\vec{m_4}}(\Betabar,\FBetabar)}_{L^\infty(\suubar)}\cdot\norm{\WDfour^{\vec{m_5}}v_4}_{L^2(\suubar)}+\norm{\WDfour^{\vec{m_4}}(\Rho,\Rhodual,\FRho,\FRhodual,\Scalar)}_{L^\infty(\suubar)}\cdot\norm{\WDfour^{\vec{m_5}}(\angv)}_{L^2(\suubar)}\right)\\
		\notag\lesssim&\left(\BAsize\frac{\delta a^{\frac{1}{2}}}{\abs{u}}\right)^{i_0+k_0}\cdot\abs{u}^2\cdot\left[\left(\BAsize\frac{\delta a^{\frac{1}{2}}}{\abs{u}^3}+\BAfluid b^{-\frac{1}{2}}\right)\cdot b^{\frac{1}{2}}+\left(\BAsize\frac{\delta a^{\frac{1}{2}}}{\abs{u}^3}+1\right)\cdot\BAfluid\right].
	\end{align}
	Similar to \eqref{PV:varvF3}, we estimate the forth term of \eqref{EQ:varvF} using Proposition \ref{PR:Weyl} and Proposition \ref{PR:Ricci0}:
	\begin{align}\label{PV:varvF4}
		&\sum\limits_{\vec{m_1}+\vec{m_2}+\vec{m_3}+\vec{m_4}+\vec{m_5}=(i-i_0-1,j,k-k_0)}\delta^{i_0+1}\abs{u}^{k_0}\norm{\WDfour^{\vec{m_1}}(\upeta+\ueta)^{k_0}\cdot\WDfour^{\vec{m_2}}\omega^{i_0}\cdot\WDfour^{\vec{m_3}}\speed^2\cdot\WDfour^{\vec{m_4}}\Riem^{I}_{434}\cdot\WDfour^{\vec{m_5}}(v_I)}_{L^2(\suubar)}\\
		\notag\lesssim&\sum\limits_{\vec{m_1}+\vec{m_2}+\vec{m_3}+\vec{m_4}+\vec{m_5}=(i-i_0-1,j,k-k_0)}\delta^{i_0+1}\abs{u}^{k_0}\norm{\WDfour^{\vec{m_1}}(\upeta+\ueta)^{k_0}\cdot\WDfour^{\vec{m_2}}\omega^{i_0}\cdot\WDfour^{\vec{m_3}}\speed^2\cdot\WDfour^{\vec{m_4}}(\Rho,\FRho,\Scalar)\cdot\WDfour^{\vec{m_5}}v_{4}}_{L^2(\suubar)}\\
		\notag&+\sum\limits_{\vec{m_1}+\vec{m_2}+\vec{m_3}+\vec{m_4}+\vec{m_5}=(i-i_0-1,j,k-k_0)}\delta^{i_0+1}\abs{u}^{k_0}\norm{\WDfour^{\vec{m_1}}(\upeta+\ueta)^{k_0}\cdot\WDfour^{\vec{m_2}}\omega^{i_0}\cdot\WDfour^{\vec{m_3}}\speed^2\cdot\WDfour^{\vec{m_4}}(\BetaA,\FBetaA)\cdot\WDfour^{\vec{m_5}}(\angv)}_{L^2(\suubar)}\\
		\notag\lesssim&\sum\limits_{\substack{\vec{m_1}+\vec{m_2}+\vec{m_3}+\vec{m_4}+\vec{m_5}=(i-i_0-1,j,k-k_0)\\\abs{\vec{m_1}+\vec{m_2}\vec{m_3}}> \frac{N-1}{2}
		}}\delta^{i_0+1}\abs{u}^{k_0}\norm{\WDfour^{\vec{m_1}}\upsi^{k_0}\cdot\WDfour^{\vec{m_2}}\psi^{i_0}\cdot\WDfour^{\vec{m_3}}\speed^2}_{L^\infty(\suubar)}\\
		\notag&\cdot\left(\norm{\WDfour^{\vec{m_4}}(\Rho,\FRho,\Scalar)}_{L^2(\suubar)}\cdot\norm{\WDfour^{\vec{m_5}}v_4}_{L^\infty(\suubar)}+\norm{\WDfour^{\vec{m_4}}(\BetaA,\FBetaA)}_{L^2(\suubar)}\cdot\norm{\WDfour^{\vec{m_5}}(\angv)}_{L^\infty(\suubar)}\right)\\
		\notag&+\sum\limits_{\substack{\vec{m_1}+\vec{m_2}+\vec{m_3}+\vec{m_4}+\vec{m_5}=(i-i_0-1,j,k-k_0)\\\abs{\vec{m_4}}> \frac{N-1}{2}
		}}\delta^{i_0+1}\abs{u}^{k_0}\norm{\WDfour^{\vec{m_1}}\upsi^{k_0}\cdot\WDfour^{\vec{m_2}}\psi^{i_0}\cdot\WDfour^{\vec{m_3}}\speed^2}_{L^\infty(\suubar)}\\
		\notag&\cdot\left(\norm{\WDfour^{\vec{m_4}}(\Rho,\FRho,\Scalar)}_{L^2(\suubar)}\cdot\norm{\WDfour^{\vec{m_5}}v_4}_{L^\infty(\suubar)}+\norm{\WDfour^{\vec{m_4}}(\BetaA,\FBetaA)}_{L^2(\suubar)}\cdot\norm{\WDfour^{\vec{m_5}}(\angv)}_{L^\infty(\suubar)}\right)\\
		\notag&+\sum\limits_{\substack{\vec{m_1}+\vec{m_2}+\vec{m_3}+\vec{m_4}+\vec{m_5}=(i-i_0-1,j,k-k_0)\\\abs{\vec{m_5}}> \frac{N-1}{2}
		}}\delta^{i_0+1}\abs{u}^{k_0}\norm{\WDfour^{\vec{m_1}}\upsi^{k_0}\cdot\WDfour^{\vec{m_2}}\psi^{i_0}\cdot\WDfour^{\vec{m_3}}\speed^2}_{L^\infty(\suubar)}\\
		\notag&\cdot\left(\norm{\WDfour^{\vec{m_4}}(\Rho,\FRho,\Scalar)}_{L^\infty(\suubar)}\cdot\norm{\WDfour^{\vec{m_5}}v_4}_{L^2(\suubar)}+\norm{\WDfour^{\vec{m_4}}(\BetaA,\FBetaA)}_{L^\infty(\suubar)}\cdot\norm{\WDfour^{\vec{m_5}}(\angv)}_{L^2(\suubar)}\right)\\
		\notag\lesssim&\left(\BAsize\frac{\delta a^{\frac{1}{2}}}{\abs{u}}\right)^{i_0+k_0}\cdot\delta\abs{u}\cdot\left[\left(\BAsize\frac{\delta a^{\frac{1}{2}}}{\abs{u}^3}+1\right)\cdot b^{\frac{1}{2}}+\left(\BAsize\frac{ a^{\frac{1}{2}}}{\abs{u}^2}+\BAfluid b^{\frac{1}{2}}\right)\cdot\BAfluid\right].
	\end{align}
	For the fifth and sixth terms of \eqref{EQ:varvF}:
	\begin{align}\label{PV:varvF5}
		&\sum\limits_{\substack{\vec{m_1}+\vec{m_2}+\vec{m_3}+\vec{m_4}=(i,j-j_0,k-k_0)\\\abs{\vec{m_3}}\geq 1
		}}\abs{u}^{j_0+k_0}\norm{\WDfour^{\vec{m_1}}(\upeta+\ueta)^{k_0}\cdot\WDfour^{\vec{m_2}}\uomega^{j_0}\cdot\WDfour^{\vec{m_3}}\speed^2\cdot\WDfour^{\vec{m_4}}\Dfour_4v_3}_{L^2(\suubar)}\\
		\notag\lesssim&\sum\limits_{\substack{\vec{m_1}+\vec{m_2}+\vec{m_3}=(i,j-j_0,k-k_0)\\\abs{\vec{m_1}}> \frac{N-1}{2}
		}}\abs{u}^{j_0+k_0}\norm{\WDfour^{\vec{m_1}}\upsi^{j_0+k_0}}_{L^2(\suubar)}\cdot\norm{\WDfour^{\vec{m_2}}\speed^2}_{L^\infty(\suubar)}\cdot\norm{\WDfour^{\vec{m_3}}\Dfour_4v_3}_{L^\infty(\suubar)}\\
		\notag&+\sum\limits_{\substack{\vec{m_1}+\vec{m_2}+\vec{m_3}=(i,j-j_0,k-k_0)\\\abs{\vec{m_2}}> \frac{N-1}{2}
		}}\abs{u}^{j_0+k_0}\norm{\WDfour^{\vec{m_1}}\upsi^{j_0+k_0}}_{L^\infty(\suubar)}\cdot\norm{\WDfour^{\vec{m_2}}\speed^2}_{L^2(\suubar)}\cdot\norm{\WDfour^{\vec{m_3}}\Dfour_4v_3}_{L^\infty(\suubar)}\\
		\notag&+\sum\limits_{\substack{\vec{m_1}+\vec{m_2}+\vec{m_3}=(i,j-j_0,k-k_0)\\\abs{\vec{m_3}+1}> \frac{N-1}{2}
		}}\abs{u}^{j_0+k_0}\norm{\WDfour^{\vec{m_1}}\upsi^{j_0+k_0}}_{L^\infty(\suubar)}\cdot\norm{\WDfour^{\vec{m_2}}\speed^2}_{L^\infty(\suubar)}\cdot\norm{\WDfour^{\vec{m_3}}\Dfour_4v_3}_{L^2(\suubar)}\\
		\notag\lesssim&\left(\BAsize\frac{\delta a^{\frac{1}{2}}}{\abs{u}}\right)^{i_0+k_0}\cdot\BAfluid^2 b^{-\frac{1}{2}}\cdot\delta^{-1}\cdot \abs{u}.
	\end{align}
	
	\begin{align}\label{PV:varvF6}
		&\sum\limits_{\substack{\vec{m_1}+\vec{m_2}+\vec{m_3}+\vec{m_4}+\vec{m_5}=(i,j-j_0,k-k_0)\\\vec{m_1},\vec{m_2},\vec{m_4}\neq(i,j-j_0,k-k_0)}}\abs{u}^{j_0+k_0}\norm{\WDfour^{\vec{m_1}}(\upeta+\ueta)^{k_0}\cdot\WDfour^{\vec{m_2}}\uomega^{j_0}\cdot\WDfour^{\vec{m_3}}\speed^2\cdot\WDfour^{\vec{m_4}}(\upchi,\omega)\cdot\WDfour^{\vec{m_5}}(v_3)}_{L^2(\suubar)}\\
		\notag\lesssim&\sum\limits_{\substack{\vec{m_1}+\vec{m_2}+\vec{m_3}+\vec{m_4}=(i,j-j_0,k-k_0)\\\abs{\vec{m_1}+\vec{m_2}}> \frac{N-1}{2}
		}}\abs{u}^{j_0+k_0}\norm{\WDfour^{\vec{m_1}}\upsi^{j_0+k_0}\cdot\WDfour^{\vec{m_2}}\psi}_{L^2(\suubar)}\cdot\norm{\WDfour^{\vec{m_3}}\speed^2}_{L^\infty(\suubar)}\cdot\norm{\WDfour^{\vec{m_4}}(v_3)}_{L^\infty(\suubar)}\\
		\notag&+\sum\limits_{\substack{\vec{m_1}+\vec{m_2}+\vec{m_3}+\vec{m_4}=(i,j-j_0,k-k_0)\\\abs{\vec{m_3}}> \frac{N-1}{2}
		}}\abs{u}^{j_0+k_0}\norm{\WDfour^{\vec{m_1}}\upsi^{j_0+k_0}\cdot\WDfour^{\vec{m_2}}\psi}_{L^\infty(\suubar)}\cdot\norm{\WDfour^{\vec{m_3}}\speed^2}_{L^2(\suubar)}\cdot\norm{\WDfour^{\vec{m_4}}(v_3)}_{L^\infty(\suubar)}\\
		\notag&+\sum\limits_{\substack{\vec{m_1}+\vec{m_2}+\vec{m_3}+\vec{m_4}=(i,j-j_0,k-k_0)\\\abs{\vec{m_4}}> \frac{N-1}{2}
		}}\abs{u}^{j_0+k_0}\norm{\WDfour^{\vec{m_1}}\upsi^{j_0+k_0}\cdot\WDfour^{\vec{m_2}}\psi}_{L^\infty(\suubar)}\cdot\norm{\WDfour^{\vec{m_3}}\speed^2}_{L^\infty(\suubar)}\cdot\norm{\WDfour^{\vec{m_4}}(v_3)}_{L^2(\suubar)}\\
		\notag\lesssim&\left(\BAsize\frac{\delta a^{\frac{1}{2}}}{\abs{u}}\right)^{i_0+k_0}\cdot \BAsize a^{\frac{1}{2}}\cdot\BAfluid b^{-\frac{1}{2}}.
	\end{align}
	For the seventh and the eighth terms of \eqref{EQ:varvF}, we take the same approach as in \eqref{PV:varvF3}-\eqref{PV:varvF4}:
	\begin{align}\label{PV:varvF7}
		&\sum\limits_{\vec{m_1}+\vec{m_2}+\vec{m_3}+\vec{m_4}+\vec{m_5}=(i,j-j_0,k-k_0-1)}\abs{u}^{j_0+k_0+1}\norm{\WDfour^{\vec{m_1}}(\upeta+\ueta)^{k_0}\cdot\WDfour^{\vec{m_2}}\omega^{i_0}\cdot\WDfour^{\vec{m_3}}\speed^2\cdot\WDfour^{\vec{m_4}}\Riem^{I}_{34A}\cdot\WDfour^{\vec{m_5}}(v_{I})}_{L^2(\suubar)}\\
		\notag\lesssim&\sum\limits_{\vec{m_1}+\vec{m_2}+\vec{m_3}+\vec{m_4}+\vec{m_5}=(i,j-j_0,k-k_0-1)}\abs{u}^{j_0+k_0+1}\norm{\WDfour^{\vec{m_1}}(\upeta+\ueta)^{k_0}\cdot\WDfour^{\vec{m_2}}\omega^{i_0}\cdot\WDfour^{\vec{m_3}}\speed^2\cdot\WDfour^{\vec{m_4}}(\BetaA,\FBetaA)\cdot\WDfour^{\vec{m_5}}(v_3)}_{L^2(\suubar)}\\
		\notag&+\sum\limits_{\vec{m_1}+\vec{m_2}+\vec{m_3}+\vec{m_4}+\vec{m_5}=(i,j-j_0,k-k_0-1)}\abs{u}^{j_0+k_0+1}\norm{\WDfour^{\vec{m_1}}(\upeta+\ueta)^{k_0}\cdot\WDfour^{\vec{m_2}}\omega^{i_0}\cdot\WDfour^{\vec{m_3}}\speed^2\cdot\WDfour^{\vec{m_4}}(\Rho,\Rhodual,\FRho,\FRhodual,\Scalar)\cdot\WDfour^{\vec{m_5}}(\angv)}_{L^2(\suubar)}\\
		\notag\lesssim&\sum\limits_{\substack{\vec{m_1}+\vec{m_2}+\vec{m_3}+\vec{m_4}+\vec{m_5}=(i,j-j_0,k-k_0-1)\\\abs{\vec{m_1}+\vec{m_2}\vec{m_3}}> \frac{N-1}{2}
		}}\abs{u}^{j_0+k_0+1}\norm{\WDfour^{\vec{m_1}}\upsi^{k_0}\cdot\WDfour^{\vec{m_2}}\psi^{i_0}\cdot\WDfour^{\vec{m_3}}\speed^2}_{L^\infty(\suubar)}\\
		\notag&\cdot\left(\norm{\WDfour^{\vec{m_4}}(\BetaA,\FBetaA)}_{L^2(\suubar)}\cdot\norm{\WDfour^{\vec{m_5}}(v_3)}_{L^\infty(\suubar)}+\norm{\WDfour^{\vec{m_4}}(\Rho,\Rhodual,\FRho,\FRhodual,\Scalar)}_{L^2(\suubar)}\cdot\norm{\WDfour^{\vec{m_5}}(\angv)}_{L^\infty(\suubar)}\right)\\
		\notag&+\sum\limits_{\substack{\vec{m_1}+\vec{m_2}+\vec{m_3}+\vec{m_4}+\vec{m_5}=(i,j-j_0,k-k_0-1)\\\abs{\vec{m_4}}> \frac{N-1}{2}
		}}\abs{u}^{j_0+k_0+1}\norm{\WDfour^{\vec{m_1}}\upsi^{k_0}\cdot\WDfour^{\vec{m_2}}\psi^{i_0}\cdot\WDfour^{\vec{m_3}}\speed^2}_{L^\infty(\suubar)}\\
		\notag&\cdot\left(\norm{\WDfour^{\vec{m_4}}(\BetaA,\FBetaA)}_{L^2(\suubar)}\cdot\norm{\WDfour^{\vec{m_5}}(v_3)}_{L^\infty(\suubar)}+\norm{\WDfour^{\vec{m_4}}(\Rho,\Rhodual,\FRho,\FRhodual,\Scalar)}_{L^2(\suubar)}\cdot\norm{\WDfour^{\vec{m_5}}(\angv)}_{L^\infty(\suubar)}\right)\\
		\notag&+\sum\limits_{\substack{\vec{m_1}+\vec{m_2}+\vec{m_3}+\vec{m_4}+\vec{m_5}=(i,j-j_0,k-k_0-1)\\\abs{\vec{m_5}}> \frac{N-1}{2}
		}}\abs{u}^{j_0+k_0+1}\norm{\WDfour^{\vec{m_1}}\upsi^{k_0}\cdot\WDfour^{\vec{m_2}}\psi^{i_0}\cdot\WDfour^{\vec{m_3}}\speed^2}_{L^\infty(\suubar)}\\
		\notag&\cdot\left(\norm{\WDfour^{\vec{m_4}}(\BetaA,\FBetaA)}_{L^\infty(\suubar)}\cdot\norm{\WDfour^{\vec{m_5}}(v_3)}_{L^2(\suubar)}+\norm{\WDfour^{\vec{m_4}}(\Rho,\Rhodual,\FRho,\FRhodual,\Scalar)}_{L^\infty(\suubar)}\cdot\norm{\WDfour^{\vec{m_5}}(\angv)}_{L^2(\suubar)}\right)\\
		\notag\lesssim&\left(\BAsize\frac{\delta a^{\frac{1}{2}}}{\abs{u}}\right)^{i_0+k_0}\cdot\abs{u}^2\cdot\left[\left(\BAsize\frac{a^{\frac{1}{2}}}{\abs{u}^2}+b^{\frac{1}{2}}\right)\cdot b^{-\frac{1}{2}}+\left(\BAsize\frac{\delta a^{\frac{1}{2}}}{\abs{u}^3}+1\right)\cdot\BAfluid\right].
	\end{align}
	
	\begin{align}\label{PV:varvF8}
		&\sum\limits_{\vec{m_1}+\vec{m_2}+\vec{m_3}+\vec{m_4}+\vec{m_5}=(i,j-j_0-1,k-k_0)}\abs{u}^{j_0+k_0+1}\norm{\WDfour^{\vec{m_1}}(\upeta+\ueta)^{k_0}\cdot\Dfour^{\vec{m_2}}\omega^{i_0}\cdot\WDfour^{\vec{m_3}}\speed^2\cdot\WDfour^{\vec{m_4}}\Riem^{I}_{334}\cdot\WDfour^{\vec{m_5}}(v_I)}_{L^2(\suubar)}\\
		\notag\lesssim&\sum\limits_{\vec{m_1}+\vec{m_2}+\vec{m_3}+\vec{m_4}+\vec{m_5}=(i,j-j_0-1,k-k_0)}\abs{u}^{j_0+k_0+1}\norm{\WDfour^{\vec{m_1}}(\upeta+\ueta)^{k_0}\cdot\WDfour^{\vec{m_2}}\omega^{i_0}\cdot\WDfour^{\vec{m_3}}\speed^2\cdot\WDfour^{\vec{m_4}}(\Rho,\FRho,\Scalar)\cdot\WDfour^{\vec{m_5}}(v_3)}_{L^2(\suubar)}\\
		\notag&+\sum\limits_{\vec{m_1}+\vec{m_2}+\vec{m_3}+\vec{m_4}+\vec{m_5}=(i,j-j_0-1,k-k_0)}\abs{u}^{j_0+k_0+1}\norm{\WDfour^{\vec{m_1}}(\upeta+\ueta)^{k_0}\cdot\WDfour^{\vec{m_2}}\omega^{i_0}\cdot\WDfour^{\vec{m_3}}\speed^2\cdot\WDfour^{\vec{m_4}}(\Betabar,\FBetabar)\cdot\WDfour^{\vec{m_5}}(\angv)}_{L^2(\suubar)}\\
		\notag\lesssim&\sum\limits_{\substack{\vec{m_1}+\vec{m_2}+\vec{m_3}+\vec{m_4}+\vec{m_5}=(i,j-j_0-1,k-k_0)\\\abs{\vec{m_1}+\vec{m_2}\vec{m_3}}> \frac{N-1}{2}
		}}\abs{u}^{j_0+k_0+1}\norm{\WDfour^{\vec{m_1}}\upsi^{k_0}\cdot\WDfour^{\vec{m_2}}\psi^{i_0}\cdot\WDfour^{\vec{m_3}}\speed^2}_{L^\infty(\suubar)}\\
		\notag&\cdot\left(\norm{\WDfour^{\vec{m_4}}(\Rho,\FRho,\Scalar)}_{L^2(\suubar)}\cdot\norm{\WDfour^{\vec{m_5}}(v_3)}_{L^\infty(\suubar)}+\norm{\WDfour^{\vec{m_4}}(\Betabar,\FBetabar)}_{L^2(\suubar)}\cdot\norm{\WDfour^{\vec{m_5}}(\angv)}_{L^\infty(\suubar)}\right)\\
		\notag&+\sum\limits_{\substack{\vec{m_1}+\vec{m_2}+\vec{m_3}+\vec{m_4}+\vec{m_5}=(i,j-j_0-1,k-k_0)\\\abs{\vec{m_4}}> \frac{N-1}{2}
		}}\abs{u}^{j_0+k_0+1}\norm{\WDfour^{\vec{m_1}}\upsi^{k_0}\cdot\WDfour^{\vec{m_2}}\psi^{i_0}\cdot\WDfour^{\vec{m_3}}\speed^2}_{L^\infty(\suubar)}\\
		\notag&\cdot\left(\norm{\WDfour^{\vec{m_4}}(\Rho,\FRho,\Scalar)}_{L^2(\suubar)}\cdot\norm{\WDfour^{\vec{m_5}}(v_3)}_{L^\infty(\suubar)}+\norm{\WDfour^{\vec{m_4}}(\Betabar,\FBetabar)}_{L^2(\suubar)}\cdot\norm{\WDfour^{\vec{m_5}}(\angv)}_{L^\infty(\suubar)}\right)\\
		\notag&+\sum\limits_{\substack{\vec{m_1}+\vec{m_2}+\vec{m_3}+\vec{m_4}+\vec{m_5}=(i,j-j_0-1,k-k_0)\\\abs{\vec{m_5}}> \frac{N-1}{2}
		}}\abs{u}^{j_0+k_0+1}\norm{\WDfour^{\vec{m_1}}\upsi^{k_0}\cdot\WDfour^{\vec{m_2}}\psi^{i_0}\cdot\WDfour^{\vec{m_3}}\speed^2}_{L^\infty(\suubar)}\\
		\notag&\cdot\left(\norm{\WDfour^{\vec{m_4}}(\Rho,\FRho,\Scalar)}_{L^\infty(\suubar)}\cdot\norm{\WDfour^{\vec{m_5}}(v_3)}_{L^2(\suubar)}+\norm{\WDfour^{\vec{m_4}}(\Betabar,\FBetabar)}_{L^\infty(\suubar)}\cdot\norm{\WDfour^{\vec{m_5}}(\angv)}_{L^2(\suubar)}\right)\\
		\notag\lesssim&\left(\BAsize\frac{\delta a^{\frac{1}{2}}}{\abs{u}}\right)^{i_0+k_0}\cdot\abs{u}^2\cdot\left[\left(\BAsize\frac{\delta a^{\frac{1}{2}}}{\abs{u}^3}+1\right)\cdot b^{-\frac{1}{2}}+\left(\BAsize\frac{ \delta a^{\frac{1}{2}}}{\abs{u}^3}+\BAfluid b^{-\frac{1}{2}}\right)\cdot\BAfluid\right].
	\end{align}
	For the last four terms of \eqref{EQ:varvF}, we use the similar approach as in \eqref{PV:varvF1}-\eqref{PV:varvF4}:
	
	\begin{align}\label{PV:varvF9}
		&\sum\limits_{\substack{\vec{m_1}+\vec{m_2}+\vec{m_3}+\vec{m_4}=(i-i_0,j-j_0,k)\\\abs{\vec{m_3}}\geq 1
		}}\delta^{i_0}\abs{u}^{j_0}\norm{\WDfour^{\vec{m_1}}\upchi^{i_0}\cdot\WDfour^{\vec{m_2}}\uchi^{j_0}\cdot\WDfour^{\vec{m_3}}\speed^2\cdot\WDfour^{\vec{m_4}}\Dfour_Bv_B}_{L^2(\suubar)}\\
		\notag\lesssim&\sum\limits_{\substack{\vec{m_1}+\vec{m_2}+\vec{m_3}+\vec{m_4}=(i-i_0,j-j_0,k)\\\abs{\vec{m_1}+\vec{m_2}}> \frac{N-1}{2}
		}}\delta^{i_0}\abs{u}^{j_0}\norm{\WDfour^{\vec{m_1}}\upsi^{k_0}\cdot\WDfour^{\vec{m_2}}\psi^{i_0}}_{L^2(\suubar)}\cdot\norm{\WDfour^{\vec{m_3}}\speed^2}_{L^\infty(\suubar)}\cdot\norm{\WDfour^{\vec{m_4}}\Dfour_Bv_B}_{L^\infty(\suubar)}\\
		\notag&+\sum\limits_{\substack{\vec{m_1}+\vec{m_2}+\vec{m_3}+\vec{m_4}=(i-i_0,j-j_0,k)\\\abs{\vec{m_3}}> \frac{N-1}{2}
		}}\delta^{i_0}\abs{u}^{j_0}\norm{\WDfour^{\vec{m_1}}\upsi^{k_0}\cdot\WDfour^{\vec{m_2}}\psi^{i_0}}_{L^\infty(\suubar)}\cdot\norm{\WDfour^{\vec{m_3}}\speed^2}_{L^2(\suubar)}\cdot\norm{\WDfour^{\vec{m_4}}\Dfour_Bv_B}_{L^\infty(\suubar)}\\
		\notag&+\sum\limits_{\substack{\vec{m_1}+\vec{m_2}+\vec{m_3}+\vec{m_4}=(i-i_0,j-j_0,k)\\\abs{\vec{m_4}+1}> \frac{N-1}{2}
		}}\delta^{i_0}\abs{u}^{j_0}\norm{\WDfour^{\vec{m_1}}\upsi^{k_0}\cdot\WDfour^{\vec{m_2}}\psi^{i_0}}_{L^\infty(\suubar)}\cdot\norm{\WDfour^{\vec{m_3}}\speed^2}_{L^\infty(\suubar)}\cdot\norm{\WDfour^{\vec{m_4}}\Dfour_Bv_B}_{L^2(\suubar)}\\
		\notag\lesssim&\left(\BAsize\frac{\delta a^{\frac{1}{2}}}{\abs{u}}\right)^{i_0+k_0}\cdot\BAfluid^2 b^{\frac{1}{2}}.
	\end{align}
	
	\begin{align}\label{PV:varvF10}
		&\sum\limits_{\substack{\vec{m_1}+\vec{m_2}+\vec{m_3}+\vec{m_4}+\vec{m_5}=(i-i_0,j-j_0,k)\\\vec{m_1},\vec{m_2},\vec{m_4}\neq(i-i_0,j-j_0,k)}}\delta^{i_0}\abs{u}^{j_0}\norm{\WDfour^{\vec{m_1}}\upchi^{i_0}\cdot\WDfour^{\vec{m_2}}\uchi^{j_0}\cdot\WDfour^{\vec{m_3}}\speed^2\cdot\WDfour^{\vec{m_4}}(\upeta_B+\ueta_B)\cdot\WDfour^{\vec{m_5}}(v_B)}_{L^2(\suubar)}\\
		\notag\lesssim&\sum\limits_{\substack{\vec{m_1}+\vec{m_2}+\vec{m_3}+\vec{m_4}+\vec{m_5}=(i-i_0,j-j_0,k)\\\abs{\vec{m_1}+\vec{m_2}}> \frac{N-1}{2}
		}}\delta^{i_0}\abs{u}^{j_0}\norm{\WDfour^{\vec{m_1}}\upsi^{k_0+1}\cdot\WDfour^{\vec{m_2}}\psi^{i_0}}_{L^2(\suubar)}\cdot\norm{\WDfour^{\vec{m_3}}\speed^2}_{L^\infty(\suubar)}\cdot\norm{\WDfour^{\vec{m_4}}(v_B)}_{L^\infty(\suubar)}\\
		\notag&+\sum\limits_{\substack{\vec{m_1}+\vec{m_2}+\vec{m_3}+\vec{m_4}+\vec{m_5}=(i-i_0,j-j_0,k)\\\abs{\vec{m_3}}> \frac{N-1}{2}
		}}\delta^{i_0}\abs{u}^{j_0}\norm{\WDfour^{\vec{m_1}}\upsi^{k_0+1}\cdot\WDfour^{\vec{m_2}}\psi^{i_0}}_{L^\infty(\suubar)}\cdot\norm{\WDfour^{\vec{m_3}}\speed^2}_{L^2(\suubar)}\cdot\norm{\WDfour^{\vec{m_4}}(v_B)}_{L^\infty(\suubar)}\\
		\notag&+\sum\limits_{\substack{\vec{m_1}+\vec{m_2}+\vec{m_3}+\vec{m_4}+\vec{m_5}=(i-i_0,j-j_0,k)\\\abs{\vec{m_4}}> \frac{N-1}{2}
		}}\delta^{i_0}\abs{u}^{j_0}\norm{\WDfour^{\vec{m_1}}\upsi^{k_0+1}\cdot\WDfour^{\vec{m_2}}\psi^{i_0}}_{L^\infty(\suubar)}\cdot\norm{\WDfour^{\vec{m_3}}\speed^2}_{L^\infty(\suubar)}\cdot\norm{\WDfour^{\vec{m_4}}(v_B)}_{L^2(\suubar)}\\
		\notag\lesssim&\left(\BAsize\frac{\delta a^{\frac{1}{2}}}{\abs{u}}\right)^{i_0+k_0+1}\cdot\BAfluid .
	\end{align}
	We bound the 11-12th term of \eqref{EQ:varvF} using Proposition \ref{PR:Weyl} and Proposition \ref{PR:Ricci0}:
	\begin{align}\label{PV:varvF11}
		&\sum\limits_{\vec{m_1}+\vec{m_2}+\vec{m_3}+\vec{m_4}+\vec{m_5}=(i-i_0,j-j_0-1,k)}\delta^{i_0}\abs{u}^{j_0+1}\norm{\WDfour^{\vec{m_1}}\upchi^{i_0}\cdot\WDfour^{\vec{m_2}}\uchi^{ j_0}\cdot\WDfour^{\vec{m_3}}\speed^2\cdot\WDfour^{\vec{m_4}}\Riem^{I}_{B3B}\cdot\WDfour^{\vec{m_5}}(v_I)}_{L^2(\suubar)}\\
		\notag\lesssim&\sum\limits_{\vec{m_1}+\vec{m_2}+\vec{m_3}+\vec{m_4}+\vec{m_5}=(i-i_0,j-j_0-1,k)}\delta^{i_0}\abs{u}^{j_0+1}\norm{\WDfour^{\vec{m_1}}(\upeta+\ueta)^{k_0}\cdot\WDfour^{\vec{m_2}}\omega^{i_0}\cdot\WDfour^{\vec{m_3}}\speed^2\cdot\WDfour^{\vec{m_4}}(\Alphabar,\FAlphabar)\cdot\WDfour^{\vec{m_5}}v_{4}}_{L^2(\suubar)}\\
		\notag&+\sum\limits_{\vec{m_1}+\vec{m_2}+\vec{m_3}+\vec{m_4}+\vec{m_5}=(i-i_0,j-j_0-1,k)}\delta^{i_0}\abs{u}^{j_0+1}\norm{\WDfour^{\vec{m_1}}(\upeta+\ueta)^{k_0}\cdot\WDfour^{\vec{m_2}}\omega^{i_0}\cdot\WDfour^{\vec{m_3}}\speed^2\cdot\WDfour^{\vec{m_4}}(\Rho,\Rhodual,\FRho,\Scalar)\cdot\WDfour^{\vec{m_5}}(v_3)}_{L^2(\suubar)}\\
		\notag&+\sum\limits_{\vec{m_1}+\vec{m_2}+\vec{m_3}+\vec{m_4}+\vec{m_5}=(i-i_0,j-j_0-1,k)}\delta^{i_0}\abs{u}^{j_0+1}\norm{\WDfour^{\vec{m_1}}(\upeta+\ueta)^{k_0}\cdot\WDfour^{\vec{m_2}}\omega^{i_0}\cdot\WDfour^{\vec{m_3}}\speed^2\cdot\WDfour^{\vec{m_4}}(\Betabar,\FBetabar)\cdot\WDfour^{\vec{m_5}}(\angv)}_{L^2(\suubar)}\\
		\notag\lesssim&\left(\BAsize\frac{\delta a^{\frac{1}{2}}}{\abs{u}}\right)^{i_0+k_0}\cdot\abs{u}^2\cdot\left[\left(\BAsize\frac{\delta^2 a}{\abs{u}^3}+b^{-1}\right)\cdot b^{\frac{1}{2}}+\left(\BAsize\frac{\delta a^{\frac{1}{2}}}{\abs{u}^3}+1\right)\cdot b^{-\frac{1}{2}}+\left(\BAsize\frac{\delta a^{\frac{1}{2}}}{\abs{u}^3}+\BAfluid b^{-\frac{1}{2}}\right)\cdot \BAfluid\right].
	\end{align}
	
	\begin{align}\label{PV:varvF12}
		&\sum\limits_{\vec{m_1}+\vec{m_2}+\vec{m_3}+\vec{m_4}+\vec{m_5}=(i-i_0-1,j-j_0,k)}\delta^{i_0+1}\abs{u}^{j_0}\norm{\WDfour^{\vec{m_1}}\upchi^{i_0}\cdot\WDfour^{\vec{m_2}}\uchi^{ j_0}\cdot\WDfour^{\vec{m_3}}\speed^2\cdot\WDfour^{\vec{m_4}}\Riem^{I}_{B4B}\cdot\WDfour^{\vec{m_5}}(v_I)}_{L^2(\suubar)}\\
		\notag\lesssim&\sum\limits_{\vec{m_1}+\vec{m_2}+\vec{m_3}+\vec{m_4}+\vec{m_5}=(i-i_0-1,j-j_0,k)}\delta^{i_0+1}\abs{u}^{j_0}\norm{\WDfour^{\vec{m_1}}(\upeta+\ueta)^{k_0}\cdot\WDfour^{\vec{m_2}}\omega^{i_0}\cdot\WDfour^{\vec{m_3}}\speed^2\cdot\WDfour^{\vec{m_4}}(\Rho,\FRho,\Scalar)\cdot\WDfour^{\vec{m_5}}v_{4}}_{L^2(\suubar)}\\
		\notag&+\sum\limits_{\vec{m_1}+\vec{m_2}+\vec{m_3}+\vec{m_4}+\vec{m_5}=(i-i_0-1,j-j_0,k)}\delta^{i_0+1}\abs{u}^{j_0}\norm{\WDfour^{\vec{m_1}}(\upeta+\ueta)^{k_0}\cdot\WDfour^{\vec{m_2}}\omega^{i_0}\cdot\WDfour^{\vec{m_3}}\speed^2\cdot\WDfour^{\vec{m_4}}(\AlphaAB,\FAlphaAB)\cdot\WDfour^{\vec{m_5}}(v_3)}_{L^2(\suubar)}\\
		\notag&+\sum\limits_{\vec{m_1}+\vec{m_2}+\vec{m_3}+\vec{m_4}+\vec{m_5}=(i-i_0-1,j-j_0,k)}\delta^{i_0+1}\abs{u}^{j_0}\norm{\WDfour^{\vec{m_1}}(\upeta+\ueta)^{k_0}\cdot\WDfour^{\vec{m_2}}\omega^{i_0}\cdot\WDfour^{\vec{m_3}}\speed^2\cdot\WDfour^{\vec{m_4}}(\BetaA,\FBetaA)\cdot\WDfour^{\vec{m_5}}(\angv)}_{L^2(\suubar)}\\
		\notag\lesssim&\left(\BAsize\frac{\delta a^{\frac{1}{2}}}{\abs{u}}\right)^{i_0+k_0}\cdot\delta\abs{u}\cdot\left[\left(\BAsize\frac{\delta a^{\frac{1}{2}}}{\abs{u}^3}+1\right)\cdot b^{\frac{1}{2}}+\left(\BAsize\frac{\delta^{-1} a^{\frac{1}{2}}}{\abs{u}}+b\right)\cdot b^{-\frac{1}{2}}+\left(\BAsize\frac{ a^{\frac{1}{2}}}{\abs{u}^2}+\BAfluid b^{\frac{1}{2}}\right)\cdot \BAfluid\right].
	\end{align}
	Combining \eqref{PV:varvF1}-\eqref{PV:varvF12}, we obtain the estimate for $\Fluid{\varv}{F}$:
	\begin{align}\label{PV:varvF}
		\norm{\Fluid{\varv}{F}}_{L^2(\suubar)}\lesssim&~\BAfluid^2\delta^{-1}b^{-\frac{1}{2}}\abs{u}+\BAsize\BAfluid a^{\frac{1}{2}} b^{-\frac{1}{2}}.
	\end{align}
	Then, the desired result \eqref{ES:EstimateFGH1} follows from \eqref{PV:varhF} and \eqref{PV:varvF}.
\end{proof}
\begin{proof}[Proof of \eqref{ES:EstimateFGH2}]
	We take the similar approach as in the proof of estimates \eqref{ES:EstimateFGH1}, where in below, we omit repetitive details.
	First of all,
	As in Proposition \ref{PR:EQVariations}, $\mathfrak{G}_{(i,j,k)}:=\Fluid{\varh}{G}+\Fluid{\varv}{G}+\Fluid{\vars}{G}$. We start with estimating terms in $\Fluid{\varh}{G}_3$. For the first and seventh terms on the RHS of \eqref{EQ:varhG}, we have:
	\begin{align}\label{PV:varhG1}
		&\sum\limits_{\substack{\vec{m_1}+\vec{m_2}+\vec{m_3}+\vec{m_4}=(i-i_0,j,k-k_0)\\\abs{\vec{m_3}}\geq1
		}}\delta^{i_0}\abs{u}^{k_0}\norm{\WDfour^{\vec{m_1}}(\upeta+\ueta)^{k_0}\cdot\WDfour^{\vec{m_2}}\omega^{i_0}\cdot\left(1,\WDfour^{\vec{m_3}}(v_4v_3)\right)\cdot\WDfour^{\vec{m_4}}\Dfour_3\lnenth}_{L^2(\suubar)}\\
		\notag\lesssim&\left(\BAsize\frac{\delta a^{\frac{1}{2}}}{\abs{u}}\right)^{i_0+k_0}\cdot\BAfluid.
	\end{align}
	For the second and eighth terms on the RHS of \eqref{EQ:varhG}:
	\begin{align}\label{PV:varhG2}
		&\sum\limits_{\substack{\vec{m_1}+\vec{m_2}+\vec{m_3}+\vec{m_4}+\vec{m_5}=(i-i_0,j,k-k_0)\\\vec{m_1},\vec{m_2},\vec{m_4}\neq(i-i_0,j,k-k_0)}}\delta^{i_0}\abs{u}^{k_0}\norm{\WDfour^{\vec{m_1}}(\upeta+\ueta)^{k_0}\cdot\WDfour^{\vec{m_2}}\omega^{i_0}\cdot\left(1,\WDfour^{\vec{m_3}}(v_4v_3)\right)\cdot\WDfour^{\vec{m_4}}(\uchi,\uomega,\upzeta)\cdot\WDfour^{\vec{m_5}}\lnenth}_{L^2(\suubar)}\\
		\notag\lesssim&\left(\BAsize\frac{\delta a^{\frac{1}{2}}}{\abs{u}}\right)^{i_0+k_0+1}\cdot\BAfluid.
	\end{align}
	The 9-10th terms of \eqref{EQ:varhG} obtain the following estimates:
	\begin{align}\label{PV:varhG3}
		&\sum\limits_{\substack{\vec{m_1}+\vec{m_2}+\vec{m_3}+\vec{m_4}=(i,j-j_0,k-k_0)\\\abs{\vec{m_3}}\geq 1
		}}\abs{u}^{j_0+k_0}\norm{\WDfour^{\vec{m_1}}(\upeta+\ueta)^{k_0}\cdot\WDfour^{\vec{m_2}}\uomega^{j_0}\cdot\WDfour^{\vec{m_3}}(v_3v_3)\cdot\WDfour^{\vec{m_4}}\Dfour_4\lnenth}_{L^2(\suubar)}\\
		\notag\lesssim&\left(\BAsize\frac{\delta a^{\frac{1}{2}}}{\abs{u}}\right)^{j_0+k_0}\cdot\BAfluid^2 b^{-1}\cdot\delta^{-1}\cdot\abs{u}.
	\end{align}
	\begin{align}\label{PV:varhG4}
		&\sum\limits_{\substack{\vec{m_1}+\vec{m_2}+\vec{m_3}+\vec{m_4}+\vec{m_5}=(i,j-j_0,k-k_0)\\\vec{m_1},\vec{m_2},\vec{m_4}\neq(i,j-j_0,k-k_0)}}\abs{u}^{j_0+k_0}\norm{\WDfour^{\vec{m_1}}(\upeta+\ueta)^{k_0}\cdot\WDfour^{\vec{m_2}}\uomega^{j_0}\cdot\WDfour^{\vec{m_3}}(v_3v_3)\cdot\WDfour^{\vec{m_4}}(\upchi,\omega)\cdot\WDfour^{\vec{m_5}}\lnenth}_{L^2(\suubar)}\\
		\notag\lesssim&\left(\BAsize\frac{\delta a^{\frac{1}{2}}}{\abs{u}}\right)^{j_0+k_0}\cdot \BAsize a^{\frac{1}{2}}\cdot\BAfluid b^{-1}.
	\end{align}
	For the last two terms in $\Fluid{\varh}{G}$, we have:
	\begin{align}\label{PV:varhG5}
		&\sum\limits_{\substack{\vec{m_1}+\vec{m_2}+\vec{m_3}+\vec{m_4}+\vec{m_5}=(i-i_0,j-j_0,k)\\\abs{\vec{m_3}+\vec{m_4}}\geq 1
		}}\delta^{i_0}\abs{u}^{j_0}\norm{\WDfour^{\vec{m_1}}\upchi^{i_0}\cdot\WDfour^{\vec{m_2}}\uchi^{j_0}\cdot\WDfour^{\vec{m_3}} (v_B)\cdot\WDfour^{\vec{m_4}}(v_3)\cdot\WDfour^{\vec{m_5}}\Dfour_B\lnenth}_{L^2(\suubar)}\\
		\notag\lesssim&\left(\BAsize\frac{\delta a^{\frac{1}{2}}}{\abs{u}}\right)^{i_0+j_0}\cdot\BAfluid^2 b^{-\frac{1}{2}} .
	\end{align}
	\begin{align}\label{PV:varhG6}
		&\sum\limits_{\substack{\vec{m_1}+\vec{m_2}+\vec{m_3}+\vec{m_4}+\vec{m_5}+\vec{m_6}=(i-i_0,j-j_0,k)\\\vec{m_1},\vec{m_2},\vec{m_5}\neq(i-i_0,j-j_0,k)}}\delta^{i_0}\abs{u}^{j_0}\norm{\WDfour^{\vec{m_1}}\upchi^{i_0}\cdot\WDfour^{\vec{m_2}}\uchi^{j_0}\cdot\WDfour^{\vec{m_3}} (v_B)\cdot\WDfour^{\vec{m_4}}(v_3)\cdot\WDfour^{\vec{m_5}}(\upeta_B+\ueta_B)\cdot\WDfour^{\vec{m_6}}\lnenth}_{L^2(\suubar)}\\
		\notag\lesssim&\left(\BAsize\frac{\delta a^{\frac{1}{2}}}{\abs{u}}\right)^{i_0+j_0+1}\cdot\BAfluid b^{-\frac{1}{2}}.
	\end{align}
	Combining \eqref{PV:varhG1}-\eqref{PV:varhG6}, we obtain the estimate for $\Fluid{\varh}{G}_3$:
	\begin{align}\label{PV:varhG}
		\norm{\Fluid{\varh}{G}_3}_{L^2(\suubar)}\lesssim&~\BAfluid^2\delta^{-1}b^{-1}\abs{u}+\BAsize\BAfluid a^{\frac{1}{2}} b^{-1}.
	\end{align}
	
	Now, we estimate terms in $\Fluid{\varv}{G}$. For the first term on the RHS of \eqref{EQ:varvG}, we have:
	\begin{align}\label{PV:varvG1}
		&\sum\limits_{\substack{\vec{m_1}+\vec{m_2}+\vec{m_3}+\vec{m_4}=(i-i_0,j,k-k_0)\\\abs{\vec{m_3}}\geq 1
		}}\delta^{i_0}\abs{u}^{k_0}\norm{\WDfour^{\vec{m_1}}(\upeta+\ueta)^{k_0}\cdot\WDfour^{\vec{m_2}}\omega^{i_0}\cdot\WDfour^{\vec{m_3}}v_4\cdot\WDfour^{\vec{m_4}}\Dfour_3v_3}_{L^2(\suubar)}\\
		\notag\lesssim&\left(\BAsize\frac{\delta a^{\frac{1}{2}}}{\abs{u}}\right)^{i_0+k_0}\cdot\BAfluid^2.
	\end{align}
	For the second term on the RHS of \eqref{EQ:varvG}:
	\begin{align}\label{PV:varvG2}
		&\sum\limits_{\substack{\vec{m_1}+\vec{m_2}+\vec{m_3}+\vec{m_4}+\vec{m_5}=(i-i_0,j,k-k_0)\\\vec{m_1},\vec{m_2},\vec{m_4}\neq(i-i_0,j,k-k_0)}}\delta^{i_0}\abs{u}^{k_0}\norm{\WDfour^{\vec{m_1}}(\upeta+\ueta)^{k_0}\cdot\WDfour^{\vec{m_2}}\omega^{i_0}\cdot\WDfour^{\vec{m_3}}v_4\cdot\WDfour^{\vec{m_4}}(\uchi,\uomega,\upzeta)\cdot\WDfour^{\vec{m_5}}(v_3)}_{L^2(\suubar)}\\
		\notag\lesssim&\left(\BAsize\frac{\delta a^{\frac{1}{2}}}{\abs{u}}\right)^{i_0+k_0+1}\cdot\BAfluid.
	\end{align}
	We bound the third term of \eqref{EQ:varvG} using Proposition \ref{PR:Weyl} and Proposition \ref{PR:Ricci0}:
	\begin{align}\label{PV:varvG3}
		&\sum\limits_{\vec{m_1}+\vec{m_2}+\vec{m_3}+\vec{m_4}+\vec{m_5}=(i-i_0,j,k-k_0-1)}\delta^{i_0}\abs{u}^{k_0+1}\norm{\WDfour^{\vec{m_1}}(\upeta+\ueta)^{k_0}\cdot\WDfour^{\vec{m_2}}\omega^{i_0}\cdot\WDfour^{\vec{m_3}}v_4\cdot\WDfour^{\vec{m_4}}\Riem^{J}_{33A}\cdot\WDfour^{\vec{m_5}}(v_{J})}_{L^2(\suubar)}\\
		\notag\lesssim&\sum\limits_{\vec{m_1}+\vec{m_2}+\vec{m_3}+\vec{m_4}+\vec{m_5}=(i-i_0,j,k-k_0-1)}\delta^{i_0}\abs{u}^{k_0+1}\norm{\WDfour^{\vec{m_1}}(\upeta+\ueta)^{k_0}\cdot\WDfour^{\vec{m_2}}\omega^{i_0}\cdot\WDfour^{\vec{m_3}}v_4\cdot\WDfour^{\vec{m_4}}(\Betabar,\FBetabar)\cdot\WDfour^{\vec{m_5}}(v_3)}_{L^2(\suubar)}\\
		\notag&+\sum\limits_{\vec{m_1}+\vec{m_2}+\vec{m_3}+\vec{m_4}+\vec{m_5}=(i-i_0,j,k-k_0-1)}\delta^{i_0}\abs{u}^{k_0+1}\norm{\WDfour^{\vec{m_1}}(\upeta+\ueta)^{k_0}\cdot\WDfour^{\vec{m_2}}\omega^{i_0}\cdot\WDfour^{\vec{m_3}}v_4\cdot\WDfour^{\vec{m_4}}(\Alphabar,\FAlphabar)\cdot\WDfour^{\vec{m_5}}(\angv)}_{L^2(\suubar)}\\
		\notag\lesssim&\left(\BAsize\frac{\delta a^{\frac{1}{2}}}{\abs{u}}\right)^{i_0+k_0}\cdot b^{\frac{1}{2}}\abs{u}^2\cdot\left[\left(\BAsize\frac{\delta a^{\frac{1}{2}}}{\abs{u}^3}+\BAfluid b^{-\frac{1}{2}}\right)\cdot b^{-\frac{1}{2}}+\left(\BAsize\frac{\delta^2 a}{\abs{u}^3}+b^{-1}\right)\cdot\BAfluid\right].
	\end{align}
	Similar to \eqref{PV:varvG3}, we estimate the forth term of \eqref{EQ:varvG} using Proposition \ref{PR:Weyl} and Proposition \ref{PR:Ricci0}:
	\begin{align}\label{PV:varvG4}
		&\sum\limits_{\vec{m_1}+\vec{m_2}+\vec{m_3}+\vec{m_4}+\vec{m_5}=(i-i_0-1,j,k-k_0)}\delta^{i_0+1}\abs{u}^{k_0}\norm{\WDfour^{\vec{m_1}}(\upeta+\ueta)^{k_0}\cdot\WDfour^{\vec{m_2}}\omega^{i_0}\cdot\WDfour^{\vec{m_3}}v_4\cdot\WDfour^{\vec{m_4}}\Riem^{J}_{334}\cdot\WDfour^{\vec{m_5}}(v_J)}_{L^2(\suubar)}\\
		\notag\lesssim&\sum\limits_{\vec{m_1}+\vec{m_2}+\vec{m_3}+\vec{m_4}+\vec{m_5}=(i-i_0-1,j,k-k_0)}\delta^{i_0+1}\abs{u}^{k_0}\norm{\WDfour^{\vec{m_1}}(\upeta+\ueta)^{k_0}\cdot\WDfour^{\vec{m_2}}\omega^{i_0}\cdot\WDfour^{\vec{m_3}}v_4\cdot\WDfour^{\vec{m_4}}(\Rho,\FRho,\Scalar)\cdot\WDfour^{\vec{m_5}}(v_3)}_{L^2(\suubar)}\\
		\notag&+\sum\limits_{\vec{m_1}+\vec{m_2}+\vec{m_3}+\vec{m_4}+\vec{m_5}=(i-i_0-1,j,k-k_0)}\delta^{i_0+1}\abs{u}^{k_0}\norm{\WDfour^{\vec{m_1}}(\upeta+\ueta)^{k_0}\cdot\WDfour^{\vec{m_2}}\omega^{i_0}\cdot\WDfour^{\vec{m_3}}v_4\cdot\WDfour^{\vec{m_4}}(\Betabar,\FBetabar)\cdot\WDfour^{\vec{m_5}}(\angv)}_{L^2(\suubar)}\\
		\notag\lesssim&\left(\BAsize\frac{\delta a^{\frac{1}{2}}}{\abs{u}}\right)^{i_0+k_0}\cdot\delta b^{\frac{1}{2}}\abs{u}\cdot\left[\left(\BAsize\frac{\delta a^{\frac{1}{2}}}{\abs{u}^3}+1\right)\cdot b^{-\frac{1}{2}}+\left(\BAsize\frac{ \delta a^{\frac{1}{2}}}{\abs{u}^3}+\BAfluid b^{-\frac{1}{2}}\right)\cdot\BAfluid\right].
	\end{align}
	For the fifth and sixth terms of \eqref{EQ:varvG}:
	\begin{align}\label{PV:varvG5}
		&\sum\limits_{\substack{\vec{m_1}+\vec{m_2}+\vec{m_3}+\vec{m_4}=(i,j-j_0,k-k_0)\\\abs{\vec{m_3}}\geq 1
		}}\abs{u}^{j_0+k_0}\norm{\WDfour^{\vec{m_1}}(\upeta+\ueta)^{k_0}\cdot\WDfour^{\vec{m_2}}\uomega^{j_0}\cdot\WDfour^{\vec{m_3}}(v_3)\cdot\WDfour^{\vec{m_4}}\Dfour_4v_3}_{L^2(\suubar)}\\
		\notag\lesssim&\left(\BAsize\frac{\delta a^{\frac{1}{2}}}{\abs{u}}\right)^{j_0+k_0}\cdot\BAfluid^2 b^{-1}\cdot\delta^{-1}\cdot \abs{u}.
	\end{align}
	
	\begin{align}\label{PV:varvG6}
		&\sum\limits_{\substack{\vec{m_1}+\vec{m_2}+\vec{m_3}+\vec{m_4}+\vec{m_5}=(i,j-j_0,k-k_0)\\\vec{m_1},\vec{m_2},\vec{m_4}\neq(i,j-j_0,k-k_0)}}\abs{u}^{j_0+k_0}\norm{\WDfour^{\vec{m_1}}(\upeta+\ueta)^{k_0}\cdot\WDfour^{\vec{m_2}}\uomega^{j_0}\cdot\WDfour^{\vec{m_3}}(v_3)\cdot\WDfour^{\vec{m_4}}(\upchi,\uomega)\cdot\WDfour^{\vec{m_5}}(v_3)}_{L^2(\suubar)}\\
		\notag\lesssim&\left(\BAsize\frac{\delta a^{\frac{1}{2}}}{\abs{u}}\right)^{j_0+k_0}\cdot \BAsize a^{\frac{1}{2}}\cdot\BAfluid b^{-1}.
	\end{align}
	For the seventh and the eighth terms of \eqref{EQ:varvG}, we have:
	\begin{align}\label{PV:varvG7}
		&\sum\limits_{\vec{m_1}+\vec{m_2}+\vec{m_3}+\vec{m_4}+\vec{m_5}=(i,j-j_0,k-k_0-1)}\abs{u}^{j_0+k_0+1}\norm{\WDfour^{\vec{m_1}}(\upeta+\ueta)^{k_0}\cdot\WDfour^{\vec{m_2}}\uomega^{i_0}\cdot\WDfour^{\vec{m_3}}(v_3)\cdot\WDfour^{\vec{m_4}}\Riem^{J}_{34A}\cdot\WDfour^{\vec{m_5}}(v_{J})}_{L^2(\suubar)}\\
		\notag\lesssim&\sum\limits_{\vec{m_1}+\vec{m_2}+\vec{m_3}+\vec{m_4}+\vec{m_5}=(i,j-j_0,k-k_0-1)}\abs{u}^{j_0+k_0+1}\norm{\WDfour^{\vec{m_1}}(\upeta+\ueta)^{k_0}\cdot\WDfour^{\vec{m_2}}\uomega^{i_0}\cdot\WDfour^{\vec{m_3}}(v_3)\cdot\WDfour^{\vec{m_4}}(\BetaA,\FBetaA)\cdot\WDfour^{\vec{m_5}}(v_3)}_{L^2(\suubar)}\\
		\notag&+\sum\limits_{\vec{m_1}+\vec{m_2}+\vec{m_3}+\vec{m_4}+\vec{m_5}=(i,j-j_0,k-k_0-1)}\abs{u}^{j_0+k_0+1}\norm{\WDfour^{\vec{m_1}}(\upeta+\ueta)^{k_0}\cdot\WDfour^{\vec{m_2}}\uomega^{i_0}\cdot\WDfour^{\vec{m_3}}(v_3)\cdot\WDfour^{\vec{m_4}}(\Rho,\Rhodual,\FRho,\FRhodual,\Scalar)\cdot\WDfour^{\vec{m_5}}(\angv)}_{L^2(\suubar)}\\
		\notag\lesssim&\left(\BAsize\frac{\delta a^{\frac{1}{2}}}{\abs{u}}\right)^{j_0+k_0}\cdot b^{-\frac{1}{2}}\abs{u}^2\cdot\left[\left(\BAsize\frac{a^{\frac{1}{2}}}{\abs{u}^2}+b^{\frac{1}{2}}\right)\cdot b^{-\frac{1}{2}}+\left(\BAsize\frac{\delta a^{\frac{1}{2}}}{\abs{u}^3}+1\right)\cdot\BAfluid\right].
	\end{align}
	
	\begin{align}\label{PV:varvG8}
		&\sum\limits_{\vec{m_1}+\vec{m_2}+\vec{m_3}+\vec{m_4}+\vec{m_5}=(i,j-j_0-1,k-k_0)}\abs{u}^{j_0+k_0+1}\norm{\WDfour^{\vec{m_1}}(\upeta+\ueta)^{k_0}\cdot\Dfour^{\vec{m_2}}\uomega^{i_0}\cdot\WDfour^{\vec{m_3}}(v_3)\cdot\WDfour^{\vec{m_4}}\Riem^{J}_{334}\cdot\WDfour^{\vec{m_5}}(v_J)}_{L^2(\suubar)}\\
		\notag\lesssim&\sum\limits_{\vec{m_1}+\vec{m_2}+\vec{m_3}+\vec{m_4}+\vec{m_5}=(i,j-j_0-1,k-k_0)}\abs{u}^{j_0+k_0+1}\norm{\WDfour^{\vec{m_1}}(\upeta+\ueta)^{k_0}\cdot\WDfour^{\vec{m_2}}\uomega^{i_0}\cdot\WDfour^{\vec{m_3}}(v_3)\cdot\WDfour^{\vec{m_4}}(\Rho,\FRho,\Scalar)\cdot\WDfour^{\vec{m_5}}(v_3)}_{L^2(\suubar)}\\
		\notag&+\sum\limits_{\vec{m_1}+\vec{m_2}+\vec{m_3}+\vec{m_4}+\vec{m_5}=(i,j-j_0-1,k-k_0)}\abs{u}^{j_0+k_0+1}\norm{\WDfour^{\vec{m_1}}(\upeta+\ueta)^{k_0}\cdot\WDfour^{\vec{m_2}}\uomega^{i_0}\cdot\WDfour^{\vec{m_3}}(v_3)\cdot\WDfour^{\vec{m_4}}(\Betabar,\FBetabar)\cdot\WDfour^{\vec{m_5}}(\angv)}_{L^2(\suubar)}\\
		\notag\lesssim&\left(\BAsize\frac{\delta a^{\frac{1}{2}}}{\abs{u}}\right)^{i_0+k_0}\cdot b^{-\frac{1}{2}}\abs{u}^2\cdot\left[\left(\BAsize\frac{\delta a^{\frac{1}{2}}}{\abs{u}^3}+1\right)\cdot b^{-\frac{1}{2}}+\left(\BAsize\frac{ \delta a^{\frac{1}{2}}}{\abs{u}^3}+\BAfluid b^{-\frac{1}{2}}\right)\cdot\BAfluid\right].
	\end{align}
	For the last four terms of \eqref{EQ:varvG}, we use the similar approach as in \eqref{PV:varvG1}-\eqref{PV:varvG4}:
	\begin{align}\label{PV:varvG9}
		&\sum\limits_{\substack{\vec{m_1}+\vec{m_2}+\vec{m_3}+\vec{m_4}=(i-i_0,j-j_0,k)\\\abs{\vec{m_3}}\geq 1
		}}\delta^{i_0}\abs{u}^{j_0}\norm{\WDfour^{\vec{m_1}}\upchi^{i_0}\cdot\WDfour^{\vec{m_2}}\uchi^{j_0}\cdot\WDfour^{\vec{m_3}} v_B\cdot\WDfour^{\vec{m_4}}\Dfour_Bv_3}_{L^2(\suubar)}\\
		\notag\lesssim&\left(\BAsize\frac{\delta a^{\frac{1}{2}}}{\abs{u}}\right)^{i_0+j_0}\cdot\BAfluid^2.
	\end{align}
	
	\begin{align}\label{PV:varvG10}
		&\sum\limits_{\substack{\vec{m_1}+\vec{m_2}+\vec{m_3}+\vec{m_4}+\vec{m_5}=(i-i_0,j-j_0,k)\\\vec{m_1},\vec{m_2},\vec{m_4}\neq(i-i_0,j-j_0,k)}}\delta^{i_0}\abs{u}^{j_0}\norm{\WDfour^{\vec{m_1}}\upchi^{i_0}\cdot\WDfour^{\vec{m_2}}\uchi^{j_0}\cdot\WDfour^{\vec{m_3}} v_B\cdot\WDfour^{\vec{m_4}}(\upeta_B+\ueta_B)\cdot\WDfour^{\vec{m_5}}(v_3)}_{L^2(\suubar)}\\
		\notag\lesssim&\left(\BAsize\frac{\delta a^{\frac{1}{2}}}{\abs{u}}\right)^{i_0+j_0+1}\cdot\BAfluid.
	\end{align}
	We bound the 11-12th term of \eqref{EQ:varvG} using Proposition \ref{PR:Weyl} and Proposition \ref{PR:Ricci0}:
	\begin{align}\label{PV:varvG11}
		&\sum\limits_{\vec{m_1}+\vec{m_2}+\vec{m_3}+\vec{m_4}+\vec{m_5}=(i-i_0,j-j_0-1,k)}\delta^{i_0}\abs{u}^{j_0+1}\norm{\WDfour^{\vec{m_1}}\upchi^{i_0}\cdot\WDfour^{\vec{m_2}}\uchi^{ j_0}\cdot\WDfour^{\vec{m_3}} v_B\cdot\WDfour^{\vec{m_4}}\Riem^{J}_{33B}\cdot\WDfour^{\vec{m_5}}(v_J)}_{L^2(\suubar)}\\
		\notag\lesssim&\sum\limits_{\vec{m_1}+\vec{m_2}+\vec{m_3}+\vec{m_4}+\vec{m_5}=(i-i_0,j-j_0-1,k)}\delta^{i_0}\abs{u}^{j_0+1}\norm{\WDfour^{\vec{m_1}}\upchi^{i_0}\cdot\WDfour^{\vec{m_2}}\uchi^{ j_0}\cdot\WDfour^{\vec{m_3}} v_B\cdot\WDfour^{\vec{m_4}}(\Alphabar,\FAlphabar)\cdot\WDfour^{\vec{m_5}}(\angv)}_{L^2(\suubar)}\\
		\notag&+\sum\limits_{\vec{m_1}+\vec{m_2}+\vec{m_3}+\vec{m_4}+\vec{m_5}=(i-i_0,j-j_0-1,k)}\delta^{i_0}\abs{u}^{j_0+1}\norm{\WDfour^{\vec{m_1}}\upchi^{i_0}\cdot\WDfour^{\vec{m_2}}\uchi^{ j_0}\cdot\WDfour^{\vec{m_3}} v_B\cdot\WDfour^{\vec{m_4}}(\Betabar,\FBetabar)\cdot\WDfour^{\vec{m_5}}(v_3)}_{L^2(\suubar)}\\
		\notag\lesssim&\left(\BAsize\frac{\delta a^{\frac{1}{2}}}{\abs{u}}\right)^{i_0+j_0}\cdot\BAfluid\abs{u}^2\cdot\left[\left(\BAsize\frac{\delta a^{\frac{1}{2}}}{\abs{u}^3}+b^{-1}\right)\cdot \BAfluid+\left(\BAsize\frac{\delta a^{\frac{1}{2}}}{\abs{u}^3}+\BAfluid b^{-\frac{1}{2}}\right)\cdot b^{-\frac{1}{2}}\right].
	\end{align}
	
	\begin{align}\label{PV:varvG12}
		&\sum\limits_{\vec{m_1}+\vec{m_2}+\vec{m_3}+\vec{m_4}+\vec{m_5}=(i-i_0-1,j-j_0,k)}\delta^{i_0+1}\abs{u}^{j_0}\norm{\WDfour^{\vec{m_1}}\upchi^{i_0}\cdot\WDfour^{\vec{m_2}}\uchi^{ j_0}\cdot\WDfour^{\vec{m_3}} v_B\cdot\WDfour^{\vec{m_4}}\Riem^{J}_{34B}\cdot\WDfour^{\vec{m_5}}(v_J)}_{L^2(\suubar)}\\
		\notag\lesssim&\sum\limits_{\vec{m_1}+\vec{m_2}+\vec{m_3}+\vec{m_4}+\vec{m_5}=(i-i_0-1,j-j_0,k)}\delta^{i_0+1}\abs{u}^{j_0}\norm{\WDfour^{\vec{m_1}}\upchi^{i_0}\cdot\WDfour^{\vec{m_2}}\uchi^{ j_0}\cdot\WDfour^{\vec{m_3}} v_B\cdot\WDfour^{\vec{m_4}}(\Rho,\Rhodual,\FRho,\FRhodual,\Scalar)\cdot\WDfour^{\vec{m_5}}(\angv)}_{L^2(\suubar)}\\
		\notag&+\sum\limits_{\vec{m_1}+\vec{m_2}+\vec{m_3}+\vec{m_4}+\vec{m_5}=(i-i_0-1,j-j_0,k)}\delta^{i_0+1}\abs{u}^{j_0}\norm{\WDfour^{\vec{m_1}}\upchi^{i_0}\cdot\WDfour^{\vec{m_2}}\uchi^{ j_0}\cdot\WDfour^{\vec{m_3}} v_B\cdot\WDfour^{\vec{m_4}}(\BetaA,\FBetaA)\cdot\WDfour^{\vec{m_5}}(v_3)}_{L^2(\suubar)}\\
		\notag\lesssim&\left(\BAsize\frac{\delta a^{\frac{1}{2}}}{\abs{u}}\right)^{i_0+j_0}\cdot\BAfluid\delta\abs{u}\cdot\left[\left(\BAsize\frac{\delta a^{\frac{1}{2}}}{\abs{u}^3}+1\right)\cdot \BAfluid+\left(\BAsize\frac{ a^{\frac{1}{2}}}{\abs{u}^2}+\BAfluid b^{\frac{1}{2}}\right)\cdot b^{-\frac{1}{2}}\right].
	\end{align}
	Combining \eqref{PV:varvG1}-\eqref{PV:varvG12}, we obtain the estimate for $\Fluid{\varv}{G}_3$:
	\begin{align}\label{PV:varvG}
		\norm{\Fluid{\varv}{G}_3}_{L^2(\suubar)}\lesssim&~\BAfluid^2\delta^{-1}b^{-1}\abs{u}+\BAsize\BAfluid a^{\frac{1}{2}} b^{-1}.
	\end{align}
	Finally, we estimate $\Fluid{\vars}{G}_3$. It is suffice to estimate first two terms in \eqref{EQ:varsG} since the rest is trivial under the contraction with $e_3$. For the first two terms in \eqref{EQ:varsG}, we have:
	\begin{align}\label{PV:varsG1}
		&\sum\limits_{\substack{\vec{m_1}+\vec{m_2}+\vec{m_3}+\vec{m_4}=(i-i_0,j,k-k_0)\\\abs{\vec{m_3}}\geq1
		}}\delta^{i_0}\abs{u}^{k_0}\norm{\WDfour^{\vec{m_1}}(\upeta+\ueta)^{k_0}\cdot\WDfour^{\vec{m_2}}\omega^{i_0}\cdot\WDfour^{\vec{m_3}}\frac{\temp}{\enth}\cdot\WDfour^{\vec{m_4}}\Dfour_3\ent}_{L^2(\suubar)}\\
		\notag\lesssim&\left(\BAsize\frac{\delta a^{\frac{1}{2}}}{\abs{u}}\right)^{i_0+k_0}\cdot\BAfluid^2.
	\end{align}
	\begin{align}\label{PV:varsG2}
		&\sum\limits_{\substack{\vec{m_1}+\vec{m_2}+\vec{m_3}+\vec{m_4}+\vec{m_5}=(i-i_0,j,k-k_0)\\\vec{m_1},\vec{m_2},\vec{m_4}\neq(i-i_0,j,k-k_0)}}\delta^{i_0}\abs{u}^{k_0}\norm{\WDfour^{\vec{m_1}}(\upeta+\ueta)^{k_0}\cdot\WDfour^{\vec{m_2}}\omega^{i_0}\cdot\WDfour^{\vec{m_3}}\frac{\temp}{\enth}\cdot\WDfour^{\vec{m_4}}(\uchi,\uomega,\upzeta)\cdot\WDfour^{\vec{m_5}}\ent}_{L^2(\suubar)}\\
		\notag\lesssim&\left(\BAsize\frac{\delta a^{\frac{1}{2}}}{\abs{u}}\right)^{i_0+k_0+1}\cdot\BAfluid.
	\end{align}
	Combining \eqref{PV:varsG1}-\eqref{PV:varsG2}, we obtain the estimate for $\Fluid{\vars}{G}_3$:
	\begin{align}\label{PV:varsG}
		\norm{\Fluid{\vars}{G}_3}_{L^2(\suubar)}\lesssim&~\BAfluid^2.
	\end{align}
	Then, the desired result \eqref{ES:EstimateFGH2} follows from \eqref{PV:varhG}, \eqref{PV:varvG} and \eqref{PV:varsG}.
\end{proof}
\begin{proof}[Proof of \eqref{ES:EstimateFGH3}]
	The proof of \eqref{ES:EstimateFGH3} is very similar to the proof of estimates \eqref{ES:EstimateFGH3}. In below, we omit repetitive details.
	First of all,
	As in Proposition \ref{PR:EQVariations}, $\mathfrak{G}_{(i,j,k)}:=\Fluid{\varh}{G}+\Fluid{\varv}{G}+\Fluid{\vars}{G}$. We start with estimating terms in $\Fluid{\varh}{G}_B$. For the seventh term on the RHS of \eqref{EQ:varhG}, we have:
	\begin{align}\label{PV:varhGA1}
		&\sum\limits_{\substack{\vec{m_1}+\vec{m_2}+\vec{m_3}+\vec{m_4}+\vec{m_5}=(i-i_0,j,k-k_0)\\\abs{\vec{m_3}+\vec{m_4}}\geq1
		}}\delta^{i_0}\abs{u}^{k_0}\norm{\WDfour^{\vec{m_1}}(\upeta+\ueta)^{k_0}\cdot\WDfour^{\vec{m_2}}\omega^{i_0}\cdot\WDfour^{\vec{m_3}}(\angv)\cdot\WDfour^{\vec{m_4}}v_4\cdot\WDfour^{\vec{m_5}}\Dfour_3\lnenth}_{L^2(\suubar)}\\
		\notag\lesssim&\left(\BAsize\frac{\delta a^{\frac{1}{2}}}{\abs{u}}\right)^{i_0+k_0}\cdot\BAfluid^2 b^{\frac{1}{2}}.
	\end{align}
	For the eighth terms on the RHS of \eqref{EQ:varhG}:
	\begin{align}\label{PV:varhGA2}
		&\sum\limits_{\substack{\vec{m_1}+\vec{m_2}+\vec{m_3}+\vec{m_4}+\vec{m_5}+\vec{m_6}=(i-i_0,j,k-k_0)\\\vec{m_1},\vec{m_2},\vec{m_5}\neq(i-i_0,j,k-k_0)}}\delta^{i_0}\abs{u}^{k_0}\norm{\WDfour^{\vec{m_1}}(\upeta+\ueta)^{k_0}\cdot\WDfour^{\vec{m_2}}\omega^{i_0}\cdot\WDfour^{\vec{m_3}}(\angv)\cdot\WDfour^{\vec{m_4}}v_4\cdot\WDfour^{\vec{m_5}}(\uchi,\uomega,\upzeta)\cdot\WDfour^{\vec{m_6}}\lnenth}_{L^2(\suubar)}\\
		\notag\lesssim&\left(\BAsize\frac{\delta a^{\frac{1}{2}}}{\abs{u}}\right)^{i_0+k_0+1}\cdot\BAfluid^2 b ^{\frac{1}{2}}.
	\end{align}
	The 9-10th terms of \eqref{EQ:varhG} obtain the following estimates:
	\begin{align}\label{PV:varhGA3}
		&\sum\limits_{\substack{\vec{m_1}+\vec{m_2}+\vec{m_3}+\vec{m_4}=(i,j-j_0,k-k_0)\\\abs{\vec{m_3}+\vec{m_4}}\geq 1
		}}\abs{u}^{j_0+k_0}\norm{\WDfour^{\vec{m_1}}(\upeta+\ueta)^{k_0}\cdot\WDfour^{\vec{m_2}}\uomega^{j_0}\cdot\WDfour^{\vec{m_3}}(\angv)\cdot\WDfour^{\vec{m_4}}(v_3)\cdot\WDfour^{\vec{m_5}}\Dfour_4\lnenth}_{L^2(\suubar)}\\
		\notag\lesssim&\left(\BAsize\frac{\delta a^{\frac{1}{2}}}{\abs{u}}\right)^{j_0+k_0}\cdot\BAfluid^2 b^{-\frac{1}{2}}\cdot\delta^{-1}\cdot\abs{u}.
	\end{align}
	\begin{align}\label{PV:varhGA4}
		&\sum\limits_{\substack{\vec{m_1}+\vec{m_2}+\vec{m_3}+\vec{m_4}+\vec{m_5}+\vec{m_6}=(i,j-j_0,k-k_0)\\\vec{m_1},\vec{m_2},\vec{m_5}\neq(i,j-j_0,k-k_0)}}\abs{u}^{j_0+k_0}\norm{\WDfour^{\vec{m_1}}(\upeta+\ueta)^{k_0}\cdot\WDfour^{\vec{m_2}}\uomega^{j_0}\cdot\WDfour^{\vec{m_3}}(\angv)\cdot\WDfour^{\vec{m_4}}(v_3)\cdot\WDfour^{\vec{m_5}}(\upchi,\omega)\cdot\WDfour^{\vec{m_6}}\lnenth}_{L^2(\suubar)}\\
		\notag\lesssim&\left(\BAsize\frac{\delta a^{\frac{1}{2}}}{\abs{u}}\right)^{j_0+k_0}\cdot \BAsize a^{\frac{1}{2}}\cdot\BAfluid^2 b^{-\frac{1}{2}}.
	\end{align}
	For the 5-6th, and last two terms in $\Fluid{\varh}{G}$, we have:
	\begin{align}\label{PV:varhGA5}
		&\sum\limits_{\substack{\vec{m_1}+\vec{m_2}+\vec{m_3}+\vec{m_4}+\vec{m_5}=(i-i_0,j-j_0,k)\\\abs{\vec{m_3}+\vec{m_4}}\geq 1
		}}\delta^{i_0}\abs{u}^{j_0}\norm{\WDfour^{\vec{m_1}}\upchi^{i_0}\cdot\WDfour^{\vec{m_2}}\uchi^{j_0}\cdot\WDfour^{\vec{m_3}}(\angv)\cdot\WDfour^{\vec{m_4}}v_B\cdot\WDfour^{\vec{m_5}}\Dfour_B\lnenth}_{L^2(\suubar)}\\
		\notag\lesssim&\left(\BAsize\frac{\delta a^{\frac{1}{2}}}{\abs{u}}\right)^{i_0+j_0}\cdot\BAfluid^2 b^{\frac{1}{2}}.
	\end{align}
	\begin{align}\label{PV:varhGA6}
		&\sum\limits_{\substack{\vec{m_1}+\vec{m_2}+\vec{m_3}+\vec{m_4}+\vec{m_5}+\vec{m_6}=(i-i_0,j-j_0,k)\\\vec{m_1},\vec{m_2},\vec{m_5}\neq(i-i_0,j-j_0,k)}}\delta^{i_0}\abs{u}^{j_0}\norm{\WDfour^{\vec{m_1}}\upchi^{i_0}\cdot\WDfour^{\vec{m_2}}\uchi^{j_0}\cdot\WDfour^{\vec{m_3}}(\angv)\cdot\WDfour^{\vec{m_4}}v_B\cdot\WDfour^{\vec{m_5}}(\upeta_B+\ueta_B)\cdot\WDfour^{\vec{m_6}}\lnenth}_{L^2(\suubar)}\\
		\notag\lesssim&\left(\BAsize\frac{\delta a^{\frac{1}{2}}}{\abs{u}}\right)^{i_0+j_0+1}\cdot\BAfluid b^{\frac{1}{2}}.
	\end{align}
	Combining \eqref{PV:varhGA1}-\eqref{PV:varhGA6}, we obtain the estimate for $\Fluid{\varh}{G}_B$:
	\begin{align}\label{PV:varhGA}
		\norm{\Fluid{\varh}{G}_B}_{L^2(\suubar)}\lesssim&~\BAfluid^2\delta^{-1}b^{-\frac{1}{2}}\abs{u}+\BAsize\BAfluid a^{\frac{1}{2}} b^{-\frac{1}{2}}.
	\end{align}
	
	Now, we estimate terms in $\Fluid{\varv}{G}_B$. For the first term on the RHS of \eqref{EQ:varvG}, we have:
	\begin{align}\label{PV:varvGA1}
		&\sum\limits_{\substack{\vec{m_1}+\vec{m_2}+\vec{m_3}+\vec{m_4}=(i-i_0,j,k-k_0)\\\abs{\vec{m_3}}\geq 1
		}}\delta^{i_0}\abs{u}^{k_0}\norm{\WDfour^{\vec{m_1}}(\upeta+\ueta)^{k_0}\cdot\WDfour^{\vec{m_2}}\omega^{i_0}\cdot\WDfour^{\vec{m_3}}v_4\cdot\WDfour^{\vec{m_4}}\Dfour_3\angv}_{L^2(\suubar)}\\
		\notag\lesssim&\left(\BAsize\frac{\delta a^{\frac{1}{2}}}{\abs{u}}\right)^{i_0+k_0}\cdot\BAfluid^2 b^{\frac{1}{2}}.
	\end{align}
	For the second term on the RHS of \eqref{EQ:varvG}:
	\begin{align}\label{PV:varvGA2}
		&\sum\limits_{\substack{\vec{m_1}+\vec{m_2}+\vec{m_3}+\vec{m_4}+\vec{m_5}=(i-i_0,j,k-k_0)\\\vec{m_1},\vec{m_2},\vec{m_4}\neq(i-i_0,j,k-k_0)}}\delta^{i_0}\abs{u}^{k_0}\norm{\WDfour^{\vec{m_1}}(\upeta+\ueta)^{k_0}\cdot\WDfour^{\vec{m_2}}\omega^{i_0}\cdot\WDfour^{\vec{m_3}}v_4\cdot\WDfour^{\vec{m_4}}(\uchi,\uomega,\upzeta)\cdot\WDfour^{\vec{m_5}}(\angv)}_{L^2(\suubar)}\\
		\notag\lesssim&\left(\BAsize\frac{\delta a^{\frac{1}{2}}}{\abs{u}}\right)^{i_0+k_0+1}\cdot\BAfluid b^{\frac{1}{2}}.
	\end{align}
	We bound the third term of \eqref{EQ:varvG} using Proposition \ref{PR:Weyl} and Proposition \ref{PR:Ricci0}:
	\begin{align}\label{PV:varvGA3}
		&\sum\limits_{\vec{m_1}+\vec{m_2}+\vec{m_3}+\vec{m_4}+\vec{m_5}=(i-i_0,j,k-k_0-1)}\delta^{i_0}\abs{u}^{k_0+1}\norm{\WDfour^{\vec{m_1}}(\upeta+\ueta)^{k_0}\cdot\WDfour^{\vec{m_2}}\omega^{i_0}\cdot\WDfour^{\vec{m_3}}v_4\cdot\WDfour^{\vec{m_4}}\Riem^{J}_{A3C}\cdot\WDfour^{\vec{m_5}}(v_{J})}_{L^2(\suubar)}\\
		\notag\lesssim&\sum\limits_{\vec{m_1}+\vec{m_2}+\vec{m_3}+\vec{m_4}+\vec{m_5}=(i-i_0,j,k-k_0-1)}\delta^{i_0}\abs{u}^{k_0+1}\norm{\WDfour^{\vec{m_1}}(\upeta+\ueta)^{k_0}\cdot\WDfour^{\vec{m_2}}\omega^{i_0}\cdot\WDfour^{\vec{m_3}}v_4\cdot\WDfour^{\vec{m_4}}(\Rho,\Rhodual,\FRho,\FRhodual,\Scalar)\cdot\WDfour^{\vec{m_5}}(v_3)}_{L^2(\suubar)}\\
		\notag&+\sum\limits_{\vec{m_1}+\vec{m_2}+\vec{m_3}+\vec{m_4}+\vec{m_5}=(i-i_0,j,k-k_0-1)}\delta^{i_0}\abs{u}^{k_0+1}\norm{\WDfour^{\vec{m_1}}(\upeta+\ueta)^{k_0}\cdot\WDfour^{\vec{m_2}}\omega^{i_0}\cdot\WDfour^{\vec{m_3}}v_4\cdot\WDfour^{\vec{m_4}}(\Betabar,\FBetabar)\cdot\WDfour^{\vec{m_5}}(\angv)}_{L^2(\suubar)}\\
		\notag&+\sum\limits_{\vec{m_1}+\vec{m_2}+\vec{m_3}+\vec{m_4}+\vec{m_5}=(i-i_0,j,k-k_0-1)}\delta^{i_0}\abs{u}^{k_0+1}\norm{\WDfour^{\vec{m_1}}(\upeta+\ueta)^{k_0}\cdot\WDfour^{\vec{m_2}}\omega^{i_0}\cdot\WDfour^{\vec{m_3}}v_4\cdot\WDfour^{\vec{m_4}}(\Alphabar,\FAlphabar)\cdot\WDfour^{\vec{m_5}}v_4}_{L^2(\suubar)}\\
		\notag\lesssim&\left(\BAsize\frac{\delta a^{\frac{1}{2}}}{\abs{u}}\right)^{i_0+k_0}\cdot b^{\frac{1}{2}}\abs{u}^2\cdot\left[\left(\BAsize\frac{\delta a^{\frac{1}{2}}}{\abs{u}^3}+1\right)\cdot b^{-\frac{1}{2}}+\left(\BAsize\frac{\delta a^{\frac{1}{2}}}{\abs{u}^3}+\BAfluid b^{-\frac{1}{2}}\right)\cdot \BAfluid+\left(\BAsize\frac{\delta^2 a}{\abs{u}^3}+b^{-1}\right)\cdot b^{\frac{1}{2}}\right].
	\end{align}
	Similar to \eqref{PV:varvGA3}, we estimate the forth term of \eqref{EQ:varvG} using Proposition \ref{PR:Weyl} and Proposition \ref{PR:Ricci0}:
	\begin{align}\label{PV:varvGA4}
		&\sum\limits_{\vec{m_1}+\vec{m_2}+\vec{m_3}+\vec{m_4}+\vec{m_5}=(i-i_0-1,j,k-k_0)}\delta^{i_0+1}\abs{u}^{k_0}\norm{\WDfour^{\vec{m_1}}(\upeta+\ueta)^{k_0}\cdot\WDfour^{\vec{m_2}}\omega^{i_0}\cdot\WDfour^{\vec{m_3}}v_4\cdot\WDfour^{\vec{m_4}}\Riem^{J}_{A34}\cdot\WDfour^{\vec{m_5}}(v_J)}_{L^2(\suubar)}\\
		\notag\lesssim&\sum\limits_{\vec{m_1}+\vec{m_2}+\vec{m_3}+\vec{m_4}+\vec{m_5}=(i-i_0-1,j,k-k_0)}\delta^{i_0+1}\abs{u}^{k_0}\norm{\WDfour^{\vec{m_1}}(\upeta+\ueta)^{k_0}\cdot\WDfour^{\vec{m_2}}\omega^{i_0}\cdot\WDfour^{\vec{m_3}}v_4\cdot\WDfour^{\vec{m_4}}(\BetaA,\FBetaA)\cdot\WDfour^{\vec{m_5}}(v_3)}_{L^2(\suubar)}\\
		\notag&+\sum\limits_{\vec{m_1}+\vec{m_2}+\vec{m_3}+\vec{m_4}+\vec{m_5}=(i-i_0-1,j,k-k_0)}\delta^{i_0+1}\abs{u}^{k_0}\norm{\WDfour^{\vec{m_1}}(\upeta+\ueta)^{k_0}\cdot\WDfour^{\vec{m_2}}\omega^{i_0}\cdot\WDfour^{\vec{m_3}}v_4\cdot\WDfour^{\vec{m_4}}\Rhodual\cdot\WDfour^{\vec{m_5}}(\angv)}_{L^2(\suubar)}\\
		\notag&+\sum\limits_{\vec{m_1}+\vec{m_2}+\vec{m_3}+\vec{m_4}+\vec{m_5}=(i-i_0-1,j,k-k_0)}\delta^{i_0+1}\abs{u}^{k_0}\norm{\WDfour^{\vec{m_1}}(\upeta+\ueta)^{k_0}\cdot\WDfour^{\vec{m_2}}\omega^{i_0}\cdot\WDfour^{\vec{m_3}}v_4\cdot\WDfour^{\vec{m_4}}(\Betabar,\FBetabar)\cdot\WDfour^{\vec{m_5}}v_4}_{L^2(\suubar)}\\
		\notag\lesssim&\left(\BAsize\frac{\delta a^{\frac{1}{2}}}{\abs{u}}\right)^{i_0+k_0}\cdot\delta b^{\frac{1}{2}}\abs{u}\cdot\left[\left(\BAsize\frac{ a^{\frac{1}{2}}}{\abs{u}^2}+\BAfluid b^{\frac{1}{2}}\right)\cdot b^{-\frac{1}{2}}+\left(\BAsize\frac{ \delta a^{\frac{1}{2}}}{\abs{u}^3}\right)\cdot\BAfluid+\left(\BAsize\frac{ \delta a^{\frac{1}{2}}}{\abs{u}^3}+\BAfluid b^{-\frac{1}{2}}\right)\cdot b^{\frac{1}{2}}\right].
	\end{align}
	For the fifth and sixth terms of \eqref{EQ:varvG}:
	\begin{align}\label{PV:varvGA5}
		&\sum\limits_{\substack{\vec{m_1}+\vec{m_2}+\vec{m_3}+\vec{m_4}=(i,j-j_0,k-k_0)\\\abs{\vec{m_3}}\geq 1
		}}\abs{u}^{j_0+k_0}\norm{\WDfour^{\vec{m_1}}(\upeta+\ueta)^{k_0}\cdot\WDfour^{\vec{m_2}}\uomega^{j_0}\cdot\WDfour^{\vec{m_3}}(v_3)\cdot\WDfour^{\vec{m_4}}\Dfour_4\angv}_{L^2(\suubar)}\\
		\notag\lesssim&\left(\BAsize\frac{\delta a^{\frac{1}{2}}}{\abs{u}}\right)^{j_0+k_0}\cdot\BAfluid^2 b^{-\frac{1}{2}}\cdot\delta^{-1}\cdot \abs{u}.
	\end{align}
	
	\begin{align}\label{PV:varvGA6}
		&\sum\limits_{\substack{\vec{m_1}+\vec{m_2}+\vec{m_3}+\vec{m_4}+\vec{m_5}=(i,j-j_0,k-k_0)\\\vec{m_1},\vec{m_2},\vec{m_4}\neq(i,j-j_0,k-k_0)}}\abs{u}^{j_0+k_0}\norm{\WDfour^{\vec{m_1}}(\upeta+\ueta)^{k_0}\cdot\WDfour^{\vec{m_2}}\uomega^{j_0}\cdot\WDfour^{\vec{m_3}}(v_3)\cdot\WDfour^{\vec{m_4}}(\upchi,\omega)\cdot\WDfour^{\vec{m_5}}(\angv)}_{L^2(\suubar)}\\
		\notag\lesssim&\left(\BAsize\frac{\delta a^{\frac{1}{2}}}{\abs{u}}\right)^{j_0+k_0}\cdot \BAsize a^{\frac{1}{2}}\cdot\BAfluid b^{-\frac{1}{2}}.
	\end{align}
	For the seventh and the eighth terms of \eqref{EQ:varvG}, we have:
	\begin{align}\label{PV:varvGA7}
		&\sum\limits_{\vec{m_1}+\vec{m_2}+\vec{m_3}+\vec{m_4}+\vec{m_5}=(i,j-j_0,k-k_0-1)}\abs{u}^{j_0+k_0+1}\norm{\WDfour^{\vec{m_1}}(\upeta+\ueta)^{k_0}\cdot\WDfour^{\vec{m_2}}\uomega^{i_0}\cdot\WDfour^{\vec{m_3}}(v_3)\cdot\WDfour^{\vec{m_4}}\Riem^{J}_{A4C}\cdot\WDfour^{\vec{m_5}}(v_{J})}_{L^2(\suubar)}\\
		\notag\lesssim&\sum\limits_{\vec{m_1}+\vec{m_2}+\vec{m_3}+\vec{m_4}+\vec{m_5}=(i,j-j_0,k-k_0-1)}\abs{u}^{j_0+k_0+1}\norm{\WDfour^{\vec{m_1}}(\upeta+\ueta)^{k_0}\cdot\WDfour^{\vec{m_2}}\uomega^{i_0}\cdot\WDfour^{\vec{m_3}}(v_3)\cdot\WDfour^{\vec{m_4}}(\AlphaAB,\FAlphaAB)\cdot\WDfour^{\vec{m_5}}(v_3)}_{L^2(\suubar)}\\
		\notag&+\sum\limits_{\vec{m_1}+\vec{m_2}+\vec{m_3}+\vec{m_4}+\vec{m_5}=(i,j-j_0,k-k_0-1)}\abs{u}^{j_0+k_0+1}\norm{\WDfour^{\vec{m_1}}(\upeta+\ueta)^{k_0}\cdot\WDfour^{\vec{m_2}}\uomega^{i_0}\cdot\WDfour^{\vec{m_3}}(v_3)\cdot\WDfour^{\vec{m_4}}(\BetaA,\FBetaA)\cdot\WDfour^{\vec{m_5}}(\angv)}_{L^2(\suubar)}\\
		\notag&+\sum\limits_{\vec{m_1}+\vec{m_2}+\vec{m_3}+\vec{m_4}+\vec{m_5}=(i,j-j_0,k-k_0-1)}\abs{u}^{j_0+k_0+1}\norm{\WDfour^{\vec{m_1}}(\upeta+\ueta)^{k_0}\cdot\WDfour^{\vec{m_2}}\uomega^{i_0}\cdot\WDfour^{\vec{m_3}}(v_3)\cdot\WDfour^{\vec{m_4}}(\Rho,\Rhodual,\FRho,\FRhodual,\Scalar)\cdot\WDfour^{\vec{m_5}}v_4}_{L^2(\suubar)}\\
		\notag\lesssim&\left(\BAsize\frac{\delta a^{\frac{1}{2}}}{\abs{u}}\right)^{j_0+k_0}\cdot b^{-\frac{1}{2}}\abs{u}^2\cdot\left[\left(\BAsize\frac{\delta^{-1}a^{\frac{3}{8}}}{\abs{u}}+b\right)\cdot b^{-\frac{1}{2}}+\left(\BAsize\frac{ a^{\frac{1}{2}}}{\abs{u}^2}+\BAfluid b^{\frac{1}{2}}\right)\cdot\BAfluid+\left(\BAsize\frac{\delta a^{\frac{1}{2}}}{\abs{u}^3}+1\right)\cdot b^{\frac{1}{2}}\right].
	\end{align}
	
	\begin{align}\label{PV:varvGA8}
		&\sum\limits_{\vec{m_1}+\vec{m_2}+\vec{m_3}+\vec{m_4}+\vec{m_5}=(i,j-j_0-1,k-k_0)}\abs{u}^{j_0+k_0+1}\norm{\WDfour^{\vec{m_1}}(\upeta+\ueta)^{k_0}\cdot\Dfour^{\vec{m_2}}\uomega^{i_0}\cdot\WDfour^{\vec{m_3}}(v_3)\cdot\WDfour^{\vec{m_4}}\Riem^{J}_{A34}\cdot\WDfour^{\vec{m_5}}(v_J)}_{L^2(\suubar)}\\
		\notag\lesssim&\sum\limits_{\vec{m_1}+\vec{m_2}+\vec{m_3}+\vec{m_4}+\vec{m_5}=(i,j-j_0-1,k-k_0)}\abs{u}^{j_0+k_0+1}\norm{\WDfour^{\vec{m_1}}(\upeta+\ueta)^{k_0}\cdot\WDfour^{\vec{m_2}}\uomega^{i_0}\cdot\WDfour^{\vec{m_3}}(v_3)\cdot\WDfour^{\vec{m_4}}(\BetaA,\FBetaA)\cdot\WDfour^{\vec{m_5}}(v_3)}_{L^2(\suubar)}\\
		\notag&+\sum\limits_{\vec{m_1}+\vec{m_2}+\vec{m_3}+\vec{m_4}+\vec{m_5}=(i,j-j_0-1,k-k_0)}\abs{u}^{j_0+k_0+1}\norm{\WDfour^{\vec{m_1}}(\upeta+\ueta)^{k_0}\cdot\WDfour^{\vec{m_2}}\uomega^{i_0}\cdot\WDfour^{\vec{m_3}}(v_3)\cdot\WDfour^{\vec{m_4}}\Rhodual\cdot\WDfour^{\vec{m_5}}(\angv)}_{L^2(\suubar)}\\
		\notag&+\sum\limits_{\vec{m_1}+\vec{m_2}+\vec{m_3}+\vec{m_4}+\vec{m_5}=(i,j-j_0-1,k-k_0)}\abs{u}^{j_0+k_0+1}\norm{\WDfour^{\vec{m_1}}(\upeta+\ueta)^{k_0}\cdot\WDfour^{\vec{m_2}}\uomega^{i_0}\cdot\WDfour^{\vec{m_3}}(v_3)\cdot\WDfour^{\vec{m_4}}(\Betabar,\FBetabar)\cdot\WDfour^{\vec{m_5}}v_4}_{L^2(\suubar)}\\
		\notag\lesssim&\left(\BAsize\frac{\delta a^{\frac{1}{2}}}{\abs{u}}\right)^{i_0+k_0}\cdot b^{-\frac{1}{2}}\abs{u}^2\cdot\left[\left(\BAsize\frac{ a^{\frac{1}{2}}}{\abs{u}^2}+\BAfluid b^{\frac{1}{2}}\right)\cdot b^{-\frac{1}{2}}+\BAsize\frac{ \delta a^{\frac{1}{2}}}{\abs{u}^3}\cdot\BAfluid+\left(\BAsize\frac{ \delta a^{\frac{1}{2}}}{\abs{u}^3}+\BAfluid b^{-\frac{1}{2}}\right)\cdot b^{\frac{1}{2}}\right].
	\end{align}
	For the last four terms of \eqref{EQ:varvG}, we use the similar approach as in \eqref{PV:varvGA1}-\eqref{PV:varvGA4}:
	\begin{align}\label{PV:varvGA9}
		&\sum\limits_{\substack{\vec{m_1}+\vec{m_2}+\vec{m_3}+\vec{m_4}=(i-i_0,j-j_0,k)\\\abs{\vec{m_3}}\geq 1
		}}\delta^{i_0}\abs{u}^{j_0}\norm{\WDfour^{\vec{m_1}}\upchi^{i_0}\cdot\WDfour^{\vec{m_2}}\uchi^{j_0}\cdot\WDfour^{\vec{m_3}} v_B\cdot\WDfour^{\vec{m_4}}\Dfour_Bv_B}_{L^2(\suubar)}\\
		\notag\lesssim&\left(\BAsize\frac{\delta a^{\frac{1}{2}}}{\abs{u}}\right)^{i_0+j_0}\cdot\BAfluid^2 b^{\frac{1}{2}}.
	\end{align}
	
	\begin{align}\label{PV:varvGA10}
		&\sum\limits_{\substack{\vec{m_1}+\vec{m_2}+\vec{m_3}+\vec{m_4}+\vec{m_5}=(i-i_0,j-j_0,k)\\\vec{m_1},\vec{m_2},\vec{m_4}\neq(i-i_0,j-j_0,k)}}\delta^{i_0}\abs{u}^{j_0}\norm{\WDfour^{\vec{m_1}}\upchi^{i_0}\cdot\WDfour^{\vec{m_2}}\uchi^{j_0}\cdot\WDfour^{\vec{m_3}} v_B\cdot\WDfour^{\vec{m_4}}(\upeta_B+\ueta_B)\cdot\WDfour^{\vec{m_5}}(\angv)}_{L^2(\suubar)}\\
		\notag\lesssim&\left(\BAsize\frac{\delta a^{\frac{1}{2}}}{\abs{u}}\right)^{i_0+j_0+1}\cdot\BAfluid^2.
	\end{align}
	We bound the 11-12th term of \eqref{EQ:varvG} using Proposition \ref{PR:Weyl} and Proposition \ref{PR:Ricci0}:
	\begin{align}\label{PV:varvGA11}
		&\sum\limits_{\vec{m_1}+\vec{m_2}+\vec{m_3}+\vec{m_4}+\vec{m_5}=(i-i_0,j-j_0-1,k)}\delta^{i_0}\abs{u}^{j_0+1}\norm{\WDfour^{\vec{m_1}}\upchi^{i_0}\cdot\WDfour^{\vec{m_2}}\uchi^{ j_0}\cdot\WDfour^{\vec{m_3}} v_B\cdot\WDfour^{\vec{m_4}}\Riem^{J}_{A3B}\cdot\WDfour^{\vec{m_5}}(v_J)}_{L^2(\suubar)}\\
		\notag\lesssim&\sum\limits_{\vec{m_1}+\vec{m_2}+\vec{m_3}+\vec{m_4}+\vec{m_5}=(i-i_0,j-j_0-1,k)}\delta^{i_0}\abs{u}^{j_0+1}\norm{\WDfour^{\vec{m_1}}\upchi^{i_0}\cdot\WDfour^{\vec{m_2}}\uchi^{ j_0}\cdot\WDfour^{\vec{m_3}} v_B\cdot\WDfour^{\vec{m_4}}(\Alphabar,\FAlphabar)\cdot\WDfour^{\vec{m_5}}v_4}_{L^2(\suubar)}\\
		\notag&+\sum\limits_{\vec{m_1}+\vec{m_2}+\vec{m_3}+\vec{m_4}+\vec{m_5}=(i-i_0,j-j_0-1,k)}\delta^{i_0}\abs{u}^{j_0+1}\norm{\WDfour^{\vec{m_1}}\upchi^{i_0}\cdot\WDfour^{\vec{m_2}}\uchi^{ j_0}\cdot\WDfour^{\vec{m_3}} v_B\cdot\WDfour^{\vec{m_4}}(\Betabar,\FBetabar)\cdot\WDfour^{\vec{m_5}}(\angv)}_{L^2(\suubar)}\\
		\notag&+\sum\limits_{\vec{m_1}+\vec{m_2}+\vec{m_3}+\vec{m_4}+\vec{m_5}=(i-i_0,j-j_0-1,k)}\delta^{i_0}\abs{u}^{j_0+1}\norm{\WDfour^{\vec{m_1}}\upchi^{i_0}\cdot\WDfour^{\vec{m_2}}\uchi^{ j_0}\cdot\WDfour^{\vec{m_3}} v_B\cdot\WDfour^{\vec{m_4}}(\Rho,\Rhodual,\FRho,\FRhodual,\Scalar)\cdot\WDfour^{\vec{m_5}}(v_3)}_{L^2(\suubar)}\\
		\notag\lesssim&\left(\BAsize\frac{\delta a^{\frac{1}{2}}}{\abs{u}}\right)^{i_0+j_0}\cdot\BAfluid b^{\frac{1}{2}}\abs{u}^2\cdot\left[\left(\BAsize\frac{\delta^2 a}{\abs{u}^3}+b^{-1}\right)\cdot b^{\frac{1}{2}}+\left(\BAsize\frac{\delta a^{\frac{1}{2}}}{\abs{u}^3}+\BAfluid b^{-\frac{1}{2}}\right)\cdot \BAfluid+\left(\BAsize\frac{\delta a^{\frac{1}{2}}}{\abs{u}^3}+1\right)\cdot b^{-\frac{1}{2}}\right].
	\end{align}
	
	\begin{align}\label{PV:varvGA12}
		&\sum\limits_{\vec{m_1}+\vec{m_2}+\vec{m_3}+\vec{m_4}+\vec{m_5}=(i-i_0-1,j-j_0,k)}\delta^{i_0+1}\abs{u}^{j_0}\norm{\WDfour^{\vec{m_1}}\upchi^{i_0}\cdot\WDfour^{\vec{m_2}}\uchi^{ j_0}\cdot\WDfour^{\vec{m_3}} v_B\cdot\WDfour^{\vec{m_4}}\Riem^{J}_{A4B}\cdot\WDfour^{\vec{m_5}}(v_J)}_{L^2(\suubar)}\\
		\notag\lesssim&\sum\limits_{\vec{m_1}+\vec{m_2}+\vec{m_3}+\vec{m_4}+\vec{m_5}=(i-i_0-1,j-j_0,k)}\delta^{i_0+1}\abs{u}^{j_0}\norm{\WDfour^{\vec{m_1}}\upchi^{i_0}\cdot\WDfour^{\vec{m_2}}\uchi^{ j_0}\cdot\WDfour^{\vec{m_3}} v_B\cdot\WDfour^{\vec{m_4}}(\Rho,\Rhodual,\FRho,\FRhodual,\Scalar)\cdot\WDfour^{\vec{m_5}}v_4}_{L^2(\suubar)}\\
		\notag&+\sum\limits_{\vec{m_1}+\vec{m_2}+\vec{m_3}+\vec{m_4}+\vec{m_5}=(i-i_0-1,j-j_0,k)}\delta^{i_0+1}\abs{u}^{j_0}\norm{\WDfour^{\vec{m_1}}\upchi^{i_0}\cdot\WDfour^{\vec{m_2}}\uchi^{ j_0}\cdot\WDfour^{\vec{m_3}} v_B\cdot\WDfour^{\vec{m_4}}(\BetaA,\FBetaA)\cdot\WDfour^{\vec{m_5}}(\angv)}_{L^2(\suubar)}\\
		\notag&+\sum\limits_{\vec{m_1}+\vec{m_2}+\vec{m_3}+\vec{m_4}+\vec{m_5}=(i-i_0-1,j-j_0,k)}\delta^{i_0+1}\abs{u}^{j_0}\norm{\WDfour^{\vec{m_1}}\upchi^{i_0}\cdot\WDfour^{\vec{m_2}}\uchi^{ j_0}\cdot\WDfour^{\vec{m_3}} v_B\cdot\WDfour^{\vec{m_4}}(\AlphaAB,\FAlphaAB)\cdot\WDfour^{\vec{m_5}}(v_3)}_{L^2(\suubar)}\\
		\notag\lesssim&\left(\BAsize\frac{\delta a^{\frac{1}{2}}}{\abs{u}}\right)^{i_0+j_0}\cdot\BAfluid b^{\frac{1}{2}}\delta\abs{u}\cdot\left[\left(\BAsize\frac{\delta a^{\frac{1}{2}}}{\abs{u}^3}+1\right)\cdot b^{\frac{1}{2}}+\left(\BAsize\frac{\delta^{-1} a^{\frac{1}{4}}}{\abs{u}}+b\right)\cdot b^{-\frac{1}{2}}+\left(\BAsize\frac{ a^{\frac{1}{2}}}{\abs{u}^2}+\BAfluid b^{\frac{1}{2}}\right)\cdot \BAfluid\right].
	\end{align}
	Combining \eqref{PV:varvGA1}-\eqref{PV:varvGA12}, we obtain the estimate for $\Fluid{\varv}{G}_B$:
	\begin{align}\label{PV:varvGA}
		\norm{\Fluid{\varv}{G}_B}_{L^2(\suubar)}\lesssim&~\BAfluid^2\delta^{-1}b^{-\frac{1}{2}}\abs{u}+\BAsize\delta^{-1}a^{\frac{1}{4}} b^{-1}\abs{u}.
	\end{align}
	Finally, we estimate $\Fluid{\vars}{G}_B$. It is suffice to estimate last two terms in \eqref{EQ:varsG} since the rest is trivial under the contraction with $e_A$. For the last two terms in \eqref{EQ:varsG}, we have:
	\begin{align}\label{PV:varsGA1}
		&\sum\limits_{\substack{\vec{m_1}+\vec{m_2}+\vec{m_3}+\vec{m_4}=(i-i_0,j-j_0,k)\\\abs{\vec{m_3}}\geq 1}}\delta^{i_0}\abs{u}^{j_0}\norm{\WDfour^{\vec{m_1}}\uchi^{j_0}\cdot\WDfour^{\vec{m_2}}\upchi^{i_0}\cdot\WDfour^{\vec{m_3}}\frac{\temp}{\enth}\cdot\WDfour^{\vec{m_4}}\angD\ent}_{L^2(\suubar)}\\
		\notag\lesssim&\left(\BAsize\frac{\delta a^{\frac{1}{2}}}{\abs{u}}\right)^{i_0+j_0}\cdot\BAfluid^2.
	\end{align}
	\begin{align}\label{PV:varsGA2}
		&\sum\limits_{\substack{\vec{m_1}+\vec{m_2}+\vec{m_3}+\vec{m_4}+\vec{m_5}=(i-i_0,j-j_0,k)\\\vec{m_1},\vec{m_2},\vec{m_4}\neq(i-i_0,j-j_0,k)}}\delta^{i_0}\abs{u}^{j_0}\norm{\WDfour^{\vec{m_1}}\uchi^{j_0}\cdot\WDfour^{\vec{m_2}}\upchi^{i_0}\cdot\WDfour^{\vec{m_3}}\frac{\temp}{\enth}\cdot\WDfour^{\vec{m_4}}(\upeta,\ueta)\cdot\WDfour^{\vec{m_5}}\ent}_{L^2(\suubar)}\\
		\notag\lesssim&\left(\BAsize\frac{\delta a^{\frac{1}{2}}}{\abs{u}}\right)^{i_0+j_0+1}\cdot\BAfluid.
	\end{align}
	Combining \eqref{PV:varsGA1}-\eqref{PV:varsGA2}, we obtain the estimate for $\Fluid{\vars}{G}_B$:
	\begin{align}\label{PV:varsGA}
		\norm{\Fluid{\vars}{G}_B}_{L^2(\suubar)}\lesssim&~\BAfluid^2.
	\end{align}
	Then, the desired result \eqref{ES:EstimateFGH2} follows from \eqref{PV:varhGA}, \eqref{PV:varvGA} and \eqref{PV:varsGA}.
\end{proof}
\begin{proof}[Proof of \eqref{ES:EstimateFGH4}]
	We take the similar approach as in the proof of estimates \eqref{ES:EstimateFGH2}-\eqref{ES:EstimateFGH3}. Therefore, in below, we omit repetitive details.
	First of all,
	As in Proposition \ref{PR:EQVariations}, $\mathfrak{G}^\alpha_{(i,j,k)}:=\Fluid{\varh}{G}^\alpha+\Fluid{\varv}{G}^\alpha+\Fluid{\vars}{G}^\alpha$. We start with estimating terms in $\Fluid{\varh}{G}_4$. For the seventh terms on the RHS of \eqref{EQ:varhG}, we have:
	\begin{align}\label{PV:varhGe41}
		&\sum\limits_{\substack{\vec{m_1}+\vec{m_2}+\vec{m_3}+\vec{m_4}=(i-i_0,j,k-k_0)\\\abs{\vec{m_3}}\geq1
		}}\delta^{i_0}\abs{u}^{k_0}\norm{\WDfour^{\vec{m_1}}(\upeta+\ueta)^{k_0}\cdot\WDfour^{\vec{m_2}}\omega^{i_0}\cdot\WDfour^{\vec{m_3}}(v_4v_4)\cdot\WDfour^{\vec{m_4}}\Dfour_3\lnenth}_{L^2(\suubar)}\\
		\notag\lesssim&\left(\BAsize\frac{\delta a^{\frac{1}{2}}}{\abs{u}}\right)^{i_0+k_0}\cdot\BAfluid^2 b.
	\end{align}
	For the eighth terms on the RHS of \eqref{EQ:varhG}:
	\begin{align}\label{PV:varhGe42}
		&\sum\limits_{\substack{\vec{m_1}+\vec{m_2}+\vec{m_3}+\vec{m_4}+\vec{m_5}=(i-i_0,j,k-k_0)\\\vec{m_1},\vec{m_2},\vec{m_4}\neq(i-i_0,j,k-k_0)}}\delta^{i_0}\abs{u}^{k_0}\norm{\WDfour^{\vec{m_1}}(\upeta+\ueta)^{k_0}\cdot\WDfour^{\vec{m_2}}\omega^{i_0}\cdot\WDfour^{\vec{m_3}}(v_4v_4)\cdot\WDfour^{\vec{m_4}}(\uchi,\uomega,\upzeta)\cdot\WDfour^{\vec{m_5}}\lnenth}_{L^2(\suubar)}\\
		\notag\lesssim&\left(\BAsize\frac{\delta a^{\frac{1}{2}}}{\abs{u}}\right)^{i_0+k_0+1}\cdot\BAfluid b.
	\end{align}
	The 3-4th, and 9-10th terms of \eqref{EQ:varhG} obtain the following estimates:
	\begin{align}\label{PV:varhGe43}
		&\sum\limits_{\substack{\vec{m_1}+\vec{m_2}+\vec{m_3}+\vec{m_4}=(i,j-j_0,k-k_0)\\\abs{\vec{m_3}}\geq 1
		}}\abs{u}^{j_0+k_0}\norm{\WDfour^{\vec{m_1}}(\upeta+\ueta)^{k_0}\cdot\WDfour^{\vec{m_2}}\uomega^{j_0}\cdot\left(1,\WDfour^{\vec{m_3}}(v_3v_4)\right)\cdot\WDfour^{\vec{m_4}}\Dfour_4\lnenth}_{L^2(\suubar)}\\
		\notag\lesssim&\left(\BAsize\frac{\delta a^{\frac{1}{2}}}{\abs{u}}\right)^{j_0+k_0}\cdot\BAfluid^2\cdot\delta^{-1}\cdot\abs{u}.
	\end{align}
	\begin{align}\label{PV:varhGe44}
		&\sum\limits_{\substack{\vec{m_1}+\vec{m_2}+\vec{m_3}+\vec{m_4}+\vec{m_5}=(i,j-j_0,k-k_0)\\\vec{m_1},\vec{m_2},\vec{m_4}\neq(i,j-j_0,k-k_0)}}\abs{u}^{j_0+k_0}\norm{\WDfour^{\vec{m_1}}(\upeta+\ueta)^{k_0}\cdot\WDfour^{\vec{m_2}}\uomega^{j_0}\cdot\left(1,\WDfour^{\vec{m_3}}(v_3v_4)\right)\cdot\WDfour^{\vec{m_4}}(\upchi,\omega)\cdot\WDfour^{\vec{m_5}}\lnenth}_{L^2(\suubar)}\\
		\notag\lesssim&\left(\BAsize\frac{\delta a^{\frac{1}{2}}}{\abs{u}}\right)^{j_0+k_0}\cdot \BAsize a^{\frac{1}{2}}\cdot\BAfluid.
	\end{align}
	For the last two terms in $\Fluid{\varh}{G}$, we have:
	\begin{align}\label{PV:varhGe45}
		&\sum\limits_{\substack{\vec{m_1}+\vec{m_2}+\vec{m_3}+\vec{m_4}=(i-i_0,j-j_0,k)\\\abs{\vec{m_3}}\geq 1
		}}\delta^{i_0}\abs{u}^{j_0}\norm{\WDfour^{\vec{m_1}}\upchi^{i_0}\cdot\WDfour^{\vec{m_2}}\uchi^{j_0}\cdot\WDfour^{\vec{m_3}}(v_Bv_4)\cdot\WDfour^{\vec{m_4}}\Dfour_B\lnenth}_{L^2(\suubar)}\\
		\notag\lesssim&\left(\BAsize\frac{\delta a^{\frac{1}{2}}}{\abs{u}}\right)^{i_0+j_0}\cdot\BAfluid^2 b ^{\frac{1}{2}}.
	\end{align}
	\begin{align}\label{PV:varhGe46}
		&\sum\limits_{\substack{\vec{m_1}+\vec{m_2}+\vec{m_3}+\vec{m_4}+\vec{m_5}=(i-i_0,j-j_0,k)\\\vec{m_1},\vec{m_2},\vec{m_4}\neq(i-i_0,j-j_0,k)}}\delta^{i_0}\abs{u}^{j_0}\norm{\WDfour^{\vec{m_1}}\upchi^{i_0}\cdot\WDfour^{\vec{m_2}}\uchi^{j_0}\cdot\WDfour^{\vec{m_3}}(v_Bv_4)\cdot\WDfour^{\vec{m_4}}(\upeta_B+\ueta_B)\cdot\WDfour^{\vec{m_5}}\lnenth}_{L^2(\suubar)}\\
		\notag\lesssim&\left(\BAsize\frac{\delta a^{\frac{1}{2}}}{\abs{u}}\right)^{i_0+j_0+1}\cdot\BAfluid^2 b^{\frac{1}{2}}.
	\end{align}
	Combining \eqref{PV:varhGe41}-\eqref{PV:varhGe46}, we obtain the estimate for $\Fluid{\varh}{G}_4$:
	\begin{align}\label{PV:varhGe4}
		\norm{\Fluid{\varh}{G}^{\alpha}\cdot (e_4)_\alpha}_{L^2(\suubar)}\lesssim&~\BAfluid^2\delta^{-1}\abs{u}+\BAsize\BAfluid a^{\frac{1}{2}}.
	\end{align}
	
	Now, we estimate terms in $\Fluid{\varv}{G}^\alpha\cdot(e_4)_{\alpha}$. For the first term on the RHS of \eqref{EQ:varvG}, we have:
	\begin{align}\label{PV:varvGe41}
		&\sum\limits_{\substack{\vec{m_1}+\vec{m_2}+\vec{m_3}+\vec{m_4}=(i-i_0,j,k-k_0)\\\abs{\vec{m_3}}\geq 1
		}}\delta^{i_0}\abs{u}^{k_0}\norm{\WDfour^{\vec{m_1}}(\upeta+\ueta)^{k_0}\cdot\WDfour^{\vec{m_2}}\omega^{i_0}\cdot\WDfour^{\vec{m_3}}v_4\cdot\WDfour^{\vec{m_4}}\Dfour_3v_4}_{L^2(\suubar)}\\
		\notag\lesssim&\left(\BAsize\frac{\delta a^{\frac{1}{2}}}{\abs{u}}\right)^{i_0+k_0}\cdot\BAfluid^2 b.
	\end{align}
	For the second term on the RHS of \eqref{EQ:varvG}:
	\begin{align}\label{PV:varvGe42}
		&\sum\limits_{\substack{\vec{m_1}+\vec{m_2}+\vec{m_3}+\vec{m_4}+\vec{m_5}=(i-i_0,j,k-k_0)\\\vec{m_1},\vec{m_2},\vec{m_4}\neq(i-i_0,j,k-k_0)}}\delta^{i_0}\abs{u}^{k_0}\norm{\WDfour^{\vec{m_1}}(\upeta+\ueta)^{k_0}\cdot\WDfour^{\vec{m_2}}\omega^{i_0}\cdot\WDfour^{\vec{m_3}}v_4\cdot\WDfour^{\vec{m_4}}(\uchi,\uomega,\upzeta)\cdot\WDfour^{\vec{m_5}}v_4}_{L^2(\suubar)}\\
		\notag\lesssim&\left(\BAsize\frac{\delta a^{\frac{1}{2}}}{\abs{u}}\right)^{i_0+k_0+1}\cdot\BAfluid b.
	\end{align}
	We bound the third term of \eqref{EQ:varvG} using Proposition \ref{PR:Weyl} and Proposition \ref{PR:Ricci0}:
	\begin{align}\label{PV:varvGe43}
		&\sum\limits_{\vec{m_1}+\vec{m_2}+\vec{m_3}+\vec{m_4}+\vec{m_5}=(i-i_0,j,k-k_0-1)}\delta^{i_0}\abs{u}^{k_0+1}\norm{\WDfour^{\vec{m_1}}(\upeta+\ueta)^{k_0}\cdot\WDfour^{\vec{m_2}}\omega^{i_0}\cdot\WDfour^{\vec{m_3}}v_4\cdot\WDfour^{\vec{m_4}}\Riem^{J}_{43A}\cdot\WDfour^{\vec{m_5}}(v_{J})}_{L^2(\suubar)}\\
		\notag\lesssim&\sum\limits_{\vec{m_1}+\vec{m_2}+\vec{m_3}+\vec{m_4}+\vec{m_5}=(i-i_0,j,k-k_0-1)}\delta^{i_0}\abs{u}^{k_0+1}\norm{\WDfour^{\vec{m_1}}(\upeta+\ueta)^{k_0}\cdot\WDfour^{\vec{m_2}}\omega^{i_0}\cdot\WDfour^{\vec{m_3}}v_4\cdot\WDfour^{\vec{m_4}}(\BetaA,\FBetaA)\cdot\WDfour^{\vec{m_5}}(v_3)}_{L^2(\suubar)}\\
		\notag&+\sum\limits_{\vec{m_1}+\vec{m_2}+\vec{m_3}+\vec{m_4}+\vec{m_5}=(i-i_0,j,k-k_0-1)}\delta^{i_0}\abs{u}^{k_0+1}\norm{\WDfour^{\vec{m_1}}(\upeta+\ueta)^{k_0}\cdot\WDfour^{\vec{m_2}}\omega^{i_0}\cdot\WDfour^{\vec{m_3}}v_4\cdot\WDfour^{\vec{m_4}}(\AlphaAB,\FAlphaAB)\cdot\WDfour^{\vec{m_5}}(\angv)}_{L^2(\suubar)}\\
		\notag\lesssim&\left(\BAsize\frac{\delta a^{\frac{1}{2}}}{\abs{u}}\right)^{i_0+k_0}\cdot b^{\frac{1}{2}}\abs{u}^2\cdot\left[\left(\BAsize\frac{a^{\frac{1}{2}}}{\abs{u}^2}+\BAfluid b^{\frac{1}{2}}\right)\cdot b^{-\frac{1}{2}}+\left(\BAsize\frac{\delta^{-1} a^{\frac{3}{8}}}{\abs{u}}+b\right)\cdot\BAfluid\right].
	\end{align}
	Similar to \eqref{PV:varvGe43}, we estimate the forth term of \eqref{EQ:varvG} using Proposition \ref{PR:Weyl} and Proposition \ref{PR:Ricci0}:
	\begin{align}\label{PV:varvGe44}
		&\sum\limits_{\vec{m_1}+\vec{m_2}+\vec{m_3}+\vec{m_4}+\vec{m_5}=(i-i_0-1,j,k-k_0)}\delta^{i_0+1}\abs{u}^{k_0}\norm{\WDfour^{\vec{m_1}}(\upeta+\ueta)^{k_0}\cdot\WDfour^{\vec{m_2}}\omega^{i_0}\cdot\WDfour^{\vec{m_3}}v_4\cdot\WDfour^{\vec{m_4}}\Riem^{J}_{434}\cdot\WDfour^{\vec{m_5}}(v_J)}_{L^2(\suubar)}\\
		\notag\lesssim&\sum\limits_{\vec{m_1}+\vec{m_2}+\vec{m_3}+\vec{m_4}+\vec{m_5}=(i-i_0-1,j,k-k_0)}\delta^{i_0+1}\abs{u}^{k_0}\norm{\WDfour^{\vec{m_1}}(\upeta+\ueta)^{k_0}\cdot\WDfour^{\vec{m_2}}\omega^{i_0}\cdot\WDfour^{\vec{m_3}}v_4\cdot\WDfour^{\vec{m_4}}(\Rho,\FRho,\Scalar)\cdot\WDfour^{\vec{m_5}}v_{4}}_{L^2(\suubar)}\\
		\notag&+\sum\limits_{\vec{m_1}+\vec{m_2}+\vec{m_3}+\vec{m_4}+\vec{m_5}=(i-i_0-1,j,k-k_0)}\delta^{i_0+1}\abs{u}^{k_0}\norm{\WDfour^{\vec{m_1}}(\upeta+\ueta)^{k_0}\cdot\WDfour^{\vec{m_2}}\omega^{i_0}\cdot\WDfour^{\vec{m_3}}v_4\cdot\WDfour^{\vec{m_4}}(\BetaA,\FBetaA)\cdot\WDfour^{\vec{m_5}}(\angv)}_{L^2(\suubar)}\\
		\notag\lesssim&\left(\BAsize\frac{\delta a^{\frac{1}{2}}}{\abs{u}}\right)^{i_0+k_0}\cdot\delta b^{\frac{1}{2}}\abs{u}\cdot\left[\left(\BAsize\frac{\delta a^{\frac{1}{2}}}{\abs{u}^3}+1\right)\cdot b^{\frac{1}{2}}+\left(\BAsize\frac{a^{\frac{1}{2}}}{\abs{u}^2}+\BAfluid b^{\frac{1}{2}}\right)\cdot\BAfluid\right].
	\end{align}
	For the fifth and sixth terms of \eqref{EQ:varvG}:
	\begin{align}\label{PV:varvGe45}
		&\sum\limits_{\substack{\vec{m_1}+\vec{m_2}+\vec{m_3}+\vec{m_4}=(i,j-j_0,k-k_0)\\\abs{\vec{m_3}}\geq 1
		}}\abs{u}^{j_0+k_0}\norm{\WDfour^{\vec{m_1}}(\upeta+\ueta)^{k_0}\cdot\WDfour^{\vec{m_2}}\uomega^{j_0}\cdot\WDfour^{\vec{m_3}}(v_3)\cdot\WDfour^{\vec{m_4}}\Dfour_4v_4}_{L^2(\suubar)}\\
		\notag\lesssim&\left(\BAsize\frac{\delta a^{\frac{1}{2}}}{\abs{u}}\right)^{j_0+k_0}\cdot\BAfluid^2 \cdot\delta^{-1}\cdot \abs{u}.
	\end{align}
	
	\begin{align}\label{PV:varvGe46}
		&\sum\limits_{\substack{\vec{m_1}+\vec{m_2}+\vec{m_3}+\vec{m_4}+\vec{m_5}=(i,j-j_0,k-k_0)\\\vec{m_1},\vec{m_2},\vec{m_4}\neq(i,j-j_0,k-k_0)}}\abs{u}^{j_0+k_0}\norm{\WDfour^{\vec{m_1}}(\upeta+\ueta)^{k_0}\cdot\WDfour^{\vec{m_2}}\uomega^{j_0}\cdot\WDfour^{\vec{m_3}}(v_3)\cdot\WDfour^{\vec{m_4}}(\upchi,\omega)\cdot\WDfour^{\vec{m_5}}v_4}_{L^2(\suubar)}\\
		\notag\lesssim&\left(\BAsize\frac{\delta a^{\frac{1}{2}}}{\abs{u}}\right)^{j_0+k_0}\cdot \BAsize a^{\frac{1}{2}}\cdot\BAfluid.
	\end{align}
	For the seventh and the eighth terms of \eqref{EQ:varvG}, we have:
	\begin{align}\label{PV:varvGe47}
		&\sum\limits_{\vec{m_1}+\vec{m_2}+\vec{m_3}+\vec{m_4}+\vec{m_5}=(i,j-j_0,k-k_0-1)}\abs{u}^{j_0+k_0+1}\norm{\WDfour^{\vec{m_1}}(\upeta+\ueta)^{k_0}\cdot\WDfour^{\vec{m_2}}\uomega^{i_0}\cdot\WDfour^{\vec{m_3}}(v_3)\cdot\WDfour^{\vec{m_4}}\Riem^{J}_{44A}\cdot\WDfour^{\vec{m_5}}(v_{J})}_{L^2(\suubar)}\\
		\notag\lesssim&\sum\limits_{\vec{m_1}+\vec{m_2}+\vec{m_3}+\vec{m_4}+\vec{m_5}=(i,j-j_0,k-k_0-1)}\abs{u}^{j_0+k_0+1}\norm{\WDfour^{\vec{m_1}}(\upeta+\ueta)^{k_0}\cdot\WDfour^{\vec{m_2}}\uomega^{i_0}\cdot\WDfour^{\vec{m_3}}(v_3)\cdot\WDfour^{\vec{m_4}}(\BetaA,\FBetaA)\cdot\WDfour^{\vec{m_5}}v_{4}}_{L^2(\suubar)}\\
		\notag&+\sum\limits_{\vec{m_1}+\vec{m_2}+\vec{m_3}+\vec{m_4}+\vec{m_5}=(i,j-j_0,k-k_0-1)}\abs{u}^{j_0+k_0+1}\norm{\WDfour^{\vec{m_1}}(\upeta+\ueta)^{k_0}\cdot\WDfour^{\vec{m_2}}\uomega^{i_0}\cdot\WDfour^{\vec{m_3}}(v_3)\cdot\WDfour^{\vec{m_4}}(\AlphaAB,\FAlphaAB)\cdot\WDfour^{\vec{m_5}}(\angv)}_{L^2(\suubar)}\\
		\notag\lesssim&\left(\BAsize\frac{\delta a^{\frac{1}{2}}}{\abs{u}}\right)^{j_0+k_0}\cdot b^{-\frac{1}{2}}\abs{u}^2\cdot\left[\left(\BAsize\frac{a^{\frac{1}{2}}}{\abs{u}^2}+b^{\frac{1}{2}}\right)\cdot b^{\frac{1}{2}}+\left(\BAsize\frac{\delta^{-1} a^{\frac{3}{8}}}{\abs{u}}+b\right)\cdot\BAfluid\right].
	\end{align}
	
	\begin{align}\label{PV:varvGe48}
		&\sum\limits_{\vec{m_1}+\vec{m_2}+\vec{m_3}+\vec{m_4}+\vec{m_5}=(i,j-j_0-1,k-k_0)}\abs{u}^{j_0+k_0+1}\norm{\WDfour^{\vec{m_1}}(\upeta+\ueta)^{k_0}\cdot\WDfour^{\vec{m_2}}\uomega^{i_0}\cdot\WDfour^{\vec{m_3}}(v_3)\cdot\WDfour^{\vec{m_4}}\Riem^{J}_{434}\cdot\WDfour^{\vec{m_5}}(v_J)}_{L^2(\suubar)}\\
		\notag\lesssim&\sum\limits_{\vec{m_1}+\vec{m_2}+\vec{m_3}+\vec{m_4}+\vec{m_5}=(i,j-j_0-1,k-k_0)}\abs{u}^{j_0+k_0+1}\norm{\WDfour^{\vec{m_1}}(\upeta+\ueta)^{k_0}\cdot\WDfour^{\vec{m_2}}\uomega^{i_0}\cdot\WDfour^{\vec{m_3}}(v_3)\cdot\WDfour^{\vec{m_4}}(\Rho,\FRho,\Scalar)\cdot\WDfour^{\vec{m_5}}v_{4}}_{L^2(\suubar)}\\
		\notag&+\sum\limits_{\vec{m_1}+\vec{m_2}+\vec{m_3}+\vec{m_4}+\vec{m_5}=(i,j-j_0-1,k-k_0)}\abs{u}^{j_0+k_0+1}\norm{\WDfour^{\vec{m_1}}(\upeta+\ueta)^{k_0}\cdot\WDfour^{\vec{m_2}}\uomega^{i_0}\cdot\WDfour^{\vec{m_3}}(v_3)\cdot\WDfour^{\vec{m_4}}(\BetaA,\FBetaA)\cdot\WDfour^{\vec{m_5}}(\angv)}_{L^2(\suubar)}\\
		\notag\lesssim&\left(\BAsize\frac{\delta a^{\frac{1}{2}}}{\abs{u}}\right)^{j_0+k_0}\cdot b^{-\frac{1}{2}}\abs{u}^2\cdot\left[\left(\BAsize\frac{\delta a^{\frac{1}{2}}}{\abs{u}^3}+1\right)\cdot b^{\frac{1}{2}}+\left(\BAsize\frac{ a^{\frac{1}{2}}}{\abs{u}^2}+\BAfluid b^{\frac{1}{2}}\right)\cdot\BAfluid\right].
	\end{align}
	For the last four terms of \eqref{EQ:varvG}, we use the similar approach as in \eqref{PV:varvGe41}-\eqref{PV:varvGe44}:
	\begin{align}\label{PV:varvGe49}
		&\sum\limits_{\substack{\vec{m_1}+\vec{m_2}+\vec{m_3}+\vec{m_4}=(i-i_0,j-j_0,k)\\\abs{\vec{m_3}}\geq 1
		}}\delta^{i_0}\abs{u}^{j_0}\norm{\WDfour^{\vec{m_1}}\upchi^{i_0}\cdot\WDfour^{\vec{m_2}}\uchi^{j_0}\cdot\WDfour^{\vec{m_3}} (v_B)\cdot\WDfour^{\vec{m_4}}\Dfour_Bv_4}_{L^2(\suubar)}\\
		\notag\lesssim&\left(\BAsize\frac{\delta a^{\frac{1}{2}}}{\abs{u}}\right)^{i_0+j_0}\cdot\BAfluid^2 b^{\frac{1}{2}}.
	\end{align}
	
	\begin{align}\label{PV:varvGe410}
		&\sum\limits_{\substack{\vec{m_1}+\vec{m_2}+\vec{m_3}+\vec{m_4}+\vec{m_5}=(i-i_0,j-j_0,k)\\\vec{m_1},\vec{m_2},\vec{m_4}\neq(i-i_0,j-j_0,k)}}\delta^{i_0}\abs{u}^{j_0}\norm{\WDfour^{\vec{m_1}}\upchi^{i_0}\cdot\WDfour^{\vec{m_2}}\uchi^{j_0}\cdot\WDfour^{\vec{m_3}} (v_B)\cdot\WDfour^{\vec{m_4}}(\upeta_B+\ueta_B)\cdot\WDfour^{\vec{m_5}}v_4}_{L^2(\suubar)}\\
		\notag\lesssim&\left(\BAsize\frac{\delta a^{\frac{1}{2}}}{\abs{u}}\right)^{i_0+j_0+1}\cdot\BAfluid b^{\frac{1}{2}}.
	\end{align}
	We bound the 11-12th term of \eqref{EQ:varvG} using Proposition \ref{PR:Weyl} and Proposition \ref{PR:Ricci0}:
	\begin{align}\label{PV:varvGe411}
		&\sum\limits_{\vec{m_1}+\vec{m_2}+\vec{m_3}+\vec{m_4}+\vec{m_5}=(i-i_0,j-j_0-1,k)}\delta^{i_0}\abs{u}^{j_0+1}\norm{\WDfour^{\vec{m_1}}\upchi^{i_0}\cdot\WDfour^{\vec{m_2}}\uchi^{ j_0}\cdot\WDfour^{\vec{m_3}} (v_B)\cdot\WDfour^{\vec{m_4}}\Riem^{J}_{43B}\cdot\WDfour^{\vec{m_5}}(v_J)}_{L^2(\suubar)}\\
		\notag\lesssim&\sum\limits_{\vec{m_1}+\vec{m_2}+\vec{m_3}+\vec{m_4}+\vec{m_5}=(i-i_0,j-j_0-1,k)}\delta^{i_0}\abs{u}^{j_0+1}\norm{\WDfour^{\vec{m_1}}\upchi^{i_0}\cdot\WDfour^{\vec{m_2}}\uchi^{ j_0}\cdot\WDfour^{\vec{m_3}} (v_B)\cdot\WDfour^{\vec{m_4}}(\Rho,\Rhodual,\FRho,\FRhodual,\Scalar)\cdot\WDfour^{\vec{m_5}}(\angv)}_{L^2(\suubar)}\\
		\notag&+\sum\limits_{\vec{m_1}+\vec{m_2}+\vec{m_3}+\vec{m_4}+\vec{m_5}=(i-i_0,j-j_0-1,k)}\delta^{i_0}\abs{u}^{j_0+1}\norm{\WDfour^{\vec{m_1}}\upchi^{i_0}\cdot\WDfour^{\vec{m_2}}\uchi^{ j_0}\cdot\WDfour^{\vec{m_3}} (v_B)\cdot\WDfour^{\vec{m_4}}(\Betabar,\FBetabar)\cdot\WDfour^{\vec{m_5}}v_4}_{L^2(\suubar)}\\
		\notag\lesssim&\left(\BAsize\frac{\delta a^{\frac{1}{2}}}{\abs{u}}\right)^{i_0+j_0}\cdot\BAfluid b^{\frac{1}{2}}\abs{u}^2\cdot\left[\left(\BAsize\frac{\delta a^{\frac{1}{2}}}{\abs{u}^3}+1\right)\cdot \BAfluid+\left(\BAsize\frac{\delta a^{\frac{1}{2}}}{\abs{u}^3}+\BAfluid b^{-\frac{1}{2}}\right)\cdot b^{\frac{1}{2}}\right].
	\end{align}
	
	\begin{align}\label{PV:varvGe412}
		&\sum\limits_{\vec{m_1}+\vec{m_2}+\vec{m_3}+\vec{m_4}+\vec{m_5}=(i-i_0-1,j-j_0,k)}\delta^{i_0+1}\abs{u}^{j_0}\norm{\WDfour^{\vec{m_1}}\upchi^{i_0}\cdot\WDfour^{\vec{m_2}}\uchi^{ j_0}\cdot\WDfour^{\vec{m_3}} (v_B)\cdot\WDfour^{\vec{m_4}}\Riem^{J}_{44B}\cdot\WDfour^{\vec{m_5}}(v_J)}_{L^2(\suubar)}\\
		\notag\lesssim&\sum\limits_{\vec{m_1}+\vec{m_2}+\vec{m_3}+\vec{m_4}+\vec{m_5}=(i-i_0-1,j-j_0,k)}\delta^{i_0+1}\abs{u}^{j_0}\norm{\WDfour^{\vec{m_1}}\upchi^{i_0}\cdot\WDfour^{\vec{m_2}}\uchi^{ j_0}\cdot\WDfour^{\vec{m_3}} (v_B)\cdot\WDfour^{\vec{m_4}}(\AlphaAB,\FAlphaAB)\cdot\WDfour^{\vec{m_5}}(\angv)}_{L^2(\suubar)}\\
		\notag&+\sum\limits_{\vec{m_1}+\vec{m_2}+\vec{m_3}+\vec{m_4}+\vec{m_5}=(i-i_0-1,j-j_0,k)}\delta^{i_0+1}\abs{u}^{j_0}\norm{\WDfour^{\vec{m_1}}\upchi^{i_0}\cdot\WDfour^{\vec{m_2}}\uchi^{ j_0}\cdot\WDfour^{\vec{m_3}} (v_B)\cdot\WDfour^{\vec{m_4}}(\BetaA,\FBetaA)\cdot\WDfour^{\vec{m_5}}v_4}_{L^2(\suubar)}\\
		\notag\lesssim&\left(\BAsize\frac{\delta a^{\frac{1}{2}}}{\abs{u}}\right)^{i_0+j_0}\cdot\BAfluid\delta b^{\frac{1}{2}}\abs{u}\cdot\left[\left(\BAsize\frac{\delta^{-1} a^{\frac{1}{4}}}{\abs{u}}+b\right)\cdot \BAfluid+\left(\BAsize\frac{ a^{\frac{1}{2}}}{\abs{u}^2}+\BAfluid b^{\frac{1}{2}}\right)\cdot b^{\frac{1}{2}}\right].
	\end{align}
	Combining \eqref{PV:varvGe41}-\eqref{PV:varvGe412}, we obtain the estimate for $\Fluid{\varv}{G}_4$:
	\begin{align}\label{PV:varvGe4}
		\norm{\Fluid{\varv}{G}_4}_{L^2(\suubar)}\lesssim&~\BAfluid^2\delta^{-1}\abs{u}+\BAsize\BAfluid \delta^{-1}a^{\frac{3}{8}} b^{-\frac{1}{2}}\cdot\abs{u}.
	\end{align}
	Finally, we estimate $\Fluid{\vars}{G}_4$. It is suffice to estimate 3-4th terms in \eqref{EQ:varsG} since the rest is trivial under the contraction with $e_4$. For the these two terms in \eqref{EQ:varsG}, we have:
	\begin{align}\label{PV:varsGe41}
		&\sum\limits_{\substack{\vec{m_1}+\vec{m_2}+\vec{m_3}+\vec{m_4}=(i,j-j_0,k-k_0)\\\abs{\vec{m_3}}\geq1
		}}\abs{u}^{j_0+k_0}\norm{\WDfour^{\vec{m_1}}(\upeta+\ueta)^{k_0}\cdot\WDfour^{\vec{m_2}}\uomega^{j_0}\cdot\WDfour^{\vec{m_3}}\frac{\temp}{\enth}\cdot\WDfour^{\vec{m_4}}\Dfour_4\ent}_{L^2(\suubar)}\\
		\notag\lesssim&\left(\BAsize\frac{\delta a^{\frac{1}{2}}}{\abs{u}}\right)^{j_0+k_0}\cdot\delta^{-1}\BAfluid^2.
	\end{align}
	\begin{align}\label{PV:varsGe42}
		&\sum\limits_{\substack{\vec{m_1}+\vec{m_2}+\vec{m_3}+\vec{m_4}+\vec{m_5}=(i,j-j_0,k-k_0)\\\abs{\vec{m_3}}\geq1
		}}\abs{u}^{j_0+k_0}\norm{\WDfour^{\vec{m_1}}(\upeta+\ueta)^{k_0}\cdot\WDfour^{\vec{m_2}}\uomega^{j_0}\cdot\WDfour^{\vec{m_3}}\frac{\temp}{\enth}\cdot\WDfour^{\vec{m_4}}\psi\cdot\WDfour^{\vec{m_5}}\ent}_{L^2(\suubar)}\\
		\notag\lesssim&\left(\BAsize\frac{\delta a^{\frac{1}{2}}}{\abs{u}}\right)^{j_0+k_0}\cdot\frac{a^{\frac{1}{2}}}{\abs{u}}\cdot\BAfluid.
	\end{align}
	Combining \eqref{PV:varsGe41}-\eqref{PV:varsGe42}, we obtain the estimate for $\Fluid{\vars}{G}_4$:
	\begin{align}\label{PV:varsGe4}
		\norm{\Fluid{\vars}{G}_4}_{L^2(\suubar)}\lesssim&~\delta^{-1}\BAfluid^2.
	\end{align}
	Then, the desired result \eqref{ES:EstimateFGH2} follows from \eqref{PV:varhGe4}, \eqref{PV:varvGe4} and \eqref{PV:varsGe4}.
\end{proof}
\begin{proof}[Proof of \eqref{ES:EstimateFGH5}]
	Estimates of $\mathfrak{H}_{(i,j,k)}$ is identical to $\Fluid{\varh}{F}$. Specifically, exactly the same argument goes through by replacing $\lnenth$ by $\ent$ in \eqref{PV:varhF1}-\eqref{PV:varhF6}.
\end{proof}


\begin{proposition}[Spacetime-norm estimates of inhomogenuous terms of equation of variations]\label{PR:EstimateofFGHtop}
	Let spacetime region $\region$ be as in Section \ref{SS:NullFrame}. Let $\BAsize,\BAfluid$ be as in Section \ref{SS:BA}. Let $\mathfrak{F}:=\mathfrak{F}_{(i,j,k)},\mathfrak{G}:=\mathfrak{G}_{(i,j,k)},\mathfrak{H}:=\mathfrak{H}_{(i,j,k)}$ be defined as in Proposition \ref{PR:EQVariations} with $i+j+k\leq N-2$ or $i+j+k\leq N$ with $k\geq 1$. Then under the bootstrap assumptions, we have the following estimates:
	\begin{subequations}\label{ES:EstimateFGHtop}
		\begin{align}
			\label{ES:EstimateFGHtop1}\norm{\mathfrak{F}_{(i,j,k)}}_{L^2(\region)}\lesssim&~\BAfluid^2\delta^{-\frac{1}{2}}b^{-\frac{1}{2}}\abs{u}^{\frac{3}{2}}+\BAsize\BAfluid \delta^{\frac{1}{2}} a^{\frac{1}{2}} b^{-\frac{1}{2}}\abs{u}^{\frac{1}{2}},\\
			\label{ES:EstimateFGHtop2}\norm{(\mathfrak{G}_{(i,j,k)})_3}_{L^2( \region)}\lesssim&~\BAfluid^2\delta^{-\frac{1}{2}}b^{-1}\abs{u}^{\frac{3}{2}}+\BAsize\BAfluid \delta^{\frac{1}{2}} a^{\frac{1}{2}} b^{-1}\abs{u}^{\frac{1}{2}},\\
			\label{ES:EstimateFGHtop3}\norm{(\mathfrak{G}_{(i,j,k)})_B}_{L^2( \region)}\lesssim&~\BAfluid^2\delta^{-\frac{1}{2}}b^{-\frac{1}{2}}\abs{u}^{\frac{3}{2}}+\BAsize\delta^{-\frac{1}{2}}a^{\frac{3}{8}} b^{-1}\abs{u}^{\frac{3}{2}},\\
			\label{ES:EstimateFGHtop4}\norm{(\mathfrak{G}_{(i,j,k)})_4}_{L^2( \region)}\lesssim&~\BAfluid^2\delta^{-\frac{1}{2}}\abs{u}^{\frac{3}{2}}+\BAsize\BAfluid\delta^{-\frac{1}{2}}a^{\frac{1}{2}} b^{-1}\abs{u}^{\frac{3}{2}},\\
			\label{ES:EstimateFGHtop5}\norm{\mathfrak{H}_{(i,j,k)}}_{L^2( \region)}\lesssim&~\BAfluid^2\delta^{-\frac{1}{2}}b^{-\frac{1}{2}}\abs{u}^{\frac{3}{2}}+\BAsize\BAfluid \delta^{\frac{1}{2}} a^{\frac{1}{2}} b^{-\frac{1}{2}}\abs{u}^{\frac{1}{2}}.
		\end{align}
	\end{subequations}
\end{proposition}
\begin{proof}[Discussion of the proof of Proposition \ref{PR:EstimateofFGHtop}]
	For the cases of $i+j+k\leq N-2$ or $i+j+k\leq N-1$ with $k\geq 1$, \eqref{ES:EstimateFGHtop} follows directly from integrating \eqref{ES:EstimateFGH} over $u,\ubar$. For the case of $i+j+k= N$ with $k\geq 1$, we use the exact same strategy as in Section \ref{SSS:ProofofBulk}. In particular, we replace every $\suubar$-norm by spacetime norm and using bootstrap assumptions \eqref{BA13}, \eqref{BA23}, \eqref{BA32} for the top-order terms (since they only have spacetime-norm bounds instead of $\suubar$-norm bounds compared to the lower-order terms). We note that, top-order terms does not enjoy a $\suubar$-norm bound is due to the energy current method for fluid and null pair energy method for the curvature terms (see Proposition \ref{PR:topL2curvature} and Section \ref{SSS:proofofcurvature} for its proof).
\end{proof}
\subsection{Energy Estimates for the Fluid Variables}\label{SS:EnergyEstimatesFluid}
\begin{proposition}[The divergence of the fluid current]\label{PR:DiveofFluidCurrent}
	Let $(\varh,\vars,\varv^\mu)$ be a solution to the equations of variation \eqref{EQ:FOE3}. Then the divergence of the fluid current $\varJ^\mu$ (defined in \eqref{EQ:EnergyCurrent}) can be expressed as follows:
	\begin{align}\label{EQ:current0}
		\Dfour_\mu\varJ^\mu=&~\Dfour_\mu \left(q^2\frac{\temp^2v^\mu}{\enth^2}\right)\cdot\vars^2+\left[\left(\frac{2q^2\temp^2}{\enth^2}\right)\vars-\frac{2\temp\varh}{\enth\speed^2}\right] \mathfrak{H}+\Dfour_\mu\left(\frac{v^\mu}{\speed^2}\right)\cdot\varh^2+\left(\frac{2\varh}{\speed^2}-\frac{2\temp\vars}{\enth\speed^2}\right)\mathfrak{F}\\
		\notag&-\Dfour_\mu\left(\frac{2\temp v^\mu}{\enth \speed^2}\right)\cdot\frac{2v^\mu}{\speed^2}\cdot\varh\vars-2\varv^\mu\cdot\Dfour_\mu\left(\frac{\temp}{\enth}\right)\cdot\vars+2\gfour_{\nu\lambda}\mathfrak{G}^\nu\varv^\lambda+\Dfour_\mu v^\mu\cdot\gfour_{\nu\lambda}\varv^\nu\varv^\lambda.
	\end{align}
\end{proposition}
\begin{proof}[Proof of Proposition \ref{PR:DiveofFluidCurrent}]
	We take divergence of $\varJ$ (defined in \eqref{EQ:EnergyCurrent}) as follows:
	\begin{align}\label{EQ:current1}
		\Dfour_\mu\varJ^\mu=&~\Dfour_\mu \left(q^2\frac{\temp^2v^\mu}{\enth^2}\right)\cdot\vars^2+2\vars\cdot\Dfour_{\mu}\vars\cdot \left(\frac{q^2\temp^2v^\mu}{\enth^2}\right)\\
		\notag&+\Dfour_\mu(\frac{v^\mu}{\speed^2})\cdot\varh^2+\frac{2\varh}{\speed^2}\cdot v^\mu\cdot\Dfour_\mu\varh\\
		\notag&-\Dfour_\mu\left(\frac{2\temp v^\mu}{\enth \speed^2}\right)\cdot\varh\vars-\frac{\temp}{\enth}\cdot\frac{2v^\mu}{\speed^2}\cdot\varh\cdot\Dfour_\mu\vars-\frac{\temp}{\enth}\cdot\frac{2v^\mu}{\speed^2}\cdot\Dfour_\mu\varh\cdot\vars\\
		\notag&+2\Dfour_\mu\varv^\mu\cdot\varh+2\varv^\mu\cdot\Dfour_\mu\varh\\
		\notag&-2\Dfour_{\mu}\varv^\mu\cdot\frac{\temp}{\enth}\cdot\vars-2\varv^\mu\cdot\Dfour_\mu\left(\frac{\temp}{\enth}\right)\cdot\vars-2\frac{\temp}{\enth}\cdot\varv^\mu\cdot\Dfour_\mu\vars\\
		\notag&+\Dfour_\mu v^\mu\cdot\gfour_{\nu\lambda} \varv^\nu\varv^\lambda+v^\mu\gfour_{\nu\lambda}\cdot \Dfour_\mu\varv^\nu\cdot\varv^\lambda+v^\mu\gfour_{\nu\lambda} \varv^\nu\cdot\Dfour_\mu\varv^\lambda
	\end{align}
	Combining \eqref{EQ:current1} with \eqref{EQ:FOE3}, we obtain the desired result \eqref{EQ:current0}.
\end{proof}
\begin{remark}
	We note that on the RHS of \eqref{EQ:current0}, $\Dfour(\varh,\vars,\varv^\alpha)$ terms are absent. That is, in \eqref{EQ:current0}, we ``gain" one derivative due to the special structure of the energy current \eqref{EQ:EnergyCurrent}. One could also consider system \eqref{EQ:FOE1}, where the energy current of the corresponding variation equations is as follows:
	\begin{align}
		\varJ^\mu=v^\mu\vars^2+\frac{1}{\speed^2(\density+\pressure)}v^\mu\varp^2+2\varv^\mu\varp+(\density+\pressure)v^\mu\gfour_{\lambda\nu}\varv^\lambda\varv^\nu.
	\end{align}
	See \cite[Chapter 1]{Christodoulou2007} for more details.
\end{remark}
\begin{lemma}[Coerciveness of the fluid current]\label{LE:Coerciveness}
	We assume the bootstrap assumptions in Section \ref{SS:BA}. Let $(\varh,\vars,\varv^\mu)$ be a solution to the equations of variation \eqref{EQ:FOE3}. Then, we have the following coercive property:
	\begin{subequations}\label{ES:J4J3}
		\begin{align}
			\label{ES:J4}-\varJ_4\geq& C\cdot\left[(-v_4)(\vars^2+\varh^2+\varv_A\varv_A)+(-v_3)\varv_4^2\right],\\
			\label{ES:J3}-\varJ_3\geq& C\cdot \left[(-v_3)(\vars^2+\varh^2+\varv_A\varv_A)+(-v_4)\varv_3^2\right],
		\end{align}
	\end{subequations}
	where $\varJ_I:=\gfour(\varJ,e_I)$, $C>0$ is a constant determined in the proof.
\end{lemma}
\begin{proof}[Proof of Lemma \ref{LE:Coerciveness}]
	We prove \eqref{ES:J3}, the proof of \eqref{ES:J4} is similar . 
	Contracting \eqref{EQ:EnergyCurrent} with $-e_3$, noting that $v_4, v_3$ are negative due to $v^\mu$ being future-directed and timelike, using the decomposition of $\gfour^{-1}$ in \eqref{EQ:decompositionofmetric} and the fact that $\gfour(\varv,v)=0$ in \eqref{EQ:gvarvv=0}, we obtain:
	\begin{align}\label{EQ:J3}
		-\varJ_3=&~-\frac{v_3}{\speed^2}\left[\varh^2+\left(\frac{q^2\temp^2\speed^2}{\enth^2}\right)\vars^2-2\frac{\temp}{\enth}\varh\vars\right]-2\varv_3(\varh-\frac{\temp}{\enth}\vars)-v_3\left(-\varv_3\varv_4+\varv_A\varv_A\right)\\
		\notag=&~-\frac{v_3}{\speed^2}\varh^2-\left(\frac{q^2\temp^2v_3}{\enth^2}\right)\vars^2-v_4\varv_3^2-v_3\varv_A\varv_A\\
		\notag&+\frac{2v_3\temp}{\speed^2\enth}\varh\vars-2\varv_3\varh+\frac{2\temp}{\enth}\varv_3\vars+2v_A\varv_3\varv_A\\
		\notag\geq&-\frac{v_3}{\speed^2}\varh^2-\left(\frac{q^2\temp^2v_3}{\enth^2}\right)\vars^2-v_4\varv_3^2-v_3\varv_A\varv_A\\
		\notag&-\abs{\frac{2v_3\temp}{\speed^2\enth}\varh\vars}-\abs{2\varv_3\varh}-\abs{\frac{2\temp}{\enth}\varv_3\vars}-\abs{2v_A\varv_3\varv_A}\\
		\notag=&~:\vec{x}^T\cdot M\cdot\vec{x},
	\end{align}
	where $I$ is the identity matrix, vectorfield $\vec{x}\in\mathbb{R}^5$ and matrix $M\in\mathbb{R}^{5\times 5}$ are as follows:
	\begin{align}
		\vec{x}:=&~
		\begin{pmatrix}
			\frac{\sqrt{-v_3}}{\speed}\abs{\varh}\\
			\frac{q\temp \sqrt{-v_3}}{\enth}\abs{\vars}\\
			\sqrt{-v_4}\abs{\varv_3}\\
			\sqrt{-v_3}\abs{\varv_1}\\
			\sqrt{-v_3}\abs{\varv_2}
		\end{pmatrix}&
		M:=&~
		\begin{pmatrix}
			1 & -\frac{1}{q\speed} & -\frac{\speed}{\sqrt{v_4v_3}}&0 & 0 \\
			-\frac{1}{q\speed} & 1 & -\frac{1}{q\sqrt{v_4v_3}}&0&0 \\
			-\frac{\speed}{\sqrt{v_4v_3}} & -\frac{1}{q\sqrt{v_4v_3}} & 1&-\frac{\abs{v_1}}{\sqrt{v_4v_3}}&-\frac{\abs{v_2}}{\sqrt{v_4v_3}} \\
			0& 0 & -\frac{\abs{v_1}}{\sqrt{v_4v_3}}&1 & 0 \\
			0 & 0 & -\frac{\abs{v_2}}{\sqrt{v_4v_3}}&0 & 1 \\
		\end{pmatrix}
	\end{align}
	We let $0<C<b^{-1}$ be a small constant (to be determined later). Using the fact that $-v_3v_4+v_Av_A=\gfour(v,v)=-1$, we derive:
	\begin{align}\label{EQ:detM}
		\det (M-CI)=\frac{1-C}{q^2\speed^2}\left\{\left[q^2\speed^2((1-\speed^2))-(1+3\speed^2)\right]+\mathcal{O}_1(C)\right\},
	\end{align}
	where $\mathcal{O}_1(C)\rightarrow0$ as $C\rightarrow0$.
	Using the bootstrap assumptions \eqref{BA},  we choose small $C>0$ and large $q>0$ such that
	\begin{align}\label{ES:qsize}
		q^2>\frac{1+3\speed^2+\mathcal{O}_1(C)}{\speed^2(1-\speed^2)}.
	\end{align} 
	Then, by \eqref{EQ:detM}, $M-CI$ is positive definite. 
	
	Then, by \eqref{EQ:J3} and bootstrap assumptions \eqref{BA},
	\begin{align}
		-\varJ_3\geq C\cdot\vec{x}^T\cdot\vec{x}\geq C\cdot\left[(-v_3)(\vars^2+\varh^2+\varv_A\varv_A)+(-v_4)\varv_3^2\right].
	\end{align}
\end{proof}
\begin{proposition}[Higher-order energy estimates for fluid variables]\label{PR:HigherOrderEnergyEstimates}
	We assume the bootstrap assumptions in Section \ref{SS:BA}. $\Energy,\uEnergy,\tilde{\Energy},\tilde{\uEnergy}$ are defined in Definition \ref{DE:ConeFluxEnergy}. Then, for solutions $(\lnenth,\ent,v^\mu)$ of \eqref{EQ:FOE2} and its variation $(\varh,\vars,\varv^\mu)$ with respect to operator $(\delta \Dfour_4)^{i}(\abs{u}\Dfour_3)^j(\abs{u}\Dfour_A)^k$ with $1\leq i+j+k\leq N$ defined in Definition \ref{DE:EnergyCurrent}, we have the following energy estimates:
	\begin{align}\label{ES:EnergyEst0}
		\Energy\left[\sqrt{-v_4}\vars,\sqrt{-v_4}\varh,\sqrt{-v_3}\varv_4,\sqrt{-v_4}\varv_A\right](u)+\uEnergy\left[\sqrt{-v_3}\vars,\sqrt{-v_3}\varh,\sqrt{-v_4}\varv_3,\sqrt{-v_3}\varv_A\right](\ubar)\lesssim&~\BAfluid^3b^{-\frac{1}{2}}+	\BAsize\BAfluid^3a^{\frac{1}{2}}b^{-1}\\
		\notag\leq&~\BAfluid^{\frac{5}{2}}b^{-\frac{1}{2}}+	\BAfluid^{\frac{5}{2}}a^{\frac{1}{2}}b^{-1}.
	\end{align}
\end{proposition}
\begin{proof}[Proof of Proposition \ref{PR:HigherOrderEnergyEstimates}]
	We apply the divergence theorem for fluid energy current $-\varJ^\mu$ on spacetime region $\region:=\{(u,\ubar,\sphere^1,\sphere^2)\}\subset[u_\infty,u]\times[0,\ubar]\times \Stwo$:
	\begin{align}\label{ES:EnergyEst1}
		\int_{\hypersur_{u}}-\lapse\varJ_4\diff\ubar\diff\spherevol+\int_{\hypersurbar_{\ubar}}-\lapse\varJ_3\diff u\diff\spherevol\leq&~\int_{\hypersur_{u_\infty}}-\lapse\varJ_4\diff\ubar\diff\spherevol+\int_{\hypersurbar_{0}}-\lapse\varJ_3\diff u\diff\spherevol+\int_{\region}\abs{\Dfour_\mu\varJ^\mu}\diff u\diff\ubar\diff\spherevol.
	\end{align}
	By Lemma \ref{LE:Coerciveness}, we have coerciveness of energy current on null hypersurfaces:
	\begin{subequations}\label{ES:EnergyEst2}
		\begin{align}
			\label{ES:EnergyEst2.1}\Energy\left[\sqrt{-v_4}\vars,\sqrt{-v_4}\varh,\sqrt{-v_3}\varv_4,\sqrt{-v_4}\varv_A\right](u)\lesssim&~ \int_{\hypersur_{u}}-\lapse\varJ_4\diff\ubar\diff\spherevol,\\
			\label{ES:EnergyEst2.2}\uEnergy\left[\sqrt{-v_3}\vars,\sqrt{-v_3}\varh,\sqrt{-v_4}\varv_3,\sqrt{-v_3}\varv_A\right](\ubar)\lesssim&~\int_{\hypersurbar_{\ubar}}-\lapse\varJ_3\diff u\diff\spherevol.
		\end{align}
	\end{subequations}
	For the first two terms on the RHS of \eqref{ES:EnergyEst1}, by initial data assumptions \eqref{Data}, we have:
	\begin{align}\label{ES:EnergyEst3}
		&\int_{\hypersur_{u_\infty}}-\lapse\varJ_4\diff\ubar\diff\spherevol+\int_{\hypersurbar_{0}}-\lapse\varJ_3\diff u\diff\spherevol\\
		\notag\lesssim&~\Energy\left[\sqrt{-v_4}\vars,\sqrt{-v_4}\varh,\sqrt{-v_3}\varv_4,\sqrt{-v_4}\varv_A\right](\uinf)+\uEnergy\left[\sqrt{-v_3}\vars,\sqrt{-v_3}\varh,\sqrt{-v_4}\varv_3,\sqrt{-v_3}\varv_A\right](0)\\
		\notag\lesssim&~\iniBAfluid^2 b^{-\frac{1}{2}}.
	\end{align} 
	Now we estimate the last terms on the RHS of \eqref{ES:EnergyEst1}. By bootstrap assumptions \eqref{BA}, combining equation \eqref{EQ:current0} and estimates \eqref{ES:EstimateFGHtop}, we have:
	\begin{align}\label{ES:EnergyEst4}
		&\int_{\region}\abs{\Dfour_\mu\varJ^\mu}\diff u\diff\ubar\diff\spherevol\\
		\notag\leq&~\norm{\Dfour_\mu \left(q^2\frac{\temp^2v^\mu}{\enth^2}\right)}_{L^\infty(\region)}\cdot\norm{\vars}^2_{L^2(\region)}\\
		\notag&+\left(\norm{\left(\frac{2q^2\temp^2}{\enth^2}\right)\vars}_{L^2(\region)}+\norm{\frac{2\temp\varh}{\enth\speed^2}}_{L^2(\region)}\right)\cdot \norm{\mathfrak{H}}_{L^2(\region)}+\norm{\Dfour_\mu\left(\frac{v^\mu}{\speed^2}\right)}_{L^\infty(\region)}\cdot\norm{\varh}^2_{L^2(\region)}\\
		\notag&+\left(\norm{\frac{2\varh}{\speed^2}}_{L^2(\region)}+\norm{\frac{2\temp\vars}{\enth\speed^2}}_{L^2(\region)}\right)\cdot\norm{\mathfrak{F}}_{L^2(\region)}\\
		\notag&+\norm{\Dfour_\mu\left(\frac{2\temp v^\mu}{\enth \speed^2}\right)}_{L^\infty(\region)}\cdot\norm{\varh}_{L^2(\region)}\cdot\norm{\vars}_{L^2(\region)}\\
		\notag&+\norm{2\Dfour_\mu\left(\frac{\temp}{\enth}\right)}_{L^\infty(\region)}\cdot\norm{\varv_4}_{L^2(\region)}\cdot\norm{\vars}_{L^\infty(\region)}\\
		\notag&+2\left[\norm{\mathfrak{G}_4}_{L^2(\region)}\cdot\norm{\varv_3}_{L^2(\region)}+\norm{\mathfrak{G}_3}_{L^2(\region)}\cdot\norm{\varv_4}_{L^2(\region)}+\norm{\mathfrak{G}_A}_{L^2(\region)}\cdot\norm{\varv_A}_{L^2(\region)}\right]\\
		\notag&+\norm{\Dfour_\mu v^\mu}_{L^\infty(\region)}\cdot\norm{\gfour_{\nu\lambda}\varv^\nu\varv^\lambda}_{L^2(\region)}\\
		\notag\lesssim&~\BAfluid^3b^{-\frac{1}{2}}+	\BAsize\BAfluid^3a^{\frac{1}{2}}b^{-1}\\
		\notag\leq&~\BAfluid^{\frac{5}{2}}b^{-\frac{1}{2}}.
	\end{align}
	Here, we use the fact that $\norm{F}_{L^2(\region)}\leq\sup_{u,\ubar}\norm{F}_{L^2(\suubar)}\cdot\delta^{\frac{1}{2}}\abs{u}^{\frac{1}{2}}$ for any tensorfield $F$.
	
	Combining \eqref{ES:EnergyEst1}-\eqref{ES:EnergyEst4}, using the size of $\iniBAfluid,\BAfluid,\BAsize$ in Section \ref{SS:InitialData}-\ref{SS:BA}, we obtain the desired result \eqref{ES:EnergyEst0}. 
\end{proof}
\begin{corollary}[Energy estimates on the null hypersurfaces for fluid variables]
	We assume the bootstrap assumptions in Section \ref{SS:BA}. Then, for solutions $(\lnenth,\ent,v^\mu)$ of \eqref{EQ:FOE2} and its variation $(\varh,\vars,\varv^\mu)$ with respect to operator $(\delta \Dfour_4)^{i}(\abs{u}\Dfour_3)^j(\abs{u}\Dfour_A)^k$ with $1\leq i+j+k\leq N$ with $k\geq 1$ defined in Definition \ref{DE:EnergyCurrent}, we have the following $L^2$ estimates hold on $\hypersur_{u}$ and $\hypersurbar_{\ubar}$:
	\begin{subequations}\label{CO:L2Fluid}
		\begin{align}
			\label{CO:L2Fluid1}\norm{\varh,\vars,\varv_A}_{L^2(\hypersurbar_{\ubar})}\lesssim&~\BAfluid^{\frac{5}{4}}+	\BAfluid^{\frac{5}{4}}a^{\frac{1}{4}}b^{-\frac{1}{4}},\\
			\label{CO:L2Fluid2}\norm{\varh,\vars,\varv_A}_{L^2(\hypersur_{u})}\lesssim&~\BAfluid^{\frac{5}{4}}b^{-\frac{1}{2}}+	\BAfluid^{\frac{5}{4}}a^{\frac{1}{4}}b^{-\frac{3}{4}},\\
			\label{CO:L2Fluid3}\norm{\varv_3}_{L^2(\hypersurbar_{\ubar})}\lesssim&~\BAfluid^{\frac{5}{4}}b^{-\frac{1}{2}}+	\BAfluid^{\frac{5}{4}}a^{\frac{1}{4}}b^{-\frac{3}{4}},\\
			\label{CO:L2Fluid4}\norm{\varv_4}_{L^2(\hypersur_{u})}\lesssim&~\BAfluid^{\frac{5}{4}}+	\BAfluid^{\frac{5}{4}}a^{\frac{1}{4}}b^{-\frac{1}{4}}.
		\end{align}
	\end{subequations}
\end{corollary}
\begin{proposition}[$L^2$-estimates for fluid variables on $\suubar$]\label{PR:L2Fluid}
	We assume the bootstrap assumptions in Section \ref{SS:BA}. Then, for solutions $(\lnenth,\ent,v^\mu)$ of \eqref{EQ:FOE2} and its variation $(\varh,\vars,\varv^\mu)$ with respect to operator $(\delta \Dfour_4)^{i}(\abs{u}\Dfour_3)^j(\abs{u}\Dfour_A)^k$ with $1\leq i+j+k\leq N-2$ or $1\leq i+j+k\leq N-1$ with $k\geq 1$, defined in Definition \ref{DE:EnergyCurrent}, we have the following $L^2$ estimates hold on $\suubar$:
	\begin{subequations}\label{ES:L2Fluid}
		\begin{align}
			\label{ES:L2Fluid1}\norm{\varh,\vars,\varv_A}_{L^2(\suubar)}\lesssim&~\BAfluid^{\frac{5}{4}}+	\BAfluid^{\frac{5}{4}}a^{\frac{1}{4}}b^{-\frac{1}{4}},\\
			\label{ES:L2Fluid2}\norm{\varv_3}_{L^2(\suubar)}\lesssim&~\BAfluid^{\frac{5}{4}}b^{-\frac{1}{2}}+	\BAfluid^{\frac{5}{4}}a^{\frac{1}{4}}b^{-\frac{3}{4}},\\
			\label{ES:L2Fluid3}\norm{\varv_4}_{L^2(\suubar)}\lesssim&~\BAfluid^{\frac{5}{4}}\delta^{-\frac{1}{2}}+	\BAfluid^{\frac{5}{4}}\delta^{-\frac{1}{2}}a^{\frac{1}{4}}b^{-\frac{1}{4}}.
		\end{align}
	\end{subequations}
	We note that \eqref{ES:L2Fluid} recovers and improves bootstrap assumption \eqref{BA32} by the choice of parameters $\delta,a,b$ in Section \ref{SS:Norms} and $\BAsize,\BAfluid$ in Section \ref{SS:BA}.
\end{proposition}
\begin{proof}[Proof of Proposition \ref{PR:L2Fluid}]
	We first prove \eqref{ES:L2Fluid1}. Recalling $\uPsi$ in Section \ref{SS:SchematicNotation}, we denote $\dot{\uPsi}\in\{\varh,\vars,\varv_A\}$. Substituting $f$ by $\norm{\uPsi}_{L^2(\suubar)}$, using the control of $\lapse$ in Proposition \ref{PR:lapse}, we integrate equation \eqref{EQ:LieuL} over $u$ to obtain:
	\begin{align}\label{ES:L2FluidProof1}
		\norm{\dot{\uPsi}}_{L^2(\suubar)}^2=&~\norm{\dot{\uPsi}}_{L^2(S_{\uinf,\ubar})}^2+\int_{\uinf}^{u}\int_{\suprimeubar}\lapse\dot{\uPsi}\Dfour_3\dot{\uPsi}+\lapse\gtr\uchi\dot{\uPsi}^2\diff\spherevol\diff\uprime\\
		\notag\lesssim&~\norm{\dot{\uPsi}}_{L^2(S_{\uinf,\ubar})}^2+\norm{\Dfour_3\dot{\uPsi}}_{L^2(\hypersurbar_{\ubar})}\norm{\dot{\uPsi}}_{L^2(\hypersurbar_{\ubar})}+\norm{\gtr\uchi}_{L^\infty(\hypersurbar_{\ubar})}\norm{\dot{\uPsi}}^2_{L^2(\hypersurbar_{\ubar})}.
	\end{align}
	Combining \eqref{ES:L2FluidProof1}, bootstrap assumptions \eqref{BA} and energy estimates \eqref{CO:L2Fluid}, we have:
	\begin{align}
		\norm{\dot{\uPsi}}_{L^2(\suubar)}\lesssim\norm{\dot{\uPsi}}_{L^2(S_{\uinf,\ubar})}+\norm{\Dfour_3\dot{\uPsi}}_{L^2(\hypersurbar_{\ubar})}\lesssim \BAfluid^{\frac{5}{4}}+	\BAfluid^{\frac{5}{4}}a^{\frac{1}{4}}b^{-\frac{1}{4}}.
	\end{align}
	Similarly, to prove \eqref{ES:L2Fluid2}, we substitute $f$ in \eqref{EQ:LieuL} by $\norm{\varv_3}_{L^2(\suubar)}$. Then, using the control of $\lapse$ in Proposition \ref{PR:lapse}, we integrate equation \eqref{EQ:LieuL} over $u$ to obtain:
	\begin{align}\label{ES:L2FluidProof2}
		\norm{\varv_3}_{L^2(\suubar)}^2=&~\norm{\varv_3}_{L^2(S_{\uinf,\ubar})}^2+\int_{\uinf}^{u}\int_{\suprimeubar}\lapse\varv_3\Dfour_3\varv_3+\lapse\gtr\uchi\varv_3^2\diff\spherevol\diff\uprime\\
		\notag\lesssim&~\norm{\varv_3}_{L^2(S_{\uinf,\ubar})}^2+\norm{\Dfour_3\varv_3}_{L^2(\hypersurbar_{\ubar})}\norm{\varv_3}_{L^2(\hypersurbar_{\ubar})}+\norm{\gtr\uchi}_{L^\infty(\hypersurbar_{\ubar})}\norm{\varv_3}^2_{L^2(\hypersurbar_{\ubar})}.
	\end{align}
	Combining \eqref{ES:L2FluidProof2}, bootstrap assumptions \eqref{BA} and energy estimates \eqref{CO:L2Fluid3}, we have:
	\begin{align}
		\norm{\varv_3}_{L^2(\suubar)}\lesssim\norm{\varv_3}_{L^2(S_{\uinf,\ubar})}+\norm{\Dfour_3\varv_3}_{L^2(\hypersurbar_{\ubar})}\lesssim\lesssim \BAfluid^{\frac{5}{4}}b^{-\frac{1}{2}}+	\BAfluid^{\frac{5}{4}}a^{\frac{1}{4}}b^{-\frac{3}{4}}.
	\end{align}
	Finally we prove \eqref{ES:L2Fluid3} by integrating \eqref{EQ:LieL} over $u$ and substituting $f$ by $\varv_4$:
	\begin{align}\label{ES:L2FluidProof3}
		\norm{\varv_4}_{L^2(\suubar)}^2=&~\norm{\varv_4}_{L^2(S_{u,0})}^2+\int_{0}^{\ubar}\int_{\suubarprime}\lapse\varv_4\Dfour_4\varv_4+\lapse\gtr\uchi\varv_4^2\diff\spherevol\diff\ubarprime\\
		\notag\lesssim&~\norm{\varv_4}_{L^2(S_{u,0})}^2+\sup_{\ubarprime\in[0,\delta]}\delta^{-\frac{1}{2}}\norm{\delta\Dfour_4\varv_4}_{L^2(\hypersur_{u})}\norm{\varv_4}_{L^2(\suubarprime)}+\sup_{\ubarprime\in[0,\delta]}\delta\norm{\gtr\upchi}_{L^\infty(\suubarprime)}\norm{\varv_4}^2_{L^2(\suubarprime)}.
	\end{align}
	Taking $\sup\limits_{\ubarprime\in[0,\delta]}$ of \eqref{ES:L2FluidProof3}, by bootstrap assumptions \eqref{BA} and energy estimates \eqref{CO:L2Fluid4}, we have:
	\begin{align}
		\sup_{\ubarprime\in[0,\delta]}\norm{\varv_4}_{L^2(\suubarprime)}\lesssim\norm{\varv_4}_{L^2(S_{u,0})}+\delta^{-\frac{1}{2}}\norm{\delta\Dfour_4\varv_4}_{L^2(\hypersur_{u})}\lesssim\delta^{-\frac{1}{2}}\left(\BAfluid^{\frac{5}{4}}+	\BAfluid^{\frac{5}{4}}\delta^{-\frac{1}{2}}a^{\frac{1}{4}}b^{-\frac{1}{4}}\right)\lesssim \BAfluid^{\frac{5}{4}}\delta^{-\frac{1}{2}}+	\BAfluid^{\frac{5}{4}}a^{\frac{1}{4}}b^{-\frac{1}{4}}.
	\end{align}
\end{proof}
\begin{corollary}[$L^\infty$-estimates for fluid variables on $\suubar$]\label{CO:LinftyFluid}
	We assume the bootstrap assumptions in Section \ref{SS:BA}. Then, for solutions $(\lnenth,\ent,v^\mu)$ of \eqref{EQ:FOE1} and its variation $(\varh,\vars,\varv^\mu)$ with respect to operator $(\delta \Dfour_4)^{i}(\abs{u}\Dfour_3)^j(\abs{u}\Dfour_A)^k$ with $1\leq i+j+k\leq N-3$ defined in Definition \ref{DE:EnergyCurrent}, we have the following $L^\infty$ estimates hold on $\suubar$:
	\begin{subequations}\label{ES:LinftyFluid1}
		\begin{align}
			\label{ES:LinftyFluid1.1}\norm{\varh,\vars,\varv_A}_{L^\infty(\suubar)}\lesssim&~\BAfluid^{\frac{5}{4}}+	\BAfluid^{\frac{5}{4}}a^{\frac{1}{4}}b^{-\frac{1}{4}},\\
			\label{ES:LinftyFluid1.2}\norm{\varv_3}_{L^\infty(\suubar)}\lesssim&~\BAfluid^{\frac{5}{4}}b^{-\frac{1}{2}}+	\BAfluid^{\frac{5}{4}}a^{\frac{1}{4}}b^{-\frac{3}{4}},\\
			\label{ES:LinftyFluid1.3}\norm{\varv_4}_{L^\infty(\suubar)}\lesssim&~\BAfluid^{\frac{5}{4}}\delta^{-\frac{1}{2}}+	\BAfluid^{\frac{5}{4}}\delta^{-\frac{1}{2}}a^{\frac{1}{4}}b^{-\frac{1}{4}}.
		\end{align}
	\end{subequations}
	Moreover, for fluid variables $(\lnenth,\ent,v_I)$, where $v_I:=\gfour(v,e_I)$ with $I=1,2,3,4$, we have:
	\begin{subequations}\label{ES:LinftyFluid2}
		\begin{align}
			\label{ES:LinftyFluid2.1}\norm{\lnenth,\ent,v_A}_{L^\infty(\suubar)}\lesssim&~\BAfluid^{\frac{5}{4}}+	\BAfluid^{\frac{5}{4}}a^{\frac{1}{4}}b^{-\frac{1}{4}},\\
			\label{ES:LinftyFluid2.2}b^{-\frac{1}{2}}\left(1-\BAfluid^{\frac{5}{4}}-	\BAfluid^{\frac{5}{4}}a^{\frac{1}{4}}b^{-\frac{1}{4}}\right)\lesssim\norm{v_3}_{L^\infty(\suubar)}\lesssim&~b^{-\frac{1}{2}}\left(1+\BAfluid^{\frac{5}{4}}+	\BAfluid^{\frac{5}{4}}a^{\frac{1}{4}}b^{-\frac{1}{4}}\right),\\
			\label{ES:LinftyFluid2.3}\delta^{-\frac{1}{2}}\left(1-\BAfluid^{\frac{5}{4}}-	\BAfluid^{\frac{5}{4}}a^{\frac{1}{4}}b^{-\frac{1}{4}}\right)\lesssim\norm{v_4}_{L^\infty(\suubar)}\lesssim&~\delta^{-\frac{1}{2}}\left(1+\BAfluid^{\frac{5}{4}}+	\BAfluid^{\frac{5}{4}}a^{\frac{1}{4}}b^{-\frac{1}{4}}\right).
		\end{align}
	\end{subequations}
	We note that \eqref{ES:LinftyFluid1} recovers and improves bootstrap assumption \eqref{BA32}, \eqref{ES:LinftyFluid2} recovers and improves bootstrap assumption \eqref{BA4} by the choice of parameters $N,\delta,a,b$ in Section \ref{SS:Norms} and $\BAsize,\BAfluid,\iniBAfluid$ in Section \ref{SS:BA}.
\end{corollary}
\begin{proof}[Proof of Corollary \ref{CO:LinftyFluid}]
	Using the control of metric $\gsphere$ in Proposition \ref{PR:gsphere} and Corollary \ref{CO:Area}, \eqref{ES:LinftyFluid1} are direct results of by Proposition \ref{PR:L2Fluid} and the Sobolev embedding $H^2\hookrightarrow L^\infty$ on $\suubar$. \eqref{ES:LinftyFluid2} follows from \eqref{ES:LinftyFluid1}, fundamental theorem of Calculus and initial data assumption \eqref{Data4}.
\end{proof}
\begin{remark}
	Here, we note that the RHS of \eqref{ES:LinftyFluid1} are $\abs{u}^2$ times of the RHS of \eqref{ES:L2Fluid} due to Sobolev embedding and the $\abs{u}$-weight on $\Dfour_A$. In this article, since $\abs{u}$ is bounded, the $\abs{u}$-weight on $\Dfour_A$ does not affect our estimates. Meanwhile, if one looks for the (mild) radiation/injection from past null infinity, then such $\abs{u}$-weight would cause a non-negligible effect.
\end{remark}

\begin{proposition}[Estimates for the Ricci curvature components]\label{PR:EstimatesforRicciCurvature}
	Let $\WangD^{\vec{m}}$ be as in Definition \ref{DE:operators} with $\vec{m}=(i,j,k)$. Ricci curvature components are defined in Definition \ref{DE:CurvatureComponents}. Under the bootstrap assumptions \eqref{BA}, we have the estimates for the derivatives of the Ricci curvature components:
	
	For $i+j+k\leq N$ with $k\geq 1$,
	\begin{subequations}\label{ES:EstimatesforRicciCurvature1}
		\begin{align}
			\sum\limits_{\abs{\vec{m}}\leq N}\norm{\WangD^{\vec{m}}\FAlphaAB}_{L^2(\hypersur_{u})}\lesssim&~  \delta^{-\frac{1}{2}},&
			\sum\limits_{\abs{\vec{m}}\leq N}\norm{\WangD^{\vec{m}}\FAlphabar}_{L^2(\hypersurbar_{\ubar})}\lesssim&~   b^{-1},\\
			\label{ES:EstimatesforRicciCurvature1.2}\sum\limits_{\abs{\vec{m}}\leq N}\norm{\WangD^{\vec{m}}\FAlphabar}_{L^2(\hypersur_{u})}\lesssim&~ \delta^{\frac{3}{2}},&
			\sum\limits_{\abs{\vec{m}}\leq N}\norm{\WangD^{\vec{m}}\FAlphaAB}_{L^2(\hypersurbar_{\ubar})}\lesssim&~   b,\\
			\label{ES:EstimatesforRicciCurvature1.3}\sum\limits_{\abs{\vec{m}}\leq N}\norm{\WangD^{\vec{m}}\FBetaA}_{L^2(\hypersur_{u})}\lesssim&~  \BAfluid^{\frac{5}{4}},&
			\sum\limits_{\abs{\vec{m}}\leq N}\norm{\WangD^{\vec{m}}\FBetabar}_{L^2(\hypersurbar_{\ubar})}\lesssim&~  \BAfluid^{\frac{5}{4}}b^{-\frac{1}{2}},\\
			\label{ES:EstimatesforRicciCurvature1.4}\sum\limits_{\abs{\vec{m}}\leq N}\norm{\WangD^{\vec{m}}\FBetabar}_{L^2(\hypersur_{u})}\lesssim&~  \BAfluid^{\frac{5}{4}}\delta,&
			\sum\limits_{\abs{\vec{m}}\leq N}\norm{\WangD^{\vec{m}}\FBetaA}_{L^2(\hypersurbar_{\ubar})}\lesssim&~   \BAfluid^{\frac{5}{4}}b^{\frac{1}{2}},\\
			\label{ES:EstimatesforRicciCurvature1.5}\sum\limits_{\abs{\vec{m}}\leq N}\norm{\WangD^{\vec{m}}(\FRho,\FRhodual,\Scalar)}_{L^2(\hypersur_{u})}\lesssim&~  b^{-\frac{1}{2}},&
			\sum\limits_{\abs{\vec{m}}\leq N}\norm{\WangD^{\vec{m}}(\FRho,\FRhodual,\Scalar)}_{L^2(\hypersurbar_{\ubar})}\lesssim&~  1.
		\end{align}
	\end{subequations}
	Moreover, on $\suubar$,	for $1\leq i+j+k\leq N-1$, we have:
	\begin{subequations}\label{ES:EstimatesforRicciCurvature2}
		\begin{align}
			\sum\limits_{1\leq\abs{\vec{m}}\leq N-1}\norm{\WangD^{\vec{m}}\FAlphaAB}_{L^2(\suubar)}\lesssim&~  \delta^{-1},&
			\sum\limits_{1\leq\abs{\vec{m}}\leq N-1}\norm{\WangD^{\vec{m}}\FAlphabar}_{L^2(\suubar)}\lesssim&~  b^{-1},\\
			\sum\limits_{1\leq\abs{\vec{m}}\leq N-1}\norm{\WangD^{\vec{m}}\FBetaA}_{L^2(\suubar)}\lesssim&~  \BAfluid^{\frac{5}{4}}\delta^{-\frac{1}{2}},&
			\sum\limits_{1\leq\abs{\vec{m}}\leq N-1}\norm{\WangD^{\vec{m}}\FBetabar}_{L^2(\suubar)}\lesssim&~  \BAfluid^{\frac{5}{4}}b^{-\frac{1}{2}},\\
			\sum\limits_{1\leq\abs{\vec{m}}\leq N-1}\norm{\WangD^{\vec{m}}(\FRho,\FRhodual,\Scalar)}_{L^2(\suubar)}\lesssim&~  1.
		\end{align}
	\end{subequations}
	
	In addition, for $1\leq i+j+k\leq N-3$, we have the following $L^\infty$-estimates on $\suubar$:
	\begin{subequations}\label{ES:EstimatesforRicciCurvature3}
		\begin{align}
			\sum\limits_{1\leq\abs{\vec{m}}\leq N-3}\norm{\WangD^{\vec{m}}\FAlphaAB}_{L^\infty(\suubar)}\lesssim&~   \delta^{-1},&
			\sum\limits_{1\leq\abs{\vec{m}}\leq N-3}\norm{\WangD^{\vec{m}}\FAlphabar}_{L^\infty(\suubar)}\lesssim&~   b^{-1},\\
			\sum\limits_{1\leq\abs{\vec{m}}\leq N-3}\norm{\WangD^{\vec{m}}\FBetaA}_{L^\infty(\suubar)}\lesssim&~  \BAfluid^{\frac{5}{4}}\delta^{-\frac{1}{2}},&
			\sum\limits_{1\leq\abs{\vec{m}}\leq N-3}\norm{\WangD^{\vec{m}}\FBetabar}_{L^\infty(\suubar)}\lesssim&~   \BAfluid^{\frac{5}{4}}b^{-\frac{1}{2}},\\
			\sum\limits_{1\leq\abs{\vec{m}}\leq N-3}\norm{\WangD^{\vec{m}}(\FRho,\FRhodual,\Scalar)}_{L^\infty(\suubar)}\lesssim&~  1.
		\end{align}
	\end{subequations}
	In particular, for $\abs{\vec{m}}=0$, we have:
	\begin{subequations}\label{ES:EstimatesforRicciCurvature4}
		\begin{align}
			\delta^{-1}(1- \BAfluid^{\frac{5}{4}})\lesssim\norm{\FAlphaAB}_{L^\infty(\suubar)}\lesssim&~ \delta^{-1}(1- \BAfluid^{\frac{5}{4}}),\\
			b^{-1}(1- \BAfluid^{\frac{5}{4}})\lesssim\norm{\FAlphabar}_{L^\infty(\suubar)}\lesssim&~ b^{-1}(1+ \BAfluid^{\frac{5}{4}}),\\
			\delta^{-\frac{1}{2}}(1- \BAfluid^{\frac{5}{4}})\lesssim\norm{\FBetaA}_{L^\infty(\suubar)}\lesssim&~ \delta^{-\frac{1}{2}}(1+ \BAfluid^{\frac{5}{4}}),\\
			b^{-\frac{1}{2}}(1- \BAfluid^{\frac{5}{4}})\lesssim\norm{\FBetabar}_{L^\infty(\suubar)}\lesssim&~ b^{-\frac{1}{2}}(1+ \BAfluid^{\frac{5}{4}}),\\
			\label{ES:EstimatesforRicciCurvature4.4}(1- \BAfluid^{\frac{5}{4}})\lesssim\norm{(\FRho,\FRhodual,\Scalar)}_{L^\infty(\suubar)}\lesssim&~ (1+ \BAfluid^{\frac{5}{4}}).
		\end{align}
	\end{subequations}
\end{proposition}
\begin{proof}[Discussion of the proof of Proposition \ref{PR:EstimatesforRicciCurvature}]
	Under the bootstrap assumptions \eqref{BA21}-\eqref{BA22}, the size of commutator terms are smaller; see Lemma \ref{LE:EstimatesofTensorDifference}. Then, simply by product rule and control of the commutator terms using Lemma \ref{LE:EstimatesofTensorDifference}, \eqref{ES:EstimatesforRicciCurvature1} is derived using \eqref{CO:L2Fluid} and Proposition \ref{PR:Ricci0}. Similarly, we obtain  \eqref{ES:EstimatesforRicciCurvature2} by \eqref{ES:L2Fluid} and \eqref{ES:EstimatesforRicciCurvature3} by \eqref{ES:LinftyFluid1}. In particular, for proving \eqref{ES:EstimatesforRicciCurvature1.2} and \eqref{ES:EstimatesforRicciCurvature1.4}, we used the bootstrap assumptions \eqref{BA4}(later improved by \eqref{ES:EstimatesforRicciCurvature4.4}) and the fact that $g(v,v)=-1$ such that, for example, $\FAlphabar=\frac{(\density+\pressure)\angv^4(\angv^2-1)^2}{v_4^2}$ on $\hypersur_{u}$.
	
	We note that, smallness $\BAfluid^{\frac{5}{4}}$ does not necessarily appears in \eqref{ES:EstimatesforRicciCurvature1}-\eqref{ES:EstimatesforRicciCurvature3}. For example,
	\begin{align}\label{EQ:D4FAlpha}
		\angD_4\FAlphaAB
		=&~\frac12\Dfour_4(\density+\pressure)v_4^2
		+(\density+\pressure)v_4\gfour(\Dfour_4v,e_4)
		-4\omega\FAlphaAB,
	\end{align}
	where the last term on the RHS of \eqref{EQ:D4FAlpha} does not enjoy a extra smallness $\BAfluid^{\frac{5}{4}}$. Nevertheless, taking $\angD_{ 3}$-derivatives would provide such smallness since $\uomega,\upeta$ are smaller (than $\BAfluid^{\frac{5}{4}}$). Hence \eqref{ES:EstimatesforRicciCurvature4} holds as one consider fundamental theorem of calculus along $e_3$-directions.
\end{proof}
\section{Below-top-order $L^2$-Estimates for Ricci Coefficients}
\begin{proposition}[Lower-order $L^2(\suubar)$-estimates for Ricci coefficients]\label{PR:L2RicciCoefficients}
	Let $\vec{m}=(i,j,k)$ with $\abs{\vec{m}}\leq N-2$ or $\abs{\vec{m}}= N-1$ and $k\geq 1$. $\tilde{\mathcal{O}}_{N,2},\mathcal{R}$ are defined in Section \ref{SS:Norms}. Let $\WangD^{\vec{m}}$ be as in Definition \ref{DE:operators} with $\vec{m}=(i,j,k)$. Under the bootstrap assumptions \eqref{BA}, we have the estimates\footnote{We note that, on the RHS of \eqref{ES:L2RicciCoefficients}-\eqref{ES:L2RicciCoefficientsImproved}, the first bound is derived from the gravity aspect, and the second bound is obtained from the fluid side. In particular, a better second bound holds since higher-order fluids enjoys a smallness $\BAfluid$ in the estimates in deriving the second bound. We do not emphasize this improvement because it is not needed in closing the bootstrap argument.} for the derivatives of the Ricci coefficients:
	\begin{subequations}\label{ES:L2RicciCoefficients}
		\begin{align}
			\label{ES:L2RicciCoefficients1}	\sum\limits_{\substack{\abs{\vec{m}}\leq N-2\\\abs{\vec{m}}= N-1,k\geq 1}}\norm{\WangD^{\vec{m}}\hchi}_{L^2(\suubar)}\lesssim&~ a^{\frac{1}{2}},\\
			\label{ES:L2RicciCoefficients2}	\sum\limits_{\substack{\abs{\vec{m}}\leq N-2\\\abs{\vec{m}}= N-1,k\geq 1}}\norm{\WangD^{\vec{m}}\gtr\upchi}_{L^2(\suubar)}\lesssim&~ \delta a+1,\\
			\label{ES:L2RicciCoefficients3}	\sum\limits_{\substack{\abs{\vec{m}}\leq N-2\\\abs{\vec{m}}= N-1,k\geq 1}}\norm{\WangD^{\vec{m}}\omega}_{L^2(\suubar)}\lesssim&~ \BAsize\delta a+1,\\
			\label{ES:L2RicciCoefficients4}	\sum\limits_{\substack{\abs{\vec{m}}\leq N-2\\\abs{\vec{m}}= N-1,k\geq 1}}\norm{\abs{u}\WangD^{\vec{m}}\upsi}_{L^2(\suubar)}\lesssim&~\frac{\delta a^{\frac{1}{2}}}{\abs{u}} \left(1+\tilde{\mathcal{O}}_{N,2}+\mathcal{R}\right).
		\end{align}
	\end{subequations}
	Additionally, we have:
	\begin{subequations}\label{ES:L2RicciCoefficientsImproved}
		\begin{align}
			\label{ES:L2RicciCoefficients5}	\sum\limits_{\substack{1\leq\abs{\vec{m}}\leq N-2\\\abs{\vec{m}}= N-1,k\geq 1}}\norm{\WangD^{\vec{m}}\left(\gtr\uchi+\frac{2}{\abs{u}}\right)}_{L^2(\suubar)}\lesssim&~\delta,\\
			\label{ES:L2RicciCoefficients6}	\sum\limits_{\substack{1\leq\abs{\vec{m}}\leq N-2\\\abs{\vec{m}}= N-1,k\geq 1}}\norm{\frac{1}{\abs{u}}\WangD^{\vec{m}}\ln\lapse}_{L^2(\suubar)}\lesssim&~ \BAsize\frac{\delta}{\abs{u}}.	
			\end{align}
		\end{subequations}
	\end{proposition}
	We prove Proposition \ref{PR:L2RicciCoefficients} in Section \ref{SS:proveRicci}. 
	
	\begin{corollary}[$L^\infty(\suubar)$-estimates for Ricci coefficients]\label{CO:LinftyRicciCoefficients}
		Let $\WangD^{\vec{m}}$ be as in Definition \ref{DE:operators} with $\vec{m}=(i,j,k)$. Under the bootstrap assumptions \eqref{BA}, we have the estimates for the derivatives of the Ricci coefficients:
		\begin{subequations}\label{ES:LinftyRicciCoefficients}
			\begin{align}
				\label{ES:LinftyRicciCoefficients1}	\sum\limits_{\abs{\vec{m}}\leq N-3}\norm{\WangD^{\vec{m}}\hchi}_{L^\infty(\suubar)}\lesssim&~ a^{\frac{1}{2}},\\
				\label{ES:LinftyRicciCoefficients2}	\sum\limits_{\abs{\vec{m}}\leq N-3}\norm{\WangD^{\vec{m}}\gtr\upchi}_{L^\infty(\suubar)}\lesssim&~ \delta a+1,\\
				\label{ES:LinftyRicciCoefficients3}	\sum\limits_{\abs{\vec{m}}\leq N-3}\norm{\WangD^{\vec{m}}\omega}_{L^\infty(\suubar)}\lesssim&~ \BAsize\delta a+1,\\
				\label{ES:LinftyRicciCoefficients4}	\sum\limits_{\abs{\vec{m}}\leq N-3}\norm{\abs{u}\WangD^{\vec{m}}\upsi}_{L^\infty(\suubar)}\lesssim&~\BAsize\frac{\delta a^{\frac{1}{2}}}{\abs{u}} \left(1+\tilde{\mathcal{O}}_{N,2}+\mathcal{R}\right).
			\end{align}
		\end{subequations}
		Additionally, we have:
		\begin{subequations}\label{ES:LinftyRicciCoefficientsMore}
			\begin{align}
				\label{ES:LinftyRicciCoefficients5}	\sum\limits_{1\leq\abs{\vec{m}}\leq N-3}\norm{\WangD^{\vec{m}}\left(\gtr\uchi+\frac{2}{\abs{u}}\right)}_{L^\infty(\suubar)}\lesssim&~ \delta,\\
				\label{ES:LinftyRicciCoefficients6}	\sum\limits_{1\leq\abs{\vec{m}}\leq N-3}\norm{\frac{1}{\abs{u}}\WangD^{\vec{m}}\ln\lapse}_{L^\infty(\suubar)}\lesssim&~ \frac{\delta}{\abs{u}}.
			\end{align}
		\end{subequations}
		We note that \eqref{ES:LinftyRicciCoefficients} recover and improve the $L^\infty$ bootstrap assumptions \eqref{BA21} since $N-3>\frac{N}{2}$.
	\end{corollary}
	\begin{proof}[Proof of Corollary \ref{CO:LinftyRicciCoefficients}]
		With the control of the induced metric $\gsphere$ in Proposition \ref{PR:gsphere} and Corollary \ref{CO:Area}, Corollary \ref{CO:LinftyRicciCoefficients} is a direct result of Proposition \ref{PR:L2RicciCoefficients} and Sobolev embedding $H^2\hookrightarrow L^\infty$ on $\suubar$.
	\end{proof}

	\subsection{Proof of Proposition \ref{PR:L2RicciCoefficients}}\label{SS:proveRicci}
	We note that, by Lemma \ref{LE:EstimatesofTensorDifference}, the size of tensorial components $\Dfour v_4,\Dfour v_3$ are the same as the size of derivatives of scalar functions $\Dfour (v_4),\Dfour(v_3)$. Then, we can use Proposition \ref{PR:L2Fluid} and Corollary \ref{CO:LinftyFluid} to estimate terms $\Dfour (v_4),\Dfour(v_3)$. In the rest of the article, we silently use Lemma \ref{LE:EstimatesofTensorDifference}.
	\begin{proof}[Proof of \eqref{ES:L2RicciCoefficients1}]
		We write the equation \eqref{EQ:angD3hchi} schematically as follows:
		\begin{align}\label{PV:hchi1}
			\angD_3\hchi+\frac{1}{2}\gtr\uchi\hchi=&~\angnabla\upeta+\psi\upsi+\upsi\upsi+\FRhodual,
		\end{align}
		where $\upsi,\psi,\angv$ are defined in Section \ref{SS:SchematicNotation}.
		Let $\vec{m}=(i,j,k)$ with $\abs{\vec{m}}\leq N-2$ or $\abs{\vec{m}}= N-1$ and $k\geq 1$. Commuting \eqref{PV:hchi1} with $\frac{1}{\abs{u}^k}\WangD^{\vec{m}}$, using \eqref{EQ:ang3ijk}-\eqref{EQ:ang3ijk2}, we obtain
		\begin{align}\label{PV:hchi2}
			&\frac{1}{\abs{u}^k}\angD_3\WangD^{\vec{m}}\hchi+\frac{1}{\abs{u}^k}\frac{k+1}{2}\gtr\uchi\cdot\WangD^{\vec{m}}\hchi\\
			\notag=&~\frac{1}{\abs{u}^k}\WangD^{\vec{m}}\angnabla\upeta+\sum\limits_{\vec{m}_1+\vec{m}_2=(i-i_0,j,k-k_0)}\frac{\delta^{i_0}}{\abs{u}^{k-k_0}}\WangD^{\vec{m}_1}\psi^{i_0+1}\cdot\WangD^{\vec{m}_2}\upsi^{k_0+1}+\sum\limits_{\vec{m}_1+\vec{m}_2=(i-i_0,j,k-k_0)}\frac{\delta^{i_0}}{\abs{u}^{k-k_0}}\WangD^{\vec{m}_1}\psi^{i_0}\cdot\WangD^{\vec{m}_2}\upsi^{k_0+2}\\
			\notag&+\sum\limits_{\vec{m}_1+\vec{m}_2=(i-i_0-1,j,k-k_0+1)}\frac{\delta^{i_0+1}}{\abs{u}^{k-k_0+1}}\WangD^{\vec{m}_1}\psi^{i_0}\cdot\WangD^{\vec{m}_2}\upsi^{k_0+2}\\
			\notag&+\sum\limits_{\vec{m}_1+\vec{m}_2+\vec{m}_3=(i-i_0,j,k-k_0)}\frac{\delta^{i_0}}{\abs{u}^{k-k_0}}\WangD^{\vec{m}_1}\psi^{i_0}\cdot\WangD^{\vec{m}_2}\upsi^{k_0}\cdot\WangD^{\vec{m}_3}\FRhodual\\
			\notag&+\sum\limits_{\vec{m}_1+\vec{m}_2+\vec{m}_3+\vec{m}_4+\vec{m}_5=(i-i_0,j,k-k_0-1)}\frac{\delta^{i_0}}{\abs{u}^{k-k_0-1}}\WangD^{\vec{m}_1}\psi^{i_0+1}\cdot\WangD^{\vec{m}_2}\upsi^{k_0}\cdot\WangD^{\vec{m}_3}\FBetabar\\
			\notag&+\sum\limits_{\vec{m}_1+\vec{m}_2+\vec{m}_3=(i-i_0,j,k-k_0-1)}\frac{\delta^{i_0}}{\abs{u}^{k-k_0-1}}\WangD^{\vec{m}_1}\psi^{i_0+1}\cdot\WangD^{\vec{m}_2}\upsi^{k_0}\cdot\WangD^{\vec{m}_3}\Betabar\\
			\notag&+\sum\limits_{\vec{m}_1+\vec{m}_2+\vec{m}_3=(i-i_0-1,j,k-k_0)}\frac{\delta^{i_0+1}}{\abs{u}^{k-k_0}}\WangD^{\vec{m}_1}\psi^{i_0+1}\cdot\WangD^{\vec{m}_2}\upsi^{k_0}\cdot\WangD^{\vec{m}_3}\Rhodual,\\
			\notag&=:\frac{1}{\abs{u}^k}\mathfrak{F}_{(\hchi)}.
		\end{align}
		We apply Proposition \ref{PR:transport2} with $\lambda_0:=\frac{k+1}{2}$ to \eqref{PV:hchi2}:
		\begin{align}\label{PV:hchi3}
			\norm{\WangD^{\vec{m}}\hchi}_{L^2(\suubar)}\lesssim&~\norm{\WangD^{\vec{m}}\hchi}_{L^2(S_{\uinf,\ubar})}+\int_{\uinf}^{u}\norm{\mathfrak{F}_{(\hchi)}}_{L^2(\suprimeubar)}\diff\uprime.
		\end{align}
		We now estimate $\int_{\uinf}^{u}\norm{\mathfrak{F}_{(\hchi)}}_{L^2(\suprimeubar)}\diff\uprime$. First, we consider the most dangerous term $\twonorms{\WangD^{\vec{m}}\angnabla\upeta}{u}{1}{\sphere}{2}{\hypersurbar_{\ubar}}$ when $\vec{m}=(i,j,k)$ with $i+j+k=N-1$ and $k\geq 1$:
		\begin{align}\label{PV:hchi4}
			\twonorms{\WangD^{\vec{m}}\angnabla\upeta}{u}{1}{\sphere}{2}{\hypersurbar_{\ubar}}\lesssim&~\norm{\frac{1}{\abs{u}}}_{L^2_u}\twonorms{\WangD^{\vec{m}_0}\upeta}{u}{2}{\sphere}{2}{\hypersurbar_{\ubar}}\\
			\notag\lesssim&~\frac{\delta^{\frac{1}{2}}a^{\frac{1}{2}}}{\abs{u}^{\frac{1}{2}}}\tilde{\mathcal{O}}_{N,2}\\
			\notag\lesssim&~\frac{\delta^{\frac{1}{2}}a^{\frac{1}{2}}}{\abs{u}^{\frac{1}{2}}}\BAsize,
		\end{align}
		where $\vec{m}_0:=(i,j,k+1)$. Then we consider the rest of the terms using bootstrap assumptions \eqref{BA}. By \eqref{ES:BoundsforProducts21},
		\begin{subequations}\label{PV:hchi5}
			\begin{align}
				\sum\limits_{\vec{m}_1+\vec{m}_2=(i-i_0,j,k-k_0)}\delta^{i_0}\twonorms{\abs{u}^{k_0}\WangD^{\vec{m}_1}\psi^{i_0+1}\cdot\WangD^{\vec{m}_2}\upsi^{k_0+1}}{u}{1}{\sphere}{2}{\hypersurbar_{\ubar}}\lesssim&~\left(\BAsize a^{\frac{1}{2}}\right)\cdot\left(\BAsize\frac{\delta a^{\frac{1}{2}}}{\abs{u}}\right)^{i_0+k_0+1},\\
				\sum\limits_{\vec{m}_1+\vec{m}_2=(i-i_0,j,k-k_0)}\delta^{i_0}\twonorms{\abs{u}^{k_0}\WangD^{\vec{m}_1}\psi^{i_0}\cdot\WangD^{\vec{m}_2}\upsi^{k_0+2}}{u}{1}{\sphere}{2}{\hypersurbar_{\ubar}}\lesssim&~\left(\BAsize\frac{\delta a^{\frac{1}{2}}}{\abs{u}}\right)^{i_0+k_0+2}.
			\end{align}
		\end{subequations}
		By bootstrap assumptions \eqref{BA}, Proposition \ref{PR:EstimatesforRicciCurvature} and \eqref{ES:BoundsforProducts2}, we have the following estimates for fluid terms:
		
		\begin{align}\label{PV:hchi6}
			&\sum\limits_{\vec{m}_1+\vec{m}_2+\vec{m}_3=(i-i_0,j,k-k_0)}\delta^{i_0}\twonorms{\abs{u}^{k_0}\WangD^{\vec{m}_1}\psi^{i_0}\cdot\WangD^{\vec{m}_2}\upsi^{k_0}\cdot\WangD^{\vec{m}_3}\FRhodual}{u}{1}{\sphere}{2}{\hypersurbar_{\ubar}}\\
			\notag\lesssim&~\sum\limits_{\substack{\vec{m}_1+\vec{m}_2+\vec{m}_3=(i-i_0,j,k-k_0)\\\abs{\vec{m}_1+\vec{m}_2}>\frac{N-1}{2}}}\delta^{i_0}\cdot\left(\twonorms{\abs{u}^{k_0}\WangD^{\vec{m}_1}\psi^{i_0}\cdot\WangD^{\vec{m}_2}\upsi^{k_0}}{u}{1}{\sphere}{2}{\hypersurbar_{\ubar}}\cdot\twonorms{\WangD^{\vec{m}_3}\FRhodual}{u}{\infty}{\sphere}{\infty}{\hypersurbar_{\ubar}}\right)\\
			\notag&+\sum\limits_{\substack{\vec{m}_1+\vec{m}_2+\vec{m}_3=(i-i_0,j,k-k_0)\\\abs{\vec{m}_1+\vec{m}_2}\leq\frac{N-1}{2}}}\delta^{i_0}\left(\twonorms{\abs{u}^{k_0}\WangD^{\vec{m}_1}\psi^{i_0}\cdot\WangD^{\vec{m}_2}\upsi^{k_0}}{u}{\infty}{\sphere}{\infty}{\hypersurbar_{\ubar}}\cdot\twonorms{\WangD^{\vec{m}_3}\FRhodual}{u}{1}{\sphere}{2}{\hypersurbar_{\ubar}}\right)\\
			\notag\lesssim&~\left(\BAsize\frac{\delta a^{\frac{1}{2}}}{\abs{u}}\right)^{i_0+k_0}\cdot\abs{u}.
		\end{align}
		\begin{align}\label{PV:hchi6.1}
			&\sum\limits_{\vec{m}_1+\vec{m}_2+\vec{m}_3=(i-i_0,j,k-k_0-1)}\delta^{i_0}\twonorms{\abs{u}^{k_0+1}\WangD^{\vec{m}_1}\psi^{i_0+1}\cdot\WangD^{\vec{m}_2}\upsi^{k_0}\cdot\WangD^{\vec{m}_3}\FBetabar}{u}{1}{\sphere}{2}{\hypersurbar_{\ubar}}\\
			\notag\lesssim&~\sum\limits_{\substack{\vec{m}_1+\vec{m}_2+\vec{m}_3=(i-i_0,j,k-k_0-1)\\\abs{\vec{m}_1+\vec{m}_2}>\frac{N-1}{2}}}\delta^{i_0}\cdot\left(\twonorms{\abs{u}^{k_0+1}\WangD^{\vec{m}_1}\psi^{i_0+1}\cdot\WangD^{\vec{m}_2}\upsi^{k_0}}{u}{1}{\sphere}{2}{\hypersurbar_{\ubar}}\cdot\twonorms{\WangD^{\vec{m}_3}\FBetabar}{u}{\infty}{\sphere}{\infty}{\hypersurbar_{\ubar}}\right)\\
			\notag&+\sum\limits_{\substack{\vec{m}_1+\vec{m}_2+\vec{m}_3=(i-i_0,j,k-k_0-1)\\\abs{\vec{m}_1+\vec{m}_2}\leq\frac{N-1}{2}}}\delta^{i_0}\cdot\left(\twonorms{\abs{u}^{k_0+1}\WangD^{\vec{m}_1}\psi^{i_0+1}\cdot\WangD^{\vec{m}_2}\upsi^{k_0}}{u}{\infty}{\sphere}{\infty}{\hypersurbar_{\ubar}}\cdot\twonorms{\WangD^{\vec{m}_3}\FBetabar}{u}{1}{\sphere}{2}{\hypersurbar_{\ubar}}\right)\\
			\notag\lesssim&~a^{\frac{1}{2}}b^{-\frac{1}{2}}\cdot\left(\BAsize\frac{\delta a^{\frac{1}{2}}}{\abs{u}}\right)^{i_0+k_0}.
		\end{align}
		By bootstrap assumptions \eqref{BA} and \eqref{ES:BoundsforProducts2}, we have the following estimates for curvature terms:
		
		\begin{align}\label{PV:hchi7.0}
			&\sum\limits_{\vec{m}_1+\vec{m}_2+\vec{m}_3=(i-i_0,j,k-k_0-1)}\delta^{i_0}\twonorms{\abs{u}^{k_0+1}\WangD^{\vec{m}_1}\psi^{i_0+1}\cdot\WangD^{\vec{m}_2}\upsi^{k_0}\cdot\WangD^{\vec{m}_3}\Betabar}{u}{1}{\sphere}{2}{\hypersurbar_{\ubar}}\\
			\notag\lesssim&~\sum\limits_{\substack{\vec{m}_1+\vec{m}_2+\vec{m}_3=(i-i_0,j,k-k_0-1)\\\abs{\vec{m}_1+\vec{m}_2}>\frac{N-1}{2}}}\delta^{i_0}\twonorms{\abs{u}^{k_0+1}\WangD^{\vec{m}_1}\psi^{i_0+1}\cdot\WangD^{\vec{m}_2}\upsi^{k_0}}{u}{1}{\sphere}{2}{\hypersurbar_{\ubar}}\cdot\twonorms{\WangD^{\vec{m}_3}\Betabar}{u}{\infty}{\sphere}{\infty}{\hypersurbar_{\ubar}}\\
			\notag&+\sum\limits_{\substack{\vec{m}_1+\vec{m}_2+\vec{m}_3=(i-i_0,j,k-k_0-1)\\\abs{\vec{m}_1+\vec{m}_2}\leq\frac{N-1}{2}}}\delta^{i_0}\twonorms{\abs{u}^{k_0+1}\WangD^{\vec{m}_1}\psi^{i_0+1}\cdot\WangD^{\vec{m}_2}\upsi^{k_0}}{u}{\infty}{\sphere}{\infty}{\hypersurbar_{\ubar}}\cdot\twonorms{\WangD^{\vec{m}_3}\Betabar}{u}{1}{\sphere}{2}{\hypersurbar_{\ubar}}\\
			\notag\lesssim&~\BAsize a^{\frac{1}{2}}\cdot\left(\BAsize\frac{\delta a^{\frac{1}{2}}}{\abs{u}}\right)^{i_0+k_0+1}.
		\end{align}
		\begin{align}\label{PV:hchi7}
			&\sum\limits_{\vec{m}_1+\vec{m}_2+\vec{m}_3=(i-i_0-1,j,k-k_0)}\delta^{i_0+1}\twonorms{\abs{u}^{k_0}\WangD^{\vec{m}_1}\psi^{i_0+1}\cdot\WangD^{\vec{m}_2}\upsi^{k_0}\cdot\WangD^{\vec{m}_3}\Rhodual}{u}{1}{\sphere}{2}{\hypersurbar_{\ubar}}\\
			\notag\lesssim&~\sum\limits_{\substack{\vec{m}_1+\vec{m}_2+\vec{m}_3=(i-i_0-1,j,k-k_0)\\\abs{\vec{m}_1+\vec{m}_2}>\frac{N-1}{2}}}\delta^{i_0+1}\twonorms{\abs{u}^{k_0}\WangD^{\vec{m}_1}\psi^{i_0+1}\cdot\WangD^{\vec{m}_2}\upsi^{k_0}}{u}{1}{\sphere}{2}{\hypersurbar_{\ubar}}\cdot\twonorms{\WangD^{\vec{m}_3}\Rhodual}{u}{\infty}{\sphere}{\infty}{\hypersurbar_{\ubar}}\\
			\notag&+\sum\limits_{\substack{\vec{m}_1+\vec{m}_2+\vec{m}_3=(i-i_0-1,j,k-k_0)\\\abs{\vec{m}_1+\vec{m}_2}\leq\frac{N-1}{2}}}\delta^{i_0+1}\twonorms{\abs{u}^{k_0}\WangD^{\vec{m}_1}\psi^{i_0+1}\cdot\WangD^{\vec{m}_2}\upsi^{k_0}}{u}{\infty}{\sphere}{\infty}{\hypersurbar_{\ubar}}\cdot\twonorms{\WangD^{\vec{m}_3}\Rhodual}{u}{1}{\sphere}{2}{\hypersurbar_{\ubar}}\\
			\notag\lesssim&~\BAsize a^{\frac{1}{2}}\left(\BAsize\frac{\delta a^{\frac{1}{2}}}{\abs{u}}\right)^{i_0+k_0+2}.
		\end{align}
		Combining \eqref{PV:hchi4}-\eqref{PV:hchi7} and the initial data assumption \eqref{Data}, by the choice of parameter $\delta,a,b$ in Section \ref{SS:Norms} and $\BAfluid,\BAsize$ in Section \eqref{SS:BA}, we obtain the desired result \eqref{ES:L2RicciCoefficients1}.
	\end{proof}
	\begin{proof}[Proof of \eqref{ES:L2RicciCoefficients2}]
		We write the equation \eqref{EQ:e4upchi} schematically:
		\begin{align}\label{PV:tracechi1}
			\angD_4\gtr\upchi=&~(\gtr\upchi)^2+\sgabs{\hchi}^2+\omega\gtr\upchi+\FAlphaAB,
		\end{align}
		Let $\vec{m}=(i,j,k)$ with $\abs{\vec{m}}\leq N-2$ or $\abs{\vec{m}}= N-1$ and $k\geq 1$. Commuting \eqref{PV:tracechi1} with $\WangD^{\vec{m}}$, using \eqref{EQ:ang4ijk}, we obtain
		\begin{align}\label{PV:tracechi2}
			& \angD_4\WangD^{\vec{m}}\gtr\upchi\\
			\notag=&~ \sum\limits_{\vec{m}_1+\vec{m}_2+\vec{m}_3=(i,j-j_0,k-k_0)}\abs{u}^{j_0+k_0}\WangD^{\vec{m}_1}\upsi^{j_0+k_0}\cdot\WangD^{\vec{m}_2}\hchi\cdot\WangD^{\vec{m}_3}\hchi\\
			\notag&+ \sum\limits_{\vec{m}_1+\vec{m}_2+\vec{m}_3=(i,j-j_0,k-k_0)}\abs{u}^{j_0+k_0}\WangD^{\vec{m}_1}\upsi^{j_0+k_0}\cdot\WangD^{\vec{m}_2}\gtr\upchi\cdot\WangD^{\vec{m}_3}\psi\\
			\notag&+ \sum\limits_{\vec{m}_1+\vec{m}_2+\vec{m}_3=(i,j-j_0,k-k_0)}\abs{u}^{j_0+k_0}\WangD^{\vec{m}_1}\upsi^{j_0+k_0+1}\cdot\WangD^{\vec{m}_2}\gtr\upchi\\
			\notag&+ \sum\limits_{\vec{m}_1+\vec{m}_2=(i,j-j_0,k-k_0)}\abs{u}^{j_0+k_0}\WangD^{\vec{m}_1}\upsi^{k_0+j_0}\cdot\WangD^{\vec{m}_2}\FAlphaAB\\
			\notag&+ \sum\limits_{\vec{m}_1+\vec{m}_2+\vec{m}_3=(i,j-j_0,k-k_0-1)}\abs{u}^{j_0+k_0+2}\WangD^{\vec{m}_1}\upsi^{k_0+j_0}\cdot\WangD^{\vec{m}_2}\gtr\upchi\cdot\WangD^{\vec{m}_3}\FBetaA\\
			\notag&+ \sum\limits_{\vec{m}_1+\vec{m}_2+\vec{m}_3=(i,j-j_0,k-k_0-1)}\abs{u}^{j_0+k_0+1}\WangD^{\vec{m}_1}\upsi^{j_0+k_0}\cdot\WangD^{\vec{m}_2}\gtr\upchi\cdot\WangD^{\vec{m}_3}\BetaA\\
			\notag&+ \sum\limits_{\vec{m}_1+\vec{m}_2+\vec{m}_3=(i,j-j_0-1,k-k_0)}\abs{u}^{j_0+k_0+1}\WangD^{\vec{m}_1}\upsi^{j_0+k_0}\cdot\WangD^{\vec{m}_2}\gtr\upchi\cdot\WangD^{\vec{m}_3}\Rhodual,\\
			\notag&=:\mathfrak{F}_{(\gtr\upchi)}.
		\end{align}
		By Proposition \ref{PR:transport1}, we have:
		\begin{align}\label{PV:tracechi3}
			\norm{\WangD^{\vec{m}}\gtr\upchi}_{L^2(\suubar)}\lesssim\norm{\WangD^{\vec{m}}\gtr\upchi}_{L^2(S_{u,0})}+\int_{0}^{\ubar}	\norm{\mathfrak{F}_{(\gtr\upchi)}}_{L^2(\suubarprime)}\diff\ubarprime.
		\end{align}
		We now estimate $\int_{0}^{\ubar}\norm{\mathfrak{F}_{(\gtr\upchi)}}_{L^2(\suubarprime)}\diff\ubarprime$. By proved result \eqref{ES:L2RicciCoefficients1}, \eqref{ES:LinftyRicciCoefficients1}, and bootstrap assumptions \eqref{BA}, for the first term on the RHS of \eqref{PV:tracechi2}:
		\begin{align}\label{PV:tracechi4}
			&\sum\limits_{\vec{m}_1+\vec{m}_2+\vec{m}_3=(i,j-j_0,k-k_0)}\abs{u}^{j_0+k_0}\twonorms{\WangD^{\vec{m}_1}\upsi^{j_0+k_0}\cdot\WangD^{\vec{m}_2}\hchi\cdot\WangD^{\vec{m}_3}\hchi}{\ubar}{1}{\sphere}{2}{\hypersur_{u}}\\
			\notag\lesssim&~\sum\limits_{\substack{\vec{m}_1+\vec{m}_2+\vec{m}_3=(i,j-j_0,k-k_0)\\\abs{\vec{m}_1}>\frac{N-1}{2}}}\delta\abs{u}^{j_0+k_0}\twonorms{\WangD^{\vec{m}_1}\upsi^{j_0+k_0}}{\ubar}{\infty}{\sphere}{2}{\hypersur_{u}}\cdot\twonorms{\WangD^{\vec{m}_2}\hchi}{\ubar}{\infty}{\sphere}{\infty}{\hypersur_{u}}\cdot\twonorms{\WangD^{\vec{m}_3}\hchi}{\ubar}{\infty}{\sphere}{\infty}{\hypersur_{u}}\\
			\notag&+\sum\limits_{\substack{\vec{m}_1+\vec{m}_2+\vec{m}_3=(i,j-j_0,k-k_0)\\\abs{\vec{m}_1}\leq\frac{N-1}{2}}}\delta\abs{u}^{j_0+k_0}\twonorms{\WangD^{\vec{m}_1}\upsi^{j_0+k_0}}{\ubar}{\infty}{\sphere}{\infty}{\hypersur_{u}}\cdot\twonorms{\WangD^{\vec{m}_2}\hchi}{\ubar}{\infty}{\sphere}{2}{\hypersur_{u}}\cdot\twonorms{\WangD^{\vec{m}_3}\hchi}{\ubar}{\infty}{\sphere}{\infty}{\hypersur_{u}}\\
			\notag\lesssim&~ \frac{\delta a}{\abs{u}}\left(\BAsize\frac{\delta a^{\frac{1}{2}}}{\abs{u}}\right)^{i_0+k_0}.
		\end{align}
		By bootstrap assumptions \eqref{BA} and Sobolev embedding $H^2\hookrightarrow L^\infty$ on $\suubar$ for $\WangD^{\vec{m}}\gtr\upchi$, for the second term on the RHS of \eqref{PV:tracechi2}, we have:
		\begin{align}\label{PV:tracechi5}
			&\sum\limits_{\vec{m}_1+\vec{m}_2+\vec{m}_3=(i,j-j_0,k-k_0)}\abs{u}^{j_0+k_0}\twonorms{\WangD^{\vec{m}_1}\upsi^{j_0+k_0}\cdot\WangD^{\vec{m}_2}\gtr\upchi\cdot\WangD^{\vec{m}_3}\psi}{\ubar}{1}{\sphere}{2}{\hypersur_{u}}\\
			\notag\lesssim&~\sum\limits_{\substack{\vec{m}_1+\vec{m}_2+\vec{m}_3=(i,j-j_0,k-k_0)\\\abs{\vec{m}_1}>\frac{N-1}{2}}}\delta\abs{u}^{j_0+k_0}\twonorms{\WangD^{\vec{m}_1}\upsi^{j_0+k_0}}{\ubar}{\infty}{\sphere}{2}{\hypersur_{u}}\cdot\twonorms{\WangD^{\vec{m}_2}\gtr\upchi}{\ubar}{\infty}{\sphere}{\infty}{\hypersur_{u}}\cdot\twonorms{\WangD^{\vec{m}_3}\psi}{\ubar}{\infty}{\sphere}{\infty}{\hypersur_{u}}\\
			\notag&+\sum\limits_{\substack{\vec{m}_1+\vec{m}_2+\vec{m}_3=(i,j-j_0,k-k_0)\\\abs{\vec{m}_2}>\frac{N-1}{2}}}\delta\abs{u}^{j_0+k_0}\twonorms{\WangD^{\vec{m}_1}\upsi^{j_0+k_0}}{\ubar}{\infty}{\sphere}{\infty}{\hypersur_{u}}\cdot\twonorms{\WangD^{\vec{m}_2}\gtr\upchi}{\ubar}{\infty}{\sphere}{2}{\hypersur_{u}}\cdot\twonorms{\WangD^{\vec{m}_3}\psi}{\ubar}{\infty}{\sphere}{\infty}{\hypersur_{u}}\\
			\notag&+\sum\limits_{\substack{\vec{m}_1+\vec{m}_2+\vec{m}_3=(i,j-j_0,k-k_0)\\\abs{\vec{m}_3}>\frac{N-1}{2}}}\delta\abs{u}^{j_0+k_0}\twonorms{\WangD^{\vec{m}_1}\upsi^{j_0+k_0}}{\ubar}{\infty}{\sphere}{\infty}{\hypersur_{u}}\cdot\twonorms{\WangD^{\vec{m}_2}\gtr\upchi}{\ubar}{\infty}{\sphere}{\infty}{\hypersur_{u}}\cdot\twonorms{\WangD^{\vec{m}_3}\psi}{\ubar}{\infty}{\sphere}{2}{\hypersur_{u}}\\
			\notag\lesssim&~ \left(\BAsize\frac{\delta a^{\frac{1}{2}}}{\abs{u}}\right)^{j_0+k_0+1}\cdot\sum\limits_{\abs{\vec{m}_2}\leq N-1}\twonorms{\WangD^{\vec{m_2}}\gtr\upchi}{\ubar}{\infty}{\sphere}{2}{\hypersur_{u}}.
		\end{align}
		By bootstrap assumptions \eqref{BA}, Proposition \ref{PR:EstimatesforRicciCurvature}, we estimate the fluid terms on the RHS of \eqref{PV:tracechi2}:
		\begin{align}\label{PV:tracechi6}
			&\sum\limits_{\vec{m}_1+\vec{m}_2=(i,j-j_0,k-k_0)}\abs{u}^{j_0+k_0}\twonorms{\WangD^{\vec{m}_1}\upsi^{j_0+k_0}\cdot\WangD^{\vec{m}_2}\FAlphaAB}{\ubar}{1}{\sphere}{2}{\hypersur_{u}}\\
			\notag\lesssim&~\sum\limits_{\substack{\vec{m}_1+\vec{m}_2=(i,j-j_0,k-k_0)\\\abs{\vec{m}_1}>\frac{N-1}{2}}}\delta\abs{u}^{j_0+k_0}\left(\twonorms{\WangD^{\vec{m}_1}\upsi^{j_0+k_0}}{\ubar}{\infty}{\sphere}{2}{\hypersur_{u}}\cdot\twonorms{\WangD^{\vec{m}_2}\FAlphaAB}{\ubar}{\infty}{\sphere}{\infty}{\hypersur_{u}}\right)\\
			\notag&+\sum\limits_{\substack{\vec{m}_1+\vec{m}_2=(i,j-j_0,k-k_0)\\\abs{\vec{m}_2}>\frac{N-1}{2}}}\delta\abs{u}^{j_0+k_0}\left(\twonorms{\WangD^{\vec{m}_1}\upsi^{j_0+k_0}}{\ubar}{\infty}{\sphere}{\infty}{\hypersur_{u}}\cdot\twonorms{\WangD^{\vec{m}_2}\FAlphaAB}{\ubar}{\infty}{\sphere}{2}{\hypersur_{u}}\right)\\
			\notag\lesssim&~  \left(\BAsize\frac{\delta a^{\frac{1}{2}}}{\abs{u}}\right)^{j_0+k_0},
		\end{align}
		In below, we also need the Sobolev embedding $H^2\hookrightarrow L^\infty$ on $\suubar$ for $\WangD^{\vec{m}}\gtr\upchi$ to obtain the estimate for the other fluid term:
		\begin{align}\label{PV:tracechi7}
			&\sum\limits_{\vec{m}_1+\vec{m}_2+\vec{m}_3=(i,j-j_0,k-k_0-1)}\abs{u}^{j_0+k_0+1}\twonorms{\WangD^{\vec{m}_1}\upsi^{k_0+j_0}\cdot\WangD^{\vec{m}_2}\gtr\upchi\cdot\WangD^{\vec{m}_3}\FBetaA}{\ubar}{1}{\sphere}{2}{\hypersur_{u}}\\
			\notag\lesssim&~\sum\limits_{\substack{\vec{m}_1+\vec{m}_2+\vec{m}_3=(i,j-j_0,k-k_0-1)\\\abs{\vec{m}_1}>\frac{N-1}{2}}}\delta\abs{u}^{j_0+k_0+1}\left(\twonorms{\WangD^{\vec{m}_1}\upsi^{j_0+k_0}}{\ubar}{\infty}{\sphere}{2}{\hypersur_{u}}\cdot\twonorms{\WangD^{\vec{m}_2}\gtr\upchi}{\ubar}{\infty}{\sphere}{\infty}{\hypersur_{u}}\cdot\twonorms{\WangD^{\vec{m}_3}\FBetaA}{\ubar}{\infty}{\sphere}{\infty}{\hypersur_{u}}\right)\\
			\notag&+\sum\limits_{\substack{\vec{m}_1+\vec{m}_2+\vec{m}_3=(i,j-j_0,k-k_0-1)\\\abs{\vec{m}_2}>\frac{N-1}{2}}}\delta\abs{u}^{j_0+k_0+1}\left(\twonorms{\WangD^{\vec{m}_1}\upsi^{j_0+k_0}}{\ubar}{\infty}{\sphere}{\infty}{\hypersur_{u}}\cdot\twonorms{\WangD^{\vec{m}_2}\gtr\upchi}{\ubar}{\infty}{\sphere}{2}{\hypersur_{u}}\cdot\twonorms{\WangD^{\vec{m}_3}\FBetaA}{\ubar}{\infty}{\sphere}{\infty}{\hypersur_{u}}\right)\\
			\notag&+\sum\limits_{\substack{\vec{m}_1+\vec{m}_2+\vec{m}_3=(i,j-j_0,k-k_0-1)\\\abs{\vec{m}_3}>\frac{N-1}{2}}}\delta\abs{u}^{j_0+k_0+1}\left(\twonorms{\WangD^{\vec{m}_1}\upsi^{j_0+k_0}}{\ubar}{\infty}{\sphere}{\infty}{\hypersur_{u}}\cdot\twonorms{\WangD^{\vec{m}_2}\gtr\upchi}{\ubar}{\infty}{\sphere}{\infty}{\hypersur_{u}}\cdot\twonorms{\WangD^{\vec{m}_3}\FBetaA}{\ubar}{\infty}{\sphere}{2}{\hypersur_{u}}\right)\\
			\notag\lesssim&~ \delta^{\frac{1}{2}}\cdot\left(\BAsize\frac{\delta a^{\frac{1}{2}}}{\abs{u}}\right)^{j_0+k_0}\cdot\sum\limits_{\abs{\vec{m}_2}\leq N-1}\twonorms{\WangD^{\vec{m_2}}\gtr\upchi}{\ubar}{\infty}{\sphere}{2}{\hypersur_{u}}.
		\end{align}
		Finally, by bootstrap assumptions \eqref{BA} and Sobolev embedding $H^2\hookrightarrow L^\infty$ on $\suubar$ for $\WangD^{\vec{m}}\gtr\upchi$, we estimate the curvature terms on the RHS of \eqref{PV:tracechi2}:
		\begin{align}\label{PV:tracechi8}
			&\sum\limits_{\vec{m}_1+\vec{m}_2+\vec{m}_3=(i,j-j_0,k-k_0-1)}\abs{u}^{j_0+k_0+1}\twonorms{\WangD^{\vec{m}_1}\upsi^{j_0+k_0}\cdot\WangD^{\vec{m}_2}\gtr\upchi\cdot\WangD^{\vec{m}_3}\BetaA}{\ubar}{1}{\sphere}{2}{\hypersur_{u}}\\
			\notag\lesssim&~\sum\limits_{\substack{\vec{m}_1+\vec{m}_2+\vec{m}_3=(i,j-j_0,k-k_0-1)\\\abs{\vec{m}_1}>\frac{N-1}{2}}}\delta\abs{u}^{j_0+k_0+1}\twonorms{\WangD^{\vec{m}_1}\upsi^{j_0+k_0}}{\ubar}{\infty}{\sphere}{2}{\hypersur_{u}}\cdot\twonorms{\WangD^{\vec{m}_2}\gtr\upchi}{\ubar}{\infty}{\sphere}{\infty}{\hypersur_{u}}\cdot\twonorms{\WangD^{\vec{m}_3}\BetaA}{\ubar}{\infty}{\sphere}{\infty}{\hypersur_{u}}\\
			\notag&+\sum\limits_{\substack{\vec{m}_1+\vec{m}_2+\vec{m}_3=(i,j-j_0,k-k_0-1)\\\abs{\vec{m}_2}>\frac{N-1}{2}}}\delta\abs{u}^{j_0+k_0+1}\twonorms{\WangD^{\vec{m}_1}\upsi^{j_0+k_0}}{\ubar}{\infty}{\sphere}{\infty}{\hypersur_{u}}\cdot\twonorms{\WangD^{\vec{m}_2}\gtr\upchi}{\ubar}{\infty}{\sphere}{2}{\hypersur_{u}}\cdot\twonorms{\WangD^{\vec{m}_3}\BetaA}{\ubar}{\infty}{\sphere}{\infty}{\hypersur_{u}}\\
			\notag&+\sum\limits_{\substack{\vec{m}_1+\vec{m}_2+\vec{m}_3=(i,j-j_0,k-k_0-1)\\\abs{\vec{m}_3}>\frac{N-1}{2}}}\delta\abs{u}^{j_0+k_0+1}\twonorms{\WangD^{\vec{m}_1}\upsi^{j_0+k_0}}{\ubar}{\infty}{\sphere}{\infty}{\hypersur_{u}}\cdot\twonorms{\WangD^{\vec{m}_2}\gtr\upchi}{\ubar}{\infty}{\sphere}{\infty}{\hypersur_{u}}\cdot\twonorms{\WangD^{\vec{m}_3}\BetaA}{\ubar}{\infty}{\sphere}{2}{\hypersur_{u}}\\
			\notag\lesssim&~ \left(\BAsize\frac{\delta a^{\frac{1}{2}}}{\abs{u}}\right)^{j_0+k_0+1}\cdot\sum\limits_{\abs{\vec{m}_2}\leq N-1}\twonorms{\WangD^{\vec{m_2}}\gtr\upchi}{\ubar}{\infty}{\sphere}{2}{\hypersur_{u}},
		\end{align}
		\begin{align}\label{PV:tracechi9}
			&\sum\limits_{\vec{m}_1+\vec{m}_2+\vec{m}_3=(i,j-j_0-1,k-k_0)}\abs{u}^{j_0+k_0+1}\twonorms{\WangD^{\vec{m}_1}\upsi^{j_0+k_0}\cdot\WangD^{\vec{m}_2}\gtr\upchi\cdot\WangD^{\vec{m}_3}\Rhodual}{\ubar}{1}{\sphere}{2}{\hypersur_{u}}\\
			\notag\lesssim&~\sum\limits_{\substack{\vec{m}_1+\vec{m}_2+\vec{m}_3=(i,j-j_0-1,k-k_0)\\\abs{\vec{m}_1}>\frac{N-1}{2}}}\delta\abs{u}^{j_0+k_0+1}\twonorms{\WangD^{\vec{m}_1}\upsi^{j_0+k_0}}{\ubar}{\infty}{\sphere}{2}{\hypersur_{u}}\cdot\twonorms{\WangD^{\vec{m}_2}\gtr\upchi}{\ubar}{\infty}{\sphere}{\infty}{\hypersur_{u}}\cdot\twonorms{\WangD^{\vec{m}_3}\Rhodual}{\ubar}{\infty}{\sphere}{\infty}{\hypersur_{u}}\\
			\notag&+\sum\limits_{\substack{\vec{m}_1+\vec{m}_2+\vec{m}_3=(i,j-j_0,k-k_0)\\\abs{\vec{m}_2}>\frac{N-1}{2}}}\delta\abs{u}^{j_0+k_0+1}\twonorms{\WangD^{\vec{m}_1}\upsi^{j_0+k_0}}{\ubar}{\infty}{\sphere}{\infty}{\hypersur_{u}}\cdot\twonorms{\WangD^{\vec{m}_2}\gtr\upchi}{\ubar}{\infty}{\sphere}{2}{\hypersur_{u}}\cdot\twonorms{\WangD^{\vec{m}_3}\Rhodual}{\ubar}{\infty}{\sphere}{\infty}{\hypersur_{u}}\\
			\notag&+\sum\limits_{\substack{\vec{m}_1+\vec{m}_2+\vec{m}_3=(i,j-j_0-1,k-k_0)\\\abs{\vec{m}_3}>\frac{N-1}{2}}}\delta\abs{u}^{j_0+k_0+1}\twonorms{\WangD^{\vec{m}_1}\upsi^{j_0+k_0}}{\ubar}{\infty}{\sphere}{\infty}{\hypersur_{u}}\cdot\twonorms{\WangD^{\vec{m}_2}\gtr\upchi}{\ubar}{\infty}{\sphere}{\infty}{\hypersur_{u}}\cdot\twonorms{\WangD^{\vec{m}_3}\Rhodual}{\ubar}{\infty}{\sphere}{2}{\hypersur_{u}}\\
			\notag\lesssim&~ \left(\BAsize\frac{\delta a^{\frac{1}{2}}}{\abs{u}}\right)^{j_0+k_0+2}\cdot\sum\limits_{\abs{\vec{m}_2}\leq N-1}\twonorms{\WangD^{\vec{m_2}}\gtr\upchi}{\ubar}{\infty}{\sphere}{2}{\hypersur_{u}}.
		\end{align}
		Combining \eqref{PV:tracechi3}-\eqref{PV:tracechi9}, we soak the $\sum\limits_{\abs{\vec{m}_2}\leq N-1}\twonorms{\WangD^{\vec{m_2}}\gtr\upchi}{\ubar}{\infty}{\sphere}{2}{\hypersur_{u}}$ terms to the left since their coefficients are $\ll 1$. Then, we obtain the desired estimates \eqref{ES:L2RicciCoefficients2} by the initial data assumption \eqref{Data} and the choice of parameter $\delta,a,b$ in Section \ref{SS:Norms}, $\BAfluid,\BAsize$ in Section \eqref{SS:BA}.
	\end{proof}
	\begin{proof}[Proof of \eqref{ES:L2RicciCoefficients3}]
		We write the equation \eqref{Lbaromega} schematically as follows:
		\begin{align}\label{PV:omega1}
			\angD_3\omega=&~\omega\upsi+\upsi\upsi+\Rho+\FRho+\Scalar,
		\end{align}
		where $\upsi,\psi$ are defined in Section \ref{SS:SchematicNotation}.
		Let $\vec{m}=(i,j,k)$ with $\abs{\vec{m}}\leq N-2$ or $\abs{\vec{m}}= N-1$ and $k\geq 1$. Commuting \eqref{PV:omega1} with $\frac{1}{\abs{u}^{k-1}}\WangD^{\vec{m}}$, using \eqref{EQ:ang3ijk}-\eqref{EQ:ang3ijk2}, we obtain
		\begin{align}\label{PV:omega2}
			&\frac{1}{\abs{u}^{k-1}}\angD_3\WangD^{\vec{m}}\omega+\frac{1}{\abs{u}^{k-1}}\frac{k}{2}\gtr\uchi\cdot\WangD^{\vec{m}}\omega\\
			\notag=&~\sum\limits_{\vec{m}_1+\vec{m}_2=(i-i_0,j,k-k_0)}\frac{\delta^{i_0}}{\abs{u}^{k-k_0-1}}\WangD^{\vec{m}_1}\omega^{i_0+1}\cdot\WangD^{\vec{m}_2}\upsi^{k_0+1}+\sum\limits_{\vec{m}_1+\vec{m}_2=(i-i_0,j,k-k_0)}\frac{\delta^{i_0}}{\abs{u}^{k-k_0-1}}\WangD^{\vec{m}_1}\omega^{i_0}\cdot\WangD^{\vec{m}_2}\upsi^{k_0+2}\\
			\notag&+\sum\limits_{\vec{m}_1+\vec{m}_2=(i-i_0-1,j,k-k_0+1)}\frac{\delta^{i_0+1}}{\abs{u}^{k-k_0}}\WangD^{\vec{m}_1}\omega^{i_0}\cdot\WangD^{\vec{m}_2}\upsi^{k_0+2}\\
			\notag&+\sum\limits_{\vec{m}_1+\vec{m}_2+\vec{m}_3=(i-i_0,j,k-k_0)}\frac{\delta^{i_0}}{\abs{u}^{k-k_0-1}}\WangD^{\vec{m}_1}\omega^{i_0}\cdot\WangD^{\vec{m}_2}\upsi^{k_0}\cdot\WangD^{\vec{m}_3}\Scalar\\
			\notag&+\sum\limits_{\vec{m}_1+\vec{m}_2+\vec{m}_3=(i-i_0,j,k-k_0)}\frac{\delta^{i_0}}{\abs{u}^{k-k_0-1}}\WangD^{\vec{m}_1}\omega^{i_0}\cdot\WangD^{\vec{m}_2}\upsi^{k_0}\cdot\WangD^{\vec{m}_3}\FRho\\
			\notag&+\sum\limits_{\vec{m}_1+\vec{m}_2+\vec{m}_3=(i-i_0,j,k-k_0-1)}\frac{\delta^{i_0+1}}{\abs{u}^{k-k_0-2}}\WangD^{\vec{m}_1}\omega^{i_0+1}\cdot\WangD^{\vec{m}_2}\upsi^{k_0}\cdot\WangD^{\vec{m}_3}\FBetabar\\
			\notag&+\sum\limits_{\vec{m}_1+\vec{m}_2+\vec{m}_3=(i-i_0,j,k-k_0)}\frac{\delta^{i_0}}{\abs{u}^{k-k_0-1}}\WangD^{\vec{m}_1}\omega^{i_0}\cdot\WangD^{\vec{m}_2}\upsi^{k_0}\cdot\WangD^{\vec{m}_3}\Rho\\
			\notag&+\sum\limits_{\vec{m}_1+\vec{m}_2+\vec{m}_3=(i-i_0,j,k-k_0-1)}\frac{\delta^{i_0}}{\abs{u}^{k-k_0-2}}\WangD^{\vec{m}_1}\omega^{i_0+1}\cdot\WangD^{\vec{m}_2}\upsi^{k_0}\cdot\WangD^{\vec{m}_3}\Betabar\\
			\notag&+\sum\limits_{\vec{m}_1+\vec{m}_2+\vec{m}_3=(i-i_0-1,j,k-k_0)}\frac{\delta^{i_0+1}}{\abs{u}^{k-k_0-1}}\WangD^{\vec{m}_1}\omega^{i_0+1}\cdot\WangD^{\vec{m}_2}\upsi^{k_0}\cdot\WangD^{\vec{m}_3}\Rhodual,\\
			\notag&=:\frac{1}{\abs{u}^{k-1}}\mathfrak{F}_{(\omega)}.
		\end{align}
		We apply Proposition \ref{PR:transport2} with $\lambda_0:=\frac{k}{2}$ to \eqref{PV:omega2}:
		\begin{align}\label{PV:omega3}
			\norm{\WangD^{\vec{m}}\omega}_{L^2(\suubar)}\lesssim&~\norm{\WangD^{\vec{m}}\omega}_{L^2(S_{\uinf,\ubar})}+\int_{\uinf}^{u}\norm{\mathfrak{F}_{(\omega)}}_{L^2(\suprimeubar)}\diff\uprime.
		\end{align}
		We now estimate $\int_{\uinf}^{u}\norm{\mathfrak{F}_{(\omega)}}_{L^2(\suprimeubar)}\diff\uprime$. Using bootstrap assumptions \eqref{BA}, Proposition \ref{PR:EstimatesforRicciCurvature}, we estimate the first two fluid terms as follows: 
		\begin{align}\label{PV:omega4}
			&\sum\limits_{\vec{m}_1+\vec{m}_2+\vec{m}_3=(i-i_0,j,k-k_0)}\delta^{i_0}\twonorms{\abs{u}^{k_0}\WangD^{\vec{m}_1}\psi^{i_0}\cdot\WangD^{\vec{m}_2}\upsi^{k_0}\cdot\WangD^{\vec{m}_3}\Scalar}{u}{1}{\sphere}{2}{\hypersurbar_{\ubar}}\\
			\notag\lesssim&~\sum\limits_{\substack{\vec{m}_1+\vec{m}_2+\vec{m}_3=(i-i_0,j,k-k_0)\\\abs{\vec{m}_1+\vec{m}_2}>\frac{N-1}{2}}}\delta^{i_0}\cdot\left(\twonorms{\abs{u}^{k_0}\WangD^{\vec{m}_1}\psi^{i_0}\cdot\WangD^{\vec{m}_2}\upsi^{k_0}}{u}{1}{\sphere}{2}{\hypersurbar_{\ubar}}\cdot\twonorms{\WangD^{\vec{m}_3}\Scalar}{u}{\infty}{\sphere}{\infty}{\hypersurbar_{\ubar}}\right)\\
			\notag&+\sum\limits_{\substack{\vec{m}_1+\vec{m}_2+\vec{m}_3=(i-i_0,j,k-k_0)\\\abs{\vec{m}_1+\vec{m}_2}\leq\frac{N-1}{2}}}\delta^{i_0}\left(\twonorms{\abs{u}^{k_0}\WangD^{\vec{m}_1}\psi^{i_0}\cdot\WangD^{\vec{m}_2}\upsi^{k_0}}{u}{\infty}{\sphere}{\infty}{\hypersurbar_{\ubar}}\cdot\twonorms{\WangD^{\vec{m}_3}\Scalar}{u}{1}{\sphere}{2}{\hypersurbar_{\ubar}}\right)\\
			\notag\lesssim&~ \left(\BAsize\frac{\delta a^{\frac{1}{2}}}{\abs{u}}\right)^{i_0+k_0}\cdot\abs{u}.
		\end{align}
		\begin{align}\label{PV:omega5}
			&\sum\limits_{\vec{m}_1+\vec{m}_2+\vec{m}_3=(i-i_0,j,k-k_0)}\delta^{i_0}\twonorms{\abs{u}^{k_0}\WangD^{\vec{m}_1}\psi^{i_0}\cdot\WangD^{\vec{m}_2}\upsi^{k_0}\cdot\WangD^{\vec{m}_3}\FRho}{u}{1}{\sphere}{2}{\hypersurbar_{\ubar}}\\
			\notag\lesssim&~\sum\limits_{\substack{\vec{m}_1+\vec{m}_2+\vec{m}_3=(i-i_0,j,k-k_0)\\\abs{\vec{m}_1+\vec{m}_2}>\frac{N-1}{2}}}\delta^{i_0}\cdot\left(\twonorms{\abs{u}^{k_0}\WangD^{\vec{m}_1}\psi^{i_0}\cdot\WangD^{\vec{m}_2}\upsi^{k_0}}{u}{1}{\sphere}{2}{\hypersurbar_{\ubar}}\cdot\twonorms{\WangD^{\vec{m}_3}\FRho}{u}{\infty}{\sphere}{\infty}{\hypersurbar_{\ubar}}\right)\\
			\notag&+\sum\limits_{\substack{\vec{m}_1+\vec{m}_2+\vec{m}_3=(i-i_0,j,k-k_0)\\\abs{\vec{m}_1+\vec{m}_2}\leq\frac{N-1}{2}}}\delta^{i_0}\left(\twonorms{\abs{u}^{k_0}\WangD^{\vec{m}_1}\psi^{i_0}\cdot\WangD^{\vec{m}_2}\upsi^{k_0}}{u}{\infty}{\sphere}{\infty}{\hypersurbar_{\ubar}}\cdot\twonorms{\WangD^{\vec{m}_3}\FRho}{u}{1}{\sphere}{2}{\hypersurbar_{\ubar}}\right)\\
			\notag\lesssim&~ \left(\BAsize\frac{\delta a^{\frac{1}{2}}}{\abs{u}}\right)^{i_0+k_0}\cdot\abs{u}.
		\end{align}
		Now, we estimate the third last term on the RHS of \eqref{PV:omega2}:
		\begin{align}\label{PV:omega6}
			&\sum\limits_{\vec{m}_1+\vec{m}_2+\vec{m}_3=(i-i_0,j,k-k_0)}\delta^{i_0}\twonorms{\abs{u}^{k_0}\WangD^{\vec{m}_1}\psi^{i_0}\cdot\WangD^{\vec{m}_2}\upsi^{k_0}\cdot\WangD^{\vec{m}_3}\Rho}{u}{1}{\sphere}{2}{\hypersurbar_{\ubar}}\\
			\notag\lesssim&~\sum\limits_{\substack{\vec{m}_1+\vec{m}_2+\vec{m}_3=(i-i_0,j,k-k_0)\\\abs{\vec{m}_1+\vec{m}_2}>\frac{N-1}{2}}}\delta^{i_0}\twonorms{\abs{u}^{k_0}\WangD^{\vec{m}_1}\psi^{i_0}\cdot\WangD^{\vec{m}_2}\upsi^{k_0}}{u}{1}{\sphere}{2}{\hypersurbar_{\ubar}}\cdot\twonorms{\WangD^{\vec{m}_3}\Rho}{u}{\infty}{\sphere}{\infty}{\hypersurbar_{\ubar}}\\
			\notag&+\sum\limits_{\substack{\vec{m}_1+\vec{m}_2+\vec{m}_3=(i-i_0,j,k-k_0)\\\abs{\vec{m}_1+\vec{m}_2}\leq\frac{N-1}{2}}}\delta^{i_0}\twonorms{\abs{u}^{k_0}\WangD^{\vec{m}_1}\psi^{i_0}\cdot\WangD^{\vec{m}_2}\upsi^{k_0}}{u}{\infty}{\sphere}{\infty}{\hypersurbar_{\ubar}}\cdot\twonorms{\WangD^{\vec{m}_3}\Rho}{u}{1}{\sphere}{2}{\hypersurbar_{\ubar}}\\
			\notag\lesssim&~\BAsize\frac{\delta a}{\abs{u}^3}\left(\BAsize\frac{\delta a^{\frac{1}{2}}}{\abs{u}}\right)^{i_0+k_0}.
		\end{align}
		Now, we estimate the rest of the terms on the RHS of \eqref{PV:omega2} by revisiting the proof of \eqref{ES:L2RicciCoefficients1}. In particular, we can derive an improved estimates compared to \eqref{PV:hchi5} and \eqref{PV:hchi7.0}-\eqref{PV:hchi7}.
		First of all, compared to \eqref{PV:hchi6.1}, we estimate the first term on the RHS of \eqref{PV:omega2} using bootstrap assumptions \eqref{BA} and Sobolev embedding $H^2\hookrightarrow L^\infty$ on $\suubar$ for $\WangD^{\vec{m}}\omega$:
		\begin{align}\label{PV:omega7}
			&\sum\limits_{\vec{m}_1+\vec{m}_2=(i-i_0,j,k-k_0)}\delta^{i_0}\twonorms{\abs{u}^{k_0}\WangD^{\vec{m}_1}\omega^{i_0+1}\cdot\WangD^{\vec{m}_2}\upsi^{k_0+1}}{u}{1}{\sphere}{2}{\hypersurbar_{\ubar}}\\
			\notag\lesssim&~\sum\limits_{\substack{\vec{m}_1+\vec{m}_2=(i-i_0,j,k-k_0)\\\abs{\vec{m}_1}>\frac{N-1}{2}}}\delta^{i_0}\twonorms{\abs{u}^{k_0}\WangD^{\vec{m}_1}\omega^{i_0+1}}{u}{1}{\sphere}{2}{\hypersurbar_{\ubar}}\cdot\twonorms{\WangD^{\vec{m}_2}\upsi^{k_0+1}}{u}{\infty}{\sphere}{\infty}{\hypersurbar_{\ubar}}\\
			\notag&+\sum\limits_{\substack{\vec{m}_1+\vec{m}_2=(i-i_0,j,k-k_0)\\\abs{\vec{m}_2}>\frac{N-1}{2}}}\delta^{i_0}\twonorms{\abs{u}^{k_0}\WangD^{\vec{m}_1}\omega^{i_0+1}}{u}{\infty}{\sphere}{\infty}{\hypersurbar_{\ubar}}\cdot\twonorms{\WangD^{\vec{m}_2}\upsi^{k_0+1}}{u}{1}{\sphere}{2}{\hypersurbar_{\ubar}}\\
			\notag\lesssim&~\left(\BAsize\frac{\delta a^{\frac{1}{2}}}{\abs{u}}\right)^{i_0+k_0+1}\cdot\sum\limits_{\abs{\vec{m}_1}\leq N-1}\twonorms{\WangD^{\vec{m}_1}\omega}{u}{\infty}{\sphere}{2}{\hypersurbar_{\ubar}}.
		\end{align}
		Then, compared to \eqref{PV:hchi7.0}-\eqref{PV:hchi7}, we estimate the last two terms on the RHS of \eqref{PV:omega2} using bootstrap assumptions \eqref{BA} and Sobolev embedding $H^2\hookrightarrow L^\infty$ on $\suubar$ for $\WangD^{\vec{m}}\omega$:
		\begin{align}\label{PV:omega8}
			&\sum\limits_{\vec{m}_1+\vec{m}_2+\vec{m}_3=(i-i_0,j,k-k_0-1)}\delta^{i_0}\twonorms{\abs{u}^{k_0+1}\WangD^{\vec{m}_1}\omega^{i_0+1}\cdot\WangD^{\vec{m}_2}\upsi^{k_0}\cdot\WangD^{\vec{m}_3}\Betabar}{u}{1}{\sphere}{2}{\hypersurbar_{\ubar}}\\
			\notag\lesssim&~\sum\limits_{\substack{\vec{m}_1+\vec{m}_2+\vec{m}_3=(i-i_0,j,k-k_0-1)\\\abs{\vec{m}_1}>\frac{N-1}{2}}}\delta^{i_0}\twonorms{\abs{u}^{k_0+1}\WangD^{\vec{m}_1}\omega^{i_0+1}}{u}{1}{\sphere}{2}{\hypersurbar_{\ubar}}\cdot\twonorms{\WangD^{\vec{m}_2}\upsi^{k_0}}{u}{\infty}{\sphere}{\infty}{\hypersurbar_{\ubar}}\cdot\twonorms{\WangD^{\vec{m}_3}\Betabar}{u}{\infty}{\sphere}{\infty}{\hypersurbar_{\ubar}}\\
			\notag&+\sum\limits_{\substack{\vec{m}_1+\vec{m}_2+\vec{m}_3=(i-i_0,j,k-k_0-1)\\\abs{\vec{m}_2}>\frac{N-1}{2}}}\delta^{i_0}\twonorms{\abs{u}^{k_0+1}\WangD^{\vec{m}_1}\omega^{i_0+1}}{u}{\infty}{\sphere}{\infty}{\hypersurbar_{\ubar}}\cdot\twonorms{\WangD^{\vec{m}_2}\upsi^{k_0}}{u}{1}{\sphere}{2}{\hypersurbar_{\ubar}}\cdot\twonorms{\WangD^{\vec{m}_3}\Betabar}{u}{\infty}{\sphere}{\infty}{\hypersurbar_{\ubar}}\\
			\notag&+\sum\limits_{\substack{\vec{m}_1+\vec{m}_2+\vec{m}_3=(i-i_0,j,k-k_0-1)\\\abs{\vec{m}_3}>\frac{N-1}{2}}}\delta^{i_0}\twonorms{\abs{u}^{k_0+1}\WangD^{\vec{m}_1}\omega^{i_0+1}}{u}{\infty}{\sphere}{\infty}{\hypersurbar_{\ubar}}\cdot\twonorms{\WangD^{\vec{m}_2}\upsi^{k_0}}{u}{\infty}{\sphere}{\infty}{\hypersurbar_{\ubar}}\cdot\twonorms{\WangD^{\vec{m}_3}\Betabar}{u}{1}{\sphere}{2}{\hypersurbar_{\ubar}}\\
			\notag\lesssim&~\left(\BAsize\frac{\delta a^{\frac{1}{2}}}{\abs{u}}\right)^{i_0+k_0+1}\cdot\sum\limits_{\abs{\vec{m}_1}\leq N-1}\twonorms{\WangD^{\vec{m}_1}\omega}{u}{\infty}{\sphere}{2}{\hypersurbar_{\ubar}}.
		\end{align}
		\begin{align}\label{PV:omega9}
			&\sum\limits_{\vec{m}_1+\vec{m}_2+\vec{m}_3=(i-i_0-1,j,k-k_0)}\delta^{i_0+1}\twonorms{\abs{u}^{k_0}\WangD^{\vec{m}_1}\omega^{i_0+1}\cdot\WangD^{\vec{m}_2}\upsi^{k_0}\cdot\WangD^{\vec{m}_3}\Rhodual}{u}{1}{\sphere}{2}{\hypersurbar_{\ubar}}\\
			\notag\lesssim&~\sum\limits_{\substack{\vec{m}_1+\vec{m}_2+\vec{m}_3=(i-i_0-1,j,k-k_0)\\\abs{\vec{m}_1>\frac{N-1}{2}}}}\delta^{i_0+1}\twonorms{\abs{u}^{k_0}\WangD^{\vec{m}_1}\omega^{i_0+1}}{u}{1}{\sphere}{2}{\hypersurbar_{\ubar}}\cdot\twonorms{\WangD^{\vec{m}_2}\upsi^{k_0}}{u}{\infty}{\sphere}{\infty}{\hypersurbar_{\ubar}}\cdot\twonorms{\WangD^{\vec{m}_3}\Rhodual}{u}{\infty}{\sphere}{\infty}{\hypersurbar_{\ubar}}\\
			\notag&+\sum\limits_{\substack{\vec{m}_1+\vec{m}_2+\vec{m}_3=(i-i_0-1,j,k-k_0)\\\abs{\vec{m}_2>\frac{N-1}{2}}}}\delta^{i_0+1}\twonorms{\abs{u}^{k_0}\WangD^{\vec{m}_1}\omega^{i_0+1}}{u}{\infty}{\sphere}{\infty}{\hypersurbar_{\ubar}}\cdot\twonorms{\WangD^{\vec{m}_2}\upsi^{k_0}}{u}{1}{\sphere}{2}{\hypersurbar_{\ubar}}\cdot\twonorms{\WangD^{\vec{m}_3}\Rhodual}{u}{\infty}{\sphere}{\infty}{\hypersurbar_{\ubar}}\\
			\notag&+\sum\limits_{\substack{\vec{m}_1+\vec{m}_2+\vec{m}_3=(i-i_0-1,j,k-k_0)\\\abs{\vec{m}_3>\frac{N-1}{2}}}}\delta^{i_0+1}\twonorms{\abs{u}^{k_0}\WangD^{\vec{m}_1}\omega^{i_0+1}}{u}{\infty}{\sphere}{\infty}{\hypersurbar_{\ubar}}\cdot\twonorms{\WangD^{\vec{m}_2}\upsi^{k_0}}{u}{\infty}{\sphere}{\infty}{\hypersurbar_{\ubar}}\cdot\twonorms{\WangD^{\vec{m}_3}\Rhodual}{u}{1}{\sphere}{2}{\hypersurbar_{\ubar}}\\
			\notag\lesssim&~\left(\BAsize\frac{\delta a^{\frac{1}{2}}}{\abs{u}}\right)^{i_0+k_0+2}\cdot\sum\limits_{\abs{\vec{m}_1}\leq N-1}\twonorms{\WangD^{\vec{m}_1}\omega}{u}{\infty}{\sphere}{2}{\hypersurbar_{\ubar}}.
		\end{align}
		Combining \eqref{PV:hchi6.1} and \eqref{PV:omega2}-\eqref{PV:omega9}, soaking the $\sum\limits_{\abs{\vec{m}_1}\leq N-1}\twonorms{\WangD^{\vec{m}_1}\omega}{u}{\infty}{\sphere}{2}{\hypersurbar_{\ubar}}$ terms by the LHS since their coefficients are $\ll 1$, we obtain the desired result \eqref{ES:L2RicciCoefficients3} by the initial data assumptions \ref{Data}, and the choice of parameter $\delta,a,b$ in Section \ref{SS:Norms}, $\BAfluid,\BAsize$ in Section \eqref{SS:BA}.
	\end{proof}
	\begin{proof}[Proof of \eqref{ES:L2RicciCoefficients4}]
		We recall $\upsi\in\{\upeta,\ueta,\huchi,\gtr\uchi+\frac{2}{\abs{u}},\uomega\}$ and $\unpsi\in\{\upeta,\huchi,\gtr\uchi+\frac{2}{\abs{u}},\uomega\}$ in Section \ref{SS:SchematicNotation}. By Proposition \ref{PR:NullStructure1}, we write the schematic transport equation as follows:
		\begin{align}\label{PV:upsi1}
			\angD_{ 4}\upsi=\angnabla(\omega,\ueta)+\psi\upsi+\upsi\upsi+\frac{1}{u}\psi+\Rho+\BetaA+\FRho+\FRhodual+\FBetaA+\Scalar.
		\end{align}
		In particular, we emphasize that $\angnabla\omega$ on the RHS of \eqref{PV:upsi1} appears only when $\upsi=\ueta$  (specifically, $\upsi=\upzeta$ or $\upsi=\angnabla\ln\lapse$).
		Let $\vec{m}=(i,j,k)$ with $\abs{\vec{m}}\leq N-2$ or $\abs{\vec{m}}= N-1$ and $k\geq 1$. Commuting \eqref{PV:tracechi1} with $\WangD^{\vec{m}}$, using \eqref{EQ:ang4ijk}, we obtain
		\begin{align}\label{PV:upsi2}
			& \angD_4\WangD^{\vec{m}}\upsi\\
			\notag=&~ \WangD^{\vec{m}}\angnabla(\omega,\ueta)+ \sum\limits_{\substack{\vec{m}_1+\vec{m}_2=(i,j-j_0,k-k_0)\\\abs{\vec{m}_2}\leq N-2}}\abs{u}^{j_0+k_0}\WangD^{\vec{m}_1}\upsi^{j_0+k_0}\cdot\WangD^{\vec{m}_2}\angnabla(\omega,\ueta)\\
			\notag&+ \sum\limits_{\vec{m}_1+\vec{m}_2=(i,j-j_0,k-k_0)}\abs{u}^{j_0+k_0}\WangD^{\vec{m}_1}\upsi^{j_0+k_0+1}\cdot\WangD^{\vec{m}_2}\psi+ \sum\limits_{\vec{m}_1=(i,j-j_0,k-k_0)}\abs{u}^{j_0+k_0}\WangD^{\vec{m}_1}\upsi^{j_0+k_0+2}\\
			\notag&+ \sum\limits_{\vec{m}_1+\vec{m}_2=(i,j-j_0,k-k_0)}\abs{u}^{j_0+k_0-1}\WangD^{\vec{m}_1}\upsi^{j_0+k_0}\cdot\WangD^{\vec{m}_2}\psi\\
			\notag&+ \sum\limits_{\vec{m}_1+\vec{m}_2=(i,j-j_0,k-k_0)}\abs{u}^{j_0+k_0}\WangD^{\vec{m}_1}\upsi^{k_0+j_0}\cdot\WangD^{\vec{m}_2}(\FRho,\FRhodual,\Scalar)\\
			\notag&+ \sum\limits_{\vec{m}_1+\vec{m}_2=(i,j-j_0,k-k_0)}\abs{u}^{j_0+k_0}\WangD^{\vec{m}_1}\upsi^{k_0+j_0}\cdot\WangD^{\vec{m}_2}\FBetaA\\
			\notag&+ \sum\limits_{\vec{m}_1+\vec{m}_2=(i,j-j_0-1,k-k_0)}\abs{u}^{j_0+k_0+1}\WangD^{\vec{m}_1}\upsi^{j_0+k_0+1}\cdot\WangD^{\vec{m}_2}\Rhodual,\\
			\notag&+ \sum\limits_{\vec{m}_1+\vec{m}_2=(i,j-j_0,k-k_0)}\abs{u}^{j_0+k_0}\WangD^{\vec{m}_1}\upsi^{j_0+k_0}\cdot\WangD^{\vec{m}_2}\Rho\\
			\notag&+ \sum\limits_{\vec{m}_1+\vec{m}_2=(i,j-j_0,k-k_0)}\abs{u}^{j_0+k_0}\WangD^{\vec{m}_1}\upsi^{j_0+k_0}\cdot\WangD^{\vec{m}_2}\BetaA\\
			\notag&=:\mathfrak{F}_{(\upsi)}.
		\end{align}
		By Proposition \ref{PR:transport1}, we have:
		\begin{align}\label{PV:upsi3}
			\norm{\WangD^{\vec{m}}\upsi}_{L^2(\suubar)}\lesssim\norm{\WangD^{\vec{m}}\upsi}_{L^2(S_{u,0})}+\int_{0}^{\ubar}	\norm{\mathfrak{F}_{(\upsi)}}_{L^2(\suubarprime)}\diff\ubarprime.
		\end{align}
		We now estimate $\int_{0}^{\ubar}\norm{\mathfrak{F}_{(\upsi)}}_{L^2(\suubarprime)}\diff\ubarprime$.
		First, we consider the top-order term $\twonorms{\WangD^{\vec{m}}\angnabla\ueta}{u}{1}{\sphere}{2}{\hypersurbar_{\ubar}}$ when $\vec{m}=(i,j,k)$ with $i+j+k=N-1$ and $k\geq 1$. Let $\vec{m}_0:=(i,j,k+1)$. By bootstrap assumptions \eqref{BA}:
		\begin{align}\label{PV:upsi4}
			\twonorms{\WangD^{\vec{m}}\angnabla\omega}{\ubar}{1}{\sphere}{2}{\hypersur_{u}}\lesssim&~\frac{\delta^{\frac{1}{2}}}{\abs{u}}\twonorms{\WangD^{\vec{m}_0}\omega}{\ubar}{2}{\sphere}{2}{\hypersur_{u}}\\
			\notag\lesssim&~\frac{\delta a^{\frac{1}{2}}}{\abs{u}}\tilde{\mathcal{O}}_{N,2},
		\end{align}
		\begin{align}\label{PV:upsi5}
			\twonorms{\WangD^{\vec{m}}\angnabla\ueta}{\ubar}{1}{\sphere}{2}{\hypersur_{u}}\lesssim&~\frac{\delta^{\frac{1}{2}}}{\abs{u}^2}\twonorms{\abs{u}\WangD^{\vec{m}_0}\ueta}{\ubar}{2}{\sphere}{2}{\hypersur_{u}}\\
			\notag\lesssim&~\frac{\delta^2 a}{\abs{u}^2},
		\end{align}
		By \eqref{ES:BoundsforProducts31}, bootstrap assumptions \eqref{BA} and proved estimates \eqref{ES:LinftyRicciCoefficients1}-\eqref{ES:LinftyRicciCoefficients3}, we have the estimates for the following terms:
		\begin{subequations}\label{PV:upsi6}
			\begin{align}
				\label{PV:upsi6.1}&\sum\limits_{\substack{\vec{m}_1+\vec{m}_2=(i,j-j_0,k-k_0)\\\abs{\vec{m}_2}\leq N-2}}\abs{u}^{j_0+k_0}\left(\twonorms{\WangD^{\vec{m}_1}\upsi^{j_0+k_0}\cdot\WangD^{\vec{m}_2}\angnabla\ueta)}{\ubar}{1}{\sphere}{2}{\hypersur_{u}}+\twonorms{\WangD^{\vec{m}_1}\upsi^{j_0+k_0+1}\cdot\WangD^{\vec{m}_2}\psi}{\ubar}{1}{\sphere}{2}{\hypersur_{u}}\right.\\
				\notag&\left.+\twonorms{\WangD^{\vec{m}_1}\upsi^{j_0+k_0+2}}{\ubar}{1}{\sphere}{2}{\hypersur_{u}}+\twonorms{\frac{1}{\abs{u}}\WangD^{\vec{m}_1}\upsi^{j_0+k_0}\cdot\WangD^{\vec{m}_2}(\omega,\gtr\upchi)}{\ubar}{1}{\sphere}{2}{\hypersur_{u}}\right)\\
				\notag\lesssim&~ \delta,\\
				\label{PV:upsi6.2}&\sum\limits_{\substack{\vec{m}_1+\vec{m}_2=(i,j-j_0,k-k_0)\\\abs{\vec{m}_2}\leq N-2}}\abs{u}^{j_0+k_0}\left(\twonorms{\WangD^{\vec{m}_1}\upsi^{j_0+k_0}\cdot\WangD^{\vec{m}_2}\angnabla\omega}{\ubar}{1}{\sphere}{2}{\hypersur_{u}}+\twonorms{\frac{1}{\abs{u}}\WangD^{\vec{m}_1}\upsi^{j_0+k_0}\cdot\WangD^{\vec{m}_2}\hchi}{\ubar}{1}{\sphere}{2}{\hypersur_{u}}\right)\\
				\notag\lesssim&~ \frac{\delta a^{\frac{1}{2}}}{\abs{u}}.
			\end{align}
		\end{subequations}
		Now we estimate the Ricci/scalar curvature terms:
		\begin{align}\label{PV:upsi7}
			&\sum\limits_{\vec{m}_1+\vec{m}_2=(i,j-j_0,k-k_0)}\abs{u}^{j_0+k_0}\twonorms{\WangD^{\vec{m}_1}\upsi^{j_0+k_0}\cdot\WangD^{\vec{m}_2}(\FRho,\FRhodual,\Scalar)}{\ubar}{1}{\sphere}{2}{\hypersur_{u}}\\
			\notag\lesssim&~\sum\limits_{\substack{\vec{m}_1+\vec{m}_2=(i,j-j_0,k-k_0)\\\abs{\vec{m}_1}>\frac{N-1}{2}}}\delta\abs{u}^{j_0+k_0}\left(\twonorms{\WangD^{\vec{m}_1}\upsi^{j_0+k_0}}{\ubar}{\infty}{\sphere}{2}{\hypersur_{u}}\cdot\twonorms{\WangD^{\vec{m}_2}(\FRho,\FRhodual,\Scalar)}{\ubar}{\infty}{\sphere}{\infty}{\hypersur_{u}}\right)\\
			\notag&\sum\limits_{\substack{\vec{m}_1+\vec{m}_2=(i,j-j_0,k-k_0)\\\abs{\vec{m}_2}>\frac{N-1}{2}}}\delta\abs{u}^{j_0+k_0}\left(\twonorms{\WangD^{\vec{m}_1}\upsi^{j_0+k_0}}{\ubar}{\infty}{\sphere}{\infty}{\hypersur_{u}}\cdot\twonorms{\WangD^{\vec{m}_2}(\FRho,\FRhodual,\Scalar)}{\ubar}{\infty}{\sphere}{2}{\hypersur_{u}}\right)\\
			\notag\lesssim&~ \delta\cdot\left(\BAsize\frac{\delta a^{\frac{1}{2}}}{\abs{u}}\right)^{j_0+k_0},
		\end{align}
		Now we estimate the largest fluid term using Proposition \ref{PR:EstimatesforRicciCurvature}:
		\begin{align}\label{PV:upsi8}
			&\sum\limits_{\vec{m}_1+\vec{m}_2=(i,j-j_0,k-k_0)}\abs{u}^{j_0+k_0}\twonorms{\WangD^{\vec{m}_1}\upsi^{j_0+k_0}\cdot\WangD^{\vec{m}_2}\FBetaA}{\ubar}{1}{\sphere}{2}{\hypersur_{u}}\\
			\notag\lesssim&~\sum\limits_{\substack{\vec{m}_1+\vec{m}_2=(i,j-j_0,k-k_0)\\\abs{\vec{m}_1}>\frac{N-1}{2}}}\delta\abs{u}^{j_0+k_0}\left(\twonorms{\WangD^{\vec{m}_1}\upsi^{j_0+k_0}}{\ubar}{1}{\sphere}{2}{\hypersur_{u}}\cdot\twonorms{\WangD^{\vec{m}_2}\FBetaA}{\ubar}{\infty}{\sphere}{\infty}{\hypersur_{u}}\right)\\
			\notag&+\sum\limits_{\substack{\vec{m}_1+\vec{m}_2=(i,j-j_0,k-k_0)\\\abs{\vec{m}_2}>\frac{N-1}{2}}}\delta\abs{u}^{j_0+k_0}\left(\twonorms{\WangD^{\vec{m}_1}\upsi^{j_0+k_0}}{\ubar}{\infty}{\sphere}{\infty}{\hypersur_{u}}\cdot\twonorms{\WangD^{\vec{m}_2}\FBetaA}{\ubar}{1}{\sphere}{2}{\hypersur_{u}}\right)\\
			\notag\lesssim&~ \delta^{\frac{1}{2}}\cdot\left(\BAsize\frac{\delta a^{\frac{1}{2}}}{\abs{u}}\right)^{j_0+k_0},
		\end{align}
		Finally, we estimate the curvature terms by bootstrap assumptions \eqref{BA} and Sobolev embedding $H^2\hookrightarrow L^\infty$ on $\suubar$ for $\WangD^{\vec{m}_2}\BetaA$:
		\begin{align}\label{PV:upsi9}
			&\sum\limits_{\vec{m}_1+\vec{m}_2=(i,j-j_0,k-k_0)}\abs{u}^{j_0+k_0}\twonorms{\WangD^{\vec{m}_1}\upsi^{j_0+k_0}\cdot\WangD^{\vec{m}_2}(\Rhodual,\Rho,\BetaA)}{\ubar}{1}{\sphere}{2}{\hypersur_{u}}\\
			\notag\lesssim&~\sum\limits_{\substack{\vec{m}_1+\vec{m}_2=(i,j-j_0,k-k_0)\\\abs{\vec{m}_1}>\frac{N-1}{2}}}\abs{u}^{j_0+k_0}\left(\twonorms{\WangD^{\vec{m}_1}\upsi^{j_0+k_0}}{\ubar}{1}{\sphere}{2}{\hypersur_{u}}\cdot\twonorms{\WangD^{\vec{m}_2}(\Rhodual,\Rho,\BetaA)}{\ubar}{\infty}{\sphere}{\infty}{\hypersur_{u}}\right)\\
			\notag&+\sum\limits_{\substack{\vec{m}_1+\vec{m}_2=(i,j-j_0,k-k_0)\\\abs{\vec{m}_2}>\frac{N-1}{2}}}\abs{u}^{j_0+k_0}\left(\twonorms{\WangD^{\vec{m}_1}\upsi^{j_0+k_0}}{\ubar}{\infty}{\sphere}{\infty}{\hypersur_{u}}\cdot\twonorms{\WangD^{\vec{m}_2}(\Rhodual,\Rho,\BetaA)}{\ubar}{1}{\sphere}{2}{\hypersur_{u}}\right)\\
			\notag\lesssim&~ \frac{\delta a^{\frac{1}{2}}}{\abs{u}}\cdot\left(\BAsize\frac{\delta a^{\frac{1}{2}}}{\abs{u}}\right)^{j_0+k_0}\cdot \mathcal{R}+\BAsize^{-1}\left(\BAsize\frac{\delta a^{\frac{1}{2}}}{\abs{u}}\right)^{j_0+k_0+2}\\
			\notag\lesssim&~\frac{\delta a^{\frac{1}{2}}}{\abs{u}}\cdot\left(1+\mathcal{R}\right).
		\end{align}
		Combining \eqref{PV:upsi2}-\eqref{PV:upsi9}, we obtain the desired result \eqref{ES:L2RicciCoefficients4} by the initial data assumptions \ref{Data}, and the choice of parameter $\delta,a,b$ in Section \ref{SS:Norms}, $\BAfluid,\BAsize$ in Section \eqref{SS:BA}.
	\end{proof}
	\begin{proof}[Proof of \eqref{ES:L2RicciCoefficients5}]
		We improve \eqref{ES:L2RicciCoefficients4} for $\WangD^{\vec{m}}\left(\gtr\uchi+\frac{2}{\abs{u}}\right)$. Compared to \eqref{PV:upsi1}, we write the schematic transport equation \eqref{Luchi} for $\gtr\uchi+\frac{2}{\abs{u}}$ as follows:
		\begin{align}\label{PV:traceuchi1}
			\angD_{ 4}\left(\gtr\uchi+\frac{2}{\abs{u}}\right)=\angnabla\ueta+\psi\upsi+\upsi\upsi+\frac{1}{u}(\omega,\gtr\upchi)+\Rho+\Scalar.
		\end{align}
		Let $\vec{m}=(i,j,k)$ with $\abs{\vec{m}}\leq N-2$ or $\abs{\vec{m}}= N-1$ and $k\geq 1$. Commuting \eqref{PV:tracechi1} with $\WangD^{\vec{m}}$, using \eqref{EQ:ang4ijk}, we obtain
		\begin{align}\label{PV:traceuchi2}
			& \angD_4\WangD^{\vec{m}}\left(\gtr\uchi+\frac{2}{\abs{u}}\right)\\
			\notag=&~ \WangD^{\vec{m}}\angnabla\ueta+ \sum\limits_{\substack{\vec{m}_1+\vec{m}_2=(i,j-j_0,k-k_0)\\\abs{\vec{m}_2}\leq N-2}}\abs{u}^{j_0+k_0}\WangD^{\vec{m}_1}\upsi^{j_0+k_0}\cdot\WangD^{\vec{m}_2}\angnabla\ueta\\
			\notag&+ \sum\limits_{\vec{m}_1+\vec{m}_2=(i,j-j_0,k-k_0)}\abs{u}^{j_0+k_0}\WangD^{\vec{m}_1}\upsi^{j_0+k_0+1}\cdot\WangD^{\vec{m}_2}\psi+ \sum\limits_{\vec{m}_1=(i,j-j_0,k-k_0)}\abs{u}^{j_0+k_0}\WangD^{\vec{m}_1}\upsi^{j_0+k_0+2}\\
			\notag&+ \sum\limits_{\vec{m}_1+\vec{m}_2=(i,j-j_0,k-k_0)}\abs{u}^{j_0+k_0-1}\WangD^{\vec{m}_1}\upsi^{j_0+k_0}\cdot\WangD^{\vec{m}_2}(\omega,\gtr\upchi)\\
			\notag&+ \sum\limits_{\vec{m}_1+\vec{m}_2=(i,j-j_0,k-k_0)}\abs{u}^{j_0+k_0}\WangD^{\vec{m}_1}\upsi^{j_0+k_0}\cdot\WangD^{\vec{m}_2}\Scalar\\
			\notag&+ \sum\limits_{\vec{m}_1+\vec{m}_2=(i,j-j_0,k-k_0-1)}\abs{u}^{j_0+k_0+1}\WangD^{\vec{m}_1}\upsi^{k_0+j_0+1}\cdot\WangD^{\vec{m}_2}\FBetaA\\
			\notag&+ \sum\limits_{\vec{m}_1+\vec{m}_2=(i,j-j_0,k-k_0)}\abs{u}^{j_0+k_0}\WangD^{\vec{m}_1}\upsi^{j_0+k_0}\cdot\WangD^{\vec{m}_2}\Rho\\
			\notag&+ \sum\limits_{\vec{m}_1+\vec{m}_2+\vec{m}_3=(i,j-j_0-1,k-k_0)}\abs{u}^{j_0+k_0+1}\WangD^{\vec{m}_1}\upsi^{j_0+k_0}\cdot\WangD^{\vec{m}_2}(\gtr\uchi+\frac{2}{\abs{u}})\cdot\WangD^{\vec{m}_3}\Rhodual,\\
			\notag&+ \sum\limits_{\vec{m}_1+\vec{m}_2+\vec{m}_3=(i,j-j_0,k-k_0-1)}\abs{u}^{j_0+k_0+1}\WangD^{\vec{m}_1}\upsi^{j_0+k_0}\cdot\WangD^{\vec{m}_2}(\gtr\uchi+\frac{2}{\abs{u}})\cdot\WangD^{\vec{m}_3}\BetaA\\
			\notag&=:\mathfrak{F}_{(\upsi)}.
		\end{align}
		By Proposition \ref{PR:transport1}, we have:
		\begin{align}\label{PV:traceuchi3}
			\norm{\WangD^{\vec{m}}\left(\gtr\uchi+\frac{2}{\abs{u}}\right)}_{L^2(\suubar)}\lesssim\norm{\WangD^{\vec{m}}\left(\gtr\uchi+\frac{2}{\abs{u}}\right)}_{L^2(S_{u,0})}+\int_{0}^{\ubar}	\norm{\mathfrak{F}_{\left(\gtr\uchi\right)}}_{L^2(\suubarprime)}\diff\ubarprime.
		\end{align}
		We now estimate $\int_{0}^{\ubar}\norm{\mathfrak{F}_{(\gtr\uchi)}}_{L^2(\suubarprime)}\diff\ubarprime$. We focus on the extra estimates compared to the proof of \eqref{ES:L2RicciCoefficients4}. For the fluid terms, we have:
		\begin{align}\label{PV:traceuchi4}
			&\sum\limits_{\vec{m}_1+\vec{m}_2=(i,j-j_0,k-k_0-1)}\abs{u}^{j_0+k_0+1}\twonorms{\WangD^{\vec{m}_1}\upsi^{j_0+k_0+1}\cdot\WangD^{\vec{m}_2}\FBetaA}{\ubar}{1}{\sphere}{2}{\hypersur_{u}}\\
			\notag\lesssim&~\sum\limits_{\substack{\vec{m}_1+\vec{m}_2=(i,j-j_0-1,k-k_0-1)\\\abs{\vec{m}_1}>\frac{N-1}{2}}}\delta\abs{u}^{j_0+k_0+1}\left(\twonorms{\WangD^{\vec{m}_1}\upsi^{j_0+k_0+1}}{\ubar}{\infty}{\sphere}{2}{\hypersur_{u}}\cdot\twonorms{\WangD^{\vec{m}_2}\FBetaA}{\ubar}{\infty}{\sphere}{\infty}{\hypersur_{u}}\right)\\
			\notag&+\sum\limits_{\substack{\vec{m}_1+\vec{m}_2=(i,j-j_0-1,k-k_0-1)\\\abs{\vec{m}_2}>\frac{N-1}{2}}}\delta\abs{u}^{j_0+k_0+1}\left(\twonorms{\WangD^{\vec{m}_1}\upsi^{j_0+k_0+1}}{\ubar}{\infty}{\sphere}{\infty}{\hypersur_{u}}\cdot\twonorms{\WangD^{\vec{m}_2}\FBetaA}{\ubar}{\infty}{\sphere}{2}{\hypersur_{u}}\right)\\
			\notag\lesssim&~  \BAsize\BAfluid\delta^{\frac{3}{2}} a^{\frac{1}{2}}\cdot\left(\BAsize\frac{\delta a^{\frac{1}{2}}}{\abs{u}}\right)^{j_0+k_0}\\
			\notag\lesssim&~ \delta^{\frac{3}{2}} a^{\frac{1}{2}}\cdot\left(\BAsize\frac{\delta a^{\frac{1}{2}}}{\abs{u}}\right)^{j_0+k_0}.
		\end{align}
		Now we estimate the second last term on the RHS of \eqref{PV:traceuchi2}. By bootstrap assumptions \eqref{BA} and Sobolev embedding $H^2\hookrightarrow L^\infty$ on $\suubar$ for $\WangD^{\vec{m}_2}\left(\gtr\uchi+\frac{2}{\abs{u}}\right)$:
		\begin{align}\label{PV:traceuchi5}
			&\sum\limits_{\vec{m}_1+\vec{m}_2+\vec{m}_3=(i,j-j_0-1,k-k_0)}\abs{u}^{k_0+j_0+1}\twonorms{\WangD^{\vec{m}_1}\upsi^{j_0+k_0}\cdot\WangD^{\vec{m}_2}\left(\gtr\uchi+\frac{2}{\abs{u}}\right)\cdot\WangD^{\vec{m}_3}\Rhodual}{\ubar}{1}{\sphere}{2}{\hypersur_{u}}\\
			\notag\lesssim&~\sum\limits_{\substack{\vec{m}_1+\vec{m}_2+\vec{m}_3=(i,j-j_0-1,k-k_0)\\\abs{\vec{m}_1}>\frac{N-1}{2}}}\abs{u}^{k_0+j_0+1}\delta\left(\twonorms{\WangD^{\vec{m}_1}\upsi^{j_0+k_0}}{\ubar}{\infty}{\sphere}{2}{\hypersur_{u}}\right.\\
			\notag&\left.\cdot\twonorms{\WangD^{\vec{m}_2}\left(\gtr\uchi+\frac{2}{\abs{u}}\right)}{\ubar}{\infty}{\sphere}{\infty}{\hypersur_{u}}\cdot\twonorms{\WangD^{\vec{m}_3}\Rhodual}{\ubar}{\infty}{\sphere}{\infty}{\hypersur_{u}}\right)\\
			\notag&+\sum\limits_{\substack{\vec{m}_1+\vec{m}_2+\vec{m}_3=(i,j-j_0-1,k-k_0)\\\abs{\vec{m}_2}>\frac{N-1}{2}}}\abs{u}^{k_0+j_0+1}\delta\left(\twonorms{\WangD^{\vec{m}_1}\upsi^{j_0+k_0}}{\ubar}{\infty}{\sphere}{\infty}{\hypersur_{u}}\right.\\
			\notag&\left.\cdot\twonorms{\WangD^{\vec{m}_2}\left(\gtr\uchi+\frac{2}{\abs{u}}\right)}{\ubar}{\infty}{\sphere}{2}{\hypersur_{u}}\cdot\twonorms{\WangD^{\vec{m}_3}\Rhodual}{\ubar}{\infty}{\sphere}{\infty}{\hypersur_{u}}\right)\\
			\notag&+\sum\limits_{\substack{\vec{m}_1+\vec{m}_2+\vec{m}_3=(i,j-j_0-1,k-k_0)\\\abs{\vec{m}_3}>\frac{N-1}{2}}}\abs{u}^{k_0+j_0+1}\delta\left(\twonorms{\WangD^{\vec{m}_1}\upsi^{j_0+k_0}}{\ubar}{\infty}{\sphere}{\infty}{\hypersur_{u}}\right.\\
			\notag&\left.\cdot\twonorms{\WangD^{\vec{m}_2}\left(\gtr\uchi+\frac{2}{\abs{u}}\right)}{\ubar}{\infty}{\sphere}{\infty}{\hypersur_{u}}\cdot\twonorms{\WangD^{\vec{m}_3}\Rhodual}{\ubar}{\infty}{\sphere}{2}{\hypersur_{u}}\right)\\
			\notag\lesssim&~\BAsize^{-1}\left(\BAsize\frac{\delta a^{\frac{1}{2}}}{\abs{u}}\right)^{j_0+k_0+2}\cdot\sum\limits_{\abs{\vec{m}_2}\leq N-1}\twonorms{\WangD^{\vec{m}_2}\left(\gtr\uchi+\frac{2}{\abs{u}}\right)}{\ubar}{\infty}{\sphere}{2}{\hypersur_{u}}.
		\end{align}
		Similarly, we estimate the last term on the RHS of \eqref{PV:traceuchi2}:
		\begin{align}\label{PV:traceuchi6}
			&\sum\limits_{\vec{m}_1+\vec{m}_2+\vec{m}_3=(i,j-j_0,k-k_0-1)}\abs{u}^{k_0+j_0}\twonorms{\WangD^{\vec{m}_1}\upsi^{j_0+k_0}\cdot\WangD^{\vec{m}_2}\left(\gtr\uchi+\frac{2}{\abs{u}}\right)\cdot\WangD^{\vec{m}_3}\BetaA}{\ubar}{1}{\sphere}{2}{\hypersur_{u}}\\
			\notag\lesssim&~\sum\limits_{\substack{\vec{m}_1+\vec{m}_2+\vec{m}_3=(i,j-j_0,k-k_0-1)\\\abs{\vec{m}_1}>\frac{N-1}{2}}}\abs{u}^{k_0+j_0}\delta\left(\twonorms{\WangD^{\vec{m}_1}\upsi^{j_0+k_0}}{\ubar}{\infty}{\sphere}{2}{\hypersur_{u}}\right.\\
			\notag&\left.\cdot\twonorms{\WangD^{\vec{m}_2}\left(\gtr\uchi+\frac{2}{\abs{u}}\right)}{\ubar}{\infty}{\sphere}{\infty}{\hypersur_{u}}\cdot\twonorms{\WangD^{\vec{m}_3}\BetaA}{\ubar}{\infty}{\sphere}{\infty}{\hypersur_{u}}\right)\\
			\notag&+\sum\limits_{\substack{\vec{m}_1+\vec{m}_2+\vec{m}_3=(i,j-j_0,k-k_0-1)\\\abs{\vec{m}_2}>\frac{N-1}{2}}}\abs{u}^{k_0+j_0}\delta\left(\twonorms{\WangD^{\vec{m}_1}\upsi^{j_0+k_0}}{\ubar}{\infty}{\sphere}{\infty}{\hypersur_{u}}\right.\\
			\notag&\left.\cdot\twonorms{\WangD^{\vec{m}_2}\left(\gtr\uchi+\frac{2}{\abs{u}}\right)}{\ubar}{\infty}{\sphere}{2}{\hypersur_{u}}\cdot\twonorms{\WangD^{\vec{m}_3}\BetaA}{\ubar}{\infty}{\sphere}{\infty}{\hypersur_{u}}\right)\\
			\notag&+\sum\limits_{\substack{\vec{m}_1+\vec{m}_2+\vec{m}_3=(i,j-j_0,k-k_0-1)\\\abs{\vec{m}_3}>\frac{N-1}{2}}}\abs{u}^{k_0+j_0}\delta\left(\twonorms{\WangD^{\vec{m}_1}\upsi^{j_0+k_0}}{\ubar}{\infty}{\sphere}{\infty}{\hypersur_{u}}\right.\\
			\notag&\left.\cdot\twonorms{\WangD^{\vec{m}_2}\left(\gtr\uchi+\frac{2}{\abs{u}}\right)}{\ubar}{\infty}{\sphere}{\infty}{\hypersur_{u}}\cdot\twonorms{\WangD^{\vec{m}_3}\BetaA}{\ubar}{\infty}{\sphere}{2}{\hypersur_{u}}\right)\\
			\notag\lesssim&~\left(\BAsize\frac{\delta a^{\frac{1}{2}}}{\abs{u}}\right)^{i_0+k_0+1}\cdot\sum\limits_{\abs{\vec{m}_2}\leq N-1}\twonorms{\WangD^{\vec{m}_2}\left(\gtr\uchi+\frac{2}{\abs{u}}\right)}{\ubar}{\infty}{\sphere}{2}{\hypersur_{u}}.
		\end{align}
		By \eqref{PV:upsi5}, \eqref{PV:upsi6.1} \eqref{PV:upsi7}, \eqref{PV:upsi9} for $\Rho$, the rest of the terms on the RHS of \eqref{PV:traceuchi2} in $L^1_uL^2_\sphere(\hypersur_{u})$-norms are bounded by $\delta$. Combining \eqref{PV:upsi5}, \eqref{PV:upsi6.1} \eqref{PV:upsi7}, \eqref{PV:upsi9} and  \eqref{PV:traceuchi2}-\eqref{PV:traceuchi5}, soaking $\sum\limits_{\abs{\vec{m}_2}\leq N-1}\twonorms{\WangD^{\vec{m}_2}\left(\gtr\uchi+\frac{2}{\abs{u}}\right)}{\ubar}{\infty}{\sphere}{2}{\hypersur_{u}}$ in \eqref{PV:traceuchi5}-\eqref{PV:traceuchi6} to the LHS of \eqref{PV:traceuchi3}, we obtain the desired result \eqref{ES:L2RicciCoefficients5}.
	\end{proof}
	\begin{proof}[Proof of \eqref{ES:L2RicciCoefficients6}]
		Let $\vec{m}=(i,j,k)$ with $\abs{\vec{m}}\leq N-2$ or $\abs{\vec{m}}= N-1$ and $k\geq 1$. Commuting $\angD_{ 4}\ln\lapse=\omega$ with $\WangD^{\vec{m}}$, using \eqref{EQ:ang4ijk}, we obtain:
		\begin{align}\label{PV:lapse1}
			& \angD_4\frac{1}{\abs{u}}\WangD^{\vec{m}}\ln\lapse\\
			\notag=&~ \frac{1}{\abs{u}}\WangD^{\vec{m}}\omega+ \sum\limits_{\vec{m}_1+\vec{m}_2+\vec{m}_3=(i,j-j_0,k-k_0)}\abs{u}^{j_0+k_0-1}\WangD^{\vec{m}_1}\upsi^{j_0+k_0}\cdot\WangD^{\vec{m}_2}\ln\lapse\cdot\WangD^{\vec{m}_3}\psi\\
			\notag&+ \sum\limits_{\vec{m}_1=(i,j-j_0,k-k_0)}\abs{u}^{j_0+k_0-1}\WangD^{\vec{m}_1}\upsi^{j_0+k_0+2}\\
			\notag&+ \sum\limits_{\vec{m}_1+\vec{m}_2=(i,j-j_0,k-k_0-1)}\abs{u}^{j_0+k_0}\WangD^{\vec{m}_1}\upsi^{k_0+j_0+1}\cdot\WangD^{\vec{m}_2}\FBetaA\\	
			\notag&+ \sum\limits_{\vec{m}_1+\vec{m}_2=(i,j-j_0,k-k_0)}\abs{u}^{j_0+k_0-1}\WangD^{\vec{m}_1}\upsi^{j_0+k_0}\cdot\WangD^{\vec{m}_2}\Rho\\
			\notag&+ \sum\limits_{\vec{m}_1+\vec{m}_2+\vec{m}_3=(i,j-j_0-1,k-k_0)}\abs{u}^{j_0+k_0}\WangD^{\vec{m}_1}\upsi^{j_0+k_0}\cdot\WangD^{\vec{m}_2}\ln\lapse\cdot\WangD^{\vec{m}_3}\Rhodual,\\
			\notag&+ \sum\limits_{\vec{m}_1+\vec{m}_2+\vec{m}_3=(i,j-j_0,k-k_0-1)}\abs{u}^{j_0+k_0-1}\WangD^{\vec{m}_1}\upsi^{j_0+k_0}\cdot\WangD^{\vec{m}_2}\ln\lapse\cdot\WangD^{\vec{m}_3}\BetaA\\
			\notag&=:\mathfrak{F}_{(\ln\lapse)}.
		\end{align}
		By Proposition \ref{PR:transport1}, we have:
		\begin{align}\label{PV:lapse2}
			\norm{\frac{1}{\abs{u}}\WangD^{\vec{m}}\ln\lapse}_{L^2(\suubar)}\lesssim\norm{\frac{1}{\abs{u}}\WangD^{\vec{m}}\ln\lapse}_{L^2(S_{u,0})}+\int_{0}^{\ubar}	\norm{\mathfrak{F}_{(\ln\lapse)}}_{L^2(\suubarprime)}\diff\ubarprime.
		\end{align}
		We now estimate $\int_{0}^{\ubar}\norm{\mathfrak{F}_{(\ln\lapse)}}_{L^2(\suubarprime)}\diff\ubarprime$.
		We consider the second term on the RHS of \eqref{PV:lapse2}. By bootstrap assumptions \eqref{BA} and Sobolev embedding $H^2\hookrightarrow L^\infty$ on $\suubar$ for $\WangD^{\vec{m}_2}\ln\lapse$:
		\begin{align}\label{PV:lapse3}
			&\sum\limits_{\vec{m}_1+\vec{m}_2+\vec{m}_3=(i,j-j_0,k-k_0)}\abs{u}^{k_0+j_0}\twonorms{\WangD^{\vec{m}_1}\upsi^{j_0+k_0}\cdot\frac{1}{\abs{u}}\WangD^{\vec{m}_2}\ln\lapse\cdot\WangD^{\vec{m}_3}\psi}{\ubar}{1}{\sphere}{2}{\hypersur_{u}}\\
			\notag\lesssim&~\left(\BAsize\frac{\delta a^{\frac{1}{2}}}{\abs{u}}\right)^{j_0+k_0+1}\cdot\sum\limits_{\abs{\vec{m}_2}\leq N-1}\twonorms{\frac{1}{\abs{u}}\WangD^{\vec{m}_2}\ln\lapse}{\ubar}{\infty}{\sphere}{2}{\hypersur_{u}}.
		\end{align}
		Applying the same approach as in \eqref{PV:traceuchi5}-\eqref{PV:traceuchi6}, we derive
		\begin{subequations}\label{PV:lapse4}
			\begin{align}
				&\sum\limits_{\vec{m}_1+\vec{m}_2+\vec{m}_3=(i,j-j_0,k-k_0-1)}\abs{u}^{k_0+j_0}\twonorms{\WangD^{\vec{m}_1}\upsi^{j_0+k_0}\cdot\frac{1}{\abs{u}}\WangD^{\vec{m}_2}\ln\lapse\cdot\WangD^{\vec{m}_3}\BetaA}{\ubar}{1}{\sphere}{2}{\hypersur_{u}}\\
				\notag\lesssim&~\left(\BAsize\frac{\delta a^{\frac{1}{2}}}{\abs{u}}\right)^{j_0+k_0+1}\cdot\sum\limits_{\abs{\vec{m}_2}\leq N-1}\twonorms{\frac{1}{\abs{u}}\WangD^{\vec{m}_2}\ln\lapse}{\ubar}{\infty}{\sphere}{2}{\hypersur_{u}},\\
				&\sum\limits_{\vec{m}_1+\vec{m}_2+\vec{m}_3=(i,j-j_0-1,k-k_0)}\abs{u}^{k_0+j_0}\twonorms{\WangD^{\vec{m}_1}\upsi^{j_0+k_0}\cdot\frac{1}{\abs{u}}\WangD^{\vec{m}_2}\ln\lapse\cdot\WangD^{\vec{m}_3}\Rhodual}{\ubar}{1}{\sphere}{2}{\hypersur_{u}}\\
				\notag\lesssim&~\BAsize^{-1}\left(\BAsize\frac{\delta a^{\frac{1}{2}}}{\abs{u}}\right)^{j_0+k_0+2}\cdot\sum\limits_{\abs{\vec{m}_2}\leq N-1}\twonorms{\frac{1}{\abs{u}}\WangD^{\vec{m}_2}\ln\lapse}{\ubar}{\infty}{\sphere}{2}{\hypersur_{u}}.
			\end{align}
		\end{subequations}
		Combining \eqref{PV:upsi6.1}, \eqref{PV:upsi8}, \eqref{PV:upsi9} for $\Rho$ and  \eqref{PV:lapse1}-\eqref{PV:lapse4}, soaking $\sum\limits_{\abs{\vec{m}_2}\leq N-1}\twonorms{\WangD^{\vec{m}_2}\left(\gtr\uchi+\frac{2}{\abs{u}}\right)}{\ubar}{\infty}{\sphere}{2}{\hypersur_{u}}$ in \eqref{PV:lapse3}-\eqref{PV:lapse4} to the LHS of \eqref{PV:traceuchi3}, we obtain the desired result \eqref{ES:L2RicciCoefficients6}.
	\end{proof}
	\section{Top-order $L^2$-Estimates for the Ricci Coefficients}
	\begin{proposition}[Top-order $L^2(\suubar)$-estimates for Ricci coefficients]\label{PR:topL2RicciCoefficients}
		$\tilde{\mathcal{O}}_{N,2},\mathcal{R}$ are defined in Section \ref{SS:Norms}. Let $\vec{m}=(i,j,k)$ with $i+j+k=N$ and $k\geq 2$. Let $\WangD^{\vec{m}}$ be as in Definition \ref{DE:operators}. Under the bootstrap assumptions \eqref{BA}, we have the estimates for the derivatives of the Ricci coefficients:
		\begin{subequations}\label{ES:topL2RicciCoefficients}
			\begin{align}
				\label{ES:topL2RicciCoefficients1}	\sum\limits_{\substack{\abs{\vec{m}}=N\\ k\geq 2}}\twonorms{\WangD^{\vec{m}}\gtr\upchi}{\ubar}{2}{\sphere}{2}{\hypersur_{u}}\lesssim&~ \delta^{\frac{1}{2}}a^{\frac{1}{2}}\cdot\left(\BAsize\delta a^{\frac{1}{2}}\right)+ \delta^{\frac{1}{2}},\\
				\label{ES:topL2RicciCoefficients2}	\sum\limits_{\substack{\abs{\vec{m}}=N\\ k\geq 2}}\twonorms{\WangD^{\vec{m}}\hchi}{\ubar}{2}{\sphere}{2}{\hypersur_{u}}\lesssim&~ \delta^{\frac{1}{2}}a^{\frac{1}{2}}\left(1+\mathcal{R}\right),\\
				\label{ES:topL2RicciCoefficients3}	\sum\limits_{\substack{\abs{\vec{m}}=N\\ k\geq 2}}\twonorms{\WangD^{\vec{m}}\omega}{\ubar}{2}{\sphere}{2}{\hypersur_{u}}\lesssim&~ \delta^{\frac{1}{2}}a^{\frac{1}{2}}\left(1+\mathcal{R}\right)+1,\\
				\label{ES:topL2RicciCoefficients4}	\sum\limits_{\substack{\abs{\vec{m}}=N\\ k\geq 2}}\norm{\abs{u}\WangD^{\vec{m}}\upeta}_{L^2_{\ubar}L_\sphere^2(\hypersur_{u})}\lesssim&~\delta^{\frac{3}{2}} a(1+\mathcal{R}),\\
				\label{ES:topL2RicciCoefficients5}	\sum\limits_{\substack{\abs{\vec{m}}=N\\ k\geq 2}}\norm{\WangD^{\vec{m}}\upeta}_{L^2_{u}L_\sphere^2(\hypersurbar_{\ubar})}\lesssim&~\delta^{\frac{1}{2}} a^{\frac{1}{2}}(1+\mathcal{R}),\\
				\label{ES:topL2RicciCoefficients6}	\sum\limits_{\substack{\abs{\vec{m}}=N\\ k\geq 2}}\norm{\abs{u}\WangD^{\vec{m}}\ueta}_{L^2_{\ubar}L_\sphere^2(\hypersur_{u})}\lesssim&~\delta^{\frac{3}{2}} a(1+\mathcal{R}).
			\end{align}
		\end{subequations}
	\end{proposition}
	We prove Proposition \ref{PR:topL2RicciCoefficients} in Section \ref{SS:provetopRicci}. 
	\subsection{Proof of Proposition \ref{PR:topL2RicciCoefficients}}\label{SS:provetopRicci}
	In this section, we only derive the estimates for the top-order terms, while the below-top-order terms can be estimated using the same method as in Section \ref{SS:proveRicci}.
	\begin{proof}[Proof of \eqref{ES:topL2RicciCoefficients1}]		
		Let $\vec{m}=(i,j,k)$ with $i+j+k=N$ and $k\geq 2$. Commuting \eqref{EQ:e4upchi} with $\abs{u}\WangD^{\vec{m}}$, using \eqref{EQ:ang4ijk}, we obtain \eqref{PV:tracechi2} multiplied by $\abs{u}$. In addition to the proof of \eqref{ES:L2RicciCoefficients2}, we need to estimate the following top-order terms using proved results in Proposition \ref{PR:EstimatesforRicciCurvature}, Corollary \ref{CO:LinftyRicciCoefficients} and bootstrap assumptions \eqref{BA}:
		\begin{subequations}\label{PV:toptracechi}
			\begin{align}
				&\twonorms{\abs{u}\cdot\psi\cdot\WangD^{\vec{m}}\psi}{\ubar}{1}{\sphere}{2}{\hypersur_{u}}\\
				\notag\lesssim&~\delta^{\frac{1}{2}}\cdot\twonorms{u\psi}{\ubar}{\infty}{\sphere}{\infty}{\hypersur_{u}}\cdot\twonorms{\WangD^{\vec{m}}\psi}{\ubar}{2}{\sphere}{2}{\hypersur_{u}}\\
				\notag\lesssim&~\BAsize\delta a,\\
				&\sum\limits_{\abs{\vec{m}_1}=N}\left(\twonorms{\abs{u}\WangD^{\vec{m}_1}\FAlphaAB}{\ubar}{1}{\sphere}{2}{\hypersur_{u}}\right)\\
				\notag\lesssim&~1,\\
				&\sum\limits_{\abs{\vec{m}_1}=N-1,k\geq 1}\twonorms{\abs{u}\cdot\gtr\upchi\cdot\WangD^{\vec{m}_1}(\BetaA,\CRhodual)}{\ubar}{1}{\sphere}{2}{\hypersur_{u}}\\
				\notag\lesssim&~\delta^{\frac{1}{2}}\cdot\twonorms{\gtr\upchi}{\ubar}{\infty}{\sphere}{\infty}{\hypersur_{u}}\cdot\twonorms{\abs{u}\WangD^{\vec{m}_1}(\BetaA,\CRhodual)}{\ubar}{2}{\sphere}{2}{\hypersur_{u}}\\
				\notag\lesssim&~\BAsize\delta a,\\
			\end{align}
		\end{subequations}
		Combining \eqref{PV:tracechi3} (with an extra $\abs
		{u}$-weight), \eqref{PV:toptracechi} and  \eqref{PV:tracechi4}-\eqref{PV:tracechi9}, we have:
		\begin{align}
			\norm{\abs{u}\WangD^{\vec{m}}\gtr\upchi}_{L^2(\suubar)}\lesssim \BAsize\delta a+1.
		\end{align}
		Taking $L^2_{\ubar}$-norm, we obtain the desired result \eqref{ES:topL2RicciCoefficients1}.
	\end{proof}
	\begin{proof}[Proof of \eqref{ES:topL2RicciCoefficients2}]
		We recall the Codazzi equation \eqref{angdivhchi} as below:
		\begin{align}\label{PV:tophchi1}
			\angdiv\hchi=\angnabla\gtr\upchi+\BetaA+\psi\upsi.
		\end{align}
		Let $\vec{m}=(i,j,k)$ with $\abs{\vec{m}}\leq N-2$ or $\abs{\vec{m}}= N-1$ and $k\geq 1$. Commuting \eqref{PV:tophchi1} with $\abs{u}\WangD^{\vec{m}}$, using \eqref{EQ:angAijk}, we obtain
		\begin{align}\label{PV:tophchi2}
			&\angdiv\WangD^{\vec{m}}\hchi\\
			\notag=&~\sum\limits_{\vec{m}_1+\vec{m}_2+\vec{m}_3=(i-i_0,j-j_0,k)}\delta^{i_0}\abs{u}^{j_0}\WangD^{\vec{m}_1}\psi^{i_0}\cdot\WangD^{\vec{m}_2}\upsi^{j_0}\cdot\WangD^{\vec{m}_3}\angnabla\gtr\upchi\\
			\notag&\sum\limits_{\vec{m}_1+\vec{m}_2=(i-i_0,j-j_0,k)}\delta^{i_0}\abs{u}^{j_0}\WangD^{\vec{m}_1}\psi^{i_0+1}\cdot\WangD^{\vec{m}_2}\upsi^{j_0+1}\\
			\notag&\sum\limits_{\vec{m}_1+\vec{m}_2+\vec{m}_3=(i-i_0-1,j-j_0,k)}\delta^{i_0+1}\abs{u}^{j_0}\WangD^{\vec{m}_1}\psi^{i_0+1}\cdot\WangD^{\vec{m}_2}\upsi^{j_0}\cdot\WangD^{\vec{m}_3}\FBetaA\\	\notag&\sum\limits_{\vec{m}_1+\vec{m}_2+\vec{m}_3=(i-i_0,j-j_0-1,k)}\delta^{i_0}\abs{u}^{j_0+1}\WangD^{\vec{m}_1}\psi^{i_0+1}\cdot\WangD^{\vec{m}_2}\upsi^{j_0}\cdot\WangD^{\vec{m}_3}\FBetabar\\
			\notag&\sum\limits_{\vec{m}_1+\vec{m}_2+\vec{m}_3=(i-i_0,j-j_0,k)}\delta^{i_0}\abs{u}^{j_0}\WangD^{\vec{m}_1}\psi^{i_0}\cdot\WangD^{\vec{m}_2}\upsi^{j_0}\cdot\WangD^{\vec{m}_3}\BetaA\\
			\notag&\sum\limits_{\vec{m}_1+\vec{m}_2+\vec{m}_3=(i-i_0,j-j_0-1,k)}\delta^{i_0}\abs{u}^{j_0+1}\WangD^{\vec{m}_1}\psi^{i_0+1}\cdot\WangD^{\vec{m}_2}\upsi^{j_0}\cdot\WangD^{\vec{m}_3}\Betabar\\
			\notag&=:\mathfrak{F}_{(\angdiv\hchi)}.
		\end{align}
		Applying elliptic estimates in Lemma \ref{LE:ellipticES} and taking $L^2_{\ubar}$-norm, we have:
		\begin{align}\label{PV:tophchi3}
			&\twonorms{\angnabla\WangD^{\vec{m}}\hchi}{\ubar}{2}{\sphere}{2}{\hypersur_{u}}\lesssim\twonorms{\mathfrak{F}_{(\angdiv\hchi)}}{\ubar}{2}{\sphere}{2}{\hypersur_{u}}.
		\end{align}
		Now we estimates terms in $\twonorms{\mathfrak{F}_{(\angdiv\hchi)}}{\ubar}{2}{\sphere}{2}{\hypersur_{u}}$.
		For the first term on the RHS of \eqref{PV:tophchi2}, by bootstrap assumptions \eqref{BA}, proved results in Proposition \ref{PR:L2RicciCoefficients}, Corollary \eqref{CO:LinftyRicciCoefficients}, Lemma \ref{LE:BoundsforProducts} and \eqref{ES:topL2RicciCoefficients1}, we obtain:	
		\begin{align}\label{PV:tophchi4}
			&\sum\limits_{\vec{m}_1+\vec{m}_2+\vec{m}_3=(i-i_0,j-j_0,k)}\delta^{i_0}\abs{u}^{j_0}\twonorms{\WangD^{\vec{m}_1}\psi^{i_0}\cdot\WangD^{\vec{m}_2}\upsi^{j_0}\cdot\WangD^{\vec{m}_3}\angnabla\gtr\upchi}{\ubar}{2}{\sphere}{2}{\hypersur_{u}}\\
			\notag\lesssim&~\left(\delta a^\frac{1}{2}\right)^{i_0}\cdot\left(\BAsize\frac{\delta a^{\frac{1}{2}}}{\abs{u}}\right)^{j_0}\cdot\left(\delta^{\frac{1}{2}}a^{\frac{1}{2}}\cdot\left(\BAsize\delta a^{\frac{1}{2}}\right)+\delta^{\frac{3}{2}} b\right)\\
			\notag\lesssim&~\delta^{\frac{1}{2}} a^\frac{1}{2}+\delta^{\frac{3}{2}} b.
		\end{align}
		By proved results in Proposition \ref{PR:L2RicciCoefficients}, Corollary \eqref{CO:LinftyRicciCoefficients}, Lemma \ref{LE:BoundsforProducts}, for the second term on the RHS of \eqref{PV:tophchi2}, we have:
		\begin{align}\label{PV:tophchi5}
			&\sum\limits_{\vec{m}_1+\vec{m}_2=(i-i_0,j-j_0,k)}\delta^{i_0}\abs{u}^{j_0}\twonorms{\WangD^{\vec{m}_1}\psi^{i_0+1}\cdot\WangD^{\vec{m}_2}\upsi^{j_0+1}}{\ubar}{2}{\sphere}{2}{\hypersur_{u}}\\
			\notag\lesssim&~\delta^{\frac{1}{2}}a^{\frac{1}{2}}\cdot\left(\delta a^\frac{1}{2}\right)^{i_0}\cdot\left(\BAsize\frac{\delta a^{\frac{1}{2}}}{\abs{u}}\right)^{j_0+1}\cdot\abs{u}^{-1}\\
			\notag\lesssim&~\delta^{\frac{1}{2}} a^\frac{1}{2}.
		\end{align}
		Now we estimate the fluid terms. By bootstrap assumptions \eqref{BA}, Proposition \ref{PR:EstimatesforRicciCurvature}, Proposition \ref{PR:L2RicciCoefficients}, Corollary \eqref{CO:LinftyRicciCoefficients}, we have the following estimates for the third and forth term on the RHS of \eqref{PV:tophchi2}:
		\begin{align}\label{PV:tophchi6}
			&\sum\limits_{\vec{m}_1+\vec{m}_2+\vec{m}_3=(i-i_0-1,j-j_0,k)}\delta^{i_0+1}\abs{u}^{j_0}\twonorms{\WangD^{\vec{m}_1}\psi^{i_0+1}\cdot\WangD^{\vec{m}_2}\upsi^{j_0}\cdot\WangD^{\vec{m}_3}\FBetaA}{\ubar}{2}{\sphere}{2}{\hypersur_{u}}\\
			\notag\lesssim&~\delta^{\frac{1}{2}}a^{\frac{1}{2}}\cdot\left(\delta a^\frac{1}{2}\right)^{i_0}\cdot\left(\BAsize\frac{\delta a^{\frac{1}{2}}}{\abs{u}}\right)^{j_0}\\
			\notag\lesssim&~\delta^{\frac{1}{2}} a^\frac{1}{2}.
		\end{align}
		\begin{align}\label{PV:tophchi7}
			&\sum\limits_{\vec{m}_1+\vec{m}_2+\vec{m}_3=(i-i_0,j-j_0-1,k)}\delta^{i_0}\abs{u}^{j_0+1}\twonorms{\WangD^{\vec{m}_1}\psi^{i_0+1}\cdot\WangD^{\vec{m}_2}\upsi^{j_0}\cdot\WangD^{\vec{m}_3}\FBetabar}{\ubar}{2}{\sphere}{2}{\hypersur_{u}}\\
			\notag\lesssim&~\delta^{\frac{1}{2}}a^{\frac{1}{2}}\cdot\left(\delta a^\frac{1}{2}\right)^{i_0}\cdot\left(\BAsize\frac{\delta a^{\frac{1}{2}}}{\abs{u}}\right)^{j_0}\cdot\abs{u}\\
			\notag\lesssim&~\delta^{\frac{1}{2}} a^\frac{1}{2}.
		\end{align}
		Lastly, we estimate the curvature terms using bootstrap assumptions \eqref{BA}, Proposition \ref{PR:L2RicciCoefficients} and Corollary \eqref{CO:LinftyRicciCoefficients} as below:
		\begin{align}\label{PV:tophchi8}
			&\sum\limits_{\vec{m}_1+\vec{m}_2+\vec{m}_3=(i-i_0,j-j_0,k)}\delta^{i_0}\abs{u}^{j_0}\twonorms{\WangD^{\vec{m}_1}\psi^{i_0}\cdot\WangD^{\vec{m}_2}\upsi^{j_0}\cdot\WangD^{\vec{m}_3}\BetaA}{\ubar}{2}{\sphere}{2}{\hypersur_{u}}\\
			\notag\lesssim&~\frac{\delta^{\frac{1}{2}}a^\frac{1}{2}}{\abs{u}^2}\mathcal{R}\cdot\left(\delta a^\frac{1}{2}\right)^{i_0}\cdot\left(\BAsize\frac{\delta a^{\frac{1}{2}}}{\abs{u}}\right)^{j_0}\\
			\notag\lesssim&~\delta^{\frac{1}{2}} a^\frac{1}{2}\mathcal{R}.
		\end{align}
		\begin{align}\label{PV:tophchi9}
			&\sum\limits_{\vec{m}_1+\vec{m}_2+\vec{m}_3=(i-i_0,j-j_0-1,k)}\delta^{i_0}\abs{u}^{j_0+1}\twonorms{\WangD^{\vec{m}_1}\psi^{i_0+1}\cdot\WangD^{\vec{m}_2}\upsi^{j_0}\cdot\WangD^{\vec{m}_3}\Betabar}{\ubar}{2}{\sphere}{2}{\hypersur_{u}}\\
			\notag\lesssim&~\delta^{\frac{1}{2}}a^{\frac{1}{2}}\left(\delta a^\frac{1}{2}\right)^{i_0}\cdot\left(\BAsize\frac{\delta a^{\frac{1}{2}}}{\abs{u}}\right)^{j_0+1}\\
			\notag\lesssim&~\delta^{\frac{1}{2}} a^\frac{1}{2}.
		\end{align}
		Combining \eqref{PV:tophchi2}-\eqref{PV:tophchi9}, we obtain the desired result \eqref{ES:topL2RicciCoefficients2}.
	\end{proof}
	\begin{proof}[Proof of \eqref{ES:topL2RicciCoefficients3}]
		We define $\omega^{\dagger}$ to be the solution to 
		\begin{align}\label{PV:topomega1}
			\angD_3\omega^{\dagger}=&~\frac{1}{2}\CRhodual,\\
			\notag\omega^{\dagger}=&~0,&
			\text{ on }\hypersur_{u_\infty}. 
		\end{align}
		Then we define 
		\begin{align}\label{PV:topomega2}
			\kappa:=\angnabla\omega+\HodgeLeftDual\angnabla\omega^{\dagger}-\frac{1}{2}\CBetaA.
		\end{align}
		Let $\vec{m}=(i,j,k)$ with $\abs{\vec{m}}\leq N-2$ or $\abs{\vec{m}}= N-1$ and $k\geq 1$. Commuting \eqref{PV:topomega1} with $\WangD^{\vec{m}}$ using \eqref{EQ:ang3ijk}, applying transport estimates in Proposition \ref{PR:transport1} and bootstrap assumptions \eqref{BA}, we have:
		\begin{align}\label{PV:topomega3}
			\sum\limits_{\abs{\vec{m}}\leq N-1}\norm{\WangD^{\vec{m}}\omega^{\dagger}}_{L^2(\suubar)}\lesssim1+\BAsize\frac{\delta a^{\frac{1}{2}}}{\abs{u}}\lesssim 1.
		\end{align}
		For convenience, we use the notation $\psi\in\{\gtr\upchi,\hchi,\omega,\omega^{\dagger}\}$ though $\omega^{\dagger}$ is smaller in size.
		We consider the transport equation for $\kappa$:
		\begin{align}\label{PV:topomega3.5}
			\angD_3\kappa+\frac{1}{2}\gtr\uchi\kappa=&~\sum_{k_1+k_2+k_3=1}\upsi^{k_1}\angnabla^{i_2}\left(\gtr\uchi+\frac{2}{\abs{u}},\huchi,\uomega\right)\angnabla^{i_3}\psi\\
			\notag&+\sum_{k_1+k_2=1}\upsi^{k_1+1}\angnabla^{i_2}(\upeta,\ueta)+\frac{1}{\abs{u}}\upsi\psi+\angnabla(\FRho,\FRhodual,\Scalar)+\frac{1}{\abs{u}}\CBetaA+\upsi(\CRho,\FRho,\Scalar)
		\end{align}
		Let $\vec{m}=(i,j,k)$ with $\abs{\vec{m}}\leq N-1$ and $k\geq 1$. Now, using \eqref{Lbaromega} and \eqref{EQ:D3CBeta}, we commuting the transport equation for $\kappa$ \eqref{PV:topomega3.5} with $\frac{1}{\abs{u}^{k-1}}\WangD^{\vec{m}}$. By \eqref{EQ:ang3ijk}-\eqref{EQ:ang3ijk2}, we obtain
		\begin{align}\label{PV:topomega4}
			&\frac{1}{\abs{u}^{k-1}}\angD_3\WangD^{\vec{m}}\kappa+\frac{1}{\abs{u}^{k-1}}\frac{k}{2}\gtr\uchi\cdot\WangD^{\vec{m}}\kappa\\
			\notag=&~\sum\limits_{\vec{m}_1+\vec{m}_2+\vec{m}_3=(i-i_0,j,k-k_0)}\frac{\delta^{i_0}}{\abs{u}^{k-k_0-1}}\WangD^{\vec{m}_1}\psi^{i_0+1}\cdot\WangD^{\vec{m}_2}\upsi^{k_0}\cdot\WangD^{\vec{m}_3}\angnabla\left(\gtr\uchi+\frac{2}{\abs{u}},\huchi,\uomega\right)\\
			\notag&+\sum\limits_{\vec{m}_1+\vec{m}_2+\vec{m}_3=(i-i_0,j,k-k_0)}\frac{\delta^{i_0}}{\abs{u}^{k-k_0-1}}\WangD^{\vec{m}_1}\omega^{i_0}\cdot\WangD^{\vec{m}_2}\upsi^{k_0+1}\cdot\WangD^{\vec{m}_3}\angnabla(\psi,\upeta,\ueta)\\		\notag&+\sum\limits_{\vec{m}_1+\vec{m}_2=(i-i_0,j,k-k_0)}\frac{\delta^{i_0}}{\abs{u}^{k-k_0}}\WangD^{\vec{m}_1}\psi^{i_0+1}\cdot\WangD^{\vec{m}_2}\upsi^{k_0+1}\\			\notag&+\sum\limits_{\vec{m}_1+\vec{m}_2+\vec{m}_3=(i-i_0,j,k-k_0)}\frac{\delta^{i_0}}{\abs{u}^{k-k_0}}\WangD^{\vec{m}_1}\omega^{i_0}\cdot\WangD^{\vec{m}_2}\upsi^{k_0}\cdot\WangD^{\vec{m}_3}\CBetaA\\
			\notag&+\sum\limits_{\vec{m}_1+\vec{m}_2+\vec{m}_3=(i-i_0,j,k-k_0)}\frac{\delta^{i_0}}{\abs{u}^{k-k_0-1}}\WangD^{\vec{m}_1}\omega^{i_0}\cdot\WangD^{\vec{m}_2}\upsi^{k_0}\cdot\WangD^{\vec{m}_3}\angnabla(\FRho,\FRhodual,\Scalar)\\
			\notag&+\sum\limits_{\vec{m}_1+\vec{m}_2+\vec{m}_3=(i-i_0,j,k-k_0)}\frac{\delta^{i_0}}{\abs{u}^{k-k_0-1}}\WangD^{\vec{m}_1}\omega^{i_0}\cdot\WangD^{\vec{m}_2}\upsi^{k_0+1}\cdot\WangD^{\vec{m}_3}(\CRho,\FRho,\Scalar)\\	
			\notag&+\sum\limits_{\vec{m}_1+\vec{m}_2+\vec{m}_3+\vec{m}_4=(i-i_0,j,k-k_0-1)}\frac{\delta^{i_0}}{\abs{u}^{k-k_0-2}}\WangD^{\vec{m}_1}\omega^{i_0+1}\cdot\WangD^{\vec{m}_2}\upsi^{k_0}\cdot\WangD^{\vec{m}_3}(\Betabar,\FBetabar)\cdot\WangD^{\vec{m}_4}\kappa\\
			\notag&+\sum\limits_{\vec{m}_1+\vec{m}_2+\vec{m}_3+\vec{m}_4=(i-i_0-1,j,k-k_0)}\frac{\delta^{i_0+1}}{\abs{u}^{k-k_0-1}}\WangD^{\vec{m}_1}\omega^{i_0+1}\cdot\WangD^{\vec{m}_2}\upsi^{k_0}\cdot\WangD^{\vec{m}_3}\Rhodual\cdot\WangD^{\vec{m}_4}\kappa,\\
			\notag&=:\frac{1}{\abs{u}^{k-1}}\mathfrak{F}_{(\kappa)}.
		\end{align}
		We apply Proposition \ref{PR:transport2} with $\lambda_0:=\frac{k}{2}$ to \eqref{PV:omega2}:
		\begin{align}\label{PV:topomega5}
			\norm{\WangD^{\vec{m}}\kappa}_{L^2_{\ubar}L_{\sphere}^2(\hypersur_{u})}\lesssim&~\norm{\WangD^{\vec{m}}\kappa}_{L^2_{\ubar}L_{\sphere}^2(\hypersur_{\uinf})}+\int_{\uinf}^{u}\norm{\mathfrak{F}_{(\kappa)}}_{L^2_{\ubar}L_{\sphere}^2(\hypersur_{\uprime})}\diff\uprime.
		\end{align}
		We now estimate $\norm{\mathfrak{F}_{(\kappa)}}_{L^1_uL^2_{\ubar}L_{\sphere}^2(\region)}\diff\uprime$. Compared to the proof of \eqref{ES:L2RicciCoefficients3}, it is suffice to only estimate the top-order terms on the RHS of \eqref{PV:topomega4}. In particular, $i_0=k_0=0$. For the first term on the RHS of \eqref{PV:topomega4}, using bootstrap assumptions \eqref{BA} and proved result \eqref{ES:L2RicciCoefficients}-\eqref{ES:LinftyRicciCoefficientsMore}, we have:
		\begin{align}\label{PV:topomega6}
			&\sum\limits_{\vec{m}_1=(i,j,k+1)}\threenorms{\frac{1}{\abs{u}}\psi\cdot\WangD^{\vec{m}_1}\left(\gtr\uchi+\frac{2}{\abs{u}},\huchi,\uomega\right)}{u}{1}{\ubar}{2}{\sphere}{2}{\region}\\
			\notag\lesssim&~\delta^{\frac{1}{2}}\cdot\sum\limits_{\vec{m}_1=(i,j,k+1)}\threenorms{\WangD^{\vec{m}_1}\left(\gtr\uchi+\frac{2}{\abs{u}},\huchi,\uomega\right)}{u}{2}{\ubar}{\infty}{\sphere}{2}{\region}\cdot\threenorms{\abs{u}\psi}{u}{\infty}{\ubar}{\infty}{\sphere}{\infty}{\region}\cdot\threenorms{u^{-2}}{u}{2}{\ubar}{\infty}{\sphere}{\infty}{\region}\\
			\notag\lesssim&~\delta^{\frac{1}{2}}\cdot\BAsize\frac{\delta a^{\frac{1}{2}}}{\abs{u}^{\frac{1}{2}}}\cdot\frac{ a^{\frac{1}{2}}}{\abs{u}^{\frac{3}{2}}}\\
			\notag\lesssim&~\BAsize\frac{\delta^{\frac{3}{2}}a}{\abs{u}^2}.
		\end{align}
		Then, for the second term on the RHS of \eqref{PV:topomega4}:
		\begin{align}\label{PV:topomega7}
			&\sum\limits_{\vec{m}_1=(i,j,k+1)}\threenorms{\frac{1}{\abs{u}}\WangD^{\vec{m}_1}\psi\cdot\upsi}{u}{1}{\ubar}{2}{\sphere}{2}{\region}\\
			\notag\lesssim&~\sum\limits_{\vec{m}_1=(i,j,k+1)}\threenorms{\WangD^{\vec{m}_1}\psi}{u}{\infty}{\ubar}{2}{\sphere}{2}{\region}\cdot\threenorms{u^2\upsi}{u}{\infty}{\ubar}{\infty}{\sphere}{\infty}{\region}\cdot\threenorms{u^{-3}}{u}{1}{\ubar}{\infty}{\sphere}{\infty}{\region}\\
			\notag\lesssim&~\BAsize\delta^{\frac{1}{2}}a^{\frac{1}{2}}\cdot\BAsize\frac{\delta a^{\frac{1}{2}}}{\abs{u}^2}\cdot\left(1+\frac{1}{\delta^{\frac{1}{2}}a^{\frac{1}{2}}}\sum\limits_{\vec{m}_1=(i,j,k+1)}\norm{\WangD^{\vec{m}_1}\omega^{\dagger}}_{L^{\infty}_{u}L^2_{\ubar}L^2_{\sphere}(\region)}\right)\\
			\notag\lesssim&~\BAsize^2\frac{\delta^{\frac{3}{2}}a}{\abs{u}^2}\left(1+\frac{1}{\delta^{\frac{1}{2}}a^{\frac{1}{2}}}\sum\limits_{\vec{m}_1=(i,j,k+1)}\norm{\WangD^{\vec{m}_1}\omega^{\dagger}}_{L^{\infty}_{u}L^2_{\ubar}L^2_{\sphere}(\region)}\right).
		\end{align}
		
		\begin{align}\label{PV:topomega8}
			&\sum\limits_{\vec{m}_1=(i,j,k+1)}\threenorms{\frac{1}{\abs{u}}\upsi\cdot\WangD^{\vec{m}_1}\left(\upeta,\ueta\right)}{u}{1}{\ubar}{2}{\sphere}{2}{\region}\\
			\notag\lesssim&~\sum\limits_{\vec{m}_1=(i,j,k+1)}\threenorms{\abs{u}\WangD^{\vec{m}_1}\left(\upeta,\ueta\right)}{u}{\infty}{\ubar}{2}{\sphere}{2}{\region}\cdot\threenorms{\abs{u}^2\upsi}{u}{\infty}{\ubar}{\infty}{\sphere}{\infty}{\region}\cdot\threenorms{u^{-4}}{u}{1}{\ubar}{\infty}{\sphere}{\infty}{\region}\\
			\notag\lesssim&~\BAsize\delta^{\frac{3}{2}}a\cdot\BAsize\frac{\delta a^{\frac{1}{2}}}{\abs{u}^{3}}\\
			\notag\lesssim&~\BAsize^2\frac{\delta^{\frac{5}{2}}a^{\frac{3}{2}}}{\abs{u}^3}.
		\end{align}
		Now we estimate the top-order term of the fourth term on the RHS of \eqref{PV:topomega4} using bootstrap assumptions \eqref{BA13} and proved result \eqref{ES:EstimatesforRicciCurvature1.4}:
		\begin{align}\label{PV:topomega9}
			&\sum\limits_{\vec{m}_1=(i,j,k)}\threenorms{\frac{1}{\abs{u}}\cdot\WangD^{\vec{m}_1}\CBetaA}{u}{1}{\ubar}{2}{\sphere}{2}{\region}\\
			\notag\lesssim&~\sum\limits_{\vec{m}_1=(i,j,k)}\threenorms{\abs{u}\WangD^{\vec{m}_1}\BetaA}{u}{\infty}{\ubar}{2}{\sphere}{2}{\region}\cdot\threenorms{u^{-2}}{u}{1}{\ubar}{\infty}{\sphere}{\infty}{\region}+\threenorms{\frac{1}{\abs{u}}\cdot\WangD^{\vec{m}_1}\FBetaA}{u}{1}{\ubar}{2}{\sphere}{2}{\region}\\
			\notag\lesssim&~\frac{\delta^{\frac{1}{2}}a^{\frac{1}{2}}}{\abs{u}}\mathcal{R}+ \frac{1}{\abs{u}}.
		\end{align}
		Then we estimate the top-order fluid term using \eqref{ES:EstimatesforRicciCurvature1.5}:
		\begin{align}\label{PV:topomega10}
			&\sum\limits_{\vec{m}_1=(i,j,k)}\threenorms{\WangD^{\vec{m}_1}\angnabla(\FRho,\FRhodual,\Scalar)}{u}{1}{\ubar}{2}{\sphere}{2}{\region}\\
			\notag\lesssim&~ \delta^{\frac{1}{2}}\abs{u}^{\frac{1}{2}}.
		\end{align}
		The remaining top-order term is the sixth term on the RHS of \eqref{PV:topomega4} for $\Rho$:
		\begin{align}\label{PV:topomega11}
			&\sum\limits_{\vec{m}_1=(i,j,k)}\threenorms{\upsi\cdot\WangD^{\vec{m}_1}\CRho}{u}{1}{\ubar}{2}{\sphere}{2}{\region}\\
			\notag\lesssim&~\delta^{\frac{1}{2}}\cdot\sum\limits_{\vec{m}_1=(i,j,k)}\threenorms{\WangD^{\vec{m}_1}\CRho}{u}{2}{\ubar}{\infty}{\sphere}{2}{\region}\cdot\threenorms{\upsi}{u}{\infty}{\ubar}{\infty}{\sphere}{\infty}{\region}\\
			\notag\lesssim&~\delta^{\frac{1}{2}}\cdot\BAsize\frac{\delta^{\frac{3}{2}}a^{\frac{3}{4}}}{\abs{u}^2}\mathcal{R}\cdot\BAsize\frac{\delta a^{\frac{1}{2}}}{\abs{u}^{2}}\\
			\notag\lesssim&~\BAsize^2\frac{\delta^{3}a^{\frac{5}{4}}}{\abs{u}^3}\mathcal{R}.
		\end{align}
		The rest of the terms on the RHS of \eqref{PV:topomega4} are below top-order, which can be estimated using the same method as in the proof of \eqref{ES:L2RicciCoefficients3}. Combining \eqref{PV:topomega4}-\eqref{PV:topomega11}, we obtain an estimate for $\WangD^{\vec{m}}\kappa$:
		\begin{align}\label{PV:topomega12}
			\threenorms{\WangD^{\vec{m}}\kappa}{u}{\infty}{\ubar}{2}{\sphere}{2}{\region}\lesssim\frac{\delta^{\frac{1}{2}}a^{\frac{1}{2}}}{\abs{u}}(1+\mathcal{R})+\BAsize^2\frac{\delta^{\frac{3}{2}}a}{\abs{u}^2}\left(1+\frac{1}{\delta^{\frac{1}{2}}a^{\frac{1}{2}}}\sum\limits_{\vec{m}_1=(i,j,k+1)}\norm{\WangD^{\vec{m}_1}\omega^{\dagger}}_{L^{\infty}_{u}L^2_{\ubar}L^2_{\sphere}(\region)}\right)+\frac{1}{\abs{u}},
		\end{align}
		Or:
		\begin{align}\label{PV:topomega13}
			\threenorms{\abs{u}\WangD^{\vec{m}}\kappa}{u}{\infty}{\ubar}{2}{\sphere}{2}{\region}\lesssim\delta^{\frac{1}{2}}a^{\frac{1}{2}}\left(1+\mathcal{R}+\frac{1}{\BAsize\delta^{\frac{1}{2}}a^{\frac{1}{2}}}\sum\limits_{\vec{m}_1=(i,j,k+1)}\norm{\WangD^{\vec{m}_1}\omega^{\dagger}}_{L^{\infty}_{u}L^2_{\ubar}L^2_{\sphere}(\region)}\right)+1.
		\end{align}
		We consider the following div-curl systems:
		\begin{subequations}\label{PV:topomega14}
			\begin{align}
				\angdiv\angnabla\omega=&~\angdiv\kappa+\frac{1}{2}\angdiv\CBetaA,&
				\angcurl\angnabla\omega=&~0,\\
				\angdiv\angnabla\omega^{\dagger}=&~0,&
				\angcurl\angnabla\omega^{\dagger}=&~\angcurl\kappa+\frac{1}{2}\angcurl\CBetaA.
			\end{align}
		\end{subequations}
		Let $\vec{m}=(i,j,k)$ with $i+j+k=N-2$. By \eqref{EQ:angAijk}, we commute \eqref{PV:topomega14} with $\WangD^{\vec{m}}$ to obtain:
		\begin{subequations}\label{PV:topomega15}
			\begin{align}
				\angdiv\WangD^{\vec{m}}\angnabla\omega=&~\sum\limits_{\vec{m_1}+\vec{m_2}+\vec{m_3}=(i-i_0,j- j_0,k)}\delta^{i_0}\abs{u}^{j_0}\WangD^{\vec{m_1}}\upchi^{i_0}\cdot\WangD^{\vec{m_2}}\uchi^{ j_0}\cdot\WangD^{\vec{m_3}}\angdiv\kappa\\
				\notag&+\sum\limits_{\vec{m_1}+\vec{m_2}+\vec{m_3}=(i-i_0,j- j_0,k)}\delta^{i_0}\abs{u}^{j_0}\WangD^{\vec{m_1}}\upchi^{i_0}\cdot\WangD^{\vec{m_2}}\uchi^{ j_0}\cdot\WangD^{\vec{m_3}}\angdiv\CBetaA\\
				\notag&+\sum\limits_{\vec{m_1}+\vec{m_2}+\vec{m_3}+\vec{m_4}=(i-i_0,j- j_0,k)}\delta^{i_0}\abs{u}^{j_0}\WangD^{\vec{m_1}}\upchi^{i_0}\cdot\WangD^{\vec{m_2}}\uchi^{ j_0}\cdot\WangD^{\vec{m_3}}(\upeta+\ueta)\cdot\WangD^{\vec{m_4}}\angnabla\omega\\
				\notag&+\sum\limits_{\vec{m_1}+\vec{m_2}+\vec{m_3}+\vec{m_4}=(i-i_0-1,j- j_0,k)}\delta^{i_0}\abs{u}^{j_0}\WangD^{\vec{m_1}}\upchi^{i_0}\cdot\WangD^{\vec{m_2}}\uchi^{ j_0}\cdot\WangD^{\vec{m_3}}(\BetaA,\FBetaA)\cdot\WangD^{\vec{m_4}}\angnabla\omega\\
				\notag&+\sum\limits_{\vec{m_1}+\vec{m_2}+\vec{m_3}+\vec{m_4}=(i-i_0,j- j_0-1,k)}\delta^{i_0}\abs{u}^{j_0}\WangD^{\vec{m_1}}\upchi^{i_0}\cdot\WangD^{\vec{m_2}}\uchi^{ j_0}\cdot\WangD^{\vec{m_3}}(\Betabar,\FBetabar)\cdot\WangD^{\vec{m_4}}\angnabla\omega,\\
				\angcurl\angnabla\omega=&~0,	
			\end{align}
		\end{subequations}
		\begin{subequations}\label{PV:topomega16}
			\begin{align}
				\angdiv\angnabla\omega^{\dagger}=&~0,\\
				\angcurl\WangD^{\vec{m}}\angnabla\omega^{\dagger}=&~\sum\limits_{\vec{m_1}+\vec{m_2}+\vec{m_3}=(i-i_0,j- j_0,k)}\delta^{i_0}\abs{u}^{j_0}\WangD^{\vec{m_1}}\upchi^{i_0}\cdot\WangD^{\vec{m_2}}\uchi^{ j_0}\cdot\WangD^{\vec{m_3}}\angcurl\kappa\\
				\notag&+\sum\limits_{\vec{m_1}+\vec{m_2}+\vec{m_3}=(i-i_0,j- j_0,k)}\delta^{i_0}\abs{u}^{j_0}\WangD^{\vec{m_1}}\upchi^{i_0}\cdot\WangD^{\vec{m_2}}\uchi^{ j_0}\cdot\WangD^{\vec{m_3}}\angcurl\CBetaA\\
				\notag&+\sum\limits_{\vec{m_1}+\vec{m_2}+\vec{m_3}+\vec{m_4}=(i-i_0,j- j_0,k)}\delta^{i_0}\abs{u}^{j_0}\WangD^{\vec{m_1}}\upchi^{i_0}\cdot\WangD^{\vec{m_2}}\uchi^{ j_0}\cdot\WangD^{\vec{m_3}}(\upeta+\ueta)\cdot\WangD^{\vec{m_4}}\angnabla\omega^{\dagger}\\
				\notag&+\sum\limits_{\vec{m_1}+\vec{m_2}+\vec{m_3}+\vec{m_4}=(i-i_0-1,j- j_0,k)}\delta^{i_0}\abs{u}^{j_0}\WangD^{\vec{m_1}}\upchi^{i_0}\cdot\WangD^{\vec{m_2}}\uchi^{ j_0}\cdot\WangD^{\vec{m_3}}(\BetaA,\FBetaA)\cdot\WangD^{\vec{m_4}}\angnabla\omega^{\dagger}\\
				\notag&+\sum\limits_{\vec{m_1}+\vec{m_2}+\vec{m_3}+\vec{m_4}=(i-i_0,j- j_0-1,k)}\delta^{i_0}\abs{u}^{j_0}\WangD^{\vec{m_1}}\upchi^{i_0}\cdot\WangD^{\vec{m_2}}\uchi^{ j_0}\cdot\WangD^{\vec{m_3}}(\Betabar,\FBetabar)\cdot\WangD^{\vec{m_4}}\angnabla\omega^{\dagger}.	
			\end{align}
		\end{subequations}
		We apply elliptic estimates Lemma \ref{LE:ellipticES} to the div-curl systems \eqref{PV:topomega15}-\eqref{PV:topomega16} and taking $L^2_{\ubar}$-norm:
		\begin{align}\label{PV:topomega17}
			\twonorms{\abs{u}^2\WangD^{\vec{m}}\angnabla^2(\omega,\omega^{\dagger})}{\ubar}{2}{\sphere}{2}{\hypersur_{u}}\lesssim\delta^{\frac{1}{2}}a^{\frac{1}{2}}\left(1+\mathcal{R}+\frac{1}{\BAsize\delta^{\frac{1}{2}}a^{\frac{1}{2}}}\sum\limits_{\vec{m}_1=(i,j,k+2)}\norm{\WangD^{\vec{m}_1}\omega^{\dagger}}_{L^{\infty}_{u}L^2_{\ubar}L^2_{\sphere}(\region)}\right)+1,
		\end{align}
		where we note that the largest term is contributed by $\angdiv\kappa,\angcurl\kappa$ terms on the RHS of \eqref{PV:topomega15}-\eqref{PV:topomega16}, while the estimates of  $\angdiv\CBetaA,\angcurl\CBetaA$ terms are in a similar fashion as \eqref{PV:topomega9} and the rest of the terms on the RHS of \eqref{PV:topomega15}-\eqref{PV:topomega16} are below top-order that can be estimated using the same approach as in the proof of \eqref{ES:L2RicciCoefficients3}. 
		Since $\BAsize$ is large, we soak $\norm{\WangD^{\vec{m}_1}\omega^{\dagger}}_{L^{\infty}_{u}L^2_{\ubar}L^2_{\sphere}(\region)}$ to the LHS of \eqref{PV:topomega17} and obtain the desired estimate \eqref{ES:topL2RicciCoefficients3}.
	\end{proof}
	\begin{proof}[Proof of \eqref{ES:topL2RicciCoefficients4}-\eqref{ES:topL2RicciCoefficients5}]
		We define the mass aspect function $\mass$ as follows:
		\begin{align}\label{DE:mass}
			\mass:=-\angdiv\upeta-\CRho.
		\end{align}
		Using \eqref{Lupeta}, \eqref{angdivhchi} and \eqref{EQ:D4Rho2}, we consider the transport equation for $\mass$:
		\begin{align}\label{PV:topupeta2}
			\angD_4\mass=&~\angdiv\FBetaA+\psi\angnabla(\upeta,\ueta)+\psi\upsi^2+\gtr\upchi\hchi\cdot\huchi+\upsi\angnabla(\gtr\upchi,\hchi)+\gtr\upchi(\CRho,\FRho,\Scalar).
		\end{align}
		Let $\vec{m}=(i,j,k)$ with $i+j+k=N-1$ and $k\geq 1$. By \eqref{EQ:ang4ijk}, we commute \eqref{PV:topupeta2} with $\abs{u}^2\WangD^{\vec{m}}$ to obtain:
		\begin{align}\label{PV:topupeta3}
			& \angD_4\abs{u}^2\WangD^{\vec{m}}\mass\\
			\notag=&~ \sum\limits_{\vec{m}_1+\vec{m}_2=(i,j-j_0,k-k_0)}\abs{u}^{j_0+k_0+2}\WangD^{\vec{m}_1}\upsi^{j_0+k_0}\cdot\WangD^{\vec{m}_2}\angdiv\FBetaA\\
			\notag&+ \sum\limits_{\vec{m}_1+\vec{m}_2+\vec{m}_3=(i,j-j_0,k-k_0)}\abs{u}^{j_0+k_0+2}\WangD^{\vec{m}_1}\upsi^{j_0+k_0}\cdot\WangD^{\vec{m}_2}\psi\cdot\WangD^{\vec{m}_3}\angnabla(\upeta,\ueta)\\
			\notag=&~ \sum\limits_{\vec{m}_1+\vec{m}_2=(i,j-j_0,k-k_0)}\abs{u}^{j_0+k_0+2}\WangD^{\vec{m}_1}\upsi^{j_0+k_0+1}\cdot\WangD^{\vec{m}_2}\angnabla(\gtr\upchi,\hchi,\upeta,\ueta)\\
			\notag&+ \sum\limits_{\vec{m}_1+\vec{m}_2=(i,j-j_0,k-k_0)}\abs{u}^{j_0+k_0+2}\WangD^{\vec{m}_1}\upsi^{j_0+k_0+1}\cdot\WangD^{\vec{m}_2}\psi^{2}\\
			\notag&++ \sum\limits_{\vec{m}_1+\vec{m}_2=(i,j-j_0,k-k_0)}\abs{u}^{j_0+k_0+2}\WangD^{\vec{m}_1}\upsi^{j_0+k_0+2}\cdot\WangD^{\vec{m}_2}\psi\\
			\notag&+ \sum\limits_{\vec{m}_1+\vec{m}_2+\vec{m}_3=(i,j-j_0,k-k_0)}\abs{u}^{j_0+k_0+2}\WangD^{\vec{m}_1}\upsi^{j_0+k_0}\cdot\WangD^{\vec{m}_2}\psi\cdot\WangD^{\vec{m}_3}(\CRho,\FRho,\Scalar)\\
			\notag&+ \sum\limits_{\vec{m}_1+\vec{m}_2+\vec{m}_3=(i,j-j_0,k-k_0-1)}\abs{u}^{j_0+k_0+3}\WangD^{\vec{m}_1}\upsi^{k_0+j_0}\cdot\WangD^{\vec{m}_2}(\BetaA,\FBetaA)\cdot\WangD^{\vec{m}_3}\mass\\
			\notag&+ \sum\limits_{\vec{m}_1+\vec{m}_2+\vec{m}_3=(i,j-j_0-1,k-k_0)}\abs{u}^{j_0+k_0+2}\WangD^{\vec{m}_1}\upsi^{j_0+k_0}\cdot\WangD^{\vec{m}_2}\Rhodual\cdot\WangD^{\vec{m}_3}\mass,\\
			\notag&=:\mathfrak{F}_{(\mass)}.
		\end{align}
		By Proposition \ref{PR:transport1}, we have:
		\begin{align}\label{PV:topupeta4}
			\norm{\abs{u}^2\WangD^{\vec{m}}\mass}_{L^2(\suubar)}\lesssim\norm{\abs{u}^2\WangD^{\vec{m}}\mass}_{L^2(S_{u,0})}+\int_{0}^{\ubar}	\norm{\mathfrak{F}_{(\mass)}}_{L^2(\suubarprime)}\diff\ubarprime.
		\end{align}
		We now estimate $\int_{0}^{\ubar}\norm{\mathfrak{F}_{(\mass)}}_{L^2(\suubarprime)}\diff\ubarprime$.
		For the top-order term of the first term on the RHS of \eqref{PV:topupeta3}, by  proved results in Proposition \ref{PR:EstimatesforRicciCurvature}, we have:
		\begin{align}\label{PV:topupeta5.0}
			&\sum\limits_{\vec{m}_1=(i,j,k)}\twonorms{\abs{u}^2\WangD^{\vec{m}_1}\angdiv\FBetaA}{\ubar}{1}{\sphere}{2}{\hypersur_{u}}\\
			\notag\lesssim&~\delta^{\frac{1}{2}}\cdot\abs{u}^2\cdot\twonorms{\WangD^{\vec{m}_1}\angdiv\FBetaA}{\ubar}{2}{\sphere}{2}{\hypersur_{u}}\\
			\notag\lesssim&~\delta^{\frac{1}{2}}\cdot\abs{u}^2.
		\end{align}
		For the top-order term of the second term on the RHS of \eqref{PV:topupeta3}, by bootstrap assumptions \eqref{BA} and proved results \eqref{ES:LinftyRicciCoefficients}, we have:
		\begin{align}\label{PV:topupeta5}
			&\sum\limits_{\vec{m}_1=(i,j,k)}\twonorms{\abs{u}^2\psi\cdot\WangD^{\vec{m}_1}\angnabla(\upeta,\ueta)}{\ubar}{1}{\sphere}{2}{\hypersur_{u}}\\
			\notag\lesssim&~\frac{\delta^{\frac{1}{2}}}{\abs{u}}\cdot\twonorms{\abs{u}\psi}{\ubar}{\infty}{\sphere}{\infty}{\hypersur_{u}}\cdot\twonorms{\abs{u}^2\cdot\WangD^{\vec{m}_1}\angnabla(\upeta,\ueta)}{\ubar}{2}{\sphere}{2}{\hypersur_{u}}\\
			\notag\lesssim&~\frac{\delta^{\frac{1}{2}}}{\abs{u}}\cdot a^{\frac{1}{2}}\cdot \BAsize\delta^{\frac{3}{2}}a\\
			\notag\lesssim&~\delta a\cdot\left(\BAsize\frac{\delta a^{\frac{1}{2}}}{\abs{u}}\right).
		\end{align}
		For the top-order term of the third term on the RHS of \eqref{PV:topupeta3}, by bootstrap assumptions \eqref{BA}:
		\begin{align}\label{PV:topupeta6}
			&\sum\limits_{\vec{m}_1=(i,j,k)}\twonorms{\abs{u}^2\upsi\cdot\WangD^{\vec{m}_1}\angnabla(\gtr\upchi,\hchi,\upeta,\ueta)}{\ubar}{1}{\sphere}{2}{\hypersur_{u}}\\
			\notag\lesssim&~\frac{\delta^{\frac{1}{2}}}{\abs{u}}\cdot\twonorms{\abs{u}^2\upsi}{\ubar}{\infty}{\sphere}{\infty}{\hypersur_{u}}\cdot\twonorms{\abs{u}\cdot\WangD^{\vec{m}_1}\angnabla(\gtr\upchi,\hchi,\upeta,\ueta)}{\ubar}{2}{\sphere}{2}{\hypersur_{u}}\\
			\notag\lesssim&~\frac{\delta^{\frac{1}{2}}}{\abs{u}}\cdot \BAsize\delta a^{\frac{1}{2}}\cdot \BAsize\delta^{\frac{1}{2}}a^{\frac{1}{2}}\\
			\notag\lesssim&~\delta a^{\frac{1}{2}}\left(\BAsize^2\frac{\delta a^{\frac{1}{2}}}{\abs{u}}\right).
		\end{align}
		Now we consider the top-order term of the sixth term for $\CRho$ on the RHS of \eqref{PV:topupeta3} using bootstrap assumptions \eqref{BA} and proved results \eqref{ES:LinftyRicciCoefficients2}, 
		\begin{align}\label{PV:topupeta7}
			&\sum\limits_{\vec{m}_1=(i,j,k)}\twonorms{\abs{u}^2\psi\cdot\WangD^{\vec{m}_1}\CRho}{\ubar}{1}{\sphere}{2}{\hypersur_{u}}\\
			\notag\lesssim&~\frac{\delta^{\frac{1}{2}}}{\abs{u}}\cdot\twonorms{\abs{u}\psi}{\ubar}{\infty}{\sphere}{\infty}{\hypersur_{u}}\cdot\twonorms{\abs{u}^2\cdot\WangD^{\vec{m}_1}\CRho}{\ubar}{2}{\sphere}{2}{\hypersur_{u}}\\
			\notag\lesssim&~\frac{\delta^{\frac{1}{2}}}{\abs{u}}\cdot a^{\frac{1}{2}}\cdot \BAsize\delta^{\frac{3}{2}}a\\
			\notag\lesssim&~\delta a\cdot\left(\BAsize\frac{\delta a^{\frac{1}{2}}}{\abs{u}}\right).
		\end{align}
		Combining \eqref{PV:topupeta3}-\eqref{PV:topupeta7}, noting that rest of the terms are below top-order, we obtain the estimates for $\mass$:
		\begin{align}
			\sum\limits_{\vec{m}=(i,j,k)}\onenorm{\abs{u}^2\WangD^{\vec{m}}\mass}{\sphere}{2}{\suubar}\lesssim&~\delta a.
		\end{align}
		Then, recalling Proposition \ref{PR:NullStructure1}, $\upeta$ satisfies the following div-curl system:
		\begin{align}\label{PV:topupeta1}
			\angdiv\upeta=&~-\mass-\CRho,&
			\angcurl\upeta=&~\CRhodual.
		\end{align}
		Let $\vec{m}=(i,j,k)$ with $i+j+k=N-1$ and $k\geq 1$. By \eqref{EQ:angAijk}, we commute \eqref{PV:topupeta1} with $\abs{u}^2\WangD^{\vec{m}}$ to obtain:
		\begin{subequations}\label{PV:topupeta8}
			\begin{align}
				\angdiv\abs{u}^2\WangD^{\vec{m}}\upeta=&~\sum\limits_{\vec{m_1}+\vec{m_2}+\vec{m_3}=(i-i_0,j- j_0,k)}\delta^{i_0}\abs{u}^{j_0+2}\WangD^{\vec{m_1}}\upchi^{i_0}\cdot\WangD^{\vec{m_2}}\uchi^{ j_0}\cdot\WangD^{\vec{m_3}}\mass\\
				\notag&+\sum\limits_{\vec{m_1}+\vec{m_2}+\vec{m_3}=(i-i_0,j- j_0,k)}\delta^{i_0}\abs{u}^{j_0+2}\WangD^{\vec{m_1}}\upchi^{i_0}\cdot\WangD^{\vec{m_2}}\uchi^{ j_0}\cdot\WangD^{\vec{m_3}}\CRho\\
				\notag&+\sum\limits_{\vec{m_1}+\vec{m_2}+\vec{m_3}=(i-i_0,j- j_0,k)}\delta^{i_0}\abs{u}^{j_0+2}\WangD^{\vec{m_1}}\upchi^{i_0}\cdot\WangD^{\vec{m_2}}\uchi^{ j_0}\cdot\WangD^{\vec{m_3}}(\upeta+\ueta)^2\\
				\notag&+\sum\limits_{\vec{m_1}+\vec{m_2}+\vec{m_3}+\vec{m_4}=(i-i_0-1,j- j_0,k)}\delta^{i_0}\abs{u}^{j_0+2}\WangD^{\vec{m_1}}\upchi^{i_0}\cdot\WangD^{\vec{m_2}}\uchi^{ j_0}\cdot\WangD^{\vec{m_3}}(\BetaA,\FBetaA)\cdot\WangD^{\vec{m_4}}\upeta\\
				\notag&+\sum\limits_{\vec{m_1}+\vec{m_2}+\vec{m_3}+\vec{m_4}=(i-i_0,j- j_0-1,k)}\delta^{i_0}\abs{u}^{j_0+2}\WangD^{\vec{m_1}}\upchi^{i_0}\cdot\WangD^{\vec{m_2}}\uchi^{ j_0}\cdot\WangD^{\vec{m_3}}(\Betabar,\FBetabar)\cdot\WangD^{\vec{m_4}}\upeta,\\
				\angcurl\abs{u}^2\angnabla\ueta=&~\sum\limits_{\vec{m_1}+\vec{m_2}+\vec{m_3}=(i-i_0,j- j_0,k)}\delta^{i_0}\abs{u}^{j_0+2}\WangD^{\vec{m_1}}\upchi^{i_0}\cdot\WangD^{\vec{m_2}}\uchi^{ j_0}\cdot\WangD^{\vec{m_3}}\CRhodual\\
				\notag&+\sum\limits_{\vec{m_1}+\vec{m_2}+\vec{m_3}=(i-i_0,j- j_0,k)}\delta^{i_0}\abs{u}^{j_0+2}\WangD^{\vec{m_1}}\upchi^{i_0}\cdot\WangD^{\vec{m_2}}\uchi^{ j_0}\cdot\WangD^{\vec{m_3}}(\upeta+\ueta)^2\\
				\notag&+\sum\limits_{\vec{m_1}+\vec{m_2}+\vec{m_3}+\vec{m_4}=(i-i_0-1,j- j_0,k)}\delta^{i_0}\abs{u}^{j_0+2}\WangD^{\vec{m_1}}\upchi^{i_0}\cdot\WangD^{\vec{m_2}}\uchi^{ j_0}\cdot\WangD^{\vec{m_3}}(\BetaA,\FBetaA)\cdot\WangD^{\vec{m_4}}\upeta\\
				\notag&+\sum\limits_{\vec{m_1}+\vec{m_2}+\vec{m_3}+\vec{m_4}=(i-i_0,j- j_0-1,k)}\delta^{i_0}\abs{u}^{j_0+2}\WangD^{\vec{m_1}}\upchi^{i_0}\cdot\WangD^{\vec{m_2}}\uchi^{ j_0}\cdot\WangD^{\vec{m_3}}(\Betabar,\FBetabar)\cdot\WangD^{\vec{m_4}}\upeta.	
			\end{align}
		\end{subequations}
		We apply elliptic estimates Lemma \ref{LE:ellipticES} to the div-curl systems \eqref{PV:topupeta8} and taking $L^2_{\ubar}$-norm:
		\begin{align}\label{PV:topupeta9}
			\onenorm{\abs{u}^2\WangD^{\vec{m}}\angnabla\upeta}{\sphere}{2}{\suubar}\lesssim\delta a+\onenorm{\abs{u}^2\WangD^{\vec{m}}(\CRho,\CRhodual)}{\sphere}{2}{\suubar},
		\end{align}
		where we note that the largest term is contributed by $\mass,\CRho,\CRhodual$ terms on the RHS of \eqref{PV:topupeta8}, while the estimates of  $\BetaA,\Betabar$ terms are in a similar fashion as \eqref{PV:topomega9} and the rest of the terms on the RHS of \eqref{PV:topupeta9} are below top-order that can be estimated using the same approach as in the proof of \eqref{ES:L2RicciCoefficients}. 
		Taking $L^2_{\ubar}$-norm and $L^2_{u}$-norm of \eqref{PV:topupeta9} respectively, we obtain the desired results \eqref{ES:topL2RicciCoefficients4}-\eqref{ES:topL2RicciCoefficients5}.
	\end{proof}
	\begin{proof}[Proof of \eqref{ES:topL2RicciCoefficients6}]
		We define the mass aspect function $\umass$ as follows:
		\begin{align}\label{DE:umass}
			\umass:=-\angdiv\ueta-\CRho.
		\end{align}
		Using \eqref{uLueta}, \eqref{angdivhuchi} and \eqref{EQ:D3Rho2}, we consider the transport equation for $\umass$:
		\begin{align}\label{PV:topueta2}
			\angD_3\umass+\frac{3}{2}\gtr\uchi\umass=&~\angdiv\FBetabar+\upsi\angnabla(\upeta,\ueta)+\psi\upsi^2+\upsi\angnabla(\gtr\uchi,\huchi)+\gtr\uchi(\FRho,\Scalar).
		\end{align}
		Let $\vec{m}=(i,j,k)$ with $i+j+k=N-1$ and $k\geq 1$. By \eqref{EQ:ang3ijk}, we commute \eqref{PV:topueta2} with $\frac{1}{\abs{u}^{k+1}}\WangD^{\vec{m}}$ to obtain:
		\begin{align}\label{PV:topueta3}
			& \angD_3\frac{1}{\abs{u}^{k+1}}\WangD^{\vec{m}}\umass+\frac{k+3}{2}\gtr\uchi\frac{1}{\abs{u}^{k+1}}\WangD^{\vec{m}}\umass\\
			\notag=&~\sum\limits_{\vec{m_1}+\vec{m_2}+\vec{m_3}=(i-i_0,j,k-k_0)}\frac{\delta^{i_0}}{\abs{u}^{k-k_0+1}}\WangD^{\vec{m_1}}\upsi^{k_0}\cdot\WangD^{\vec{m_2}}\omega^{i_0}\cdot\WangD^{\vec{m_3}}\angdiv\FBetabar\\
			\notag&+\sum\limits_{\vec{m_1}+\vec{m_2}+\vec{m_3}=(i-i_0,j,k-k_0)}\frac{\delta^{i_0}}{\abs{u}^{k-k_0+1}}\WangD^{\vec{m_1}}\upsi^{k_0+1}\cdot\WangD^{\vec{m_2}}\omega^{i_0}\cdot\WangD^{\vec{m_3}}\angnabla(\gtr\uchi,\huchi,\upeta,\ueta)\\
			\notag&+\sum\limits_{\vec{m_1}+\vec{m_2}=(i-i_0,j,k-k_0)}\frac{\delta^{i_0}}{\abs{u}^{k-k_0+1}}\WangD^{\vec{m_1}}\upsi^{k_0+2}\cdot\WangD^{\vec{m_2}}\psi^{i_0+1}\\
			\notag&+\sum\limits_{\vec{m_1}+\vec{m_2}+\vec{m_3}=(i-i_0,j,k-k_0)}\frac{\delta^{i_0}}{\abs{u}^{k-k_0+1}}\WangD^{\vec{m_1}}\upsi^{k_0+1}\cdot\WangD^{\vec{m_2}}\omega^{i_0}\cdot\WangD^{\vec{m_3}}(\CRho,\FRho,\Scalar)\\
			\notag&+\sum\limits_{\vec{m_1}+\vec{m_2}+\vec{m_3}+\vec{m_4}=(i-i_0,j,k-k_0-1)}\frac{\delta^{i_0}}{\abs{u}^{k-k_0}}\WangD^{\vec{m_1}}(\upeta+\ueta)^{k_0}\cdot\WangD^{\vec{m_2}}\omega^{i_0}\cdot\WangD^{\vec{m_3}}(\Betabar,\FBetabar)\cdot\WangD^{\vec{m_4}}\umass\\
			\notag&+\sum\limits_{\vec{m_1}+\vec{m_2}+\vec{m_3}+\vec{m_4}=(i-i_0-1,j,k-k_0)}\frac{\delta^{i_0+1}}{\abs{u}^{k-k_0+1}}\WangD^{\vec{m_1}}(\upeta+\ueta)^{k_0}\cdot\WangD^{\vec{m_2}}\omega^{i_0}\cdot\WangD^{\vec{m_3}}\Rhodual\cdot\WangD^{\vec{m_4}}\umass\\
			\notag&=:\frac{1}{\abs{u}^{k+2}}\mathfrak{F}_{(\umass)}.
		\end{align}
		We apply Proposition \ref{PR:transport2} with $\lambda_0:=\frac{k+3}{2}$ to \eqref{PV:omega2}:
		\begin{align}\label{PV:topueta4}
			\norm{\abs{u}\WangD^{\vec{m}}\umass}_{L^2(\suubar)}\lesssim&~\norm{\abs{u}\WangD^{\vec{m}}\umass}_{L^2(S_{\uinf,\ubar})}+\int_{\uinf}^{u}\norm{\mathfrak{F}_{(\umass)}}_{L^2(\suprimeubar)}\diff\uprime.
		\end{align}
		Taking $L^2_{\ubar}$-norm of \eqref{PV:topueta4}, we have:
		\begin{align}
			\norm{\abs{u}\WangD^{\vec{m}}\umass}_{L^2(\hypersur_u)}\lesssim&~\norm{\abs{u}\WangD^{\vec{m}}\umass}_{L^2(\hypersur_{\uinf})}+\norm{\mathfrak{F}_{(\umass)}}_{L^2_{\ubar}L^1_uL^2_{\sphere}(\region)}.
		\end{align}
		We now estimate $\norm{\mathfrak{F}_{(\umass)}}_{L^2_{\ubar}L^1_uL^2_{\sphere}(\region)}$.
		For the top-order term of the first term on the RHS of \eqref{PV:topueta3}, by  proved results in Proposition \ref{PR:EstimatesforRicciCurvature}, we have:
		\begin{align}\label{PV:topueta5.0}
			&\sum\limits_{\vec{m}_1=(i,j,k)}\threenorms{\abs{u}\WangD^{\vec{m}_1}\angdiv\FBetabar}{\ubar}{2}{u}{1}{\sphere}{2}{\region}\\
			\notag\lesssim&~\norm{\abs{u}}_{L^1_u}\cdot\threenorms{\WangD^{\vec{m}_1}\angdiv\FBetaA}{\ubar}{2}{u}{\infty}{\sphere}{2}{\region}\\
			\notag\lesssim&~\delta\cdot\abs{u}^{2}.
		\end{align}
		For the top-order term of the second term on the RHS of \eqref{PV:topupeta3}, by bootstrap assumptions \eqref{BA} and proved results \eqref{ES:LinftyRicciCoefficients}, we have:
		\begin{align}\label{PV:topueta5}
			&\sum\limits_{\vec{m}_1=(i,j,k)}\threenorms{\abs{u}\upsi\cdot\WangD^{\vec{m}_1}\angnabla(\upeta,\ueta)}{\ubar}{2}{u}{1}{\sphere}{2}{\region}\\
			\notag\lesssim&~\norm{u^{-3}}_{L^1_u}\threenorms{\abs{u}^2\upsi}{\ubar}{\infty}{u}{\infty}{\sphere}{\infty}{\region}\cdot\threenorms{\abs{u}\cdot\WangD^{\vec{m}_1}\angnabla(\upeta,\ueta)}{\ubar}{2}{u}{\infty}{\sphere}{2}{\region}\\
			\notag\lesssim&~\frac{1}{\abs{u}^2}\cdot \BAsize\delta a^{\frac{1}{2}}\cdot \BAsize\delta^{\frac{3}{2}}a\\
			\notag\lesssim&~\frac{\delta^{\frac{3}{2}}a}{\abs{u}} \cdot\left(\BAsize^2\frac{\delta a^{\frac{1}{2}}}{\abs{u}}\right),
		\end{align}
		\begin{align}\label{PV:topueta6}
			&\sum\limits_{\vec{m}_1=(i,j,k)}\threenorms{\abs{u}\upsi\cdot\WangD^{\vec{m}_1}\angnabla(\gtr\uchi,\huchi)}{\ubar}{2}{u}{1}{\sphere}{2}{\region}\\
			\notag\lesssim&~\delta^{\frac{1}{2}}\norm{u^{-\frac{5}{2}}}_{L^2_u}\threenorms{\abs{u}^2\upsi}{\ubar}{\infty}{u}{\infty}{\sphere}{\infty}{\region}\cdot\threenorms{\abs{u}^{\frac{1}{2}}\cdot\WangD^{\vec{m}_1}\angnabla(\gtr\uchi,\huchi)}{\ubar}{\infty}{u}{2}{\sphere}{2}{\region}\\
			\notag\lesssim&~\delta^{\frac{1}{2}}\cdot\frac{1}{\abs{u}^{2}}\cdot \BAsize\delta a^{\frac{1}{2}}\cdot \BAsize\delta a^{\frac{1}{2}}\\
			\notag\lesssim&~\frac{\delta^{\frac{3}{2}}a^{\frac{1}{2}}}{\abs{u}} \cdot\left(\BAsize^2\frac{\delta a^{\frac{1}{2}}}{\abs{u}}\right),
		\end{align}
		Now we consider the top-order term of the fourth term for $\CRho$ on the RHS of \eqref{PV:topupeta3} using bootstrap assumptions \eqref{BA} and proved results \eqref{ES:LinftyRicciCoefficients2}, 
		\begin{align}\label{PV:topueta7}
			&\sum\limits_{\vec{m}_1=(i,j,k)}\threenorms{\abs{u}\upsi\cdot\WangD^{\vec{m}_1}\CRho}{\ubar}{2}{u}{1}{\sphere}{2}{\region}\\
			\notag\lesssim&~\norm{\abs{u}^{-4}}_{L^1_u}\threenorms{\abs{u}^2\upsi}{\ubar}{\infty}{u}{\infty}{\sphere}{\infty}{\region}\cdot\threenorms{\abs{u}^{2}\cdot\WangD^{\vec{m}_1}\CRho}{\ubar}{2}{u}{\infty}{\sphere}{2}{\region}\\
			\notag\lesssim&~\frac{1}{\abs{u}^{3}}\cdot \BAsize\delta a^{\frac{1}{2}}\cdot \BAsize\delta^{\frac{3}{2}} a\\
			\notag\lesssim&~\frac{\delta^{\frac{3}{2}}a}{\abs{u}^2} \cdot\left(\BAsize^2\frac{\delta a^{\frac{1}{2}}}{\abs{u}}\right),
		\end{align}
		Combining \eqref{PV:topueta3}-\eqref{PV:topueta7}, noting that rest of the terms are below top-order, we obtain the estimates for $\mass$:
		\begin{align}
			\sum\limits_{\vec{m}=(i,j,k)}\twonorms{\abs{u}\WangD^{\vec{m}}\umass}{\ubar}{2}{\sphere}{2}{\hypersur_{u}}\lesssim&~\frac{\delta^{\frac{3}{2}} a}{\abs{u}}.
		\end{align}
		Then, recalling Proposition \ref{PR:NullStructure1}, $\ueta$ satisfies the following div-curl system:
		\begin{align}\label{PV:topueta1}
			\angdiv\ueta=&~-\umass-\CRho,&
			\angcurl\upeta=&~-\CRhodual.
		\end{align}
		Let $\vec{m}=(i,j,k)$ with $i+j+k=N-1$ with $k\geq 1$. By \eqref{EQ:angAijk}, we commute \eqref{PV:topupeta1} with $\abs{u}\WangD^{\vec{m}}$ to obtain:
		\begin{subequations}\label{PV:topueta8}
			\begin{align}
				\angdiv\abs{u}\WangD^{\vec{m}}\ueta=&~\sum\limits_{\vec{m_1}+\vec{m_2}+\vec{m_3}=(i-i_0,j- j_0,k)}\abs{u}^2\WangD^{\vec{m_1}}\upchi^{i_0}\cdot\WangD^{\vec{m_2}}\uchi^{ j_0}\cdot\WangD^{\vec{m_3}}\umass\\
				\notag&+\sum\limits_{\vec{m_1}+\vec{m_2}+\vec{m_3}=(i-i_0,j- j_0,k)}\abs{u}\WangD^{\vec{m_1}}\upchi^{i_0}\cdot\WangD^{\vec{m_2}}\uchi^{ j_0}\cdot\WangD^{\vec{m_3}}\CRho\\
				\notag&+\sum\limits_{\vec{m_1}+\vec{m_2}+\vec{m_3}=(i-i_0,j- j_0,k)}\abs{u}\WangD^{\vec{m_1}}\upchi^{i_0}\cdot\WangD^{\vec{m_2}}\uchi^{ j_0}\cdot\WangD^{\vec{m_3}}(\upeta+\ueta)^2\\
				\notag&+\sum\limits_{\vec{m_1}+\vec{m_2}+\vec{m_3}+\vec{m_4}=(i-i_0-1,j- j_0,k)}\abs{u}\WangD^{\vec{m_1}}\upchi^{i_0}\cdot\WangD^{\vec{m_2}}\uchi^{ j_0}\cdot\WangD^{\vec{m_3}}(\BetaA,\FBetaA)\cdot\WangD^{\vec{m_4}}\ueta\\
				\notag&+\sum\limits_{\vec{m_1}+\vec{m_2}+\vec{m_3}+\vec{m_4}=(i-i_0,j- j_0-1,k)}\abs{u}\WangD^{\vec{m_1}}\upchi^{i_0}\cdot\WangD^{\vec{m_2}}\uchi^{ j_0}\cdot\WangD^{\vec{m_3}}(\Betabar,\FBetabar)\cdot\WangD^{\vec{m_4}}\ueta,\\
				\angcurl\abs{u}\angnabla\ueta=&~\sum\limits_{\vec{m_1}+\vec{m_2}+\vec{m_3}=(i-i_0,j- j_0,k)}\abs{u}\WangD^{\vec{m_1}}\upchi^{i_0}\cdot\WangD^{\vec{m_2}}\uchi^{ j_0}\cdot\WangD^{\vec{m_3}}\CRhodual\\
				\notag&+\sum\limits_{\vec{m_1}+\vec{m_2}+\vec{m_3}=(i-i_0,j- j_0,k)}\abs{u}\WangD^{\vec{m_1}}\upchi^{i_0}\cdot\WangD^{\vec{m_2}}\uchi^{ j_0}\cdot\WangD^{\vec{m_3}}(\upeta+\ueta)^2\\
				\notag&+\sum\limits_{\vec{m_1}+\vec{m_2}+\vec{m_3}+\vec{m_4}=(i-i_0-1,j- j_0,k)}\abs{u}\WangD^{\vec{m_1}}\upchi^{i_0}\cdot\WangD^{\vec{m_2}}\uchi^{ j_0}\cdot\WangD^{\vec{m_3}}(\BetaA,\FBetaA)\cdot\WangD^{\vec{m_4}}\ueta\\
				\notag&+\sum\limits_{\vec{m_1}+\vec{m_2}+\vec{m_3}+\vec{m_4}=(i-i_0,j- j_0-1,k)}\abs{u}\WangD^{\vec{m_1}}\upchi^{i_0}\cdot\WangD^{\vec{m_2}}\uchi^{ j_0}\cdot\WangD^{\vec{m_3}}(\Betabar,\FBetabar)\cdot\WangD^{\vec{m_4}}\ueta.	
			\end{align}
		\end{subequations}
		We apply elliptic estimates Lemma \ref{LE:ellipticES} to the div-curl systems \eqref{PV:topupeta8} and taking $L^2_{\ubar}$-norm:
		\begin{align}\label{PV:topueta9}
			\twonorms{\abs{u}^2\WangD^{\vec{m}}\angnabla\upeta}{\ubar}{2}{\sphere}{2}{\hypersur_{u}}\lesssim\delta^{\frac{3}{2}} a+\twonorms{\abs{u}^2\WangD^{\vec{m}}(\CRho,\CRhodual)}{\ubar}{2}{\sphere}{2}{\hypersur_{u}},
		\end{align}
		where we note that the largest term is contributed by $\umass,\CRho,\CRhodual$ terms on the RHS of \eqref{PV:topupeta8}, while the estimates of  $\BetaA,\Betabar$ terms are in a similar fashion as \eqref{PV:topomega9} and the rest of the terms on the RHS of \eqref{PV:topueta9} are below top-order that can be estimated using the same approach as in the proof of \eqref{ES:L2RicciCoefficients}.  We have obtained the desired results \eqref{ES:topL2RicciCoefficients6}.
	\end{proof}
	
	\section{Estimates for the Weyl Tensor Components}\label{S:Curvature}
	\subsection{Integration by Parts Formula}\label{SS:IntegrationByParts}
	\begin{lemma}\cite[Proposition 8.1]{AnLuk2017}.\label{LE:Int1}
		Let $F$ and $G$ be $\suubar$-tensorfields. We have the following integration identity:
		\begin{align}
			&\int_{\region}F\cdot\angD_4G\diff\gvol+\int_{\region}G\cdot\angD_4F\diff\gvol\\
			\notag=&~\int_{\hypersurbar_{\ubar}}F\cdot G\diff u\diff\spherevol-\int_{\hypersurbar_{0}}F\cdot G\diff u\diff\spherevol+\int_{\region}(2\omega-\gtr\upchi)F\cdot G\diff \gvol.
		\end{align}
	\end{lemma}
	\begin{lemma}\cite[Proposition 8.2]{AnLuk2017}.\label{LE:Int2}
		Let $F$ be $\suubar$-tangent $r$-tensorfield and $G$ be $\suubar$-tangent $r-1$-tensorfield. We have the following integration identity:
		\begin{align}
			\int_{\region}F\cdot\angnabla G\diff\gvol+\int_{\region}G\cdot\angdiv F\diff\gvol=-\int_{\region}(\upeta+\ueta)\cdot F\cdot G\diff \gvol.
		\end{align}
	\end{lemma}
	\begin{lemma}\cite[Proposition 8.3]{AnLuk2017}.\label{LE:Int3}
		Let $F$ be a $\suubar$-tensorfield. Let $\lambda_1=2(\lambda_0-\frac{1}{2})$. We have the following integration identity:
		\begin{align}
			&2\int_{\region}\abs{u}^{2\lambda_1}F\cdot(\angD_3+\lambda_0\gtr\uchi)F\diff\gvol\\
			\notag=&~\int_{\hypersur_{u}}\abs{u}^{2\lambda_1}\sgabs{F}^2\diff \ubar\diff\spherevol-\int_{\hypersur_{\uinf}}\abs{u}^{2\lambda_1}\sgabs{F}^2\diff \ubar\diff\spherevol+\int_{\region}\abs{u}^{2\lambda_1}\Chfour\cdot\sgabs{F}^2\diff \gvol,
		\end{align}
		where $\Chfour$ follows the estimate:
		\begin{align}
			\abs{\Chfour}\lesssim\BAsize\frac{\delta a^{\frac{1}{2}}}{\abs{u}^2}.
		\end{align}
	\end{lemma}
	\begin{proposition}\cite[Section 7.2]{AnAthanasiou2022}.\label{PR:Int4}
		Let $\vec{m}=(i,j,k)$ with $i+j+k=N-1$ with $k\geq 1$. Let $\WangD^{\vec{m}}$ be as in Definition \ref{DE:operators}. Let $F$ and $G$ be $\suubar$-tensorfields. Let $s(F)$ be the ``signature\footnote{In the analysis, we substitute $F,G$ by Weyl tensor components, and signature is designed to count the number difference in $e_3$ and $e_4$ derivatives (of the metric) that the Weyl tensor component has; see \cite{KlainermanRodnianski2012} and \cite{AnAthanasiou2022}.} function" of $F$. We assume pair $(\WangD^{\vec{m}}F,\WangD^{\vec{m}}G)$ satisfies the following:
		\begin{subequations}\label{EQ:Int4}
			\begin{align}
				\angD_3\WangD^{\vec{m}}F+\left(\frac{1+k}{2}+s(F)\right)\gtr\uchi\WangD^{\vec{m}}F-\mathcal{D}\WangD^{\vec{m}}G=&~P_{(F,G;\vec{m})},\\
				\angD_4\WangD^{\vec{m}}G-\HodgeLeftDual\mathcal{D}\WangD^{\vec{m}}F=&~Q_{(F,G;\vec{m})},
			\end{align}
		\end{subequations}
		where $\mathcal{D}$ is an angular derivative operator on $\suubar$ and $\HodgeLeftDual\mathcal{D}$ is its $L^2$-adjoint from a Hodge system. Then, we have the following estimates:
		\begin{align}\label{ES:Int4}
			&\twonorms{\abs{u}^{2s(F)}\WangD^{\vec{m}}F}{\ubar}{2}{\sphere}{2}{\hypersur_{u}}^2+\twonorms{\abs{u}^{2s(F)}\WangD^{\vec{m}}G}{u}{2}{\sphere}{2}{\hypersurbar_{\ubar}}^2\\
			\notag\lesssim&~\twonorms{\abs{u}^{2s(F)}\WangD^{\vec{m}}F}{\ubar}{2}{\sphere}{2}{\hypersur_{\uinf}}^2+\twonorms{\abs{u}^{2s(F)}\WangD^{\vec{m}}G}{u}{2}{\sphere}{2}{\hypersurbar_{0}}^2\\
			\notag&+\int_{\region}\abs{u}^{4s(F)}\abs{\WangD^{\vec{m}}F}\cdot \abs{P_{(F,G;\vec{m})}}\diff\gvol+\int_{\region}\abs{u}^{4s(F)}\abs{\WangD^{\vec{m}}G}\cdot \abs{Q_{(F,G;\vec{m})}}\diff\gvol.
		\end{align}
	\end{proposition}
	\subsection{Energy Estimates for the Renormalized Weyl Tensor Components}
	\begin{proposition}[Below top-order energy estimates for the renormalized Weyl tensor components]\label{PR:L2curvature}
		Let renormalized Weyl tensor components $\CAlphaAB,\CBetaA,\CRho,\CRhodual,\CBetabar,\CAlphabar$ be defined as in Definition \ref{DE:Renormalization}. Let $\vec{m}=(i,j,k)$ with $i+j+k=N-2$. Let $\WangD^{\vec{m}}$ be as in Definition \ref{DE:operators}. Under the initial data assumptions \eqref{Data} and bootstrap assumptions \eqref{BA}, we have the following top-order energy estimates for the renormalized Weyl tensor components\footnote{Similar to Proposition \ref{PR:L2Fluid} and Proposition \ref{PR:L2RicciCoefficients}, the first term on the RHS of \eqref{ES:L2curvature} represents the largest term and the second term represents the largest impact from the fluid.}:
		\begin{subequations}\label{ES:L2curvature}
			\begin{align}\label{ES:L2curvature1}
				\sum\limits_{\abs{\vec{m}}\leq N-2}\norm{\WangD^{\vec{m}}\CAlphaAB}_{L^2(\suubar)}&\lesssim \delta^{-1}a^{\frac{1}{4}}(1+\tilde{\mathcal{R}}_{N-1,2})+ \delta^{-1},\\
				\label{ES:L2curvature2}\sum\limits_{\abs{\vec{m}}\leq N-2}\norm{\abs{u}\WangD^{\vec{m}}\CBetaA}_{L^2(\suubar)}&\lesssim \delta^{\frac{1}{2}}a^{\frac{1}{2}}(\delta^{-\frac{1}{2}}+\tilde{\mathcal{R}}_{N-1,2})+ b^{\frac{1}{2}},\\
				\label{ES:L2curvature3}\sum\limits_{\abs{\vec{m}}\leq N-2}\norm{\abs{u}^2\WangD^{\vec{m}}(\CRho,\CRhodual,\CBetabar)}_{L^2(\suubar)}&\lesssim \delta^{\frac{3}{2}} a(\delta^{-\frac{1}{2}}a^{-\frac{1}{2}}+\tilde{\mathcal{R}}_{N-1,2})+ b^{-\frac{1}{2}},\\
				\label{ES:L2curvature4}\sum\limits_{\abs{\vec{m}}\leq N-2}\norm{\abs{u}^{2}\WangD^{\vec{m}}\CAlphabar}_{L^2(\suubar)}&\lesssim \delta^{2}a \left(a^{-\frac{1}{2}}+\tilde{\mathcal{R}}_{N-1,2}\right)+ b^{-1}.
			\end{align}
		\end{subequations}
	\end{proposition}
	We prove Proposition \ref{PR:L2curvature} in Section \ref{SSS:proofofcurvature0}. 
	\begin{corollary}[Below top-order energy estimates for the renormalized Weyl tensor components]\label{CO:Linftycurvature}
		Let renormalized Weyl tensor components $\CAlphaAB,\CBetaA,\CRho,\CRhodual,\CBetabar,\CAlphabar$ be defined as in Definition \ref{DE:Renormalization}. Let $\vec{m}=(i,j,k)$ with $i+j+k=N-2$. Let $\WangD^{\vec{m}}$ be as in Definition \ref{DE:operators}. Under the initial data assumptions \eqref{Data} and bootstrap assumptions \eqref{BA}, we have the following top-order energy estimates for the renormalized Weyl tensor components\footnote{Similar to Proposition \ref{PR:L2Fluid} and Proposition \ref{PR:L2RicciCoefficients}, the first term on the RHS of \eqref{ES:Linftycurvature} represents the largest term and the second term represents the largest impact from the fluid.}:
		\begin{subequations}\label{ES:Linftycurvature}
			\begin{align}\label{ES:Linftycurvature1}
				\sum\limits_{\abs{\vec{m}}\leq N-4}\norm{\WangD^{\vec{m}}\CAlphaAB}_{L^\infty(\suubar)}&\lesssim \delta^{-1}a^{\frac{1}{4}}(1+\tilde{\mathcal{R}}_{N-1,2})+ \delta^{-1},\\
				\label{ES:Linftycurvature2}\sum\limits_{\abs{\vec{m}}\leq N-4}\norm{\abs{u}\WangD^{\vec{m}}\CBetaA}_{L^\infty(\suubar)}&\lesssim \delta^{\frac{1}{2}}a^{\frac{1}{2}}(\delta^{-\frac{1}{2}}+\tilde{\mathcal{R}}_{N-1,2})+ b^{\frac{1}{2}},\\
				\label{ES:Linftycurvature3}\sum\limits_{\abs{\vec{m}}\leq N-4}\norm{\abs{u}^2\WangD^{\vec{m}}(\CRho,\CRhodual,\CBetabar)}_{L^\infty(\suubar)}&\lesssim \delta^{\frac{3}{2}} a(\delta^{-\frac{1}{2}}a^{-\frac{1}{2}}+\tilde{\mathcal{R}}_{N-1,2})+ b^{-\frac{1}{2}},\\
				\label{ES:Linftycurvature4}\sum\limits_{\abs{\vec{m}}\leq N-4}\norm{\abs{u}^{2}\WangD^{\vec{m}}\CAlphabar}_{L^\infty(\suubar)}&\lesssim \delta^{2}a \left(a^{-\frac{1}{2}}+\tilde{\mathcal{R}}_{N-1,2}\right)+ b^{-1}.
			\end{align}
		\end{subequations}
	\end{corollary}
	\begin{proposition}[Top-order energy estimates for the renormalized Weyl tensor components]\label{PR:topL2curvature}
		Let renormalized Weyl tensor components $\CAlphaAB,\CBetaA,\CRho,\CRhodual,\CBetabar,\CAlphabar$ be defined as in Definition \ref{DE:Renormalization}. Let $\vec{m}=(i,j,k)$ with $i+j+k=N-1$ with $k\geq 1$. Let $\WangD^{\vec{m}}$ be as in Definition \ref{DE:operators}. Under the initial data assumptions \eqref{Data} and bootstrap assumptions \eqref{BA}, we have the following top-order energy estimates for the renormalized Weyl tensor components:
		\begin{subequations}\label{ES:topL2curvature}
			\begin{align}\label{ES:topL2curvature1}
				\sum\limits_{\substack{\abs{\vec{m}}=N-1\\ k\geq 1}}\left(\norm{\WangD^{\vec{m}}\CAlphaAB}_{L^2(\hypersur_{u})}+\norm{\WangD^{\vec{m}}\CBetaA}_{L^2(\hypersurbar_{\ubar})}\right)&\lesssim \delta^{\frac{1}{2}} a^{\frac{3}{4}}+ \left(\delta^{-1} b^{-\frac{1}{2}}+a^{\frac{1}{2}}\right),\\
				\label{ES:topL2curvature2}\sum\limits_{\substack{\abs{\vec{m}}=N-1\\ k\geq 1}}\left(\norm{\abs{u}\WangD^{\vec{m}}\CBetaA}_{L^2(\hypersur_{u})}+\norm{\abs{u}\WangD^{\vec{m}}(\CRho,\CRhodual)}_{L^2(\hypersurbar_{\ubar})}\right)&\lesssim \delta^{\frac{3}{2}} a^{\frac{3}{2}}+ \left( b^{-\frac{1}{2}}+\delta^{\frac{1}{2}}\right),\\
				\label{ES:topL2curvature3}\sum\limits_{\substack{\abs{\vec{m}}=N-1\\ k\geq 1}}\left(\norm{\abs{u}^2\WangD^{\vec{m}}(\CRho,\CRhodual)}_{L^2(\hypersur_{u})}+\norm{\abs{u}^2\WangD^{\vec{m}}\CBetabar}_{L^2(\hypersurbar_{\ubar})}\right)&\lesssim \delta^{2} a^{\frac{3}{2}}+ b^{-\frac{1}{2}},\\
				\label{ES:topL2curvature4}\sum\limits_{\substack{\abs{\vec{m}}=N-1\\ k\geq 1}}\left(\norm{\abs{u}^3\WangD^{\vec{m}}\CBetabar}_{L^2(\hypersur_{u})}+\norm{\abs{u}^3\WangD^{\vec{m}}\CAlphabar}_{L^2(\hypersurbar_{\ubar})}\right)&\lesssim \delta^{\frac{5}{2}} a^{\frac{3}{2}}+ \delta^{\frac{1}{2}}b^{-\frac{1}{2}}.
			\end{align}
		\end{subequations}
	\end{proposition}
	We prove Proposition \ref{PR:topL2curvature} in Section \ref{SSS:proofofcurvature}.
	\subsubsection{Proof of Proposition \ref{PR:L2curvature}}\label{SSS:proofofcurvature0}
	In this section, we only highlight the analysis that are not present in Section \ref{SS:proveRicci}.
	\begin{proof}[Proof of \eqref{ES:L2curvature1}]
		Let $\vec{m}=(i,j,k)$ with $i+j+k\leq N-2$. Let $\WangD^{\vec{m}}$ be as in Definition \ref{DE:operators}. We commute equations \eqref{EQ:D3CAlpha} with $\frac{1}{\abs{u}^{k}}\WangD^{\vec{m}}$ using \eqref{EQ:ang3ijk}-\eqref{EQ:ang3ijk2}, we obtain
		\begin{align}\label{PV:alpha2}
			&\frac{1}{\abs{u}^k}\angD_3\WangD^{\vec{m}}\CAlphaAB+\frac{1}{\abs{u}^k}\frac{k+1}{2}\gtr\uchi\cdot\WangD^{\vec{m}}\CAlphaAB\\
			\notag=&~\sum\limits_{\vec{m_1}+\vec{m_2}+\vec{m_3}=(i-i_0,j,k-k_0)}\frac{\delta^{i_0}}{\abs{u}^{k-k_0}}\WangD^{\vec{m_1}}\upsi^{k_0}\cdot\WangD^{\vec{m_2}}\psi^{i_0}\cdot\WangD^{\vec{m_3}}\angnabla\CBetaA\\
			\notag&+\sum\limits_{\vec{m_1}+\vec{m_2}+\vec{m_3}=(i-i_0,j,k-k_0)}\frac{\delta^{i_0}}{\abs{u}^{k-k_0}}\WangD^{\vec{m_1}}\upsi^{k_0+1}\cdot\WangD^{\vec{m_2}}\psi^{i_0}\cdot\WangD^{\vec{m_3}}(\CAlphaAB,\CBetaA,\FAlphaAB)\\
			\notag&+\sum\limits_{\vec{m_1}+\vec{m_2}+\vec{m_3}=(i-i_0,j,k-k_0)}\frac{\delta^{i_0}}{\abs{u}^{k-k_0}}\WangD^{\vec{m_1}}\upsi^{k_0}\cdot\WangD^{\vec{m_2}}\psi^{i_0+1}\cdot\WangD^{\vec{m_3}}(\CRho,\CRhodual,\Scalar,\FRho,\FRhodual)\\
			\notag&+\sum\limits_{\vec{m_1}+\vec{m_2}+\vec{m_3}=(i-i_0,j,k-k_0)}\frac{\delta^{i_0}}{\abs{u}^{k-k_0}}\WangD^{\vec{m_1}}\upsi^{k_0}\cdot\WangD^{\vec{m_2}}\psi^{i_0}\cdot\WangD^{\vec{m_3}}(\angnabla\FBetaA,\angD_4\Scalar)\\
			\notag&+\sum\limits_{\vec{m_1}+\vec{m_2}+\vec{m_3}+\vec{m_4}=(i-i_0-1,j,k-k_0+1)}\frac{\delta^{i_0+1}}{\abs{u}^{k-k_0+1}}\WangD^{\vec{m_1}}\upsi^{k_0}\cdot\WangD^{\vec{m_2}}\psi^{i_0}\cdot\WangD^{\vec{m_3}}\upzeta\cdot\WangD^{\vec{m_4}}\CAlphaAB\\
			\notag&+\sum\limits_{\vec{m_1}+\vec{m_2}+\vec{m_3}+\vec{m_4}=(i-i_0-1,j,k-k_0)}\frac{\delta^{i_0+1}}{\abs{u}^{k-k_0}}\WangD^{\vec{m_1}}\upsi^{k_0}\cdot\WangD^{\vec{m_2}}\psi^{i_0}\cdot\WangD^{\vec{m_3}}\Rhodual\cdot\WangD^{\vec{m_4}}\CAlphaAB\\
			\notag&+\sum\limits_{\vec{m_1}+\vec{m_2}+\vec{m_3}+\vec{m_4}=(i-i_0,j,k-k_0-1)}\frac{\delta^{i_0}}{\abs{u}^{k-k_0-1}}\WangD^{\vec{m_1}}\upsi^{k_0}\cdot\WangD^{\vec{m_2}}\psi^{i_0}\cdot\WangD^{\vec{m_3}}(\CBetabar,\FBetabar)\cdot\WangD^{\vec{m_4}}\CAlphaAB\\
			\notag&=:\frac{1}{\abs{u}^k}\mathfrak{F}_{(\CAlphaAB)}.
		\end{align}
		We apply Proposition \ref{PR:transport2} with $\lambda_0:=\frac{k+1}{2}$ to \eqref{PV:alpha2}:
		\begin{align}\label{PV:alpha3}
			\norm{\WangD^{\vec{m}}\CAlphaAB}_{L^2(\suubar)}\lesssim&~\norm{\WangD^{\vec{m}}\CAlphaAB}_{L^2(S_{\uinf,\ubar})}+\int_{\uinf}^{u}\norm{\mathfrak{F}_{(\CAlphaAB)}}_{L^2(\suprimeubar)}\diff\uprime.
		\end{align}
			We now estimate $\norm{\mathfrak{F}_{(\CAlphaAB)}}_{L^1_uL^2_{\sphere}(\hypersurbar_{\ubar})}$. For the first term on the RHS of \eqref{PV:alpha2}, by bootstrap assumptions \eqref{BA}, we have:
			\begin{align}\label{PV:alpha4}
				&\sum\limits_{\vec{m}_1+\vec{m}_2+\vec{m}_3=(i-i_0,j,k-k_0)}\delta^{i_0}\twonorms{\abs{u}^{k_0}\WangD^{\vec{m}_1}\psi^{i_0}\cdot\WangD^{\vec{m}_2}\upsi^{k_0}\cdot\WangD^{\vec{m}_3}\angnabla\CBetaA}{u}{1}{\sphere}{2}{\hypersurbar_{\ubar}}\\
				\notag\lesssim&~\sum\limits_{\substack{\vec{m}_1+\vec{m}_2+\vec{m}_3=(i-i_0,j,k-k_0)\\\abs{\vec{m}_1+\vec{m}_2}>N-3}}\delta^{i_0}\cdot\left(\twonorms{\abs{u}^{k_0}\WangD^{\vec{m}_1}\psi^{i_0}\cdot\WangD^{\vec{m}_2}\upsi^{k_0}}{u}{1}{\sphere}{2}{\hypersurbar_{\ubar}}\cdot\twonorms{\WangD^{\vec{m}_3}\angnabla\CBetaA}{u}{\infty}{\sphere}{\infty}{\hypersurbar_{\ubar}}\right)\\
				\notag&+\sum\limits_{\substack{\vec{m}_1+\vec{m}_2+\vec{m}_3=(i-i_0,j,k-k_0)\\\abs{\vec{m}_1+\vec{m}_2}\leq N-3}}\delta^{i_0}\left(\twonorms{\abs{u}^{k_0}\WangD^{\vec{m}_1}\psi^{i_0}\cdot\WangD^{\vec{m}_2}\upsi^{k_0}}{u}{\infty}{\sphere}{\infty}{\hypersurbar_{\ubar}}\cdot\twonorms{\WangD^{\vec{m}_3}\angnabla\CBetaA}{u}{1}{\sphere}{2}{\hypersurbar_{\ubar}}\right)\\
				\notag\lesssim&~\left(\BAsize\frac{\delta a^{\frac{1}{2}}}{\abs{u}}\right)^{i_0+k_0}\cdot\delta^{-\frac{1}{2}}a^{\frac{1}{4}}(1+\tilde{\mathcal{R}}_{N-1,2}).
			\end{align}
			For the second term on the RHS of \eqref{PV:alpha2}, by Sobolev embedding $H^2\hookrightarrow L^{\infty}$ for $\CAlphaAB$ on $\suubar$, bootstrap assumptions \eqref{BA} and proved results in Proposition \ref{PR:EstimatesforRicciCurvature}, we have:
			\begin{align}\label{PV:alpha4.2}
				&\sum\limits_{\vec{m}_1+\vec{m}_2+\vec{m}_3=(i-i_0,j,k-k_0)}\delta^{i_0}\twonorms{\abs{u}^{k_0}\WangD^{\vec{m}_1}\psi^{i_0}\cdot\WangD^{\vec{m}_2}\upsi^{k_0+1}\cdot\WangD^{\vec{m}_3}(\CAlphaAB,\CBetaA,\FAlphaAB)}{u}{1}{\sphere}{2}{\hypersurbar_{\ubar}}\\
				\notag\lesssim&~\sum\limits_{\substack{\vec{m}_1+\vec{m}_2+\vec{m}_3=(i-i_0,j,k-k_0)\\\abs{\vec{m}_1+\vec{m}_2}>N-3}}\delta^{i_0}\cdot\left[\twonorms{\abs{u}^{k_0}\WangD^{\vec{m}_1}\psi^{i_0}\cdot\WangD^{\vec{m}_2}\upsi^{k_0+1}}{u}{1}{\sphere}{2}{\hypersurbar_{\ubar}}\right.\\
				\notag&\left.\cdot\left(\twonorms{\WangD^{\vec{m}_3}(\CBetaA,\FAlphaAB)}{u}{\infty}{\sphere}{\infty}{\hypersurbar_{\ubar}}+\twonorms{\WangD^{\vec{m}_3}\angnabla^2\CAlphaAB}{u}{\infty}{\sphere}{2}{\hypersurbar_{\ubar}}\right)\right]\\
				\notag&+\sum\limits_{\substack{\vec{m}_1+\vec{m}_2+\vec{m}_3=(i-i_0,j,k-k_0)\\\abs{\vec{m}_1+\vec{m}_2}\leq N-3}}\delta^{i_0}\left(\twonorms{\abs{u}^{k_0}\WangD^{\vec{m}_1}\psi^{i_0}\cdot\WangD^{\vec{m}_2}\upsi^{k_0+1}}{u}{\infty}{\sphere}{\infty}{\hypersurbar_{\ubar}}\cdot\twonorms{\WangD^{\vec{m}_3}(\CAlphaAB,\CBetaA,\FAlphaAB)}{u}{1}{\sphere}{2}{\hypersurbar_{\ubar}}\right)\\
				\notag\lesssim&~\left(\BAsize\frac{\delta a^{\frac{1}{2}}}{\abs{u}}\right)^{i_0+k_0+1}\cdot\sum\limits_{\abs{\vec{m}}\leq N-2}\norm{\WangD^{\vec{m}}\CAlphaAB}_{L^2(\suubar)}+   b\cdot\abs{u}^{\frac{1}{2}}\cdot\left(\BAsize\frac{\delta a^{\frac{1}{2}}}{\abs{u}}\right)^{i_0+k_0+1}.
			\end{align}
			For the fourth term on the RHS of \eqref{PV:alpha2}, by bootstrap assumptions \eqref{BA} and proved results in Proposition \ref{PR:EstimatesforRicciCurvature}, we have:
			\begin{align}\label{PV:alpha5}
				&\sum\limits_{\vec{m}_1+\vec{m}_2+\vec{m}_3=(i-i_0,j,k-k_0)}\delta^{i_0}\twonorms{\abs{u}^{k_0}\WangD^{\vec{m}_1}\psi^{i_0}\cdot\WangD^{\vec{m}_2}\upsi^{k_0}\cdot\WangD^{\vec{m}_3}(\angnabla\FBetaA,\angD_{4}\Scalar)}{u}{1}{\sphere}{2}{\hypersurbar_{\ubar}}\\
				\notag\lesssim&~\sum\limits_{\substack{\vec{m}_1+\vec{m}_2+\vec{m}_3=(i-i_0,j,k-k_0)\\\abs{\vec{m}_1+\vec{m}_2}>N-3}}\delta^{i_0}\cdot\left(\twonorms{\abs{u}^{k_0}\WangD^{\vec{m}_1}\psi^{i_0}\cdot\WangD^{\vec{m}_2}\upsi^{k_0}}{u}{1}{\sphere}{2}{\hypersurbar_{\ubar}}\cdot\twonorms{\WangD^{\vec{m}_3}(\angnabla\FBetaA,\angD_{4}\Scalar)}{u}{\infty}{\sphere}{\infty}{\hypersurbar_{\ubar}}\right)\\
				\notag&+\sum\limits_{\substack{\vec{m}_1+\vec{m}_2+\vec{m}_3=(i-i_0,j,k-k_0)\\\abs{\vec{m}_1+\vec{m}_2}\leq N-3}}\delta^{i_0}\left(\twonorms{\abs{u}^{k_0}\WangD^{\vec{m}_1}\psi^{i_0}\cdot\WangD^{\vec{m}_2}\upsi^{k_0}}{u}{\infty}{\sphere}{\infty}{\hypersurbar_{\ubar}}\cdot\twonorms{\WangD^{\vec{m}_3}(\angnabla\FBetaA,\angD_{4}\Scalar)}{u}{1}{\sphere}{2}{\hypersurbar_{\ubar}}\right)\\
				\notag\lesssim&~\BAfluid \delta^{-1}\cdot\abs{u}^{\frac{1}{2}}\cdot\left(\BAsize\frac{\delta a^{\frac{1}{2}}}{\abs{u}}\right)^{i_0+k_0}.
			\end{align}
			For the last two terms on the RHS of \eqref{PV:alpha2}, using Sobolev embedding $H^2\hookrightarrow L^{\infty}$ for $\CAlphaAB$ on $\suubar$ and proved results in Proposition \ref{PR:EstimatesforRicciCurvature}, we obtain:
			\begin{align}\label{PV:alpha6}
				&\sum\limits_{\vec{m}_1+\vec{m}_2+\vec{m}_3+\vec{m}_4=(i-i_0-1,j,k-k_0)}\delta^{i_0+1}\twonorms{\abs{u}^{k_0}\WangD^{\vec{m}_1}\psi^{i_0}\cdot\WangD^{\vec{m}_2}\upsi^{k_0}\cdot\WangD^{\vec{m}_3}\Rhodual\cdot\WangD^{\vec{m}_4}\CAlphaAB}{u}{1}{\sphere}{2}{\hypersurbar_{\ubar}}\\
				\notag\lesssim&~\sum\limits_{\substack{\vec{m}_1+\vec{m}_2+\vec{m}_3+\vec{m}_4=(i-i_0-1,j,k-k_0)\\\abs{\vec{m}_1+\vec{m}_2}>N-3}}\delta^{i_0+1}\cdot\left(\twonorms{\abs{u}^{k_0}\WangD^{\vec{m}_1}\psi^{i_0}\cdot\WangD^{\vec{m}_2}\upsi^{k_0}}{u}{1}{\sphere}{2}{\hypersurbar_{\ubar}}\right.\\
				\notag&\left.\cdot\twonorms{\WangD^{\vec{m}_3}\Rhodual}{u}{\infty}{\sphere}{\infty}{\hypersurbar_{\ubar}}\cdot\twonorms{\WangD^{\vec{m}_4}\angnabla^2\CAlphaAB}{u}{\infty}{\sphere}{2}{\hypersurbar_{\ubar}}\right)\\
				\notag&+\sum\limits_{\substack{\vec{m}_1+\vec{m}_2+\vec{m}_3+\vec{m}_4=(i-i_0-1,j,k-k_0)\\\abs{\vec{m}_3}\geq N-3}}\delta^{i_0+1}\left(\twonorms{\abs{u}^{k_0}\WangD^{\vec{m}_1}\psi^{i_0}\cdot\WangD^{\vec{m}_2}\upsi^{k_0}}{u}{\infty}{\sphere}{\infty}{\hypersurbar_{\ubar}}\right.\\
				\notag&\left.\cdot\twonorms{\WangD^{\vec{m}_3}\Rhodual}{u}{1}{\sphere}{2}{\hypersurbar_{\ubar}}\cdot\twonorms{\WangD^{\vec{m}_4}\angnabla^2\CAlphaAB}{u}{\infty}{\sphere}{2}{\hypersurbar_{\ubar}}\right)\\		\notag&+\sum\limits_{\substack{\vec{m}_1+\vec{m}_2+\vec{m}_3+\vec{m}_4=(i-i_0-1,j,k-k_0)\\\abs{\vec{m}_4}\geq N-3}}\delta^{i_0+1}\left(\twonorms{\abs{u}^{k_0}\WangD^{\vec{m}_1}\psi^{i_0}\cdot\WangD^{\vec{m}_2}\upsi^{k_0}}{u}{\infty}{\sphere}{\infty}{\hypersurbar_{\ubar}}\right.\\
				\notag&\left.\cdot\twonorms{\WangD^{\vec{m}_3}\Rhodual}{u}{\infty}{\sphere}{\infty}{\hypersurbar_{\ubar}}\cdot\twonorms{\WangD^{\vec{m}_4}\CAlphaAB}{u}{1}{\sphere}{2}{\hypersurbar_{\ubar}}\right)\\
				\notag\lesssim&~\left(\BAsize\frac{\delta a^{\frac{1}{2}}}{\abs{u}}\right)^{i_0+k_0+2}\cdot\sum\limits_{\abs{\vec{m}}\leq N-2}\norm{\WangD^{\vec{m}}\CAlphaAB}_{L^2(\suubar)}.
			\end{align}
			\begin{align}\label{PV:alpha7}
				&\sum\limits_{\vec{m}_1+\vec{m}_2+\vec{m}_3+\vec{m}_4=(i-i_0-1,j,k-k_0-1)}\delta^{i_0}\twonorms{\abs{u}^{k_0+1}\WangD^{\vec{m}_1}\psi^{i_0}\cdot\WangD^{\vec{m}_2}\upsi^{k_0}\cdot\WangD^{\vec{m}_3}(\CBetabar,\FBetabar)\cdot\WangD^{\vec{m}_4}\CAlphaAB}{u}{1}{\sphere}{2}{\hypersurbar_{\ubar}}\\
				\notag\lesssim&~\sum\limits_{\substack{\vec{m}_1+\vec{m}_2+\vec{m}_3+\vec{m}_4=(i-i_0,j,k-k_0-1)\\\abs{\vec{m}_1+\vec{m}_2}>N-3}}\delta^{i_0}\cdot\left(\twonorms{\abs{u}^{k_0+1}\WangD^{\vec{m}_1}\psi^{i_0}\cdot\WangD^{\vec{m}_2}\upsi^{k_0}}{u}{1}{\sphere}{2}{\hypersurbar_{\ubar}}\right.\\
				\notag&\left.\cdot\twonorms{\WangD^{\vec{m}_3}(\CBetabar,\FBetabar)}{u}{\infty}{\sphere}{\infty}{\hypersurbar_{\ubar}}\cdot\twonorms{\WangD^{\vec{m}_4}\angnabla^2\CAlphaAB}{u}{\infty}{\sphere}{2}{\hypersurbar_{\ubar}}\right)\\
				\notag&+\sum\limits_{\substack{\vec{m}_1+\vec{m}_2+\vec{m}_3+\vec{m}_4=(i-i_0,j,k-k_0-1)\\\abs{\vec{m}_3}\geq N-3}}\delta^{i_0}\left(\twonorms{\abs{u}^{k_0+1}\WangD^{\vec{m}_1}\psi^{i_0}\cdot\WangD^{\vec{m}_2}\upsi^{k_0}}{u}{\infty}{\sphere}{\infty}{\hypersurbar_{\ubar}}\right.\\
				\notag&\left.\cdot\twonorms{\WangD^{\vec{m}_3}(\CBetabar,\FBetabar)}{u}{1}{\sphere}{2}{\hypersurbar_{\ubar}}\cdot\twonorms{\WangD^{\vec{m}_4}\angnabla^2\CAlphaAB}{u}{\infty}{\sphere}{2}{\hypersurbar_{\ubar}}\right)\\		\notag&+\sum\limits_{\substack{\vec{m}_1+\vec{m}_2+\vec{m}_3+\vec{m}_4=(i-i_0,j,k-k_0-1)\\\abs{\vec{m}_4}\geq N-3}}\delta^{i_0}\left(\twonorms{\abs{u}^{k_0+1}\WangD^{\vec{m}_1}\psi^{i_0}\cdot\WangD^{\vec{m}_2}\upsi^{k_0}}{u}{\infty}{\sphere}{\infty}{\hypersurbar_{\ubar}}\right.\\
				\notag&\left.\cdot\twonorms{\WangD^{\vec{m}_3}(\CBetabar,\FBetabar)}{u}{\infty}{\sphere}{\infty}{\hypersurbar_{\ubar}}\cdot\twonorms{\WangD^{\vec{m}_4}\CAlphaAB}{u}{1}{\sphere}{2}{\hypersurbar_{\ubar}}\right)\\
				\notag\lesssim&~\left(\BAsize\frac{\delta a^{\frac{1}{2}}}{\abs{u}}\right)^{i_0+k_0}\cdot\left(\BAsize\frac{\delta a^{\frac{1}{2}}}{\abs{u}}+  b^{-\frac{1}{2}}\right)\cdot\sum\limits_{\abs{\vec{m}}\leq N-2}\norm{\WangD^{\vec{m}}\CAlphaAB}_{L^2(\suubar)}.
			\end{align}
			The rest of the terms are identical to the proofs in Section \ref{SS:proveRicci}, e.g. \eqref{PV:hchi6}-\eqref{PV:hchi7} for the third term on the RHS of \eqref{PV:alpha2}. Combining \eqref{PV:alpha2}-\eqref{PV:alpha7}, soaking $\sum\limits_{\abs{\vec{m}}\leq N-2}\norm{\WangD^{\vec{m}}\CAlphaAB}_{L^2(\suubar)}$ terms on the RHS of \eqref{PV:alpha4.2}, \eqref{PV:alpha6}, \eqref{PV:alpha7} to the LHS, we obtain the desired result \eqref{ES:L2curvature1}.
		\end{proof}
		\begin{proof}[Proof of \eqref{ES:L2curvature2}]
			Let $\vec{m}=(i,j,k)$ with $i+j+k\leq N-2$. Let $\WangD^{\vec{m}}$ be as in Definition \ref{DE:operators}. We commute equations \eqref{EQ:D3CBeta} with $\frac{1}{\abs{u}^{k}}\WangD^{\vec{m}}$ using \eqref{EQ:ang3ijk}-\eqref{EQ:ang3ijk2}, we obtain
			\begin{align}\label{PV：beta2}
				&\frac{1}{\abs{u}^{k}}\angD_3\WangD^{\vec{m}}\CBetaA+\frac{1}{\abs{u}^{k}}\frac{k+2}{2}\gtr\uchi\cdot\WangD^{\vec{m}}\CBetaA\\
				\notag=&~\sum\limits_{\vec{m_1}+\vec{m_2}+\vec{m_3}=(i-i_0,j,k-k_0)}\frac{\delta^{i_0}}{\abs{u}^{k-k_0}}\WangD^{\vec{m_1}}\upsi^{k_0}\cdot\WangD^{\vec{m_2}}\psi^{i_0}\cdot\WangD^{\vec{m_3}}\angnabla(\CRho,\CRhodual)\\
				\notag&+\sum\limits_{\vec{m_1}+\vec{m_2}+\vec{m_3}=(i-i_0,j,k-k_0)}\frac{\delta^{i_0}}{\abs{u}^{k-k_0}}\WangD^{\vec{m_1}}\upsi^{k_0+1}\cdot\WangD^{\vec{m_2}}\psi^{i_0}\cdot\WangD^{\vec{m_3}}(\CBetaA,\FBetaA,\CRho,\CRhodual,\Scalar,\FRho,\FRhodual)\\
				\notag&+\sum\limits_{\vec{m_1}+\vec{m_2}+\vec{m_3}=(i-i_0,j,k-k_0)}\frac{\delta^{i_0}}{\abs{u}^{k-k_0}}\WangD^{\vec{m_1}}\upsi^{k_0}\cdot\WangD^{\vec{m_2}}\psi^{i_0+1}\cdot\WangD^{\vec{m_3}}(\FBetabar,\CBetabar)\\
				\notag&+\sum\limits_{\vec{m_1}+\vec{m_2}+\vec{m_3}=(i-i_0,j,k-k_0)}\frac{\delta^{i_0}}{\abs{u}^{k-k_0}}\WangD^{\vec{m_1}}\upsi^{k_0}\cdot\WangD^{\vec{m_2}}\psi^{i_0}\cdot\WangD^{\vec{m_3}}\angnabla(\FRho,\FRhodual,\Scalar)\\
				\notag&+\sum\limits_{\vec{m_1}+\vec{m_2}=(i-i_0,j,k-k_0+1)}\frac{\delta^{i_0}}{\abs{u}^{k-k_0+1}}\WangD^{\vec{m_1}}\upsi^{k_0+1}\cdot\WangD^{\vec{m_2}}\psi^{i_0+1}\\
				\notag&+\sum\limits_{\vec{m_1}+\vec{m_2}+\vec{m_3}+\vec{m_4}=(i-i_0-1,j,k-k_0+1)}\frac{\delta^{i_0+1}}{\abs{u}^{k-k_0+1}}\WangD^{\vec{m_1}}\upsi^{k_0}\cdot\WangD^{\vec{m_2}}\psi^{i_0}\cdot\WangD^{\vec{m_3}}\upzeta\cdot\WangD^{\vec{m_4}}\CBetaA\\
				\notag&+\sum\limits_{\vec{m_1}+\vec{m_2}+\vec{m_3}+\vec{m_4}=(i-i_0-1,j,k-k_0)}\frac{\delta^{i_0+1}}{\abs{u}^{k-k_0}}\WangD^{\vec{m_1}}\upsi^{k_0}\cdot\WangD^{\vec{m_2}}\psi^{i_0}\cdot\WangD^{\vec{m_3}}\Rhodual\cdot\WangD^{\vec{m_4}}\CBetaA\\
				\notag&+\sum\limits_{\vec{m_1}+\vec{m_2}+\vec{m_3}+\vec{m_4}=(i-i_0,j,k-k_0-1)}\frac{\delta^{i_0}}{\abs{u}^{k-k_0-1}}\WangD^{\vec{m_1}}\upsi^{k_0}\cdot\WangD^{\vec{m_2}}\psi^{i_0}\cdot\WangD^{\vec{m_3}}(\CBetabar,\FBetabar)\cdot\WangD^{\vec{m_4}}\CBetaA\\
				\notag&=:\frac{1}{\abs{u}^{k+1}}\mathfrak{F}_{(\CBetaA)}.
			\end{align}
			We apply Proposition \ref{PR:transport2} with $\lambda_0:=\frac{k+2}{2}$ to \eqref{PV：beta2}:
			\begin{align}\label{PV：beta3}
				\norm{\abs{u}\WangD^{\vec{m}}\CBetaA}_{L^2(\suubar)}\lesssim&~\norm{\abs{u}\WangD^{\vec{m}}\CBetaA}_{L^2(S_{\uinf,\ubar})}+\int_{\uinf}^{u}\norm{\mathfrak{F}_{(\CBetaA)}}_{L^2(\suprimeubar)}\diff\uprime.
			\end{align}
				We now estimate $\norm{\mathfrak{F}_{(\CBetaA)}}_{L^1_uL^2_{\sphere}(\hypersurbar_{\ubar})}$. For the first term on the RHS of \eqref{PV：beta2}, by bootstrap assumptions \eqref{BA}, we have:
				\begin{align}\label{PV：beta4}
					&\sum\limits_{\vec{m}_1+\vec{m}_2+\vec{m}_3=(i-i_0,j,k-k_0)}\delta^{i_0}\twonorms{\abs{u}^{k_0+1}\WangD^{\vec{m}_1}\psi^{i_0}\cdot\WangD^{\vec{m}_2}\upsi^{k_0}\cdot\WangD^{\vec{m}_3}\angnabla(\CRho,\CRhodual)}{u}{1}{\sphere}{2}{\hypersurbar_{\ubar}}\\
					\notag\lesssim&~\sum\limits_{\substack{\vec{m}_1+\vec{m}_2+\vec{m}_3=(i-i_0,j,k-k_0)\\\abs{\vec{m}_1+\vec{m}_2}>N-3}}\delta^{i_0}\cdot\left(\twonorms{\abs{u}^{k_0}\WangD^{\vec{m}_1}\psi^{i_0}\cdot\WangD^{\vec{m}_2}\upsi^{k_0}}{u}{1}{\sphere}{2}{\hypersurbar_{\ubar}}\cdot\twonorms{\abs{u}\WangD^{\vec{m}_3}\angnabla(\CRho,\CRhodual)}{u}{\infty}{\sphere}{\infty}{\hypersurbar_{\ubar}}\right)\\
					\notag&+\sum\limits_{\substack{\vec{m}_1+\vec{m}_2+\vec{m}_3=(i-i_0,j,k-k_0)\\\abs{\vec{m}_1+\vec{m}_2}\leq N-3}}\delta^{i_0}\left(\twonorms{\abs{u}^{k_0}\WangD^{\vec{m}_1}\psi^{i_0}\cdot\WangD^{\vec{m}_2}\upsi^{k_0}}{u}{\infty}{\sphere}{\infty}{\hypersurbar_{\ubar}}\cdot\twonorms{\abs{u}^2\WangD^{\vec{m}_3}\angnabla(\CRho,\CRhodual)}{u}{2}{\sphere}{2}{\hypersurbar_{\ubar}}\cdot\norm{\abs{u}^{-1}}_{L^2_u}\right)\\
					\notag\lesssim&~\left(\BAsize\frac{\delta a^{\frac{1}{2}}}{\abs{u}}\right)^{i_0+k_0}\cdot\frac{\delta^{\frac{1}{2}}a^{\frac{1}{2}}}{\abs{u}^{\frac{1}{2}}}(1+\tilde{\mathcal{R}}_{N-1,2}).
				\end{align}
				For the second term and the fifth term on the RHS of \eqref{PV：beta2}, by Sobolev embedding $H^2\hookrightarrow L^{\infty}$ for $\CBetaA$ on $\suubar$, bootstrap assumptions \eqref{BA} and proved results in Proposition \ref{PR:EstimatesforRicciCurvature}, Proposition \ref{PR:L2RicciCoefficients} and Corollary \ref{ES:LinftyRicciCoefficients} we have:
				\begin{align}\label{PV：beta4.2}
					&\sum\limits_{\vec{m}_1+\vec{m}_2+\vec{m}_3=(i-i_0,j,k-k_0)}\delta^{i_0}\twonorms{\abs{u}^{k_0+1}\WangD^{\vec{m}_1}\psi^{i_0}\cdot\WangD^{\vec{m}_2}\upsi^{k_0+1}\cdot\WangD^{\vec{m}_3}(\CBetaA,\FBetaA,\CRho,\CRhodual,\Scalar,\FRho,\FRhodual,\angnabla\psi)}{u}{1}{\sphere}{2}{\hypersurbar_{\ubar}}\\
					\notag\lesssim&~\sum\limits_{\substack{\vec{m}_1+\vec{m}_2+\vec{m}_3=(i-i_0,j,k-k_0)\\\abs{\vec{m}_1+\vec{m}_2}>N-3}}\delta^{i_0}\cdot\left[\twonorms{\abs{u}^{k_0+1}\WangD^{\vec{m}_1}\psi^{i_0}\cdot\WangD^{\vec{m}_2}\upsi^{k_0+1}}{u}{1}{\sphere}{2}{\hypersurbar_{\ubar}}\right.\\
					\notag&\left.\cdot\left(\twonorms{\WangD^{\vec{m}_3}(\FBetaA,\CRho,\CRhodual,\Scalar,\FRho,\FRhodual,\angnabla\psi)}{u}{\infty}{\sphere}{\infty}{\hypersurbar_{\ubar}}+\twonorms{\WangD^{\vec{m}_4}\angnabla^2\CBetaA}{u}{\infty}{\sphere}{2}{\hypersurbar_{\ubar}}\right)\right]\\
					\notag&+\sum\limits_{\substack{\vec{m}_1+\vec{m}_2+\vec{m}_3=(i-i_0,j,k-k_0)\\\abs{\vec{m}_1+\vec{m}_2}\leq N-3}}\delta^{i_0}\left(\twonorms{\abs{u}^{k_0+1}\WangD^{\vec{m}_1}\psi^{i_0}\cdot\WangD^{\vec{m}_2}\upsi^{k_0+1}}{u}{\infty}{\sphere}{\infty}{\hypersurbar_{\ubar}}\cdot\twonorms{\WangD^{\vec{m}_3}(\CBetaA,\FBetaA,\CRho,\CRhodual,\Scalar,\FRho,\FRhodual,\angnabla\psi)}{u}{1}{\sphere}{2}{\hypersurbar_{\ubar}}\right)\\
					\notag\lesssim&~\left(\BAsize\frac{\delta a^{\frac{1}{2}}}{\abs{u}}\right)^{i_0+k_0+1}\cdot\left\{\sum\limits_{\abs{\vec{m}}\leq N-2}\norm{\abs{u}\WangD^{\vec{m}}\CBetaA}_{L^2(\suubar)}+a^{\frac{1}{2}}+   b^{\frac{1}{2}}\right\}.
				\end{align}
				For the third term and the fifth term on the RHS of \eqref{PV：beta2},  by bootstrap assumptions \eqref{BA} and proved results in Proposition \ref{PR:EstimatesforRicciCurvature}, Proposition \ref{PR:L2RicciCoefficients} and Corollary \ref{ES:LinftyRicciCoefficients} we have:
				\begin{align}\label{PV：beta4.3}
					&\sum\limits_{\vec{m}_1+\vec{m}_2+\vec{m}_3=(i-i_0,j,k-k_0)}\delta^{i_0}\twonorms{\abs{u}^{k_0+1}\WangD^{\vec{m}_1}\psi^{i_0+1}\cdot\WangD^{\vec{m}_2}\upsi^{k_0}\cdot\WangD^{\vec{m}_3}(\CBetabar,\FBetabar)}{u}{1}{\sphere}{2}{\hypersurbar_{\ubar}}\\
					\notag\lesssim&~\sum\limits_{\substack{\vec{m}_1+\vec{m}_2+\vec{m}_3=(i-i_0,j,k-k_0)\\\abs{\vec{m}_1+\vec{m}_2}>N-3}}\delta^{i_0}\cdot\left(\twonorms{\abs{u}^{k_0+1}\WangD^{\vec{m}_1}\psi^{i_0+1}\cdot\WangD^{\vec{m}_2}\upsi^{k_0}}{u}{1}{\sphere}{2}{\hypersurbar_{\ubar}}\cdot\twonorms{\WangD^{\vec{m}_3}(\CBetabar,\FBetabar)}{u}{\infty}{\sphere}{\infty}{\hypersurbar_{\ubar}}\right)\\
					\notag&+\sum\limits_{\substack{\vec{m}_1+\vec{m}_2+\vec{m}_3=(i-i_0,j,k-k_0)\\\abs{\vec{m}_1+\vec{m}_2}\leq N-3}}\delta^{i_0}\left(\twonorms{\abs{u}^{k_0+1}\WangD^{\vec{m}_1}\psi^{i_0+1}\cdot\WangD^{\vec{m}_2}\upsi^{k_0}}{u}{\infty}{\sphere}{\infty}{\hypersurbar_{\ubar}}\cdot\twonorms{\WangD^{\vec{m}_3}(\CBetabar,\FBetabar)}{u}{1}{\sphere}{2}{\hypersurbar_{\ubar}}\right)\\
					\notag\lesssim&~\left(\BAsize\frac{\delta a^{\frac{1}{2}}}{\abs{u}}\right)^{i_0+k_0}\cdot a^{\frac{1}{2}}\cdot\left(\frac{\delta a^{\frac{1}{2}}}{\abs{u}}+   b^{-\frac{1}{2}}\cdot\abs{u}^2\right).
				\end{align}
				For the fourth term on the RHS of \eqref{PV：beta2}, by bootstrap assumptions \eqref{BA} and proved results in Proposition \ref{PR:EstimatesforRicciCurvature}, we have:
				\begin{align}\label{PV：beta5}
					&\sum\limits_{\vec{m}_1+\vec{m}_2+\vec{m}_3=(i-i_0,j,k-k_0)}\delta^{i_0}\twonorms{\abs{u}^{k_0+1}\WangD^{\vec{m}_1}\psi^{i_0}\cdot\WangD^{\vec{m}_2}\upsi^{k_0}\cdot\WangD^{\vec{m}_3}\angnabla(\FRho,\FRhodual,\Scalar)}{u}{1}{\sphere}{2}{\hypersurbar_{\ubar}}\\
					\notag\lesssim&~\sum\limits_{\substack{\vec{m}_1+\vec{m}_2+\vec{m}_3=(i-i_0,j,k-k_0)\\\abs{\vec{m}_1+\vec{m}_2}>N-3}}\delta^{i_0}\cdot\left(\twonorms{\abs{u}^{k_0+1}\WangD^{\vec{m}_1}\psi^{i_0}\cdot\WangD^{\vec{m}_2}\upsi^{k_0}}{u}{1}{\sphere}{2}{\hypersurbar_{\ubar}}\cdot\twonorms{\WangD^{\vec{m}_3}\angnabla(\FRho,\FRhodual,\Scalar)}{u}{\infty}{\sphere}{\infty}{\hypersurbar_{\ubar}}\right)\\
					\notag&+\sum\limits_{\substack{\vec{m}_1+\vec{m}_2+\vec{m}_3=(i-i_0,j,k-k_0)\\\abs{\vec{m}_1+\vec{m}_2}\leq N-3}}\delta^{i_0}\left(\twonorms{\abs{u}^{k_0+1}\WangD^{\vec{m}_1}\psi^{i_0}\cdot\WangD^{\vec{m}_2}\upsi^{k_0}}{u}{\infty}{\sphere}{\infty}{\hypersurbar_{\ubar}}\cdot\twonorms{\WangD^{\vec{m}_3}\angnabla(\FRho,\FRhodual,\Scalar)}{u}{1}{\sphere}{2}{\hypersurbar_{\ubar}}\right)\\
					\notag\lesssim&~ \abs{u}^{2}\cdot\left(\BAsize\frac{\delta a^{\frac{1}{2}}}{\abs{u}}\right)^{i_0+k_0}.
				\end{align}
				For the last two terms on the RHS of \eqref{PV：beta2}, using Sobolev embedding $H^2\hookrightarrow L^{\infty}$ for $\CBetaA$ on $\suubar$ and proved results in Proposition \ref{PR:EstimatesforRicciCurvature}, we obtain:
				\begin{align}\label{PV：beta6}
					&\sum\limits_{\vec{m}_1+\vec{m}_2+\vec{m}_3+\vec{m}_4=(i-i_0-1,j,k-k_0)}\delta^{i_0+1}\twonorms{\abs{u}^{k_0+1}\WangD^{\vec{m}_1}\psi^{i_0}\cdot\WangD^{\vec{m}_2}\upsi^{k_0}\cdot\WangD^{\vec{m}_3}\Rhodual\cdot\WangD^{\vec{m}_4}\CBetaA}{u}{1}{\sphere}{2}{\hypersurbar_{\ubar}}\\
					\notag\lesssim&~\sum\limits_{\substack{\vec{m}_1+\vec{m}_2+\vec{m}_3+\vec{m}_4=(i-i_0-1,j,k-k_0)\\\abs{\vec{m}_1+\vec{m}_2}>N-3}}\delta^{i_0+1}\cdot\left(\twonorms{\abs{u}^{k_0+1}\WangD^{\vec{m}_1}\psi^{i_0}\cdot\WangD^{\vec{m}_2}\upsi^{k_0}}{u}{1}{\sphere}{2}{\hypersurbar_{\ubar}}\right.\\
					\notag&\left.\cdot\twonorms{\WangD^{\vec{m}_3}\Rhodual}{u}{\infty}{\sphere}{\infty}{\hypersurbar_{\ubar}}\cdot\twonorms{\WangD^{\vec{m}_4}\angnabla^2\CBetaA}{u}{\infty}{\sphere}{2}{\hypersurbar_{\ubar}}\right)\\
					\notag&+\sum\limits_{\substack{\vec{m}_1+\vec{m}_2+\vec{m}_3+\vec{m}_4=(i-i_0-1,j,k-k_0)\\\abs{\vec{m}_3}\geq N-3}}\delta^{i_0+1}\left(\twonorms{\abs{u}^{k_0+1}\WangD^{\vec{m}_1}\psi^{i_0}\cdot\WangD^{\vec{m}_2}\upsi^{k_0}}{u}{\infty}{\sphere}{\infty}{\hypersurbar_{\ubar}}\right.\\
					\notag&\left.\cdot\twonorms{\WangD^{\vec{m}_3}\Rhodual}{u}{1}{\sphere}{2}{\hypersurbar_{\ubar}}\cdot\twonorms{\WangD^{\vec{m}_4}\angnabla^2\CBetaA}{u}{\infty}{\sphere}{2}{\hypersurbar_{\ubar}}\right)\\		\notag&+\sum\limits_{\substack{\vec{m}_1+\vec{m}_2+\vec{m}_3+\vec{m}_4=(i-i_0-1,j,k-k_0)\\\abs{\vec{m}_4}\geq N-3}}\delta^{i_0+1}\left(\twonorms{\abs{u}^{k_0+1}\WangD^{\vec{m}_1}\psi^{i_0}\cdot\WangD^{\vec{m}_2}\upsi^{k_0}}{u}{\infty}{\sphere}{\infty}{\hypersurbar_{\ubar}}\right.\\
					\notag&\left.\cdot\twonorms{\WangD^{\vec{m}_3}\Rhodual}{u}{\infty}{\sphere}{\infty}{\hypersurbar_{\ubar}}\cdot\twonorms{\WangD^{\vec{m}_4}\CBetaA}{u}{1}{\sphere}{2}{\hypersurbar_{\ubar}}\right)\\
					\notag\lesssim&~\left(\BAsize\frac{\delta a^{\frac{1}{2}}}{\abs{u}}\right)^{i_0+k_0+2}\cdot\sum\limits_{\abs{\vec{m}}\leq N-2}\norm{\abs{u}\WangD^{\vec{m}}\CBetaA}_{L^2(\suubar)}.
				\end{align}
				\begin{align}\label{PV：beta7}
					&\sum\limits_{\vec{m}_1+\vec{m}_2+\vec{m}_3+\vec{m}_4=(i-i_0-1,j,k-k_0-1)}\delta^{i_0}\twonorms{\abs{u}^{k_0+2}\WangD^{\vec{m}_1}\psi^{i_0}\cdot\WangD^{\vec{m}_2}\upsi^{k_0}\cdot\WangD^{\vec{m}_3}(\CBetabar,\FBetabar)\cdot\WangD^{\vec{m}_4}\CBetaA}{u}{1}{\sphere}{2}{\hypersurbar_{\ubar}}\\
					\notag\lesssim&~\sum\limits_{\substack{\vec{m}_1+\vec{m}_2+\vec{m}_3+\vec{m}_4=(i-i_0,j,k-k_0-1)\\\abs{\vec{m}_1+\vec{m}_2}>N-3}}\delta^{i_0}\cdot\left(\twonorms{\abs{u}^{k_0+2}\WangD^{\vec{m}_1}\psi^{i_0}\cdot\WangD^{\vec{m}_2}\upsi^{k_0}}{u}{1}{\sphere}{2}{\hypersurbar_{\ubar}}\right.\\
					\notag&\left.\cdot\twonorms{\WangD^{\vec{m}_3}(\CBetabar,\FBetabar)}{u}{\infty}{\sphere}{\infty}{\hypersurbar_{\ubar}}\cdot\twonorms{\WangD^{\vec{m}_4}\angnabla^2\CBetaA}{u}{\infty}{\sphere}{2}{\hypersurbar_{\ubar}}\right)\\
					\notag&+\sum\limits_{\substack{\vec{m}_1+\vec{m}_2+\vec{m}_3+\vec{m}_4=(i-i_0,j,k-k_0-1)\\\abs{\vec{m}_3}\geq N-3}}\delta^{i_0}\left(\twonorms{\abs{u}^{k_0+2}\WangD^{\vec{m}_1}\psi^{i_0}\cdot\WangD^{\vec{m}_2}\upsi^{k_0}}{u}{\infty}{\sphere}{\infty}{\hypersurbar_{\ubar}}\right.\\
					\notag&\left.\cdot\twonorms{\WangD^{\vec{m}_3}(\CBetabar,\FBetabar)}{u}{1}{\sphere}{2}{\hypersurbar_{\ubar}}\cdot\twonorms{\WangD^{\vec{m}_4}\angnabla^2\CBetaA}{u}{\infty}{\sphere}{2}{\hypersurbar_{\ubar}}\right)\\		\notag&+\sum\limits_{\substack{\vec{m}_1+\vec{m}_2+\vec{m}_3+\vec{m}_4=(i-i_0,j,k-k_0-1)\\\abs{\vec{m}_4}\geq N-3}}\delta^{i_0}\left(\twonorms{\abs{u}^{k_0+2}\WangD^{\vec{m}_1}\psi^{i_0}\cdot\WangD^{\vec{m}_2}\upsi^{k_0}}{u}{\infty}{\sphere}{\infty}{\hypersurbar_{\ubar}}\right.\\
					\notag&\left.\cdot\twonorms{\WangD^{\vec{m}_3}(\CBetabar,\FBetabar)}{u}{\infty}{\sphere}{\infty}{\hypersurbar_{\ubar}}\cdot\twonorms{\WangD^{\vec{m}_4}\CBetaA}{u}{1}{\sphere}{2}{\hypersurbar_{\ubar}}\right)\\
					\notag\lesssim&~\left(\BAsize\frac{\delta a^{\frac{1}{2}}}{\abs{u}}\right)^{i_0+k_0}\cdot\left(\BAsize\frac{\delta a^{\frac{1}{2}}}{\abs{u}}+  b^{-\frac{1}{2}}\cdot\abs{u}\right)\cdot\sum\limits_{\abs{\vec{m}}\leq N-2}\norm{\abs{u}\WangD^{\vec{m}}\CBetaA}_{L^2(\suubar)}.
				\end{align}
				The rest of the terms are identical to the proofs in Section \ref{SS:proveRicci}. Combining \eqref{PV：beta2}-\eqref{PV：beta7}, soaking $\sum\limits_{\abs{\vec{m}}\leq N-2}\norm{\abs{u}\WangD^{\vec{m}}\CBetaA}_{L^2(\suubar)}$ terms on the RHS of \eqref{PV：beta4.2}, \eqref{PV：beta6}, \eqref{PV：beta7} to the LHS, we obtain the desired result \eqref{ES:L2curvature2}.
			\end{proof}
			\begin{proof}[Proof of \eqref{ES:L2curvature3}]
				Let $\vec{m}=(i,j,k)$ with $i+j+k\leq N-2$. Let $\WangD^{\vec{m}}$ be as in Definition \ref{DE:operators}. We commute equations \eqref{EQ:D3Rho2}, \eqref{EQ:D3Rhodual2} with $\frac{1}{\abs{u}^{k}}\WangD^{\vec{m}}$ and commute \eqref{EQ:D3CBetabar} with $\frac{1}{\abs{u}^{k+1}}\WangD^{\vec{m}}$,using \eqref{EQ:ang3ijk}-\eqref{EQ:ang3ijk2}, we obtain
				\begin{align}\label{PV：rhorhodualbetabar2}
					&\frac{1}{\abs{u}^{k}}\angD_3\WangD^{\vec{m}}(\CRho,\CRhodual)+\frac{1}{\abs{u}^{k}}\frac{k+3}{2}\gtr\uchi\cdot\WangD^{\vec{m}}(\CRho,\CRhodual)\\
					\notag=&~\sum\limits_{\vec{m_1}+\vec{m_2}+\vec{m_3}=(i-i_0,j,k-k_0)}\frac{\delta^{i_0}}{\abs{u}^{k-k_0}}\WangD^{\vec{m_1}}\upsi^{k_0}\cdot\WangD^{\vec{m_2}}\psi^{i_0}\cdot\WangD^{\vec{m_3}}\angnabla\CBetabar\\
					\notag&+\sum\limits_{\vec{m_1}+\vec{m_2}+\vec{m_3}=(i-i_0,j,k-k_0)}\frac{\delta^{i_0}}{\abs{u}^{k-k_0}}\WangD^{\vec{m_1}}\upsi^{k_0+1}\cdot\WangD^{\vec{m_2}}\psi^{i_0}\cdot\WangD^{\vec{m_3}}(\FBetabar,\Scalar,\FRho)\\
					\notag&+\sum\limits_{\vec{m_1}+\vec{m_2}+\vec{m_3}=(i-i_0,j,k-k_0)}\frac{\delta^{i_0}}{\abs{u}^{k-k_0}}\WangD^{\vec{m_1}}\upsi^{k_0}\cdot\WangD^{\vec{m_2}}\psi^{i_0}\cdot\WangD^{\vec{m_3}}\angnabla\FBetabar\\
					\notag&+\sum\limits_{\vec{m_1}+\vec{m_2}+\vec{m_3}=(i-i_0,j,k-k_0)}\frac{\delta^{i_0}}{\abs{u}^{k-k_0}}\WangD^{\vec{m_1}}\upsi^{k_0}\cdot\WangD^{\vec{m_2}}\psi^{i_0+1}\cdot\WangD^{\vec{m_3}}\angnabla\upsi\\
					\notag&+\sum\limits_{\vec{m_1}+\vec{m_2}+\vec{m_3}=(i-i_0,j,k-k_0+1)}\frac{\delta^{i_0}}{\abs{u}^{k-k_0+1}}\WangD^{\vec{m_1}}\upsi^{k_0+2}\cdot\WangD^{\vec{m_2}}\psi^{i_0}\cdot\WangD^{\vec{m_3}}(\psi,\upsi)\\
					\notag&+\sum\limits_{\vec{m_1}+\vec{m_2}+\vec{m_3}+\vec{m_4}=(i-i_0-1,j,k-k_0+1)}\frac{\delta^{i_0+1}}{\abs{u}^{k-k_0+1}}\WangD^{\vec{m_1}}\upsi^{k_0}\cdot\WangD^{\vec{m_2}}\psi^{i_0}\cdot\WangD^{\vec{m_3}}\upzeta\cdot\WangD^{\vec{m_4}}(\CRho,\CRhodual)\\
					\notag&+\sum\limits_{\vec{m_1}+\vec{m_2}+\vec{m_3}+\vec{m_4}=(i-i_0-1,j,k-k_0)}\frac{\delta^{i_0+1}}{\abs{u}^{k-k_0}}\WangD^{\vec{m_1}}\upsi^{k_0}\cdot\WangD^{\vec{m_2}}\psi^{i_0}\cdot\WangD^{\vec{m_3}}\Rhodual\cdot\WangD^{\vec{m_4}}(\CRho,\CRhodual)\\
					\notag&+\sum\limits_{\vec{m_1}+\vec{m_2}+\vec{m_3}+\vec{m_4}=(i-i_0,j,k-k_0-1)}\frac{\delta^{i_0}}{\abs{u}^{k-k_0-1}}\WangD^{\vec{m_1}}\upsi^{k_0}\cdot\WangD^{\vec{m_2}}\psi^{i_0}\cdot\WangD^{\vec{m_3}}(\CBetabar,\FBetabar)\cdot\WangD^{\vec{m_4}}(\CRho,\CRhodual)\\
					\notag&=:\frac{1}{\abs{u}^{k+2}}\mathfrak{F}_{((\CRho,\CRhodual))}.
				\end{align}
				\begin{align}\label{PV：rhorhodualbetabar2.1}
					&\frac{1}{\abs{u}^{k+1}}\angD_3\WangD^{\vec{m}}\CBetabar+\frac{1}{\abs{u}^{k+1}}\frac{k+4}{2}\gtr\uchi\cdot\WangD^{\vec{m}}\CBetabar\\
					\notag=&~\sum\limits_{\vec{m_1}+\vec{m_2}+\vec{m_3}=(i-i_0,j,k-k_0)}\frac{\delta^{i_0}}{\abs{u}^{k-k_0+1}}\WangD^{\vec{m_1}}\upsi^{k_0}\cdot\WangD^{\vec{m_2}}\psi^{i_0}\cdot\WangD^{\vec{m_3}}\angnabla\CAlphabar\\
					\notag&+\sum\limits_{\vec{m_1}+\vec{m_2}+\vec{m_3}=(i-i_0,j,k-k_0)}\frac{\delta^{i_0}}{\abs{u}^{k-k_0+1}}\WangD^{\vec{m_1}}\upsi^{k_0+1}\cdot\WangD^{\vec{m_2}}\psi^{i_0}\cdot\WangD^{\vec{m_3}}(\CAlphabar,\FAlphabar,\FBetabar,\CBetabar)\\
					\notag&+\sum\limits_{\vec{m_1}+\vec{m_2}+\vec{m_3}+\vec{m_4}=(i-i_0-1,j,k-k_0+1)}\frac{\delta^{i_0+1}}{\abs{u}^{k-k_0+2}}\WangD^{\vec{m_1}}\upsi^{k_0}\cdot\WangD^{\vec{m_2}}\psi^{i_0}\cdot\WangD^{\vec{m_3}}\upzeta\cdot\WangD^{\vec{m_4}}\CBetabar\\
					\notag&+\sum\limits_{\vec{m_1}+\vec{m_2}+\vec{m_3}+\vec{m_4}=(i-i_0-1,j,k-k_0)}\frac{\delta^{i_0+1}}{\abs{u}^{k-k_0+1}}\WangD^{\vec{m_1}}\upsi^{k_0}\cdot\WangD^{\vec{m_2}}\psi^{i_0}\cdot\WangD^{\vec{m_3}}\Rhodual\cdot\WangD^{\vec{m_4}}\CBetabar\\
					\notag&+\sum\limits_{\vec{m_1}+\vec{m_2}+\vec{m_3}+\vec{m_4}=(i-i_0,j,k-k_0-1)}\frac{\delta^{i_0}}{\abs{u}^{k-k_0}}\WangD^{\vec{m_1}}\upsi^{k_0}\cdot\WangD^{\vec{m_2}}\psi^{i_0}\cdot\WangD^{\vec{m_3}}(\CBetabar,\FBetabar)\cdot\WangD^{\vec{m_4}}\CBetabar\\
					\notag&=:\frac{1}{\abs{u}^{k+3}}\mathfrak{F}_{(\CBetabar)}.
				\end{align}
				We apply Proposition \ref{PR:transport2} with $\lambda_0:=\frac{k+3}{2}$,  to \eqref{PV：rhorhodualbetabar2}, and with $\lambda_0:=\frac{k+4}{2}$ to \eqref{PV：rhorhodualbetabar2.1}:
				\begin{align}\label{PV：rhorhodualbetabar3}
					\norm{\abs{u}^2\WangD^{\vec{m}}(\CRho,\CRhodual,\CBetabar)}_{L^2(\suubar)}\lesssim&~\norm{\abs{u}^2\WangD^{\vec{m}}(\CRho,\CRhodual,\CBetabar)}_{L^2(S_{\uinf,\ubar})}+\int_{\uinf}^{u}\norm{\mathfrak{F}_{(\CRho,\CRhodual)},\mathfrak{F}_{(\CBetabar)}}_{L^2(\suprimeubar)}\diff\uprime.
				\end{align}
					We now estimate $\norm{\mathfrak{F}_{(\CRho,\CRhodual)},\mathfrak{F}_{(\CBetabar)}}_{L^1_uL^2_{\sphere}(\hypersurbar_{\ubar})}$. For the first term on the RHS of \eqref{PV：rhorhodualbetabar2} and the first term on the RHS of \eqref{PV：rhorhodualbetabar2.1}, by bootstrap assumptions \eqref{BA}, we have:
					\begin{align}\label{PV：rhorhodualbetabar4}
						&\sum\limits_{\vec{m}_1+\vec{m}_2+\vec{m}_3=(i-i_0,j,k-k_0)}\delta^{i_0}\twonorms{\abs{u}^{k_0+2}\WangD^{\vec{m}_1}\psi^{i_0}\cdot\WangD^{\vec{m}_2}\upsi^{k_0}\cdot\WangD^{\vec{m}_3}\angnabla(\CAlphabar,\CBetabar)}{u}{1}{\sphere}{2}{\hypersurbar_{\ubar}}\\
						\notag\lesssim&~\sum\limits_{\substack{\vec{m}_1+\vec{m}_2+\vec{m}_3=(i-i_0,j,k-k_0)\\\abs{\vec{m}_1+\vec{m}_2}>N-3}}\delta^{i_0}\cdot\left(\twonorms{\abs{u}^{k_0+2}\WangD^{\vec{m}_1}\psi^{i_0}\cdot\WangD^{\vec{m}_2}\upsi^{k_0}}{u}{1}{\sphere}{2}{\hypersurbar_{\ubar}}\cdot\twonorms{\WangD^{\vec{m}_3}\angnabla(\CAlphabar,\CBetabar)}{u}{\infty}{\sphere}{\infty}{\hypersurbar_{\ubar}}\right)\\
						\notag&+\sum\limits_{\substack{\vec{m}_1+\vec{m}_2+\vec{m}_3=(i-i_0,j,k-k_0)\\\abs{\vec{m}_1+\vec{m}_2}\leq N-3}}\delta^{i_0}\left(\twonorms{\abs{u}^{k_0}\WangD^{\vec{m}_1}\psi^{i_0}\cdot\WangD^{\vec{m}_2}\upsi^{k_0}}{u}{\infty}{\sphere}{\infty}{\hypersurbar_{\ubar}}\cdot\twonorms{\abs{u}^3\WangD^{\vec{m}_3}\angnabla(\CAlphabar,\CBetabar)}{u}{2}{\sphere}{2}{\hypersurbar_{\ubar}}\cdot\norm{\abs{u}^{-1}}_{L^2_u}\right)\\
						\notag\lesssim&~\left(\BAsize\frac{\delta a^{\frac{1}{2}}}{\abs{u}}\right)^{i_0+k_0}\cdot\frac{\delta^{\frac{3}{2}} a}{\abs{u}^{\frac{1}{2}}}\left(\frac{\delta^{-\frac{1}{2}}a^{-\frac{1}{2}}}{\abs{u}^{\frac{3}{2}}}+\tilde{\mathcal{R}}_{N-1,2}\right).
					\end{align}
					For the last two terms on the RHS of \eqref{PV：rhorhodualbetabar2} and \eqref{PV：rhorhodualbetabar2.1}, using Sobolev embedding $H^2\hookrightarrow L^{\infty}$ for $\CRho,\CRhodual,\CBetabar$ on $\suubar$ and proved results in Proposition \ref{PR:EstimatesforRicciCurvature}, we obtain:
					\begin{align}\label{PV：rhorhodualbetabar6}
						&\sum\limits_{\vec{m}_1+\vec{m}_2+\vec{m}_3+\vec{m}_4=(i-i_0-1,j,k-k_0)}\delta^{i_0+1}\twonorms{\abs{u}^{k_0+2}\WangD^{\vec{m}_1}\psi^{i_0}\cdot\WangD^{\vec{m}_2}\upsi^{k_0}\cdot\WangD^{\vec{m}_3}\Rhodual\cdot\WangD^{\vec{m}_4}(\CRho,\CRhodual,\CBetabar)}{u}{1}{\sphere}{2}{\hypersurbar_{\ubar}}\\
						\notag\lesssim&~\sum\limits_{\substack{\vec{m}_1+\vec{m}_2+\vec{m}_3+\vec{m}_4=(i-i_0-1,j,k-k_0)\\\abs{\vec{m}_1+\vec{m}_2}>N-3}}\delta^{i_0+1}\cdot\left(\twonorms{\abs{u}^{k_0+2}\WangD^{\vec{m}_1}\psi^{i_0}\cdot\WangD^{\vec{m}_2}\upsi^{k_0}}{u}{1}{\sphere}{2}{\hypersurbar_{\ubar}}\right.\\
						\notag&\left.\cdot\twonorms{\WangD^{\vec{m}_3}\Rhodual}{u}{\infty}{\sphere}{\infty}{\hypersurbar_{\ubar}}\cdot\twonorms{\WangD^{\vec{m}_4}\angnabla^2(\CRho,\CRhodual,\CBetabar)}{u}{\infty}{\sphere}{2}{\hypersurbar_{\ubar}}\right)\\
						\notag&+\sum\limits_{\substack{\vec{m}_1+\vec{m}_2+\vec{m}_3+\vec{m}_4=(i-i_0-1,j,k-k_0)\\\abs{\vec{m}_3}\geq N-3}}\delta^{i_0+1}\left(\twonorms{\abs{u}^{k_0+2}\WangD^{\vec{m}_1}\psi^{i_0}\cdot\WangD^{\vec{m}_2}\upsi^{k_0}}{u}{\infty}{\sphere}{\infty}{\hypersurbar_{\ubar}}\right.\\
						\notag&\left.\cdot\twonorms{\WangD^{\vec{m}_3}\Rhodual}{u}{1}{\sphere}{2}{\hypersurbar_{\ubar}}\cdot\twonorms{\WangD^{\vec{m}_4}\angnabla^2(\CRho,\CRhodual,\CBetabar)}{u}{\infty}{\sphere}{2}{\hypersurbar_{\ubar}}\right)\\		\notag&+\sum\limits_{\substack{\vec{m}_1+\vec{m}_2+\vec{m}_3+\vec{m}_4=(i-i_0-1,j,k-k_0)\\\abs{\vec{m}_4}\geq N-3}}\delta^{i_0+1}\left(\twonorms{\abs{u}^{k_0+2}\WangD^{\vec{m}_1}\psi^{i_0}\cdot\WangD^{\vec{m}_2}\upsi^{k_0}}{u}{\infty}{\sphere}{\infty}{\hypersurbar_{\ubar}}\right.\\
						\notag&\left.\cdot\twonorms{\WangD^{\vec{m}_3}\Rhodual}{u}{\infty}{\sphere}{\infty}{\hypersurbar_{\ubar}}\cdot\twonorms{\WangD^{\vec{m}_4}(\CRho,\CRhodual,\CBetabar)}{u}{1}{\sphere}{2}{\hypersurbar_{\ubar}}\right)\\
						\notag\lesssim&~\left(\BAsize\frac{\delta a^{\frac{1}{2}}}{\abs{u}}\right)^{i_0+k_0+2}\cdot\sum\limits_{\abs{\vec{m}}\leq N-2}\norm{\abs{u}^2\WangD^{\vec{m}}(\CRho,\CRhodual,\CBetabar)}_{L^2(\suubar)}.
					\end{align}
					\begin{align}\label{PV：rhorhodualbetabar7}
						&\sum\limits_{\vec{m}_1+\vec{m}_2+\vec{m}_3+\vec{m}_4=(i-i_0-1,j,k-k_0-1)}\delta^{i_0}\twonorms{\abs{u}^{k_0+3}\WangD^{\vec{m}_1}\psi^{i_0}\cdot\WangD^{\vec{m}_2}\upsi^{k_0}\cdot\WangD^{\vec{m}_3}(\CBetabar,\FBetabar)\cdot\WangD^{\vec{m}_4}(\CRho,\CRhodual,\CBetabar)}{u}{1}{\sphere}{2}{\hypersurbar_{\ubar}}\\
						\notag\lesssim&~\sum\limits_{\substack{\vec{m}_1+\vec{m}_2+\vec{m}_3+\vec{m}_4=(i-i_0,j,k-k_0-1)\\\abs{\vec{m}_1+\vec{m}_2}>N-3}}\delta^{i_0}\cdot\left(\twonorms{\abs{u}^{k_0+3}\WangD^{\vec{m}_1}\psi^{i_0}\cdot\WangD^{\vec{m}_2}\upsi^{k_0}}{u}{1}{\sphere}{2}{\hypersurbar_{\ubar}}\right.\\
						\notag&\left.\cdot\twonorms{\WangD^{\vec{m}_3}(\CBetabar,\FBetabar)}{u}{\infty}{\sphere}{\infty}{\hypersurbar_{\ubar}}\cdot\twonorms{\WangD^{\vec{m}_4}\angnabla^2(\CRho,\CRhodual,\CBetabar)}{u}{\infty}{\sphere}{2}{\hypersurbar_{\ubar}}\right)\\
						\notag&+\sum\limits_{\substack{\vec{m}_1+\vec{m}_2+\vec{m}_3+\vec{m}_4=(i-i_0,j,k-k_0-1)\\\abs{\vec{m}_3}\geq N-3}}\delta^{i_0}\left(\twonorms{\abs{u}^{k_0+3}\WangD^{\vec{m}_1}\psi^{i_0}\cdot\WangD^{\vec{m}_2}\upsi^{k_0}}{u}{\infty}{\sphere}{\infty}{\hypersurbar_{\ubar}}\right.\\
						\notag&\left.\cdot\twonorms{\WangD^{\vec{m}_3}(\CBetabar,\FBetabar)}{u}{1}{\sphere}{2}{\hypersurbar_{\ubar}}\cdot\twonorms{\WangD^{\vec{m}_4}\angnabla^2(\CRho,\CRhodual,\CBetabar)}{u}{\infty}{\sphere}{2}{\hypersurbar_{\ubar}}\right)\\		\notag&+\sum\limits_{\substack{\vec{m}_1+\vec{m}_2+\vec{m}_3+\vec{m}_4=(i-i_0,j,k-k_0-1)\\\abs{\vec{m}_4}\geq N-3}}\delta^{i_0}\left(\twonorms{\abs{u}^{k_0+3}\WangD^{\vec{m}_1}\psi^{i_0}\cdot\WangD^{\vec{m}_2}\upsi^{k_0}}{u}{\infty}{\sphere}{\infty}{\hypersurbar_{\ubar}}\right.\\
						\notag&\left.\cdot\twonorms{\WangD^{\vec{m}_3}(\CBetabar,\FBetabar)}{u}{\infty}{\sphere}{\infty}{\hypersurbar_{\ubar}}\cdot\twonorms{\WangD^{\vec{m}_4}(\CRho,\CRhodual,\CBetabar)}{u}{1}{\sphere}{2}{\hypersurbar_{\ubar}}\right)\\
						\notag\lesssim&~\left(\BAsize\frac{\delta a^{\frac{1}{2}}}{\abs{u}}\right)^{i_0+k_0}\cdot\left(\BAsize\frac{\delta a^{\frac{1}{2}}}{\abs{u}}+  b^{-\frac{1}{2}}\cdot\abs{u}^2\right)\cdot\sum\limits_{\abs{\vec{m}}\leq N-2}\norm{\abs{u}^2\WangD^{\vec{m}}(\CRho,\CRhodual,\CBetabar)}_{L^2(\suubar)}.
					\end{align}
					The rest of the terms are identical to the proofs in Section \ref{SS:proveRicci} and the previous proof in this section. Combining \eqref{PV：rhorhodualbetabar2}-\eqref{PV：rhorhodualbetabar7}, soaking $\sum\limits_{\abs{\vec{m}}\leq N-2}\norm{\WangD^{\vec{m}}(\CRho,\CRhodual,\CBetabar)}_{L^2(\suubar)}$ terms on the RHS of \eqref{PV：rhorhodualbetabar6}, \eqref{PV：rhorhodualbetabar7} to the LHS, we obtain the desired result \eqref{ES:L2curvature3}.
				\end{proof}
				\begin{proof}[Proof of \eqref{ES:L2curvature4}]
					Let $\vec{m}=(i,j,k)$ with $i+j+k\leq N-2$. Let $\WangD^{\vec{m}}$ be as in Definition \ref{DE:operators}. We commute equations \eqref{EQ:D4CAlphabar} with $\WangD^{\vec{m}}$ using \eqref{EQ:ang4ijk}, we obtain:
					\begin{align}\label{PV:Alphabar2}
						& \WangD_4\WangD^{\vec{m}}\CAlphabar\\
						\notag=~& \sum\limits_{\vec{m}_1+\vec{m}_2=(i,j-j_0,k-k_0)}\abs{u}^{j_0+k_0}\WangD^{\vec{m}_1}\upsi^{j_0+k_0}\cdot\WangD^{\vec{m}_2}\angnabla\CBetabar\\	\notag&+\sum\limits_{\vec{m}_1+\vec{m}_2+\vec{m}_3=(i,j-j_0,k-k_0)}\abs{u}^{j_0+k_0}\WangD^{\vec{m}_1}\upsi^{j_0+k_0}\cdot\WangD^{\vec{m}_2}(\psi,\upsi)\cdot\WangD^{\vec{m}_3}\CAlphabar\\
						\notag&+\sum\limits_{\vec{m}_1+\vec{m}_2+\vec{m}_3=(i,j-j_0,k-k_0)}\abs{u}^{j_0+k_0}\WangD^{\vec{m}_1}\upsi^{j_0+k_0}\cdot\WangD^{\vec{m}_2}\psi\cdot\WangD^{\vec{m}_3}\FAlphabar\\
						\notag&+\sum\limits_{\vec{m}_1+\vec{m}_2=(i,j-j_0,k-k_0)}\abs{u}^{j_0+k_0}\WangD^{\vec{m}_1}\upsi^{j_0+k_0}\cdot\WangD^{\vec{m}_2}(\angnabla\FBetabar,\angD_3\Scalar)\\
						\notag&+\sum\limits_{\vec{m}_1+\vec{m}_2=(i,j-j_0,k-k_0)}\abs{u}^{j_0+k_0}\WangD^{\vec{m}_1}\upsi^{j_0+k_0+1}\cdot\WangD^{\vec{m}_3}(\CBetabar,\Rho,\Rhodual,\FRhodual,\Scalar,\FRho)\\
						\notag&+ \sum\limits_{\vec{m}_1+\vec{m}_2+\vec{m}_3=(i,j-j_0,k-k_0-1)}\abs{u}^{j_0+k_0+1}\WangD^{\vec{m}_1}\upsi^{k_0+j_0}\cdot\WangD^{\vec{m}_2}\CAlphabar\cdot\WangD^{\vec{m}_3}(\CBetaA,\FBetaA)\\
						\notag&+ \sum\limits_{\vec{m}_1+\vec{m}_2+\vec{m}_3=(i,j-j_0-1,k-k_0)}\abs{u}^{j_0+k_0+1}\WangD^{\vec{m}_1}\upsi^{j_0+k_0}\cdot\WangD^{\vec{m}_2}\CAlphabar\cdot\WangD^{\vec{m}_3}\Rhodual,\\
						\notag&=:\mathfrak{F}_{(\CAlphabar)}.
					\end{align}
					By Proposition \ref{PR:transport1}, we have:
					\begin{align}\label{PV:Alphabar3}
						\norm{\WangD^{\vec{m}}\CAlphabar}_{L^2(\suubar)}\lesssim\norm{\WangD^{\vec{m}}\CAlphabar}_{L^2(S_{u,0})}+\int_{0}^{\ubar}	\norm{\mathfrak{F}_{(\CAlphabar)}}_{L^2(\suubarprime)}\diff\ubarprime.
					\end{align}
					We now estimate $\int_{\uinf}^{u}\norm{\mathfrak{F}_{(\CAlphabar)}}_{L^2(\suprimeubar)}\diff\uprime$.
					By bootstrap assumptions \eqref{BA}, Proposition \ref{PR:L2RicciCoefficients} and Corollary \ref{CO:LinftyRicciCoefficients}, for the first term on the RHS of \eqref{PV:Alphabar2}, we have:
					\begin{align}\label{PV:Alphabar4}
						&\sum\limits_{\vec{m}_1+\vec{m}_2=(i,j-j_0,k-k_0)}\abs{u}^{j_0+k_0}\twonorms{\WangD^{\vec{m}_1}\upsi^{j_0+k_0}\cdot\WangD^{\vec{m}_2}\angnabla\CBetabar}{\ubar}{1}{\sphere}{2}{\hypersur_{u}}\\
						\notag\lesssim&~\sum\limits_{\substack{\vec{m}_1+\vec{m}_2=(i,j-j_0,k-k_0)\\\abs{\vec{m}_1}>\frac{N-2}{2}}}\delta\abs{u}^{j_0+k_0}\left(\twonorms{\WangD^{\vec{m}_1}\upsi^{j_0+k_0}}{\ubar}{\infty}{\sphere}{2}{\hypersur_{u}}\cdot\twonorms{\WangD^{\vec{m}_2}\angnabla\CBetabar}{\ubar}{\infty}{\sphere}{\infty}{\hypersur_{u}}\right)\\
						\notag&+\sum\limits_{\substack{\vec{m}_1+\vec{m}_2=(i,j-j_0,k-k_0)\\\abs{\vec{m}_2}>\frac{N-2}{2}}}\delta^{\frac{1}{2}}\abs{u}^{j_0+k_0}\left(\twonorms{\WangD^{\vec{m}_1}\upsi^{j_0+k_0}}{\ubar}{\infty}{\sphere}{\infty}{\hypersur_{u}}\cdot\twonorms{\WangD^{\vec{m}_2}\angnabla\CBetabar}{\ubar}{2}{\sphere}{2}{\hypersur_{u}}\right)\\
						\notag\lesssim&~\frac{\delta^{2}a}{\abs{u}^2} \left(a^{-\frac{1}{2}}+\tilde{\mathcal{R}}_{N-1,2}\right)\cdot\left(\BAsize\frac{\delta a^{\frac{1}{2}}}{\abs{u}}\right)^{j_0+k_0},
					\end{align}
					By bootstrap assumptions \eqref{BA} and Sobolev embedding $H^2\hookrightarrow L^\infty$ on $\suubar$ for $\WangD^{\vec{m}}\CAlphabar$, for the second term on the RHS of \eqref{PV:Alphabar2}, we have:
					\begin{align}\label{PV:Alphabar5}
						&\sum\limits_{\vec{m}_1+\vec{m}_2+\vec{m}_3=(i,j-j_0,k-k_0)}\abs{u}^{j_0+k_0}\twonorms{\WangD^{\vec{m}_1}\upsi^{j_0+k_0}\cdot\WangD^{\vec{m}_2}\CAlphabar\cdot\WangD^{\vec{m}_3}(\psi,\upsi)}{\ubar}{1}{\sphere}{2}{\hypersur_{u}}\\
						\notag\lesssim&~\sum\limits_{\substack{\vec{m}_1+\vec{m}_2+\vec{m}_3=(i,j-j_0,k-k_0)\\\abs{\vec{m}_1}>\frac{N-2}{2}}}\delta\abs{u}^{j_0+k_0}\twonorms{\WangD^{\vec{m}_1}\upsi^{j_0+k_0}}{\ubar}{\infty}{\sphere}{2}{\hypersur_{u}}\cdot\twonorms{\WangD^{\vec{m}_2}\CAlphabar}{\ubar}{\infty}{\sphere}{\infty}{\hypersur_{u}}\cdot\twonorms{\WangD^{\vec{m}_3}(\psi,\upsi)}{\ubar}{\infty}{\sphere}{\infty}{\hypersur_{u}}\\
						\notag&+\sum\limits_{\substack{\vec{m}_1+\vec{m}_2+\vec{m}_3=(i,j-j_0,k-k_0)\\\abs{\vec{m}_2}>\frac{N-2}{2}}}\delta\abs{u}^{j_0+k_0}\twonorms{\WangD^{\vec{m}_1}\upsi^{j_0+k_0}}{\ubar}{\infty}{\sphere}{\infty}{\hypersur_{u}}\cdot\twonorms{\WangD^{\vec{m}_2}\CAlphabar}{\ubar}{\infty}{\sphere}{2}{\hypersur_{u}}\cdot\twonorms{\WangD^{\vec{m}_3}(\psi,\upsi)}{\ubar}{\infty}{\sphere}{\infty}{\hypersur_{u}}\\
						\notag&+\sum\limits_{\substack{\vec{m}_1+\vec{m}_2+\vec{m}_3=(i,j-j_0,k-k_0)\\\abs{\vec{m}_3}>\frac{N-2}{2}}}\delta\abs{u}^{j_0+k_0}\twonorms{\WangD^{\vec{m}_1}\upsi^{j_0+k_0}}{\ubar}{\infty}{\sphere}{\infty}{\hypersur_{u}}\cdot\twonorms{\WangD^{\vec{m}_2}\CAlphabar}{\ubar}{\infty}{\sphere}{\infty}{\hypersur_{u}}\cdot\twonorms{\WangD^{\vec{m}_3}(\psi,\upsi)}{\ubar}{\infty}{\sphere}{2}{\hypersur_{u}}\\
						\notag\lesssim&~ \left(\BAsize\frac{\delta a^{\frac{1}{2}}}{\abs{u}}\right)^{j_0+k_0+1}\cdot\sum\limits_{\abs{\vec{m}_2}\leq N-2}\twonorms{\WangD^{\vec{m_2}}\CAlphabar}{\ubar}{\infty}{\sphere}{2}{\hypersur_{u}}.
					\end{align}
					By Proposition \ref{PR:EstimatesforRicciCurvature} and bootstrap assumptions \eqref{BA}, for the third term on the RHS of \eqref{PV:Alphabar2}, we have:
					\begin{align}\label{PV:Alphabar6}
						&\sum\limits_{\vec{m}_1+\vec{m}_2+\vec{m}_3=(i,j-j_0,k-k_0)}\abs{u}^{j_0+k_0}\twonorms{\WangD^{\vec{m}_1}\upsi^{j_0+k_0}\cdot\WangD^{\vec{m}_2}\FAlphabar\cdot\WangD^{\vec{m}_3}\psi}{\ubar}{1}{\sphere}{2}{\hypersur_{u}}\\
						\notag\lesssim&~\sum\limits_{\substack{\vec{m}_1+\vec{m}_2+\vec{m}_3=(i,j-j_0,k-k_0)\\\abs{\vec{m}_1}>\frac{N-2}{2}}}\delta\abs{u}^{j_0+k_0}\twonorms{\WangD^{\vec{m}_1}\upsi^{j_0+k_0}}{\ubar}{\infty}{\sphere}{2}{\hypersur_{u}}\cdot\twonorms{\WangD^{\vec{m}_2}\FAlphabar}{\ubar}{\infty}{\sphere}{\infty}{\hypersur_{u}}\cdot\twonorms{\WangD^{\vec{m}_3}\psi}{\ubar}{\infty}{\sphere}{\infty}{\hypersur_{u}}\\
						\notag&+\sum\limits_{\substack{\vec{m}_1+\vec{m}_2+\vec{m}_3=(i,j-j_0,k-k_0)\\\abs{\vec{m}_2}>\frac{N-2}{2}}}\delta\abs{u}^{j_0+k_0}\twonorms{\WangD^{\vec{m}_1}\upsi^{j_0+k_0}}{\ubar}{\infty}{\sphere}{\infty}{\hypersur_{u}}\cdot\twonorms{\WangD^{\vec{m}_2}\FAlphabar}{\ubar}{\infty}{\sphere}{2}{\hypersur_{u}}\cdot\twonorms{\WangD^{\vec{m}_3}\psi}{\ubar}{\infty}{\sphere}{\infty}{\hypersur_{u}}\\
						\notag&+\sum\limits_{\substack{\vec{m}_1+\vec{m}_2+\vec{m}_3=(i,j-j_0,k-k_0)\\\abs{\vec{m}_3}>\frac{N-2}{2}}}\delta\abs{u}^{j_0+k_0}\twonorms{\WangD^{\vec{m}_1}\upsi^{j_0+k_0}}{\ubar}{\infty}{\sphere}{\infty}{\hypersur_{u}}\cdot\twonorms{\WangD^{\vec{m}_2}\FAlphabar}{\ubar}{\infty}{\sphere}{\infty}{\hypersur_{u}}\cdot\twonorms{\WangD^{\vec{m}_3}\psi}{\ubar}{\infty}{\sphere}{2}{\hypersur_{u}}\\
						\notag\lesssim&~ \left(\BAsize\frac{\delta a^{\frac{1}{2}}}{\abs{u}}\right)^{j_0+k_0+1}\cdot b^{-1}
					\end{align}
					In below, using Proposition \ref{PR:EstimatesforRicciCurvature} and bootstrap assumptions \eqref{BA}, we also need the Sobolev embedding $H^2\hookrightarrow L^\infty$ on $\suubar$ for $\WangD^{\vec{m}}\CAlphabar$ to obtain the estimate for the second last term on the RHS of \eqref{PV:Alphabar2}:
					\begin{align}\label{PV:Alphabar7}
						&\sum\limits_{\vec{m}_1+\vec{m}_2+\vec{m}_3=(i,j-j_0,k-k_0-1)}\abs{u}^{j_0+k_0+1}\twonorms{\WangD^{\vec{m}_1}\upsi^{k_0+j_0}\cdot\WangD^{\vec{m}_2}\CAlphabar\cdot\WangD^{\vec{m}_3}(\CBetaA,\FBetaA)}{\ubar}{1}{\sphere}{2}{\hypersur_{u}}\\
						\notag\lesssim&~\sum\limits_{\substack{\vec{m}_1+\vec{m}_2+\vec{m}_3=(i,j-j_0,k-k_0-1)\\\abs{\vec{m}_1}>\frac{N-2}{2}}}\delta\abs{u}^{j_0+k_0+1}\left(\twonorms{\WangD^{\vec{m}_1}\upsi^{j_0+k_0}}{\ubar}{\infty}{\sphere}{2}{\hypersur_{u}}\cdot\twonorms{\WangD^{\vec{m}_2}\CAlphabar}{\ubar}{\infty}{\sphere}{\infty}{\hypersur_{u}}\cdot\twonorms{\WangD^{\vec{m}_3}(\CBetaA,\FBetaA)}{\ubar}{\infty}{\sphere}{\infty}{\hypersur_{u}}\right)\\
						\notag&+\sum\limits_{\substack{\vec{m}_1+\vec{m}_2+\vec{m}_3=(i,j-j_0,k-k_0-1)\\\abs{\vec{m}_2}>\frac{N-2}{2}}}\delta\abs{u}^{j_0+k_0+1}\left(\twonorms{\WangD^{\vec{m}_1}\upsi^{j_0+k_0}}{\ubar}{\infty}{\sphere}{\infty}{\hypersur_{u}}\cdot\twonorms{\WangD^{\vec{m}_2}\CAlphabar}{\ubar}{\infty}{\sphere}{2}{\hypersur_{u}}\cdot\twonorms{\WangD^{\vec{m}_3}(\CBetaA,\FBetaA)}{\ubar}{\infty}{\sphere}{\infty}{\hypersur_{u}}\right)\\
						\notag&+\sum\limits_{\substack{\vec{m}_1+\vec{m}_2+\vec{m}_3=(i,j-j_0,k-k_0-1)\\\abs{\vec{m}_3}>\frac{N-2}{2}}}\delta\abs{u}^{j_0+k_0+1}\left(\twonorms{\WangD^{\vec{m}_1}\upsi^{j_0+k_0}}{\ubar}{\infty}{\sphere}{\infty}{\hypersur_{u}}\cdot\twonorms{\WangD^{\vec{m}_2}\CAlphabar}{\ubar}{\infty}{\sphere}{\infty}{\hypersur_{u}}\cdot\twonorms{\WangD^{\vec{m}_3}(\CBetaA,\FBetaA)}{\ubar}{\infty}{\sphere}{2}{\hypersur_{u}}\right)\\
						\notag\lesssim&~ \left( \delta^{\frac{1}{2}}+\BAsize\frac{\delta a^{\frac{1}{2}}}{\abs{u}^2}\right)\cdot\left(\BAsize\frac{\delta a^{\frac{1}{2}}}{\abs{u}}\right)^{j_0+k_0}\cdot\sum\limits_{\abs{\vec{m}_2}\leq N-2}\twonorms{\WangD^{\vec{m_2}}\CAlphabar}{\ubar}{\infty}{\sphere}{2}{\hypersur_{u}}.
					\end{align}
					Similarly, for the last term on the RHS of \eqref{PV:Alphabar2}:
					\begin{align}\label{PV:Alphabar9}
						&\sum\limits_{\vec{m}_1+\vec{m}_2+\vec{m}_3=(i,j-j_0-1,k-k_0)}\abs{u}^{j_0+k_0+1}\twonorms{\WangD^{\vec{m}_1}\upsi^{j_0+k_0}\cdot\WangD^{\vec{m}_2}\CAlphabar\cdot\WangD^{\vec{m}_3}\Rhodual}{\ubar}{1}{\sphere}{2}{\hypersur_{u}}\\
						\notag\lesssim&~\sum\limits_{\substack{\vec{m}_1+\vec{m}_2+\vec{m}_3=(i,j-j_0-1,k-k_0)\\\abs{\vec{m}_1}>\frac{N-2}{2}}}\delta\abs{u}^{j_0+k_0+1}\twonorms{\WangD^{\vec{m}_1}\upsi^{j_0+k_0}}{\ubar}{\infty}{\sphere}{2}{\hypersur_{u}}\cdot\twonorms{\WangD^{\vec{m}_2}\CAlphabar}{\ubar}{\infty}{\sphere}{\infty}{\hypersur_{u}}\cdot\twonorms{\WangD^{\vec{m}_3}\Rhodual}{\ubar}{\infty}{\sphere}{\infty}{\hypersur_{u}}\\
						\notag&+\sum\limits_{\substack{\vec{m}_1+\vec{m}_2+\vec{m}_3=(i,j-j_0,k-k_0)\\\abs{\vec{m}_2}>\frac{N-2}{2}}}\delta\abs{u}^{j_0+k_0+1}\twonorms{\WangD^{\vec{m}_1}\upsi^{j_0+k_0}}{\ubar}{\infty}{\sphere}{\infty}{\hypersur_{u}}\cdot\twonorms{\WangD^{\vec{m}_2}\CAlphabar}{\ubar}{\infty}{\sphere}{2}{\hypersur_{u}}\cdot\twonorms{\WangD^{\vec{m}_3}\Rhodual}{\ubar}{\infty}{\sphere}{\infty}{\hypersur_{u}}\\
						\notag&+\sum\limits_{\substack{\vec{m}_1+\vec{m}_2+\vec{m}_3=(i,j-j_0-1,k-k_0)\\\abs{\vec{m}_3}>\frac{N-2}{2}}}\delta\abs{u}^{j_0+k_0+1}\twonorms{\WangD^{\vec{m}_1}\upsi^{j_0+k_0}}{\ubar}{\infty}{\sphere}{\infty}{\hypersur_{u}}\cdot\twonorms{\WangD^{\vec{m}_2}\CAlphabar}{\ubar}{\infty}{\sphere}{\infty}{\hypersur_{u}}\cdot\twonorms{\WangD^{\vec{m}_3}\Rhodual}{\ubar}{\infty}{\sphere}{2}{\hypersur_{u}}\\
						\notag\lesssim&~ \left(\BAsize\frac{\delta a^{\frac{1}{2}}}{\abs{u}}\right)^{j_0+k_0+2}\cdot\sum\limits_{\abs{\vec{m}_2}\leq N-2}\twonorms{\WangD^{\vec{m_2}}\CAlphabar}{\ubar}{\infty}{\sphere}{2}{\hypersur_{u}}.
					\end{align}
					The rest of the terms can be controlled using the methods identical to the proofs in Section \ref{SS:proveRicci}. In particular, the fourth term and the fifth term on the RHS of \eqref{PV:Alphabar2} can be estimated using the same approach\footnote{We note that $\BetaA$ term in \eqref{PV:upsi9} is absent in the fourth term and the fifth term on the RHS of \eqref{PV:Alphabar2}. Hence we obtain a smaller upper bound $\delta\BAsize^{-1}\left(\BAsize\frac{\delta a^{\frac{1}{2}}}{\abs{u}}\right)^{i_0+k_0+2}$ for the term $\sum\limits_{\vec{m}_1+\vec{m}_2=(i,j-j_0,k-k_0)}\abs{u}^{j_0+k_0}\WangD^{\vec{m}_1}\upsi^{j_0+k_0+1}\cdot\WangD^{\vec{m}_3}(\CBetabar,\Rho,\Rhodual)$.} as in \eqref{PV:upsi7} and \eqref{PV:upsi9}. Combining \eqref{PV:Alphabar3}-\eqref{PV:Alphabar9}, we soak the $\sum\limits_{\abs{\vec{m}_2}\leq N-2}\twonorms{\WangD^{\vec{m}}\CAlphabar}{\ubar}{\infty}{\sphere}{2}{\hypersur_{u}}$ terms to the left since their coefficients are $\ll 1$. Then, we obtain the desired estimates \eqref{ES:L2curvature4} by the initial data assumption \eqref{Data} and the choice of parameter $\delta,a,b$ in Section \ref{SS:Norms}, $\BAfluid,\BAsize$ in Section \eqref{SS:BA}.
				\end{proof}
				\subsubsection{Proof of Proposition \ref{PR:topL2curvature}}\label{SSS:proofofcurvature}
				\begin{proof}[Proof of \eqref{ES:topL2curvature1}]
					Let $\vec{m}=(i,j,k)$ with $i+j+k\leq N-1$ with $k\geq 1$. Let $\WangD^{\vec{m}}$ be as in Definition \ref{DE:operators}. We commute equations \eqref{EQ:D3CAlpha} and \eqref{EQ:D4CBeta} with $\WangD^{\vec{m}}$ using \eqref{EQ:ang3ijk}-\eqref{EQ:ang4ijk}. Let $F:=\CAlphaAB$, $G:=\CBetaA$ and
					
					\begin{align}
						\label{PV:AlphaBeta1}P_{(\CAlphaAB,\CBetaA;\vec{m})}=&~\sum\limits_{\vec{m_1}+\vec{m_2}+\vec{m_3}=(i-i_0,j,k-k_0)}\delta^{i_0}\abs{u}^{k_0}\WangD^{\vec{m_1}}\upsi^{k_0+1}\cdot\WangD^{\vec{m_2}}\psi^{i_0}\cdot\WangD^{\vec{m_3}}(\CAlphaAB,\CBetaA,\FAlphaAB)\\
						\notag&+\sum\limits_{\vec{m_1}+\vec{m_2}+\vec{m_3}=(i-i_0,j,k-k_0)}\delta^{i_0}\abs{u}^{k_0}\WangD^{\vec{m_1}}\upsi^{k_0}\cdot\WangD^{\vec{m_2}}\psi^{i_0+1}\cdot\WangD^{\vec{m_3}}(\CRho,\CRhodual,\Scalar,\FRho,\FRhodual)\\
						\notag&+\sum\limits_{\vec{m_1}+\vec{m_2}+\vec{m_3}=(i-i_0,j,k-k_0)}\delta^{i_0}\abs{u}^{k_0}\WangD^{\vec{m_1}}\upsi^{k_0}\cdot\WangD^{\vec{m_2}}\psi^{i_0}\cdot\WangD^{\vec{m_3}}(\angnabla\FBetaA,\angD_4\Scalar)\\
						\notag&+\sum\limits_{\vec{m_1}+\vec{m_2}+\vec{m_3}+\vec{m_4}=(i-i_0-1,j,k-k_0+1)}\delta^{i_0+1}\abs{u}^{k_0-1}\WangD^{\vec{m_1}}\upsi^{k_0}\cdot\WangD^{\vec{m_2}}\psi^{i_0}\cdot\WangD^{\vec{m_3}}\upzeta\cdot\WangD^{\vec{m_4}}\CAlphaAB\\
						\notag&+\sum\limits_{\vec{m_1}+\vec{m_2}+\vec{m_3}+\vec{m_4}=(i-i_0-1,j,k-k_0)}\delta^{i_0+1}\abs{u}^{k_0}\WangD^{\vec{m_1}}\upsi^{k_0}\cdot\WangD^{\vec{m_2}}\psi^{i_0}\cdot\WangD^{\vec{m_3}}\Rhodual\cdot\WangD^{\vec{m_4}}\CAlphaAB\\
						\notag&+\sum\limits_{\vec{m_1}+\vec{m_2}+\vec{m_3}+\vec{m_4}=(i-i_0,j,k-k_0-1)}\delta^{i_0}\abs{u}^{k_0+1}\WangD^{\vec{m_1}}\upsi^{k_0}\cdot\WangD^{\vec{m_2}}\psi^{i_0}\cdot\WangD^{\vec{m_3}}(\CBetabar,\FBetabar)\cdot\WangD^{\vec{m_4}}\CAlphaAB,	
					\end{align}
					\begin{align}
						\label{PV:AlphaBeta2}Q_{(\CAlphaAB,\CBetaA;\vec{m})}=&~\sum\limits_{\vec{m_1}+\vec{m_2}+\vec{m_3}=(i,j- j_0,k-k_0)}\abs{u}^{j_0+k_0}\WangD^{\vec{m_1}}\upsi^{k_0+j_0}\cdot\WangD^{\vec{m_2}}\psi\cdot\WangD^{\vec{m_3}}(\CBetaA,\FBetaA)\\
						\notag&+\sum\limits_{\vec{m_1}+\vec{m_2}=(i,j- j_0,k-k_0)}\abs{u}^{j_0+k_0}\WangD^{\vec{m_1}}\upsi^{k_0+j_0+1}\cdot\WangD^{\vec{m_2}}(\CAlphaAB,\FAlphaAB)\\
						\notag&+\sum\limits_{\vec{m_1}+\vec{m_2}+\vec{m_3}+\vec{m_4}=(i,j- j_0,k-k_0-1)}\abs{u}^{j_0+k_0+1}\WangD^{\vec{m_1}}\upsi^{j_0+k_0}\cdot\WangD^{\vec{m_2}}(\CBetaA,\FBetaA)\cdot\WangD^{\vec{m_3}}\CBetaA\\
						\notag&+\sum\limits_{\vec{m_1}+\vec{m_2}+\vec{m_3}=(i,j- j_0-1,k-k_0+1)}\abs{u}^{j_0+k_0}\WangD^{\vec{m_1}}\upsi^{j_0+k_0}\cdot\WangD^{\vec{m_2}}\upzeta\cdot\WangD^{\vec{m_3}}\CBetaA\\
						\notag&+\sum\limits_{\vec{m_1}+\vec{m_2}+\vec{m_3}=(i,j- j_0-1,k-k_0)}\abs{u}^{j_0+k_0+1}\WangD^{\vec{m_1}}\upsi^{j_0+k_0}\cdot\WangD^{\vec{m_2}}\Rhodual\cdot\angD^{\vec{m_3}}\CBetaA.
					\end{align}
					We apply Proposition \ref{PR:Int4} with $s(\CAlphaAB)=0$ to obtain
					\begin{align}\label{PV:AlphaBeta3}
						&\twonorms{\WangD^{\vec{m}}\CAlphaAB}{\ubar}{2}{\sphere}{2}{\hypersur_{u}}^2+\twonorms{\WangD^{\vec{m}}\CBetaA}{u}{2}{\sphere}{2}{\hypersurbar_{\ubar}}^2\\
						\notag\lesssim&~\twonorms{\WangD^{\vec{m}}\CAlphaAB}{\ubar}{2}{\sphere}{2}{\hypersur_{\uinf}}^2+\twonorms{\WangD^{\vec{m}}\CBetaA}{u}{2}{\sphere}{2}{\hypersurbar_{0}}^2+\int_{\region}\abs{\WangD^{\vec{m}}\CAlphaAB}\cdot \abs{P_{(\CAlphaAB,\CBetaA;\vec{m})}}\diff\gvol+\int_{\region}\abs{\WangD^{\vec{m}}\CBetaA}\cdot \abs{Q_{(\CAlphaAB,\CBetaA;\vec{m})}}\diff\gvol.
					\end{align}
					We focus on estimating $\threenorms{P_{(\CAlphaAB,\CBetaA;\vec{m})}}{u}{1}{\ubar}{2}{\sphere}{2}{\region}$ and $\threenorms{Q_{(\CAlphaAB,\CBetaA;\vec{m})}}{u}{2}{\ubar}{2}{\sphere}{2}{\region}$ since 
					\begin{align}\label{PV:AlphaBeta4}
						&\int_{\region}\abs{\WangD^{\vec{m}}\CAlphaAB}\cdot \abs{P_{(\CAlphaAB,\CBetaA;\vec{m})}}\diff\gvol+\int_{\region}\abs{\WangD^{\vec{m}}\CBetaA}\cdot \abs{Q_{(\CAlphaAB,\CBetaA;\vec{m})}}\diff\gvol\\
						\notag\lesssim&~\threenorms{\WangD^{\vec{m}}\CAlphaAB}{u}{\infty}{\ubar}{2}{\sphere}{2}{\region}\threenorms{P_{(\CAlphaAB,\CBetaA;\vec{m})}}{u}{1}{\ubar}{2}{\sphere}{2}{\region}+\delta^{\frac{1}{2}}\threenorms{\WangD^{\vec{m}}\CBetaA}{u}{2}{\ubar}{\infty}{\sphere}{2}{\region}\threenorms{Q_{(\CAlphaAB,\CBetaA;\vec{m})}}{u}{2}{\ubar}{2}{\sphere}{2}{\region}.
					\end{align}
					For the terms involving $\CAlphaAB$, namely, the first term and the last three terms on the RHS of \eqref{PV:AlphaBeta1}, using Proposition \ref{PR:EstimatesforRicciCurvature} for $\FBetabar$, bootstrap assumptions \eqref{BA} and Sobolev embedding $H^2\hookrightarrow L^\infty$ on $\suubar$ for $\WangD^{\vec{m}}\CAlphaAB$, we have:
					\begin{align}\label{PV:AlphaBeta5}
						&\sum\limits_{\vec{m}_1+\vec{m}_2+\vec{m}_3+\vec{m}_4=(i-i_0,j,k-k_0)}\delta^{i_0}\threenorms{\abs{u}^{k_0}\WangD^{\vec{m}_1}\psi^{i_0}\cdot\WangD^{\vec{m}_2}\upsi^{k_0}\cdot\left(\WangD^{\vec{m}_3}\upsi,\delta\WangD^{\vec{m}'_3}\Rhodual,\abs{u}\WangD^{\vec{m}'_3}\CBetabar,\abs{u}\WangD^{\vec{m}'_3}\FBetabar\right)\cdot\WangD^{\vec{m}_4}\CAlphaAB}{u}{1}{\ubar}{2}{\sphere}{2}{\region}\\
						\notag\lesssim&~\sum\limits_{\substack{\vec{m}_1+\vec{m}_2+\vec{m}_3+\vec{m}_4=(i-i_0,j,k-k_0)\\\abs{\vec{m}_4}<\frac{N-1}{2}}}\delta^{i_0}\threenorms{\abs{u}^{k_0}\WangD^{\vec{m}_1}\psi^{i_0}\cdot\WangD^{\vec{m}_2}\upsi^{k_0}\cdot\left(\WangD^{\vec{m}_3}\upsi,\delta\WangD^{\vec{m}'_3}\Rhodual,\abs{u}\WangD^{\vec{m}'_3}\CBetabar,\abs{u}\WangD^{\vec{m}'_3}\FBetabar\right)}{u}{1}{\ubar}{\infty}{\sphere}{2}{\region}\\
						\notag&\cdot\threenorms{\WangD^{\vec{m}_4}\CAlphaAB}{u}{\infty}{\ubar}{2}{\sphere}{\infty}{\region}\\
						\notag&+\sum\limits_{\substack{\vec{m}_1+\vec{m}_2+\vec{m}_3+\vec{m}_4=(i-i_0,j,k-k_0)\\\abs{\vec{m}_4}\geq\frac{N-1}{2}}}\delta^{i_0}\threenorms{\abs{u}^{k_0}\WangD^{\vec{m}_1}\psi^{i_0}\cdot\WangD^{\vec{m}_2}\upsi^{k_0}\cdot\left(\WangD^{\vec{m}_3}\upsi,\delta\WangD^{\vec{m}'_3}\Rhodual,\abs{u}\WangD^{\vec{m}'_3}\CBetabar,\abs{u}\WangD^{\vec{m}'_3}\FBetabar\right)}{u}{1}{\ubar}{\infty}{\sphere}{\infty}{\region}\\
						\notag&\cdot\threenorms{\WangD^{\vec{m}_4}\CAlphaAB}{u}{\infty}{\ubar}{2}{\sphere}{2}{\region}\\
						\notag\lesssim&~\left(\BAsize\frac{\delta a^{\frac{1}{2}}}{\abs{u}}\right)^{i_0+k_0}\left(\BAsize\frac{\delta a^{\frac{1}{2}}}{\abs{u}}+ b^{-\frac{1}{2}}\right)\cdot\sum\limits_{\abs{\vec{m}_4}\leq N-1}\threenorms{\WangD^{\vec{m}_4}\CAlphaAB}{u}{\infty}{\ubar}{2}{\sphere}{2}{\region},
					\end{align}
					where $\abs{\vec{m}_3'}=\abs{\vec{m}_3}-1$.
					Next, we estimate the rest terms on the RHS of \eqref{PV:AlphaBeta1}. For the terms only depends on gravity, using bootstrap assumptions \eqref{BA}:
					\begin{align}\label{PV:AlphaBeta6}
						&\sum\limits_{\vec{m}_1+\vec{m}_2+\vec{m}_3+\vec{m}_4=(i-i_0,j,k-k_0)}\delta^{i_0}\threenorms{\abs{u}^{k_0}\WangD^{\vec{m}_1}\psi^{i_0}\cdot\WangD^{\vec{m}_2}\upsi^{k_0}\cdot\left(\WangD^{\vec{m}_3}\upsi\cdot\WangD^{\vec{m}_4}\CBetaA,\WangD^{\vec{m}_3}\psi\cdot\WangD^{\vec{m}_4}(\CRho,\CRhodual)\right)}{u}{1}{\ubar}{2}{\sphere}{2}{\region}\\
						\notag\lesssim&~\sum\limits_{\substack{\vec{m}_1+\vec{m}_2+\vec{m}_3=(i-i_0,j,k-k_0)\\\abs{\vec{m}_3}<\frac{N-1}{2}}}\delta^{i_0}\left(\threenorms{\abs{u}^{k_0-1}\WangD^{\vec{m}_1}\psi^{i_0}\cdot\WangD^{\vec{m}_2}\upsi^{k_0+1}}{u}{1}{\ubar}{2}{\sphere}{2}{\region}\cdot\threenorms{\abs{u}\WangD^{\vec{m}_3}\CBetaA}{u}{\infty}{\ubar}{\infty}{\sphere}{\infty}{\region}\right.\\
						\notag&\left.+\threenorms{\abs{u}^{k_0-2}\WangD^{\vec{m}_1}\psi^{i_0+1}\cdot\WangD^{\vec{m}_2}\upsi^{k_0}}{u}{1}{\ubar}{2}{\sphere}{2}{\region}\cdot\threenorms{\abs{u}^2\WangD^{\vec{m}_3}(\CRho,\CRhodual)}{u}{\infty}{\ubar}{\infty}{\sphere}{\infty}{\region}\right)\\
						\notag&+\sum\limits_{\substack{\vec{m}_1+\vec{m}_2+\vec{m}_3=(i-i_0,j,k-k_0)\\\abs{\vec{m}_3}\geq\frac{N-1}{2}}}\delta^{i_0}\left(\threenorms{\abs{u}^{k_0}\WangD^{\vec{m}_1}\psi^{i_0}\cdot\WangD^{\vec{m}_2}\upsi^{k_0+1}}{u}{2}{\ubar}{\infty}{\sphere}{\infty}{\region}\cdot\norm{\abs{u}^{-1}}_{L^2_u}\cdot\threenorms{\abs{u}\WangD^{\vec{m}_3}\CBetaA}{u}{\infty}{\ubar}{2}{\sphere}{2}{\region}\right.\\
						\notag&\left.+\threenorms{\abs{u}^{k_0}\WangD^{\vec{m}_1}\psi^{i_0+1}\cdot\WangD^{\vec{m}_2}\upsi^{k_0}}{u}{2}{\ubar}{\infty}{\sphere}{\infty}{\region}\cdot\norm{\abs{u}^{-2}}_{L^2_u}\cdot\threenorms{\abs{u}^2\WangD^{\vec{m}_3}(\CRho,\CRhodual)}{u}{\infty}{\ubar}{2}{\sphere}{2}{\region}\right)\\
						\notag\lesssim&~\left(\BAsize\frac{\delta a^{\frac{1}{2}}}{\abs{u}}\right)^{i_0+k_0}\left(\BAsize\frac{\delta a^{\frac{1}{2}}}{\abs{u}}+1\right)\cdot\BAsize\delta^{\frac{1}{2}}a^{\frac{1}{2}}.
					\end{align}
					For the terms that depends on fluids, using bootstrap assumptions \eqref{BA} and Proposition \ref{PR:EstimatesforRicciCurvature}:
					\begin{align}\label{PV:AlphaBeta7}
						&\sum\limits_{\vec{m}_1+\vec{m}_2+\vec{m}_3+\vec{m}_4=(i-i_0,j,k-k_0)}\delta^{i_0}\threenorms{\abs{u}^{k_0}\WangD^{\vec{m}_1}\psi^{i_0}\cdot\WangD^{\vec{m}_2}\upsi^{k_0}\cdot\left(\WangD^{\vec{m}_3}\upsi\cdot\WangD^{\vec{m}_4}\FAlphaAB,\WangD^{\vec{m}_3}\psi\cdot\WangD^{\vec{m}_4}(\Scalar,\FRho,\FRhodual)\right)}{u}{1}{\ubar}{2}{\sphere}{2}{\region}\\
						\notag\lesssim&~\sum\limits_{\substack{\vec{m}_1+\vec{m}_2+\vec{m}_3=(i-i_0,j,k-k_0)\\\abs{\vec{m}_3}<\frac{N-1}{2}}}\delta^{i_0}\left(\threenorms{\abs{u}^{k_0}\WangD^{\vec{m}_1}\psi^{i_0}\cdot\WangD^{\vec{m}_2}\upsi^{k_0+1}}{u}{1}{\ubar}{2}{\sphere}{2}{\region}\cdot\threenorms{\WangD^{\vec{m}_3}\FAlphaAB}{u}{\infty}{\ubar}{\infty}{\sphere}{\infty}{\region}\right.\\
						\notag&\left.+\threenorms{\abs{u}^{k_0}\WangD^{\vec{m}_1}\psi^{i_0+1}\cdot\WangD^{\vec{m}_2}\upsi^{k_0}}{u}{1}{\ubar}{2}{\sphere}{2}{\region}\cdot\threenorms{\WangD^{\vec{m}_3}(\Scalar,\FRho,\FRhodual)}{u}{\infty}{\ubar}{\infty}{\sphere}{\infty}{\region}\right)\\
						\notag&+\sum\limits_{\substack{\vec{m}_1+\vec{m}_2+\vec{m}_3=(i-i_0,j,k-k_0)\\\abs{\vec{m}_3}\geq\frac{N-1}{2}}}\delta^{i_0}\left(\threenorms{\abs{u}^{k_0}\WangD^{\vec{m}_1}\psi^{i_0}\cdot\WangD^{\vec{m}_2}\upsi^{k_0+1}}{u}{1}{\ubar}{\infty}{\sphere}{\infty}{\region}\cdot\threenorms{\WangD^{\vec{m}_3}\FAlphaAB}{u}{\infty}{\ubar}{2}{\sphere}{2}{\region}\right.\\
						\notag&\left.+\threenorms{\abs{u}^{k_0}\WangD^{\vec{m}_1}\psi^{i_0+1}\cdot\WangD^{\vec{m}_2}\upsi^{k_0}}{u}{1}{\ubar}{\infty}{\sphere}{\infty}{\region}\cdot\threenorms{\WangD^{\vec{m}_3}(\Scalar,\FRho,\FRhodual)}{u}{\infty}{\ubar}{2}{\sphere}{2}{\region}\right)\\
						\notag\lesssim&~\left(\BAsize\frac{\delta a^{\frac{1}{2}}}{\abs{u}}\right)^{i_0+k_0}\BAsize  a^{\frac{1}{2}}.
					\end{align}
					\begin{align}\label{PV:AlphaBeta8}
						&\sum\limits_{\vec{m}_1+\vec{m}_2+\vec{m}_3=(i-i_0,j,k-k_0)}\delta^{i_0}\threenorms{\abs{u}^{k_0}\WangD^{\vec{m}_1}\psi^{i_0}\cdot\WangD^{\vec{m}_2}\upsi^{k_0}\cdot\WangD^{\vec{m}_3}(\angnabla\FBetaA,\angD_{ 4}\Scalar)}{u}{1}{\ubar}{2}{\sphere}{2}{\region}\\
						\notag\lesssim&~\sum\limits_{\substack{\vec{m}_1+\vec{m}_2+\vec{m}_3=(i-i_0,j,k-k_0)\\\abs{\vec{m}_3}<\frac{N-1}{2}}}\delta^{i_0}\threenorms{\abs{u}^{k_0}\WangD^{\vec{m}_1}\psi^{i_0}\cdot\WangD^{\vec{m}_2}\upsi^{k_0}}{u}{1}{\ubar}{2}{\sphere}{2}{\region}\cdot\threenorms{\WangD^{\vec{m}_3}(\angnabla\FBetaA,\angD_{ 4}\Scalar)}{u}{\infty}{\ubar}{\infty}{\sphere}{\infty}{\region}\\
						\notag&+\sum\limits_{\substack{\vec{m}_1+\vec{m}_2+\vec{m}_3=(i-i_0,j,k-k_0)\\\abs{\vec{m}_3}\geq\frac{N-1}{2}}}\delta^{i_0}\threenorms{\abs{u}^{k_0}\WangD^{\vec{m}_1}\psi^{i_0}\cdot\WangD^{\vec{m}_2}\upsi^{k_0}}{u}{1}{\ubar}{\infty}{\sphere}{\infty}{\region}\cdot\threenorms{\WangD^{\vec{m}_3}(\angnabla\FBetaA,\angD_{ 4}\Scalar)}{u}{\infty}{\ubar}{2}{\sphere}{2}{\region}\\
						\notag\lesssim&~\left(\BAsize\frac{\delta a^{\frac{1}{2}}}{\abs{u}}\right)^{i_0+k_0}  \frac{\delta^{-1} b^{-\frac{1}{2}}}{\abs{u}}.
					\end{align}
					Now, we estimate the terms on the RHS of \eqref{PV:AlphaBeta2}. For the terms with $\CBetaA$ present, using Proposition \ref{PR:EstimatesforRicciCurvature} for $\FBetabar$, bootstrap assumptions \eqref{BA}, we have:
					\begin{align}\label{PV:AlphaBeta9}
						&\sum\limits_{\vec{m}_1+\vec{m}_2+\vec{m}_3=(i,j-j_0,k-k_0)}\threenorms{\abs{u}^{j_0+k_0}\WangD^{\vec{m}_1}\upsi^{j_0+k_0}\cdot\left(\WangD^{\vec{m}_2}(\psi,\upsi),\abs{u}\WangD^{\vec{m}'_2}(\CBetaA,\FBetaA,\Rhodual)\right)\cdot\WangD^{\vec{m}_3}(\CBetaA,\FBetaA)}{u}{2}{\ubar}{2}{\sphere}{2}{\region}\\
						\notag\lesssim&~\sum\limits_{\substack{\vec{m}_1+\vec{m}_2+\vec{m}_3=(i,j-j_0,k-k_0)\\\abs{\vec{m}_3}<\frac{N-1}{2}}}\threenorms{\abs{u}^{j_0+k_0}\WangD^{\vec{m}_1}\upsi^{j_0+k_0}\cdot\left(\WangD^{\vec{m}_2}(\psi,\upsi),\abs{u}\WangD^{\vec{m}_2}(\CBetaA,\FBetaA,\Rhodual)\right)}{u}{2}{\ubar}{2}{\sphere}{2}{\region}\\
						\notag&\cdot\threenorms{\WangD^{\vec{m}_3}(\CBetaA,\FBetaA)}{u}{\infty}{\ubar}{\infty}{\sphere}{\infty}{\region}\\
						\notag&+\sum\limits_{\substack{\vec{m}_1+\vec{m}_2+\vec{m}_3=(i,j-j_0,k-k_0)\\\abs{\vec{m}_3}\geq\frac{N-1}{2}}}\threenorms{\abs{u}^{j_0+k_0}\WangD^{\vec{m}_1}\upsi^{j_0+k_0}\cdot\left(\WangD^{\vec{m}_2}(\psi,\upsi),\abs{u}\WangD^{\vec{m}'_2}(\CBetaA,\FBetaA,\Rhodual)\right)}{u}{2}{\ubar}{\infty}{\sphere}{\infty}{\region}\\
						\notag&\cdot\threenorms{\WangD^{\vec{m}_3}(\CBetaA,\FBetaA)}{u}{\infty}{\ubar}{2}{\sphere}{2}{\region}\\
						\notag\lesssim&~\left(\BAsize\frac{\delta a^{\frac{1}{2}}}{\abs{u}}\right)^{j_0+k_0}\left(\BAsize\frac{\delta^{\frac{1}{2}}a}{\abs{u}^{\frac{3}{2}}}+ \frac{a^{\frac{1}{2}}}{\abs{u}^{\frac{1}{2}}}+ \delta^{-\frac{1}{2}}\right),
					\end{align}
					where $\abs{\vec{m}_2'}=\abs{\vec{m}_2}-1$.
					Finally, we estimate the second term on the RHS of \eqref{PV:AlphaBeta2} using Proposition \ref{PR:EstimatesforRicciCurvature}, bootstrap assumptions \eqref{BA}, and Corollary \ref{CO:Linftycurvature}:
					\begin{align}\label{PV:AlphaBeta10}
						&\sum\limits_{\vec{m}_1+\vec{m}_2=(i,j-j_0,k-k_0)}\threenorms{\abs{u}^{j_0+k_0}\WangD^{\vec{m}_1}\upsi^{j_0+k_0+1}\cdot\WangD^{\vec{m}_2}(\CAlphaAB,\FAlphaAB)}{u}{2}{\ubar}{2}{\sphere}{2}{\region}\\
						\notag\lesssim&~\sum\limits_{\substack{\vec{m}_1+\vec{m}_2=(i,j-j_0,k-k_0)\\\abs{\vec{m}_2}<\frac{N-1}{2}}}\threenorms{\abs{u}^{j_0+k_0}\WangD^{\vec{m}_1}\upsi^{j_0+k_0+1}}{u}{2}{\ubar}{2}{\sphere}{2}{\region}\cdot\threenorms{\WangD^{\vec{m}_2}(\CAlphaAB,\FAlphaAB)}{u}{\infty}{\ubar}{\infty}{\sphere}{\infty}{\region}\\
						\notag&+\sum\limits_{\substack{\vec{m}_1+\vec{m}_2=(i,j-j_0,k-k_0)\\\abs{\vec{m}_2}\geq\frac{N-1}{2}}}\threenorms{\abs{u}^{j_0+k_0}\WangD^{\vec{m}_1}\upsi^{j_0+k_0+1}}{u}{2}{\ubar}{\infty}{\sphere}{\infty}{\region}\cdot\threenorms{\WangD^{\vec{m}_2}(\CAlphaAB,\FAlphaAB)}{u}{\infty}{\ubar}{2}{\sphere}{2}{\region}\\
						\notag\lesssim&~\left(\BAsize\frac{\delta a^{\frac{1}{2}}}{\abs{u}}\right)^{j_0+k_0}\left(\frac{\BAsize\delta^{\frac{1}{2}}a^{\frac{3}{4}}}{\abs{u}^{\frac{3}{2}}}+ \frac{\delta^{\frac{1}{2}}a^{\frac{1}{2}}}{\abs{u}^{\frac{1}{2}}}\right),
					\end{align}
					Combining \eqref{PV:AlphaBeta1}-\eqref{PV:AlphaBeta10}, summing \eqref{PV:AlphaBeta3} over $\abs{\vec{m}}\leq N-1$, we soak the $\sum\limits_{\abs{\vec{m}_4}\leq N-1}\threenorms{\WangD^{\vec{m}_4}\CAlphaAB}{u}{\infty}{\ubar}{2}{\sphere}{2}{\region}$ terms in \eqref{PV:AlphaBeta5} to the LHS of \eqref{PV:AlphaBeta3} since their coefficients are $\ll 1$. Then, we obtain the desired estimates \eqref{ES:topL2curvature1} by the initial data assumption \eqref{Data} and the choice of parameter $\delta,a,b$ in Section \ref{SS:Norms}, $\BAfluid,\BAsize$ in Section \eqref{SS:BA}.
				\end{proof}
				\begin{proof}[Proof of \eqref{ES:topL2curvature2}]
					Let $\vec{m}=(i,j,k)$ with $i+j+k\leq N-1$ with $k\geq 1$. Let $\WangD^{\vec{m}}$ be as in Definition \ref{DE:operators}. We commute equations \eqref{EQ:D3CBeta}, \eqref{EQ:D4Rho2}, \eqref{EQ:D4Rhodual2} with $\WangD^{\vec{m}}$ using \eqref{EQ:ang3ijk}-\eqref{EQ:ang4ijk}. Let $F:=\CBetaA$, $G:=\CRho+\CRhodual$ and
					\begin{align}
						\label{PV:BetaRho1}P_{(\CBetaA,\CRho+\CRhodual;\vec{m})}=&~\sum\limits_{\vec{m_1}+\vec{m_2}=(i-i_0,j,k-k_0+1)}\delta^{i_0}\abs{u}^{k_0}\WangD^{\vec{m_1}}\upsi^{k_0+1}\cdot\WangD^{\vec{m_2}}\psi^{i_0+1}\\
						\notag&+\sum\limits_{\vec{m_1}+\vec{m_2}+\vec{m_3}=(i-i_0,j,k-k_0)}\delta^{i_0}\abs{u}^{k_0}\WangD^{\vec{m_1}}\upsi^{k_0+1}\cdot\WangD^{\vec{m_2}}\psi^{i_0}\cdot\WangD^{\vec{m_3}}(\CBetaA,\FBetaA,\Rho,\Rhodual,\Scalar,\FRho,\FRhodual)\\
						\notag&+\sum\limits_{\vec{m_1}+\vec{m_2}+\vec{m_3}=(i-i_0,j,k-k_0)}\delta^{i_0}\abs{u}^{k_0}\WangD^{\vec{m_1}}\upsi^{k_0}\cdot\WangD^{\vec{m_2}}\psi^{i_0+1}\cdot\WangD^{\vec{m_3}}(\CBetabar,\FBetabar)\\
						\notag&+\sum\limits_{\vec{m_1}+\vec{m_2}+\vec{m_3}=(i-i_0,j,k-k_0)}\delta^{i_0}\abs{u}^{k_0}\WangD^{\vec{m_1}}\upsi^{k_0}\cdot\WangD^{\vec{m_2}}\psi^{i_0}\cdot\WangD^{\vec{m_3}}\angnabla(\Scalar,\FRho,\FRhodual)\\
						\notag&+\sum\limits_{\vec{m_1}+\vec{m_2}+\vec{m_3}+\vec{m_4}=(i-i_0-1,j,k-k_0+1)}\delta^{i_0+1}\abs{u}^{k_0-1}\WangD^{\vec{m_1}}\upsi^{k_0}\cdot\WangD^{\vec{m_2}}\psi^{i_0}\cdot\WangD^{\vec{m_3}}\upzeta\cdot\WangD^{\vec{m_4}}\CBetaA\\
						\notag&+\sum\limits_{\vec{m_1}+\vec{m_2}+\vec{m_3}+\vec{m_4}=(i-i_0-1,j,k-k_0)}\delta^{i_0+1}\abs{u}^{k_0}\WangD^{\vec{m_1}}\upsi^{k_0}\cdot\WangD^{\vec{m_2}}\psi^{i_0}\cdot\WangD^{\vec{m_3}}\Rhodual\cdot\WangD^{\vec{m_4}}\CBetaA\\
						\notag&+\sum\limits_{\vec{m_1}+\vec{m_2}+\vec{m_3}+\vec{m_4}=(i-i_0,j,k-k_0-1)}\delta^{i_0}\abs{u}^{k_0+1}\WangD^{\vec{m_1}}\upsi^{k_0}\cdot\WangD^{\vec{m_2}}\psi^{i_0}\cdot\WangD^{\vec{m_3}}(\CBetabar,\FBetabar)\cdot\WangD^{\vec{m_4}}\CBetaA,	
					\end{align}
					\begin{align}
						\label{PV:BetaRho2}Q_{(\CBetaA,\CRho+\CRhodual;\vec{m})}=&~\sum\limits_{\vec{m_1}+\vec{m_2}+\vec{m_3}=(i,j- j_0,k-k_0)}\abs{u}^{j_0+k_0}\WangD^{\vec{m_1}}\upsi^{j_0+k_0+1}\cdot\WangD^{\vec{m_2}}\psi\cdot\WangD^{\vec{m_3}}(\psi,\upsi)\\
						\notag&+\sum\limits_{\vec{m_1}+\vec{m_2}+\vec{m_3}=(i,j- j_0,k-k_0)}\abs{u}^{j_0+k_0}\WangD^{\vec{m_1}}\upsi^{j_0+k_0}\cdot\WangD^{\vec{m_2}}\psi\cdot\WangD^{\vec{m_3}}\angnabla\upsi\\
						\notag&+\sum\limits_{\vec{m_1}+\vec{m_2}+\vec{m_3}=(i,j- j_0,k-k_0)}\abs{u}^{j_0+k_0}\WangD^{\vec{m_1}}\upsi^{j_0+k_0}\cdot\WangD^{\vec{m_2}}\psi\cdot\WangD^{\vec{m_3}}(\Scalar,\FRho,\CRho,\CRhodual)\\
						\notag&+\sum\limits_{\vec{m_1}+\vec{m_2}+\vec{m_3}=(i,j- j_0,k-k_0)}\abs{u}^{j_0+k_0}\WangD^{\vec{m_1}}\upsi^{j_0+k_0+1}\cdot\WangD^{\vec{m_2}}(\CBetaA,\FBetaA)\\
						\notag&+\sum\limits_{\vec{m_1}+\vec{m_2}=(i,j- j_0,k-k_0)}\abs{u}^{j_0+k_0}\WangD^{\vec{m_1}}\upsi^{k_0+j_0}\cdot\WangD^{\vec{m_2}}\angnabla\FBetaA\\
						\notag&+\sum\limits_{\vec{m_1}+\vec{m_2}+\vec{m_3}=(i,j- j_0,k-k_0-1)}\abs{u}^{j_0+k_0+1}\WangD^{\vec{m_1}}\upsi^{j_0+k_0}\cdot\WangD^{\vec{m_3}}(\CBetaA,\FBetaA)\cdot\WangD^{\vec{m_3}}(\CRho,\CRhodual)\\
						\notag&+\sum\limits_{\vec{m_1}+\vec{m_2}+\vec{m_3}=(i,j- j_0-1,k-k_0+1)}\abs{u}^{j_0+k_0}\WangD^{\vec{m_1}}\upsi^{j_0+k_0}\cdot\WangD^{\vec{m_2}}\upzeta\cdot\WangD^{\vec{m_3}}(\CRho,\CRhodual)\\
						\notag&+\sum\limits_{\vec{m_1}+\vec{m_2}+\vec{m_3}=(i,j- j_0-1,k-k_0)}\abs{u}^{j_0+k_0+1}\WangD^{\vec{m_1}}\upsi^{j_0+k_0}\cdot\WangD^{\vec{m_2}}\Rhodual\cdot\angD^{\vec{m_3}}(\CRho,\CRhodual).
					\end{align}
					We apply Proposition \ref{PR:Int4} with $s(\CBetaA)=\frac{1}{2}$ to obtain
					\begin{align}\label{PV:BetaRho3}
						&\twonorms{\abs{u}\WangD^{\vec{m}}\CBetaA}{\ubar}{2}{\sphere}{2}{\hypersur_{u}}^2+\twonorms{\abs{u}\WangD^{\vec{m}}(\CRho,\CRhodual)}{u}{2}{\sphere}{2}{\hypersurbar_{\ubar}}^2\\
						\notag\lesssim&~\twonorms{\abs{u}\WangD^{\vec{m}}\CBetaA}{\ubar}{2}{\sphere}{2}{\hypersur_{\uinf}}^2+\twonorms{\abs{u}\WangD^{\vec{m}}(\CRho,\CRhodual)}{u}{2}{\sphere}{2}{\hypersurbar_{0}}^2\\
						\notag&+\int_{\region}\abs{u}^2\abs{\WangD^{\vec{m}}\CBetaA}\cdot \abs{P_{(\CBetaA,\CRho+\CRhodual;\vec{m})}}\diff\gvol+\int_{\region}\abs{u}^2\abs{\WangD^{\vec{m}}(\CRho,\CRhodual)}\cdot \abs{Q_{(\CBetaA,\CRho+\CRhodual;\vec{m})}}\diff\gvol.
					\end{align}
					We focus on estimating $\threenorms{\abs{u}P_{(\CBetaA,\CRho+\CRhodual;\vec{m})}}{u}{1}{\ubar}{2}{\sphere}{2}{\region}$ and $\threenorms{\abs{u}Q_{(\CBetaA,\CRho+\CRhodual;\vec{m})}}{u}{2}{\ubar}{2}{\sphere}{2}{\region}$ since 
					\begin{align}\label{PV:BetaRho4}
						&\int_{\region}\abs{u}^2\abs{\WangD^{\vec{m}}\CBetaA}\cdot \abs{P_{(\CBetaA,\CRho+\CRhodual;\vec{m})}}\diff\gvol+\int_{\region}\abs{u}^2\abs{\WangD^{\vec{m}}(\CRho,\CRhodual)}\cdot \abs{Q_{(\CBetaA,\CRho+\CRhodual;\vec{m})}}\diff\gvol\\
						\notag\lesssim&~\threenorms{\abs{u}\WangD^{\vec{m}}\CBetaA}{u}{\infty}{\ubar}{2}{\sphere}{2}{\region}\threenorms{\abs{u}P_{(\CBetaA,\CRho+\CRhodual;\vec{m})}}{u}{1}{\ubar}{2}{\sphere}{2}{\region}+\delta^{\frac{1}{2}}\threenorms{\abs{u}\WangD^{\vec{m}}(\CRho,\CRhodual)}{u}{2}{\ubar}{\infty}{\sphere}{2}{\region}\threenorms{\abs{u}Q_{(\CBetaA,\CRho+\CRhodual;\vec{m})}}{u}{2}{\ubar}{2}{\sphere}{2}{\region}.
					\end{align}
					For the terms involving $\CBetaA$, namely, the second term and the last three terms on the RHS of \eqref{PV:BetaRho1}, similar to \eqref{PV:AlphaBeta5}, using Proposition \ref{PR:EstimatesforRicciCurvature} for $\FBetabar$, bootstrap assumptions \eqref{BA} and Sobolev embedding $H^2\hookrightarrow L^\infty$ on $\suubar$ for $\WangD^{\vec{m}}\CBetaA$, we have:
					\begin{align}\label{PV:BetaRho5}
						&\sum\limits_{\vec{m}_1+\vec{m}_2+\vec{m}_3+\vec{m}_4=(i-i_0,j,k-k_0)}\delta^{i_0}\threenorms{\abs{u}^{k_0+1}\WangD^{\vec{m}_1}\psi^{i_0}\cdot\WangD^{\vec{m}_2}\upsi^{k_0}\cdot\left(\WangD^{\vec{m}_3}\upsi,\delta\WangD^{\vec{m}'_3}\Rhodual,\abs{u}\WangD^{\vec{m}'_3}\CBetabar,\abs{u}\WangD^{\vec{m}'_3}\FBetabar\right)\cdot\WangD^{\vec{m}_4}\CBetaA}{u}{1}{\ubar}{2}{\sphere}{2}{\region}\\
						\notag\lesssim&~\sum\limits_{\substack{\vec{m}_1+\vec{m}_2+\vec{m}_3+\vec{m}_4=(i-i_0,j,k-k_0)\\\abs{\vec{m}_4}<\frac{N-1}{2}}}\delta^{i_0}\threenorms{\abs{u}^{k_0}\WangD^{\vec{m}_1}\psi^{i_0}\cdot\WangD^{\vec{m}_2}\upsi^{k_0}\cdot\left(\WangD^{\vec{m}_3}\upsi,\delta\WangD^{\vec{m}'_3}\Rhodual,\abs{u}\WangD^{\vec{m}'_3}\CBetabar,\abs{u}\WangD^{\vec{m}'_3}\FBetabar\right)}{u}{1}{\ubar}{\infty}{\sphere}{2}{\region}\\
						\notag&\cdot\threenorms{\abs{u}\WangD^{\vec{m}_4}\CBetaA}{u}{\infty}{\ubar}{2}{\sphere}{\infty}{\region}\\
						\notag&+\sum\limits_{\substack{\vec{m}_1+\vec{m}_2+\vec{m}_3+\vec{m}_4=(i-i_0,j,k-k_0)\\\abs{\vec{m}_4}\geq\frac{N-1}{2}}}\delta^{i_0}\threenorms{\abs{u}^{k_0}\WangD^{\vec{m}_1}\psi^{i_0}\cdot\WangD^{\vec{m}_2}\upsi^{k_0}\cdot\left(\WangD^{\vec{m}_3}\upsi,\delta\WangD^{\vec{m}'_3}\Rhodual,\abs{u}\WangD^{\vec{m}'_3}\CBetabar,\abs{u}\WangD^{\vec{m}'_3}\FBetabar\right)}{u}{1}{\ubar}{\infty}{\sphere}{\infty}{\region}\\
						\notag&\cdot\threenorms{\abs{u}\WangD^{\vec{m}_4}\CBetaA}{u}{\infty}{\ubar}{2}{\sphere}{2}{\region}\\
						\notag\lesssim&~\left(\BAsize\frac{\delta a^{\frac{1}{2}}}{\abs{u}}\right)^{i_0+k_0}\left(\BAsize\frac{\delta a^{\frac{1}{2}}}{\abs{u}}+ b^{-\frac{1}{2}}\right)\cdot\sum\limits_{\abs{\vec{m}_4}\leq N-1}\threenorms{\abs{u}\WangD^{\vec{m}_4}\CBetaA}{u}{\infty}{\ubar}{2}{\sphere}{2}{\region},
					\end{align}
					where $\abs{\vec{m}_3'}=\abs{\vec{m}_3}-1$.
					Next, we estimate the rest terms on the RHS of \eqref{PV:BetaRho1}. For the first three terms, using bootstrap assumptions \eqref{BA} and Proposition \ref{PR:EstimatesforRicciCurvature}:
					\begin{align}\label{PV:BetaRho6}
						&\sum\limits_{\vec{m}_1+\vec{m}_2+\vec{m}_3+\vec{m}_4=(i-i_0,j,k-k_0)}\delta^{i_0}\\
						\notag&\cdot\threenorms{\abs{u}^{k_0+1}\WangD^{\vec{m}_1}\psi^{i_0}\cdot\WangD^{\vec{m}_2}\upsi^{k_0}\cdot\left(\WangD^{\vec{m}_3}\upsi\cdot\WangD^{\vec{m}_4}\left(\FBetaA,\CRho,\CRhodual,\Scalar,\FRho,\FRhodual\psi\right),\WangD^{\vec{m}_3}\psi\cdot\WangD^{\vec{m}_4}\left(\CBetabar,\FBetabar\right)\right)}{u}{1}{\ubar}{2}{\sphere}{2}{\region}\\
						\notag\lesssim&~\sum\limits_{\substack{\vec{m}_1+\vec{m}_2+\vec{m}_3=(i-i_0,j,k-k_0)\\\abs{\vec{m}_3}<\frac{N-1}{2}}}\delta^{i_0}\left(\threenorms{\abs{u}^{k_0+1}\WangD^{\vec{m}_1}\psi^{i_0}\cdot\WangD^{\vec{m}_2}\upsi^{k_0+1}}{u}{1}{\ubar}{2}{\sphere}{2}{\region}\cdot\threenorms{\WangD^{\vec{m}_3}\left(\FBetaA,\CRho,\CRhodual,\Scalar,\FRho,\FRhodual,\psi\right)}{u}{\infty}{\ubar}{\infty}{\sphere}{\infty}{\region}\right.\\
						\notag&\left.+\threenorms{\abs{u}^{k_0+1}\WangD^{\vec{m}_1}\psi^{i_0+1}\cdot\WangD^{\vec{m}_2}\upsi^{k_0}}{u}{1}{\ubar}{2}{\sphere}{2}{\region}\cdot\threenorms{\WangD^{\vec{m}_3}(\CBetabar,\FBetabar)}{u}{\infty}{\ubar}{\infty}{\sphere}{\infty}{\region}\right)\\
						\notag&+\sum\limits_{\substack{\vec{m}_1+\vec{m}_2+\vec{m}_3=(i-i_0,j,k-k_0)\\\abs{\vec{m}_3}\geq\frac{N-1}{2}}}\delta^{i_0}\left(\threenorms{\abs{u}^{k_0-1}\WangD^{\vec{m}_1}\psi^{i_0}\cdot\WangD^{\vec{m}_2}\upsi^{k_0+1}}{u}{2}{\ubar}{\infty}{\sphere}{\infty}{\region}\cdot\threenorms{\abs{u}^2\WangD^{\vec{m}_3}(\FBetaA,\CRho,\CRhodual,\Scalar,\FRho,\FRhodual,\psi)}{u}{\infty}{\ubar}{2}{\sphere}{2}{\region}\right.\\
						\notag&\left.+\threenorms{\abs{u}^{k_0-2}\WangD^{\vec{m}_1}\psi^{i_0+1}\cdot\WangD^{\vec{m}_2}\upsi^{k_0}}{u}{2}{\ubar}{\infty}{\sphere}{\infty}{\region}\cdot\threenorms{\abs{u}^3\WangD^{\vec{m}_3}(\CBetabar,\FBetabar)}{u}{\infty}{\ubar}{2}{\sphere}{2}{\region}\right)\\
						\notag\lesssim&~\left(\BAsize\frac{\delta a^{\frac{1}{2}}}{\abs{u}}\right)^{i_0+k_0}\cdot\left(\BAsize \delta^{\frac{1}{2}}a^{\frac{1}{2}}+\BAsize\frac{\delta^{\frac{3}{2}} a}{\abs{u}^{\frac{3}{2}}}\right).
					\end{align}
					For the fourth terms that depends on the highest order of fluids, using bootstrap assumptions \eqref{BA} and Proposition \ref{PR:EstimatesforRicciCurvature}:
					\begin{align}\label{PV:BetaRho8}
						&\sum\limits_{\vec{m}_1+\vec{m}_2+\vec{m}_3=(i-i_0,j,k-k_0)}\delta^{i_0}\threenorms{\abs{u}^{k_0+1}\WangD^{\vec{m}_1}\psi^{i_0}\cdot\WangD^{\vec{m}_2}\upsi^{k_0}\cdot\WangD^{\vec{m}_3}\angnabla(\Scalar,\FRho,\FRhodual)}{u}{1}{\ubar}{2}{\sphere}{2}{\region}\\
						\notag\lesssim&~\sum\limits_{\substack{\vec{m}_1+\vec{m}_2+\vec{m}_3=(i-i_0,j,k-k_0)\\\abs{\vec{m}_3}<\frac{N-1}{2}}}\delta^{i_0}\threenorms{\abs{u}^{k_0+1}\WangD^{\vec{m}_1}\psi^{i_0}\cdot\WangD^{\vec{m}_2}\upsi^{k_0}}{u}{1}{\ubar}{2}{\sphere}{2}{\region}\cdot\threenorms{\WangD^{\vec{m}_3}\angnabla(\Scalar,\FRho,\FRhodual)}{u}{\infty}{\ubar}{\infty}{\sphere}{\infty}{\region}\\
						\notag&+\sum\limits_{\substack{\vec{m}_1+\vec{m}_2+\vec{m}_3=(i-i_0,j,k-k_0)\\\abs{\vec{m}_3}\geq\frac{N-1}{2}}}\delta^{i_0}\threenorms{\abs{u}^{k_0+1}\WangD^{\vec{m}_1}\psi^{i_0}\cdot\WangD^{\vec{m}_2}\upsi^{k_0}}{u}{1}{\ubar}{\infty}{\sphere}{\infty}{\region}\cdot\threenorms{\WangD^{\vec{m}_3}\angnabla(\Scalar,\FRho,\FRhodual)}{u}{\infty}{\ubar}{2}{\sphere}{2}{\region}\\
						\notag\lesssim&~\left(\BAsize\frac{\delta a^{\frac{1}{2}}}{\abs{u}}\right)^{i_0+k_0}  b^{-\frac{1}{2}}.
					\end{align}
					Now, we estimate the terms on the RHS of \eqref{PV:BetaRho2}. For the below-top-order terms, using Proposition \ref{PR:EstimatesforRicciCurvature}, bootstrap assumptions \eqref{BA}, we have:
					\begin{align}\label{PV:BetaRho9}
						&\sum\limits_{\vec{m}_1+\vec{m}_2+\vec{m}_3=(i,j-j_0,k-k_0)}\threenorms{\abs{u}^{j_0+k_0+1}\WangD^{\vec{m}_1}\upsi^{j_0+k_0}\cdot\left(\WangD^{\vec{m}_2}(\psi,\upsi),\abs{u}\WangD^{\vec{m}'_2}(\CBetaA,\FBetaA,\Rhodual)\right)\cdot\WangD^{\vec{m}_3}(\CRho,\CRhodual,\Scalar,\FRho,\psi\cdot\upsi)}{u}{2}{\ubar}{2}{\sphere}{2}{\region}\\
						\notag\lesssim&~\sum\limits_{\substack{\vec{m}_1+\vec{m}_2+\vec{m}_3=(i,j-j_0,k-k_0)\\\abs{\vec{m}_3}<\frac{N-1}{2}}}\threenorms{\abs{u}^{j_0+k_0+1}\WangD^{\vec{m}_1}\upsi^{j_0+k_0}\cdot\left(\WangD^{\vec{m}_2}(\psi,\upsi),\abs{u}\WangD^{\vec{m}_2}(\CBetaA,\FBetaA,\Rhodual)\right)}{u}{2}{\ubar}{2}{\sphere}{2}{\region}\\
						\notag&\cdot\threenorms{\WangD^{\vec{m}_3}(\CRho,\CRhodual,\Scalar,\FRho,\psi\cdot\upsi)}{u}{\infty}{\ubar}{\infty}{\sphere}{\infty}{\region}\\
						\notag&+\sum\limits_{\substack{\vec{m}_1+\vec{m}_2+\vec{m}_3=(i,j-j_0,k-k_0)\\\abs{\vec{m}_3}\geq\frac{N-1}{2}}}\threenorms{\abs{u}^{j_0+k_0-1}\WangD^{\vec{m}_1}\upsi^{j_0+k_0}\cdot\left(\WangD^{\vec{m}_2}(\psi,\upsi),\abs{u}\WangD^{\vec{m}'_2}(\CBetaA,\FBetaA,\Rhodual)\right)}{u}{2}{\ubar}{\infty}{\sphere}{\infty}{\region}\\
						\notag&\cdot\threenorms{\abs{u}^2\WangD^{\vec{m}_3}(\CRho,\CRhodual,\Scalar,\FRho,\psi\cdot\upsi)}{u}{\infty}{\ubar}{2}{\sphere}{2}{\region}\\
						\notag\lesssim&~\left(\BAsize\frac{\delta a^{\frac{1}{2}}}{\abs{u}}\right)^{j_0+k_0}\left(\BAsize\frac{\delta a}{\abs{u}^{\frac{3}{2}}}+ \frac{a^{\frac{1}{2}}b^{-\frac{1}{2}}}{\abs{u}^{\frac{1}{2}}}+ \delta^{-\frac{1}{2}}b^{-\frac{1}{2}}\right),
					\end{align}
					where $\abs{\vec{m}_2'}=\abs{\vec{m}_2}-1$.
					Finally, we estimate the top-order terms, namely, the second and fifth terms on the RHS of \eqref{PV:BetaRho2} using Proposition \ref{PR:EstimatesforRicciCurvature}, bootstrap assumptions \eqref{BA}:
					\begin{align}\label{PV:BetaRho10}
						&\sum\limits_{\vec{m}_1+\vec{m}_2=(i,j-j_0,k-k_0)}\threenorms{\abs{u}^{j_0+k_0+1}\WangD^{\vec{m}_1}\upsi^{j_0+k_0}\cdot\WangD^{\vec{m}_2}(\psi\cdot\angnabla\upsi,\angnabla\FBetaA)}{u}{2}{\ubar}{2}{\sphere}{2}{\region}\\
						\notag\lesssim&~\sum\limits_{\substack{\vec{m}_1+\vec{m}_2=(i,j-j_0,k-k_0)\\\abs{\vec{m}_2}<\frac{N-1}{2}}}\threenorms{\abs{u}^{j_0+k_0+1}\WangD^{\vec{m}_1}\upsi^{j_0+k_0}}{u}{2}{\ubar}{2}{\sphere}{2}{\region}\cdot\threenorms{\WangD^{\vec{m}_2}(\psi\cdot\angnabla\upsi,\angnabla\FBetaA)}{u}{\infty}{\ubar}{\infty}{\sphere}{\infty}{\region}\\
						\notag&+\sum\limits_{\substack{\vec{m}_1+\vec{m}_2=(i,j-j_0,k-k_0)\\\abs{\vec{m}_2}\geq\frac{N-1}{2}}}\threenorms{\abs{u}^{j_0+k_0+1}\WangD^{\vec{m}_1}\upsi^{j_0+k_0}}{u}{2}{\ubar}{\infty}{\sphere}{\infty}{\region}\cdot\threenorms{\WangD^{\vec{m}_2}(\psi\cdot\angnabla\upsi,\angnabla\FBetaA)}{u}{\infty}{\ubar}{2}{\sphere}{2}{\region}\\
						\notag\lesssim&~\left(\BAsize\frac{\delta a^{\frac{1}{2}}}{\abs{u}}\right)^{j_0+k_0}\left(\frac{\BAsize\delta a}{\abs{u}^{\frac{5}{2}}}+ \right),
					\end{align}
					Combining \eqref{PV:BetaRho1}-\eqref{PV:BetaRho10}, summing \eqref{PV:BetaRho3} over $\abs{\vec{m}}\leq N-1$, we soak the $\sum\limits_{\abs{\vec{m}_4}\leq N-1}\threenorms{\WangD^{\vec{m}_4}\CBetaA}{u}{\infty}{\ubar}{2}{\sphere}{2}{\region}$ terms in \eqref{PV:BetaRho5} to the LHS of \eqref{PV:BetaRho3} since their coefficients are $\ll 1$. Then, we obtain the desired estimates \eqref{ES:topL2curvature2} by the initial data assumption \eqref{Data} and the choice of parameter $\delta,a,b$ in Section \ref{SS:Norms}, $\BAfluid,\BAsize$ in Section \eqref{SS:BA}.
				\end{proof}
				\begin{proof}[Proof of \eqref{ES:topL2curvature3}]
					Let $\vec{m}=(i,j,k)$ with $i+j+k\leq N-1$ with $k\geq 1$. Let $\WangD^{\vec{m}}$ be as in Definition \ref{DE:operators}. We commute equations \eqref{EQ:D4CBetabar}, \eqref{EQ:D3Rho2}, \eqref{EQ:D3Rhodual2} with $\WangD^{\vec{m}}$ using \eqref{EQ:ang3ijk}-\eqref{EQ:ang4ijk}. Let $F:=\CRho+\CRhodual$, $G:=\CBetabar$, and
					\begin{align}
						\label{PV:RhoBetabar1}P_{(\CRho+\CRhodual,\CBetabar;\vec{m})}=&~\sum\limits_{\vec{m_1}+\vec{m_2}=(i-i_0,j,k-k_0)}\delta^{i_0}\abs{u}^{k_0}\WangD^{\vec{m_1}}\upsi^{k_0+3}\cdot\WangD^{\vec{m_2}}\psi^{i_0}\\
						\notag&+\sum\limits_{\vec{m_1}+\vec{m_2}+\vec{m_3}=(i-i_0,j,k-k_0)}\delta^{i_0}\abs{u}^{k_0}\WangD^{\vec{m_1}}\upsi^{k_0+1}\cdot\WangD^{\vec{m_2}}\psi^{i_0}\cdot\WangD^{\vec{m_3}}(\CBetabar,\FBetabar,\Scalar,\FRho)\\
						\notag&+\sum\limits_{\vec{m_1}+\vec{m_2}+\vec{m_3}=(i-i_0,j,k-k_0)}\delta^{i_0}\abs{u}^{k_0}\WangD^{\vec{m_1}}\upsi^{k_0+1}\cdot\WangD^{\vec{m_2}}\psi^{i_0}\cdot\WangD^{\vec{m_3}}\angnabla\upsi\\
						\notag&+\sum\limits_{\vec{m_1}+\vec{m_2}+\vec{m_3}=(i-i_0,j,k-k_0)}\delta^{i_0}\abs{u}^{k_0}\WangD^{\vec{m_1}}\upsi^{k_0}\cdot\WangD^{\vec{m_2}}\psi^{i_0}\cdot\WangD^{\vec{m_3}}\angnabla\FBetabar\\
						\notag&+\sum\limits_{\vec{m_1}+\vec{m_2}+\vec{m_3}=(i-i_0,j,k-k_0)}\delta^{i_0}\abs{u}^{k_0}\WangD^{\vec{m_1}}\upsi^{k_0+1}\cdot\WangD^{\vec{m_2}}\psi^{i_0}\cdot\WangD^{\vec{m_3}}(\CRho,\CRhodual)\\
						\notag&+\sum\limits_{\vec{m_1}+\vec{m_2}+\vec{m_3}+\vec{m_4}=(i-i_0-1,j,k-k_0+1)}\delta^{i_0+1}\abs{u}^{k_0-1}\WangD^{\vec{m_1}}\upsi^{k_0}\cdot\WangD^{\vec{m_2}}\psi^{i_0}\cdot\WangD^{\vec{m_3}}\upzeta\cdot\WangD^{\vec{m_4}}(\CRho,\CRhodual)\\
						\notag&+\sum\limits_{\vec{m_1}+\vec{m_2}+\vec{m_3}+\vec{m_4}=(i-i_0-1,j,k-k_0)}\delta^{i_0+1}\abs{u}^{k_0}\WangD^{\vec{m_1}}\upsi^{k_0}\cdot\WangD^{\vec{m_2}}\psi^{i_0}\cdot\WangD^{\vec{m_3}}\Rhodual\cdot\WangD^{\vec{m_4}}(\CRho,\CRhodual)\\
						\notag&+\sum\limits_{\vec{m_1}+\vec{m_2}+\vec{m_3}+\vec{m_4}=(i-i_0,j,k-k_0-1)}\delta^{i_0}\abs{u}^{k_0+1}\WangD^{\vec{m_1}}\upsi^{k_0}\cdot\WangD^{\vec{m_2}}\psi^{i_0}\cdot\WangD^{\vec{m_3}}(\CBetabar,\FBetabar)\cdot\WangD^{\vec{m_4}}(\CRho,\CRhodual),	
					\end{align}
					\begin{align}
						\label{PV:RhoBetabar2}Q_{(\CRho+\CRhodual,\CBetabar;\vec{m})}=&~\sum\limits_{\vec{m_1}+\vec{m_2}=(i,j- j_0,k-k_0+1)}\abs{u}^{j_0+k_0}\WangD^{\vec{m_1}}\upsi^{j_0+k_0+1}\cdot\WangD^{\vec{m_2}}\psi\\
						\notag&+\sum\limits_{\vec{m_1}+\vec{m_2}+\vec{m_3}=(i,j- j_0,k-k_0)}\abs{u}^{j_0+k_0}\WangD^{\vec{m_1}}\upsi^{j_0+k_0}\cdot\WangD^{\vec{m_2}}\psi\cdot\WangD^{\vec{m_3}}(\CBetabar,\FBetabar)\\
						\notag&+\sum\limits_{\vec{m_1}+\vec{m_2}=(i,j- j_0,k-k_0)}\abs{u}^{j_0+k_0}\WangD^{\vec{m_1}}\upsi^{j_0+k_0+1}\cdot\WangD^{\vec{m_2}}(\CBetaA,\FBetaA,\Scalar,\FRho,\CRho,\CRhodual)\\
						\notag&+\sum\limits_{\vec{m_1}+\vec{m_2}=(i,j- j_0,k-k_0)}\abs{u}^{j_0+k_0}\WangD^{\vec{m_1}}\upsi^{k_0+j_0}\cdot\WangD^{\vec{m_2}}\angnabla(\Scalar,\FRho,\FRhodual)\\
						\notag&+\sum\limits_{\vec{m_1}+\vec{m_2}+\vec{m_3}=(i,j- j_0,k-k_0-1)}\abs{u}^{j_0+k_0+1}\WangD^{\vec{m_1}}\upsi^{j_0+k_0}\cdot\WangD^{\vec{m_3}}(\CBetaA,\FBetaA)\cdot\WangD^{\vec{m_3}}\CBetabar\\
						\notag&+\sum\limits_{\vec{m_1}+\vec{m_2}+\vec{m_3}=(i,j- j_0-1,k-k_0+1)}\abs{u}^{j_0+k_0}\WangD^{\vec{m_1}}\upsi^{j_0+k_0}\cdot\WangD^{\vec{m_2}}\upzeta\cdot\WangD^{\vec{m_3}}\CBetabar\\
						\notag&+\sum\limits_{\vec{m_1}+\vec{m_2}+\vec{m_3}=(i,j- j_0-1,k-k_0)}\abs{u}^{j_0+k_0+1}\WangD^{\vec{m_1}}\upsi^{j_0+k_0}\cdot\WangD^{\vec{m_2}}\Rhodual\cdot\angD^{\vec{m_3}}\CBetabar.
					\end{align}
					We apply Proposition \ref{PR:Int4} with $s(\CRho+\CRhodual)=1$ to obtain
					\begin{align}\label{PV:RhoBetabar3}
						&\twonorms{\abs{u}^2\WangD^{\vec{m}}(\CRho,\CRhodual)}{\ubar}{2}{\sphere}{2}{\hypersur_{u}}^2+\twonorms{\abs{u}^2\WangD^{\vec{m}}\CBetabar}{u}{2}{\sphere}{2}{\hypersurbar_{\ubar}}^2\\
						\notag\lesssim&~\twonorms{\abs{u}^2\WangD^{\vec{m}}(\CRho,\CRhodual)}{\ubar}{2}{\sphere}{2}{\hypersur_{\uinf}}^2+\twonorms{\abs{u}^2\WangD^{\vec{m}}\CBetabar}{u}{2}{\sphere}{2}{\hypersurbar_{0}}^2\\
						\notag&+\int_{\region}\abs{u}^4\abs{\WangD^{\vec{m}}(\CRho,\CRhodual)}\cdot \abs{P_{(\CRho+\CRhodual,\CBetabar;\vec{m})}}\diff\gvol+\int_{\region}\abs{u}^4\abs{\WangD^{\vec{m}}\CBetabar}\cdot \abs{Q_{(\CRho+\CRhodual,\CBetabar;\vec{m})}}\diff\gvol.
					\end{align}
					We focus on estimating $\threenorms{\abs{u}^2P_{(\CRho+\CRhodual,\CBetabar;\vec{m})}}{u}{1}{\ubar}{2}{\sphere}{2}{\region}$ and $\threenorms{\abs{u}^2Q_{(\CRho+\CRhodual,\CBetabar;\vec{m})}}{u}{2}{\ubar}{2}{\sphere}{2}{\region}$ since 
					\begin{align}\label{PV:RhoBetabar4}
						&\int_{\region}\abs{u}^4\abs{\WangD^{\vec{m}}(\CRho,\CRhodual)}\cdot \abs{P_{(\CRho+\CRhodual,\CBetabar;\vec{m})}}\diff\gvol+\int_{\region}\abs{u}^4\abs{\WangD^{\vec{m}}\CBetabar}\cdot \abs{Q_{(\CRho+\CRhodual,\CBetabar;\vec{m})}}\diff\gvol\\
						\notag\lesssim&~\threenorms{\abs{u}^2\WangD^{\vec{m}}(\CRho,\CRhodual)}{u}{\infty}{\ubar}{2}{\sphere}{2}{\region}\threenorms{\abs{u}^2P_{(\CRho+\CRhodual,\CBetabar;\vec{m})}}{u}{1}{\ubar}{2}{\sphere}{2}{\region}+\delta^{\frac{1}{2}}\threenorms{\abs{u}^2\WangD^{\vec{m}}\CBetabar}{u}{2}{\ubar}{\infty}{\sphere}{2}{\region}\threenorms{\abs{u}^2Q_{(\CRho+\CRhodual,\CBetabar;\vec{m})}}{u}{2}{\ubar}{2}{\sphere}{2}{\region}.
					\end{align}
					For the terms involving $(\CRho,\CRhodual)$, namely, the last four terms on the RHS of \eqref{PV:RhoBetabar1}, using Proposition \ref{PR:EstimatesforRicciCurvature} for $\FBetabar$, bootstrap assumptions \eqref{BA} and Sobolev embedding $H^2\hookrightarrow L^\infty$ on $\suubar$ for $\WangD^{\vec{m}}(\CRho,\CRhodual)$, we have:
					\begin{align}\label{PV:RhoBetabar5}
						&\sum\limits_{\vec{m}_1+\vec{m}_2+\vec{m}_3+\vec{m}_4=(i-i_0,j,k-k_0)}\delta^{i_0}\threenorms{\abs{u}^{k_0+2}\WangD^{\vec{m}_1}\psi^{i_0}\cdot\WangD^{\vec{m}_2}\upsi^{k_0}\cdot\left(\WangD^{\vec{m}_3}\upsi,\delta\WangD^{\vec{m}'_3}\Rhodual,\abs{u}\WangD^{\vec{m}'_3}\CBetabar,\abs{u}\WangD^{\vec{m}'_3}\FBetabar\right)\cdot\WangD^{\vec{m}_4}(\CRho,\CRhodual)}{u}{1}{\ubar}{2}{\sphere}{2}{\region}\\
						\notag\lesssim&~\sum\limits_{\substack{\vec{m}_1+\vec{m}_2+\vec{m}_3+\vec{m}_4=(i-i_0,j,k-k_0)\\\abs{\vec{m}_4}<\frac{N-1}{2}}}\delta^{i_0}\threenorms{\abs{u}^{k_0}\WangD^{\vec{m}_1}\psi^{i_0}\cdot\WangD^{\vec{m}_2}\upsi^{k_0}\cdot\left(\WangD^{\vec{m}_3}\upsi,\delta\WangD^{\vec{m}'_3}\Rhodual,\abs{u}\WangD^{\vec{m}'_3}\CBetabar,\abs{u}\WangD^{\vec{m}'_3}\FBetabar\right)}{u}{1}{\ubar}{\infty}{\sphere}{2}{\region}\\
						\notag&\cdot\threenorms{\abs{u}^2\WangD^{\vec{m}_4}(\CRho,\CRhodual)}{u}{\infty}{\ubar}{2}{\sphere}{\infty}{\region}\\
						\notag&+\sum\limits_{\substack{\vec{m}_1+\vec{m}_2+\vec{m}_3+\vec{m}_4=(i-i_0,j,k-k_0)\\\abs{\vec{m}_4}\geq\frac{N-1}{2}}}\delta^{i_0}\threenorms{\abs{u}^{k_0}\WangD^{\vec{m}_1}\psi^{i_0}\cdot\WangD^{\vec{m}_2}\upsi^{k_0}\cdot\left(\WangD^{\vec{m}_3}\upsi,\delta\WangD^{\vec{m}'_3}\Rhodual,\abs{u}\WangD^{\vec{m}'_3}\CBetabar,\abs{u}\WangD^{\vec{m}'_3}\FBetabar\right)}{u}{1}{\ubar}{\infty}{\sphere}{\infty}{\region}\\
						\notag&\cdot\threenorms{\abs{u}^2\WangD^{\vec{m}_4}(\CRho,\CRhodual)}{u}{\infty}{\ubar}{2}{\sphere}{2}{\region}\\
						\notag\lesssim&~\left(\BAsize\frac{\delta a^{\frac{1}{2}}}{\abs{u}}\right)^{i_0+k_0}\left(\BAsize\frac{\delta a^{\frac{1}{2}}}{\abs{u}}+ b^{-\frac{1}{2}}\right)\cdot\sum\limits_{\abs{\vec{m}_4}\leq N-1}\threenorms{\abs{u}^2\WangD^{\vec{m}_4}(\CRho,\CRhodual)}{u}{\infty}{\ubar}{2}{\sphere}{2}{\region},
					\end{align}
					where $\abs{\vec{m}_3'}=\abs{\vec{m}_3}-1$.
					Next, we estimate the rest terms on the RHS of \eqref{PV:RhoBetabar1}. For the below-top-order terms, namely, the first two terms, using bootstrap assumptions \eqref{BA} and Proposition \ref{PR:EstimatesforRicciCurvature}:
					\begin{align}\label{PV:RhoBetabar6}
						&\sum\limits_{\vec{m}_1+\vec{m}_2+\vec{m}_3=(i-i_0,j,k-k_0)}\delta^{i_0}\cdot\threenorms{\abs{u}^{k_0+2}\WangD^{\vec{m}_1}\psi^{i_0}\cdot\WangD^{\vec{m}_2}\upsi^{k_0+1}\cdot\left(\WangD^{\vec{m}_3}\upsi^2,\WangD^{\vec{m}_3}\left(\CBetabar,\FBetabar,\Scalar,\FRho\right)\right)}{u}{1}{\ubar}{2}{\sphere}{2}{\region}\\
						\notag\lesssim&~\sum\limits_{\substack{\vec{m}_1+\vec{m}_2+\vec{m}_3=(i-i_0,j,k-k_0)\\\abs{\vec{m}_3}<\frac{N-1}{2}}}\delta^{i_0}\left(\threenorms{\abs{u}^{k_0+2}\WangD^{\vec{m}_1}\psi^{i_0}\cdot\WangD^{\vec{m}_2}\upsi^{k_0+1}}{u}{1}{\ubar}{2}{\sphere}{2}{\region}\cdot\threenorms{\WangD^{\vec{m}_3}\upsi^2}{u}{\infty}{\ubar}{\infty}{\sphere}{\infty}{\region}\right.\\
						\notag&\left.+\threenorms{\abs{u}^{k_0+2}\WangD^{\vec{m}_1}\psi^{i_0}\cdot\WangD^{\vec{m}_2}\upsi^{k_0+1}}{u}{1}{\ubar}{2}{\sphere}{2}{\region}\cdot\threenorms{\WangD^{\vec{m}_3}\left(\CBetabar,\FBetabar,\Scalar,\FRho\right)}{u}{\infty}{\ubar}{\infty}{\sphere}{\infty}{\region}\right)\\
						\notag&+\sum\limits_{\substack{\vec{m}_1+\vec{m}_2+\vec{m}_3=(i-i_0,j,k-k_0)\\\abs{\vec{m}_3}\geq\frac{N-1}{2}}}\delta^{i_0}\left(\threenorms{\abs{u}^{k_0-1}\WangD^{\vec{m}_1}\psi^{i_0}\cdot\WangD^{\vec{m}_2}\upsi^{k_0+1}}{u}{2}{\ubar}{\infty}{\sphere}{\infty}{\region}\cdot\threenorms{\abs{u}^2\WangD^{\vec{m}_3}\upsi^2}{u}{\infty}{\ubar}{2}{\sphere}{2}{\region}\right.\\
						\notag&\left.+\threenorms{\abs{u}^{k_0-1}\WangD^{\vec{m}_1}\psi^{i_0}\cdot\WangD^{\vec{m}_2}\upsi^{k_0+1}}{u}{2}{\ubar}{\infty}{\sphere}{\infty}{\region}\cdot\threenorms{\abs{u}^3\WangD^{\vec{m}_3}\left(\CBetabar,\FBetabar,\Scalar,\FRho\right)}{u}{\infty}{\ubar}{2}{\sphere}{2}{\region}\right)\\
						\notag\lesssim&~\left(\BAsize\frac{\delta a^{\frac{1}{2}}}{\abs{u}}\right)^{i_0+k_0}\cdot\left(\BAsize \delta^{\frac{3}{2}} a^{\frac{1}{2}}+\BAsize\frac{\delta^{2} a}{\abs{u}^{\frac{5}{2}}}\right).
					\end{align}
					For the top-order terms, using bootstrap assumptions \eqref{BA} and Proposition \ref{PR:EstimatesforRicciCurvature}:
					\begin{align}\label{PV:RhoBetabar8}
						&\sum\limits_{\vec{m}_1+\vec{m}_2+\vec{m}_3=(i-i_0,j,k-k_0)}\delta^{i_0}\threenorms{\abs{u}^{k_0+2}\WangD^{\vec{m}_1}\psi^{i_0}\cdot\WangD^{\vec{m}_2}\upsi^{k_0}\cdot\WangD^{\vec{m}_3}(\upsi\cdot\angnabla\upsi,\FBetabar)}{u}{1}{\ubar}{2}{\sphere}{2}{\region}\\
						\notag\lesssim&~\sum\limits_{\substack{\vec{m}_1+\vec{m}_2+\vec{m}_3=(i-i_0,j,k-k_0)\\\abs{\vec{m}_3}<\frac{N-1}{2}}}\delta^{i_0}\threenorms{\abs{u}^{k_0+2}\WangD^{\vec{m}_1}\psi^{i_0}\cdot\WangD^{\vec{m}_2}\upsi^{k_0}}{u}{1}{\ubar}{2}{\sphere}{2}{\region}\cdot\threenorms{\WangD^{\vec{m}_3}(\upsi\cdot\angnabla\upsi,\FBetabar)}{u}{\infty}{\ubar}{\infty}{\sphere}{\infty}{\region}\\
						\notag&+\sum\limits_{\substack{\vec{m}_1+\vec{m}_2+\vec{m}_3=(i-i_0,j,k-k_0)\\\abs{\vec{m}_3}\geq\frac{N-1}{2}}}\delta^{i_0}\threenorms{\abs{u}^{k_0+2}\WangD^{\vec{m}_1}\psi^{i_0}\cdot\WangD^{\vec{m}_2}\upsi^{k_0}}{u}{1}{\ubar}{\infty}{\sphere}{\infty}{\region}\cdot\threenorms{\WangD^{\vec{m}_3}(\upsi\cdot\angnabla\upsi,\FBetabar)}{u}{\infty}{\ubar}{2}{\sphere}{2}{\region}\\
						\notag\lesssim&~\left(\BAsize\frac{\delta a^{\frac{1}{2}}}{\abs{u}}\right)^{i_0+k_0}\left(\BAsize\frac{\delta^2a}{\abs{u}^{\frac{3}{2}}}+  b^{-\frac{1}{2}}\right).
					\end{align}
					Now, we estimate the terms on the RHS of \eqref{PV:RhoBetabar2}. For the below-top-order terms, using Proposition \ref{PR:EstimatesforRicciCurvature}, bootstrap assumptions \eqref{BA}, we have:
					\begin{align}\label{PV:RhoBetabar8.5}
						&\sum\limits_{\vec{m}_1+\vec{m}_2=(i,j-j_0,k-k_0)}\threenorms{\abs{u}^{j_0+k_0+2}\WangD^{\vec{m}_1}\upsi^{j_0+k_0+1}\cdot\WangD^{\vec{m}_2}(\CRho,\CRhodual,\Scalar,\FRho,\psi,\CBetaA,\FBetaA)}{u}{2}{\ubar}{2}{\sphere}{2}{\region}\\
						\notag\lesssim&~\sum\limits_{\substack{\vec{m}_1+\vec{m}_2=(i,j-j_0,k-k_0)\\\abs{\vec{m}_2}<\frac{N-1}{2}}}\threenorms{\abs{u}^{j_0+k_0+2}\WangD^{\vec{m}_1}\upsi^{j_0+k_0+1}}{u}{2}{\ubar}{2}{\sphere}{2}{\region}\cdot\threenorms{\WangD^{\vec{m}_2}(\CRho,\CRhodual,\Scalar,\FRho,\psi,\CBetaA,\FBetaA)}{u}{\infty}{\ubar}{\infty}{\sphere}{\infty}{\region}\\
						\notag&+\sum\limits_{\substack{\vec{m}_1+\vec{m}_2=(i,j-j_0,k-k_0)\\\abs{\vec{m}_2}\geq\frac{N-1}{2}}}\threenorms{\abs{u}^{j_0+k_0}\WangD^{\vec{m}_1}\upsi^{j_0+k_0+1}}{u}{2}{\ubar}{\infty}{\sphere}{\infty}{\region}\cdot\threenorms{\abs{u}^2\WangD^{\vec{m}_2}(\CRho,\CRhodual,\Scalar,\FRho,\psi,\CBetaA,\FBetaA)}{u}{\infty}{\ubar}{2}{\sphere}{2}{\region}\\
						\notag\lesssim&~\left(\BAsize\frac{\delta a^{\frac{1}{2}}}{\abs{u}}\right)^{j_0+k_0}\left(\BAsize\frac{\delta a}{\abs{u}^{\frac{1}{2}}}+ \delta^{\frac{1}{2}}a^{\frac{1}{2}}\right),
					\end{align}
					\begin{align}\label{PV:RhoBetabar9}
						&\sum\limits_{\vec{m}_1+\vec{m}_2+\vec{m}_3=(i,j-j_0,k-k_0)}\threenorms{\abs{u}^{j_0+k_0+2}\WangD^{\vec{m}_1}\upsi^{j_0+k_0}\cdot(\WangD^{\vec{m}_2}\psi,\abs{u}\WangD^{\vec{m}'_2}(\CBetaA,\FBetaA,\Rhodual))\cdot\WangD^{\vec{m}_3}(\CBetabar,\FBetabar)}{u}{2}{\ubar}{2}{\sphere}{2}{\region}\\
						\notag\lesssim&~\sum\limits_{\substack{\vec{m}_1+\vec{m}_2+\vec{m}_3=(i,j-j_0,k-k_0)\\\abs{\vec{m}_3}<\frac{N-1}{2}}}\threenorms{\abs{u}^{j_0+k_0+2}\WangD^{\vec{m}_1}\upsi^{j_0+k_0}\cdot(\WangD^{\vec{m}_2}\psi,\abs{u}\WangD^{\vec{m}'_2}(\CBetaA,\FBetaA,\Rhodual))}{u}{2}{\ubar}{2}{\sphere}{2}{\region}\\
						\notag&\cdot\threenorms{\WangD^{\vec{m}_3}(\CBetabar,\FBetabar)}{u}{\infty}{\ubar}{\infty}{\sphere}{\infty}{\region}\\
						\notag&+\sum\limits_{\substack{\vec{m}_1+\vec{m}_2+\vec{m}_3=(i,j-j_0,k-k_0)\\\abs{\vec{m}_3}\geq\frac{N-1}{2}}}\threenorms{\abs{u}^{j_0+k_0-1}\WangD^{\vec{m}_1}\upsi^{j_0+k_0}\cdot(\WangD^{\vec{m}_2}\psi,\abs{u}\WangD^{\vec{m}'_2}(\CBetaA,\FBetaA,\Rhodual))}{u}{2}{\ubar}{\infty}{\sphere}{\infty}{\region}\\
						\notag&\cdot\threenorms{\abs{u}^3\WangD^{\vec{m}_3}(\CBetabar,\FBetabar)}{u}{\infty}{\ubar}{2}{\sphere}{2}{\region}\\
						\notag\lesssim&~\left(\BAsize\frac{\delta a^{\frac{1}{2}}}{\abs{u}}\right)^{j_0+k_0}\cdot\left(\BAsize\delta a+ a^{\frac{1}{2}}b^{-\frac{1}{2}}+ \delta^{-\frac{1}{2}}b^{-\frac{1}{2}}\right),
					\end{align}
					where $\abs{\vec{m}_2'}=\abs{\vec{m}_2}-1$.
					Finally, we estimate the top-order terms, namely, the fourth terms on the RHS of \eqref{PV:RhoBetabar2} using Proposition \ref{PR:EstimatesforRicciCurvature}, bootstrap assumptions \eqref{BA}:
					\begin{align}\label{PV:RhoBetabar10}
						&\sum\limits_{\vec{m}_1+\vec{m}_2=(i,j-j_0,k-k_0)}\threenorms{\abs{u}^{j_0+k_0+2}\WangD^{\vec{m}_1}\upsi^{j_0+k_0}\cdot\WangD^{\vec{m}_2}\angnabla(\Scalar,\FRho,\FRhodual)}{u}{2}{\ubar}{2}{\sphere}{2}{\region}\\
						\notag\lesssim&~\sum\limits_{\substack{\vec{m}_1+\vec{m}_2=(i,j-j_0,k-k_0)\\\abs{\vec{m}_2}<\frac{N-1}{2}}}\threenorms{\abs{u}^{j_0+k_0+2}\WangD^{\vec{m}_1}\upsi^{j_0+k_0}}{u}{2}{\ubar}{2}{\sphere}{2}{\region}\cdot\threenorms{\WangD^{\vec{m}_2}\angnabla(\Scalar,\FRho,\FRhodual)}{u}{\infty}{\ubar}{\infty}{\sphere}{\infty}{\region}\\
						\notag&+\sum\limits_{\substack{\vec{m}_1+\vec{m}_2=(i,j-j_0,k-k_0)\\\abs{\vec{m}_2}\geq\frac{N-1}{2}}}\threenorms{\abs{u}^{j_0+k_0+2}\WangD^{\vec{m}_1}\upsi^{j_0+k_0}}{u}{2}{\ubar}{\infty}{\sphere}{\infty}{\region}\cdot\threenorms{\WangD^{\vec{m}_2}\angnabla(\Scalar,\FRho,\FRhodual)}{u}{\infty}{\ubar}{2}{\sphere}{2}{\region}\\
						\notag\lesssim&~\left(\BAsize\frac{\delta a^{\frac{1}{2}}}{\abs{u}}\right)^{j_0+k_0}\cdot b^{-\frac{1}{2}},
					\end{align}
					Combining \eqref{PV:RhoBetabar1}-\eqref{PV:RhoBetabar10}, summing \eqref{PV:RhoBetabar3} over $\abs{\vec{m}}\leq N-1$, we soak the $\sum\limits_{\abs{\vec{m}_4}\leq N-1}\threenorms{\WangD^{\vec{m}_4}(\CRho,\CRhodual)}{u}{\infty}{\ubar}{2}{\sphere}{2}{\region}$ terms in \eqref{PV:RhoBetabar5} to the LHS of \eqref{PV:RhoBetabar3} since their coefficients are $\ll 1$. Then, we obtain the desired estimates \eqref{ES:topL2curvature3} by the initial data assumption \eqref{Data} and the choice of parameter $\delta,a,b$ in Section \ref{SS:Norms}, $\BAfluid,\BAsize$ in Section \eqref{SS:BA}.
				\end{proof}
				\begin{proof}[Proof of \eqref{ES:topL2curvature4}]
					Let $\vec{m}=(i,j,k)$ with $i+j+k\leq N-1$ with $k\geq 1$. Let $\WangD^{\vec{m}}$ be as in Definition \ref{DE:operators}. We commute equations \eqref{EQ:D4CAlphabar}, \eqref{EQ:D3CBetabar} with $\WangD^{\vec{m}}$ using \eqref{EQ:ang3ijk}-\eqref{EQ:ang4ijk}. Let $F:=\CBetabar$, $G:=\CAlphabar$, and
					\begin{align}
						\label{PV:BetabarAlphabar1}P_{(\CBetabar,\CAlphabar;\vec{m})}=&~\sum\limits_{\vec{m_1}+\vec{m_2}+\vec{m_3}=(i-i_0,j,k-k_0)}\delta^{i_0}\abs{u}^{k_0}\WangD^{\vec{m_1}}\upsi^{k_0+1}\cdot\WangD^{\vec{m_2}}\psi^{i_0}\cdot\WangD^{\vec{m_3}}(\CAlphabar,\CBetabar,\FAlphabar,\FBetabar)\\
						\notag&+\sum\limits_{\vec{m_1}+\vec{m_2}+\vec{m_3}+\vec{m_4}=(i-i_0-1,j,k-k_0+1)}\delta^{i_0+1}\abs{u}^{k_0-1}\WangD^{\vec{m_1}}\upsi^{k_0}\cdot\WangD^{\vec{m_2}}\psi^{i_0}\cdot\WangD^{\vec{m_3}}\upzeta\cdot\WangD^{\vec{m_4}}\CBetabar\\
						\notag&+\sum\limits_{\vec{m_1}+\vec{m_2}+\vec{m_3}+\vec{m_4}=(i-i_0-1,j,k-k_0)}\delta^{i_0+1}\abs{u}^{k_0}\WangD^{\vec{m_1}}\upsi^{k_0}\cdot\WangD^{\vec{m_2}}\psi^{i_0}\cdot\WangD^{\vec{m_3}}\Rhodual\cdot\WangD^{\vec{m_4}}\CBetabar\\
						\notag&+\sum\limits_{\vec{m_1}+\vec{m_2}+\vec{m_3}+\vec{m_4}=(i-i_0,j,k-k_0-1)}\delta^{i_0}\abs{u}^{k_0+1}\WangD^{\vec{m_1}}\upsi^{k_0}\cdot\WangD^{\vec{m_2}}\psi^{i_0}\cdot\WangD^{\vec{m_3}}(\CBetabar,\FBetabar)\cdot\WangD^{\vec{m_4}}\CBetabar,	
					\end{align}
					\begin{align}
						\label{PV:BetabarAlphabar2}Q_{(\CBetabar,\CAlphabar;\vec{m})}=&~\sum\limits_{\vec{m_1}+\vec{m_2}+\vec{m_3}=(i,j- j_0,k-k_0)}\abs{u}^{j_0+k_0}\WangD^{\vec{m_1}}\upsi^{j_0+k_0}\cdot\WangD^{\vec{m_2}}\psi\cdot\WangD^{\vec{m_3}}(\CAlphabar,\FAlphabar)\\
						\notag&+\sum\limits_{\vec{m_1}+\vec{m_2}=(i,j- j_0,k-k_0)}\abs{u}^{j_0+k_0}\WangD^{\vec{m_1}}\upsi^{j_0+k_0+1}\cdot\WangD^{\vec{m_2}}(\CBetabar,\Scalar,\FRho,\FRhodual,\CRho,\CRhodual)\\
						\notag&+\sum\limits_{\vec{m_1}+\vec{m_2}=(i,j- j_0,k-k_0)}\abs{u}^{j_0+k_0}\WangD^{\vec{m_1}}\upsi^{k_0+j_0}\cdot\WangD^{\vec{m_2}}(\angnabla\FBetabar,\angD_{ 3}\Scalar)\\
						\notag&+\sum\limits_{\vec{m_1}+\vec{m_2}+\vec{m_3}=(i,j- j_0,k-k_0-1)}\abs{u}^{j_0+k_0+1}\WangD^{\vec{m_1}}\upsi^{j_0+k_0}\cdot\WangD^{\vec{m_3}}(\CBetaA,\FBetaA)\cdot\WangD^{\vec{m_3}}\CAlphabar\\
						\notag&+\sum\limits_{\vec{m_1}+\vec{m_2}+\vec{m_3}=(i,j- j_0-1,k-k_0+1)}\abs{u}^{j_0+k_0}\WangD^{\vec{m_1}}\upsi^{j_0+k_0}\cdot\WangD^{\vec{m_2}}\upzeta\cdot\WangD^{\vec{m_3}}\CAlphabar\\
						\notag&+\sum\limits_{\vec{m_1}+\vec{m_2}+\vec{m_3}=(i,j- j_0-1,k-k_0)}\abs{u}^{j_0+k_0+1}\WangD^{\vec{m_1}}\upsi^{j_0+k_0}\cdot\WangD^{\vec{m_2}}\Rhodual\cdot\angD^{\vec{m_3}}\CAlphabar.
					\end{align}
					We apply Proposition \ref{PR:Int4} with $s(\CBetabar)=\frac{3}{2}$ to obtain
					\begin{align}\label{PV:BetabarAlphabar3}
						&\twonorms{\abs{u}^3\WangD^{\vec{m}}\CBetabar}{\ubar}{2}{\sphere}{2}{\hypersur_{u}}^2+\twonorms{\abs{u}^3\WangD^{\vec{m}}\CAlphabar}{u}{2}{\sphere}{2}{\hypersurbar_{\ubar}}^2\\
						\notag\lesssim&~\twonorms{\abs{u}^3\WangD^{\vec{m}}\CBetabar}{\ubar}{2}{\sphere}{2}{\hypersur_{\uinf}}^2+\twonorms{\abs{u}^3\WangD^{\vec{m}}\CAlphabar}{u}{2}{\sphere}{2}{\hypersurbar_{0}}^2\\
						\notag&+\int_{\region}\abs{u}^6\abs{\WangD^{\vec{m}}\CBetabar}\cdot \abs{P_{(\CBetabar,\CAlphabar;\vec{m})}}\diff\gvol+\int_{\region}\abs{u}^6\abs{\WangD^{\vec{m}}\CAlphabar}\cdot \abs{Q_{(\CBetabar,\CAlphabar;\vec{m})}}\diff\gvol.
					\end{align}
					We focus on estimating $\threenorms{\abs{u}^3P_{(\CBetabar,\CAlphabar;\vec{m})}}{u}{1}{\ubar}{2}{\sphere}{2}{\region}$ and $\threenorms{\abs{u}^3Q_{(\CBetabar,\CAlphabar;\vec{m})}}{u}{2}{\ubar}{2}{\sphere}{2}{\region}$ since 
					\begin{align}\label{PV:BetabarAlphabar4}
						&\int_{\region}\abs{u}^6\abs{\WangD^{\vec{m}}\CBetabar}\cdot \abs{P_{(\CBetabar,\CAlphabar;\vec{m})}}\diff\gvol+\int_{\region}\abs{u}^6\abs{\WangD^{\vec{m}}\CAlphabar}\cdot \abs{Q_{(\CBetabar,\CAlphabar;\vec{m})}}\diff\gvol\\
						\notag\lesssim&~\threenorms{\abs{u}^3\WangD^{\vec{m}}\CBetabar}{u}{\infty}{\ubar}{2}{\sphere}{2}{\region}\threenorms{\abs{u}^3P_{(\CBetabar,\CAlphabar;\vec{m})}}{u}{1}{\ubar}{2}{\sphere}{2}{\region}+\delta^{\frac{1}{2}}\threenorms{\abs{u}^3\WangD^{\vec{m}}\CAlphabar}{u}{2}{\ubar}{\infty}{\sphere}{2}{\region}\threenorms{\abs{u}^3Q_{(\CBetabar,\CAlphabar;\vec{m})}}{u}{2}{\ubar}{2}{\sphere}{2}{\region}.
					\end{align}
					For the terms involving $\CBetabar$, using Proposition \ref{PR:EstimatesforRicciCurvature} for $\FBetabar$, bootstrap assumptions \eqref{BA} and Sobolev embedding $H^2\hookrightarrow L^\infty$ on $\suubar$ for $\WangD^{\vec{m}}\CBetabar$, we have:
					\begin{align}\label{PV:BetabarAlphabar5}
						&\sum\limits_{\vec{m}_1+\vec{m}_2+\vec{m}_3+\vec{m}_4=(i-i_0,j,k-k_0)}\delta^{i_0}\threenorms{\abs{u}^{k_0+3}\WangD^{\vec{m}_1}\psi^{i_0}\cdot\WangD^{\vec{m}_2}\upsi^{k_0}\cdot\left(\WangD^{\vec{m}_3}\upsi,\delta\WangD^{\vec{m}'_3}\Rhodual,\abs{u}\WangD^{\vec{m}'_3}\CBetabar,\abs{u}\WangD^{\vec{m}'_3}\FBetabar\right)\cdot\WangD^{\vec{m}_4}\CBetabar}{u}{1}{\ubar}{2}{\sphere}{2}{\region}\\
						\notag\lesssim&~\sum\limits_{\substack{\vec{m}_1+\vec{m}_2+\vec{m}_3+\vec{m}_4=(i-i_0,j,k-k_0)\\\abs{\vec{m}_4}<\frac{N-1}{2}}}\delta^{i_0}\threenorms{\abs{u}^{k_0}\WangD^{\vec{m}_1}\psi^{i_0}\cdot\WangD^{\vec{m}_2}\upsi^{k_0}\cdot\left(\WangD^{\vec{m}_3}\upsi,\delta\WangD^{\vec{m}'_3}\Rhodual,\abs{u}\WangD^{\vec{m}'_3}\CBetabar,\abs{u}\WangD^{\vec{m}'_3}\FBetabar\right)}{u}{1}{\ubar}{\infty}{\sphere}{2}{\region}\\
						\notag&\cdot\threenorms{\abs{u}^3\WangD^{\vec{m}_4}\CBetabar}{u}{\infty}{\ubar}{2}{\sphere}{\infty}{\region}\\
						\notag&+\sum\limits_{\substack{\vec{m}_1+\vec{m}_2+\vec{m}_3+\vec{m}_4=(i-i_0,j,k-k_0)\\\abs{\vec{m}_4}\geq\frac{N-1}{2}}}\delta^{i_0}\threenorms{\abs{u}^{k_0}\WangD^{\vec{m}_1}\psi^{i_0}\cdot\WangD^{\vec{m}_2}\upsi^{k_0}\cdot\left(\WangD^{\vec{m}_3}\upsi,\delta\WangD^{\vec{m}'_3}\Rhodual,\abs{u}\WangD^{\vec{m}'_3}\CBetabar,\abs{u}\WangD^{\vec{m}'_3}\FBetabar\right)}{u}{1}{\ubar}{\infty}{\sphere}{\infty}{\region}\\
						\notag&\cdot\threenorms{\abs{u}^3\WangD^{\vec{m}_4}\CBetabar}{u}{\infty}{\ubar}{2}{\sphere}{2}{\region}\\
						\notag\lesssim&~\left(\BAsize\frac{\delta a^{\frac{1}{2}}}{\abs{u}}\right)^{i_0+k_0}\left(\BAsize\frac{\delta a^{\frac{1}{2}}}{\abs{u}}+ b^{-\frac{1}{2}}\right)\cdot\sum\limits_{\abs{\vec{m}_4}\leq N-1}\threenorms{\abs{u}^3\WangD^{\vec{m}_4}\CBetabar}{u}{\infty}{\ubar}{2}{\sphere}{2}{\region},
					\end{align}
					where $\abs{\vec{m}_3'}=\abs{\vec{m}_3}-1$.
					Next, we estimate the rest terms on the RHS of \eqref{PV:BetabarAlphabar1}. Using bootstrap assumptions \eqref{BA} and Proposition \ref{PR:EstimatesforRicciCurvature}, we have:
					\begin{align}\label{PV:BetabarAlphabar6}
						&\sum\limits_{\vec{m}_1+\vec{m}_2+\vec{m}_3=(i-i_0,j,k-k_0)}\delta^{i_0}\cdot\threenorms{\abs{u}^{k_0+3}\WangD^{\vec{m}_1}\psi^{i_0}\cdot\WangD^{\vec{m}_2}\upsi^{k_0+1}\cdot\WangD^{\vec{m}_3}\left(\FBetabar,\CAlphabar,\FAlphabar\right)}{u}{1}{\ubar}{2}{\sphere}{2}{\region}\\
						\notag\lesssim&~\sum\limits_{\substack{\vec{m}_1+\vec{m}_2+\vec{m}_3=(i-i_0,j,k-k_0)\\\abs{\vec{m}_3}<\frac{N-1}{2}}}\delta^{i_0}\left(\threenorms{\abs{u}^{k_0+3}\WangD^{\vec{m}_1}\psi^{i_0}\cdot\WangD^{\vec{m}_2}\upsi^{k_0+1}}{u}{1}{\ubar}{2}{\sphere}{2}{\region}\cdot\threenorms{\WangD^{\vec{m}_3}\left(\FBetabar,\CAlphabar,\FAlphabar\right)}{u}{\infty}{\ubar}{\infty}{\sphere}{\infty}{\region}\right)\\
						\notag&+\sum\limits_{\substack{\vec{m}_1+\vec{m}_2+\vec{m}_3=(i-i_0,j,k-k_0)\\\abs{\vec{m}_3}\geq\frac{N-1}{2}}}\delta^{i_0}\left(\threenorms{\abs{u}^{k_0}\WangD^{\vec{m}_1}\psi^{i_0}\cdot\WangD^{\vec{m}_2}\upsi^{k_0+1}}{u}{2}{\ubar}{\infty}{\sphere}{\infty}{\region}\cdot\threenorms{\abs{u}^3\WangD^{\vec{m}_3}\left(\FBetabar,\CAlphabar,\FAlphabar\right)}{u}{\infty}{\ubar}{2}{\sphere}{2}{\region}\right)\\
						\notag\lesssim&~\left(\BAsize\frac{\delta a^{\frac{1}{2}}}{\abs{u}}\right)^{i_0+k_0}\cdot\left(\BAsize \delta^{2} a^{\frac{1}{2}}+\BAsize\frac{\delta^{\frac{5}{2}} a^{\frac{3}{2}}}{\abs{u}^{\frac{3}{2}}}\right).
					\end{align}
					Now, we estimate the terms on the RHS of \eqref{PV:BetabarAlphabar2}. For the gravity terms and the below-top-order fluid terms, using Proposition \ref{PR:EstimatesforRicciCurvature}, bootstrap assumptions \eqref{BA}, we have:
					\begin{align}\label{PV:BetabarAlphabar8.5}
						&\sum\limits_{\vec{m}_1+\vec{m}_2=(i,j-j_0,k-k_0)}\threenorms{\abs{u}^{j_0+k_0+3}\WangD^{\vec{m}_1}\upsi^{j_0+k_0+1}\cdot\WangD^{\vec{m}_2}(\CRho,\CRhodual,\Scalar,\FRho,\FRhodual,\CBetabar)}{u}{2}{\ubar}{2}{\sphere}{2}{\region}\\
						\notag\lesssim&~\sum\limits_{\substack{\vec{m}_1+\vec{m}_2=(i,j-j_0,k-k_0)\\\abs{\vec{m}_2}<\frac{N-1}{2}}}\threenorms{\abs{u}^{j_0+k_0+3}\WangD^{\vec{m}_1}\upsi^{j_0+k_0+1}}{u}{2}{\ubar}{2}{\sphere}{2}{\region}\cdot\threenorms{\WangD^{\vec{m}_2}(\CRho,\CRhodual,\Scalar,\FRho,\FRhodual,\CBetabar)}{u}{\infty}{\ubar}{\infty}{\sphere}{\infty}{\region}\\
						\notag&+\sum\limits_{\substack{\vec{m}_1+\vec{m}_2=(i,j-j_0,k-k_0)\\\abs{\vec{m}_2}\geq\frac{N-1}{2}}}\threenorms{\abs{u}^{j_0+k_0}\WangD^{\vec{m}_1}\upsi^{j_0+k_0+1}}{u}{2}{\ubar}{\infty}{\sphere}{\infty}{\region}\cdot\threenorms{\abs{u}^3\WangD^{\vec{m}_2}(\CRho,\CRhodual,\Scalar,\FRho,\FRhodual,\CBetabar)}{u}{\infty}{\ubar}{2}{\sphere}{2}{\region}\\
						\notag\lesssim&~\left(\BAsize\frac{\delta a^{\frac{1}{2}}}{\abs{u}}\right)^{j_0+k_0}\left(\BAsize\frac{\delta^{\frac{3}{2}} a}{\abs{u}^{\frac{1}{2}}}+ \delta a^{\frac{1}{2}}b^{-\frac{1}{2}}\right),
					\end{align}
					\begin{align}\label{PV:BetabarAlphabar9}
						&\sum\limits_{\vec{m}_1+\vec{m}_2+\vec{m}_3=(i,j-j_0,k-k_0)}\threenorms{\abs{u}^{j_0+k_0+3}\WangD^{\vec{m}_1}\upsi^{j_0+k_0}\cdot(\WangD^{\vec{m}_2}\psi,\abs{u}\WangD^{\vec{m}'_2}(\CBetaA,\FBetaA,\Rhodual))\cdot\WangD^{\vec{m}_3}(\CAlphabar,\FAlphabar)}{u}{2}{\ubar}{2}{\sphere}{2}{\region}\\
						\notag\lesssim&~\sum\limits_{\substack{\vec{m}_1+\vec{m}_2+\vec{m}_3=(i,j-j_0,k-k_0)\\\abs{\vec{m}_3}<\frac{N-1}{2}}}\threenorms{\abs{u}^{j_0+k_0+3}\WangD^{\vec{m}_1}\upsi^{j_0+k_0}\cdot(\WangD^{\vec{m}_2}\psi,\abs{u}\WangD^{\vec{m}'_2}(\CBetaA,\FBetaA,\Rhodual))}{u}{2}{\ubar}{2}{\sphere}{2}{\region}\\
						\notag&\cdot\threenorms{\WangD^{\vec{m}_3}(\CAlphabar,\FAlphabar)}{u}{\infty}{\ubar}{\infty}{\sphere}{\infty}{\region}\\
						\notag&+\sum\limits_{\substack{\vec{m}_1+\vec{m}_2+\vec{m}_3=(i,j-j_0,k-k_0)\\\abs{\vec{m}_3}\geq\frac{N-1}{2}}}\threenorms{\abs{u}^{j_0+k_0}\WangD^{\vec{m}_1}\upsi^{j_0+k_0}\cdot(\WangD^{\vec{m}_2}\psi,\abs{u}\WangD^{\vec{m}'_2}(\CBetaA,\FBetaA,\Rhodual))}{u}{2}{\ubar}{\infty}{\sphere}{\infty}{\region}\\
						\notag&\cdot\threenorms{\abs{u}^3\WangD^{\vec{m}_3}(\CAlphabar,\FAlphabar)}{u}{\infty}{\ubar}{2}{\sphere}{2}{\region}\\
						\notag\lesssim&~\left(\BAsize\frac{\delta a^{\frac{1}{2}}}{\abs{u}}\right)^{j_0+k_0}\cdot\left(\BAsize\delta a+ a^{\frac{1}{2}}b^{-1}+ \delta^{-\frac{1}{2}}b^{-1}\right),
					\end{align}
					where $\abs{\vec{m}_2'}=\abs{\vec{m}_2}-1$.
					Finally, we estimate the top-order terms, namely, the third terms on the RHS of \eqref{PV:BetabarAlphabar2} using Proposition \ref{PR:EstimatesforRicciCurvature}, bootstrap assumptions \eqref{BA}:
					\begin{align}\label{PV:BetabarAlphabar10}
						&\sum\limits_{\vec{m}_1+\vec{m}_2=(i,j-j_0,k-k_0)}\threenorms{\abs{u}^{j_0+k_0+3}\WangD^{\vec{m}_1}\upsi^{j_0+k_0}\cdot\WangD^{\vec{m}_2}(\angnabla\FBetabar,\angD_{ 3}\Scalar)}{u}{2}{\ubar}{2}{\sphere}{2}{\region}\\
						\notag\lesssim&~\sum\limits_{\substack{\vec{m}_1+\vec{m}_2=(i,j-j_0,k-k_0)\\\abs{\vec{m}_2}<\frac{N-1}{2}}}\threenorms{\abs{u}^{j_0+k_0+3}\WangD^{\vec{m}_1}\upsi^{j_0+k_0}}{u}{2}{\ubar}{2}{\sphere}{2}{\region}\cdot\threenorms{\WangD^{\vec{m}_2}(\angnabla\FBetabar,\angD_{ 3}\Scalar)}{u}{\infty}{\ubar}{\infty}{\sphere}{\infty}{\region}\\
						\notag&+\sum\limits_{\substack{\vec{m}_1+\vec{m}_2=(i,j-j_0,k-k_0)\\\abs{\vec{m}_2}\geq\frac{N-1}{2}}}\threenorms{\abs{u}^{j_0+k_0+3}\WangD^{\vec{m}_1}\upsi^{j_0+k_0}}{u}{2}{\ubar}{\infty}{\sphere}{\infty}{\region}\cdot\threenorms{\WangD^{\vec{m}_2}(\angnabla\FBetabar,\angD_{ 3}\Scalar)}{u}{\infty}{\ubar}{2}{\sphere}{2}{\region}\\
						\notag\lesssim&~\left(\BAsize\frac{\delta a^{\frac{1}{2}}}{\abs{u}}\right)^{j_0+k_0}\cdot b^{-\frac{1}{2}},
					\end{align}
					Combining \eqref{PV:BetabarAlphabar1}-\eqref{PV:BetabarAlphabar10}, summing \eqref{PV:BetabarAlphabar3} over $\abs{\vec{m}}\leq N-1$, we soak the $\sum\limits_{\abs{\vec{m}_4}\leq N-1}\threenorms{\WangD^{\vec{m}_4}\CBetabar}{u}{\infty}{\ubar}{2}{\sphere}{2}{\region}$ terms in \eqref{PV:BetabarAlphabar5} to the LHS of \eqref{PV:BetabarAlphabar4} since their coefficients are $\ll 1$. Then, we obtain the desired estimates \eqref{ES:topL2curvature3} by the initial data assumption \eqref{Data} and the choice of parameter $\delta,a,b$ in Section \ref{SS:Norms}, $\BAfluid,\BAsize$ in Section \eqref{SS:BA}.
				\end{proof}
				\section{The Formation of Trapped Surfaces}\label{S:Formation}
				\subsection{Corollaries of the Control of Geometry}
				\begin{lemma}[Improved control of lapse function $\lapse$]\label{LE:improvedlapse}
					Under the initial data assumptions in Section \ref{SS:InitialData} and estimates for $\omega$ as in \eqref{ES:LinftyRicciCoefficients3}, we have:
					\begin{align}\label{ES:improvedlapse1}
						\norm{\lapse^{-1}-1}_{L^\infty(\suubar)}\lesssim\frac{\delta a^{\frac{1}{2}}}{\abs{u}}.
					\end{align}
					Moreover,
					\begin{align}\label{ES:improvedlapse2}
						\norm{\lapse-1}_{L^\infty(\suubar)}\lesssim\frac{\delta a^{\frac{1}{2}}}{\abs{u}}.
					\end{align}
				\end{lemma}
				\begin{proof}[Proof of Lemma \ref{LE:improvedlapse}]
					The proof of \eqref{ES:improvedlapse1} is similar to the proof of proof of Proposition \ref{PR:lapse}. By \eqref{ES:LinftyRicciCoefficients3}, we integrate the equation \eqref{EQ:omegauomega} along $\hypersur_{u}$ and use the fact that $\lapse=1$ on $\hypersurbar_{0}$ to obtain the desired result. 
					
					\eqref{ES:improvedlapse2} follows directly from \eqref{ES:improvedlapse1} and the fact that $\norm{\lapse-1}_{L^\infty(\suubar)}\leq \norm{\lapse}_{L^\infty(\suubar)}\norm{\lapse^{-1}-1}_{L^\infty(\suubar)}$
				\end{proof}
				\begin{lemma}[Estimate of angular coordinate function $\Lambda^A$]\label{LE:EstimateofLambda}
					For $A=1,2$, the angular coordinate function $\Lambda^A$ defined in \eqref{DE:Coordinatese3e4} satisfies the following estimate:
					\begin{align}\label{ES:EstimateofLambda}
						\norm{\Lambda^A}_{L^\infty(\suubar)}\lesssim\frac{\delta^2a^\frac{1}{2}}{\abs{u}^2}.
					\end{align}
				\end{lemma}
				\begin{proof}[Proof of Lemma \ref{LE:EstimateofLambda}]
					By Proposition \ref{PR:commutator}, we have
					\begin{align}\label{ES:EstimateofLambda1}
						[\lapse e_3,\lapse e_4]=4\lapse^2\upzeta\cdot\angnabla.
					\end{align}
					\eqref{ES:EstimateofLambda1} and \eqref{DE:Coordinatese3e4} implies that
					\begin{align}\label{ES:EstimateofLambda2}
						\frac{\p\Lambda}{\p\ubar}=-4\lapse^2\upzeta.
					\end{align}
					Integrating \eqref{ES:EstimateofLambda2} over $\ubar$, using Corollary \ref{CO:LinftyRicciCoefficients}, \eqref{ES:improvedlapse2} and the fact that $\Lambda^A=0$ on $\hypersurbar_{0}$, we obtain the desired result \eqref{ES:EstimateofLambda}.
				\end{proof}
				\subsection{First Trapped Surface Formation Argument}
				\begin{lemma}[Estimates of $\abs{u}^2\sgabs{\hchi}^2$]\label{LE:Estimatesforuhchi}
					Under the initial data assumptions \eqref{Data}, \eqref{AS:SmallonHbar}, \eqref{AS:LargenessofData} and \eqref{AS:LargenessofData1}, we have the following bound for $\abs{\uts}^2\sgabs{\hchi}^2$:
					\begin{align}\label{ES:Estimatesforuhchi}
						\int_{0}^{\delta}\abs{\uts}^2\sgabs{\hchi(\uts,\ubarprime,\sphere^1,\sphere^2)}^2\diff \ubarprime&>4\abs{\uts},&
						\text{for all } (\sphere^1,\sphere^2)&\in\Stwo.
					\end{align}
				\end{lemma}
				\begin{proof}[Proof of Lemma \ref{LE:Estimatesforuhchi}]
					Contracting \eqref{EQ:angD3hchi} with $\hchi$, we obtain:
					\begin{align}
						\label{PV:Estimatesforuhchi1}\frac{1}{2}e_3\left(\sgabs{\hchi}^2\right)+\frac{1}{2}(\gtr\uchi)\sgabs{\hchi}^2-2\uomega\sgabs{\hchi}^2=&~\left(-\frac{1}{2}\gtr\upchi\huchi+\angnabla\tensorproduct\upeta+\upeta\tensorproduct\upeta+\FRhodual\right)\cdot\hchi.
					\end{align}
					By \eqref{EQ:omegauomega}, multiplying \eqref{PV:Estimatesforuhchi1} by $2\lapse^2$, we re-write \eqref{PV:Estimatesforuhchi2} as follows:
					\begin{align}
						\label{PV:Estimatesforuhchi2}e_3\left(\lapse^2\sgabs{\hchi}^2\right)+\lapse^2(\gtr\uchi)\sgabs{\hchi}^2=&~2\lapse^2\left(-\frac{1}{2}\gtr\upchi\huchi+\angnabla\tensorproduct\upeta+\upeta\tensorproduct\upeta+\FRhodual\right)\cdot\hchi.
					\end{align}
					Using coordinates expression \eqref{DE:Coordinatese3e4}, multiplying \eqref{PV:Estimatesforuhchi2} by $u^2\lapse$ reads:
					\begin{align}
						\label{PV:Estimatesforuhchi3}\frac{\p}{\p u}\left(u^2\lapse^{2}\sgabs{\hchi}^2\right)=&~2\abs{u}^2\lapse^3\left(-\frac{1}{2}\gtr\upchi\huchi+\angnabla\tensorproduct\upeta+\upeta\tensorproduct\upeta+\FRhodual\right)\cdot\hchi-\abs{u}^2\Lambda^A\frac{\p}{\p\theta^A}\left(\lapse^{2}\sgabs{\hchi}^2\right)\\
						\notag&-\abs{u}^2\lapse\left(\gtr\uchi+\frac{2}{\abs{u}}\right)\left(\lapse^{2}\sgabs{\hchi}^2\right)+2\abs{u}(\lapse-1)\left(\lapse^{2}\sgabs{\hchi}^2\right).
					\end{align}
					Now we estimate the terms on the RHS of \eqref{PV:Estimatesforuhchi3}. By Corollary \ref{CO:LinftyRicciCoefficients}, Proposition \ref{PR:EstimatesforRicciCurvature}, Lemma \ref{LE:improvedlapse}, and Lemma \ref{LE:EstimateofLambda}, we have the following bounds for the terms on the RHS of \eqref{PV:Estimatesforuhchi3} as follows: 
					\begin{align}\label{PV:Estimatesforuhchi4}
						\norm{2\abs{u}^2\lapse^3\left(-\frac{1}{2}\gtr\upchi\huchi+\angnabla\tensorproduct\upeta+\upeta\tensorproduct\upeta+\FRhodual\right)\cdot\hchi}_{L^\infty(\suubar)}&\lesssim\abs{u}^2\cdot\frac{a^{\frac{1}{2}}}{\abs{u}}\cdot\left(\delta^2 a^{\frac{3}{2}}+ \right)\lesssim a^{\frac{1}{2}},\\
						\norm{\abs{u}^2\Lambda^A\frac{\p}{\p\theta^A}\left(\lapse^{2}\sgabs{\hchi}^2\right)}_{L^\infty(\suubar)}&\lesssim\abs{u}^2\cdot\frac{\delta^2 a^{\frac{1}{2}}}{\abs{u}^2}\cdot\frac{a}{\abs{u}^2}\lesssim \frac{\delta^2 a^{\frac{3}{2}}}{\abs{u}^2},\\
						\norm{\abs{u}^2\lapse\left(\gtr\uchi+\frac{2}{\abs{u}}\right)\left(\lapse^{2}\sgabs{\hchi}^2\right)}_{L^\infty(\suubar)}&\lesssim\abs{u}^2\cdot\delta\cdot\frac{a}{\abs{u}^2}\lesssim \delta a,\\
						\norm{2\abs{u}(\lapse-1)\left(\lapse^{2}\sgabs{\hchi}^2\right)}_{L^\infty(\suubar)}&\lesssim\abs{u}\cdot\frac{\delta a^{\frac{1}{2}}}{\abs{u}}\cdot\frac{a}{\abs{u}^2}\lesssim \frac{\delta a^{\frac{3}{2}}}{\abs{u}^2}.
					\end{align}
					Combining \eqref{PV:Estimatesforuhchi3}-\eqref{PV:Estimatesforuhchi4}, we have:
					\begin{align}
						\label{PV:Estimatesforuhchi5}\abs{\frac{\p}{\p u}\left(u^2\lapse^{2}\sgabs{\hchi}^2\right)}\lesssim&~a^{\frac{1}{2}}.
					\end{align}
					Integrate \eqref{PV:Estimatesforuhchi3} over $u$ then over $\ubar$ and using \eqref{PV:Estimatesforuhchi5}, for any $(\sphere^1,\sphere^2)\in\Stwo$, we have:
					\begin{align}\label{PV:Estimatesforuhchi6}
						\int_{0}^{\delta}\abs{\uinf}^2\lapse^{2}\sgabs{\hchi}^2(\uinf,\ubarprime,\sphere^1,\sphere^2)\diff\ubarprime-C\delta a^{\frac{1}{2}}&\leq \int_{0}^{\delta}\abs{\uts}^2\lapse^{2}\sgabs{\hchi}^2(\uts,\ubarprime,\sphere^1,\sphere^2)\diff\ubarprime\\
						\notag&\leq \int_{0}^{\delta}\abs{\uinf}^2\lapse^{2}\sgabs{\hchi}^2(\uinf,\ubarprime,\sphere^1,\sphere^2)\diff\ubarprime+C\delta a^{\frac{1}{2}},
					\end{align}
					where $C$ is a universal constant ($C\ll\delta^{-1}a^{-\frac{1}{2}}$).
					Combining \eqref{PV:Estimatesforuhchi6} and the assumption on the largeness of $\hchi$ in \eqref{AS:LargenessofData}, we obtain the desired result \eqref{ES:Estimatesforuhchi}.
				\end{proof}
				\subsubsection{Proof of the first case in Theorem \ref{TH:Main}}
				$\gtr\uchi<0$ in $\region$ is directly given by the estimates \eqref{ES:LinftyRicciCoefficients5}. By data prescription in Section \ref{SS:CharData} and assumptions \eqref{AS:SmallonHbar}, \eqref{AS:LargenessofData}, $\gtr\upchi>0$ on $\hypersur_{\uinf}\cup\hypersurbar_{0}$. 
				
				Now, provided assumption \eqref{AS:LargenessofData1}, we show that $\gtr\upchi<0$ everywhere on $S_{\uts,\delta}$. Multiplying equation \eqref{EQ:e4upchi} by $\lapse^{-1}$, using \eqref{EQ:omegauomega}, we have:
				\begin{align}\label{PV:Main1}
					e_4\left(\lapse^{-1}\gtr\upchi\right)=&~-\frac{\lapse^{-1}}{2}(\gtr\upchi)^2-\lapse^{-1}\sgabs{\hchi}^2-2\lapse^{-1}\FAlphaAB.
				\end{align}
				Using the coordinates expression \eqref{DE:Coordinatese3e4}, \eqref{PV:Main1} reads:
				\begin{align}\label{PV:Main2}
					\frac{\p}{\p\ubar}\left(\lapse^{-1}\gtr\upchi\right)=&~-\frac{1}{2}(\gtr\upchi)^2-\sgabs{\hchi}^2-2\FAlphaAB.
				\end{align}
				Integrating \eqref{PV:Main2} over $\ubar$, by the assumption on $\gtr\upchi$ on $\hypersurbar_{0}$ as in \eqref{AS:SmallonHbar}, Lemma \ref{LE:improvedlapse}, and Lemma \ref{LE:Estimatesforuhchi}, for any $(\sphere^1,\sphere^2)\in\Stwo$, we derive
				\begin{align}
					\lapse^{-1}\gtr\upchi(\uts,\delta,\sphere^1,\sphere^2)&\leq \lapse^{-1}\gtr\upchi(\uts,0,\sphere^1,\sphere^2)-\int_{0}^{\delta}\left(\sgabs{\hchi}^2+2\FAlphaAB\right)(\uts,\ubarprime,\sphere^1,\sphere^2)\diff\ubarprime\\
					\notag&< \frac{3}{\abs{\uts}}-\frac{4}{\abs{\uts}}+\mathcal{O}(\delta a^{\frac{1}{2}})
					\\
					\notag&<0.
				\end{align}
				\subsection{Second Trapped Surface Formation Argument}
				We consider the region $\mathcal{D}_{\uinf,\delta}^{\uts,2\delta}$ bounded by $\hypersur_{\uinf},\hypersur_{\uts},\hypersurbar_{0}$ and $\mathcal{T}:=\{\ubar=\left(1+2(u-\uinf)\right)\delta,u\in[\uinf,\uts]\}$. That is
				\begin{align}
					\mathcal{D}_{\uinf,\delta}^{\uts,2\delta}=\left\{(u,\ubar,\sphere^1,\sphere^2)|u\in[\uinf,\uts], \ubar\in\left[0,\left(1+2(u-\uinf)\right)\delta\right]\right\}.
				\end{align}
				$\mathcal{T}$ is a timelike manifold (generated by $\frac{1}{\sqrt{8\delta}}(2\delta e_4+e_3)$). The unit outer normal $\mathbf{S}$ of $\mathcal{T}$ is as follows:
				\begin{align}
					\mathbf{S}=\frac{1}{\sqrt{8\delta}}\left(2\delta e_4-e_3\right).
				\end{align}
				\begin{lemma}[Estimates of $\abs{u}^2\cdot2\FAlphaAB$]\label{LE:EstimatesforuFAlpha}
					Under the initial data assumptions \eqref{Data}, \eqref{AS:SmallonHbar}, \eqref{AS:LargenessofData} and \eqref{AS:LargenessofData2}, we have the following bound for $\abs{\uts}^2\FAlphaAB$:
					\begin{align}\label{ES:EstimatesforuFAlpha}
						\int_{0}^{\delta+2(\uts-\uinf)\delta}\abs{\uts}^22\FAlphaAB(\uts,\ubarprime,\sphere^1,\sphere^2)\diff \ubarprime&>4\abs{\uts},&
						\text{for all } (\sphere^1,\sphere^2)&\in\Stwo.
					\end{align}
				\end{lemma}
				\begin{proof}[Proof of Lemma \ref{LE:EstimatesforuFAlpha}]
					We apply divergence theorem on \eqref{EQ:DivCurrent} in the region $\mathcal{D}_{\uinf,\delta}^{\uts,2\delta}$:
					\begin{align}
						0=&~\int_{\mathcal{D}_{\uinf,\delta}^{\uts,2\delta}}\Dfour_\alpha I^\alpha \diff\gvol\\
						\notag=&~-\int_{\hypersur_{\uts}\cap\{\ubar\in[0,\delta+2(\uts-\uinf)\delta]\}}\gfour(I,e_4)\diff\spherevol\diff \ubar\\
						\notag&+\int_{\hypersur_{\uinf}\cap\{\ubar\in[0,\delta]\}}\gfour(I,e_4)\diff\spherevol\diff \ubar+\int_{\hypersurbar_{0}}\gfour(I,e_3)\diff\spherevol\diff u+\int_{\mathcal{T}}\gfour(I,\mathbf{S}) \diff\gvol.
					\end{align}
					Express $I^\alpha:=nv^\alpha$, we have:
					\begin{align}\label{EQ:2ndTrapped1}
						\int_{\hypersur_{\uts}\cap\{\ubar\in[0,\delta+2(\uts-\uinf)\delta]\}}n v_4\diff\spherevol\diff \ubar=\int_{\hypersur_{\uinf}\cap\{\ubar\in[0,\delta]\}}n v_4\diff\spherevol\diff \ubar+\int_{\hypersurbar_{0}}nv_3\diff\spherevol\diff u+\int_{\mathcal{T}}n v_\alpha \mathbf{S}^\alpha \diff\gvol.
					\end{align}
					By Corollary \ref{CO:LinftyFluid} and the fact that $\delta b\leq \abs{\uts} \leq 1$, the integrand in the last term of the \eqref{EQ:2ndTrapped1} is strictly negative. Then, \eqref{EQ:2ndTrapped1} gives the following estimate
					\begin{align}\label{EQ:2ndTrapped2}
						\int_{\hypersur_{\uts}\cap\{\ubar\in[0,\delta+2(\uts-\uinf)\delta]\}}n v_4\diff\spherevol\diff \ubar\leq \int_{\hypersur_{\uinf}\cap\{\ubar\in[0,\delta]\}}n v_4\diff\spherevol\diff \ubar.
					\end{align}
					Recall $v_4<0$. By \eqref{ES:LinftyFluid2.3}, we have $\abs{b^{\frac{1}{2}}v_4(\uinf,\ubar,\sphere^1,\sphere^2)-b^{\frac{1}{2}}v_4(\uts,\ubar,\sphere^1,\sphere^2)}\leq \BAfluid$. Thus \eqref{EQ:2ndTrapped2} implies
					\begin{align}\label{EQ:2ndTrapped3}
						\int_{\hypersur_{\uts}\cap\{\ubar\in[0,\delta+2(\uts-\uinf)\delta]\}}n\diff\spherevol\diff \ubar \geq \frac{7}{8} \int_{\hypersur_{\uinf}\cap\{\ubar\in[0,\delta]\}}n\diff\spherevol\diff \ubar.
					\end{align}
					By estimate of $\lnenth$ in \eqref{ES:LinftyFluid2.1}, relations between $\lnenth,n$ and $\density+\pressure$ in \eqref{DE:enthalpy}, \eqref{DE:logenthalpy}, we have
					\begin{align}\label{EQ:2ndTrapped4}
						\int_{\hypersur_{\uts}\cap\{\ubar\in[0,\delta+2(\uts-\uinf)\delta]\}}(\density+\pressure)v_4^2\diff\spherevol\diff \ubar \geq &~\frac{3}{4} \int_{\hypersur_{\uinf}\cap\{\ubar\in[0,\delta]\}}(\density+\pressure)v_4^2\diff\spherevol\diff \ubar,&
						\text{for all }&(\sphere^1,\sphere^2)\in\Stwo.
					\end{align}
					Using the expression of $\FAlphaAB$ in Proposition \ref{PR:Ricci0}, since the difference of $\FAlphaAB$ on $\suubar$ is negligible (indeed, $\abs{\FAlphaAB(u,\ubar,\sphere^1,\sphere^2)-\FAlphaAB(u,\ubar,(\sphere')^1,(\sphere')^2)}\leq \BAfluid \abs{\uinf}$ for all $u,\ubar,\sphere,\sphere'$), \eqref{EQ:2ndTrapped4} implies the following:
					\begin{align}\label{EQ:2ndTrapped5}
						\abs{\uts}^2\int_0^{\delta+2(\uts-\uinf)\delta}\FAlphaAB(\uts,\ubarprime,\sphere^1,\sphere^2)\diff \ubarprime\geq \frac{1}{2}	\abs{\uinf}^2\int_0^{\delta}\FAlphaAB(\uinf,\ubarprime,\sphere^1,\sphere^2)\diff \ubarprime.
					\end{align}
					Combining \eqref{EQ:2ndTrapped5} with the assumption on the initial data \eqref{AS:LargenessofData2}, we have:
					\begin{align}
						\int_0^{\delta+2(\uts-\uinf)\delta}2\FAlphaAB(\uts,\ubarprime,\sphere^1,\sphere^2)\diff \ubarprime>   \frac{4}{\abs{\uts}}.
					\end{align}
				\end{proof}
				\subsubsection{Proof of the second case in Theorem \ref{TH:Main}}
				$\gtr\uchi<0$ in $\region$ is directly given by the estimates \eqref{ES:LinftyRicciCoefficients5}. By data prescription in Section \ref{SS:CharData} and assumptions \eqref{AS:SmallonHbar}, \eqref{AS:LargenessofData}, $\gtr\upchi>0$ on $\hypersur_{\uinf}\cup\hypersurbar_{0}$. 
				
				Now, provided assumption \eqref{AS:LargenessofData2}, we show that $\gtr\upchi<0$ everywhere on $S_{\uts,(1+2(\uts-\uinf))\delta}$. Multiplying equation \eqref{EQ:e4upchi} by $\lapse^{-1}$, using \eqref{EQ:omegauomega}, we have:
				\begin{align}\label{PV:Main21} 
					e_4\left(\lapse^{-1}\gtr\upchi\right)=&~-\frac{\lapse^{-1}}{2}(\gtr\upchi)^2-\lapse^{-1}\sgabs{\hchi}^2-2\lapse^{-1}\FAlphaAB.
				\end{align}
				Using coordinates expression \eqref{DE:Coordinatese3e4}, \eqref{PV:Main1} reads:
				\begin{align}\label{PV:Main22}
					\frac{\p}{\p\ubar}\left(\lapse^{-1}\gtr\upchi\right)=&~-\frac{1}{2}(\gtr\upchi)^2-\sgabs{\hchi}^2-2\FAlphaAB.
				\end{align}
				Integrating \eqref{PV:Main22} over $\ubar$, by assumption on $\gtr\upchi$ on $\hypersurbar_{0}$ as in \eqref{AS:SmallonHbar}, Lemma \ref{LE:improvedlapse}, and Lemma \ref{LE:EstimatesforuFAlpha}, for any $(\sphere^1,\sphere^2)\in\Stwo$, we derive
				\begin{align}
					\lapse^{-1}\gtr\upchi(\uts,\delta+2(\uts-\uinf)\delta,\sphere^1,\sphere^2)&\leq \lapse^{-1}\gtr\upchi(\uts,0,\sphere^1,\sphere^2)-\int_{0}^{\delta+2(\uts-\uinf)\delta}\left(\sgabs{\hchi}^2+2\FAlphaAB\right)(\uts,\ubarprime,\sphere^1,\sphere^2)\diff\ubarprime\\
					\notag&< \frac{3}{\abs{\uts}}-\frac{4}{\abs{\uts}}+\mathcal{O}(\delta a^{\frac{1}{2}})
					\\
					\notag&<0.
				\end{align}
				
				\appendix
				\section*{Appendices}
				\addcontentsline{toc}{section}{Appendices}
				\renewcommand{\thesubsection}{\Alph{subsection}}
				\subsection{Scaling of the Null Vectorfields}
				For convenience, we list our null vectorfields with different scaling as follows:
				\begin{subequations}
					\begin{align}
						\Lgeo^\alpha=&~-(\gfour^{-1})^{\alpha\beta}\p_\beta u,&
						\uLgeo^\alpha=&~-(\gfour^{-1})^{\alpha\beta}\p_\beta \ubar,\\	
						e_4^\alpha=&~-\lapse(\gfour^{-1})^{\alpha\beta}\p_\beta u,&
						e_3^\alpha=&~-\lapse(\gfour^{-1})^{\alpha\beta}\p_\beta \ubar,\\
						\Lunit^\alpha=&~-\lapse^2(\gfour^{-1})^{\alpha\beta}\p_\beta u,&
						\uLunit^\alpha=&~-\lapse^2(\gfour^{-1})^{\alpha\beta}\p_\beta \ubar.
					\end{align}
				\end{subequations}
				\subsection{Notations}\label{AP:notations}
				In this appendix, we gather notations that we use throughout the article.
				\begin{center}

					\begin{tabular}{|c|c|} 
						\hline
						Symbol & Reference  \\ 
						\hline\hline
						$\lesssim,\p,\pfour,\antisymmetic_{\alpha\beta\gamma\delta}$ & Section \ref{SS:Notations} \\
						\hline
						$\hypersurbar_{\ubar},\gfour,\Dfour,\Riem_{\alpha\beta\gamma\delta},\Ricfour{\alpha}{\beta},\Scalar$ & Section \ref{SS:despacetime} \\
						\hline
						$\sphereproject,\angnabla,\angD$ & Section \ref{SS:despacetime} \\
						\hline
						$v^\alpha,\density,\pressure,n,\temp,\enth,\ent,\speed,\speed_0$ & Section \ref{SSS:thebasicfluidvariables} \\
						\hline
						$\vort,\lnenth,\average{H}$ & Section \ref{SSS:v-ortho} \\
						\hline
						$T_{\mu\nu}$ & Section \ref{SS:EinsteinEuler} \\
						\hline
						$u,\ubar,\sphere^1,\sphere^2,\uptau,\rgeo$ & Section \ref{SS:NullFrame} \\
						\hline
						$\Lgeo,\uLgeo,e_4,e_3,e_A,\Lunit,\uLunit$ & Section \ref{SS:NullFrame} \\
						\hline
						$\lapse,\gsphere, \hypersur,\hypersurbar, v_I$ & Section \ref{SS:NullFrame} \\
						\hline
						$\tensorproduct,\wedge$ & Section \ref{SSS:basicnotations} \\
						\hline
						$\deform{V}$ & Definition \ref{DE:DeformationTensor} \\
						\hline
						$\upchi,\uchi,\upeta,\ueta,\omega,\uomega,\upzeta$ & Definition \ref{D:Defconnectioncoeff}\\
						\hline 
						$\gtr\upchi,\hchi,\gtr\uchi,\huchi$ & Definition \ref{D:Defconnectioncoeff}\\
						\hline 
						$\Weyl_{\alpha\beta\gamma\delta}$ & Definition \ref{DE:Weyl} \\
						\hline
						$\HodgeLeftDual\xi$ & Definition \ref{DE:HodgeDual} \\
						\hline
						$\AlphaAB,\Alphabar,\BetaA,\Betabar,\Rho,\Rhodual,\GaussCurvature$ & Definition \ref{DE:CurvatureComponents} \\
						\hline
						$\FAlphaAB,\FAlphabar,\FBetaA,\FBetabar,\FRho,\FRhodual$ & Definition \ref{DE:CurvatureComponents}\\
						\hline
						$\DWeyl,\DWeylStar,\DWeylBianchi,\DWeylBianchiStar$ & Definition \ref{DE:DerivativeofWeyl} \\
						\hline
						$\angRicfour$ & Definition \ref{DE:ProjectedRicci} \\
						\hline
						$\CAlphaAB,\CAlphabar,\CBetaA,\CBetabar,\CRho,\CRhodual$ & Definition \ref{DE:Renormalization} \\
						\hline
						
					\end{tabular}
					\begin{tabular}{|c|c|} 
						\hline
						Symbol & Reference  \\ 
						\hline\hline	
						$\vproject$ & Section \ref{SS:EulerEquations}\\
						\hline
						$\psi,\upsi,\upRiem,\uRiem,\variables,\uvariable$ & Section \ref{SS:SchematicNotation} \\	
						\hline
						$\Dfour^{\vec{m}},\angD^{\vec{m}},\WDfour^{\vec{m}},\WangD^{\vec{m}}$ & Definition \ref{DE:operators}\\
						\hline
						$N,a,b,\delta$ & Section \ref{SS:Norms}\\
						\hline
						$\mathcal{R}_{i,j,k,\infty},\mathcal{R}_{i,j,k,2},\mathcal{R}$ & Section \ref{SS:Norms}\\
						\hline
						$\mathcal{O}_{i,j,k,\infty},\mathcal{O}_{i,j,k,2},\tilde{\mathcal{O}}_{N,2}$ & Section \ref{SS:Norms}\\
						\hline
						$\mathcal{E}_{i,j,k,\infty},\mathcal{E}_{i,j,k,2},\tilde{\mathcal{E}}_{N,2}$ & Section \ref{SS:Norms}\\
						\hline
						$\iniBAsize,\iniBAfluid$ & Section \ref{SS:InitialData}\\
						\hline
						$\BAsize,\BAfluid$ & Section \ref{SS:BA}\\
						\hline
						$\Lie$ & Lemma \ref{LE:LieandCovariant}\\
						\hline
						$\varv,\varh,\vars,\varJ$ & Definition \ref{DE:EnergyCurrent}\\
						\hline
						$\mathfrak{F},\mathfrak{G}^\alpha,\mathfrak{H}$ & Section \ref{SS:CommutedEulerEq} \\
						\hline
						$\omega^{\dagger},\kappa,\uomega^{\dagger},\underline{\kappa},\mass,\umass$ & Section \ref{SS:provetopRicci} \\
						\hline
						
					\end{tabular}
				\end{center}
				\bibliography{bib.bib}
				\bibliographystyle{plain} 
			\end{document}